# Justice Through the Lens of Calculus

## FRAMING NEW POSSIBILITIES FOR DIVERSITY, EQUITY, AND INCLUSION

# Introduction to the Editorial Team

### Matthew Voigt
*Clemson University*

Matthew Voigt (He, Him, His) is an Assistant Professor in Engineering and Science Education at Clemson University. I grew up in a rural farming community in central Minnesota, and I am a proud first-generation college student and Queer STEM education researcher. My disciplinary background is in applied mathematics. I earned my Ph.D. in Mathematics and Science Education from San Diego State University and the University of California San Diego. I enjoy traveling, cooking, playing strategy board games, and hiking in National Parks.

### Jess Ellis Hagman
*Colorado State University – Fort Collins*

Jess Ellis Hagman (She, Her, Hers) is an Associate Professor in the Department of Mathematics at Colorado State University in Fort Collins. My area of research is undergraduate mathematics education, and my work is driven by the goal of dramatically increasing the number and diversity of people who thrive and succeed in undergraduate mathematics - especially introductory mathematics courses that often function as a barrier for STEM intending students. The ideas shared in this volume have dramatically informed my approach to research and practice as a member of a mathematics department, and I am grateful to all contributors. In addition to research, mentoring, and teaching, I enjoy spending time with my family in nature as much as possible.

### Jessica Gehrtz
*University of Texas – San Antonio*

Jessica Gehrtz (She, Her, Hers) is an Assistant Professor in the Department of Mathematics at the University of Texas at San Antonio. After teaching middle school and high school mathematics for several years, I returned to school and earned my M.S. in Mathematics and my Ph.D. in Mathematics Education from Colorado State University. My research is focused on key stakeholders – students, instructors, and departments – in the adoption and sustained use of evidence-based teaching practices at the collegiate level. This volume, and the collaboration with authors and editors, has impacted my personal reflection, teaching, and research, and for that, I am grateful. I enjoy traveling, reading, hiking, and spending time outdoors.

### Brea Ratliff
*Auburn University*

Brea Ratliff is the immediate Past-President of the Benjamin Banneker Association, Inc. (BBA), and has taught mathematics at the elementary, middle, high school, and collegiate levels. I am the founder and Educational Strategist of Me to the Power of Three, a consulting company which specializes in developing inquiry-based programs and resources to support the teaching and learning of mathematics at all levels. I am currently pursuing a Ph.D. in secondary mathematics education at Auburn University, where I have been recognized as a 2021 recipient of the Presidential Graduate Opportunity Program fellowship. My research interests include equity in mathematics education, teacher education and leadership, and community/youth development.

### Nathan Alexander
*Morehouse College*

Nathan Alexander is an Assistant Professor in the Department of Mathematics at Morehouse College and faculty director of the Quantitative Histories Workshop. My work integrates mathematical and statistical modeling with the computational social sciences to speak to critical traditions in education that seek to uncover and dismantle historical systems of power and privilege. I have worked as a teacher for more than 15 years in secondary, two-year, four-year, and alternative educational institutions, such as prison systems and community-based organizations. I received my Ph.D. from Columbia University in mathematics and education where I supported youth educational programs in mathematics in Harlem and Bedford Stuyvesant, and I double majored in mathematics and sociology as a Pogue Scholar at the University of North Carolina at Chapel Hill. I enjoy backpacking trips to Central and South America, cooking, and playing tennis.

**Rachel Levy**

*North Carolina State University*

Rachel Levy is Executive Director of North Carolina State's Data Science Academy and Professor of Mathematics. I was the 2020-21 American Mathematical Society Congressional Fellow and served as Deputy Executive Director of the Mathematical Association of America from 2018-2020. I served as Vice President for Education for the Society for Industrial and Applied Mathematics and was Professor of Mathematics and Associate Dean for Faculty Development at Harvey Mudd College. I have taught mathematics and English at the middle school, high school and undergraduate levels and mathematics and mathematics education at the graduate level. I would like to thank my family, students, colleagues, and communities for all the ways they help me grow. I love mountain biking in the forest and boating in the ocean.

# Table of Contents

# *Part 1: Introduction and Framing*

# Building an Evolving Framework: A Clarion Call / Manifesto

**Gregory V. Larnell**
*University of Illinois - Chicago*

A volume such as this one represents both an abiding concern and an urgent call to the field: No longer is it possible (if ever, at least not in good conscience) to look beyond longstanding, abysmal, systemic patterns of exclusion attributable to the curricular trajectory to and through calculus—patterns that stubbornly, chronically, yet persistently adhere and contribute to identity-based disparities within the broader society. The call this time, however, is even broader and includes the following corollary: *The responsibility for redressing these deep inequities should not be relegated to a small, impassioned group of reformers but shared by the entire community of faculty members and academic researchers who steward the undergraduate mathematics education corridor.* Simply put, the responsibility for redressing injustice belongs to all of us.

As Kilpatrick (2011) argues, we need more research and more researchers within this community toward "seeing our work whole" (p. 173)—that is, toward seeing the whole of undergraduate mathematics education collectively and fully serve all learners and, ultimately, toward contributing to the advancement of our society. This volume represents considerable, yet ultimately incremental, progress toward that goal.

This volume also represents an invitation—an invitation to engage. No matter what your particular area of focus may be, it is necessary to consider how your work can connect to and expand the imperatives of greater diversity, equity, and inclusion. The authors of this volume are working toward this goal and offer all of us evidence that this work is not just possible but vital and crucial.

Justice through the lens of calculus as an evolving framework seeks to unify a vision for calculus programs that are diverse, equitable, and inclusive while critically examining how they are situated within current power structures and levels of enactment. Drawing on the joint vision from TODOS: Mathematics for ALL[1] and the National Council of Supervisors of Mathematics (NCSM), "A social justice stance requires a systemic approach that includes fair and equitable teaching practices, high expectations for all students, access to rich, rigorous, and relevant mathematics, and strong family/community relationships to promote positive mathematics learning and achievement. Equally important, a social justice stance interrogates and challenges the roles power, privilege, and oppression play in the current unjust system of mathematics education—and in society as a whole" (todos-math.org). Toward that end, our collective work must expand and include fervent attention to how our students and colleagues are experiencing not only mathematical content and processes writ large but also critically analyzing longstanding legacies of exclusion, tracking, and gatekeeping practices; status-laden hierarchies; and a variety of implicit practices that all-too-often inveigle individuals to shed parts of their identities in order to belong.

---

[1] TODOS: Mathematics for ALL is an international professional organization that advocates for equity and excellence in mathematics education for ALL students - in particular, Latina/o students

# Introduction to this Volume

**Matthew Voigt, Rachel Levy, Jess Ellis Hagman, Jessica Gehrtz, Brea Ratliff, Nathan Alexander**

## 1. *Preface*

The idea for creating this volume emerged from an NSF-funded research project entitled "Progress through Calculus" (DUE I-USE #1430540). The focus of the Progress through Calculus (PtC) project was to examine the Precalculus through Calculus II course sequence and associated supports in United States (U.S.) colleges and universities. The project consisted of two phases of data collection and analysis. In the first phase, a national census survey was sent to mathematics departments offering a graduate degree (master's or Ph.D.) in mathematics. The survey gathered information to better understand the characteristics of successful calculus programs. In the second phase of the project, in-depth longitudinal case studies were conducted at 12 colleges and universities. The case studies investigated models of the Precalculus through Calculus II sequence, their implementations, and their impact on student outcomes. More details about the entire project are located at: https://maa.org/ptc

As a result of this research, Editors Hagman, Voigt, and Gehrtz formed a thematic research team examining issues of diversity, equity, and inclusion (DEI) across the national census survey data and the 12 case study universities. It became clear that while many of the members of mathematics departments valued issues of DEI, most did not yet have actionable ideas or strategies for addressing these ideas locally within their departments. Additionally, many departments pointed towards broader university-wide programs for addressing issues of DEI, with only a select few having local initiatives within the purview of the mathematics department.

*As such, the Editors saw a clear need to gather a multitude of ideas, works in progress, and creative solutions to systemically and centrally address topics of diversity, equity, and inclusion in mathematics programs.*

Through generous support from the National Science Foundation, we were able to address this need through the development and creation of a resource attending to DEI in Calculus programs that could be made available to the broader mathematics community. This volume began with a call for individual case studies (see figure 1) from math departments attending to DEI issues, which became the 30 case studies of this volume. The case study submissions were generally written by mathematics faculty engaged in teaching and administering the Calculus programs.

Upon analysis of the case studies, we identified salient topics within the submissions and invited Mathematics Educators to author cross-cutting thematic chapters. We asked this group of thematic authors to read the case studies, situate them in extant literature, and create a future vision for a more diverse, equitable, inclusive, and justice-oriented field of mathematics. The thematic authors also pose questions to the readers of this volume to allow for further exploration and insights into local contexts. In addition to the cross-cutting thematic authors, we recruited an author team to bring in student voices and another author team to share relevant data to problematize how we measure DEI efforts and to situate this volume in its historical period. The authors and Editors met several times to choose and discuss the cross-cutting theme chapters of interest that were relevant to the ideas discussed in the case studies.

**Scope:** The MAA Notes Volume on Diverse Equitable and Inclusive (DEI) Issues in Calculus Programs encourages a broad array of submissions that highlight issues of DEI in introductory mathematics programs with special attention to precalculus, differential calculus, and integral calculus and surrounding departmental programs to support students in these courses. For this volume we are soliciting "illustrative case studies" that showcase ways in which departments and instructors are attending to promoting diverse introductory mathematics programs, achieving, or monitoring equitable student outcomes and experiences, and promoting inclusive teaching practices. Achieving and promoting DEI issues in introductory mathematics programs is not an easy undertaking so we encourage submissions of models in progress, discussions of potential obstacles, challenges, and what departments/instructors have done to overcome barriers to address these issues. We also encourage collaborations between mathematics department members and people outside the department involved in programs, such as individuals in administrative positions or working with student support centers. Submissions might address (but are not limited to) the following topics: a) How coordination can support fairness and also justness b) How placement procedures can value multiple ways to demonstrate readiness c) How professional development (for faculty, instructors, and GTAs) can address DEI d) How changes to the curriculum respond to an increasingly diverse student population e) How departments/universities collect and use data to inform changes related to DEI f) How centers or programs operate to support inclusion and student success in STEM.

*Figure 1. Call for Case Studies Included in this volume*

We want to acknowledge our own struggles as we curated a collection of voices, the concurrent privilege and burden of doing this work, and all the power dynamics that can come into play. Our goal was mindful, respectful, and collaborative work. At the same time, during the development of this volume, we were challenged to think about our positionality and how to ensure the presence of diverse voices within this volume. Even with the best intentions of promoting diversity, equity, and inclusion, our push to develop the volume needed to be balanced and paused at times to ensure that multiple perspectives were being given space. As such, we began a re-envisioning of our project. We altered authorship and editorial teams, included student voices in authorship, and explicitly reframed and changed how we presented the case studies and chapters in this book to promote an anti-deficit framing.

Our work on this volume occurred during some extraordinary moments of national focus on race, violence, disinformation, and the disruption caused by a global pandemic. Many of us experienced personal trauma and loss while we were in the process of creating this work. The Editors feel grateful that this community has been a source of hope, support, and friendship. We hope that you will find ideas, solace, and discomfort in this book as you engage on a journey of supporting issues of diversity, equity, and inclusion in your own context.

## 2. *Acknowledgments and Dedication*

The Editors would like to acknowledge support from the National Science Foundation (DUE-1432381), the Mathematical Association of America, and the American Mathematical Society. In addition, we would like to thank the MAA Notes Volume Editor and Reviewers, the Progress through Calculus research team, and the case study and chapter author teams. We would also like to thank each other as well as our families, friends, and colleagues who have walked this road with us. A special thank you to Kiera Edwards for her expertise in preparing the grant supplement, Victoria Barron for her support in updating references and citations, and Destinee Cooper for her assistance in formatting and reviewing the volume.

We dedicate this volume to anyone who has received the message that mathematics was not for them. We dedicate this volume to anyone who felt they did not belong in mathematics. And finally, we dedicate this volume to you, the reader, for taking the time and energy to engage with these issues.

### 3. *How to Use this Book*

Our goal was for this book to be of value to many different stakeholders, including faculty, administrators, evaluators, students, and mentors. The book is structured into three main parts: *Introduction and framing*, *Cross-cutting thematic*, and *Illustrative case studies*.

The *Introduction and framing* part contains general information about the book structure, a guiding manifesto, and cross-cutting inquiries to frame your exploration of DEI issues. This part also contains Chapter 1 (A Critical Examination of Undergraduate Degree Completion Data Problematizing Institutional Diversity Measures in STEM), which uses national data to provoke a critical examination of undergraduate degrees in STEM and how we assess and think about issues related to DEI. This chapter also helps to contextualizes this volume within the current historical period.

The *Cross-cutting thematic* part contains six chapters addressing various topics as they arise within the case studies and one chapter that leverages student voices to reflect on the volume. Chapter 2 (Developing the Next Generation of Change Agents in College Mathematics Instruction) broadly examines professional development through the lens of change agents of instruction. Chapter 3 (Three Frameworks for Advanced and Theoretically-Driven Analyses of Classroom Practices) examines classroom practices in the case studies by drawing on three frameworks from mathematics education. Chapter 4 (The Role of Role Models in Addressing Diversity, Equity, and Inclusion) examines how role models are leveraged in the case study work across departments. Chapter 5 (Course Redesign: Pathways Towards Transformation) explores course pathways and course redesign efforts in calculus programs to support DEI. Chapter 6 (Leveraging Identity and Language to Promote DEI in Undergraduate Mathematics) investigates the role of language and identity in promoting DEI. Chapter 7 (Reflecting on the Role of Assessment in Placement and Evaluating Student Learning in Calculus)Looks at the ways alternative assessments are used to promote DEI. The culminating Chapter 8 (Reflections From a Student Advisory Committee to Deepen Equity in STEM Classrooms) elevates the voices of marginalized STEM students who react to the thematic chapters and paint their vision for a more diverse, equitable, and inclusive field of mathematics.

Please note that the cross-cutting thematic chapters represent several different levels of scholarly inquiry. Some may require the reading of reference material to fully engage the concepts. Some refer to many or all of the case studies, while others are more focused on a few ideas. Our hope is that the ideas in this volume will help readers wrestle with and discuss their own practices with colleagues.

The *illustrative case study* part contains 30 case studies from authors across the United States that address local issues of diversity, equity, and inclusion. The case studies are not presented in particular order or grouping as they address many related topics. Here, we also wish to acknowledge an important limitation of the volume: cases in this section come from the United States perspective, despite issues of equity being global in both form and nature.

As Editors, we made a conscious decision to provide a glossary of terms, but not to standardize language across the volume. This decision enabled authors to select choices that were most suited to their needs and communities. There may be disagreements and tensions with what language is used in different contexts and settings, but this provides an authentic lens into the equity work that is occurring across a variety of institutions.

As you read, you may choose to start with some of the case studies and then turn to a broader thematic chapter. Alternatively, you may choose to read the thematic chapters and then dive into relevant case studies. Within the thematic chapters, case studies are referred to in APA style format with Author names and then (CS #) – for CASE STUDY #. The numeric case study numeral is hyperlinked and will aid navigation between the theme chapters and case studies.

The following table presents an overview of which case studies are specifically mentioned in each of the thematic chapters.

*Table 1*
*Tabular organizer of the illustrative case studies and cross-cutting thematic chapters.*

| | | | **Cross-Cutting Thematic Chapters** | | | | | | |
|---|---|---|---|---|---|---|---|---|---|
| **Case Studies** | **#** | **Case Study Title** | 2- Professional Development | 3- Classroom practices | 4- Role models | 5- Course redesign | 6 – Language and Identity | 7- Assessment | 8– Student Voices |
| | 1 | Equity in Grading Systems: Moving away from "fair" towards transparency and inclusion in coordinated Calculus courses | | X | | X | | X | X |
| | 2 | Using a Multi-Tiered and Data-Driven Approach to Support Student Access and Success in Introductory Mathematics at a Large, Hispanic-Serving Institution | | | | X | | X | |
| | 3 | Increasing Inclusion in Large Enrollment, Uniform Math Courses: Instructor Training and Course Assessment | X | | | X | | X | X |
| | 4 | Promoting Racial Equity in the STEM Math Pathway in Community Colleges | | | | X | X | X | X |
| | 5 | Stretching Calculus: A Yearlong Calculus Class that Stretches Minds | | X | X | X | | X | |
| | 6 | Redefining success: A multifaceted approach to assessing departmental change in first year mathematics/statistics | | X | | X | | X | X |
| | 7 | Addressing Diverse Equitable and Inclusive Issues in a Four-year Private Hispanic Serving Institution | | X | X | X | X | | |
| | 8 | Supporting Underrepresented Minority Students in STEM Through In-Class Peer Tutoring | | X | X | | X | | X |
| | 9 | A Case Study of an Emerging Scholars Program at a Predominately White Institution: Analyzing the Success of Underrepresented Students in an Inquiry-Based Calculus Sequence | | X | X | X | X | | |
| | 10 | Towards student-ready mathematics departments: Creating mathematics placement messaging within an equity-minded framework | | X | | | X | X | |
| | 11 | Instructional Practices that Foster Success in Underrepresented Populations | | X | X | X | | | |
| | 12 | Creating Space for Student Agency to Support Success in Online Pre-Calculus | | X | | X | | | X |
| | 13 | Coordinating Courses Below Calculus: Creating a Community of Practice | X | | | X | | X | |
| | 14 | A College Algebra Intervention to Address Power Dynamics and Promote Students' Reasoning | X | X | | X | | | |
| | 15 | Collaboratively Re-envisioning Calculus for the Modern Student | | X | | X | | | |
| | 16 | Peer-Led Team Learning in Mathematics: An Effort to Address Diversity and Inclusion Through Learning and Leadership | X | X | X | | X | | |
| | 17 | Creating Humane Exam Structures in College Calculus Courses | | X | | | X | X | X |
| | 18 | Supporting the STEM Pathway at Appalachian | | X | | X | | X | |
| | 19 | "I Wish I Could Say 'You Should Not be Here'": An Analysis of Instructors' and Students' Contrasting Perceptions of a Racialized and Gendered Gatekeeping Practice in Calculus | | X | | | X | | X |
| | 20 | Enhancing Academic Capital to Improve Equity and Inclusion in First Semester Calculus | | X | X | X | | X | |
| | 21 | A Multi-pronged Approach to Closing Opportunity Gaps in Calculus I | X | X | | X | | | |
| | 22 | A Focus on DEI for Students Yields DEI for Instructors | X | X | | X | X | | |
| | 23 | Undermining the "Bad at Math" Identity: Overhauling Differential Calculus to Empower Marginalized Students | | X | | X | X | X | X |
| | 24 | The Effective Thinking Calculus Project at UT Austin | | X | | X | | X | |
| | 25 | Leveraging Classroom Data to Promote Equitable Instruction | X | X | | | X | | |
| | 26 | CALCULUS: Crossing the Bridge to Success in STEM | | | | X | | | |
| | 27 | Shifting Pre-Calculus from Gatekeeper to Gateway at a Hispanic-Serving Institution | X | | X | | X | | |
| | 28 | Nivelando el Campo Educativo / Leveling the Educational Field | | | X | X | X | | |
| | 29 | First-generation students' experience of entering a STEM major | | X | X | | X | | X |
| | 30 | Inclusive Calculus: Active, Applied, and Inquiry-based Learning with Linked Workshop and First Year Seminar | | X | X | X | | | |

## Cross-cutting inquiries: Moving Towards a Systems-critical Practice

**Nathan Alexander, Brea Ratliff, Rachel Levy, Matthew Voigt, Jess Ellis Hagman, Jessica Gehrtz**

Building a framework around justice at the level of one's daily practices requires not only an evolving set of actions but also an ever-evolving and explicit intentionality toward diversity, equity, and inclusion. Each effort, diverse in focus and goal, is likely to cut across various levels of practice and result in different levels of influence, power, and privilege. The dynamic and interactive nature of these influences led us (the editors) to identify cross-cutting inquiries that help guide readers in their reflections. We have identified a set of leading inquiries that cut across the volume's themes and various points of entry that any particular submission may focus on. These entry points are 1) individual and self-work, 2) classroom-level work, 3) departmental- or division-level work, 4) university-level or cross-institutional work, 5) broader policy and research efforts, and 6) mathematical content areas. We see each point of entry, or "level," as neither clean-cut nor confined and likely to overlap in multiple ways. More broadly, however, the leading inquiries can support one's practical transition to develop a critical perspective that challenges the status quo. Or, more simply, a shift in the lens from a focus on individual experiences towards an acknowledgment of the interconnected and historical systems of unequal power and privilege known as systemic inequalities.

***Leading inquiry #1 (individual-level and self-work): How do I engage in the work of being critically reflexive to support justice-making?***

Change requires efforts across multiple levels. However, many would argue that individual and "self" work is the most important component of generating change, especially when that change is aligned with one's position or power. Building on this idea, offered by James Baldwin, who challenged us to reconsider our own beliefs and identities in the wake of injustice, change requires deep and critical development at the individual level. We hope that this leading inquiry provides opportunities for everyone to reflect on how shifts in individual beliefs, behaviors, and actions toward supporting diversity, equity, and inclusion have come as a result of critical reflexivity. Moreover, how have individual reflections surfaced new pathways toward justice in the local and global contexts? Additional questions to consider:

- What data are available to guide my future efforts, and what data need to be collected?
- How can I help people to move forward from where they are?
- What feedback mechanisms exist so that DEI becomes a process inherent to my professional growth?
- What helps people take ownership and responsibility for the need to make changes?

***Leading inquiry #2 (classroom-level work): How do I engage in the work of supporting justice-making and change in my classroom?***

Mathematics instructors have an opportunity to leverage their power and influence in the classroom to promote justice. Au, Bigelow, and Karp (2007) note that "classrooms can be places of hope, where students and teachers gain glimpses of the kind of society we could live in and where students learn academic and critical skills needed to make it a reality." It is important to note that justice in the classroom is not achieved in a singular way or through a specific means of instruction. Rather, it is an ongoing journey that is started by acknowledging the injustice within the system and then working to build a safe and encouraging space for every student's voice to be valued and respected (even when that does not sound like the instructor's voice or the discipline's voice). How can instructors be critically reflexive to interrogate the structures, instructional materials, and teaching practices toward supporting diversity, equity, and inclusion? Additional questions to consider:

- How do/can/might classroom-based efforts move us toward broader visions of justice in calculus?
- What are the priorities for my classroom?
- What data are available to guide our future efforts, and what data need to be collected?
- How can we help people to move forward from where they are?
- What feedback mechanisms exist so that DEI becomes a process of continuous improvement within my classroom?

***Leading inquiry #3 (department- and division-level work): How are mathematics departments and divisions positioned to move us toward change in support of multiple visions of justice-making?***

Mathematics departments represent a community of practitioners that embody and represent the field of mathematics. As such, they have inherent power in supporting and challenging normative practices within the field. At the same time, departments are not monolithic entities and are composed of individuals with different identities, beliefs, positionality, and visions for mathematics education. How can departments (or department members) be critically reflexive to interrogate the structures and programs toward supporting diversity, equity, and inclusion? Additional questions to consider:

- How can we engage ourselves and our colleagues within the department in a continuous state of growth and development?
- What are the priorities for our department/introductory mathematics program?
- What data are available to guide our future efforts, and what data need to be collected?
- How can we help departments move forward from where they are currently?
- What feedback mechanisms exist within the department so that DEI becomes a process of continuous improvement?

***Leading inquiry #4 (university-level or cross-institutional and community-level work): How are institutions positioned to move us toward change in support of multiple visions of justice-making as it relates to introductory mathematics programs?***

Institutions represent a community of diverse positions, visions and objectives, beliefs, positionality, and identities working together to support undergraduate (and graduate) education. While only a subset of this community is focused on mathematics education, mathematics is integral (pun intended) to every degree and therefore carries power as a field. How can institutions (and the relevant members of the institutions) be critically reflexive to interrogate the structure and programs toward supporting diversity, equity, and inclusion within introductory mathematics programs? Additional questions to consider:

- How can we engage ourselves and our colleagues across the institution in a continuous state of growth and development?

- What are the priorities for our institution?
- What data are available to guide our future efforts, and what data need to be collected?
- How can we help entities and the people within them to move forward from where they are?
- How can we understand not just what happened in the case studies, but also WHY the various institutions made the choices that they did?
- What feedback mechanisms exist so that DEI becomes a process of continuous improvement within our institution?

***Leading Inquiry #5 (broader visions, policy-level efforts): What broader vision statements and policy-level efforts should be enacted to support systemic justice efforts in mathematics?***
The long-term success of implementing these changes at any level requires alignment with local, state, and national policies to increase students' opportunities in postsecondary mathematics and improve continuity across institutions/education systems. To ensure the sustainability of these efforts, a broader vision statement of justice-oriented mathematics should be developed by an inclusive and diverse audience. How can institutions engage in equitable communication and collaboration with students, community members, legislators, and other key stakeholders to establish new structures and practices promoting the immediate need for justice-oriented policies that directly impact postsecondary mathematics? Additional questions to consider:

- What are possible justice-oriented, collective, and long-term visions for DEI in calculus, specifically, and mathematics, more broadly, as it relates to research and public policy?
- What is the work needed to push DEI towards the type of change that we seek in mathematics education?
- In an age of injustice, how can DEI efforts move us toward structural and systemic changes in mathematics?

***Leading Inquiry #6 (content, curriculum): How can choices of course topics and content support change in support of multiple visions of justice-making?***
The choice of mathematical content, how courses are sequenced, and what forms of knowledge are valued within mathematical spaces require a critical lens to support the inclusion of all learners. How can we re-envision the selection and treatment of topics included in a calculus course so that the course engages students in issues related to DEI?

- Who are presented as the doers of mathematics? How might the work of contemporary mathematicians be included in the course to introduce students to a diverse group of mathematicians and mathematical approaches?
- Are there ways for students to influence the pace, topics, and applications of the course?
- What will be the role of computing in the course?
- Will students engage in calculus using the tools they will use in their life and jobs?
- What interdisciplinary opportunities might be possible that engage students in topics of interest and relevance to them?

As Editors, we strongly encourage readers to revisit these inquiries again, and then again, as you consider the various components and your framing of the work taking place in this volume

# 1 A Critical Examination of Undergraduate Degree Completion Data Problematizing Institutional Diversity Measures in STEM

**Naneh Apkarian**
*Arizona State University*

I, Dr. Naneh Apkarian, am currently an Assistant Professor at Arizona State University in the School of Mathematical and Statistical Sciences. I hold two degrees (BA & MA) in mathematics and a doctorate in mathematics education; my field is undergraduate mathematics education research. I identify as a woman and use she/her pronouns. Ethnically I identify as Armenian-American with cultural ties to the Middle East where much of my paternal family was born; on census forms I select "white" or "prefer not to disclose" when asked about my race. My father immigrated to the USA for college while my mother's family has lived in the USA for generations. I am not a first-generation college student; many of my family members hold advanced degrees and some are in academia. These experiences necessarily impact my perception and interpretation of the world; I acknowledge these and strive to navigate my research, practice, and lived experience responsibly.

As we set forth this compilation of efforts to create more equitable educational systems in undergraduate STEM, we look to the future. We aspire toward educational excellence that is inclusive – because if it is not, then we cannot call it excellence. It is our fervent hope that our words will seem passé sooner rather than later. We know that our field's understandings of systemic inequity, biased and exclusionary structures in education, and how to counteract these features of society and higher education will advance. In light of that, I offer up this chapter as a time capsule to contextualize what we, in the rest of the book, describe as current efforts in support of diversity, equity, and inclusion. In addition to a broad snapshot of the current situation regarding STEM degree completion, I describe the limitations of the dataset and consider some aspects of students' lived experience which are hidden from similar analyses.

Many of the efforts toward inclusion described in this volume are motivated, at least in part, by oft-repeated statistics about populations which are "underrepresented in STEM." This is an ill-defined term which means many different things to different people, and which obscures the cultural structures and barriers to entry and sustained participation in STEM; it can also result in a reduction of inclusion and equity to numerically representative diversity (Williams, 2017). The most common question and report in this vein is, "what proportion of people who […] belong to [each demographic group]." This formulation puts the activity (e.g., declare a major in STEM) or outcome (e.g., pass calculus) first, and the people and their experiences second. Furthermore, this framing can convey images of a finite pie – when one group increases their share, another's share must decrease. In an effort to focus our attention primarily on the *people* we are considering, and to minimize zero-sum interpretations, I instead ask,

> *What proportion of each (available) race-gender group of BA/BS degree completers obtain their degree in a STEM field?*

This chapter privileges intersectional race-gender identities and considers the distribution of each group's bachelor's degrees instead of facts about numeric representation or statistical diversity. I present a brief analysis of publicly available data from postsecondary institutions in the US, revealing different patterns in degree field across different groups of students. The race-gender categories available for this analysis are limited in scope, ignoring many identities and nuances which are not part of US federal race, ethnicity, or gender designations and preclude deep critical or intersectional analysis (Covarrubias & Vélez, 2013). While I offer up figures and some slight interpretation, the goals of this chapter are primarily to provoke *your* thoughts and make explicit some of the assumptions and blind spots which are frequently present in conversations related to diversity, equity, and inclusion in the context of postsecondary STEM. I hope that this chapter will serve to (a) contextualize the initiatives presented in other chapters in this volume; (b) allow for some quantitative assessment of future progress; and (c) spur you to assess your local context with a critical eye toward understanding students' pursuit and completion of STEM degrees.

## 1. *Understanding & Framing the Analysis*

Instead of a comprehensive analysis of extant literature, as is done in some other chapters in this volume, here I leverage an analysis of institutional postsecondary data to contextualize the volume *and* to provoke critical examination of how such data is collected and interpreted. I also discuss some of what is left out by this analytic approach and urge those on the ground to do better. The data presented later in this chapter comes from the Integrated Postsecondary Education Data System (IPEDS), a database of publicly available information from the U.S. Department of Education's Institute of Education Sciences (IES NCES, 2020). Specifically, I use data reported by 2057 not-for-profit, not entirely online, degree-granting postsecondary institutions about students who earned a baccalaureate (BA/BS) degree in the 2018-19 academic year, with indicators of race, gender, and degree field (see Appendix A for more detail about the institutions).

## 2. *Critical & Intersectional Approach to Interpretation*

The data I use in this chapter was *not* collected using a critical or intersectional lens, however to the extent possible I draw on aspects of Critical Race Quantitative Intersectionality (CRQI):

> *Critical race quantitative intersectionality is an explanatory framework and methodological approach that utilizes quantitative methods to account for the material impact of race and racism at its intersection with other forms of subordination and works toward identifying and challenging oppression at this intersection in hopes of achieving social justice for students of color, their families, and their communities (Covarrubias & Vélez, 2013, p. 276)*

In line with some of the principles of CRQI, in this chapter I:

- Use intersectional race-gender identity markers to focus on the experiences of different groups of people
- Do not assume that statistics about a group allow us to assume the trajectory of any individual
- Reject the idea of race and/or gender as a causal variable in considering differential patterns within the data
- Acknowledge that these blunt groupings homogenize heterogeneous populations and hide intra-group differences

- Ask questions aimed at (along with the rest of this volume) addressing injustice and seeking transformation of postsecondary STEM

### 3. *Which degrees, and why?*

As the bulk of this volume discusses undergraduate mathematics programs at four-year-or-above institutions, I distinguish between majors (i.e., degree fields) based primarily on their (probable) mathematical requirements. This allows me to present analysis about outcomes for students most impacted by the efforts discussed in the other chapters of this volume. There is no consistent agreement about what is and is not a STEM degree, but I have adopted three categories using Classification of Instructional Program (CIP) codes (more detail in Appendix B).

The first, *pSTEM*, includes engineering (CIP 14), mathematics and statistics (CIP 27), and physical sciences (CIP 40) degrees; most of the majors subsumed by these headings require a full year of single-variable calculus and some additional courses beyond (e.g., multivariable calculus, differential equations, linear algebra). The second, *sSTEM*, includes computer and information sciences (CIP 11) and biological and biomedical sciences (CIP 26); most degrees within these codes require multiple mathematics courses, but there is more variation in the number and nature of the course sequences expected (e.g., statistics, discrete mathematics). Other degrees, categorized as *non-STEM* here, generally have fewer and more variable requirements for mathematics courses as part of the major track. This distinction is helpful when considering which students are impacted most by efforts to improve particular course sequences, but it can also be linked to STEM hierarchies such as a broad cultural tendency to refer to the more mathematics-intensive disciplines as “hard” or more valuable than others; this strategy of considering mathematics as the gatekeeper of the sciences can be damaging.

### 4. *Demographic categorizations & commentary*

Identity is not something which can be easily defined, and when dealing with individuals it should not be reduced to a few generic terms. However, in considering patterns related to systemic inequity that is based on limited understandings of identity, we find ourselves in a bind. To reduce students to a group based on single words for “race” and/or “gender” is to deny their complex experience and to gloss over extensive variation in lived experience; to deny that these kinds of categorization are relevant for their lived experiences can obscure the systematic aspects of discrimination; and large national-scale data sets are collected using problematic and limited federal designations. Data at our disposal for this analysis uses the same categories as the US census, which is collected without a critical lens and (at a minimum) homogenizes information about heterogeneous groups (Covarrubias & Vélez, 2013). I report on patterns using available information about race/ethnicity, gender, and STEM degree completion – but I do not discuss these outcomes as inherent to a particular group, nor do I present these as useful for understanding the experiences, past present or future, of any individual.

With all these caveats, the race/ethnicity groupings I use from the IPEDS database are presented in Table 1, along with estimates of the proportion of the US population aged 18-24 who fall into each category. This is presented as contextual information about what we know about the traditional college-going age population of the US in the late 2010s.

| Race/Ethnicity Category | US Population (18-24) |
|---|---|
| American Indian and/or Alaska Native | 0.8% |
| Asian | 5% |
| Black and/or African American | 13% |
| Hispanic and/or Latinx | 15% |
| Native Hawaiian and/or Pacific Islander | 0.2% |
| White | 55% |
| Two-or-more | 4% |
| Other (unknown & non-resident) | 6% |

*Table 1: Descriptions of the US census race/ethnicity categories and estimates for the proportion of the traditionally college-age US population which fall into each category.*

There are myriad reasons why these "race/ethnicity" categories are problematic, limiting, and reductionist; experts have written extensively on the topics, but we only briefly mention a few which relate particularly to this analysis and discussion.

In the IPEDS data set, the category Hispanic and/or Latinx includes all students who report their ethnicity as Hispanic and/or Latinx, regardless of how they report their race, an example of quantitative homogenization of a heterogeneous group. The conglomeration of so many ethnicities and identities into "Asian" is another colossal challenge when trying to understand racialized stereotypes and barriers. When this data was collected, people were instructed to select "Asian" if they had origins in the original peoples of "the Far East, Southeast Asia, or the Indian subcontinent." This includes over 20 ethnicities and obscures differential access and participation in higher education (Nguyen et al., 2015). Finally, students were asked to self-identify their race using limited categories, which means that analyses of the data will "miss discrimination based on racial appearance by relying on data that instead capture racial self-identification, a different aspect of race that often does not correspond" (Roth, 2010, p. 1288). This racial mismatch is generally attributable to varied (and evolving) conceptualizations of race which can vary substantially across cultures, contexts, and time (Brown, 2020; Gupta, 2020; McKenney & Bennett, 1994; Roth, 2010). Because, of course, race is a social construct.

The social construct of gender is another aspect of identity which is relevant to a student's experience. Today, the U.S. census continues to ask about biological sex and provide the binary options "female" and "male" on their surveys; the IPEDS database retains the binary but uses the gender terms "man" and "woman." Projections from census data and other surveys and polls estimate that 51% of the U.S. population age 18-24 are women and 49% are men – nonresponse or "gender unknown" are not provided with the dataset, effectively erasing the existence of those who do not ascribe to this binary. The utter lack of large data sets using expansive notions of gender is disheartening. Two estimates suggests that (at least) 0.6% of the US population identify as transgender (perhaps with higher rates among younger Americans) but these did not account for genderfluid, genderqueer, two-spirit, or other non-binary gender identities (Crissman et al., 2017; Flores et al., 2016). Growing understanding of gender identity as a social construct with material implications for the lived experience is leading people to call for additional data collection by federal agencies (e.g., Freeman et al., 2018) and scholars (e.g., American Psychological Association, 2016) but these have not, as of yet, gained much traction.

## 5. *Race-gender groups: Intersection & omission*

From here on, this chapter considers students' race-gender groups. As per CRQI, this allows us to consider "the impact of race and racism at its intersection with other forms of subordination" (Covarrubias & Vélez, 2013, p. 276). In this case, we consider race and racism at its intersection with gender and sexism. This is in part because race-gender groups are the most visible in terms of data accessibility and existing research. For example, documentation of racialized-gendered gatekeeping and classroom mechanisms as they serve to propagate exclusionary views of "who belongs" in mathematics and impact students' decisions to major in STEM (Leyva et al., 2020; Rainey et al., 2018). Some dimensions which I do not include here are socio-economic status, sexual orientation, citizenship status, geographic location, disability status, and age; in part because the data is not readily available for such an analysis. The analysis therefore loses sight of some of the heterogeneity and differential experiences of students within these race-gender groups.

### 5.1. *Questions for the reader:*

- What other intersectional demographic categories *ought* we to be able to report on?
- What "forms of subordination" or systems of oppression do those categories relate to?
- What other categories have you considered when assessing local phenomena?
- What will you start looking at now?

## 6. *Some Numbers*

Data used in this chapter is from the 2018-19 academic year (the most recent year at the time of writing for which data was available and not impacted by the COVID-19 pandemic). As data is presented, I ask that you pause to reflect on these numbers and your reaction to them. Are you surprised, pleasantly or unpleasantly? Are you comparing these numbers to some ideal, and where did the ideal come from? These are some of the questions we grapple with as we assess current evidence of inequities and attempt to measure progress in the coming years.

## 7. *Aggregate Degrees Awarded*

BA/BS degree completion in the 2018-19 academic year data reported by 2057 postsecondary institutions is presented in Table 2. I distinguish between US residents and non-residents for a few reasons. These include: (a) non-residents include international students temporarily relocating to the United States while obtaining a degree; (b) I am interested in race/ethnicity questions in the context of US society; (c) the cultural experiences of non-residents in regard to race and gender may differ in important ways from those who reside permanently in the US; (d) it is not actually a marker of race or ethnicity, and so is not appropriate in the same analysis; and (e) the distribution of degrees varies significantly from that among US residents.

| | US Residents | Non-residents | Total |
|---|---|---|---|
| Engineering, Mathematics, Statistics, Physical Sciences (pSTEM) | 162,483 | 20,644 | 183,127 |
| Biological, Biomedical, Computer Science (sSTEM) | 188,404 | 11,616 | 200,020 |
| Non-STEM | 1,435,801 | 67,299 | 1,503,100 |
| **Total** | **1,786,688** | **99,559** | **1,886,247** |

*Table 2: BA/BS degrees awarded to US residents and non-residents in the 2018-19 academic year. We differentiate two kinds of STEM degrees: "pSTEM" refers to degrees awarded under CIP codes 14 (Engineering), 27 (Mathematics & Statistics), and 40 (Physical Sciences); "sSTEM" refers to degrees awarded under CIP codes 11 (Computer & Information Sciences) and 26 (Biological & Biomedical Sciences).*

The numbers in Table 2 include those US residents whose race/ethnicity is listed as "unknown," so as to provide an indicator of the total number of degrees awarded in various disciplines; in the following, I omit data corresponding to "race unknown," as it is not comparable to data from groups with identified racial/ethnic markers. Furthermore, "unknown" often results from an individual's decisions not to declare their racial identity and can be considered as a way of opting out of analyses like this one. We do include data from those categorized as "multiracial," i.e., those who select more than one racial/ethnic category on official documents, despite the many obvious shortcomings.

### 7.1. *Questions for the reader:*

- Do any of these numbers, or the relative size of these numbers, surprise you? Why?
- The vast majority (80%) of degrees obtained by US residents are in non-STEM majors. Does this seem appropriate? Why or why not?
- What rate of STEM degrees should we aspire to? Why?
- Among US residents, fewer degrees are awarded in the pSTEM fields (engineering, mathematics, and physical sciences) than in sSTEM (biological, biomedical, computer, and information sciences). Should we aim to change this ratio? Why or why not?

## 8. *BA/BS Degrees within Race-Gender Groups*

This section will provide a degree breakdown similar to that in Table 2, disaggregated by race-gender groups. Before reading the results, pause to consider what you expect to see, and why.

- If you assume that a particular race-gender group will have a higher or lower rate of attaining pSTEM or sSTEM degrees, is that…
  - Because you ascribe to certain stereotypes?
  - Because you've read about similar data elsewhere?
  - Because you recognize systemic racialized and/or gendered barriers to access and acceptance in STEM?
- What systemic gendered and/or racialized barriers to participation and continuation in STEM do you expect to see reflected when this data is disaggregated by race-gender groupings? What size impact do you expect to see?
- What variation in the proportions (pSTEM, sSTEM, non-STEM) across demographic groups do you think is acceptable? Why?

Note that this measure is calculated out of the total number of undergraduate degrees granted on a national scale, it does not reflect outcomes at particular institutions. These numbers also do not reflect students' original intentions (i.e., how many students initially intended to major in STEM), nor incorporate the numbers of students who leave college without degrees. This analysis provides a blunt, aggregate metric of the proportion of undergraduate degrees – for certain demographic subgroups – which are STEM degrees. And yet – from these tables, it is clear that there are racialized and gendered patterns in STEM degree completion – we can say little else from this data alone. Table 3 presents the rates at which each population obtains degrees in pSTEM, sSTEM, and non-STEM degrees. After the table, the same information is presented in pie charts (Figure 1).

| | | BA/BS Degrees Awarded 2018-19 AY | | | |
|---|---|---|---|---|---|
| Race | Gender | Row Total | Engineering, Mathematics, Statistics, Physical Sciences | Biological, Biomedical, Computer Science | Non-STEM Fields |
| American Indian and/or Alaska Native | Men | 3,020 | 10.93% (330) | 10.13% (306) | 78.94% (2,384) |
| | Women | 4,990 | 2.99% (149) | 5.85% (292) | 91.16% (4,549) |
| Asian | Men | 64,536 | 21.09% (13,613) | 26.02% (16,795) | 52.88% (34,128) |
| | Women | 76,577 | 8.46% (6,478) | 18.73% (14,345) | 72.81% (55,754) |
| Black | Men | 61,182 | 7.82% (4,786) | 11.58% (7,086) | 80.6% (49,310) |
| | Women | 107,093 | 2.46% (2,635) | 7.87% (8,430) | 89.67% (96,028) |
| Hispanic and/or Latinx | Men | 101,078 | 13.97% (14,117) | 12.01% (12,139) | 74.02% (74,822) |
| | Women | 157,162 | 3.55% (5,573) | 7.54% (11,848) | 88.92% (139,741) |
| Multiracial | Men | 27,680 | 15.6% (4,318) | 14.57% (4,032) | 69.83% (19,330) |
| | Women | 39,681 | 5.41% (2,146) | 9.85% (3,907) | 84.75% (33,628) |
| Native Hawaiian and/or Pacific Islander | Men | 1,570 | 10.45% (164) | 12.93% (203) | 76.62% (1,203) |
| | Women | 2,000 | 2.95% (59) | 8.4% (168) | 88.65% (1,773) |
| White | Men | 473,971 | 15.88% (75,277) | 11.9% (56,423) | 72.21% (342,271) |
| | Women | 609,635 | 4.62% (28,143) | 7.64% (46,596) | 87.74% (534,896) |
| *Aggregate* | | *1,730,175* | *9.12% (157,788)* | *10.55% (182,570)* | *80.33% (1,389,817)* |

*Table 3. Percentages (counts) of BA/BS degrees awarded in pSTEM, sSTEM, and non-STEM fields to students of different race-gender groups. Percentages are of row totals (i.e., 8.46% of Asian Women's BA/BS degrees are awarded in pSTEM majors).*

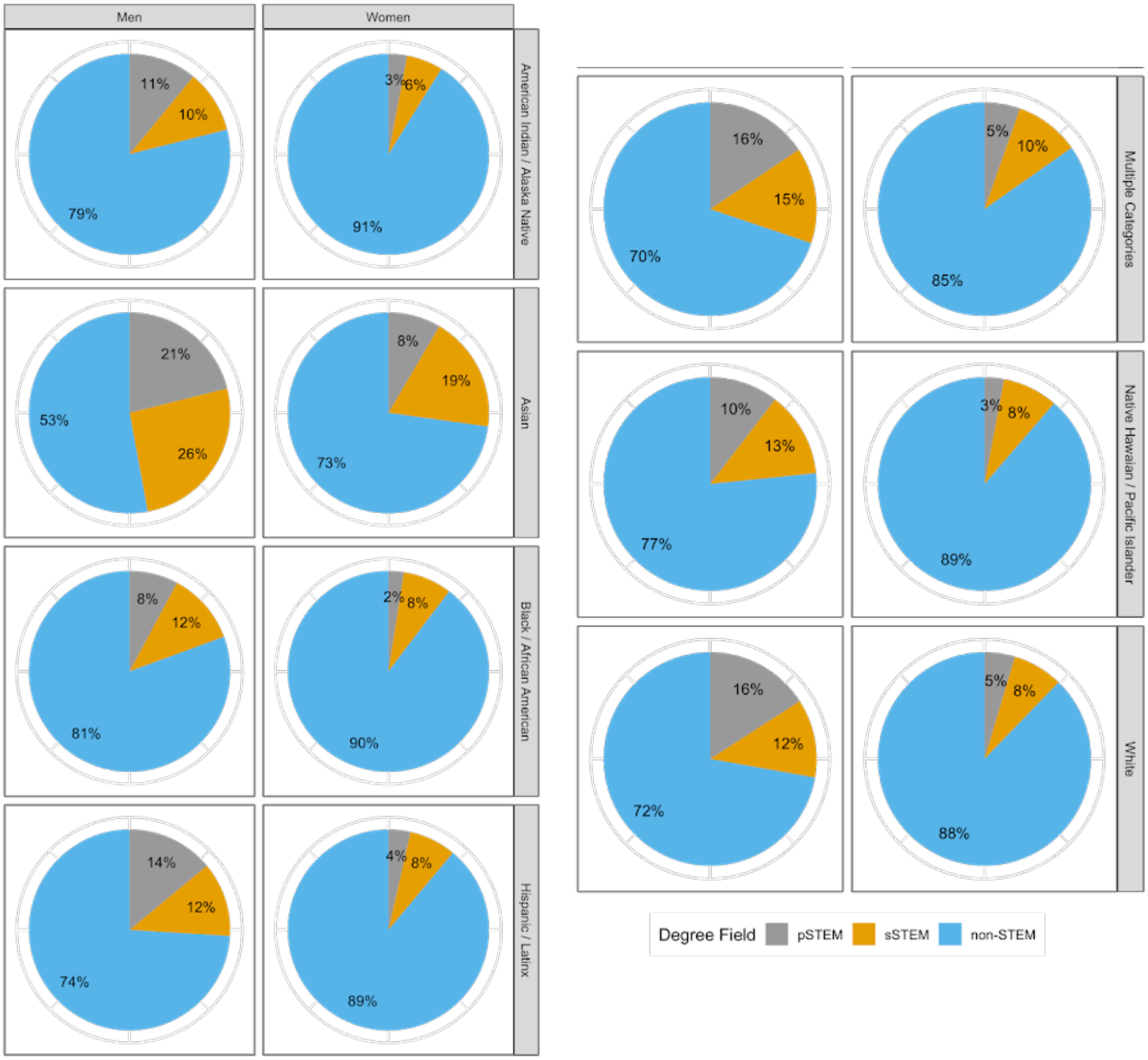

*Figure 1. Pie charts. Each cell describes the degrees awarded to a particular demographic group (race/ethnicity and gender), broken down into pSTEM (engineering, math, physical sciences), sSTEM (biological and computer sciences), and non-STEM (all other). The section dedicated to pSTEM starts at the 12 o'clock position and opens clockwise (gray); non-STEM starts at the 12 o'clock position and opens counterclockwise (blue); the remaining portion (orange) is dedicated to sSTEM degrees. Proportions are the same as those presented in Table 3.*

I expect that you have, in reading Table 3 and Figure 1, noticed some of the racialized and gendered patterns in this data and the wide range of rates of STEM degree attainment. Across race-gender populations, the rates of non-STEM degrees vary from 53% to 91%. Rates of pSTEM degrees vary from 2% to 21%. Rates of sSTEM degrees vary from 6% to 26%. However, the numbers do not speak for themselves (Covarrubias & Vélez, 2013) and so I offer some commentary and additional questions.

**9.** ***Deliberate choices***

The choices to use row totals (Table 3) and pie charts (Figure 1) were deliberate. By focusing on within-group patterns, we hope to minimize discussion of numerical, compositional, or representational diversity within STEM. I also hope to minimize comparative gap-gazing (Gutiérrez, 2008) by using graphical representations (pie charts) which are not ideal for comparisons (Tufte, 2001). I present these pie charts normalized to population size, to avoid implications that "small" groups are less worthy of study or attention *and* to avoid reinforcing status hierarchies because of subconscious "bigger is better" biases. Focusing on within-group patterns is also to avoid zero-sum thinking about STEM degrees. Had I displayed the proportion of STEM degrees awarded to each race-gender group, it would appear that at least one group's participation in STEM must decrease if another group's participation in STEM increases. This is not how it works – there is plenty of science (and science careers) to go around (PCAST, 2012).

**10.** ***(Even more) sobering thoughts and questions***

The data presented in this chapter show unequal outcomes for different race-gender groups. If they showed equal (or nearly equal) rates of STEM degree completion, that would not in and of itself indicate equity. Compositional diversity is one thing, but it does not equate to respectful treatment or inclusion. Consider, for example, that there are large numbers of Asian students in STEM fields, but these students are subject to damaging "model minority" myths which are dehumanizing (Shah, 2019) and lead to stereotyping and discrimination (Trytten et al., 2012). For those subject to negative stereotypes about their group's ability or propensity to succeed in STEM, participating in STEM may require *stereotype management*, wherein students expend additional energy not only recognizing stereotype threat and discriminatory comments but redefining their identities and/or success – an additional burden which is not carried by those who are assumed to be competent (McGee & Martin, 2011). Classroom mechanisms (e.g., interactions, dialogue, perceptions) are racialized and gendered in ways which contribute to stereotype visibility, STEM identity, and feelings of belonging/inclusion (Leyva et al., 2020); these are also not captured in this chapter's analysis. CHAPTER 7 of this volume incorporates students' voices to add nuance and a human touch to the stories behind these numbers.

This chapter considers student groups' propensity to complete an undergraduate degree in STEM. But this is not the end of the path to STEM careers, and that road is further littered with racialized and gendered obstacles. For example, Castro and Collins (2021) report on Asian-American women's experiences pursuing doctorate degrees in STEM. In these dominantly white and male-normed environments, participants reported gendered and racialized microaggressions, stereotypes, harassment, and hostile environments which led to (among other outcomes) re-examination and re-negotiation of their own identities. Though participants were still in STEM, they reported stories of others who left the field due to the exclusive culture and pressures to assimilate and accept mistreatment (Castro & Collins, 2021). It has been pointed out in other contexts that diversity and inclusion (and an increasing list of associated buzzwords) often mask demands for assimilation and participation in problematic structures (Martin, 2019). Discrimination, pressure, isolation, and pressure to conform continue throughout postdoctoral professions as well and lead to trauma, exits, and career changes (e.g., Cheney & Shattuck, 2020; Goins, 2017; Harmon, 2019).

10.1. *Questions for the reader:*

- How will you avoid interpreting Table 3 and Figure 1 in ways that reinforce hierarchies and/or stereotypes about students from different race-gender groups?
- How do you think some of the missing aspects of student context (e.g., age; sexual orientation; ability status) might change these results?
- How do you expect the numbers in Table 3 and Figure 1 compare to students' initial interest (i.e., initially intended major)?
- How can you begin to remake STEM spaces to be more inclusive and supportive?
    - Can you do so while allowing students to maintain and celebrate their identities, rather than assimilate into a dominant image of what mathematicians and scientists look like?

**11. *Applications to Local Context***

The effects of intersecting and pervasive systems of racism and sexism as they relate to STEM impact the rates at which members of race-gender groups complete their degrees in STEM. It is not race and/or gender which determine likelihood to major in STEM, but the disproportionate and differential impact of racialized-gendered stereotypes, gatekeeping, barriers, and exclusion within and around STEM. These are systemic problems which require systemic dismantling at multiple levels.

This must be addressed and must be addressed at every turn and in every way possible. This includes (but is certainly not limited to) admissions, placement processes, graduation requirements, assessments of learning, classroom practice, program support, local culture, individual interactions, messaging, and disciplinary culture.

One of our goals is to spur you to investigate students' pursuit and completion of undergraduate STEM degrees in contexts you have access to, and to do so with a critical eye. We have offered up one approach to doing this but have noted extensive limitations and problematic aspects to the approach here. We encourage you to explore students' experiences along with STEM-related outcomes, and to do so with a critical and intersectional lens to the extent possible. Some questions to ask yourself (and others you might work with) when undertaking such investigations:

- Can you access local data about students' major intentions and rates of switching/persisting across race-gender groups? Can you incorporate other dimensions of student identity and context? If not, why not? Can anyone?
    - How can you cooperate with other institutional units to capture, interpret, and act on this information?
- How can you capture additional aspects of students' lived experiences, beyond outcome measures and numeric representation?
    - How can this information inform program reform to improve students' experiences in mathematics classes? In the rest of their academic experience?
    - How can you (and others) use this information to remake exclusive spaces into inclusive ones?

- Does the institution, and any program within the institution, support whole students in their pursuit of fulfilling educational experiences?
- In reflecting on your data and analysis, consider the ways in which it can:
    - Help **you** learn and guide your future actions to support change
    - Support program assessment and contribute to sustained improvement without settling for marginal numeric changes
    - Engage colleagues within and across disciplines and institutions to grow and develop into more humane spaces
    - Contribute to structural work aimed at creating justice out of unjust systems

## 12. *Parting Thoughts*

Other chapters in this volume describe 30 current (or recent) case studies of efforts to address inequitable experiences and/or unequal outcomes in calculus, which is a foundational course for STEM majors and a key gatekeeper. These are starting points, conceived of and executed at a time when inequities are stark and painfully obvious despite limited and non-critical data. Most of the initiatives reviewed in this text leverage locally accessible levers, taking advantage of opportunities at hand to make some change possible. Thematic chapters review subsets of these cases with an eye toward initiatives based on developing change agents (CHAPTER 2 ), classroom practices (CHAPTER 3 ) , role models (CHAPTER 4 ), course design (CHAPTER 5 ), language use and identity (CHAPTER 6 ), and assessment (CHAPTER 7 ) . The final thematic chapter (CHAPTER 8 ) provides the perspective of those who matter most: STEM students themselves.

We hope that you will consider the lessons learned from this volume and take steps to make change in your own local context. In doing so, we urge you to go beyond what is presented here in considering intersecting axes of inequity and imagining strategies for change.

### 12.1. *Final questions for reflection and exploration:*

- In what ways are intersecting systems of hierarchy, discrimination, and oppression present in postsecondary STEM education, and STEM more broadly?
- To what extent are you complicit in upholding those systems?
- What can you do to dismantle those systems and enlist allies in that process?
- How will you know if progress is being made?

## 14. *Appendix A: Data Sources*

All data used in this work is publicly available. Institutional data comes from IPEDS, US demographic data comes from the US Census databases and projections/analyses of that data.

The following criteria were used to select institutions from the available IPEDS database. These selection criteria reflect questions about STEM BA/BS degree attainment at the types of institutions open to all students in the United States. Using IPEDS selection terminology, for those seeking to replicate or update this analysis:

*IPEDS "First Look" Universe Member*

- Report data to IPEDS
- Open to the public
- Participant in federal financial aid programs (Title IV Participant)

*Sector*

- Public, 4-year or above
- Private, not-for-profit, 4-year or above

*Other*

- Degree-granting institution
- Not all programs offered completely via distance education

In all, 2057 institutions met the selection criteria **and** reported the relevant data for the 2018-19 academic year. These institutions are located in all 50 US states and the District of Columbia. Of these, 733 (36%) are public institutions, the other 1324 (64%) are private not-for-profit institutions (though in general, there are more students enrolled in public institutions). There are 84 Historically Black Colleges and Universities, 16 Tribal Colleges, and 87 Land Grant institutions in the sample. Of the 2057 schools, 535 (26%) award no degrees higher than a bachelor's degree, while the others offer one or more graduate degrees. Using the 2018 Basic Carnegie Classifications (collapsed for simplicity):

| Institution Type | Count |
|---|---|
| Associates | 12 |
| Baccalaureate | 638 |
| Masters | 601 |
| Doctoral | 390 |
| Special Focus | 369 |
| Tribal College | 16 |
| Other/No Class | 31 |
| *Total* | *2057* |

## 15. *Appendix B: CIP Codes*

*Classification of Instructional Program* (CIP) codes are used by the National Center for Education Statistics to categorize programs of study across institutions (see https://nces.ed.gov/ipeds/cipcode/ for exhaustive lists). While many agencies and groups define STEM using CIP codes, there is not much agreement on what counts as STEM. For simplicity, we restricted our categorizations here to the two-digit top-level CIP codes, further distinguishing based on whether a majority of undergraduate majors within each code require the P2C2 sequence or a more varied set of mathematics courses as a prerequisite to completion. Below we present the two- and four-digit CIP codes included in this study (note that there is a finer-grain six-digit categorization which we omit for space):

Table B.1 Four-digit subcodes of the three two-digit CIP codes categorized in this chapter as pSTEM.

| | |
|---|---|
| **14. Engineering** | |
| 14.01 Engineering, General | 14.24 Ocean Engineering |
| 14.02 Aerospace, Aeronautical and Astronautical Engineering | 14.25 Petroleum Engineering |
| 14.03 Agricultural/Biological Engineering and Bioengineering | 14.27 Systems Engineering |
| 14.04 Architectural Engineering | 14.28 Textile Sciences and Engineering |
| 14.05 Biomedical/Medical Engineering | 14.32 Polymer/Plastics Engineering |
| 14.06 Ceramic Sciences and Engineering | 14.33 Construction Engineering |
| 14.07 Chemical Engineering | 14.34 Forest Engineering |
| 14.08 Civil Engineering | 14.35 Industrial Engineering |
| 14.09 Computer Engineering, General. | 14.36 Manufacturing Engineering |
| 14.10 Electrical, Electronics and Communications Engineering | 14.37 Operations Research. |
| 14.11 Engineering Mechanics. | 14.38 Surveying Engineering |
| 14.12 Engineering Physics. | 14.39 Geological/Geophysical Engineering |
| 14.13 Engineering Science. | 14.40 Paper Science and Engineering |
| 14.14 Environmental/Environmental Health Engineering | 14.41 Electromechanical Engineering |
| 14.18 Materials Engineering | 14.42 Mechatronics, Robotics, and Automation Engineering |
| 14.19 Mechanical Engineering | 14.43 Biochemical Engineering |
| 14.20 Metallurgical Engineering | 14.44 Engineering Chemistry |
| 14.21 Mining and Mineral Engineering | 14.45 Biological/Biosystems Engineering |
| 14.22 Naval Architecture and Marine Engineering | 14.99 Engineering, Other |
| 14.23 Nuclear Engineering | |
| **27. Mathematics and Statistics** | |
| 27.01 Mathematics | 27.05 Statistics |
| 27.03 Applied Mathematic | 27.99 Mathematics and Statistics, Other |
| **40. Physical Sciences** | |
| 40.01 Physical Sciences | 40.06 Geological and Earth Sciences / Geosciences |
| 40.02 Astronomy and Astrophysics | 40.08 Physics |
| 40.04 Atmospheric Sciences and Meteorology | 40.10 Materials Science |
| 40.05 Chemistry | 40.99 Physical Sciences, Other |

Table B.2 Four-digit subcodes of the three two-digit CIP codes categorized in this chapter as sSTEM.

| **11. Computer and Information Sciences & Support Services** | |
|---|---|
| 11.01 Computer and Information Sciences, General | 11.07 Computer Science |
| 11.02 Computer Programming | 11.08 Computer Software and Media Applications |
| 11.03 Data Processing | 11.09 Computer Systems Networking and Telecommunications |
| 11.04 Information Science/Studies | 11.10 Computer/Information Technology Administration and Management |
| 11.05 Computer Systems Analysis | 11.99 Computer and Information Sciences and Support Services, Other |
| 11.06 Data Entry/Microcomputer Applications | |

| **26. Biological and Biomedical Sciences** | |
|---|---|
| 26.01 Biology, General. | 26.10 Pharmacology and Toxicology |
| 26.02 Biochemistry, Biophysics and Molecular Biology | 26.11 Biomathematics and Bioinformatics. |
| 26.03 Botany/Plant Biology | 26.12 Biotechnology |
| 26.04 Cell/Cellular Biology and Anatomical Sciences. | 26.13 Ecology, Evolution, Systematics, and Population Biology |
| 26.05 Microbiological Sciences and Immunology | 26.14 Molecular Medicine. |
| 26.07 Zoology/Animal Biology | 26.15 Neurobiology and Neurosciences. |
| 26.08 Genetics. | 26.99 Biological and Biomedical Sciences, Other. |
| 26.09 Physiology, Pathology and Related Sciences. | |

# Part 2: Cross-Cutting Thematic Chapters

# 2 Developing the Next Generation of Change Agents in College Mathematics Instruction

**Shandy Hauk**
*San Francisco State University*

Shandy Hauk <hauk@sfsu.edu> grew up in and around Los Angeles. Since finishing a mathematics Ph.D. at the University of California, Irvine, and a mathematics education post-doc at Arizona State, she has specialized in intercultural competence and post-secondary mathematics education research and development. She is a faculty member in Mathematics at San Francisco State University.

**Natasha Speer**
*University of Maine*

Natasha Speer <natasha.speer@maine.edu> grew up in Brooklyn, New York and Eastport, Maine. She is a faculty member at The University of Maine, affiliated with the Department of Mathematics and Statistics and the Maine Center for Research in STEM Education. Her work includes design of materials and programs for teaching-focus professional development as well as research on the teaching of undergraduate mathematics.

## 1. *Pre-Reading Reflection Questions*

Bring to mind the teaching-related professional learning opportunities that are available to novice instructors in your department (e.g., graduate students, new faculty, undergraduate learning assistants). These would include structured professional development courses, workshops, and coordination meetings. Make a few notes about what you know related to the questions:

1. What *skills*, *knowledge*, and *orientations to working with others* (e.g., students, colleagues) will novice instructors need in order to realize intended teaching-related goals?
2. In what ways are the opportunities designed so that novices *experience the same kinds of instructional practices in professional learning* as they are expected to enact in classrooms?
3. What information do you or colleagues gather to provide *data about how novice instructors are creating student learning opportunities reflective of department/program goals*?
4. What might novices learn by being a part of your program that will help *prepare them to seek out and pursue opportunities to enact change later in their careers*?

## 2. *Introduction*

The case studies in this volume capture a diverse and exciting collection of examples of change. Some are changes in content, others in processes. All occur within and across various (sub)systems of post-secondary mathematics. The values behind, and goals of, the changes differ. Yet, the stories of change share something: the need for professional learning by instructors. In most case studies, the instructors who were involved included novices such as undergraduate learning assistants, graduate student instructors, and early career faculty.

Making progress in justice, equity, diversity, and inclusion requires systemic change and the development of novice college mathematics instructor professional knowledge is a critical subsystem.

The mathematics community is fortunate to include the authors of the case studies, who have taken up the challenge of improving learning for current undergraduate students. But, from where will the next generation of such folk come? What can be done now to ensure the mathematics community has people who are equipped to address future obstacles and leverage future opportunities in equitable and just post-secondary mathematics teaching, the particulars of which may be unknowable now? In other words, what does it mean for professional learning opportunities to support novices to become effective instructors now and also join the ranks of the next generation of change agents? As noted in Larnell's opening manifesto in this volume (BUILDING AN EVOLVING FRAMEWORK: A CLARION CALL / MANIFESTO), the answers require a multi-threaded effort across personal reflections, curriculum development, cycles of instructional preparation, implementation, and evaluation, department- and division-level accountability, within and across-institution continuous improvement work, and broad policy efforts to support large scale change within and across subsystems.

### 2.1. *Looking Back: Incremental Change*

The history of educational reform is replete with examples of how excellent efforts, aimed at only some subsystems, can fail. Without attention to instructional development needs, some beautifully designed curricula have fallen short of their potential to improve student learning. Conversely, efforts to have instructors adopt a particular practice, such as collaborative group work, have been hampered by a scarcity of group-worthy tasks. To succeed as an instructional change agent, one needs awareness and knowledge about subsystems, understanding and anticipation about how those subsystems connect, interact, and influence one another, and skill at seizing opportunities to enact change as well as capacity to generate new learning from each professional encounter as preparation for future change-focused endeavors.

### 2.2. *Looking Ahead: Generative Change*

As was noted more than 20 years ago by Franke and colleagues (1998),

> *Self-sustaining, generative change does not involve acquiring a set of procedures to implement with fidelity; rather it frequently entails teachers making changes in their basic epistemological perspectives, their knowledge of what it means to learn, as well as their conceptions of classroom practice. It means conceptualizing teacher change in terms of teachers becoming ongoing learners (p. 67).*

Knowledge about equitable and just post-secondary teaching and learning is rapidly growing. Today's novice instructors will be better equipped to navigate future demands and shape tomorrow's departments and universities through self-sustaining, generative change if their professional learning is informed by the growing knowledge base, including the lessons today's change agents shared in the case studies.

This chapter is not a primer on how to design professional learning about teaching for novice instructors (for that, see, e.g., College Mathematics Instructor Development Source, 2021;

Council of Graduate Schools, 2021; Deshler et al., 2015; Saichai & Theisen, 2020 and references therein). Rather, this report digs into two dimensions of professional growth related to diversity, equity, and inclusion (DEI). We contend that the development of needed DEI-rich knowledge depends on two skills for interacting with ideas and people: *decentering* and *interconnecting*. The aim of this chapter is to describe and illustrate the role of these dimensions in professional development for novice college mathematics instructors. The report closes with next steps and resources for purposeful attention to decentering and interconnecting as vital components in developing novices as future change agents in college mathematics instruction.

## 3. *Decentering as a Professional Skill*

Around the world, calculus courses and mathematics programs in most post-secondary institutions are built on an instructor-centered model. This approach has been effectively self-sustaining for many decades. The practice of lecturing has been passed on from generation to generation of college teachers through personal classroom experiences and through graduate school training rooted in curricula that preserve lecturing as the status quo. Now, however, it is clear that an instructor-centered approach is not universally effective or appropriate (see, e.g., Abell et al., 2018; Bressoud et al., 2015; Freeman et al., 2014; Laursen et al., 2014).

The apprenticeship of observation is powerful. People tend to teach the way they were taught (Lortie, 1975). It is important for novice instructors to experience teaching that models and provides touchpoints in their efforts to teach differently, more equitably. For example, professional learning opportunities can be offered in ways that model equitable instruction. Thus, novices can refer to how they have *recently* been taught (e.g., in professional development) to contrast with the power-culture-driven, instructor-centered experiences that likely make up the bulk of their histories as learners.

### 3.1. *What is Decentering?*

Successful implementations of the kinds of instructional practices described throughout this volume call for teaching that elicits and utilizes student contributions. Now, the volume authors argue, the expectation is that instructors facilitate discussions to which students contribute their thinking and voices. This kind of *instructional decentering is, at its most basic, the act of seeing from someone else's point of view* and has historical roots in the work of Piaget (1955). It means engaging with other people as one of the *participants* in interaction, rather than as the center of interaction.

In decentering, instructional attention is on uncovering, understanding, and expanding on what students know and do to include novel, non-standard, and standard mathematical ideas and methods (Carlson et al., 2007; Rahman, 2018; Teuscher et al., 2016). Being self-aware and facilitating self-aware learning by students are the focus (instead of attention and authority being vested largely or solely in the instructor).

### 3.2. *Development of Decentering Skills*

Decentering requires attention to other people as (potentially) different from oneself, noticing nuances in similarity and difference between one's own perspectives or experiences and those of others. In its most developed forms, instruction that leverages decentering also bridges across similarities and differences in formulating in-the-moment responses to situations.

Decentering depends on a variety of individual instructor factors (e.g., self-knowledge, goals, orientations, beliefs, psycho-social challenges). Professional development can provide opportunities for instructors to build skill in decentering, along with other facets of cross- or intercultural competence. There is a developmental continuum for decentering: from an ethno-centric view that everyone is like me to an ethno-relative view, that any person (including me) is like and unlike every other person in identifiable and valuable ways (Bennett, 2004).

In particular, learners in instructor-<u>de</u>centered instruction are encouraged to see themselves as decision-makers in the classroom. Novice instructors need to learn how to initiate, sustain, and manage undergraduates as participants in student-centered learning. This includes instructors learning about many things from a student-centered perspective, such as content, curriculum, and assessment (Bok, 2009), communication and interaction (e.g., related to classroom authority or socio-political factors, Gutiérrez, 2009, 2013; Winter & Yackel, 2000), as well as how to learn in and from instruction itself (Speer & Hald, 2008). Learning to elicit student thinking and learning how to shape instruction based on that thinking is the foundation on which generative change is built (Franke et al., 1998).

To be a change agent in the future, instructors will need decentering skills and habits in and outside the classroom. Decentering will allow them to engage in a cycle of learning from their experiences so that they can change themselves, be aware of the experiences of others, and formulate responses that are tuned to circumstances and challenges that might arise.

### 3.3. *Illustrating the Idea of Decentering*

Decentering as an instructor means inviting and scaffolding students to contribute to mathematically dense conversations in the classroom. The case studies share a variety of approaches to supporting instructional decentering.

Consider the course redesign for calculus 1 described in the Oliver et al. case study (CS 21 ). Their multi-pronged approach was rooted in an instructional tool – a "Dynamic Calendar" of resources available to instructors that included materials designed to promote student engagement, tools for formative assessment, and videos for use in a flipped classroom format. Oliver and colleagues' report (CS 21 ) includes glimpses into the professional learning by instructors in a local Mathematical Community of Practice. Monthly community meetings included dissemination of ideas by leaders that were relevant to instructional decentering such as: using student names (e.g., through nametags on desks or online tools in remote teaching), creating multiple mechanisms for student contribution (e.g., writing and speaking), and orchestrating multiple small conversations (e.g., (re)assigning people to groups). These are all examples of acknowledging students and creating opportunities for instructors to make student thinking "visible." Thus, though primarily about their Dynamic Calendar and examining the potential downstream effects of its use as evidenced by student grades, Oliver et al. (CS 21 ) provide an example of an important step in building skill at noticing how the thinking of students is similar to or different from an instructor's own: creating the opportunity in one's classroom to hear and see student thinking.

In a different vein, the Soto et al. case study (CS 22 ) centers on an investigation of instructor learning (rather than undergraduate student learning). Soto et al. (CS 22 ), created and then examined the impact of a professional course for instructors. The course aimed to build instructor knowledge and lived-experience of classroom community, rich tasks, and learner-centered teaching. The authors note that modeling target practices is important and not only did the course encourage decentering by instructors, the providers of the course also engaged in decentering as they taught the instructors. The positive impact on instructors (as learners) mirrored those intended for undergraduates. Participating instructors reported an increase in valuable, useful knowledge, a shift in identity, and a shift in perceptions of power to include students as essential in making decisions in and for classroom activity. These, the case study authors argued, were a consequence of the designers/providers of the professional learning engaging in learner-centered approaches with the instructor-participants.

Two other case studies report similar efforts, though with undergraduate learning assistants (LAs). Liou-Mark et al. (CS 16 ) and Villa et al. (CS 27 ) discussed components in the professional development for undergraduate peer leaders that were consistent with decentering (e.g., have peer-leaders attempt the mathematics course learning modules and come together to compare similarities and differences in problem-solving strategies). The curricular activities were designed to showcase a diversity of solution methods. Both Liou-Mark et al. (CS 16 ) and Villa et al. (CS 27 ) report ideas consistent with decentered practices (e.g., cooperative learning and being inclusive) and some of the approaches in Villa et al. (CS 27 ) have potential for increasing linguistic and socio-cultural supports for students. There is, however, one subtle difference between the two in the preparation for decentering of novice LAs that may interfere with decentering. By having LAs orchestrate a consensus among students on a publicly shared solution method as well as a final answer, the LA-led solution documentation might not actually represent a strategy used by any student in the room. This could constrain opportunities to celebrate a diversity of strategies and the documented solution may not represent how any actual student thought (except, perhaps the LA; Villa et al. (CS 27 )).

In the Jensen-Vallin case study (CS 13 ), the department implemented a program to provide coordination and curriculum development for pre-calculus (taught through paired co-requisite and main courses). Meetings of course instructors included discussions of content, teaching practices, and pacing. As in Oliver et al. (CS 21 ) and Soto et al. (CS 22 ), the changes were intended to build a community of practice among instructors, with attention to diversity and inclusion. Unlike the other two case studies, Jensen-Vallin (CS 13 ) reported re-purposing an existing tool to support instructors to decenter themselves within their community (i.e., a truncation of a model from the Dana Center). That is, instructors were asked to notice how they were similar to and different from each other.

In the case study by Johnson et al. (CS 14 ), a new curricular component (Techtivities) was introduced into a college algebra course. Teaching guides and professional development included attention to the redistribution of power in the classroom, so that instructors were not seen as sole experts and student contributions were encouraged. Suggestions for instructor self-awareness to decenter also appeared as annotations in the teaching guide for the Techtivities (e.g., Figure 2, Johnson et al. (CS 14 )). That is, both the professional development and the written supports for use of the curriculum were intended to be student-centered and were explicit about instructor decentering.

In Bennett et al. (CS 3 ), a two-part effort, also in pre-calculus, is described. The authors report on trying to influence instruction at their large, historically white university by changing curriculum in the pre-calculus courses to rely on mastery-based assessment for 50% of the course grade and by adding an inclusive teaching workshop with follow-up meetings to the preparation of new instructors. The professional development program in Bennett et al. (CS 3 ) had been in operation for decades and already had individual, department, and administrative supports focused on inquiry-based instruction. In the case study work, professional learning for novice instructors was broadened to provide decentering activities, including a skit-based experience and a thread throughout about how students may differ from the instructor (i.e., explicit attention to decentering). At the same time, as the authors note, the program served an inequitable system. Introducing mastery-based exam structures perturbed that system in a way that was student-centered, offering students more opportunities to attempt assessment. The bulk of the grade in the course was, however, still determined by proctored assessments (nine “mastery assessments” and two exams).

## 4. *Interconnecting as a Professional Skill*

While decentering is awareness from within the perspectives of others, *interconnecting uses meta-awareness to make connections across perspectives and contexts*. This can occur at many levels and grain sizes. Linking across contexts is essential in developing and nurturing coalitions, an essential component of local and systemic change (Kotter, 2012). The case studies and other essays in this volume describe many instructional practices (e.g., in Figure 1, inside the disk labeled INSTRUCTION) and several describe content aspects (Figure 1, MATH). A few case studies address what happens from the perspective of the third layer in Figure 1, PROFESSIONAL DEVELOPMENT. The previous discussion of decentering referenced several of these (Johnson et al. (CS 14 ), Oliver et al. (CS 21 ), Bennett et al. (CS 3 ), Liou-Mark et al. (CS 16 ), Soto et al. (CS 22 ), and Villa et al. (CS 27 )). This chapter, for those involved in providing professional development, is a contribution to the outermost region in Figure 1, LEADERSHIP DEVELOPMENT.

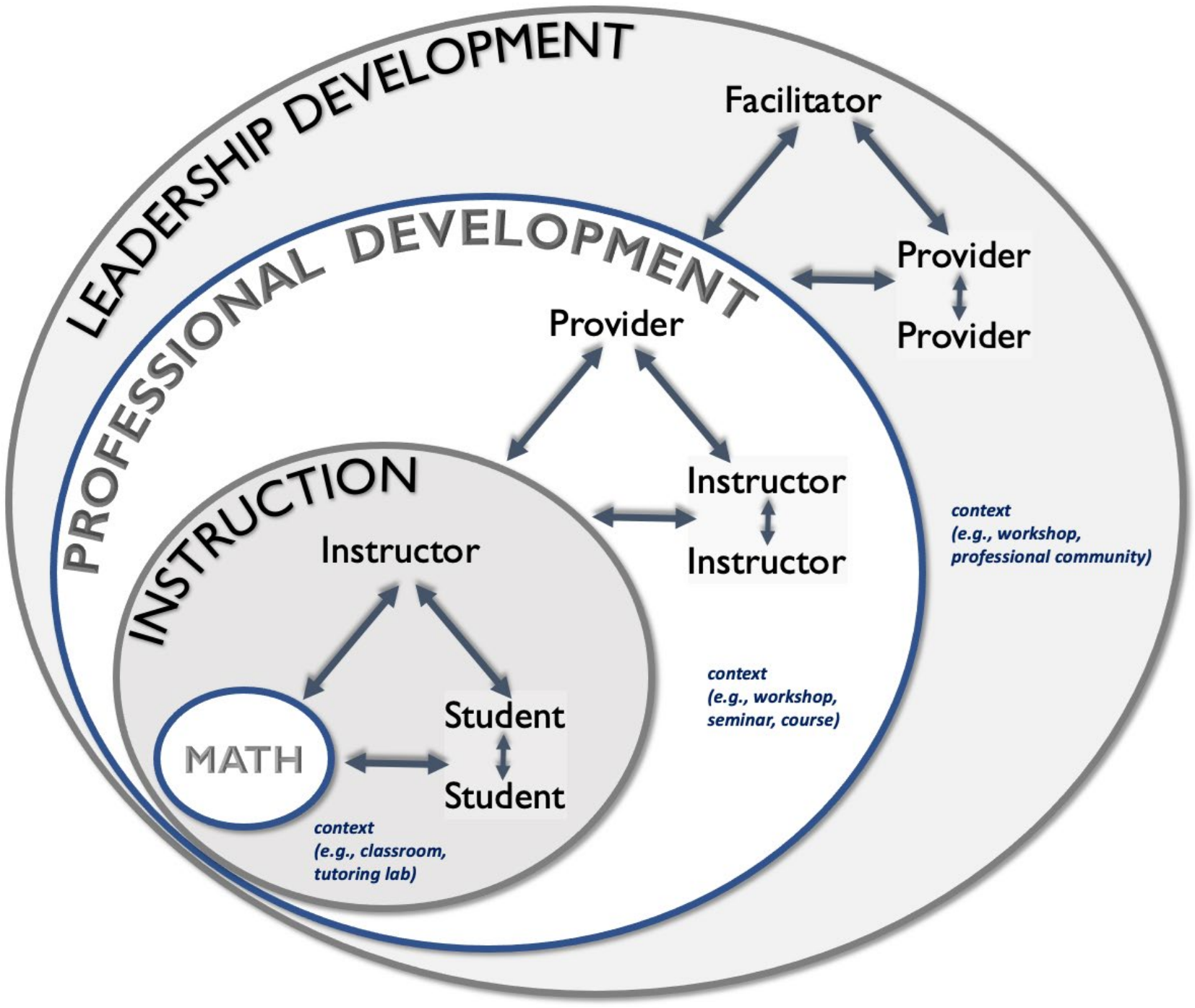

Figure 1. Nested structure of professional learning (adapted from Carroll & Mumme, 2007).

### 4.1. *Interconnecting: Examples and Non-Examples*

As an example, consider interconnecting rooted in the concept of derivative (MATH in Figure 1). If an instructor is aiming to help students learn about the idea of a derivative, the instructional goals for that are influenced by what the instructor may understand about students' knowledge of slope, ratio, and change. Interconnecting by the instructor involves noticing how students' conceptions may support or constrain the way learning progresses. That is, the instructor considers what is happening at the Student <—> Student node in INSTRUCTION (Figure 1), where student thinking might include $m$ in $y = mx + b$, previous experience with unit rate, experiences that discretize change (e.g., compare slope at point A to slope at point B), or treating change as covariational. Instructors develop interconnecting skills by knowing these student conceptions as well as the dynamics of communicating about them in a multi-contributor, student-centered, context (the arrows in the INSTRUCTION region). Instructors also consider and link to students' thinking in selecting what type of learning supports might be useful (e.g., group work on some problems, viewing a pre-class video for a topic). Instructors connect across and prioritize the mathematical and contextual factors, to decide what is instructionally useful. In this part of the example, the interconnecting is all within the INSTRUCTION region of Figure 1.

A *non-example of interconnecting* would be an instructor deciding to use collaborative group work as an instructional approach without also having a set of problems designed to support student thinking and learning in a group setting (e.g., the instructor might not yet know or connect the nuances of group-worthy tasks to the instructional context of group work). Another non-example of interconnecting would be an instructor asserting a single approach (e.g., their own preferred method) to solving a differentiation problem using a limit quotient. Learning opportunities for students are reduced if an instructor relies only on one way of thinking.

Notice that "Instructor" and "Student<—> Student" are separate nodes in Figure 1. Rather than INSTRUCTION containing only Student and MATH with an arrow between them, the figure makes explicit that the instructor is one component of interaction and students interacting with each other is another. With interconnecting meta-awareness, instructor decisions about what problems to choose and how to elicit and respond to student thinking during class are informed by what is known (personally and more generally) about the diverse ways students make sense of limit and limit quotient, lesson designs that elicit student thinking, and classroom activity structures to support every student in constructing mathematically robust understandings. Also, as diagrammed in the INSTRUCTION region in Figure 1, as students are experiencing the interactions with the content and each other, the instructor is attending to both the mathematics and the students, making meta-aware decisions about how to design and support those interactions. Novices learn about these interacting connections in and from their teaching and through professional development.

When someone embarks on *providing* professional development to support instructors to learn to teach, they are taking on the challenge of thinking at another level out, a meta-meta-awareness of connections is required. Providers of professional development think about how to have novice instructors interact with ideas of mathematics teaching (and each other) in ways that will help those instructors think about how to support students to interact productively with mathematical ideas (and each other).

Continuing the derivative example in the PROFESSIONAL DEVELOPMENT region of Figure 1, for the provider there are additional considerations about what instructors will need to know and do if they are to create the desired learning opportunities for students. For example, novice instructors might need assistance to build a clear idea of what learning goals for the derivative are. Identifying learning goals might occur in various professional activities (e.g., as a step in lesson design, or as a provided component in a set of curricular materials). Novice instructors may also need opportunities to learn about, have practice with, and connect across the mathematics of limit, limit quotient, and derivative along with instructional approaches that may be particularly effective (e.g., students working on a group-worthy task in which limits, rates, and limits of ratios are compared and contrasted). This, too, may occur through various professional learning tasks that are designed for use with novice college instructors (e.g., in an interactive activity, or in a teaching guide).

Also needed by novice instructors– as they notice and structure their understanding of a web of interconnections– is guidance about how to orchestrate mathematical and instructional ideas while decentering. This includes building knowledge about working in racially, ethnically, and linguistically diverse classrooms (e.g., a reading about student funds of knowledge and how to leverage those in teaching; González, et al., 2011) as well as creating and maintaining socio-mathematical and social norms in the classroom (e.g., a professional learning activity about what to do the first day of class to begin setting norms). These all make up the “content” of the professional development (in Figure 1, this includes the MATH, information about students and instructors as well as context knowledge from research and practice about the interaction arrows among MATH, Instructor, and Student). Also, part of information that novice instructors interconnect is the learning about teaching they encounter in hallway conversations with colleagues and other informal interactions (Latulippe, 2009).

In particular, interconnecting includes instructors thinking about students’ thinking about mathematics. In an analogous fashion, those who provide professional learning opportunities must concern themselves with an additional level of interconnecting: providers think about how instructors are thinking about how students are thinking about mathematics.

### 4.2. *Development of Interconnecting Skills*

Like decentering, skill at interconnecting can be purposefully developed. In particular, development begins from a self-focused denial of differences (e.g., no connections are needed since only MATH matters). Skills for interconnecting will further develop through a tendency to polarize and focus on mathematics as disconnected from human interaction (e.g., there is one best or right way to solve every problem). From there, development can progress to a search for universals, connections that are compressed into a single process, but not multiple interconnected processes (e.g., there are “objective” or “mastery based” ways that are universally applicable to assess *all* students, and grades become the essential element of interaction for instructor and student, disconnected from students’ mathematical funds of knowledge). With time and intentional development, one can learn more about mathematical ideas, contexts, and human interactions and reflect on teaching with greater attention to nuance (e.g., learning about implicit bias or micro aggression and suddenly noticing it in every word problem in the text), but how to use this knowledge to improve opportunities to learn, classroom climate, and interactions with others remains elusive. At its most developed, interconnecting is adaptive — networks of people and their interactions can be anticipated (enough) that teaching serves the needs of the people, and networks of people, in and outside of the room.

Interconnecting is important for change agents because they need to know how people, policies, and perspectives function in and across interacting systems (e.g., in and beyond those shown in Figure 1). The knowledge from that is extremely valuable when advocating for a change in one part of a (sub)system: one can anticipate how change in one place will cause or necessitate change in another. Such anticipation supports planful thinking about who the stakeholders are in the change and how to contact and be involved with them.

### 4.3. *An Illustration of Interconnecting*

Most of the case studies that describe professional learning report on groups involving four or more instructors. One-on-one or small group coaching is also a form of professional development. Stone-Johnstone et al. (CS 25 ) is a case study of an instructor who was working to identify and respond to equity related instructional goals. Decentering and interconnecting were the key tools of the small professional group in this case study. Unique to Stone-Johnstone et al. (CS 25 ) was the use of multiple classroom observations with the EQUIP method for documenting participation in classroom interaction (Reinholz et al., 2020a). The use of the method supported the case instructor to notice his assumptions about similarity, difference, and participation (decentering). Discussions by the professional group, based on the observation data, generated and strengthened interconnectedness among ideas and across awareness and meta-awareness for equitable and inclusive instructional moves. There was, quite literally, a non-instructor, non-student, “other” in the room (the observer). And the EQUIP method of tracking interaction provided evidence of what the other had observed about the mathematics (Figure 1, MATH), the classroom interactions (Figure 1, arrows between instructor and students, student and student), and the birds-eye view of the classroom (Figure 1, all of INSTRUCTION).

## 5. *Decentering, Interconnecting, Data, and Developing Agents for Change*

Take any case study where the department has implemented an innovative program and instructional-development questions arise: What should be happening for novice instructors so they can teach in the targeted ways? And, later, when a person graduates from Innovative U. and lands a job at Status Quo College, how does that person – who is a relative novice in college mathematics instruction and departmental politics – initiate and maintain productive exchanges to improve teaching? The answers to both questions include building social and management skills, particularly for interacting with people and policies seen as "other" in various ways (e.g., by decentering and interconnecting). For change that disrupts the academic mathematics status quo, skills must be rooted in an understanding of the norms and values of the status quo, how those are different from and similar to new norms and values (e.g., those anchored in justice and anti-racism), how to enact change from the former to the latter, and how to determine the nature of the success of the change effort (and start a new cycle of change based on it).

The case studies and chapters in this volume raise awareness that things do not look or work the same everywhere. They also provide examples of what it is like to engage in, struggle with, and monitor progress in efforts as change agents. Many case study endeavors were based on local evidence of success (e.g., course or department data on student achievement). Eliciting thinking from many students in the classroom allows instructional decisions to be informed by multiple perspectives. In a similar fashion, eliciting information from many sources at the next layer out (e.g., for a course with multiple sections or a department program) enables understanding and decision-making to be informed by a diverse and more comprehensive set of ideas. End-of-term grades are only one form of readily accessible local data (like hearing from only one student in the classroom). In taking a systems approach to change, useful data for determining need and success are generated in intentional and inclusive ways, from across diverse stakeholder groups. Identifying stakeholder groups happens when change agents decenter, look outside themselves and the voices of the status quo.

### 5.1. *Data-driven Interconnecting and Decentering*

It is clear from the variety in approach and the nature of successes reported in the case studies that data-driven decision-making is valuable. This mirrors what is known from the research and practice literatures (Laursen, 2019). The case studies have done a good job with rationales prompted by student data. However, more can be learned through the kinds of data that are largely absent in the case reports: instructor and implementation data. Now is the time to interconnect across contexts. Examples of how to do that as part of the professional development of novice instructors include activities in which novices and/or providers:

- Gather observational data about the nature of classroom questions and answers (e.g., Stone-Johnstone et al. (CS 25 ) used EQUIP).
- Data-mine learning management systems for evidence of equitable and inclusive instruction (e.g., an audit of time/contributions to discussion fora broken down by student demographics, or a review with feedback to instructors of course sites in learning management systems using a rubric; Baldwin et al., 2018).
- Conduct surveys of instructors about uses of various practices and instructor interactional experiences of teaching – including experiences of racialized or sexualized or gendered interactions (Sue, et al., 2011); repeat these types of data gathering in and through the professional development workshops/courses/experiences, where instructors are the learners.

Some tools for such data collection exist already (see the Resources Appendix). However, next steps for every department include moving into data gathering at the broader levels in Figure 1 and a group-level decentering that includes identifying and inviting outside experts to support connected knowledge growth about each other as thinkers and doers of instruction (Reinholz et al., 2020b). Such expertise will help programs create professional learning opportunities that are informed by research on equitable and just teaching development (including ideas learned in other disciplines, e.g., in biology, Gormally et al., 2016).

### 5.2. *Generating the Future*

Important for developing the next generation of instructional change agents is that an effective agent of change has self-aware communication and advocacy skills. Change or "reform" efforts can fail because people do not know about, or do not know to, pay attention to interconnected-ness and the complexity that comes from decentering. For example, it can be fatal to a change effort to focus attention on one thing (e.g., "all that is needed to reform how this course is taught is changing the textbook" or "student-centered instruction will fix the problem").

Just as mathematical expertise is needed to identify and describe what mathematical ideas to focus on (Figure 1, MATH) and additional expertise in how people think and learn must be added to create effective instruction (Figure 1, INSTRUCTION), yet another type of expertise must be tapped for professional learning about college instruction that is generative (Figure 1, PROFESSIONAL DEVELOPMENT). That expertise includes developing, modeling, and practicing the use of decentering and interconnecting in the design, delivery, and evaluation of professional development. While mathematics departments include experts in MATH, and many have access to experts in INSTRUCTION (e.g., those knowledgeable about research and development in post-secondary curriculum design and instruction), few mathematics departments already have members with expertise in decentering and interconnecting, associated intercultural competence development, or modeling these, related to how adults gain professional knowledge about equitably and inclusively teaching a diversity of other adults.

As the cases demonstrate, many members of the undergraduate mathematics community have taken on roles as change agents, responding to needs and shaping efforts to improve teaching and learning. There may be additional benefits to the community if, in parallel with these efforts, we aim to develop and refine theory. The history of mathematics education includes many examples of the productive interplay between empirical research and theory development (e.g., studies of problem-solving and theory about meta-cognition, studies of teachers and theory about mathematical knowledge for teaching). Having additional theory-development around the characteristics, roles, knowledge, and skills of effective change agents would certainly be welcomed. That can further inform the theory-building efforts we have begun to lay out in this report.

We have asserted that decentering and interconnecting are key to success as a change agent. We based this assertion on examination of what organizational change entails (e.g., Deszca et al., 2019; Kotter, 2012). The research and practice communities will likely benefit if these claims are investigated in conjunction with organizational change efforts. Parallel efforts to create theory and enact change can lead to the accumulation of knowledge which can, in turn, inform the next cycle of efforts.

As is true in the broader literature, the case stories demonstrate the importance of collaboration by a group of change agents, not all of whom are mathematicians (Laursen & Austin, 2020; McShannon & Hynes, 2005; Saichai & Theisen, 2020; Theobald et al., 2020). A corollary of acquiring skill at decentering is that it prepares one to be a participant in a collective effort (e.g., with a classroom full of students, with a department full of colleagues, with a cross-professions team). With attention to decentering and interconnecting, the next generation of change agents will be equipped to participate in collective action.

### 6. *Summary*

The case studies provide examples of how today's change agents addressed issues of diversity, equity, and inclusion by utilizing multiple strands of knowledge. At the classroom level, for example, the work entailed (among other things) making design decisions for curriculum and instruction that took into consideration the mathematics to be learned, the learning opportunities that were well-suited to that mathematics, the people doing the learning, and how these were connected and informed one another. The case stories describe many teaching strategies that might be learned in professional development. These include using learner names (Oliver et al. (CS 21 ).), learner-participation in decisions in and for classroom activity (Soto et al. (CS 22 )), eliciting and celebrating diversity in strategies (Liou-Mark et al. (CS 16 ); Jensen-Vallin (CS 13 ); Villa et al. (CS 27 )), providing structure within which learners have choices about what and how to focus intellectual effort (Johnson et al. (CS 14 ); Bennett et al. (CS 3 )), and purposeful work to notice and respond to interconnections consequential in equity and inclusion (Stone-Johnstone et al. (CS 25 )). Similarly, change at the course or department level was more apt to be successful if approached with information in mind about the different subsystems involved and with knowledge of how those subsystems connected to and interacted with one another.

The challenges and opportunities that the mathematics community will face in the future are unknown. However, the community can approach professional development for novice instructors in ways that build capacity in decentered information-gathering (e.g., about student thinking) and in engaging in the intellectual work of using information to understand connections among subsystems. The instructional practices of eliciting student thinking (to generate information) and making decisions based on that (interconnecting) are a mini-version or microcosm of what is needed to design and enact change at any level of the educational system.

Professional development for novice instructors is an opportunity to develop data-driven cycles of awareness, use, and monitoring of success for decentering and interconnecting. Also, it is an opportunity for novices like graduate students to imagine and build purposefully towards broad future-self goals – addressing questions like "How will I define my professional success as an instructor? How will instructional success be connected to my other professional goals for research and contributions to my professional community? How will I measure that and be accountable for it?" These purposes are served when activities in professional learning environments go beyond readings and presentations about teaching, when the environment is itself clearly and explicitly modeling target practices and demonstrating how to assess progress in ways that leverage decentering and interconnecting. Professional development that engages novices in building skill at self-sustaining, generative change as professionals is the ground in which agency for change is seeded and nurtured.

## 7. *Acknowledgements*

This material is based upon work supported by the National Science Foundation under Grant Nos. (DUE-142381, DUE-1625215, and DUE-2021139). Any opinions, findings, and conclusions or recommendations expressed in this material are those of the authors and do not necessarily reflect the views of the National Science Foundation.

## 9. *Appendix*

The resources in this appendix are offered to the reader as a sampler of currently available supports for monitoring progress (evaluation), curricula for use in professional development of novice instructors, and DEI-focused decentering and interconnecting.

### 9.1. *Evaluation: Monitoring Progress Towards Goals in Professional Development*

Four interconnecting activities underpin successful strategies to transform instruction (e.g., professional development offered to novice college teachers; Laursen & Austin, 2020; pp. 201-217). Each of the steps calls for decentering and interconnecting:

1. Identify and frame challenges in consultation with others (e.g., with stakeholders such as undergraduate tutors, novice instructors, as well as experts in noticing and addressing difference who may not come from within the department).
2. Analyze the context in consultation with others (e.g., make sure the picture of the contextual landscape is detailed and attends to justice, equity, diversity, and inclusion).
3. Choose strategies and interventions that are responsive to the context. This involves two pieces:
   (a) including strategies that have demonstrated success (discussed at length in the chapters in this book and the collegiate mathematics education literature) and
   (b) creating a relevant change model (Reinholz & Andrews, 2020). As providers develop a change model, they need to ask themselves these questions:
      - What are the overall goals that the program seeks to achieve?
      - What strategies or interventions do we (providers) believe will help to achieve each goal? Why do we believe these strategies will work? How do various strategies complement each other? How might they conflict? How do we resolve values conflicts?
      - How will we know if we have been successful? What short-term and longer-term indicators will show progress toward a goal?
4. Implement the plan and monitor progress in a feedback loop for continued improvement.

Additional ideas related to monitoring progress can be found in the American Association for the Advancement of Science (2013) report on describing and measuring STEM teaching practices.

9.2. *Curricula for Professional Development of Novice College Mathematics Instructors*

The College Mathematics Instructor Development Source (CoMInDS; 2021) is an NSF-funded project done in collaboration with the Mathematical Association of America. Since its inception in 2014, CoMInDS has offered instructional resources and professional learning opportunities to those who provide teaching-focused professional development to novice college mathematics instructors. The collection of resources (accessible via connect.maa.org) include activities and assignments to help novice instructors develop knowledge and skills for teaching mathematics. The collection includes materials that providers can use to help novice instructors gain facility with the types of practices featured in this volume, including engaging students, eliciting students' many ways of thinking, facilitating collaborative group work, and having courageous conversations about race, equity, justice, diversity, and inclusion. Also included are providers' syllabi for teaching seminars and profiles of novice instructor development programs from a variety of different institutions. New resources are regularly added to the CoMInDS collection.

9.3. *Decentering and Interconnecting*

Specific resources for expertise for decentering and/or connecting at the departmental, institutional, regional, and national levels include:

- The Professional and Organizational Development (POD) Network in Higher Education
  https://podnetwork.org/
- Aspire: The National Alliance for Inclusive and Diverse STEM Faculty
  https://www.aspirealliance.org
- CIRTL: The Center for the Integration of Research, Teaching and Learning
  https://www.cirtl.net/
- Intercultural competence development. Contact the Intercultural Development Inventory group with a brief description of a project (e.g., workshop, professional learning community, department development) and they will send recommendations for licensed Qualified Administrator consultants suited to the project.
  https://idiinventory.com/about-us/

Specific readings that share the perspectives of people of color include:

- Racial Dialogues: Challenges Faculty of Color Face in the Classroom (Sue et al., 2011).
- Inside Higher Ed's series on The Black Experience in Higher Education (e.g., The Souls of Black Professors, Flaherty, 2020).
- Anti-Racism Training, recordings and list of consultants (Academics for Black Survival and Wellness, 2020).
- Racelighting: A Prevalent Version of Gaslighting Facing People of Color (Wood & Harris, 2021).

**10. *Take Away Messages***

- Professional development has to be part of whatever the endeavor is unless all you are trying to do is reproduce what has already been happening (i.e., suppress change).
- The chapter proposes attention to two processes as useful: decentering and interconnecting.
    - Instructional decentering is, at its most basic, the act of seeing from someone else's point of view and has historical roots in the work of Piaget (1955). It means engaging with other people as a participant in interaction, rather than as the center of interaction.
    - While decentering is awareness from within the perspectives of others, interconnecting uses meta-awareness to make connections across perspectives and contexts. The example in the chapter is anchored in the interconnections among mathematics, classroom instruction, and professional development.
- Professional development for novice instructors is an opportunity to develop data-driven cycles of awareness, use, and monitoring of success for decentering and interconnecting. Also, it is an opportunity for novices like graduate students to imagine and build purposefully towards broad future-self goals – addressing questions like "How will I define my professional success as an instructor? How will instructional success be connected to my other professional goals for research and contributions to my professional community? How will I measure that and be accountable for it?" These purposes are served when activities in professional learning environments go beyond readings and presentations about teaching, when the environment is itself clearly and explicitly modeling target practices and demonstrating how to assess progress in ways that leverage decentering and interconnecting.
- Save some professional energy. Designing a program and choosing curriculum for novice instructor development is an important part of the work. Equally important are modeling equitable and inclusive practices with instructors-as-learners during professional development and situating the professional learning in the larger systems of course, department, university, and communities.
- Novices need opportunities to try, stumble, and learn how to keep moving and keep learning in and through teaching when the ground under them changes.

# 3 Three Frameworks for Advanced and Theoretically-Driven Analyses of Classroom Practices

**Brian Katz**
*California State University – Long Beach*

Brian P Katz (BK) is faculty in Mathematics Education at CSULB. BK is passionate about interactions between inquiry, epistemology, identity, authority, and justice as both a scholar and teacher, especially in the context of preparing teachers to lead student-centered and rehumanizing classrooms of their own. BK is Associate and Communications Editor with PRIMUS, a former Chair of IBL SIGMAA, a current member of the executive committee of SIGMAA RUME, an MAA textbook author, co-Editor of the AMS inclusion/exclusion blog, Chair of the MAA SoCal-NV Section, and part of the design and facilitation PRODUCT professional development team. BK supports two of the best cats in the world and loves to sing.

This volume contains a large collection of case studies from mathematics educators about efforts to make their work and institutions more diverse, equitable, inclusive, and just. My goal with this chapter is to help you make sense of and begin to leverage these case studies, with a focus on classroom practices. For the purposes of this chapter, I will use classroom practices as it is used in the MAA Instructional Practices Guide (Abell et al., 2018), as part of a trio with Design Practices and Assessment Practices making up Instructional Practices. These three collections of practices are interrelated and overlapping, but I will focus on the classroom practices as the live interactions with students because design and assessment themes appear strongly in other cross-cutting theme chapters.

I plan to approach this goal in two ways. First, I will summarize three theoretical frameworks for thinking about equity and justice in mathematics classrooms. I believe that these frameworks will help you be more effective in stepping back from the particular details within individual case studies as you think about transferability to your context, and I believe that these frameworks will supporting you in reading these case studies critically, asking about the implicit assumptions and wondering about what might be missing. In each framework section, I will make connections to illuminating examples of the framework ideas in the collection of case studies. These frameworks are explicitly cited in several of the case studies as motivation and theoretical framing, but this chapter will develop them more explicitly as frameworks. That said, my overview of these frameworks cannot replace reading the original works, but I hope to communicate a heuristic understanding. Second, I will engage more directly the case studies that bring classroom practices to the foreground by discussing common themes and by analyzing a small number as in-depth case study analyses.

This chapter is focused on the classroom level leading inquiry: *How do I engage in the work of supporting justice-making and change in my classroom?*

## 1. *Frameworks*

Each of these sections is focused on a framework from a core reference, by Gutiérrez (2017), Gutiérrez (2019), and Boaler (2006) respectively. I will summarize them and connect them to the case studies in this volume, but I also hope you engage the original sources, so I have separated these references out from the other supporting references at the end of this chapter.

The first framework (Section 1.1) focuses on a bridge between the broader societal context and instructor choices in the classroom. The second framework (Section 1.2) will offer a vision of mathematics as we hope students will experience it in our courses. And the third framework (Section 1.3) will describe some teacher facilitation moves and classroom structures that support this work, though there are also moves and structures embedded in the other two frameworks as well. Rather than framing this as getting progressively more concrete, I suggest that you view context (Section 1.1), vision and values (Section 1.2), and strategies (Section 1.3) as dimensions of every decision you make when working for more equitable classroom practices. However, if you are feeling overwhelmed while reading, you could read the subsections in reverse order or start by reading the case study analyses (Section 2) or the summary table (Section 3.1) as examples first before returning to the frameworks in general.

### 1.1. *Context: Nepantla and Political Conocimiento (Gutiérrez)*

Classroom practices are part of a larger educational and cultural system, and this first pair of frameworks developed by Rochelle Gutiérrez help us interpret equity and justice work in this context. This section is derived from her paper "Political Conocimiento for Teaching Mathematics: Why teachers need it and how to develop it" (Gutiérrez, 2017).

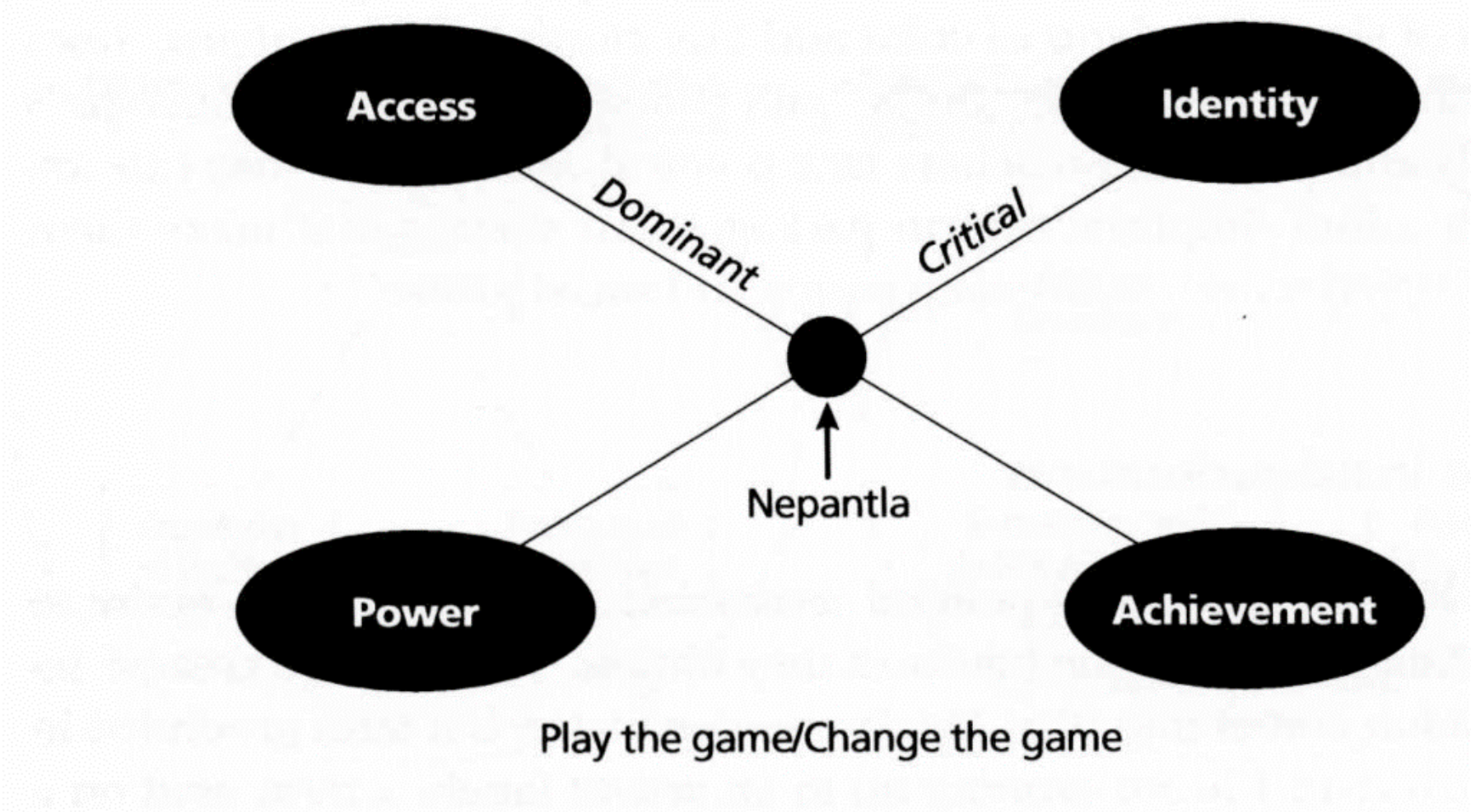

**Figure 2.2** Dimensions of equity/learning.

Gutiérrez challenges us to develop a careful and sophisticated definition of equity, especially in light of the many ways that this term gets co-opted and demonized in our society. She points out that many discussions of equity in mathematics classrooms focus on themes of *Access*, who is present in these spaces, and *Achievement*, how are those present succeeding in these spaces. These two ideas are already complex. Discussions of access can use metaphors of pipelines, suggesting that students are leaking out of the system, rather than being forced out, and discussions of preparation can focus blame on individual students while ignoring the systems that failed to offer rich learning experiences to some students. Gutiérrez refers to these two dimensions as the dominant axis in this model because they center the powerful educational systems as they are. Discussions centering these dimensions often leave topics such as the nature and value of assessment unexamined, because educators from this powerful system have historically focused on these dimensions, and because (as a result) these dimensions align with preparing students to participate economically in society (Gutiérrez, 2009).

Gutiérrez introduces two other dimensions for discussion of equity: *Identity*, students' pasts, their cultural and linguistic resources, and how they see themselves and are seen by others, and *Power*, students' agency in making choices and changes. She refers to these two dimensions as the critical axis because they center student perspectives, often offer counter-narratives to dominant framings, and help build critical citizens who can "change the game" (Gutiérrez, 2009). In particular, Gutiérrez advocates playing the game of helping students succeed as framed by the dominant axis while also changing the game in ways that are guided by the critical axis. Gutiérrez labels the intersection of these two axes with the term "Nepantla" to indicate the way educators live with the tensions between these axes, walking the "middle path" with a foot in each. Note that the term "axis" here is not meant to imply that one must choose between power and identity or between access and achievement, nor is it meant to construct a two-dimensional space. Instead, I suggest you view each axis as a lens made from a strongly connected pair of aspects of justice work.

Many of the case studies in this volume focus on access and achievement. We certainly don't want students to be driven out of mathematics, but I encourage you to look deeper than this goal, which can treat students from privileged groups as neutral and invisible while focusing on perceived deficits in other groups. How are we pushing back against narrow definitions of achievement? How are we valuing the strengths that students bring with them? How are we helping students achieve the goals they value?

Several case studies explicitly engage aspects of identity. Cross-cutting CHAPTER 6 discusses the ways that classes incorporate students' skills with multiple languages as resources for learning. Multiple case studies attend to potential stereotype threat (e.g., Soto et al., (CS 22 )), whereby students can feel that their value and identity is threatened if they are seen as a member of a marginalized group who is not excelling in math. Kilty et al. (CS 15 ) explicitly names the challenge that many students in Calculus have to manage perceptions that "'all' the other students seem to 'already know' Calculus (although, often, this "knowing" of Calculus can be more of a perception, than a reality)". And multiple case studies incorporate discussion of growth mindsets. I believe that growth mindset interventions can help students identify as learners and see themselves as having more power and agency, but I am also sensitive to the ways that growth mindset can suggest the responsibility for learning rests solely on the student's mindset, leading to individual students being blamed for not learning while their teachers ignore structural forces working against the student's learning (TODOS, 2020; Gross-Loh, 2016). As Gutiérrez has pointed out, this shallow version of the mindset discourse is a new version of the toxic American bootstrap myth. I think DiGregorio et al. (CS 10 ) is a nice example of a productive use of growth mindset; this case study discusses a revision to a placement exam and procedure that was previously encouraging a harmful and dehumanizing fixed mindset.

Johnson et al. (CS 14 ) explicitly engages the dimension of power. This case study starts from the observation that a reliance on "compliant answer finding" locates all classroom power with instructors. Instead, these authors "argue that strategic implementation of mathematical tasks to promote reasoning, rather than answer finding, can serve as a catalyst to redistribute power in the classroom, to afford students more agency in their mathematics learning".

Gutiérrez goes further to develop the concept of *political conocimiento*, which is part of teachers' knowledge with students and communities. Older models of teacher knowledge emphasize content knowledge, pedagogical knowledge, and knowledge of students. Political conocimiento, a term that intentionally uses the Spanish verb for knowing that emphasizes knowledge through connections rather than abstracted knowledge, replaces knowledge "of" students with knowledge "with" students and communities and adds political knowledge as a key component of teacher knowledge. This term thus means teachers' political knowledge with students and communities. This framework goes further to explicitly center the students who are framed as "others" in many dominant narratives and to acknowledge the relationships between educational systems and the histories of the people involved.

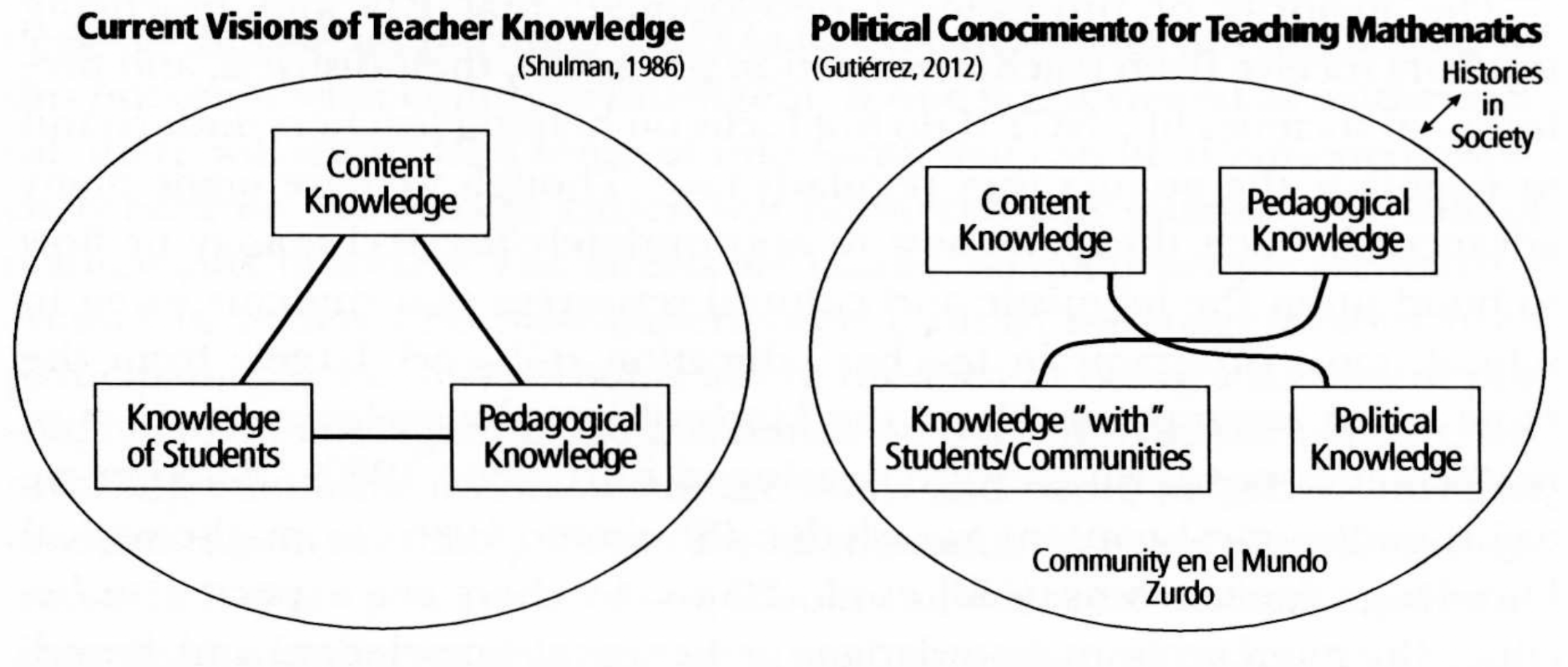

**Figure 2.1** Teacher knowledge.

Gutiérrez argues that teachers need political conocimiento because all mathematics teaching is political, especially teaching that is trying to change the game. Teachers need to be able to push back when they are pressured to do inequitable and dehumanizing things, and they need to find ways to be creatively insubordinate that align with goals for equity and justice as they manage demands made of them by the system.

Many of the case studies discuss reform within a system, especially the recent movement to remove courses that do not carry college credit as barriers for students. I will leave a deeper discussion of this theme to the cross-cutting CHAPTER 5 "*Course Redesign: Pathways towards transformation*" on pathways and course redesign, but I encourage you to look at these through the lens of questions asked from political conocimiento. Which opportunities were taken to be creative within systems and which were not? How was shared knowledge of the challenges developed with students? Change efforts are strategic and involve leveraging limited resources; I am not suggesting that compromises are a problem, but I do believe we need to look carefully for patterns in the impacts of those compromises.

Several case studies hint at potential resistance to change, either through institutional inertia or through adherence to a lecture paradigm. It is now clear that faculty hold many strong beliefs that student-centered pedagogies and a focus on equity are at odds with content goals, but it is also clear that these fears are not realized (e.g., Johnson et al., 2018; Rasmussen & Kwon, 2007;

Yoshinobu & Jones, 2012). Similarly, it is also clear that students often report learning less when their instructors use active learning approaches while they are actually learning more (Deslauriers et al., 2019). These factors can make it dangerous for some instructors to pursue reform efforts, a danger that is often built into faculty evaluation systems, which amplifies these dangers for faculty without security of employment or who do not present with the most privileged identities. And while I do not see evidence in the case studies that the authors experienced retaliation from within their institutions or outside, retaliation against reform is growing in terms of frequency and ferocity (Gutiérrez, 2017b). I encourage you to consider how the authors of the case studies address the challenges of making change. Finally, some case studies bring explicit critiques of systems into the classroom with students. For example, Oliver et al. (CS 21 ) includes a DEI Statement that acknowledges identities and identifies policies that might seem fair that nevertheless lead to inequities.

### 1.1. *Vision and Values: Rehumanizing Mathematics (Gutiérrez)*

This section focuses on a different framework, also developed by Rochelle Gutiérrez, called *Rehumanizing Mathematics*. I encourage you to listen to Gutiérrez describe these ideas in her own words; the core reference for this section (Gutiérrez, 2019) is a presentation she gave to a higher education mathematics audience at the Park City Mathematics Institute (http://projects.ias.edu/pcmi/hstp/sum2019/photos/gutierrez/) that builds on an earlier paper (Gutiérrez, 2018).

Gutiérrez is very careful with the wording in the name of this framework. First, she points out that mathematics is already and always has been a human activity, so while many of us have experienced dehumanization through mathematics, the work is REhumanizing it, not humanizing it from an inhuman starting point. And second, she points out that in choosing the present progressive participle rehumanizING, she is asserting that this is an active and ongoing process; mathematics will never be permanently rehumanizED.

The rehumanizing mathematics framework offers eight dimensions of the work. These dimensions interact and overlap; their ordering in the diagram below is not significant, and the dimensions are more and less salient in various applications of this framework. Furthermore, Gutiérrez does not claim that these eight dimensions are a complete taxonomy of what is needed. Here are short interpretations of these dimensions as questions for the listener/reader with references to examples from the collection of case studies.

*Positioning/Participation:* How does membership and participation in this classroom look? What status differences and hierarchies exist, and how do they relate to those from the broader culture? How do teachers position students as legitimate and competent members of the class and disciplinary communities?

Several case studies describe students developing class norms or group commitments (e.g., Byrne (CS 5 )); I think this helps position students with agency while also making the modes of participation less opaque. Similarly, several case studies include classroom interventions explicitly intended to help students see themselves as members of a STEM community and doers of math/science (e.g., Wagnon et al. (CS 29 )), though many are less explicit about positioning students as knowers.

*Cultures/Histories:* How are students' pasts valued and leveraged in this classroom? What resources and funds of knowledge do students bring with them to the classroom? How are various cultural approaches to mathematics incorporated?

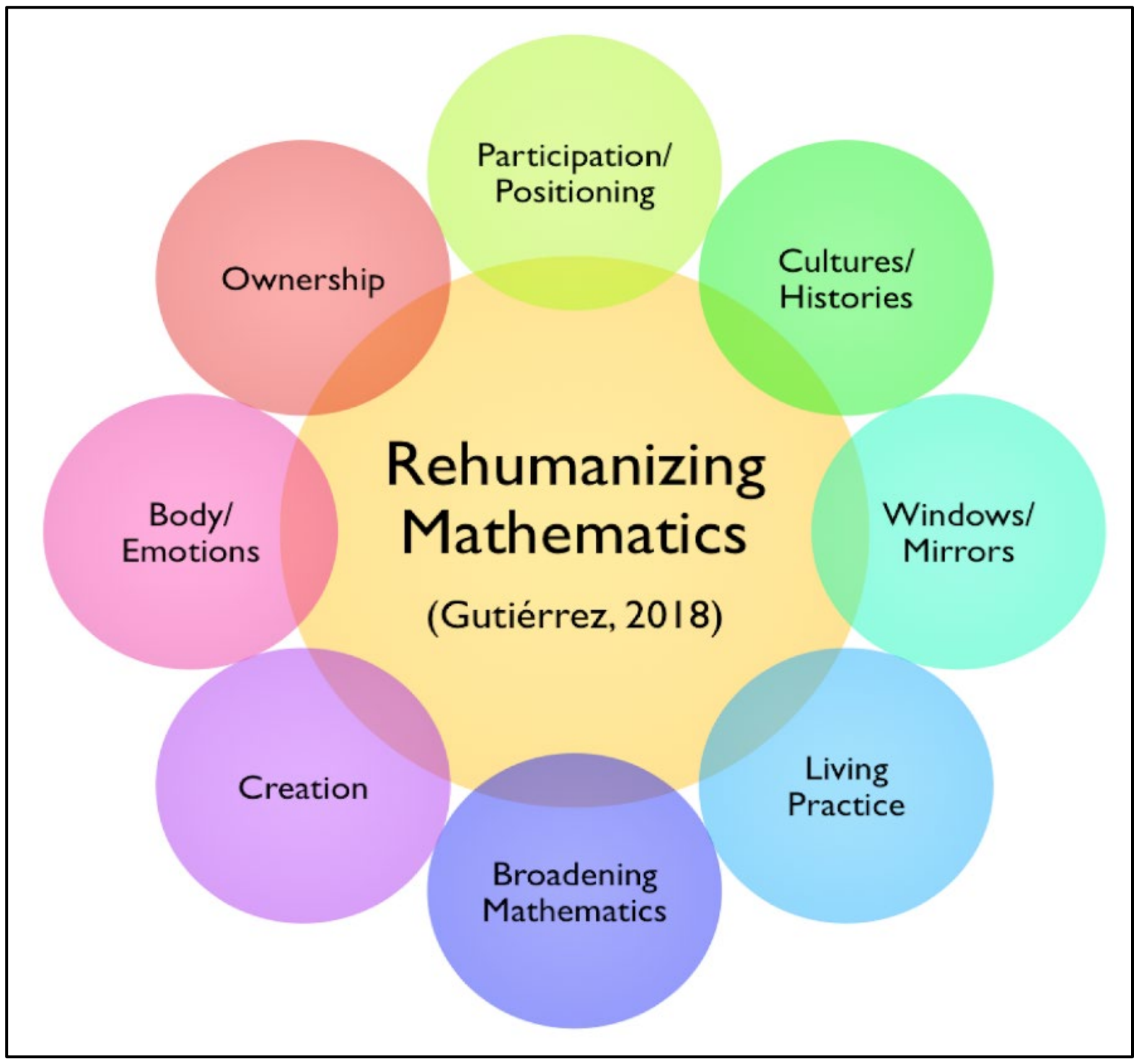

*Figure 2.3 Rehumanizing Mathematics Diagram Reproduced from (Gutierrez, 2018).*

All of the case studies are pushing back against a cultural narrative that some students are "bad at math" and need to be weeded or tracked out of collegiate math classrooms (e.g., Stacy (CS 23 )). The case studies that focus on multilingual students are doing important work in this dimension by framing students' cultures as strengths (see CHAPTER 6 "Leveraging Identity and Language to Promote DEI in Undergraduate Mathematics"). This collection of case studies is much more likely to critique western academic traditions (Gutiérrez, 2017c; Brown, 2015) than to center other ways of doing mathematics. Many of these case studies do work on enculturating students into STEM rather than changing that culture (e.g., Oliver et al. (CS 21 )).

*Windows/Mirrors:* How do students see themselves and others represented in the curriculum? How is mathematics used as a path toward understanding their own experiences and developing empathy for the experiences of their peers, especially the experiences of those historically minoritized in mathematics classrooms?

In terms of course content, this dimension could include moving away from tasks that treat human labor as an abstracted commodity (as is often done in optimization problems) and moving toward tasks that help students seek safety and wellness for all or grapple with issues such as voter suppression. Several of the case studies incorporate messy projects that have potential to focus on issues of local importance to students (e.g., Zobitz et al. (CS 30 )). In terms of processes, many of these case studies do advocate a shift away from forced consensus and toward dignity for students. These elements allow students to bring more of themselves into the classroom, but as a discipline we have a lot more work to do to make mathematics a place in which we learn about ourselves and other and expand our empathy.

*Living Practice:* In what ways is mathematics experienced as an ongoing, human endeavor?

Many students come to college having only experienced mathematics as an inert body of knowledge, not a field with questions that is always growing (e.g., Schoenfeld, 1988). The case studies that incorporate messy application projects are helping students see that there are new questions to ask and that modeling involves lots of ongoing human choices. Several other case studies mention teaching with inquiry. While this is not explicit in the case studies, for me, the term inquiry refers to a classroom discourse that is guided by questions that are being asked (or at least made sense of) from the learners' perspectives. In other words, I think that if executed well, inquiry experiences can help students view mathematical knowledge as something built by humans in response to their questions.

*Broadening Mathematics:* How does this classroom incorporate aspects of mathematics outside of abstract symbolism, emphasis on the general cases, and the narrow geometry/algebra/calculus pathway? How are qualitative forms of mathematics valued?

This is a challenging dimension for a volume focused on calculus courses, because that focus centers this narrow pathway. I believe that we, as a disciplinary community, need to challenge ourselves to work on broadening mathematics in calculus courses because so many students enter and leave these classes thinking that understanding calculus means that they can compute derivatives quickly. Many of these case studies do incorporate group and team learning explicitly in ways that broaden mathematics. These collaborative settings challenge the vision of mathematics as work done in isolation by a "genius", who produces something symbolic fully formed. Instead, students are able to see that mathematics is done in community and often starts with qualitative approaches focused on meaning.

*Creation:* How are students in this classroom engaged in inventing ideas that are new (at least to them)? What new forms of mathematics are likely to appear in this classroom?

In contrast with the dimension about histories, which seeks to incorporate more diverse forms of mathematics that already exist, this dimension seeks to support students in generating new mathematics (and/or new-to-them mathematics). All learning is in some sense an act of creation of structure and connections, but this dimension challenges us to make sure that this learning goes beyond asking students to rebuild others' thinking. I think this dimension has a reciprocal relationship with viewing mathematics as a living practice because generating new mathematics demonstrates that it is living, while believing that it is a living practice may be important for students to understand when we encourage them to create. As a result, the clearest example of this dimension in the case studies are also those modeling projects discussed above; for more ideas, see the GAIMME second edition report (Garfunkel et al., 2016).

*Body/Emotions:* How does this classroom acknowledge and leverage the fact that mathematics is done by people, with bodies and emotions?

Many of the case studies (e.g., Chang & Chen. (CS 7 )) acknowledge the strong emotions of shame and anxiety that many students feel in mathematics classrooms, especially students also managing stereotype threat. This dimension also encourages us to attend to the ways that humans think with their bodies, a perspective called embodied cognition, and to value intuition and emotions as resources for doing mathematics rather than as barriers or distractions. Some of the case study authors are known for their other writings about embodied cognition, even if it is not explicit in many case studies. For more, see Soto (2018) and Barnes and Libertini (2013, 2018).

*Ownership:* In what ways does this classroom make space for students to do mathematics for their own purposes? How are students' joy and desire to play supported?

Many of the case studies focus on student populations that are intending to pursue work in STEM, so it is reasonable to see their calculus course broadly as supporting their goals, but we know that many of these students see calculus as disconnected from those goals and as a barrier between them and the path they hope to follow (e.g., Ellis et al., 2016). The case studies that incorporate student choice for projects are making some space for this dimension, and the case studies that use student mathematical autobiographies (see CHAPTER 4 ) are communicating to students that their goals matter in the course environment.

### 1.2. *Strategies: Complex Instruction+ (Boaler)*

This section describes our final framework for thinking about equitable classroom practices. The ideas are based on a paper by Jo Boaler (2006), and they build on an earlier framework called complex instruction, so I will refer to this third framework as *Complex Instruction+*. The two previous frameworks emphasize the context and goals for classroom practices; this framework emphasizes smaller-scale teacher moves. Canner et al. (CS 6 ) discusses complex instruction explicitly.

*Multidimensionality:* How does this classroom incorporate multiple ways to be successful?

Multidimensional classrooms acknowledge that mathematics is a complex, community activity that requires multiple skills and perspectives, and that it is expected that different members of that community contribute differently to the overall success. Case studies that valued open-ended tasks without a singular correct solution as well as those that explicitly value communication and collaboration as mathematical skills are supporting multidimensionality.

*Roles:* How is group work structured? How does the work get distributed so that each learner has something important to do?

Students rarely come to our classrooms knowing how to collaborate effectively. Furthermore, we know that unstructured group work can devolve into unproductive learning environments and can lead to marginalizing experiences for students (e.g., Johnson et al., 2020). Several case studies mention roles as part of their strategy, though they do not go into detail about the roles. Several other case studies incorporate peer learning support; while this is a role between learner and expert (e.g., Davis (CS 8 )), I think readers will find it productive to read those case studies through this lens.

*Assigning Competence:* How is the competence of students made visible in the classroom?

This dimension is related to the ways that students are positioned in the classroom from the rehumanizing mathematics framework. Teachers are doing important work when they make explicit statements about the value of students' ideas because it communicates information about productive ways to do mathematics and because it can counteract the ways that some students are not positioned as competent. I think the term “assigning” is problematic because it suggests that the teacher is arbitrarily labeling a student as competent. Instead, teachers should illuminate students’ competencies and should engineer situations in which those competencies come to the surface.

*Student Responsibility:* In what ways are the students responsible for each other and collective success?

Boaler gives two examples. First, when teachers ask students to revoice and explain peer ideas in their own words, teachers are supporting the notion that members of the classroom community are responsible for trying to understand how their peers are thinking and for trying to communicate in a way that makes sense to those peers. Similarly, teachers can ask a student to report for a group, leaving and returning to the group if the student needed group support; this communicates that groups should progress by making sure that everyone shares understanding rather than moving on when only some group members understand. This kind of responsibility is not clear in the case studies, likely because of the grain-size of these pieces of writing. Second, students can be more formally responsible to each other in group work with group assessments. Deshler et al. (CS 9 ) discusses a course with formal groups leading to student responsibility, and MacAurthur (CS 17 ) discusses an effort to rehumanize group exams. Akin et al. (CS 1 ) discusses a situation in which a grading system encouraged competition rather than shared responsibility.

*Justification/Reasoning:* How do students participate in the work of reasoning and justifying mathematical assertions?

Not only does this dimension position students as knowers of mathematics, focusing on justification and reasoning makes the classroom space more multidimensional and democratizes access to the work of mathematics (Buell et al., 2017), shifting away from a focus on which individuals do or do not already hold certain knowledge and toward collectively deciding if the support for a claim is convincing. Johnson et al. (CS 14 ) is about moving away from answer-seeking and Liou-Mark et al. (CS 16 ) about moving peer tutors away from answer-giving. Mingus et al. (CS 20 ) includes an intervention that helps students understand the qualities of a justification.

*High Expectations:* How does this course communicate the belief that all students can and should achieve success with challenging tasks?

At a structural level, the case studies are clear that the goal is to reimagine their courses so that all students are supported in meeting high standards and that there is no tension between equity and these standards. These beliefs are critical for faculty to continue offering rich learning experiences for students. And while we don’t see much about the impact of these beliefs on individual students in these case studies, students can clearly tell the difference between environments that seek to support their brilliance (Martin, 2019) and those that assume they will fail. Zobitz et al. (CS 30 ) contains explicit discussion of high expectations.

*Effort over Ability:* How is the work of learning valued and viewed as changeable rather than innate?

Stacy (CS 23 ) and others push back against the harmful narrative that some people are just "bad at math". This dimension connects strongly with growth mindset, but I will discuss that more in the final dimension.

*Learning Practices:* How are productive practices for learning made visible and named in this classroom?

I believe that every course needs to help students learn about how teaching and learning work. Discussions of growth mindset are often used to open discussions about learning, such as in Byrne (CS 5 ). Other chapters including Mawhinney et al. (CS 18 ) and Wagnon et al. (CS 29 ) focus on self-regulated learning, an aspect of metacognition. Fuller et al. (CS 11 ) connects justification and learning practices as a core tenet of the course design. In addition to general principles about learning, I expect that many of these courses include strategies specific to mathematics, such as checking a claim with a familiar example, generating examples and non-examples, adding assumptions to help get started, and translating ideas between representations to consider them differently.

Boaler concludes with a powerful observation. The students experiencing this kind of instruction were learning to treat each other in respectful ways through these mathematics teaching practices. In contrast to the myth that mathematics is independent of equity, this framework asserts that mathematics has particular ways to help students develop relational equity skills.

## 2. *Case Study Analyses*

My assertion at the start of this chapter was that the frameworks above would help you analyze and make meaning with the case studies. Now I'll model that with four short analyses of individual case studies.

### 2.1. *"The Effective Thinking Calculus Project at UT Austin"*

Starbird et al. (CS 24 ) describes a pair of linked courses, one Calculus I and one entitled "Effective Thinking" focused on having explicit discussions about thinking strategies.

The initial motivation for this course focuses on access to higher education in Texas and the achievement of the students who do have access, and the case study concludes by reporting encouraging achievement outcomes in terms of completion rates in this course and grades in subsequent courses. The body of this case study asserts that all students "can view mathematics as part of who they are", which I interpret as an identity goal. This pair of courses represents a creative approach to working within Texas law and to leveraging existing campus structures to support students.

The ideals in this innovation align well with the rehumanizing mathematics dimensions of creation and living practice, as the students are asked to participate in the production of mathematics, though students might experience this as recreating others' mathematics. That said, the goal of asking students to tell the story of the creation of an idea is very interesting. I think this could be quite powerful for seeing mathematics as a living discipline and possibly a sense of ownership, but I think there's also a chance that the emphasis on these concepts being natural might undermine that goal, depending on how it's implemented.

The strongest alignment between this case study and the three frameworks is the focus on learning practices. The point of the paired course is to make space for explicit discussions of effective thinking and practices, and I think that finding a way to make a whole course's worth of space that supports the other learning goals is a really exciting piece in this case study. The course also includes student responsibility structures, but the focus is on higher-stakes assignments rather than on learning tasks distributed across the whole course environment.

### 2.2. *"I Wish I Could Say 'You Should Not be Here'": An Analysis of Instructors' and Students' Contrasting Perceptions of a Racialized and Gendered Gatekeeping Practice in Calculus*

McNeill et al. (CS 19 ) describes interviews with students and faculty around a generalized episode from a Calculus II course in which the instructor suggests to the whole class that students who cannot complete some computations quickly should consider dropping back a course or not taking Calculus II. The instructors approached this episode with assumptions about mathematics teaching being socially neutral and hence read this episode as helpful and supportive, while the students of color felt that this episode made them feel unwelcome in the course or mathematics, more stressed in ways similar to stereotype threat, and discouraged about their learning even when they persist.

This case study discusses the dominant axis of access and achievement in mathematics classrooms in depth, with both the students and instructors concerned about student achievement in the course. The students also bring up the negative impacts on their identities, as being told they did not belong in math or science, and the discouragement as an impact on their sense of power. The instructors take an approach that might be called knowledge of students when they assert that they know best who should or should not be in their courses, while the authors of the case study push for more knowledge with students in the context of histories in society to resist these narratives and center the student experiences and goals with more nuance.

The episode is dehumanizing for the students. They feel positioned as unable to participate, they wonder if this is why they don't see themselves represented in the disciplinary community, they have to manage intense emotions that are not supporting learning, and perhaps they express frustration at a very narrow vision of mathematics as presented by their instructors.

We can also analyze the episode using the complex instruction+ framework. This episode suggests a unidimensional classroom that only values computational speed and not reasoning, and the students certainly feel positioned as incompetent. While the instructors may feel they are expressing high expectations, these expectations appear to be missing the belief that all students can rise to those expectations and the focus on supporting that learning. The students also experience this episode as reinforcing the idea that learning outcomes are innate rather than related to effort, and there is no room in that framing for a discussion of learning practices.

This case study also engages LEADING INQUIRY #1, about instructor critical self-awareness, and I think the larger project from which this case study was taken could be used for powerful professional development interventions in which instructors learn about marginalizing experiences from students' perspectives.

### 2.3. *Creating Space for Student Agency to Support Success in Online Pre-Calculus*

Golden et al. (CS 12 ) describes a support course for an online PreCalculus course, focusing on adaptations made while teaching the course based on feedback that it was not meeting students' needs. Initially, the course was created to address an identified achievement gap for structurally disadvantaged students.

The existing online course structure (as 1-credit, self-paced modules) may have empowered some students with agency and freedom over their educational choices, but for others it could be exacerbating negative emotions about both mathematics and online learning, perhaps by conceptualizing the courses as remediation that could be delivered to individuals without regard to the fact that they are people with bodies and emotions or the fact that some students are positioned as far from competent and disconnected from support resources.

The course leadership designed an in-person supplement to the online modules to build a class community that focused on developing productive mathematical identities and incorporated reflective journal elements that would highlight learning practices. However, the students let the design team know that the course was not meeting its intended goals, in part because they experienced it as putting extra work between them and their achievement goals. The design team learned that their initial plan was developed from the perspective of knowledge of students, rather than with them. They sought student perspectives on the issue at the time that the need for change was identified, and they continued to leverage the journals to gather ongoing insight into the students' needs. The revised course was better aligned with the needs of the students. Moreover, the team learned that they need to continue making decisions with the students in future sections of this course rather than assuming they now had complete knowledge of their new students based on the pilot group.

The initial design of this supplement course approached nepantla, the intersection of the dominant and critical axes, by trying to support students in both their achievement and identity, but students found the initial top-down approach to this goal counterproductive. Starting with the course correction, the focus shifted to making sure students had power and agency in their educational choices, in support of their achievement goals.

### 2.4. *Leveraging Classroom Data to Promote Equitable Instruction*

Stone-Johnstone et al. (CS 25 ) describes the use of a classroom observation tool (EQUIP) to support an instructor's reflection on participation in the classroom. The first author ran a faculty learning community in which the second author participated, in part by recording and analyzing his teaching through the lens of the EQUIP tool, co-created by the third author.

This case study focuses on what the second author learned about himself as a teacher through this reflective analysis. The classroom data showed gaps and other patterns in participation, and the instructor learned that his sense that all students participated approximately equally was not accurate. The data showed that many student contributions in class were short, which contrasted with the instructor's intention and belief that students were participating at high levels with long contributions that explain their thinking. Learning this from the data-driven reflection led him to adopt more intentional strategies to facilitate high-level contributions from students. In particular, this shift in strategies appears to have been coordinated with a shift away from asking students questions for which they are expected to have right answers and toward their strategies.

This case study focuses on how the different kinds of questions instructors ask students lead to patterns of participation. These patterns position students (as a group) in relation to the discipline and the power of the instructor and position individual students as competent or not in the context of the class. These questions also have a strong impact on whether the classroom discourse is focused on reasoning and justification or more inert knowledge. The example transcript shows how these questions can support a multidimensional classroom in which process and effort have value and offers a model for an instructor illuminating a student's competence.

This case study also engages Leading Inquiry #1, about instructor critical self-awareness, which was the focus of the professional learning community that generated this case study.

### 3. *Conclusion*

I'd like to finish this chapter with a few questions for both myself and the readers.

*What is missing?* It's hard to know what you don't see. Theoretical frameworks help articulate what we expect to see, which makes it easier to notice both when expected elements were missing (and patterns therein) and when included elements were not anticipated. I identified some aspects of the frameworks above that are not well represented in these case studies. Stepping back a little further: what kind of conversation about math and social justice is this volume encoding? I see lots of examples of efforts to make the classroom community a more diverse, inclusive, equitable, and perhaps just place through more just pedagogical choices. I also see some examples of efforts to develop mathematics in the context of questions about justice through applications, though this could be a much stronger theme for Calculus courses. I see very little about critiquing the history, philosophy, and culture of the discipline of mathematics as an effort to make mathematics more just. And I see very little about using the context of mathematics to teach justice skills, such as ethics and empathy. These theoretical lenses help me notice these areas for future work.

Flipping this question around, the case studies illuminated a critical theme that is not overt enough in these three frameworks: classroom community. Community is certainly related to knowledge with students, positioning, student responsibility, and other aspects of these frameworks, but many of the authors of these case studies are explicitly trying to build communities in which students feel that they belong, are connected, and are supported. This focus on community is a critical element of the work of improving teaching and learning. Similarly, these frameworks all take an anti-deficit, asset-based stance toward supporting students (Adiredja, 2019), but this axiom is not as clear as I'd like in my summaries. Finally, these frameworks were not developed specifically for the context of Calculus, which is the focus of this volume, so they may not identify the distinctive aspects of these reforms related to Calculus. I encourage readers to consider the role of community, the positioning of students as having assets or deficits, and the specific context of Calculus courses as they engage the case studies.

*How is this cross-cutting chapter connected to the others?* This chapter focused on classroom practices is most closely aligned with LEADING INQUIRY #2: "How do I engage in the work of supporting justice-making and change in my classroom?". However, I think that the analysis and frameworks in this chapter show how the classroom-level work is done on a foundation of critical self-reflection (LEADING INQUIRY #1) and in the context of the educational system (LEADING INQUIRIES #3-6). And just like in the MAA Instructional Practices Guide, classroom practices are intricately linked to design and assessment practices.

*What can I do next?* I hope that this chapter will help you read the case studies in this volume through frameworks that support you in making meaning from them and support you in reading them critically. These case studies should also inspire you to improve your own work. Beyond this volume, I encourage you to engage the original sources for the ideas from Gutiérrez (2017, 2019) and Boaler (2006). I selected the three core references for this chapter carefully to be accessible to a broad audience of mathematics educators so that I could encourage you to learn more about these ideas from the original authors.

*What questions should I be asking myself?* This chapter is filled with questions that I hope you will ask yourself frequently, but I will highlight one three-part question that I think is at the core of all of this work. I have learned this question from the work of Rochelle Gutiérrez and Belin Tsinnajinnie: *How is equity being defined, by whom, and in service of what?* Are we centering the experiences and voices of our most marginalized students in support of their goals, or are we trying to bring students into our vision of mathematics in service of our own agendas?

### 3.1. *Summary of Frameworks and Practices*

I have organized the three frameworks in this chapter around questions that instructors can ask themselves as they read the case studies or as they work on making their courses more equitable and inclusive. Here are these questions collected together.

Nepantla and Political Conocimiento (Gutiérrez, 2017):

- What does equity mean here, and who gets to define it?
- What about the educational system is taken as immutable? Is the critical axis engaged?
- How are student perspectives about their identity and power incorporated?
- How are the authors being creative about dealing with constraints and resistance? How do the authors stand in solidarity with students?
- How are the contexts and histories engaged, and how are the needs of the most marginalized centered?
- Other key concepts: access, achievement, identity, power, nepantla, critical and dominant axes

Rehumanizing Mathematics (Gutiérrez, 2019):

- *Positioning/Participation:* How does membership and participation in this classroom look? What status differences and hierarchies exist, and how do they relate to those from the broader culture? How do teachers position students as legitimate and competent members of the class and disciplinary communities?
- *Cultures/Histories:* How are students' pasts valued and leveraged in this classroom? What resources and funds of knowledge do students bring with them to the classroom? How are various cultural approaches to mathematics incorporated?
- *Windows/Mirrors:* How do students see themselves and others represented in the curriculum? How is mathematics used as a path toward understanding their own experiences and developing empathy for the experiences of their peers, especially the experiences of those historically minoritized in mathematics classrooms?
- *Living Practice:* In what ways is mathematics experienced as an ongoing, human endeavor?
- *Broadening Mathematics:* How does this classroom incorporate aspects of mathematics outside of abstract symbolism, emphasis on the general cases, and the narrow geometry/algebra/calculus pathway? How are qualitative forms of mathematics valued?
- *Creation:* How are students in this classroom engaged in inventing ideas that are new (at least to them)? What new forms of mathematics are likely to appear in this classroom?
- *Body/Emotions:* How does this classroom acknowledge and leverage the fact that mathematics is done by people, with bodies and emotions?
- *Ownership:* In what ways does this classroom make space for students to do mathematics for their own purposes? How are students' joy and desire to play supported?

Complex Instruction+ (Boaler, 2006):

- *Multidimensionality:* How does this classroom incorporate multiple ways to be successful?
- *Roles:* How is group work structured? How does the work get distributed so that each learner has something important to do?
- *Assigning Competence:* How is the competence of students made visible in the classroom?
- *Student Responsibility:* In what ways are the students responsible for each other and collective success?
- *Justification/Reasoning:* How do students participate in the work of reasoning and justifying mathematical assertions?
- *High Expectations:* How does this course communicate the belief that all students can and should achieve success with challenging tasks?
- *Effort over Ability:* How is the work of learning valued and viewed as changeable rather than innate?
- *Learning Practices:* How are productive practices for learning made visible and named in this classroom?

This final table lists some of the most prominent examples of practices for individual instructors, including examples of these practices from the case studies in this volume.

| Practices | Some Case Studies containing examples or discussion of this practice |
|---|---|
| Promote reasoning rather than answer finding; position students as knowers and creators of knowledge; emphasize mathematics as a human endeavor | Johnson et al. (CS 14 ); Liou-Mark et al. (CS 16 ); Mingus et al. (CS 20 ); Starbird et al. (CS 24 ); Stone-Johnstone et al. (CS 25 ) |
| Acknowledge the ways that our education system is failing to support students; make space for and validate student critiques of the system and the associated emotions | Oliver et al. (CS 21 ); Stacy (CS 23 ); Chang & Chen (CS 7 ); McNeill et al. (CS 9 ) |
| Approach students as having resources for success rather than deficits that need to be fixed; design course activities to support and value multiple pathways to success; make sure that student competence is illuminated in class discussions; express both high expectations and the belief that each student can meet and exceed them; frame learning around effort rather than predetermined by some innate ability | CHAPTER 5 ; Zobitz et al. (CS 30 ); Stacy (CS 23 ) |
| Have students develop a set of group commitments (or norms) about their rights and responsibilities in the classroom; return to these commitments consistently, especially during challenging episodes | Byrne (CS 5 ) |
| Support students in growing their affirming identities as learners and members of a STEM community | DiGregorio et al. (CS 10 ); Wagnon et al. (CS 29 ); Starbird et al. (CS 24 ) |
| Seek to understand student motivation for enrolling; connect course activities with students' goals; work to ensure that students see themselves and their peers represented in the curriculum | Zobitz et al. (CS 30 ) |
| Scaffold effective group collaborations; teach collaboration as a skill | Johnson et al. (CS 14 ); Davis (CS8 ) |
| Design for collective responsibility rather than individualism or competition | Deshler et al. (CS 9 ); MacAurthur (CS 17 ); Akin et al. (CS 1 ) |
| Incorporate explicit discussions of productive mathematical moves and learning practices; incorporate metacognition and help students self-regulate learning | Byrne (CS 5 ); Mawhinney et al. (CS 18 ); Wagnon et al. (CS 29 ); Fuller et al. (CS 11 ); Starbird et al. (CS 24 ) |
| Start from the assumption that faculty are in positions of power, so the impact of our actions is often hidden from our view; do not generalize from the experience of successful faculty members to the experiences of all students; seek data and student perspectives about the classroom experience and adapt based on what you learn | McNeill et al. (CS 19 ); Golden et al. (CS 12 ); Stone-Johnstone et al. (CS 25 ) |
| Accept that no action or inaction from a position of power is neutral and seek to use power to support students not being supported by our systems; learn about histories of oppression and ongoing hierarchies that impact the classroom to better understand patterns in our courses | McNeill et al. (CS 19 ); Stone-Johnstone et al. (CS 25 ) |

## 4. ***References***

### 4.1. *Core References*

### 4.2. *Supporting References*

# 4 The Role of Role Models in Addressing Diversity, Equity, and Inclusion

**Jessica R. Gladstone**
*Virginia Commonwealth University*

Jessica R. Gladstone is a postdoctoral research fellow at Virginia Commonwealth University. Her research investigates the individual and social contexts that facilitate achievement motivation and engagement in STEM.

**Hortensia Soto**
*Colorado State University*

Hortensia Soto is a professor of mathematics at Colorado State University. Her research centers on the teaching and learning of undergraduate mathematics, where she adopts an embodied cognition lens. She frequently delivers professional development that addresses DEI issues for K-16 teachers.

## 1. *Introduction*

In this chapter we begin by briefly summarizing literature on role models in STEM (science, technology, engineering, mathematics) disciplines. We then illustrate how some authors of these case studies made use of role models in their efforts to address diversity, equity and inclusion (DEI). Although one might generally think of role models as instructors or graduate teaching assistants, these are not the types of role models addressed in the case studies. This could be due to the fact that STEM departments continue to struggle to attract and retain graduate students and faculty from underrepresented groups. In these case studies authors detail how they integrated undergraduate students to serve as learning assistants (LA) and thus, as role models for other undergraduates. As part of this summary, we synthesize how LAs were integrated into courses and served as role models as well as the benefits and challenges with integrating LAs into the mathematics classroom. Thus, it is possible that the departments are attempting to leverage the resources that they do have, i.e., students, in order to provide near-peer role models for other students. In the last section, we offer suggestions for training LAs and instructors to serve as role models. We also suggest cost effective ways to implement role models that extend beyond implementing undergraduates into the mathematics classroom because LA implementation can be time consuming and costly. Furthermore, we advocate for institutions to create mechanisms that support and encourage faculty to attend to DEI issues, by finding creative ways to offer role models for our students.

## 2. *Role Models Literature*

One promising way to make college mathematics classrooms more inclusive and in turn address issues of diversity and equity is to introduce role models. In fact, popular news outlets (e.g., Dean, 2014; Levere, 2018) have promoted the notion of using role models to help decrease inequalities in STEM classrooms. **A role model can be broadly defined as someone who is a successful exemplar and can impact a person's motivation, achievement, and sense of belonging** (Gladstone & Cimpian, 2021). In the mathematics classroom there are a number of people who can be considered role models for students. These include, but are not limited to, LAs, peer mentors, tutors, and instructors.

Research has shown that the use of role models as an intervention is effective in several ways (Cheryan et al., 2011; Johnson et al., 2019; Marx & Roman, 2002; Shapiro et al., 2013; Stout et al., 2011). Role models who are perceived to be competent, perceived to be similar to the student, or whose success is perceived to be attainable by the student can have a positive impact on students. Furthermore, research shows that when role models are effectively implemented in STEM classrooms there is the potential for increased student motivation (Dennehy & Dasgupta, 2017; Shin et al., 2016), achievement (Herrman et al., 2016; Krämer et al., 2016; Marx & Roman, 2002), STEM identity (Gilbert, 2015; Shin et al., 2016), and feelings of belonging in STEM classrooms (Dennehy & Dasgupta, 2017; Johnson et al., 2019). This is especially true when the role model and student share similar social identities, that are generally viewed as underrepresented identities in STEM fields (Johnson et al., 2019). Thus, the use of role models can be beneficial in terms of promoting diversity, equity, and inclusion in the mathematics classroom. Bandura's (2001) work shows that a student who observes a role model succeed on a task may then come to believe that they can also succeed on the task; thus, students come to value competence. Further, students who perceive themselves to be similar either demographically or psychologically (e.g., shared values) to the role model may be more likely to be influenced by the role model. For example, students can better assess their own chances of success if they perceive that the potential role model has had similar experiences as the students (Schunk & DiBenedetto, 2020; Schunk & Usher, 2019). Finally, if students perceive the role model's success as unattainable then it can lead to thoughts of "she can do that, but I can't." Therefore, it is important that role models describe their path to success in terms of controllable factors such as effort rather than innate talent or luck, and other strategies that promote a growth mindset (Dweck, 2006).

However, because the study of mathematics is a highly stereotyped domain that tends to exclude women and students of color (Murphy et al., 2007; Steele, 1997) it is critical that role models not heighten feelings of stereotype threat or of not belonging. It is also vital that neither the role model nor the students experience a decrease in motivation. Thus, it is important that careful consideration be made when implementing role model interventions and that role models be adequately trained before working with students. We believe this is true whether the role model is or is not from a marginalized group; however, it is important to consider the social identity of the role model and the students served. Therefore, training and implementing role models must extend beyond showing students how to appropriately solve mathematics problems. We argue that training of role models should include messaging about how to help students recognize that they can be successful in mathematics and that they belong and are welcome in the mathematics classroom. Furthermore, it is important that role models know how to best share their own mathematical journey, including their struggles.

Several of the authors of these case studies indicate that they made use of LAs in their classrooms and that these LAs served as role models because they are successful exemplars. In the following section we summarize how the LAs were implemented into the classroom, how the LAs were trained, benefits of LAs to students and to the LAs, and challenges that can accompany the use of LAs. We note that although several of the case studies illustrate content and pedagogical training of the LAs, we offer suggestions of how to train LAs to serve as role models.

## 3. *Role Models as Illustrated in Case Studies*

The case studies showcased in this volume primarily illustrate how LAs can serve as role models for other students in the classroom ranging from college algebra to abstract algebra and even for K-12 students (Fuller et al. (CS 11 ), Mingus et al. (CS 20 ), Villalobos et al. (CS 28 )). Learning assistants are generally described as paid undergraduate students who completed a particular course and performed well enough in the course to help facilitate discussions or activities with the guidance of a faculty member. Such facilitations generally occur in the classroom but they can also occur outside of class. Authors of these case studies use different terminology to refer to a learning assistant, including peer mentor (Wagnon & Hubbard (CS 29 ), Zobitz et. al (CS 30 )), learning community mentor (Byrne et al. (CS 5 )), peer tutor (Davis (CS 8 )), peer leader (Liou-Mark et al. (CS 16 )), and undergraduate teaching assistants (Deshler et al. (CS 9 ), Villa et al. (CS 27 )). In an effort to honor the authors' work, we will use their term when we reference their case study.

### 3.1. *Summary of case studies*

Appropriate and intentional training is necessary in order for LAs to be effective role models. We found that most programs described in the case studies integrate LA training that focuses on content and pedagogy. For example, in an effort to help Hispanic students, African American students, first-generation college students, and students from a low socio-economic status family, Liou-Mark et al. (CS 16 ) describe how their institution created a weekly extra hour problem-solving session for some mathematics courses to be facilitated by a peer leader. Such sessions were available for calculus prep courses (e.g., College Algebra, Trigonometry, and Precalculus) and for Calculus I and II. In order to serve as a peer leader, students completed weekly modules related to each problem-solving session before the semester began. After they completed the modules, the students attended a 2-hour content orientation where they discussed solutions to the tasks, prior knowledge required for the tasks, anticipated student errors and questions, and alternative solutions to the tasks. First time peer leaders also registered for a one-credit course that focused on pedagogical training. This included learning about facilitation strategies, learning theories, stereotype threat, diversity, and inclusion. As part of this pedagogical training the peer leaders also learned how to integrate ice-breakers, facilitate group work, encourage students to share their work, and create norms for their sessions.

While such training is valuable, it could be advantageous if the instructors also supervised the LAs as they worked with the students to document how the LAs explained mathematical concepts using their new-found teaching strategies and how they attended to stereotype threat, diversity, and inclusion while working with the students. Given that some faculty members lack training in attending to these concepts when teaching, we suspect it can be especially challenging for first time LAs. Institutions who do not have the ability to increase their credit requirements, without negatively impacting students' financial aid, might struggle to implement these strategies. This might require negotiating other university requirements, which could encourage such institutions to adopt some of these practices.

It has typically been difficult for students who are African American, Hispanic, first-generation college, or Pell-grant recipients to gain academic capital; to support those students and others, Mingus et.al. (CS 20 ) also incorporated LAs into most of their Calculus I classes. The authors selected the LAs based on their gender and ethnicity that best represented the students whom they were trying to serve. Besides receiving pedagogical training, the LAs met with the course instructor each week. The authors describe their LAs as bridges between the students and the instructor, as mentors, as cheerleaders, and as someone who has to be able to co-learn with the students. As such, at the weekly meetings with the instructor the LAs shared the students' experiences with the content during class or outside of class. This gave the instructor a better sense of how well students were responding to the content and how to modify future lessons accordingly. The LAs also provided the students with good study habits, time-management strategies, insights into homework expectations, advice on how to develop their relationship with the instructor, information regarding course and instructor selection, and course registration. The LAs also facilitated review sessions before an exam at times and at locations that were convenient for the students.

Strengths of this case study are the weekly debriefing meetings between the instructor and the LA, the presence of both the instructor and the LA during class, and opportunities for students to gain academic capital from the LAs. This model certainly gives instructors an opportunity to observe and assess the LA- student interactions which could also be discussed during the weekly debriefings. Given that this model requires quite a bit of extra effort from the instructors, it is possible that instructors who teach at research one institutions may find the extra work burdensome. As such, we recommend that institutions implement strategies to encourage buy-in from faculty to address DEI issues, such as counting this work towards tenure and promotion.

While the above case studies' central focus was the implementation of LAs, other studies briefly introduced the notion of LAs. For example, Fuller et al. (CS 11 ) integrated LAs into their active learning-based calculus course aimed to support their large Hispanic and African-American population. As part of their duties, the LAs met with their instructor each week and were also instrumental in relaying ideas between the students and the instructor. Similarly, Deshler et al. (CS 9 ) describe how they implemented undergraduate teaching assistants into their calculus sequence to primarily facilitate group work during class. Villalobos et al. (CS 28 ), describe a unique case where a majority of their LAs were bilingual in Spanish and English. As such they facilitated the recitations for the Calculus I course which was offered in Spanish. This helped put the students at ease because they could discuss concepts in their native tongue with the LA.

Although most of the case studies discuss integration of LAs into their calculus preparatory or calculus courses, Davis (CS 8 ) is the only case study author to discuss integration of peer tutors in upper division inquiry-based classrooms including linear algebra, proof writing courses, and abstract algebra. Davis' case study (CS 8 ) is a description of her attempt to support Native Hawaiians and Pacific Islander students. With National Science Foundation funds, she hired students from traditionally underrepresented groups (e.g., women, Hispanics, Native

Hawaiians) to serve as peer tutors. Prior to the start of the semester the peer tutors met with their faculty mentor to discuss expectations and duties. This was also an opportunity for the peer mentors to request a copy of the text, advanced copies of notes, access to Blackboard, and solutions to homework assignments. Because a majority of the peer tutors are mathematics education majors, they are familiar with facilitating problem-solving, asking guiding questions, and ensuring that they are not doing the work for the students. As such, these peer tutors do not have additional training.

Fuller et al. (CS 11 ) acknowledged that sometimes students struggle in inquiry-based calculus courses and that both the instructor and the LAs are not equipped to address these challenges. As such they suggest professional development training for both faculty and LAs, so that they learn to foster best practices as students work collaboratively. Villalobos et al. (CS 28 ) provide a rich description of the professional development that they offer their undergraduate teaching assistants, who tend to be the same age, ethnicity, and bilingual in Spanish and English as the Precalculus students who they teach. The professional development consists of a 3-hour workshop centered on cooperative learning and inclusion along with their theoretical underpinnings. The facilitators of the workshop teach about cooperative learning by modeling individual accountability, group reflections, and social skills. Learning about inclusion begins by having the teaching assistants reflect on their own high school experiences and by recognizing that their best experiences occurred when they felt a sense of belonging. These conversations prepare the teaching assistants to learn about communities of practice (Lave, 1991; Wenger, 1998) and how to create such a community with the students who they will teach.

### 3.2. *Benefits and challenges*

Those who integrate LAs into their classrooms claim that LAs serve as role models especially when they look like the students who they are serving, help students navigate the hidden curriculum, assist students with research opportunities, alleviate students' fear and embarrassment in asking questions, foster a sense of belonging, and are easier to understand (Chang & Chen (CS 7 ), Davis (CS 8 ), Liou-Mark et al. (CS 16 ), Mingus et al. (CS 20 ), Zobitz et al. (CS 30 )). These are all mechanisms that can help support students from underrepresented groups to be successful in their mathematics courses. As an example, Chang and Chen (CS 7 ) describe their efforts to address DEI at a private Hispanic-serving institution by redesigning their College Algebra course to include embedded peer tutors who were mostly women majoring in STEM. The authors found that College Algebra courses with a peer tutor had a higher A, B, C (i.e., grades typically higher than a 70%) rate of 62.66% compared to the courses without a peer tutor, whose A, B, C rate was 38.30%. As such, it seems that LAs can also help students be successful in a course and thus, see themselves as future LAs. Similarly, Davis' (CS 8 ) survey data from Calculus II and Linear Algebra students indicate that students felt that their peer tutor encouraged them to continue in STEM, inspired them to graduate, and served as a role model.

Besides benefits to the students, there are also benefits to the LAs. Liou-Mark et al. (CS 16 ) indicate that students who serve as a peer leader can experience an increased level of confidence as a leader, reinforcement of course content, and a feeling of influence because they serve as a role model. The peer leaders also appreciate that once they are trained in the pedagogical component they can serve as a peer leader in other courses as long as they have been successful in that course. Davis (CS 8 ) also documents that their peer tutors expressed an improvement in their mathematics skills, a sense of enjoyment in helping others, and an increase in their confidence level. More importantly, these peer tutors were more determined to complete their mathematics major or add it to their existing major. Davis (CS 8 ) adds that serving as a peer tutor empowers women, Hispanic students, and Native Hawaiian students because they gain status as experts and as role models, which may in turn empower students who are from similar marginalized groups. Mingus et al. (CS 20 ) emphasize that their LAs also experienced a sense of empowerment because they were instrumental in helping the instructor re-negotiate classroom norms with the students, especially those related to assignments that required higher cognitive demand.

As mentioned above a majority of the case studies related to role models focused on LAs, but Villalobos et al. (CS 28 ) also describe efforts to support K-12 Latino students. As part of their outreach efforts, their university brings in K-12 students to engage in hands-on STEM activities. Student staff, who are generally Latinos, facilitate the activities in English, but they can deliver the activities and answer questions in Spanish as need be. As such, undergraduates can serve as role models for their peers as an LA or for K-12 students as part of outreach activities. Villalobos et al. (CS 28 ) also offers panels of scientists for their undergraduates to learn about career and research opportunities in STEM. The panelists are generally women or Latinos; this allows a space for undergraduates to have their own role models. As such these case studies illustrate how mathematics departments can move "toward change in support of multiple visions of justice-making" as suggested by Larnell (this volume) in the Manifesto.

Although several authors document benefits to the LAs, Byrne et al. (CS 5 ) delineate some of the downfalls with integrating LAs. First, there is the financial stress that this may place on departments. Many of these cases studies relied on grant funding or work study in order to launch their efforts, thus, there is a concern regarding sustainability of integrating LAs into courses. Besides the extra cost, some of Byrne et al.'s (CS 5 ) learning community mentors experienced added stress trying to meet the needs of the course and their own course obligations. Some of the stress was due to the extra work completing course materials and to facilitating the extended two-hour calculus sessions for 30 students, without appropriate training. We speculate that integrating LAs into mathematics courses could also result in LA-student conflicts where some students do not view the LA as an authority figure, as knowledgeable enough, or as a peer. This is problematic, especially if the LA is a student from an underrepresented group. Such conflict can be demoralizing for the LA who is supposed to serve as a role model. LA-student conflict might also arise if the LA is from a majority group and the students in the course are from underrepresented groups. Conflict could arise because the students are unable to identify with the LA and vice-versa. Another potential challenge not directly addressed in the case studies is the fact that undergraduate students may not know what it means to serve as a role model, which we elaborate on below.

## 4. *Recommendations*

As can be seen in the case studies reviewed above, there are a number of ways in which undergraduates can serve as LAs and potentially as role models for other students in college mathematics classrooms. While we applaud those making efforts to provide role models for students via LAs, we are concerned that much of the burden or privilege falls on the students and not the faculty. We strongly believe that faculty also need to be involved in serving as role models. Thus, while pedagogical and content training is important in order for LAs to help students with course material, **LAs and faculty also need training about how to serve as a role model**. In this section we argue that knowing how to create a sense of community could help to intentionally cultivate future role models. We also advocate for mathematics departments to have role models that extend beyond undergraduate students and propose ways in which this might be achieved. After all, as stated by Larnell in this volume (BUILDING AN EVOLVING FRAMEWORK: A CLARION CALL / MANIFESTO), "The responsibility for redressing these deep inequities should not be relegated to a small impassioned group of reformers but shared by the entire community of faculty members and academic researchers who steward the undergraduate mathematics education trajectory. Putting it plainly, the responsibility for it is on all of us."

One might wonder how someone from an underrepresented group suddenly becomes entrusted to be a role model for others if they have never had a role model. We argue that this begins with those who have privilege and access to mathematics. In the mathematics classroom this can be the instructor because a role model is a successful exemplar who can impact a person's motivation and achievement (Gladstone & Cimpian, 2021). For example, in their case study Wagnon and Hubbard (CS 29 ), describe how one of their participants, David, attributed his perception that he could be successful in a STEM career to his junior high school teacher. In junior high, David had no aspirations of attending college; in fact, his plan was to quit high school, get a job, and help support his family. But, David's teacher saw potential in him and invited him to enroll in a college readiness program. The fact that someone - a person of authority - took enough interest in him helped David transform his self-image and motivated him to be fully engaged in his studies and with his peers. Such engagement led to a sense of belonging and David's persistence to graduate with a STEM degree. As such, David's high school teacher might be his first STEM role model, and the teacher might not have been from an underrepresented group. This might lead us to also ask, how do instructors learn to be role models?

Given the systemic exclusion of certain groups, it can sometimes make it difficult to have role models from marginalized groups. Therefore, role models who look like our students aren't always available, but we argue that role models should humanize all students and creating a sense of belonging can be a first step to achieve this goal. Having a sense of belonging is one of the four characteristics of a sense of community. The other three characteristics are feeling nurtured, feeling a sense of influence on others, and experiencing an emotional connection (McMillan & Chavis, 1986). As such efforts to create a sense of community and to rehumanize mathematics (Gutiérrez, 2018), requires "shifting the role of authority from teacher/text to other students, … creating mathematics, … conjuring feelings of joy, … [and through] ownership [of the content]" (p. 5) and this work must be intentional. Thus, we recommend that anyone who serves as a role model acknowledge these efforts and convey to students the rationale behind these efforts. In particular, LA training should acknowledge that beyond knowing the content and knowing how to explain the content**, LAs need to be able to help the students feel** **welcomed, nurtured, valued, and emotionally connected to the rest of the class.** These conversations are important so that LAs do not become role models through osmosis, but rather through purposeful actions and conversations. Such tasks are not always easy for faculty and thus, we should expect that some LAs might struggle executing these practices, especially because these practices invite one to be vulnerable. This is an invitation that some may quickly decline, especially populations from underrepresented groups because vulnerability can lead to embarrassment or may be a threat to their personhood. Furthermore, such invitations can appear as simply exclusionary for those with marginalized identities. This is why it is important for those with power to serve as role models also. On the other hand, we believe that exhibiting vulnerability is part of serving as a positive role model because it is part of humanizing our experiences with mathematics. As such, acknowledging these practices might also help students recognize what it means to be a role model and begin to critically reflect on how to support justice-making. Furthermore, instructors and LAs who make such efforts may ease the disruption of current mathematical norms which privilege certain students and exclude others (Adiredja & Andrews-Larson, 2017; Gutiérrez, 2017).

As Liou-Mark et al. (CS 16 ) convey, such training requires institutional support for appropriate professional development for the faculty and for students who may serve as LAs. There are also costs associated with hiring students to be role models in some form (i.e., LAs, peer tutors, teaching assistants, etc.) and extra required resources such as time to train role models, space for extra study sessions, and extra texts. Although we strongly support the notion that students should be paid for their time as role models, we understand that this is not always possible and that resources for training are not always available. Therefore, in what follows we provide low- or no-cost alternatives for integrating role models into mathematics classrooms.

One popular and successful way to present low- or no-cost role models to students is via text or videos. In fact, a recent meta-analysis by Lawner and colleagues (2019) found that videos or printed materials of role models are just as effective as in-person interactions. For example, Shin and colleagues (2016) had undergraduate students read biographies of successful women and successful non-European American men in STEM who obtained success in STEM through hard work as opposed to a perception of natural talent. The researchers found that reading these role models' biographies increased students' interest in STEM, identification with STEM, and self-efficacy regarding STEM. In particular they found that women who read the biographies of other women had an increase in their perceived identity compatibility with STEM. Thus, the effects of having students read short biographies of role models in STEM is very promising, cost effective, and does not require training of students. Using texts or videos to present role models is cost effective because they can be widely distributed and they can be presented to large numbers of students at different times and locations. They also limit any additional burdens that role models from underrepresented groups may experience as they become more in demand. This is especially true with role models who belong to traditionally underrepresented groups in STEM.

Although presenting students with role models via text or videos has shown to successfully impact students' experiences with STEM (Shin et al., 2016), one needs to carefully consider how the role models are portrayed and described in the text or video. Research shows that such texts and videos must strike the right amount of balance regarding the role models' competence. It seems that if the role models are described as extremely competent, then there is the potential of demotivating students (Woodcock, 2012; Ziegler & Stoeger, 2008) because they do not believe they can attain similar success and levels of competence. As with the use of LAs, role models portrayed in texts or videos should also be portrayed as meaningfully similar to the students (e.g., the role model worked hard for their success rather than being naturally talented). Additionally, it is important that the role model's success is described in a way that students will perceive as attainable. It is also beneficial if the text offers suggestions for how students can be as successful as the described role models, while also acknowledging that feelings of struggle are normal. Such acknowledgements can help humanize the STEM learning process. As such, texts and videos should be versatile enough so that students can identify with someone portrayed in these media. The book, *Living Proof: Stories of Resilience along the Mathematical Journey*, published by the MAA and AMS is a good resource with a variety of personal stories. The figure below illustrates four recommendations for implementing effective role models in STEM (see Gladstone & Cimpian, 2021). In the next section we provide a description of Soto's (one of the authors of this chapter) efforts to provide role models for her students.

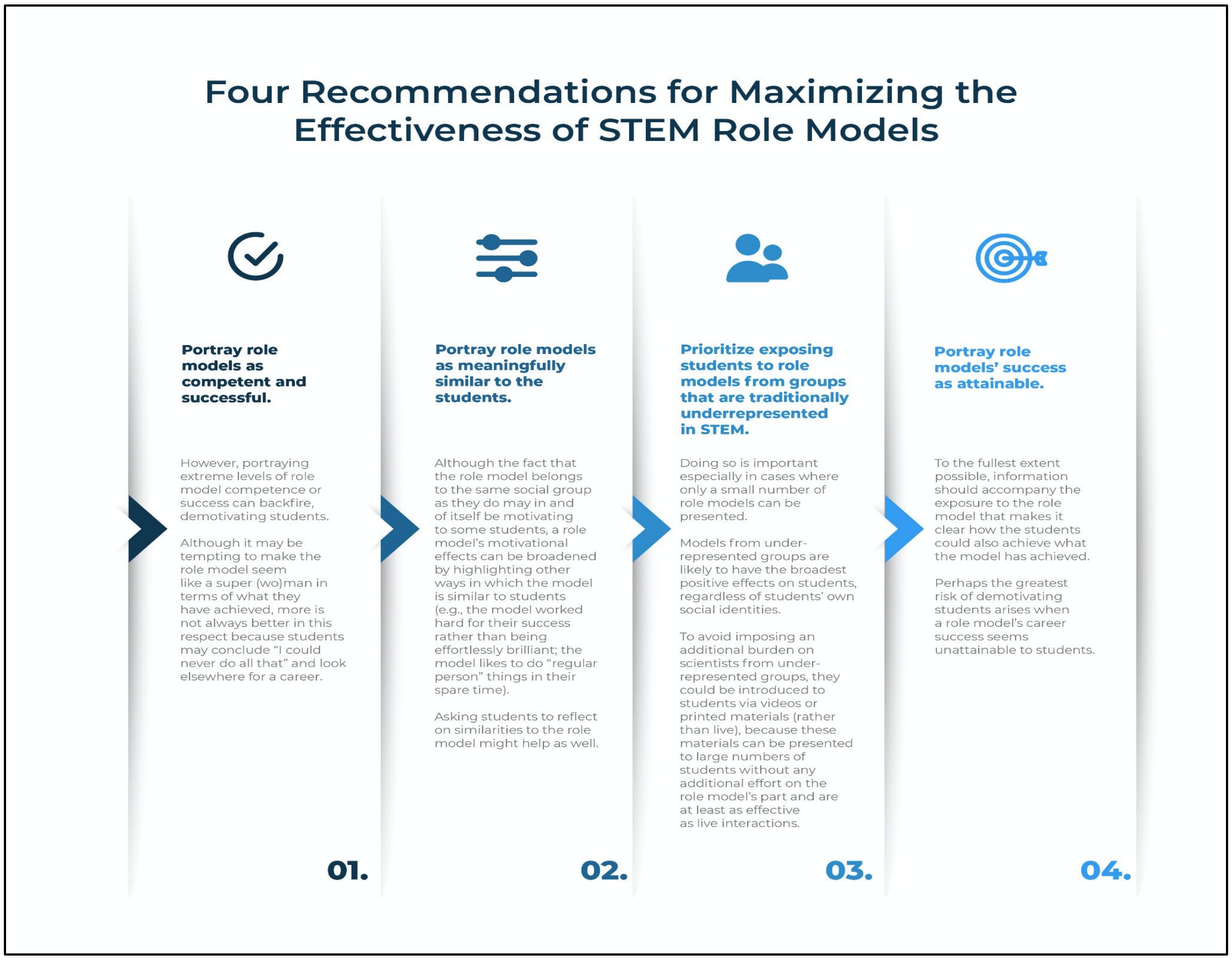

*Figure 1. Recommendations as cited in Gladstone and Cimpian (p. 16, 2021). High-resolution versions of this infographic in .eps, .jpg, .pdf, and .png formats are freely available at https://doi.org/10.6084/m9.fgshare.c.5681674.v1*

## 5. *Personal Efforts to Provide Role Models*

It is well established that part of transforming our perceptions about self or others requires reflection. Thus, reflecting must be part of introducing students to text or videos regarding role models. Such reflections can require students to share how their mathematical experiences are the same or different as the role models, how they can identify with the role models, how they hope to be similar to the role models in the future, etc. For example, Van Camp et al. (2019) indicate that students who reflect on their similarities to the role model may also have increased feelings of belonging. This can contribute to feeling that one is part of a community - in this case as part of the mathematics community. Soto, one of the authors of this chapter, recently integrated the *Living Proof* book into her History of Mathematics course, which primarily consisted of preservice mathematics teachers. Soto asked her students to read and select 2-3 stories and then write a paper where the students shared why they selected the stories, how they were surprised by the stories, how they related to the mathematicians depicted in the stories, or how they experienced a visceral feeling as they read the stories. The students were also asked to explain how such knowledge and reflection was useful to them as a future teacher and mathematician. Finally, Soto asked her students to write their own story about their mathematical journey. This assignment was a first step to help students become aware of some of the struggles that some mathematicians experience and to introduce them to mathematicians from underrepresented groups. In her reflection one student commented, "When mathematics became hard for me, I felt that I didn't belong in the mathematics community. Laura highlights her struggles and how they have helped her to form a 'growth mindset' so she can enjoy working through the struggle. This is something that resonates strongly with me. Just as with Laura, struggling through mathematics has made me better at it than I ever would've been had I never struggled and it has helped me to develop a growth mindset that I will teach to my students."

Soto also leveraged online teaching and invited mathematicians from underrepresented groups (women, Hispanic, and African American) to give a guest lecture each month. The purpose of this activity was for her students to hear about the mathematical contributions of these mathematicians and how they connected to the history of mathematics. Too many times we ask mathematicians from underrepresented groups to only share their struggles, which can be inspiring but also burdensome. Having them share their mathematics expertise is a way to showcase members of underrepresented groups as the "doers of mathematics" as suggested by Larnell in this volume (BUILDING AN EVOLVING FRAMEWORK: A CLARION CALL / MANIFESTO). Besides inviting monthly mathematicians, Soto also had students conduct a virtual interview with a mathematician from an underrepresented group regarding their professional and personal endeavors including successes and challenges. As part of the assignment the students again had to reflect on their experience with this interaction. Frequently students never interact with mathematicians from underrepresented groups, especially at institutions that primarily serve white populations. Assignments such as these can help students become aware of DEI issues, learn how to mentor students from underrepresented groups even if the students aren't themselves from an underrepresented group, and possibly learn how to serve as a role model for others. Furthermore, such teaching practices may prevent stereotype threat and support justice-making in the classroom. For example, in his reflection related to this project a student remarked, “I think my interview ... was extremely enlightening on the importance of diversity and inclusivity in mathematics. It is something I have been aware of for quite a while, but never something I have internalized and truly thought about past surface level. The fact that there are only 12 Native Americans with Ph.D.’s in mathematics is an obvious indicator that there is not enough diversity in mathematics, and there are systemic issues prevalent that mitigate the ability for people to be able to become mathematicians and even succeed in mathematics at any level.”

## 6. *Conclusion*

Although the case studies portrayed in this volume primarily portray how undergraduates can serve as role models for other undergraduates, we believe all of us including administrators, faculty, graduate students, undergraduates, and K-12 students need good role models who motivate us and help us believe that we can do and be better human beings. Furthermore, we believe that **we are all capable of serving as role models for others**. From our experiences, we have found that role models can appear in our life in many different ways even if the role models do not look like us or have not had similar experiences as us. We have found that effective role models impacted us professionally and personally. They have motivated us to succeed in our courses, encouraged us to pursue graduate degrees, and inspired us to be role models for others. In general, **effective role models model good character**. This character can manifest as the ability to create a space where others feel welcome, feel that their ideas and their humanity are valued, feel nurtured, or where there is an opportunity to develop an emotional connection with the role model. Different people attend to these characteristics in a variety of ways based on their power, their privilege, and most importantly their **willingness to be a change agent**. For example, administrators can serve as change agents because they have the power to hire faculty who are from the same underrepresented groups as the students. Faculty and graduate teaching assistants have the ability to create an emotional bond with students by interacting with students during class or during office hours, by sharing personal aspects of their life (such as their mathematical journey), or by simply pronouncing the students’ names correctly. Faculty can also demonstrate that they value students’ ideas by encouraging them to apply to graduate school or inviting them to partake in the faculty’s research projects. Such gestures are acts of kindness that affirm social inclusion and allow students to identify as a member of the classroom community (Estrada et al., 2018), which in turn can support diversity and equity.

We invite the reader to reflect on the following questions.

- What actions might you take to encourage others to serve as role models and to serve as a role model yourself?
- How might you integrate role models into your courses and assess their effectiveness?
- How might you take into consideration the individual students in your courses and the role models that would be most effective for the students?
- How might you engage your colleagues in your department to develop and implement training to humanize the subject of mathematics?

# 5 Course Redesign: Pathways Towards Transformation

**Geillan Aly**
*Compassionate Math*

Geillan Aly is the founder and CEO of Compassionate Math, a company founded to transform mathematics teaching and learning for teachers and students. A former assistant professor of mathematics, Aly received a Master's degree in mathematics and a Ph.D. in teaching and teacher education from the University of Arizona. Her research focuses on affective learning factors among developmental mathematics students. Underlying her work is a dedication to equity and social justice in mathematics education and a firm belief that learning mathematics is both an emotional and cognitive endeavor. Geillan enjoys traveling, seeing live music, and is an avid chef.

*"And the day came when the risk it took to remain tight inside the bud was more painful than the risk it took to blossom."*

At some point, faculty in a mathematics department come to the realization that they are not serving their students well. A critical moment comes when major course redesign or redevelopment is necessary. Rather than making small modifications, transforming a course comes from a radical shift in perspective. These changes may be due to advances in pedagogical research, a change in stakeholders' needs, changing demographics, or external pressures. The case studies in this volume share department's approaches to these changes, and this thematic chapter provides an organizing framework for understanding the myriad motivations, responses, and results.[2]

For example, after the University of Texas at Austin was legally required to accept any Texas high school student who graduated in the top 10 percent (now 6 percent) of their class, the university had a more diverse pool of students, both demographically and in terms of their preparation. As a result, they faced a sudden need to revamp mathematics courses to support a body of students with a more diverse background (Jensen-Vallin et al. (CS 13 ); Starbird et al. (CS 24 )). Similarly, when developmental coursework was reduced or eliminated by mandate, the California State University system and universities in Texas had to develop co-requisite courses or institute changes to support students (Benken et al. (CS 2 ); Canner et al. (CS 6 ); Jensen-Vallin et al. (CS 13 )). In these cases, these decrees provided an opportunity to evaluate and redesign classes to better serve those with the most needs. In other cases, signs that change is needed may also be more subtle and come from within. A member of the faculty at one department of mathematics was shocked and prompted to make changes when a student remarked that the environment in STEM classes was not welcoming to underserved students (Stacy (CS 23 )).

The field of mathematics education is increasingly acknowledging that lecture-based mathematics courses and curricula uphold a culture of White supremacy and marginalize underserved students0F[3]. Mathematics education evolved from focusing on mathematical content to more research on equity, power, and affect (Gutiérrez, 2013; Lubienski, 2002; Martin et al., 2010). These and other scholars' works have become more mainstream while highlighting and exploring the existence and consequences of an inequitable mathematics classroom. In turn, equity-directed instructional practices have been developed and incorporated into the mathematics classroom (Gutstein, 2016; Gutstein & Peterson, 2005; Rubel, 2017). As a result, mathematics departments across the country at various types of institutions have redesigned their mathematics courses to better serve a given proportion of their students1F[4]. As the realities of our inequitable society become more stark, it's hard to ignore that fundamental changes should be made in the mathematics classroom in order to provide for a more equitable education.

[2] Because every case study reflects a change, this thematic chapter brings together many more case studies than in other chapters, and is thus longer than other chapters.

Undergoing such significant changes is a daunting task and the consequences of making course changes could be severe if redesign outcomes do not meet expectations. Students' education and future are at stake, especially for those majoring in or working in a career focused on mathematics (such as STEM or business programs). Moreover, students who have struggled in previous classes may attribute their low performance in a badly-designed course to their own inabilities and not the structure of the course. This could reinforce students' previously held negative beliefs about their ability or trigger feelings of imposter syndrome which can undermine their sense of belonging in college. Professional careers can be jeopardized for tenure-track faculty who are involved in program changes that do not work as proposed. Furthermore, overhauling a course takes significant time, energy, and emotional dedication. Junior faculty could later be judged that their efforts should have been better spent focusing on scholarship, especially if redesign is not immediately successful. The overall reputation of a department may also be tarnished. Other departments may not be comfortable with intended changes.

Changes to courses and programs intended to create a more inclusive experience may be uncomfortable for stakeholders. The majority of mathematics department chairs, faculty, and those who would redesign or implement changes likely do not share the same backgrounds or experiences with the students they intend to support. While mathematics departments are made of mostly White or Asian male faculty[5], students are not as homogeneous a group. Being members of a hegemonic group and implementing changes to better serve a given underserved population requires one to come to terms with their own privilege, recognize the inherent biases in a system that made them successful, and consider an unfamiliar path from the one which brought them success. Such vulnerability and empathy requires significant effort.

---

[3] To define a culture of White supremacy, I refer to Ansley's (1997, p. 592) writing, "[By] 'white supremacy' I do not mean to allude only to the self-conscious racism of white supremacist hate groups. I refer instead to a political, economic, and cultural system in which whites overwhelmingly control power and material resources, conscious and unconscious ideas of white superiority and entitlement are widespread, and relations of white dominance and non-white subordination are daily reenacted across a broad array of institutions and social settings." To further explore and understand how White supremacy is a part of the fabric of our society and systemic in multiple facets of our day-to-day living see Mills (2003).

[4] Although most of the drives towards equity in these case studies focus on BIPOC students, this is not always the case. Some changes were designed to support BIPOCs, others focused on overlapping but non-identical groups: first-generation, or low income students. Women were also a population of focus to drive changes. Each case study discusses their population of interest. I recognize that most implemented policies were designed to support a specific population but that population is not always specified so the reader can focus on the larger process of mindfully implementing solutions to support your population(s) with the greatest needs.

[5] Data gathered from NCES in 2003 demonstrates that females, Blacks, Hispanics, Native Americans, and Pacific Islanders are underrepresented as faculty when compared to US Census data (NCES, 2019). It is important to note that The American

The price for not changing may be even greater. Traditional lecture-based classes do not favor underserved students (Palmer, Maramba & Dancy, 2011). As a group, instructors have demonstrated biases towards such students throughout their educational lives (Sue, Lin, Torino, Capodilupo, & Rivera, 2009). Textbooks are written so that the language and examples may not be relevant or accessible to diverse populations (Abedi & Lord, 2001). Mathematics curriculum has roots in Western and European academic traditions which excluded women and people of non-European descent (Burton, 1990; Joseph, 1987, 2010). For example, the credit of Khayyam quadrilaterals is given to Saccheri, an Italian who translated Khayyam's work from Persian to Latin 600 years after Khayyam wrote it. It is therefore not surprising that many students do not see themselves in their studies of mathematics. It is clear that changes should be made to the traditional mathematics classroom to embrace and include more students.

When transforming a course, there are endless possibilities. Choosing a textbook, deciding on whether and how to include a computer program for students' use, structuring class meetings, altering assessments, providing students with outside support, and rethinking curriculum are just some of the many decisions that can be considered. As such, after recognizing the need for change, inertia can easily replace momentum after countless department meetings to address these decisions. This chapter is designed to support faculty and administrators considering revamping mathematics courses to more equitably meet students' needs. This analysis outlines ten areas of consideration based upon the collective wisdom of this collection's case studies. In this chapter, I outline a ten-phase framework to transform a course to foster a more equitable environment; each can be considered individually or in tandem with others. My hope is to inspire, not overwhelm and paralyze, the reader.

The corresponding case studies focusing on Calculus courses demonstrate that changes made to each course varied and were based on addressing specific needs or challenges faced by the faculty or a particular group of students while building on the myriad strengths of the institution and its members. Each case leveraged only one or two phases in the outlined framework. Overall, the process by which these departments reconceptualized their Calculus courses was geared towards their department's strengths and students' specific needs. In many cases the redeveloped Calculus courses were a result of creative and bold initiatives to reconceptualize what teaching Calculus means and why it matters.

Equity is not a destination, it is a journey (Gutiérrez, 2012, 2017). This journey requires faculty and administrators to take steps in that direction. Transforming mathematics classes to foster an equitable learning environment requires reflection and an openness among the faculty to

---

Mathematical Society publishes an annual profile of mathematical sciences and departments at four-year colleges and universities in the United States. This AMS profile provides a gender breakdown of faculty but does not parse the information by race or ethnic background. It is critical that such data be regularly gathered and analyzed by AMS to recognize the stark inequities in our field as an initial step towards achieving a more equitable mathematics community.

come to terms with and recognize their potential biases, privileges, and how these phenomena can impact their students' learning environment.

**1. *Intentionally search for equity***

Before any significant work can be undertaken to alter mathematics courses and make them more equitable, the intent behind the course redesign efforts must be considered. What specific group are you looking to benefit from the redesigned course? Blanket "for all" approaches can be problematic in that they often carry the assumption that success for students of color means being less Black, Latinx, or Native American and adopting more White mannerisms and values (Martin, 2009). Second, working to support "all" students and working from a colorblind approach has a tendency to reinforce Whiteness. This premise is well argued and a fundamental assumption outlined in the accompanying vignette of the *Transitioning Learners to Calculus in Community Colleges* research project designed to use an equity-focused lens to identify successful mathematics programs in community colleges (Burn et al. (CS 4 )). Furthermore, there are problematic consequences to considering mathematics to be race and culture neutral (Ball et al., 2005; Fasheh, 1982; Nasir et al., 2008). Finally, critical race theorists warn against the argument that intent based on self-serving motivations should be examined so that benefits to the hegemonic groups become a consequence, rather than a motivator, of acting (Basile & Lopez, 2015; Delgado & Stefancic, 2011).

The motivation to alter a course to support a group of students can affect the direction a new course takes. Attending to equity in STEM has been categorized by the motivational perspective of the stakeholders (Basile & Lopez, 2015). The leaky STEM pipeline3F[6] can be seen as a problem because students who are fundamentally interested in the sciences do not always pursue their passions or interests because they do not see themselves as a full participant of that community, and thus we need to increase diversity in STEM to support minoritized students' passions. From a different perspective, the leaky STEM pipeline may be a problem because demographic shifts project fewer engineers who can compete with China, and thus we need to increase diversity in STEM to compete globally. These two perspectives, both aimed at increasing minoritized students in STEM, can result in vastly different course designs aimed at supporting the goal of increased diversity in STEM. Consider the differences between two hypothetical changes made to statistics courses to increase diversity in the major. In the first, the course could establish a fundamental assumption that quantitative literacy is a civil right and that being more statistically literate can help shape students' future social engagement by giving students the ability to understand and critique their world. Such a perspective would center the rich history and contributions of marginalized communities and foster discussion on what it means to be statistically literate or draw valid conclusions that oppose the hegemonic narrative. Students in such a class can analyze policing records as a basis to reallocate city budgets to potentially better serve the community. The other course would be designed with the fundamental assumption that the United States needs to graduate more STEM majors in order to maintain geopolitical hegemony. This course could provide more attention or support to increase target students' grades.

[6] The leaky STEM pipeline is often used to describe attrition among STEM majors. However, this metaphor does not account for the critical perspective that the STEM pipeline leaks due to societal inequities. The recognition that this "leakage" is not random or that there is an implication that students are passive in their education is also not generally attributed to this metaphor. Filters or membranes may provide a more appropriate perspective. For more, refer to the following (Blickenstaff, 2005; Herzig, 2004)

The reasoning and motivation behind the intent to redesign a mathematics course to be more equitable can carry a commitment to the well-being of a group of students or be based on short-sighted self-interest. This can ultimately shape the course and the degree to which this endeavor is successful and can help rehumanize mathematics (Goffney et al., 2018; Gutiérrez, 2017). Furthermore, motivation based on supporting the curiosity and desires of the individual can positively affect larger national interests.

Analyzing student grades in mathematics classes at your institution is a straightforward approach to determining whether certain groups of students disproportionately appear in the DFW group. The achievement gap is a natural catalyst which can trigger stakeholders to redesign a course. However, "gap gazing" is problematic and has been demonstrated to perpetuate deficit thinking and a reliance on surface-level changes to content rather than reflecting and investigating deeper issues that contribute to inequitable outcomes among and between students (Adiredja & Andrews-Larson, 2017; Gutiérrez, 2008; Gutiérrez & Dixon-Román, 2010; Harper, 2010). Furthermore, the notion of an achievement gap versus an opportunity gap, which recognizes that individual achievement is related to the societal opportunities available to an individual, endorses the current inequitable educational system and places the burden of closing that gap on those who are most affected (Oliver et al. (CS 21 )). If an "achievement gap" exists in your department, this is an indication that changes should be made, however it is critical to further investigate the source of these differences. The existence of any gap is an indication that the department is not meeting the needs of a specific group of students, but there is no indication of what the root causes may be.

When Centre College and Southwestern University recognized that an increasingly diverse student body was not proportionately represented in their Calculus passing rates, further investigation determined that students from minority or first-generation backgrounds were struggling with unfamiliar content in their Calculus class. In contrast, students from high-performing high schools, continuing generation students, and White and Asian students were more likely to be familiar with the content and had often previously taken Calculus courses in high school (Kilty et al. (CS 15 ); Stacy (CS 23 )). The experienced students had an advantage over the inexperienced students leaving the latter feeling marginalized. This difference was mitigated at Centre College and Southwestern University by modifying the order in which course content was presented so that all students were wrestling with unfamiliar ideas. Furthermore, students learned to use the open-source software R and RStudio as their technological tools – software that was unfamiliar to almost all of the students. Both changes provided a more equitable learning experience where all students were challenged with new ideas regardless of their ability to afford a graphing calculator or their previous experiences and preparation (Kilty et al. (CS 15 )).

An honest accounting of the community's beliefs and assumptions can be transformative. When efforts to change at a Midwestern institution emphasized the power dynamics in the classroom, professional development workshops were instituted to de-emphasize finding answers, help instructors recognize their implicit biases, and address the misguided notion that students just needed to work harder to be successful. In this case, there was a recognition that the responsibility for student success and engagement came from the instructors' practices and was not based on students' abilities, attitudes, or any other deficit-based perspective (Johnson et al. (CS 14 )). Faculty at Centre College and Southwestern University presumed that traditional methods of teaching mathematics weed out BIPOC, Pell-eligible, and first-generation students. Thus, the entire Calculus sequence was revamped to provide a more equitable learning environment (Kilty et al. (CS 15 )).

Evaluating enrollment in Calculus classes may point to changes that need to occur in earlier classes. Mercy College is an HSI with 70% female student enrollment (Chang & Chen (CS 7 )). Many declared STEM majors at Mercy College were not meeting prerequisites to enroll in Calculus classes. There was also a low sense of community in this almost entirely commuter-based school. This very specific need was addressed by incorporating emporium classes[7], a peer tutoring program, and course coordination to foster a sense of community among the students. College Algebra also acted as a weeding-out course, attributing low enrollment to a lack of preparation. There was also a high turnover ratio and low sense of engagement among the adjunct faculty who taught most of the lower-level mathematics classes. To support students' mathematical knowledge and prepare them to succeed in Calculus, basic math classes were redesigned to make up for the lack of attention, support, and resources from their high school preparation. Sonoma State University faculty recognized that STEM majors were out of reach for underserved students since they would have often needed to complete at least one non-credit developmental mathematics class and the college-level prerequisites to qualify for Calculus (Byrne et al. (CS 5 )). To make STEM classes more accessible, a year-long stretch Calculus course was developed. Students were expected to meet the same level of rigor required in other classes and concepts from Calculus were introduced at the very start of the semester. Similarly, a stretch Calculus class was offered to students at Washington College and placement tests were better enforced (Stacy (CS 23 )).

Talking directly with students is a valuable way to understand the challenges they face and how to best address them. The Precalculus requirement at Colorado State University is fulfilled through five one-credit online courses. This gives the students agency to choose which courses they need to complete, and flexibility in completing the courses. An analysis of student completion showed that STEM students who identified as BIPOC, first-generation, Pell-eligible, and female needed support. This group was targeted specifically due to being systemically underserved and lacking a sense of belonging in STEM. Faculty conjectured that the online-only Precalculus courses exacerbated the anxieties and challenges experienced by these students and created a face-to-face supplemental section to support these students. (Golden et al. (CS 12 )) At Washington College University, when a female black student stated at a workshop that the science building and her mathematics classes were not welcoming, a faculty member investigated and found that students enrolled in Calculus with different abilities and those who were less prepared felt less confident and were alienated (Stacy (CS 23 )). The course was

7 Emporium classes are held in computer labs where students receive their primary instruction and assessments from an online computer program that incorporates text or video lessons with questions to test students' knowledge.

revised into stretch Calculus courses which focused on developing students' academic capital in a comfortable and supportive setting with lessons devoted to active learning. It should be noted, however, that stretch Calculus classes may have a buy-in challenge. Students may worry that a year-long class will hinder their progress through a STEM major. Administrators may also have challenges in scheduling these types of courses. Working with colleagues from the initial conception of the plan can increase support and cooperation among fellow faculty members and administrators.

Though opportunity gaps can be seen in marginalized students, other populations at your institution may also statistically underperform relative to their peers. Transfer students may have a difficult time navigating new systems and may not receive the same level of support as first-year students, such as at Augsburg University (Zobitz et al. (CS 30 )). These students' sense of belonging or existing obligations may prevent them from connecting with other students. Vulnerable populations may require specific programs to foster success and a community spirit.

**2. *Determine the problem that needs to be addressed***

Taking the time to analyze the specific needs of the department and students may provide insight into the types of changes that need to be implemented. Once a target population is chosen, it may be helpful to reflect and ask, "What is the equity question that needs to be addressed? What motivates making improvements in the curriculum? How should this problem be addressed?"

The opportunity gap may not be content related, but could be addressed by considering whether the norms of classes foster an inclusive environment for certain students. At University of Texas at Austin, students' mathematical challenges were viewed as stemming from two roots: strategies for thinking about mathematics and students' sense of belonging in the mathematics classroom (Starbird et al. (CS 24 )). Thus, two side-by-side classes were developed to address each challenge in tandem. Recognizing that these were the challenges experienced by students made the response easier to design. In this case, two courses were offered to students, one on effective thinking that included mathematical proofs, the other an inquiry-based Calculus course designed to foster a sense of belonging. The opportunity gap at Western Michigan University was attributed to barriers associated with possessing academic capital, the set of skills and assets needed to navigate and succeed in educational systems and professional organizations (Mingus et al. (CS 20 )). A task force charged with addressing this issue decided to coordinate instruction among sections through a redesigned course that included changes to content and provided support structures for students. At Florida International University, low passing rates among Latinx and Black students were attributed to traditional lecture-based classes and the artificial norms and boundaries which appeared to alienate and hamper the success of students unaccustomed to the practices in middle-class, suburban, White schools (Fuller et al. (CS 11 )). This problem was addressed by implementing a studio model in Calculus classes where students continuously developed their mathematical knowledge using their experiences through inquiry-based learning.

Course redesign may also necessitate instructing faculty on how to recognize issues of equity and social justice that could be addressed in their classrooms. California State University

Channel Islands developed a professional development program aimed at building a classroom community and creating and implementing rich mathematical tasks in a student-centered environment (Soto et al. (CS 22 )). These aspects were designed to create a more equitable learning environment to address the inequitable STEM retention and graduation rates.

Power dynamics are an example of a dimension of equity which can often be addressed by departmental changes (Johnson et al. (CS 14 )). Power dynamics determine who has a voice or agency, the validity of that voice, and what defines legitimate mathematics in the classroom. The recognition that the only opportunities students had to express their mathematical reasoning was by answering the instructor's answers led to a desire to equitably redistribute power in the classroom so that students had more mathematical agency.

### 3. *Take advantage of opportunities and inspiration*

A catalyst for change can come from a variety of directions. Opportunities to improve classes provide the chance to make changes that may have previously seemed impossible due to institutional inertia or a desire to not push established norms or boundaries. Seizing opportunities for small-scale, "safe" changes can build a foundation for more significant change in the future.

Duke University's grading system in Calculus was historic and storied (Akin & Viel (CS 1 )). Grades in Calculus classes were curved across each course's multi-section common final exam and among each section. In each course, both individual sections and all sections in a given course as a whole were curved based on preset benchmarks related to final exam scores. Lower-level classes had lower-benchmarks. Thus, students in lower-level classes, in which more minoritized students enrolled, received lower grades than they may have expected based upon their experiences and efforts in their individual course section. That is, some students who thought they were doing "well," received lower scores than they anticipated based on the curve. Once again, minoritized students were further marginalized. First, their grades did not reflect their own efforts entirely; they were partially dependent on how well other students performed on the common final exam. Second, these lower grades likely demotivated students and once again had them question their sense of belonging.

When new faculty were hired at Duke University, their energy brought a sense of change to the department. Supported by a new department chair who encouraged experimentation in classes, and a university initiative to improve introductory classes, the ten-year-old grading scheme was changed. The changes to this system resulted in students in lower-level courses receiving grades which more accurately reflect their efforts and level of knowledge rather than on their relative final exam score. Although changes to Duke University's Calculus grading scheme were being slowly implemented, COVID-19 and the subsequent university policy to allow students to choose Pass/Not Pass options provided further opportunities for change.

Grants and funding can provide support to foster course changes. Funding for faculty to develop new lessons and curriculum, compensate adjuncts for time spent in professional development, or to purchase in-class manipulatives allows a department to invest in changes that otherwise may be unaffordable. Grants have opened the door for many course changes, such as at California State University, East Bay, where funding helped incorporate and enact active learning techniques in the classroom (Oliver et al. (CS 21 )), and at West Virginia University where an NSF grant (secured in conjunction with nine other institutions) helped create programs to better support underserved students (Deshler et al. (CS 9 )). At Mercy College, a grant from the Department of Education provided seed funding for their peer tutoring program (Chang & Chen (CS 7 )). The program was later sustained through federal work-study positions, a course fee that covers this and other aspects of the course, and support from the College after the program demonstrated significant success. Small institutional grants can also help evaluate course redesign (Canner et al. (CS 6 )). Augsburg University faculty received a grant so mathematics faculty could work with faculty from other STEM fields to help revitalize the Calculus content (Zobitz et al. (CS 30 )).

Interactions with faculty at other universities can provide inspiration for ideas. The resulting collaborations can also help mitigate trepidation associated with implementing large-scale course changes. The faculty at Florida International University were inspired to redesign their Calculus program after seeing results of other projects and the successes observed in other disciplines (Fuller et al. (CS 11 )). After seeing success at a neighboring institution, Centre College and Southwestern University were inspired to work together and revamp their Calculus program and secured a grant to develop collaborative faculty workshops to re-envision courses (Kilty et al. (CS 15 )). A symposium of institutions in the region to redesign mathematics courses for engagement crystalized three major challenges at Mercy College that needed to be remedied to increase enrollment in Calculus (Chang & Chen (CS 7 )). Faculty at Lamar University consulted with the Dana Center at the University of Texas at Austin to help shape regular faculty meetings to develop a community of practice (Jensen-Vallin et al. (CS 13 )). A regional grant among California State University Channel Islands and three other universities helped fund a professional development program to teach instructors how to better address equity and inclusion in their classroom (Soto et al. (CS 22 )). At Sonoma State University, a workshop on mindset interventions to help change students' perspective on their own success and their mathematical identity helped faculty design interventions for use in their revamped stretch Calculus program (Byrne et al. (CS 5 )).

The campus community and student body can also be a significant basis from which new opportunities can develop. At the University of Texas Rio Grande Valley, the department is establishing a dual-language certificate in mathematics where mathematics classes are offered in both English and Spanish (Villalobos et al. (CS 28 )). This program was based on the recognition that many students are bilingual and that their knowledge of Spanish can be an academic, and not just personal asset. This allows students to be fully bilingual, bicultural, and biliterate. Opportunities abound for ways to support the level of change you believe would better support your students.

### 4. *Identify your theoretical framework*

Once the general problem comes into shape, a theoretical framework to address the problem can provide guidance about how to move forward. Building consensus on guiding principles to address existing educational inequities can help re-envision the intended course and let interested faculty buy into the project, such as was done at Centre College and Southwestern University (Kilty et al. (CS 15 )). A theoretical framework can provide a roadmap to the desired changes. The Deep Theory of Change layering framework developed by Johnson outlines how small incremental changes to existing structures and policies can slowly alter the power dynamics in College Algebra so that student agency is a more central component of the classroom (Johnson et al. (CS 14 )). Small changes accumulate over time to radically transform the teaching and learning dynamics while at the same time, the small nature of each individual change minimizes disruption. Outlining fundamental assumptions about how to approach the problems can serve as a guideline to changes to come and maintain focus on the initial goals of revising a course, such as was done at Sonoma State University, The University of Texas at Austin, and the Transitioning Learners to Calculus in Community Colleges case sites (Burn et al. (CS 4 ); Byrne et al. (CS 5 ); Starbird et al. (CS 24 )). Aside from frameworks to provide general guidance for course revision, theoretical frameworks associated with defining equity and fostering an equitable learning environment have been used to shape redesign efforts.

Defining equity provides an agreed basis from which stakeholders can build upon. In particular, these practices should not privilege one group over another (Adiredja & Andrews-Larson, 2017; Gutiérrez, 2009). California State University Channel Islands' professional development program was framed by and paid specific attention to defining equity (Soto et al. (CS 22 )). Equitable teaching practices were interpreted as those that mitigated systemic differences in how students experienced or gained from educational opportunities. The research project *Transitioning Learners to Calculus in Community Colleges* aimed at studying practices which help African American, Latinx, Native American, and Southeast Asian succeed. To do so, they relied on an anti-deficit framework (Harper, 2010) and practices that supported men of color in community colleges (Wood et al., 2015). Rochelle Gutiérrez' (Gutiérrez, 2009, 2012) dual axis framework defines equitable outcomes in terms of a dominant axis of access and achievement[8]. Additionally, a critical axis of power and identity has been used to define equitable outcomes or as a framework to shape course redesign and evaluation in several case studies, such as at California State Universities at East Bay and Monterey, and a third university located in the Midwestern United States (Canner et al. (CS 6 ); Johnson et al. (CS 14 ); Oliver et al. (CS 21 ))

Pedagogical frameworks, particularly those which promote equity in the classroom, can be used to redesign courses to create opportunities for engagement and participation among traditionally marginalized students. Equitable teaching practices at California State University Channel Islands were envisioned through building a classroom community, student-centered teaching, and creating and implementing rich mathematical tasks (Soto et al. (CS 22 )). These three goals were based on seven teaching practices of a successful Calculus program (Bressoud & Rasmussen, 2015) and the positive, equitable results which come from a framework focusing on social inclusion in the classroom (Estrada et al., 2018). Many challenges arise out of group work, but Complex Instruction (Cohen et al., 1995) can provide a framework to help make this powerful pedagogical tool succeed and ensure an equitable learning environment. This framework was adopted at California State University, Monterey Bay (Canner et al. (CS 6 )). Supplementing Complex Instruction was the Reading Apprenticeship framework (Schoenbach et al., 2012). The goal of Reading Apprenticeship is to help students become readers in a given discipline.

---

8 For a more detailed explanation of Gutiérrez' work see "A College Algebra Intervention to Address Power Dynamics and Promote Students' Reasoning", pages XXX-XXX in this volume.

A guiding framework can help ensure that any changes made are deliberate, necessary, and fall within the scope of the department's vision. Working from a guiding framework can also provide a better means of assessing and troubleshooting the efficacy of the implemented changes.

### 5. ***Take a multi-pronged approach***

There is rarely a single-handed approach to revamping a mathematics course to support underserved students. Many mathematics faculty did not take a unilateral view of what changes were necessary to foster the level of changes they felt were necessary and instead relied on multiple approaches to achieve their goals.

At California State University, East Bay, a three-pronged approach was outlined to help support BIPOC students in Calculus I: course redesign and curation, opt-in lessons and coordinated calendar, and the cultivation of a robust community of practice (Oliver et al. (CS 21 )). When faculty at University of Michigan saw that fewer minority students enrolled in higher-level Calculus courses, they developed a two-dimensional program to revamp their classes (Bennett et al. (CS 3 )). First, an intensive week-long training program along with ongoing support was provided to instructors every year. Second, assessments were revised to account for mastery learning. Faculty at Florida International University (Fuller et al. (CS 11 )) incorporated an active learning-based approach, culturally relevant pedagogy, and group work into their classrooms. They also developed a program of near peer support, where students were supported by peers who were only slightly advanced in their mathematics courses. In Augsburg University, revisions to Calculus focused in two general areas, classroom practices and pedagogy, and community building (Zobitz et al. (CS 30 )). The content was revamped to focus on application problems from other departments. Class size was deliberately decreased so that IBL exploratory problems could be done in small groups. Weekly Calculus workshops to work on harder problems and build community were initiated. A first-year seminar was linked to a Calculus class. Connections with alumni in STEM fields were fostered. Finally, added support was given to transfer students from community colleges in this multi-dimensional endeavor. Sonoma State University's stretch Calculus program incorporated mindset interventions, active learning through group work, and peer support (Byrne et al. (CS 5 )). The department at California State University, Long Beach (Benken et al. (CS 2 )) redesigned their Calculus course, developed an intrusive advising program (an early alert program and mandatory academic advising if necessary), rethought their placement policies, and developed a summer program in response to inequitable success rates and the elimination of developmental courses. Mathematics content instruction was intertwined with helping students develop their growth mindset, metacognition, STEM identity, and a sense of belonging.

Multi-pronged approaches are necessary for target populations that are not monolithic where a single approach can be reductive or insufficient to meet students' needs. Multiple modes of support help students find a program which can best meet their needs, while building on their strengths. More importantly, recognizing the varied needs of students brings a newfound level of compassion to teaching.

### 6. *Think beyond the math*

It is tempting to think that revising a syllabus, changing the textbook, or implementing new student-centered lessons may be enough to support students. However, oftentimes affective, non-content based factors must be considered to fully realize and address how the structure of the course supports or undermines specific students.

Mathematics is an intimidating subject. Undergraduates may be haunted by past failures or negative experiences in previous mathematics classes. Consciously building in "kindness cues" and "welcomeness to engage" practices to promote social inclusion is a recognition and a step towards rehumanizing the mathematics classroom and can be a major step towards promoting equity while affirming the worth of the student on a micro and macro level (Burn et al. (CS 4 ); Soto et al. (CS 22 )). Likewise, students greatly benefit from specific practices or actions that work to deliberately foster a growth mindset in the classroom such as instituting norms to normalize making mistakes so that mistakes are a part of learning, such as at several case universities in this volume (Benken et al. (CS 2 ); Bennett et al. (CS 3 ); Burn et al. (CS 4 ); Byrne et al. (CS 5 ); Canner et al. (CS6 ); Soto et al. (CS 22 ); Stacy (CS 23 )). Changing students' perceptions can be a fundamental objective when redesigning classes. At University of Texas at Austin, helping students see themselves as producers of knowledge rather than consumers of knowledge was the goal of both the Elements of Effective Thinking class and the linked Calculus class (Starbird et al. (CS 24 )). To increase persistence in Calculus and decrease attrition in STEM, Centre College and Southwestern University consciously worked to improve student efficacy in mathematics and to increase their understanding of the usefulness of mathematical skills (Kilty et al. (CS 15 )).

This desire to focus on auxiliary factors can be student driven. Student feedback about the supplemental course created to support students showed that students were not initially satisfied with the course structure (Golden et al. (CS 12 )). Initial feedback showed that students wanted a course to help them succeed in the program, not a course that focused on their mathematical development. In response to this feedback, the course was changed so that the teaching goals were geared towards students' needs: teaching test-taking strategies, different problem-solving methods, and how to use a graphing calculator. Content was not emphasized; the skills that support academic success were. Similarly, the effective thinking course at University of Texas at Austin (Starbird et al. ((CS 24 )) was designed to help students develop high-level study skills such as how to effectively learn from one's mistakes in the Calculus class. The Summer Bridge program at St. Joseph's University was not strictly focused on helping students work through a set curriculum (Terry (CS 26 )). Instead, the Summer Bridge program helped students manage the heavy workload, adjust their perceptions, and increase their self-efficacy. Various student skills were addressed such as how to take notes, how to write lab reports, and how to understand faculty expectations. A workshop was also included that was designed to empower students to be more knowledgeable about university resources and services; knowing more about university services can help students manage their self-efficacy and socio-emotional skills to contribute to their overall academic success. Such topics are particularly beneficial for first-generation students and students who come from under-served schools.

The long-term career benefits of a college degree can be a strong motivator for students. However, student motivation can decrease if the connection between courses and careers is tenuous or not understood. Augsburg University linked a first-year seminar to Calculus I to help students understand how Calculus skills can lead to financial success (Zobitz et al. (CS 30 )). Recognizing that students "can't be what they can't see", BIPOC, LGBTQIA, and female alumni in STEM fields were recruited into a mentoring program. This helped students see the direct benefits of their classes[9]. This program also discussed implicit biases in and out of mathematics, intercultural communication, financial aid, study abroad opportunities, advising, and mental health.

Along with reconsidering content or course structure, the major framework used to design the stretch Calculus course at Sonoma State University focused more on student affect than on content knowledge (Byrne et al. (CS 5 )). This was based on the assumption that when factoring for success, existing mathematical knowledge comes second to mindset, confidence, and a sense of belonging; mathematical knowledge can be built over time if these other factors are in place. The framework for the course redesign was based on agency, authority, identity, and equitable access to content. These dimensions assumed that access and mindset contribute more to mathematical success than previously measured mathematical skills; content came second to affect. Affect-oriented perspectives can be interspersed throughout a mathematics course such as by providing students with a weekly checklist that outlines required tasks and provides suggestions to enhance learning (Mingus et al. (CS 20 )). These types of mechanisms can help support students' development of their academic capital.

Attending to affective factors which contribute to success should be considered within a larger equity framework and should be carefully selected so that activities that are intended to be helpful do not become counterproductive. For example, at Appalachian State University, students traditionally responded to prompts such as "I attend every class", "I ask questions during class", "I read the course textbook", and "I read my course notes". These prompts were intended to help students recognize that they are responsible for their own progress. However, these prompts may trigger a sense of shame or disengagement with BIPOC, low-income, rural, and first-generation students because the student's responses may be for reasons beyond his or her control. Traditionally underserved students may not have the privilege of affording the textbook or attending every class. Cultural factors such as a deference to authority, a students' lived experience with racism, or microaggressions can make the vulnerable act of asking questions difficult (Burn et al. (CS 4 )). Such activities should be designed to support students' sense of agency while recognizing the larger societal factors they are working against.

9 For more on the benefits of role models see the chapter "A Focus on DEI for Students Yields DEI for Instructors" pages XXX-XXX

## 7. *Build a community to help inspire and thrive*

The drive towards enhancing diversity is not just about righting sociohistorical wrongs, it is also about calling upon the strength of multiple perspectives to make larger advances towards a goal. We are smarter together and greater strides are made when multiple ideas, experiences, and knowledge bases come together to tackle a problem. A diverse group of individuals unencumbered by bias or asymmetrical power relationships are more creative, innovative, and productive (Smith-Doerr, Alegria, & Sacco, 2017; Wooley, Chabris, Pentland, Hashmi, & Malone, 2010). Thus, faculty looking to revamp their courses should look to others for guidance, advice, and input. Furthermore, if the overall goal is to support nonhegemonic students, faculty would benefit from learning to listen to, work with, and implement the suggestions of others who can provide valuable insights they may not have previously considered.

### 7.1. *Build Communities Among Faculty and other Stakeholders*

The larger university community can be a strong resource to tap as you endeavor to change your mathematics classes. As many mathematics classes, and in particular Calculus classes, are intended to serve other departments, these departments should be consulted to provide input into specific course objectives, programming, and to help provide real-world applications and problems for students (Kilty et al. (CS 15 ); Terry (CS 26 ); Zobitz et al. (CS 30 )). Look to your institution's central resource for teaching to help develop revitalize pedagogical practices or to more objectively evaluate your program (Bennett et al. (CS 3 ); Canner et al. (CS 6 ); Chang & Chen (CS 7 ); Fuller et al. (CS 11 )). Aligning with and educating academic advisors on the affordances of specific programs or special sections can increase student enrollment and participation (Mawhinney et al. (CS 18 )). Your Office of Institutional Effectiveness and administrative support staff can help assess the effectiveness of policies on a longitudinal basis and potentially raise the profile of related endeavors. This in turn may lead to generating more institutional support and resources (Benken et al. (CS 2 ); Canner et al. (CS 6 ); Chang & Chen (CS 7 ); Oliver et al. (CS 21 )). Bringing together a team of stakeholders, instructors, and experts from the mathematics department and the school of education can maximize individual investment and minimize objections to proposed changes (Johnson et al. (CS 14 )).

Developing a community of practice strengthens and supports the pedagogical skills of participating mathematics faculty (Fuller et al. (CS 11 ); Jensen-Vallin et al. (CS 13 )). When instructors work collectively, such as through weekly meetings, everyone has the opportunity to contribute to the process and feel invested in the changes. This also gives adjunct professors and instructors the opportunity to share concerns and suggestions, giving voice to those who often teach students but may have limited agency in the structure of their courses. At Lamar University, community interaction has been strengthened through social teas and luncheons with topical themes that are aimed at improving teaching, learning technology, and teaching towards equity (Jensen-Vallin et al. (CS 13 )). Participants in California State University Channel Islands' professional development program gave adjunct professors who felt invisible and voiceless a supportive community in which they could interact and learn (Soto et al. (CS 22 )). Having new or inexperienced instructors co-teach allowed instructors to master new teaching techniques and develop camaraderie (Byrne et al. (CS 5 )).

When multiple sections of a course are offered, it is more difficult to standardize students' learning experience. Differences in instructor styles and expectations are certainly expected. However, it is only when these differences result in significant inconsistencies such as rigor or covered content, that students can be subject to an inequitable educational experience (Akin & Viel (CS 1 ); Chang & Chen (CS 7 ); Jensen-Vallin et al. (CS 13 ); Mingus et al. (CS 20 )). Inequitable class structures may also result in an inequitable workload among the instructors. For example, one instructor may be favored because students have a better classroom experience with them or because they have lower grading expectations. Streamlining and synchronizing course content is an opportunity to make the course more uniform among sections. Coordinating a course can also provide an opportunity to develop a community among faculty if they are a part of the coordination effort by co-creating assessments, rubrics, and coordinating pacing (Mingus et al. (CS 20 )). Variation across instructors provides an opportunity to collaborate, share best practices, and present a cohesive course to students. A dynamic calendar tracking course content, which can be edited easily by instructors and contain links to resources associated with each lesson, encourages instructors to coordinate schedules and lessons across different sections in a class (Oliver et al. (CS 21 )). Materials specific to fostering equity in the classroom and ways to introduce "just in time" review topics (where prerequisite content is reviewed as needed) can also be highlighted and shared. Creating a repository of materials for instructors allows for more consistency across instructors and is an easy way to reference and share resources (Chang & Chen (CS 7 ); Johnson et al. (CS 14 ); Oliver et al. (CS 21 )). This system works particularly well for inexperienced instructors or adjuncts who may not have the time to develop their own resources, be familiar with the nuances of the institution's specific needs and content or meet with faculty or course coordinators in their department.

Regular meetings among instructors can be used to discuss issues of pedagogy, assessment, pacing, design and use of activities, how best to provide effective support for instructors, and share resources (Bennett et al. (CS 3 ); Byrne et al. (CS 5 ); Chang & Chen (CS 7 ); Golden et al. (CS 12 ); Oliver et al. (CS 21 )). These meetings are an opportunity to explore issues of equity including opportunity gaps in classes with multiple sections or differences in instructor grading criteria, and to discuss common readings on inclusive teaching and general strategies around creating a more equitable learning environment. These meetings can also allow instructors to celebrate their successes and support one another through challenges.

### 7.2. *Build community between faculty and students and among students*

There are various types of effective pedagogical practices in mathematics such as complex instruction, project-based learning, inquiry-based learning, and other forms of student-centered instruction. One unifying characteristic in how these practices compare with an instructor-centered lecture is the bilateral nature of the student-instructor relationship. Traditionally, an instructor lectures *at* students and students' responses are in passive forms of communication such as written homework, quizzes, and exams. In the aforementioned alternatives, the instructor works *with* the student to nurture their curiosity and foster the development of their knowledge.

One of the most common changes that was made in the Calculus classes in the above case studies is the inclusion of a program or aspect which recognized the importance of including the students as an active member of the classroom community. Giving students a voice, space, and agency to control their experience is a major component of many redesign endeavors. This may be as “simple” as redeveloping content so that classroom content is presented in a more student-centered manner, or as “radical” as giving students a larger role in building the overall mathematics community.

Student satisfaction and input can provide an indication that change is necessary. At times, students can provide insight into problematic policies or situations. At Duke University, students indicated that the grading scheme was both non-transparent and inequitable because final grades did not accurately reflect each student’s individual knowledge or abilities. This left students feeling demoralized and discouraged from pursuing their intended STEM majors (Akin & Viel (CS 1 )). Mentoring students and listening to them can illuminate some of the challenges and misunderstandings that students have which can result in students believing that they are not capable of completing a STEM degree (Terry (CS 26 )). Personal interactions provide frank discussions on some of the limitations experienced by students. The resulting summer bridge program was designed to address problems that contributed to BIPOC students abandoning their STEM major. These conversations can initiate significant change. In another aforementioned example, a student was the catalyst for change when faculty realized that some students did not feel welcome in their mathematics class (Stacy (CS 23 )). Student feedback can also help divert resources to address students’ actual learning needs rather than simply focusing on mathematical remediation (Golden et al. (CS 12 )). The department gave students agency through a bottom-up approach to change.

Messaging is an important component to community building. Emphasizing the creative aspects of mathematical practices can increase the self-efficacy of someone who does not believe they are a “math person” but can nevertheless succeed in mathematics based on their creative abilities (Stacy (CS 23 )). Demonstrating students’ worth demonstrates that the students are valued and are a part of the community—a critical factor in STEM attrition (Miller, Williams, & Silberstein, 2019). Small, conscious acts such as learning students’ names and encouraging them to use chat features during online meetings, make students comfortable and can significantly increase their sense of belonging (Oliver et al. (CS 21 )). Social events such as holding a tea can foster community (Soto et al. (CS22 )). In one case, students and teaching assistants are invited to a professor’s home (Starbird et al. (CS 24 )). Faculty can also help students navigate the hidden curriculum4F[10] (Stacy (CS 23 )). At the University of Texas at Austin, the students’ ability to succeed is a fundamental belief for instructors. Students in their program are explicitly and repeatedly encouraged to major in mathematics and struggling students are actively supported (Starbird et al. (CS 24 )). Professional development around student power dynamics also can help instructors reflect on how subtle messaging can have significant ramifications (Johnson et al. (CS 14 )). One such exercise, which helps instructors consider which students they connected with, was designed to demonstrate how intentional and unintentional actions could affect students' sense of belonging, achievement, and persistence in mathematics. In other cases of professional development, instructors can be made aware of how their actions engage or alienate students (Oliver et al. (CS 21 )). If an instructor calls on the same students, they deny others an opportunity to participate, illustrating how actions such as cold-calling should be discouraged.

[10] The Hidden Curriculum is the concept that students learn more than the formal content when in a classroom. There are unwritten, unofficial, and at times unintended ideas, values, and perspectives that students learn. To further understand the hidden curriculum, see Portelli (1993). To see how the hidden curriculum is manifested in the classroom see Anyon (1980).

Students can also feel valued if they are looked to as a source of knowledge and support. Programs that incorporate peer learning in and out of the classroom help build community and show the peer tutors that they are a valued resource (Benken et al. (CS 2 ); Chang & Chen (CS 7 ); Fuller et al. (CS 11 ); Mingus et al. (CS 20 )). Peers can also support students' emotional well-being by acting as a role model and by helping others navigate their institution's hidden curriculum (Mingus et al. (CS 20 ); Zobitz et al. (CS 30 )). There is also the added benefit that the tutors also gain academically from their assignments (Leung, 2019). The stretch Calculus class at Washington College assigns each student a learning objective from the course. Each student is then responsible for reviewing this objective with the class in the manner of their choice (Stacy (CS 23 )). Culturally responsive teaching likewise underscores the personal worth of a student's knowledge base and experience and makes mathematical content more relevant, accessible, and meaningful (Fuller et al. (CS 11 )).

A welcoming environment can take various forms. Recognizing that there is a lack of community in an HSI composed mostly of commuter students can provide insight into the types of small changes which can make significant differences to students' engagement and achievement (Benken et al. (CS 2 ); Chang & Chen (CS 7 )). Other forms of messaging can be more subtle such as incorporating kindness cues (Soto et al. (CS 22 )), cultivating a sense of pride in students' linguistic and cultural heritage and acknowledging that being bilingual is an asset and a source of strength can be very beneficial (Villalobos et al. (CS 28 )). A welcoming environment can also recruit students to support programs. Specific language that lets students know they are included in a support program can be helpful. For example, inclusion based on a "nomination" process indicates that their presence is important and honored. Negatively-toned and deficit-based language about placement or progress should be eliminated in lieu of messaging that highlights benefits and the positive aspects of programs (Deshler et al. (CS 9 ); Mawhinney et al. (CS 18 )).

Having students consistently work in groups can help advance mathematical knowledge while developing a sense of community (Deshler et al. (CS 9 ); Fuller et al. (CS 11 ); Mawhinney et al. (CS 18 ); Mingus et al. (CS 20 ); Zobitz et al. (CS 30 )). Community can be further strengthened when students are assessed as a group and must rely on one another to be successful (Starbird et al. (CS 24 )). Fostering community can also be a long-term endeavor such as at the University of West Virginia where one instructor dedicated themselves to teach all four Calculus courses in a sequence to a cohort (Deshler et al. (CS 9 )). These faculty also provided small acts of meaningful support such as helping students with registration changes in their first semester, removing challenges faced by new students who are unfamiliar with how to navigate registration systems. Policies which hinder the development of community should also be reconsidered, as was done at Duke University when evaluations determined that the curved grading scheme increased competition among students and decreased student camaraderie and sense of belonging (Akin & Viel (CS 1 )).

The mathematics classroom can also be a way to help students feel welcome in the larger campus community. Students can be required to complete assignments where they interact with various services on campus, such as the university writing center and library (Starbird et al. (CS 24 )). Such purposeful acts teach students that success comes from taking advantage of all opportunities and available resources.

### 8. *Rethink assessments*

When courses undergo significant redesign, traditional assessments such as quizzes and exams may no longer be appropriate and may undermine the type of cultural change that may be necessary to appropriately support and educate students. Scaling or curving exams also disadvantages women and students of color (Akin & Viel (CS 1 ); Bennett et al. (CS 3)). In addition, the COVID-19 pandemic demonstrated that relying on traditional assessments is problematic when faculty need to be more nimble and respond to fast-changing situations. When mathematical reasoning is emphasized over answer finding, concerns over academic honesty may not be a byproduct of online assessments (Johnson et al. (CS 14)).

Writing assignments are one form of assessment that can help students clarify their thinking. Students can be asked to write papers explaining central ideas and key concepts in Calculus, or they can outline attempts to solve a specific puzzle or problem rather than actually finding a solution (Mingus et al. (CS 20 ); Stacy (CS 23 ); Starbird et al. (CS 24 )). Students can also create videos explaining Calculus concepts and can complete a final project with both a written and oral component. Washington College has students create a comprehensive review guide as a final course project (Stacy (CS 23 )).

Mastery-based assessments give students multiple opportunities to demonstrate their skills and knowledge (Akin & Viel (CS 1 ); Bennett et al. (CS 3 ); Byrne et al. (CS 5 ); Stacy (CS 23 )). Such exams provide students an opportunity to review difficult ideas and reassess their skills (Starbird et al. (CS 24 )). They have an added benefit of developing a growth mindset among students, which can reduce stereotype threat and counterbalance the effect of lower socioeconomic status on achievement (Aronson et al., 2002; Claro et al., 2016). Traditional assessments can be modified to be more equitable, particularly among multi-section classes. At University of Michigan, providing more frequent, mastery-based assessments gave instructors and students regular feedback on their progress so they could make more informed decisions about whether to remain or drop a class (Bennett et al. (CS 3 )). Furthermore, the mastery-based assessments were available for students to take over a two-week period. This provided students with more flexibility and was responsive to their individual and personal needs, providing for a more equitable learning environment.

Among the regular discussions at California State University, East Bay, it was discovered that instructors' policy on allowing notes during exams varied (Oliver et al. (CS 21 )). This inequity was resolved when a policy was established to allow students a 3x5 index card on all assessments. Not only does this provide an equitable policy across classes, but it encourages instructors to include questions which require conceptual rather than procedural knowledge.

Including more low-stakes assessments can also reduce anxiety and provide students with more regular feedback (Bennett et al. (CS 3 ); Oliver et al. (CS 21 )). Overall, all assessments should be evaluated to ensure they align closely with course objectives and content (Bennett et al. (CS 3 )).

Placement exams can have significant effects on students' long-term successes in mathematics. Placement into a developmental or lower-level college course can make a Calculus class, and therefore a STEM major, unattainable. Incorrectly placing a student into a course lower than appropriate, particularly for students of color, can affect students' sense of belonging; conversely, good placement challenges the student and can positively influence their success (Burn et al. (CS 4 )). Multiple options for placement can provide a more accurate reflection of a student's ability (Benken et al. (CS 2 ); Burn et al. (CS 4 )). Nuanced placement procedures take time to develop and implement and may require last-minute revisions to course schedules. Thus, training and dedicated time in the summer is necessary to effectively ensure student placement (Benken et al. (CS 2 )). Similarly, providing students with opportunities to "test out" of developmental classes at the start of the semester gives them more agency and can positively affect their self-efficacy (Burn et al. (CS 4 )). Ensuring that the placement exam accurately reflects the needs of the department is also crucial. Higher cut-offs for placement exams can prevent students with weaker algebra and precalculus skills from being overwhelmed in Calculus (Mawhinney et al. (CS 18 )). Interventions which can provide topical or just in time support can expedite a student's progress by refreshing and meeting standards for specific topics rather than requiring the student to take a semester-length course (Benken et al. (CS 2 ); Mawhinney et al. (CS 18 )). At Appalachian State University the department redeveloped their policy for placing students to ensure that every student was taking a class that challenged them appropriately. When a student enrolled in Calculus 1, he or she had access to several mechanisms for support including a corequisite course (Mawhinney et al. (CS 18 )).

### 1. *Rethinking everything*

Non-Euclidean geometry developed from questioning the necessity of the parallel postulate. Netflix's rise can be attributed to questioning the assumption that late fees are a necessary part of movie rentals. One factor that led to the 2008 financial crisis was the assumption that real estate never lost value. Revolutionary changes sometimes result from questioning the reality with which we are presented and the assumptions that are embedded in that reality. Questioning all aspects of how we teach mathematics should underlie any reform efforts, particularly those intended to foster a more equitable and inclusive learning environment. It may not be obvious that the smallest traditions or practices may reinforce a culture of White supremacy or make someone feel excluded, nevertheless this is a reality with which we must all reckon. Considering whether the most mundane factors contribute or hinder student success and making radical changes can make a revised mathematics class truly transformative. The ramifications of such changes can have profound effects on students, faculty, and administrators. For example, altering the infamous grading system at Duke led to a more equitable learning environment. In this case, equity was the result, rather than the catalyst, of questioning the existence of a grading tradition and then making changes to the system.

Some departments have stretched their courses over multiple semesters to provide academic and affective support for students (Benken et al. (CS 2 ); Byrne et al. (CS 5 ); Stacy (CS 23 )). At West Virginia University (Deshler et al. (CS 9 )) class meetings were held for longer periods of time over fewer days to provide a better setting for inquiry-based learning. Although the revised meeting times resulted in more time in the classroom, students appreciated meeting less often because it freed up their schedule for other classes.

Keeping in mind that mathematics is often a service subject that prepares students to succeed in STEM and other disciplines, some departments (Kilty et al. (CS 15 ); Zobitz et al. (CS 30 )) have eliminated a focus on rote algebra and paired their course down to essential concepts and skills. This led to students being more successful by eliminating the burden brought about by difficult algebraic manipulations. They could focus more fully on fundamental concepts and ideas. Essential concepts necessary for success in STEM classes or the longer-term needs of the students were kept and explored more deeply. Those who ultimately majored in mathematics were able to learn concepts, skills, and theorems at a later point in time.

Another challenge that departments often contend with is what to do with students who enroll in an introductory level course and do not complete the semester or need support during the semester. At Appalachian State University an alternative course was developed for those who started in Calculus but needed to withdraw during the semester (Mawhinney et al. (CS 18 )). Rather than drop and repeat the course the next semester, students enrolled in a fallback course offered in the middle of the semester. This helped those students focus on algebraic techniques that would be useful when they repeated the Calculus class. Courses designed to supplement or support courses could start during the semester or could include voluntary review sessions held at the same time other days of the week (Golden et al. (CS 12 )).

Reflecting on what is valued and what is not in a mathematics classroom can determine how inclusive the classroom is for students. English language dominance in a mathematics classroom is not necessary if multiple languages and modes of communication among students can foster a sense of community and enhance their ability to better understand mathematical concepts (Villalobos et al. (CS 28 )).

The physical classroom environment is also a potential source for reform. At West Virginia University (Deshler et al. (CS 9 )), the Emerging Scholars Program secured a dedicated classroom and study space. The room was renovated to be more conducive to group activities by incorporating moveable furniture, large whiteboards, markers, and erasers for students to use in class. Washington College's stretch Calculus class is currently held in a beautiful, historic house with comfortable couches, snacks, and other amenities not usually found in a classroom (Stacy (CS 23 )). Dedicated spaces for studying mathematics and tutoring can foster a stronger sense of belonging while at the same time helping develop mathematical knowledge. These spaces can also be used to help students prepare for their next-level mathematics course (Burn et al. (CS 4 )).

Some departments re-envisioned their courses by doubling down on the things they thought would be most helpful. Several departments changed course schedules to provide instructors with the time they need to explore content in a more meaningful way or veer away from a traditional Calculus curriculum (Benken et al. (CS 2 ); Byrne et al. (CS 5 ); Kilty et al. (CS 15 )). In such a case, it is beneficial to ensure that faculty concerns are expressed and addressed, as occurred when the topic of limits was moved to a later course in a Calculus sequence (Kilty et al. (CS 15 )). Some faculty were uncomfortable with such a radical transformation, but the eventual consensus was that changing when the topic appeared did not compromise rigor but provided greater context and deeper understanding.

### 9. *Revise revise revise*

A final note should be remembered when revising classes to foster a more equitable learning environment. The more transformative the proposed changes, the more likely it is that revisions to the original format will be needed. Assess how changes to course structures align with the developed theoretical framework or address the outlined problem established at the beginning of the process. From there, expect that revisions and changes to the initial plan will have to be considered. Faculty at Augsburg University regularly adopted active learning into their classes, valued the scholarship of teaching and learning, and consistently developed their knowledge based on issues of equity and social justice. Nevertheless, demographic shifts among the student body and anecdotal statements on BIPOC, transfer, and first-generation students' anxieties over their mathematics classes led to implementing five significant changes to their Calculus classes (Zobitz et al. (CS 30 )). Changes were both proactive and reactive to students' needs. These faculty members were open to students' feedback that more work was needed to foster a more inclusive and supportive learning environment. The changes outlined by California State University, Long Beach (Benken et al. (CS 2 )) have been a part of an ongoing effort to redesign lower-level mathematics classes. Evaluating passing rates and longitudinal analysis of student success resulted in offering a stretched version of Precalculus. Patience and compassion are necessary for students who are trying to master a difficult subject, but faculty and administrators undergoing this important endeavor also need patience and compassion.

The effectiveness of changes in course design and content to ensure more equitable outcomes can be evaluated using multiple sources of data (Benken et al. (CS 2 ); Canner et al. (CS 6 )). Demographic data and student grades can be compiled from your Office of Institutional Effectiveness. Interviews and focus groups, attitude surveys given at the beginning and end of a semester, and an experience survey given to students to gain insight into how course changes affect students' attitudes and their sense of belonging provide rich context for more nuanced results to identify areas of improvement. Triangulating evaluative data can also help motivate students. When surveyed students indicated that they did not enjoy the active learning component of their classes or think it was useful, their higher grades provided a counternarrative which was shared with subsequent cohorts (Canner et al. (CS 6 )). Furthermore, the evaluation process can help uncover positive unintended consequences of various programs. At Mercy College, students who were recruited to the peer tutoring program and worked in classes supporting students found the experience so positive, they enrolled into a 5-year Master's Degree program in STEM education (Chang & Chen (CS 7 )). Thus, a secondary problem of diversifying the teaching force is also being addressed through course redesign.

In many of the accompanying case studies, initial ideas and programs were revised once feedback demonstrated that the policy was not effective. Struggling students at the University of Texas at Austin were advised and offered the opportunity to transfer from Calculus to Precalculus (Starbird et al. (CS 24 )). Most students were not willing to leave the class and remained in Calculus. However, these students ultimately succeeded in the class with a C or better. The choice to transfer was thereafter eliminated. Instructor meetings among mathematics faculty at Lamar University initially weren't successful due to low attendance which prevented the intended standardization of pacing across sections (Jensen-Vallin et al. (CS 13 )). In the following term, meetings had preset scheduled times, norms, expectations, and a template agenda to ensure regular, smooth-running meetings. Similarly low attendance in supplemental classes stopped being a problem when attendance was linked to a small component of the grade (Golden et al. (CS 12 )).

Evaluating how a course helps students' achievement may provide limited insight into changes that need to be made. At Mercy College's emporium classes, students who did not own a computer used a smartphone. Compared to computers, the phones had limited capabilities and many distractions (Chang & Chen (CS 7 )). A lab fee was instituted to give students access to computers. Additionally, many students worked full-time and did not have time to work independently in the flipped class format. Thus, the format of the class was modified so it would be a quasi-flipped class. Longitudinal data suggest some areas of improvement. In particular, students who were successful in the emporium format seemed to struggle in subsequent courses, while others preferred conventional lecturing over the emporium model. These policies were instituted to foster a community among students, however the reliance on computer-centered learning may limit students' abilities to construct mathematical knowledge rather than respond to rote procedural questions.

Consider that changes do not need to be immediate and radical. Smaller, more careful changes may be more suited to a department's culture. At Duke, revising the grading scheme in Calculus was a slow, deliberative and iterative process (Akin & Viel (CS 1 )). Departmental structures encouraged incremental change and inertia. Experienced instructors and teaching assistants were reluctant to make changes. However, with time and patience the grading system has been slowly changing. At Washington College, senior faculty were reluctant to consider making drastic changes to the mathematics curriculum (Stacy (CS 23 )). The lack of faculty buy-in was overcome when students showed significant support for the changes proposed. Although the senior faculty had reservations about making the changes, they were open to changing their minds once they saw the level of support the students had for these course changes.

### 2. *The Risk to Blossom*

The opening quote, attributed to Anaïs Nin[11], states how the pain of inaction can be overwhelming. Remaining ambivalent to the needs of students and expecting students to bend to an archaic system may be easier and safer for you the faculty member, but it harms students and maintains systemic inequities that must be overhauled if we are to improve our society. Inaction ensures that cohort after cohort of students enroll and only a few advance. Inaction results in students continuing to hate mathematics. Inaction maintains the myth of being "a math person." Inaction muddies the beauty and magnificence of mathematics for so many. Inaction upholds the school to prison pipeline. Inaction prevents others with new perspectives from participating and expanding fields of research. Inaction allows misogyny, racism (systemic and overt), and a culture of White supremacy to be further entrenched.

[11] Although this quote is often attributed to Nin, the reference seems to be circular, and the author has not found a text written by Nin containing these words. An online blog contains a post where Elizabeth Appell claims to have written the poem "Risk" in 1979 to inspire and motivate students at John F. Kennedy University in Orinda. (http://anaisninblog.skybluepress.com/2013/03/who-wrote-risk-is-the-mystery-solved/)

There is another choice, the challenging and more difficult path of change. It is important to recognize that transformative changes are not easy. These changes are personal; they require a need to reckon with a reality that may be difficult to face. These changes are difficult; they require time, dedication, and strong relationships. These changes are slow; reflecting on the needs of students, planning changes, and deliberating over outcomes requires patience. Every department, every student body, every institution has its own set of challenges and affordances which should be the basis to which course transformation should flow. Changes made to support student learning may not be drastic, but they may be significant in students' overall mathematical success.

Although the pain that comes with taking the risk to blossom may be less, blossoming invites danger; blossoming elicits fear. Mitigate that fear by outlining a clear plan of transformation. Mitigate that fear by ensuring administrative support and through careful communication with other departments. Mitigate that fear by pausing to honestly reflect and reassess the transformative process. Mitigate that fear by allowing you and your department to be vulnerable, to recognize a change in the status quo, and to welcome it. Mitigate that fear by recognizing that vulnerability leads to innovation. That innovation will lead to a blossoming.

### 3. *Reflect and plan*

Once you are motivated by action, it's best to take a moment and begin to reflect on how any of these phases can help you achieve your goal. Here are a few questions to consider when you begin this transformative process.

Intentionally search for equity

- What assumptions do you hold about students' abilities? Are those assumptions valid?
- Which groups of students are most marginalized in your mathematics classroom? Remember, marginalization may not be related to grades and outcomes.

Determine the problem that needs to be addressed

- Aside from outcomes or STEM attrition rates, what are your expectations for improvement?
- What traditions or practices are held in your department? How may they contribute to making some students feel marginalized?

Take advantage of opportunities and inspiration

- What are your students' strengths that you can use to foster their mathematical growth?
- What experiences do new and visiting faculty bring that you can build upon?

Identify your theoretical framework

- What researchers or visionaries are sources for inspiration? Consider inspiration from areas outside your field such as K-12 teaching, business, or self-help.

- What is your timeline of implementation? What outcomes do you hope to achieve by then?

Take a multi-pronged approach

- List ideas to target three different areas for improvement in your department or courses. How could each idea also help improve the other two areas?
- What other benefits do you see from these ideas both in and beyond the mathematics classroom?

Think beyond the math

- What sociocultural factors can help make your students more successful? How can you weave those skills or beliefs into your classes?
- How can you support every student's "whole being"?

Build a community to help inspire and thrive

- What are some of your department's and institution's strengths that you can build upon?
- What can local leaders and community members do to help support your students' mathematical knowledge?

Rethink assessments

- If you replaced the word "test" with "assessment," "accounting", or another word, how would your "test" change?
- If a traditional test is necessary, is it necessary to limit students' time, their use of resources, or the number of opportunities they have to demonstrate their knowledge? How can you modify a traditional test to more fully assess your students' knowledge?

Rethink everything

- In your students' (and your) ideal world, with unlimited resources and support, what does a successful mathematics class look like? How can you bring this vision to life?
- What strengths/barriers exist to foster/prevent your ideal implementation and outcomes?

Revise revise revise

- What changes or improvements do you anticipate to occur after one semester or year? Were your expectations realistic? What modifications do you need to consider?
- What reactions or suggestions do your students have after changes have been made?

# 6 Leveraging Identity and Language to Promote DEI in Undergraduate Mathematics

**Ernesto Daniel Calleros**
*San Diego State University*
*University of California San Diego*

Ernesto Daniel Calleros is a Ph.D. candidate in the Mathematics and Science Education joint doctoral program at San Diego State University and the University of California San Diego. He earned an M.A. in Mathematics from Rice University. His research interests include undergraduate mathematics education with a focus on linear algebra. He also investigates ways to make mathematics linguistically accessible to all students, especially multilingual learners. Prior to the Ph.D. program, Ernesto taught college-level mathematics.

**William Zahner**
*San Diego State University*

William Zahner, Ph.D., is an Associate Professor in the department of Mathematics and Statistics at San Diego State University. Zahner teaches mathematics courses for prospective secondary teachers and is the principal investigator of research projects focused on language diversity and mathematics. Prior to his work as a faculty member, Zahner was a high school mathematics teacher in Chuuk, Federated States of Micronesia and San Jose, California.

## 1. *Introduction*

The case studies in this volume describe efforts to transform calculus at a wide variety of institutions. Transformation is needed because recruiting and retaining students reflecting the diversity of the US population into mathematics (and STEM, which often has calculus as a gateway) continues to be problematic in the United States. Research shows that choosing STEM majors and careers is correlated with developing a positive identity toward STEM disciplines (Godwin et al., 2016; Hazari et al., 2010). Hence, increased STEM major participation and retention rates are linked to positive identity formation within a given STEM field (Carlone & Johnson, 2007). This leads to the question, what is identity, and how does one form an identity in relation to an academic discipline?

In this chapter, we will begin with a brief literature review of the work on identity in mathematics education, including a discussion of different conceptualizations of identity and their implications. This review will highlight how students' lived experiences related to gender, race, language, and other socially meaningful categories may affect students' identity formation toward mathematics (or STEM more broadly). We will then showcase how the concept of identity was leveraged in the case studies in this volume to transform the learning experiences of students in calculus courses with an eye toward fostering diversity, equity, and inclusion (DEI).

In this showcase, we will devote a particular segment to instructors' uses of language in the math classroom, since we see language as being related to identity in two important ways. First, language is often bound up with culture and national origin, and speakers of a language (particularly a minoritized language in the US context) may feel an affinity with other speakers of that language. Second, all people, including speakers of the "same" language, use particular words, phrases, and ways of using words (discourse practices) to signal membership with particular communities (Gee, 1996). For example, at our university, the tutors in the mathematics learning center have developed ways of talking about their job (Bjorkman, 2019). The use of specialized terms (e.g., "good student" for the students who came to the tutoring center every day, or "member of the club" for non-tutors who hung out in the center) was one way for the tutors to form a discourse community and develop an identity within mathematics. We will close the chapter with recommendations for incorporating identity work to make students' calculus learning experiences more diversified, equitable, and inclusive.

We note as a preface to this chapter that the language and words that educators use to describe students' socio-cultural identities matters. For example, some research studies and policy documents use the word "Hispanic," while others use "Latino/a", "Latin@", "Latinx", or "Latin*" (Salinas, 2020). In the case studies in this volume, the authors' choice of terms to describe students might reflect regional differences, demographic labels used in policy documents, or an intentional socio-political stance related to identity. In our discussions, we use the terms that appear in each case study as we describe the case, but we acknowledge that terminology can be a contested space.

## 2. *Literature on Identity*

Intuitively, identity is intended to capture some essential aspect of who one is (Langer-Osuna & Esmonde, 2017). There has been extensive research on identity in mathematics education, most grounded in K-12 settings (see Langer-Osuna & Esmonde, 2017 for a synthesis) and more recently extending to the undergraduate mathematics setting (e.g., Larnell, 2016; Leyva, 2017). Yet, as Langer-Osuna and Esmonde (2017) note, literature authors have not come to an agreement on one single way to define the construct of identity. There are at least two different ways to consider identity: (a) ***individual*** (developed by an individual relative to a social setting) and (b) ***membership*** (based on membership in a socially salient group such as a gender, racial, or language group). When thinking of mathematics, we sometimes refer to one's mathematical identity, informally, an individual's self-concept in relation to mathematics. A person's mathematical identity is an example of an individual identity, whereas, for example, being "Black" (as a racial identity) is a membership identity. Research shows that students' mathematical identities are closely related to their membership identities (Barwell et al., 2017; Langer-Osuna & Esmonde, 2017; Martin et al., 2017). In this section, we will discuss different conceptualizations of mathematical identity and their implications, as well as unpack the relationship between students' membership identities and student's mathematical identities and experiences.

Researchers have drawn on four conceptualizations of mathematical identity: *psychoanalytic, narrative, positioning, and poststructural* (Langer-Osuna & Esmonde, 2017). A **psychoanalytic** view posits that math identities develop through coping with emotions, fears, and desires associated with mathematical experiences. A **narrative** approach views math identity as making sense of one's and others' mathematical experiences and reifying them as stories of success, failure, belonging, or distance. In terms of **positioning**, math identities develop through the statements made about oneself and others in mathematical settings. Finally, a **poststructural** perspective of identity sees math identities as being largely influenced by broader power establishments, such as education policy and cultural discourses. Although these conceptualizations of identity arose as ways for researchers to investigate people's mathematical identity formation, in this chapter we will suggest ways in which these conceptualizations can be leveraged by individuals and other classroom, departmental, or institutional agents of change to shape students' mathematical identities and learning experiences.

Each conceptualization has different implications for math identity formation (Langer-Osuna & Esmonde, 2017). A psychoanalytic perspective notes how crucial it is to consider the emotions (e.g., anxiety and fear) evoked in the learners in relation to mathematics or a mathematics-related context. A narrative approach emphasizes the benefit to inquire about the ways students describe experiences, whether students and others have experienced certain past events in or out of the classroom as positive or negative, and whether students have constructed these experiences as stories of belonging or distance. A positioning lens marks the importance of attending to how interactions in learning environments may position students as more or less mathematically competent, and how such interactions may make certain students feel included or excluded from the local or broader math community. Finally, a poststructural view highlights the need to consider broader relations of power, such as the regulation of higher education by policy makers, and how that may promote or hinder students' positive mathematical identity formation.

Synthesizing past research, there are certain classroom structures that are more likely to lead to positive STEM identity formation in each STEM field. For example, project-based learning environments that incorporate real-world problems with contexts related to students' lives may promote students' "rightful presence" in STEM (Tan et al., 2018), and active learning spaces that extend mathematical authority to students may help construct students' identities as mathematically capable (Cobb et al., 2009). Yet, thinking of how a "general student" develops STEM identity is not sufficient for determining an intervention for promoting a positive STEM identity, as research indicates that STEM identities and learning experiences are co-constructed along with other salient membership identities such as race, social status, and language background (Barwell et al., 2017; Langer-Osuna & Esmonde, 2017; Martin et al., 2017). This highlights the importance of considering students' membership identities in the context of mathematics education.

Racial identity and racialized experiences have been a major focus of attention in past mathematics education research (e.g., Larnell, 2016; Martin et al., 2017; Nasir & Shah, 2011; Stinson, 2008). Researchers have conceptualized race in multiple ways, including as a categorical variable (or identifier) in achievement gap studies and as a racial identity -- a socially negotiated self-concept associated with membership in a racial group, with economic, legal, and political ramifications (Martin et al., 2017). In the literature, one may find several racial overarching narratives that are consequential for students' educational journeys. These narratives include unfounded statements that African Americans are less capable of succeeding in math (Stinson, 2008) and that Asians are good at math (Nasir & Shah, 2011). Yet, although these discourses circulate widely and shape students' identity construction, the link between racialized overarching narratives and individuals' mathematical identities is not deterministic. Some students internalize the racialized overarching narratives about mathematics and there is a strong link between their racialized and mathematical identities. Yet, other students actively resist racialized narratives linked to their mathematical identities (Larnell, 2016; Martin et al., 2017).

There are longstanding and persistent inequities in STEM related to gender. In addition to documenting inequities in, for example, degree attainment and career persistence, authors in the classroom-focused literature in mathematics education have documented how women have been positioned as not mathematically capable or fit for doing mathematics due to their gender. For example, Ernest, Reinholz and Shah (2019) documented a group discussion during an inquiry-oriented undergraduate math class in which a female student encountered sexist remarks that insinuated women are for looking pretty and being in the kitchen rather than for doing mathematics. Scholars note that much of this research has drawn on binary perspectives of gender, generally leaving out conversations about gender queer identities (Langer-Osuna & Esmonde, 2017). As for research on LGBTQ+ in STEM, Mattheis et al. (2019) reports that queer students' conversations about gender and sexuality in STEM are frequently silenced due to assumptions of heteronormativity, and this leads to complications in forming a STEM identity.

The research on language in the context of mathematics teaching and learning has emphasized that language is both a cultural tool for expressing ideas and mathematical sensemaking and a sociopolitical tool for being recognized as a certain kind of someone in a mathematically-related context (Barwell et al., 2017). For example, Garza (2018) showcases a seventh-grade classroom with Latinx bilingual students in which a teacher used translanguaging practices, which integrated Spanish and English and their associated ways of saying and doing things (e.g., linguistic markers), to allow students to use their entire linguistic repertoire for developing mathematics biliteracy. On the other hand, LópezLeiva and Khisty (2014) demonstrate the sociopolitical role of language in a study situated in one elementary mathematics classroom. They report how a predominantly Spanish-speaking Latina student named Elsita was positioned as less capable mathematically through the teacher's use of language. For instance, during a group discussion, the teacher switched from English to Spanish particularly when she questioned Elsita's mathematical understanding in an annoyed tone. Furthermore, this interaction occurred in a classroom where a group peer had demanded "English only," which positioned Elsita's language as inferior. While there is a relatively large body of research on language and mathematics in K-12 contexts (de Araujo et al., 2018), the research on linguistic diversity in undergraduate mathematics classrooms is relatively sparse, so more research at this level is needed. Given this need, we devote a significant portion of Section 3 of this chapter to uses of language and their relationship to identity in the context of college mathematics.

Adding to the complex relationship between students' membership identities and students' mathematical learning experiences, teachers' membership and individual (e.g., mathematical teaching) identities can also play a critical role in shaping this relationship. For example, considerable research has documented the positive effect that matching students with a teacher of the same race or ethnicity as them can have on students' mathematical achievement (e.g., Egalite et al., 2015; Redding, 2019). In particular, Egalite et al. (2015) reported that low-performing Black and white students in grades 3 through 10 appear to benefit academically from being assigned a teacher that matches their racial group. The authors attribute this result to multiple factors, including potentially greater shared cultural understanding that may aid teachers in becoming more effective role models or advocates as well as selecting more culturally relevant tools for teaching, learning, and mentoring students in the same racial group as the teachers. Another related factor may be the teachers' attitudes towards the students. This claim was supported by Redding (2019)'s finding, based on a comprehensive review of the literature, that teachers' perceptions of the mathematical achievement of Black (respectively, Latinx) students tended to be higher when the teachers were also Black (respectively, Latinx). But matching membership identities between students and teachers is not the only way teachers' identities can influence students' mathematical experiences. For example, as noted in Soto et al.'s study (CS 22 ), Estrada et al. (2018) indicate that increasing teachers' self-efficacy in the social aspect of the mathematics classroom (e.g., teaching practices that promote social inclusion) can lead to increased student engagement and connection to the classroom community, and ultimately promote more positive student mathematical identity formation. This result points to the importance of building positive teacher identities to shape students' mathematical identities and learning experiences.

## 3. *Identity as Illustrated in Case Studies*

Nearly all of the case studies in this volume raise the issue of identity in relation to the teaching and learning of calculus, though identity may not be the primary focus of each study or the intervention described. Given this volume's focus on diversity, equity, and inclusion, students' various identities were discussed and made salient for locating and addressing inequities in undergraduate mathematics education. For example, as illustrated in Table 1, there are case studies that describe analyses or interventions focused on students' racial or ethnic identity, gender identity, linguistic identity, identity as a first-generation college student, identity due to socioeconomic status, and even students' mathematical (or STEM) identity.

In this section, we characterize the ways in which case study authors drew on the concept of identity to promote DEI in their classrooms and institutions. First, we describe how the case studies framed students' membership identities as well as the benefits and limitations of these framings for promoting DEI. Next, we illustrate how the case studies leveraged different conceptualizations of identity to promote a positive mathematical identity. We end with a subsection focused on the role of language in identity formation as illustrated in the case studies.

*Table 1*
*Variety of Students' Identities Made Salient in this Volume's Case Studies*

| Case Studies | Racial or ethnic | Gender | Linguistic | First-generation college student | Socio-economic status | Math (or STEM) |
|---|---|---|---|---|---|---|
| Burn et al. (CS 4 ) | X | X | | | | X |
| Chang & Chen (CS 7 ) | X | X | | X | | X |
| Davis (CS 8 ) | X | | | | | X |
| Deshler et al. (CS 9 ) | X | | | X | | X |
| DiGregorio et al. (CS 10 ) | X | | | X | X | X |
| Liou-Mark et al. (CS 16 ) | X | X | | X | X | X |
| MacArthur (CS 17 ) | X | X | | | | X |
| McNeill et al. (CS 19 ) | X | X | | | | X |
| Soto et al. (CS 22 ) | X | | | | X | X |
| Stacy (CS 23 ) | X | X | | X | | X |
| Stone-Johnstone et al. (CS 25 ) | X | X | | | | X |
| Villa et al. (CS 27 ) | X | | X | X | | X |
| Villalobos et al. (CS 28 ) | X | | X | | | X |
| Wagnon & Hubbard (CS 29 ) | | | | X | | X |

### 3.1. *Uses of Students' Membership Identities*

The case studies highlighted students' membership identities (e.g., racial, gendered, and linguistic identities) with different depths and to different effects. Similar to conceptualizations of race (Martin et al. 2017), two main ways the case studies framed students' membership identities were: (a) as a **categorical variable or label** used to identify students from different social groups, or (b) as a **self-concept or lived experience** with shared cultural tools and legal, economic, and political negotiations. For example, in their case study, Chang and Chen (CS 7 ) analyzed whether a math course redesign with additional features led to more equitable outcomes for students from certain groups. To assess the effectiveness of this intervention, the authors conducted a quantitative study comparing the performance of different gender groups (female versus male) and different ethnicity groups (Asian, Blacks, Hispanic, and white), leveraging the students' demographic information as categorical variables. An example of a case study that framed the students' linguistic identities as a self-concept or lived experience was the study by Villalobos et al. (CS 28 ). The intervention in this study acknowledged the linguistic assets and sociopolitical experiences of bilingual students from a US-Mexico border region, where Spanish-English bilingualism is more the norm rather than the exception and where Spanish may be seen as inferior to English. To send the message that being bilingual is a valuable skill, the math department at this study's institution offered several undergraduate math courses, including Calculus, in Spanish along with a dual language certificate. Additionally, building on students' bilingualism, these courses allowed students to speak both English and Spanish.

Each case study's framing of students' membership identities generally corresponded to different views of equity, each with its own benefits and limitations. As referenced in MacArthur's case study (CS 17 ), Gutiérrez (2018) distinguished two ways to think of equity: (a) as a **destination** (a static outcome after an event), or (b) as a **process** (a dynamic phenomenon during an event). The case studies that used the categorical variable or label view of identity generally framed equity as a destination -- usually as a small achievement gap, as in Cheng and Chan's case study (CS 7 ). On the other hand, the studies that used the self-concept or lived experience view of membership identity tended to frame equity as a process. For example, Villalobos et al.'s case study (CS 28 ) considered equity in terms of the learning dynamics of the classroom in relation to the students' bilingual identities. They note that in one group session, one student announced his mathematical discovery to the class in English and then explained it to his group in Spanish, spreading his excitement and eureka moment to his group peers. Not all case studies' conceptualizations of equity fell into only one of the two views of equity. For example, in Stone-Johnstone et al.'s study (CS 25 ), the teacher considered the (destination) patterns of participation as a tool to reflect on the process of teaching (e.g., calling or not on certain students to participate in class) relative to teacher's biases about students with particular membership identities. Such conceptualizations of equity incorporated both (even if partially) a process view and a destination view of equity.

As can be seen in the examples above, two benefits of the destination view of equity are that it provides a simple and easy way to assess and compare equity outcomes among different membership identities, and that it assesses *whether* a certain event achieved equitable outcomes for groups of students with different membership identities. However, one major limitation is that the destination view fails to explain the process by which students with certain membership identities achieve certain outcomes. This can be problematic because reporting static outcomes (e.g., achievement gaps) may reinforce existing deficit narratives. On the other hand, two benefits of the process view of equity are that it highlights *how* an event achieved (or not) equity for different social identities, and emphasizes the ongoing effort needed to sustain equity across time and places (Gutiérrez, 2018). Yet, one limitation is that the process view is more complex and difficult to implement for assessing equity among students from a large set of membership identities.

### 3.2. *Uses of Mathematical Identity Conceptualizations*

Several of the case studies in this volume leveraged the construct of mathematical identity to inform the design or analysis of interventions. The ways these case studies approached mathematical identity spanned all four conceptualizations of identity (psychoanalytic, narrative, positioning, and poststructural) discussed earlier, in Section 2 of this chapter.

Some case studies leveraged multiple conceptualizations of mathematical identity, allowing them to reap the benefits afforded by each conceptualization. For example, we see evidence that MacArthur (CS 17 ) drew on the psychoanalytic, narrative, and positioning views of identity to design more humane exam structures. More specifically, MacArthur split the exams into two portions, one to be completed by each student with support from a group of peers and another to be completed by each student alone. One factor that motivated MacArthur to move from pure solo exams to exams with a group portion was her aim to reduce test anxiety. One could argue then that her inspiration was afforded by her psychoanalytic conception of identity, which considers that mathematical identities develop as students cope with unconscious forces such as anxiety. MacArthur also leveraged this view of identity to assess the effectiveness of her intervention, noting that some students reported feeling less anxious during the exam as a result of including a group portion. By asking students to provide written survey responses and commentary, MacArthur also drew on the narrative perspective of identity. Through these students' narratives, MacArthur found that students reported stronger feelings of belonging with the new testing structure compared with traditional exams, rendering the use of a group exam structure worth keeping. Hence, MacArthur leveraged one of the benefits of the narrative approach of identity, which highlights *whether* a certain classroom structure like a group exam is experienced by students holistically in a positive way. Finally, one could see how MacArthur may also have made use of a positioning view of identity. More specifically, in the analysis of students' comments, MacArthur also attended to students' descriptions of negative interactions. With this focus in mind, MacArthur noticed that two white women mentioned being talked over in their group (interactional experience) as the negative side effect of group exams. So, MacArthur decided to develop group norms and pay closer attention to how students are grouped in order to better serve women. Hence, one could see how MacArthur may have reaped one benefit from the positioning view of identity, which has the power to highlight *how* interactions during a particular activity or classroom structure (such as a group exam) shape students' experiences.

For additional examples of how the case studies leveraged the construct of mathematical identity, consider Davis' (CS 8 ) and Burn et al.'s (CS 4 ) studies. Davis (CS 8 ) drew on the narrative view of mathematical identity by attending to students' comments about their mathematical experiences in relation to an in-class peer tutoring program. In their stories, students mentioned that having peer tutors showed that "people (ladies in particular for this class) can make it through STEM degree plans" (Davis (CS 8 ), p. 5), giving a positive indication of effectiveness of the intervention in addressing DEI for women. On the other hand, Burn et al.'s case study (CS 4 ) applied a poststructural view of mathematical identity. Prior to suggesting ways to support the engagement of African American, Latinx, Native American, and Southeast Asian students through practicing welcomeness to engage, Burn et al. (CS 4 ) reflected on what may lead certain students not to ask for help or engage in mathematics learning. To do this, the case study authors attended to broader power structures (e.g., racism and overarching narratives related to race, gender, or language). For instance, they noted that masculine identity, model minority stereotype threat, and students' prior experiences with racism can prevent students from seeking mathematical help. Drawing on a poststructural view of identity by way of identifying broader structures of power, may have productively directed the authors towards considering or analyzing instructional ways to make students' engagement feel welcomed and desired.

One additional facet of identity related to mathematics is the connection between teachers' identities and students' mathematical identity. The case study by Soto et al. (CS 22 ) shows one way to impact students' mathematical identity through implementing a professional (PD) workshop. To explain the effectiveness of the workshop, the case study authors built on Adiredja and Andrew-Larson's (2017) framework, which argued for targeting students' knowledge, identity, and power to promote equity. The workshop targeted those same three components, but on teachers instead of directly on students. In terms of teachers' knowledge, the workshop engaged the PD participants in reading and reflecting on articles about equity, high-cognitive demand tasks, sense of community, stereotype threat, and student-centered learning. For example, the instructors were supported in learning about kindness clues of social inclusion to promote equity in the mathematics classroom (Estrada, Eroy-Reveles, & Matsui, 2018). The case study authors reported that the fellows came to feel more knowledgeable about implementing rich tasks and about the positive impact of this teaching approach on students' mathematical identity. To target teachers' identities, the workshop promoted a sense of community by creating a space where fellows felt safe and valued while sharing their teaching struggles and efforts. The case study authors shared how one fellow felt grateful for the sense of community the workshop gave her, in contrast to the limited camaraderie she experienced prior to joining the PD due to her position as a part-time faculty. Another fellow also reported feeling a strong sense of belonging and higher self-efficacy with the community of fellows. Finally, the workshop targeted the teachers' power. After the PD, the teachers reported feeling more empowered to make decisions in their own teaching, including advocating for student-centered teaching even when formal power structures discouraged it (e.g., when administrators pressured faculty to teach too high an amount of material). This case study successfully highlights that "increasing DEI for our students requires addressing DEI for our faculty" (Soto et al. (CS 22 )).

### 3.3. *Uses of Language in Relation to Students' Identities*

One aspect of identity that played an important role in six of the case studies in this volume is language, in the two ways described in the introduction: (a) language as **bound up with culture and national origin**, and (b) language as **discourse practices** that people use to signal membership within socio-cultural communities. Two of the case studies address language in the first sense. In four additional studies, the use of language is related more to the second sense.

Villalobos et al. (CS 28 ) and Villa et al. (CS 27 ) address language in the first sense. In the case by Villalobos et al. (CS 28 ), the members of the mathematics department designed mathematics classes in Spanish and implemented them bilingually (Spanish-English) for undergraduates at their university located along the US-Mexico border. These faculty members were using language in the first sense of creating cohesion among speakers of the same named language and using that as a resource for teaching mathematics. This program was designed to build on the prominent culture of bilingualism in this region. Through this intervention, the authors highlighted the potential for creating a linguistically heterogeneous mathematics department, and they showed how this stance opened new opportunities for bilingual students to develop identities in the context of undergraduate mathematics classes. In a related case study, Villa et al. (CS 27 ) found that Latinx students and teaching assistants in precalculus workshop classes fluidly engaged in translanguaging, even when there was not an official bilingual program in place. The innovative program showcased by Villalobos et al. (CS 28 ) and the observation from Villa et al. (CS 27 ) provide models for faculty in other universities, particularly Hispanic-Serving Institutions (HSIs), to consider how the mathematics department can be aligned with the university's HSI identity.

In the case studies by DiGregorio et al. (CS 10 ) and Deshler et al. (CS 9 ), we see another use of language in relation to the construction of mathematical identities in an undergraduate mathematics program. In this case, the use of language is related more to the concept of language as discourse practices. Specifically, in DiGregorio et al. (CS 10 ), the designers of the mathematics program's placement process shifted the name of their mathematics placement test to the mathematics placement *tool*. This seemingly minor change in terminology was intended to signal that the placement process was not a "test" in the sense of a summative judgment, but rather, a process with the goal of matching students with their classes. Similarly, in Deshler et al. (CS 9 ), the faculty who developed an Emerging Scholars Program (ESP) -- a program to support Black, Hispanic, Native American, or Hawaiian/Pacific Islander students -- called the invitation to join the program a "nomination." These shifts in terminology are related to the discourse of school mathematics. In school mathematics discourse, "tests" invoke notions such as measurement, precision, and judgment. We know from research on stereotype threat that "tests" can provoke anxiety and lead to diminished performance by students who perceive a threat to their membership identity group (Steele, 1997; Steele & Arouson, 1995). Similarly, for the ESP, the authors of Deshler et al. (CS 9 ) describe their intention of using the word nomination to signal that the invitation to join the program was a special opportunity, rather than a penalty. In both cases, with the redesign of the placement process and renaming the selection process into ESP, DiGregorio et al. (CS 10 ) and Deshler et al. (CS 9 ) sought to transform students' perception of selection and placement processes and also, perhaps, to humanize those processes.

Finally, two additional case studies highlight how instructors' use of language in calculus classrooms might lead students to feel included or excluded from the local mathematics classroom community, or even the mathematics community at large. In particular, the case study of the "course drop event" by McNeill et al. (CS 19 ) highlights how a teacher's use of language might provoke racialized feelings of exclusion among students. These feelings may be especially deleterious for students from non-dominant racial ethnic or gender groups. McNeill et al.'s case study (CS 19 ) focuses on language in the sense of "discourse." Within college mathematics, it is common for instructors to mark certain types of problems as "easy" or "review." This type of talk is an example of a Discourse practice in the sense that Gee (1996) used the word Discourse (with a capital D). For people who are part of the discourse community, these words, phrases and ways of using language signal affiliation and membership within the community. However, for students who are outside of the core community, being present in the use of this discourse may be a form of symbolic violence. Students may feel like they are being excluded or reminded of their marginal status within the community. This feeling is compounded again when students identify with a historically excluded racial, ethnic, language, or gender group. In the case study by McNeill and colleagues (CS 19 ), we see how students hear the implicit messages of these classroom discourse events and interpret these messages through their membership identities, including racial and gender identities. This case study highlights for us (as faculty in STEM departments) the importance of instructors attending to how they use language and languages in their lectures and how their use of language may serve to exclude students.

On the other hand, Stacy (CS 23 ) noted how a mathematics instructor (herself) used language to undermine the “not good at math” identity. To achieve this, Stacy highlighted for students that math is a creative effort rather than a result of innate ability. In doing this, Stacy engaged students in conversations around the nature of the “math person” identity. By framing math success as a creative process, Stacy highlighted that mathematical identity is not biologically predetermined, but rather something that could be built through hard work. This discourse might invite more students to feel included in the mathematical community, whereas a discourse around innate mathematical ability (e.g., “this is easy”) might make some students feel excluded.

These six cases highlight or suggest interventions at various levels (e.g., individual, classroom, departmental, and institutional) regarding the use of language that is inclusive of all students. As an individual, a teacher could reflect on the framing of particular events by considering students’ racialized and gendered lenses aligning with their lived experiences. As an agent in the classroom, a teacher could frame math success as a creative process rather than innate ability. Another way a teacher could leverage language in the classroom is to encourage students to use their full linguistic repertoire (e.g., through translanguaging) to learn mathematics. In addition, mathematics departments could build on students’ linguistic assets, for example by offering courses in Spanish, and the departments could coordinate with their university to offer a dual-language certificate as one way to empower students (and their families) whose linguistic assets may be currently undervalued in political, economic, and sociopolitical ways. Additional ways departments and institutions could encourage discourses of inclusion are by carefully naming and framing certain structures and processes, such as the nomination of participants for particular programs and the mathematics placement tool.

## 4. *Recommendations*

As illustrated in the case studies, there are several ways to leverage the construct of identity for promoting DEI in college mathematics classrooms. In this section, we provide recommendations for addressing DEI using identity and language. In alignment with the cross-cutting inquiries of this volume, we speak to the reader at multiple levels (as an individual and as an agent of change in a classroom, department, and institution). With this diverse audience in mind, we elaborate on three aspects of identity that can be leveraged to promote DEI in college mathematics: (a) students’ membership identities, (b) conceptualizations of mathematical identity, and (c) uses of language.

One recommendation is to **leverage your students’ membership identities**. As an individual, we encourage you to reflect on your conceptualization of students’ membership identities. We invite you to acknowledge that students are more than identifiers by considering the economic, legal, and political implications and cultural assets associated with these identities. For example, consider ways in which events may be experienced as microaggressions by students from certain racial, gendered, or linguistic groups, and also learn about your students' assets. We also invite you to reflect on your conceptualization of equity. It may be useful to think of equity, not only as a destination, but also as a process (Gutiérrez, 2018). This may help you expand your assessment of equity beyond achievement gap studies that, without process identification, may inadvertently reify deficit narratives about students with certain membership identities. As an agent of a classroom, consider using the cultural assets associated with students' membership identities as tools for teaching and learning. For instance, consider encouraging students to use their full linguistic repertoire through allowing them to use translanguaging practices. As an agent in a department or institution, consider using students' assets for curriculum design (e.g., by designing classes in Spanish or identifying applications of interest to diverse students) and for student empowerment (e.g., by offering a dual-language certificate that incorporates STEM as well as language classes).

Our second suggestion is to **leverage your conceptualization of mathematical identity**. As an individual, we encourage you to reflect on your conceptualization of identity and associated benefits and limitations. It may help you to consider the four conceptualizations of identity (psychoanalytic, narrative, positioning, or poststructural) that were provided in the literature review above and illustrated through the case studies. As an agent of a classroom, department, or institution, we invite you to leverage different conceptualizations to design or analyze interventions for promoting a positive identity in your students in relation to mathematics. For example, we elaborated on how one mathematics instructor (MacArthur (CS 17 )) may have leveraged the psychoanalytic, narrative, and positioning conceptualizations of identity to create a more humane exam structure (in this case, a group exam portion) and assess its effectiveness as an intervention. Similarly, you could select a classroom structure you might want to improve. You could then ask yourself questions or take approaches that leverage the benefits of different conceptualizations of mathematical identity. For instance, in alignment with the psychoanalytic view of identity, you may ask yourself: What fears and desires may be at play in this classroom structure and how may I mitigate them? After implementing the change in structure, you could leverage a narrative perspective of identity by requesting student comments about their mathematical experiences and see whether they develop stories of belonging or distance.

Our final recommendation is to **leverage your use of language to promote a positive identity** in students in relation to mathematics. As an individual, we encourage you to reflect on your use of language by considering the two ways illustrated in the chapter. One way is associated with choosing a language that is reflective of students' culture and national origin, and the second way is related to choosing particular words, phrases, and discourse practices that signal membership with particular communities. In thinking about your language, a critical self-examination might reveal your hidden assumptions about language and your uses of exclusionary discourse practices in your class (e.g., the "course drop" event in McNeill et al.'s case study ((CS 19 )). We then invite you to choose to use language that is inclusive of all students. The end of the subsection 3.3 provides examples of how you can achieve this at the individual, classroom, departmental, and institutional levels.

## 5. *Concluding Remarks*

Connecting to Larnell's Manifesto in this volume (BUILDING AN EVOLVING FRAMEWORK: A CLARION CALL / MANIFESTO) and the leading inquiries (CROSS-CUTTING INQUIRIES: MOVING TOWARDS A SYSTEMS-CRITICAL PRACTICE) that guide this book, change must include more than individual faculty members raising their personal awareness about the roles of identity and language in calculus learning (although such individual changes are needed). Rather, change will require collective action from broad and diverse groups of stakeholders in each college or university, working at multiple levels. We encourage readers of this chapter to form coalitions to continue unpacking questions related to DEI and identity that were highlighted in this chapter. Additionally, we invite you to work with peers and to critically analyze local practices at your institution around calculus. Ask, how do current practices support equity or lead to inequitable outcomes? After engaging in reflective processes, take action. With collective attention to how calculus course structures position students, groups of faculty members, department leaders, and colleagues from across your campus can create learning environments in which more students of varied identities feel welcomed and succeed in undergraduate mathematics.

We invite the reader to reflect on the following questions:

- What actions might you take to encourage others to value and build on their students' membership identities and to do so yourself with your students' identities?
- How might you develop interventions to promote positive mathematical identities in your courses and assess their effectiveness?
- How might you be more inclusive of all students with your language uses?
- How might you take into consideration individual variations among students classified with similar demographic information?

# 7 Reflecting on the Role of Assessment in Placement and Evaluating Student Learning in Calculus

**Valerie Cheathon**

*I Am STEM ABLE*

Valerie Cheathon received her master’s in Applied Mathematics from WSU. Her university mathematics education inspired her to take time off to understand faculty perspectives around STEM education and assessment while observing their pedagogical strategies to refine her own methods of teaching mathematics in a quarter and semester system. This year she initiated a research opportunity with I Am STEM ABLE that she hopes will transition her back to graduate school to acquire a Ph.D. in Mathematics Education.

**Nicholas Fortune**

*Western Kentucky University*

Nicholas Fortune is an Assistant Professor of Mathematics Education at Western Kentucky University. He received his Ph.D. in Mathematics Education from North Carolina State University and his masters and bachelors in Applied Mathematics from Rensselaer Polytechnic Institute. His research centers around instructional change in undergraduate mathematics and how mathematics faculty can collaborate on pedagogy and student thinking to support their instructional change and incorporate inclusive pedagogical practices.

**William Hall**

*Washington State University*

William Hall is an Assistant Professor of Mathematics Education at Washington State University. His research interest is the role of identity in the teaching and learning of mathematics. Will earned undergraduate and MST degrees from UMaine before a Ph.D. in mathematics education from North Carolina State University. He is a proud first-generation college student and different learner. He tries to aim his teaching, research, and service goals/responsibilities towards making mathematics a better place.

## 1. *Introduction*

All math assessments are of human design. Every person, assessment creator and taker alike, is embedded within various social structures and subject to various social norms and customs. Each quiz, exam, project, portfolio, essay assignment, etc., while serving to measure students’ understanding, is inherently a cultural activity. Assumptions regarding how much mathematics one should be able to complete, to what degree of “correctness” that mathematics must be demonstrated, and under what time and resource constraints that knowledge is demonstrated, are not solely the province of mathematics, but very much of culture.

People experience assessment differently. Completing a high-stakes exam, for example, is an in situ cognitive construction for each test-taker, subject to what they are experiencing at the very moment of completing the assessment in addition to their mathematical knowledge. Results of such an exam, like all measures, serve as proxies for a person’s knowledge. Those of us who take on roles in undergraduate mathematics education, have the responsibility and duty to continuously reflect on the following question: **Do our assessment practices and policies fairly and accurately measure what we claim they measure?** Reflecting on this question and actively investigating local practices is a necessary step in beginning and/or sustaining an equitable and effective calculus program. This responsibility stems in large part from the fact that measurement is not inherently benign.

In fact, the act of measuring a person's mathematical knowledge leaves a mark. Larnell (2016) illustrates how non-credit bearing remedial mathematics courses and an online placement exam serve as racialized experiences in undergraduate mathematics. One person, Vanessa, shared that she did not try to answer many of the questions on her placement exam, "I just pressed next like on all of them because I knew that I wanted to be in the lowest math. I knew that I wasn't that good in math...I didn't want to come into college, um, like behind" (Larnell, 2016, p. 254). Clearly, the assessment and placement process was not benign in Vanessa's college mathematics experience. Had she been subject to a different placement system, one not dependent on the results of a single exam alone, she would have been placed differently. Vanessa had taken and done well in precalculus and AP Statistics/Probability, graduating near the top of her class (Larnell, 2016).

Unfortunately, exactly these kinds of practices have been celebrated for their purported objectivity. Au (2013) writes, "under the assumption that standardised tests provide fair and objective measurement of individuals, such testing seemingly held the promise that every test taker is offered a fair and equal shot at educational, social, and economic achievement. Problems like racism and class privilege are thus supposedly ameliorated through testing" (p. 13). This view relies upon a meritocratic perspective on mathematics assessment and achievement. Au claims that with such a perspective, "the low achievement on standardised tests of working-class people, non-white populations, and some immigrant groups can then be simply and neatly attributed to the failure of individual students, individual groups, or individual cultures, and not attributed to existing structural inequalities" (p. 13). Improperly assigned objectivity to a given assessment thus becomes the tool of the oppressor.

Extremists use purportedly objective assessments to justify causing harm to people. Results of IQ tests were used alongside other measures to justify forcibly sterilizing more than 60,000 people in the U.S. during the 20th century (National Public Radio [NPR], 2016). In 1927, the U.S. Supreme Court upheld the right of a state to forcibly sterilize those deemed unfit. Justice Oliver Wendell Holmes, Jr. wrote the summary for the 8-1 decision: "three generations of imbeciles are enough." Carrie Buck was forcibly sterilized. Her infant daughter Vivian was classified as feeble-minded by eugenics researchers for "not looking at a coin held in front of her face" (NPR, 2019). Author Paul Lombardo writes, "although the Buck decision set the stage for more than sixty thousand operations in the United States and was cited at the Nuremberg trials in defense of Nazi sterilization experiments, it has never been overturned" (Lombardo, 2009).

Our intentions are not to draw parallels between the suffering endured during these atrocities and calculus assessments today. But we also do not want to marginalize the impact of using allegedly objective assessments for placement into and progress through calculus. Assessment practices have evolved in the past and we should continue to reflect on what assessments we use and how we use them, particularly in calculus - one of the narrowest gateways in undergraduate studies. University admissions are doing this work. In reflecting on the systemic bias that has been demonstrated within standardized assessments (e.g., SAT/ACT/GRE), many universities are no longer requiring them (Jaschik, 2020). In fact, more than 1,400 U.S. four-year colleges have dropped their requirements for these exams (National Center for Fair and Open Testing, 2021). Only through such reflection will we be able to provide more equitable mathematics pathways.

In this chapter, we reflect on the work of faculty and change agents who have re-imagined various aspects of either how students' mathematical understanding is assessed and/or the structure and policies of assessment in university calculus. We begin with a reflection on the role of assessment in placement into calculus programs followed by ways mathematics instructors are reimagining assessment practices in evaluating student learning of calculus. Throughout the chapter, we reflect on the role of assessment within the context of this volume's leading inquiries: individual level, classroom level, department level, university level, and policy level. As authors of this chapter, we acknowledge that we are products of and have played our roles in perpetuating many normative calculus assessment and placement techniques. Therefore, we position ourselves with the reader, as learners who want to make the practice, assessment, and structures of calculus more equitable, fair, and just.

## 2. ***The Role of Assessment in Placement***

As members of the mathematics community, we understand the importance of students entering the classroom equipped with mastery of prior mathematical concepts, providing some assurances towards the students' success in their current mathematics course. A variety of math placement philosophies describe placement as a tool that will enable students to have a positive and productive university mathematical journey. There are several models that mathematics departments use to place students into calculus. These placement models vary in structure and requirements for students with one approach being a single exam built by a department or a third-party publisher software (e.g., ALEKS) while another being a multi-measure approach that considers these exams, ACT/SAT scores, high school or previous GPA, and/or previous mathematics coursework. In this section, we discuss what some universities from the case studies are doing to build towards more equitable placement practices, such as examples in multi-measure placement. We also briefly touch on a corequisite model that some of these universities are using to support students after they have been placed into a mathematics course. We conclude the section with guiding and reflective questions.

### 2.1. *Multi-Measure Placement Practices*

> *It is time to support all students who are excited about STEM careers, especially those who have been structurally disadvantaged, and ensure that their mathematics placement does not discourage these interests, but rather supports them in thriving in their STEM pursuits (DiGregorio et al., CS 10).*

Assessments designed for placement in undergraduate mathematics should not be one dimensional (Hsu & Bressoud, 2015; Scott-Clayton, 2012). Placement based on a single timed assessment gives an incomplete view of students' understanding. Further, inaccurate and/or unfair placement appends time and debt to students' path to degree. Published in *PRIMUS*, DiGregorio and Hagman (2020) found that students from historically-marginalized groups in mathematics (e.g., first-generation college students, low-income, and Students of Color) have high anxiety associated with mathematics placement exams, perceived the placement exam as very high-stakes, and associated the exams with fixed mindsets rather than growth mindsets.

These populations deserve our ingenuity and effort to craft a plan that positions them for success, rather than sidelining them or weeding them out. DiGregorio and Hagman (2020) proposed an asset framed mathematics placement process as a starting point for reform. Suggestions include refraining from terminology such as “exam or test” and calling it a “process,” messaging that the department supports students through the placement process, refraining from using only one score (e.g., ACT/SAT) for placement, incorporating motivation and self-efficacy into placement, and providing online preparation materials that can be engaged in a collaborative manner. DiGregorio and colleagues (CS 10 ) described their implementations of their recommendations in this volume.

Another approach included reviewing the data analytics from first-year students’ performance on standardized tests such as the SAT or ACT and directing them towards self-guided instruction tools like ALEKS (Assessment and Learning in Knowledge Spaces) to prepare them for the upcoming semester (Benken & Slowinski (CS 2 )). Similarly, Mawhinney et al. (CS 18 ), created a multi-measure approach to placement. While raising the score required to place into Calculus I, students were provided an “alternate entry” that prepared them for Calculus I without requiring an entire semester of preparation. If students still found it difficult to pass the Calculus I entry exam, entrance could still be earned by completing a majority of a self-paced online course. “This alternate entry allows students a more timely path to Calculus I, along with tools and background knowledge needed to be calculus ready” (Mawhinney et al., (CS 18 )). This is an example of how university mathematics departments are rethinking both placement and prerequisites. This is referred to as a corequisite model, one in which the qualification bar is raised alongside significant real-time support for any student who needs it.

### 2.2. *Reimagining Placement and Prerequisites*

Reimagining placement is much more than a suggestion or innovative idea. It is a responsibility. Educational structures, particularly those that serve as barriers to entry for certain areas of society like assessments in mathematics do, are built upon and reinforce social structures. Recall Vanessa. Larnell (2016) describes her and other Black students’ experiences with incorrect placement and racialized identity threat in non-credit-bearing remedial courses. Many of these structures remain undisturbed from days when many federal, state, and educational policies were anti-Black and anti-woman - some are calcified from when slavery was legal but teaching a slave to read was not. As highlighted by one of the case studies:

> *Despite that equity and opportunity gaps in mathematics are framed around race/ethnicity, efforts to close gaps tend to take a race-neutral approach that “lifts all ships” but does little to advance racial equity in mathematics. Colorblind approaches reinforce whiteness as normative in mathematics, creating barriers to racialized student engagement, learning, and attainment in mathematics (Battey, 2013; Battey & Leyva, 2016; Davis & Martin, 2008; Martin, 2009; Gutiérrez, 2012, 2017; Palmer & Wood, 2013; Wood & Palmer, 2014; as cited in Burn et al., (CS 4 )).*

Ensuring that the needs of the systematically oppressed and vulnerable students are met naturally envelopes the needs of the entire population.

Corequisite models serve as a current alternative to existing semester-long prerequisites. Corequisite models are those that provide just-in-time support to students in foundational mathematics courses. They often include one additional credit-hour where the student explores missing prerequisite mathematics content. They can also support in developing students' mathematical confidence, providing a safe place where being wrong is a welcomed part of learning and offering review sessions that show diversified perspectives of mathematics concepts, thus accommodating a variety of learner types. Importantly, research has shown corequisite models work (Logue et al., 2019). In a randomized controlled trial, Logue et al. (2019) found that the "corequisite group not only demonstrated significantly higher quantitative course pass rates but also success in many other disciplines, as well as significantly higher graduation rates" (p. 294). Corequisite models are a refreshing replacement to self-guided learning tools or semester-long courses that set a student back without more personal analysis of their needs (Canner et al., (CS 6 )). Mawhinney et al. (CS 18 ) describe offering such a corequisite model as an option to students, emphasizing to them the technical and emotional benefits of taking part in such an offering:

> *Not only do students engage in tasks that increase their understanding of concepts, but they are actively doing mathematics in a safe space with classmates and with faculty who believe that they can be successful in mathematics (Mawhinney et al., (CS 18 )).*

The inclusivity and connection that a corequisite model provides can enhance the mathematical learning environment that mathematics departments strive to provide. They can also help redefine the culture of mathematics and mathematics teaching at your university. "We are working to create a community of practice, where we share information on the content and pedagogy, and search for ways to provide better support for all students in a first-year mathematics course" (Jensen-Vallin, (CS 13 )).

We can look beyond traditional means of mathematics course placement and make room for a holistic means of mathematics development for our students. Solid footing on a mathematical pathway is important for every student. Once their path has been securely and adequately set, our responsibility transitions from how we place students to how we evaluate their learning.

### 2.3. *Guiding and Reflective Question on the Role of Assessment in Placement*

- How does your department currently place students into calculus courses?
- What role do exams play in placing students into calculus? Are there other ways for students to place into calculus besides an exam?
- In what ways can your department evaluate and/or alter the role of assessment in your placement policies?
- How do/could placement efforts relate to diversity, equity, and inclusion efforts in your department?

## 3. *The Role of Assessment in Evaluating Student Learning*

Assessment in mathematics, at the individual and classroom level, is often in the form of exams and quizzes that are time and resource limited (Burn & Mesa, 2015; Hillel, 2001). As part of the Characteristics of Successful Programs in College Calculus study (Bressoud et al., 2015) conducted by the Mathematical Association of America (MAA), Burn and Mesa (2015) found that only 25% of faculty reported to have projects as part of their calculus courses. 36% of faculty give three exams and 47% of faculty give short quizzes in some class sessions (Burn & Mesa, 2015). This reflects what Hillel (2001) described 14 years ago. These assessment tools are used to evaluate student learning, often with confidence. "Mathematics faculty tend to have a high degree of confidence in exams and quizzes as indicators of student learning in part because assessing whether students have mastered concepts and principles is presumed to be more easily quantifiable in mathematics" (Lattuca & Stark, 2009, as cited in Burn & Mesa, 2015, p. 51). However, decades of education research have argued that to evaluate student learning, we cannot do it in a vacuum (Feldman, 2018). Our assessment practices must align with our pedagogical practice, consider the language that we use, be aware of our students' diverse backgrounds, and more.

The case study authors in this volume have done work to address these recommendations at the individual, classroom, and department level. Namely, they have worked to modify timed exam policies (through exam revisions and group exams), provide alternatives to timed exams (portfolios, projects, and oral exams), and reform how they transparently calculate grades. These modifications often focus on inherent structural inequities in how we assess student learning. We conclude the section with guiding and reflective questions.

### 3.1. *Reimagining Timed Exam Policies: Exam Revisions and Group Exams*

According to research from K-12 education, one assessment technique that offers the potential for more equitable grading is make-up work, often in the form of exam revisions (Feldman, 2018). Another assessment technique, recently more commonly seen in undergraduate mathematics, is group exams (Kinnear, 2020). Two case studies in this volume addressed implementing exam revision policies and two discussed how they utilized group exams in calculus.

Byrne et al. (CS 5 ) implemented a yearlong Calculus I course designed specifically for students who would have been placed into remedial courses. This was to move away from a system that placed students in a course that did not contribute to the credits needed to graduate, but did cost tuition. The population in such courses was disproportionately underrepresented and racially minoritized and students from low-income households. The transition to a yearlong Calculus I course was implemented for these students to "provide better access and pathways to the entire group of poorly-served students" (Byrne et al., (CS 5 )). For the first semester of the yearlong course, students would receive a grade of Credit or No Credit. There were two semester exams that both had opportunities, and importantly expectations, to submit revisions. The second semester, however, was more traditional with "more emphasis on graded work that aims to prepare students for expectations in Calculus II" (Byrne et al., (CS 5 )).

Similarly, Stacy (CS 23 ) also created a yearlong stretch course in an effort to address low pass rates in Differential Calculus after they noticed that successfully passing Precalculus was not indicating success in Differential Calculus. Stacy (CS 23 ) provided specific details on their exam revision policy.

> *Students are graded qualitatively, with each problem earning A, B, C, D, F, or 0. After the exam, students have a chance to correct their work to raise their grade up to one letter grade per problem. To have their grade raised, they must correctly redo the problem, and explain any errors they made originally. For example, if a student earns an F on each problem, but works very hard on their corrections, they could earn a D on the exam. Thus, a student would hypothetically earn an F on every exam problem and still pass the class through hard work. This is deliberate. I want to lower the stress of exams and encourage students to dig deep into the material.*

The concept of a stretch course addresses many assessment concerns. However, the enactment of the policy is important to reflect upon. Specifically, this stretch course "was offered for any student who felt either academically or emotionally unprepared for Differential Calculus" (Stacy, (CS 23 )), and while the university is a Predominantly White Institution (PWI), roughly half of the students in the course in 2019-2020 were Black women. How we place students, or how we ask them to self-evaluate and self-place, can perpetuate existing structural inequities experienced by people historically minoritized in mathematics, in this case Black women (Larnell, 2016).

The second common assessment reform is utilizing group exams (Kinnear, 2020). Kinnear (2020) investigated two staged exams in the context of undergraduate mathematics. In Stage 1, students would complete exams individually, and in Stage 2, they formed groups to complete the exam again. While Kinnear (2020) found that two-stage exams do not have significant impact on longer-term learning of mathematics, they argued that "instructors may still wish to use them to emphasize a collaborative classroom pedagogy" (p. 1). Moreover, more robust studies are needed on group exam effectiveness. In this volume, two case studies discussed how they utilized group exams.

Starbird et al. (CS 24 ) implemented a two-stage exam structure similar to Kinnear (2020). The first hour of the exam was individual. During the second hour of the exam, students were randomly placed in groups of three. Exam questions were either identical or variations to problems found on the individual part of the exam. Each student earned the same grade on the collaborative piece, which accounted for 30%-40% of a student's total exam grade. Starbird et al. (CS 24 ) "view[ed] those concentrated discussions during the group portion of the tests as among the most effective learning moments that the students experience."

MacArthur (CS 17 ) implemented a group exam structure, with the explicit goal of rehumanizing mathematics (Gutiérrez, 2018) by "creating a humane learning environment, produc[ing] higher exam scores, improv[ing] student confidence and increas[ing] students' sense of belonging" (p. 1). This case study focused particularly on how Black, Latinx, and Indigenous students, "who were the intended targets of the rehumanizing framework, and women experienced this exam change, given their historical underrepresentation in college calculus courses." MacArthur (CS 17 ) provides specific details on the group exam structure. In contrast to Starbird et al. (CS 24 ),

MacArthur (CS 17 ) implemented the group portion of the exam prior to the individual portion of the exam. The initial exam grouping was random but subsequent groupings took into account "student needs, personalities, gender, anxieties, etc." MacArthur found that Black, Latinx, Indigenous, or women students approved of the group exam structure because it felt more humane, supported collaboration and cooperation, and mimicked what they might see in their future careers. These findings are encouraging, but MacArthur concludes with powerful and important caution:

> *It's prudent to recognize that simply changing exam structures, without the foundation of implementing equitable and inclusive pedagogy during class, will likely not have similar success for learning outcomes and belongingness of calculus students.*

Overall, group exams, and their inherent collaborative nature, can provide a more authentic opportunity for what it means to do mathematics. The important take-away is that the structures, rules, and policies that are in place for the group exam need to be carefully looked at and defined to ensure equitable experiences for students. Again, no one measure will give us an accurate view by itself.

### 3.2. *Alternatives to Timed Exams: Portfolios, Projects, and Oral Exams*

There are ways to assess student learning of calculus that are not seated, timed exams. These include but are not limited to portfolios, projects, and oral exams, and several case study authors describe using such alternatives. Byrne et al. (CS 5 ) discussed the use of a group project that concluded with a final written report and presentation. While the use of writing in calculus is not common (Cooley, 2002), portfolios and group projects afford students the opportunity to demonstrate their understanding of mathematics in different ways. As Burn and Mesa (2015) pointed out, these types of assessments can seem harder to quantify (in comparison to sit-down, timed assessments) but they promote higher order thinking. Moreover, portfolios, group projects, and oral exams more accurately reflect the work that academics, mathematicians, engineers, biologists, etc. all do because they exhibit collaborative complex and real-world phenomena. Boedigheimer et al. (2015) "found that this alternative assessment gives the instructor considerable insight into students' [mathematics] understanding, and, in general, students value the opportunity to demonstrate their abilities via this kind of assessment" (p. 99).

These structures have potential to serve as more inclusive practices, to provide a sense of belonging, but implementing them without considerations of power dynamics is likely only to be a different version of assessment that still marginalizes the historically marginalized in mathematics.

### 3.3. *Grading Schemes and Transparency*

Another important piece of assessment is transparency in grade determination and how much that grade is worth in terms of a final course grade (i.e., the weight and importance of some items over others). Two case studies discussed how their departments took steps to improve transparency and reform how final course grades were calculated. Akin & Viel (CS 1 ) discussed how their grading scheme was reformed in coordinated Calculus courses because of a growing concern that their previous system had unintended inequitable consequences. They were motivated by equity (students expressing different experiences), transparency (students were unclear of their standing in the class before the final exam), and inclusion (students were being evaluated in comparison to other students). Their department was first using a block grading system for coordinated calculus courses (see their case study for details). They shifted to a mastery grading system with no exams and a resubmission cycle, where students were getting regular feedback about their progress in the course.

Similarly, Bennett et al. (CS 3 ) discussed shifting assessment of student learning to mastery-based grading. This shift was supported, importantly, by department faculty and the university at large. They shifted from a grading scheme where 95% of a student's grade came from three timed high-stakes assessments to a grading scheme where 50% of a student's course grade comes from repeatable mastery assessments, 10% comes from work done in class and class assignments, and 40% comes from two timed exams.

Grading schemes are often structures that faculty set individually, or for coordinated courses may be set at the course or department level. Research has shown the adverse impacts of high stakes assessment on minoritized students (e.g., Madeaus & Clarke, 2001). By reforming grading schemes, we directly (and indirectly) lower the stakes of assessments for our calculus students. Mastery based grading is a frequent topic of discussion in mathematics education. Oftentimes, its implementation (and thus success) varies widely. For a deeper look into mastery-based grading, read the *Special Issue on Implementing Mastery Grading in the Undergraduate Mathematics Classroom* in *PRIMUS* (Campbell et al., 2020).

### 3.4. *Guiding and Reflective Questions on the Role Assessment in Evaluating Student Learning*

- How do faculty in your department currently assess student understanding of calculus?
- Are there guidelines and expectations for fair treatment of students by instructors in assessment practices?
- How do you gather feedback on your assessment practices inside and outside the classroom?
- In what ways can your department make reforms in evaluating student learning?
- How do assessment efforts relate to diversity, equity, and inclusion efforts in your department?

## 4. *Concluding Thoughts*

Case studies highlighted in this chapter have proposed several ways to reimagine the role of assessment in placement and evaluating student learning. In terms of placing students, multi-measure approaches are crucial. Additionally, using a corequisite model provides necessary support for students. In terms of evaluating student learning, many reimagined policies occurred. Faculty are using exam revisions, groups exams, and portfolios, as well as reforming how they calculate overall course grades. The case studies all had similar starting points. In most cases, the impetus for change was improving pass rates in calculus. The force behind that impetus sometimes differed. Sometimes, the pass rates for "all" students were not good. Sometimes, authors noted that certain groups of students were not succeeding compared to others (e.g., Mingus et al. (CS 20 ), this volume, explicitly pointed out Black or multiracial students and students who were first-generation college or Pell-eligible). Additionally, some authors specifically were targeting the experiences of certain groups of students (e.g., MacArthur, (CS 17 )). We commend the case study authors, and their local institutions, for the initiatives they have undertaken.

We are encouraged by those who are centering the lived experiences of minoritized and underrepresented students in their revisions of assessment. Our intention with this chapter is to push this conversation forward. We hazard against the approach of reforming assessment and placement with a 'better for all' framework and reinforce calls for educators to center the experiences of marginalized and underrepresented students in all assessment reimaginings. In 2003, Dr. Danny Martin discussed the inherent problems with "mathematics for all" rhetoric and action. Similarly, we argue that we cannot approach this assessment reimagining simply from the point of view that it will make assessment better for everyone. Mathematics spaces are dominated by whiteness (Battey & Leyva, 2016; Martin, 2013) and students have different racialized and gendered mathematics experiences (Leyva et al., 2021). Further, the male superiority myth (Leyva, 2017) has dominated for years, only recently being pushed back against, with, for example, counter-stories of resilience of Black women in mathematics (Leyva, 2021). Initiatives must be explicit in their intention to benefit marginalized students in mathematics. If initiatives are only targeted to improve calculus, in our case, assessment, for all, then we are continuing a meritocratic perspective on mathematics assessment and achievement.

We return to the discussion of the "assumption that standardised tests provide fair and objective measurement of individuals" (Au, 2013, p. 13). Such a perspective on "the low achievement on standardised tests of working-class people, non-white populations, and some immigrant groups can then be simply and neatly attributed to the failure of individual students, individual groups, or individual cultures, and not attributed to existing structural inequalities" (Au, 2013, p. 13). There are many avenues for reimagining our placement and assessment practices. We argue that faculty cannot look at their individual, classroom, department, university, or policy level changes as if they are trying to make assessment practice inherently fair. **We always have to explore the implications of our attempted reimagining. We do not solve this problem in one attempt, we must acknowledge that this reflection is continually part of the work.**

Our intent with this chapter is not to provide all the various roles of assessment and placement one might try; we refer the reader to the phenomenal work of the individual faculty and mathematics departments as a whole in the case studies in this volume. Rather, our intent with this chapter was to provide a spark for the reader, an initial point from which to start in making calculus more equitable. We leave you with these final guiding questions:

- How can we re-imagine the role of assessment for placement into calculus to better support students through calculus?

- How can we re-imagine the role of assessment in evaluating student learning of the concepts and practices of calculus?
- What do we care about that our students learn in calculus?
- What do we care about that our students learn *from* calculus?
- Are the changes we make specifically focused on improving the experiences of marginalized groups in mathematics?
- As stated in the framing of this volume: Does our work "build a framework … around the concepts of justice and calculus"?

# 8 Reflections From a Student Advisory Committee to Deepen Equity in STEM Classrooms

**Nyote J Calixte**

*Elon University*

Nyote J. Calixte, Ph.D. is an Assistant Professor of Chemistry at Elon University who serves as steadfast proponent for Science, Technology, Engineering, and Math (STEM) endeavors and is an aspiring diversity, equity, and inclusion in STEM leader. She earned a Ph.D. in chemistry from Louisiana State University and a B.S. degree in chemistry with a mathematics minor from Xavier University of Louisiana. Dr. Calixte's research projects focus on the use of student voice as a tool for institutional STEM curricular reform as well as demonstrating how metacognition and Bloom's Taxonomy can be employed in general chemistry classrooms to improve student performance. Dr. Calixte also is part of a grant writing team at Elon focused on transforming introductory STEM courses to promote student sense of belonging.

**Ricardo Mattias**

*Independent Researcher*

Ricardo Matias is an Associate Director at Reach Higher, former First Lady Michelle Obama's education access initiative, where he focuses on improving education access for historically underrepresented students. Prior to this role, Ricardo was a Business Analyst at McKinsey & Company, where he primarily served higher education institutions in addition to other public and social sector clients.

**Kristel Black**

*Independent Researcher*

Kristel Black is a former post-baccalaureate student from University of Pennsylvania and former clinical research coordinator in the Gynecologic Oncology department at Penn Medicine. Kristel is currently pursuing her MD at the UNC School of Medicine as a Kenan Urban Primary Care Scholar.

**Lesley Ihionkhan**

*Independent Researcher*

Lesley Ihionkhan is currently a full-time MBA student and Rath Foundation Merit Scholar at The University of Chicago Booth School of Business. At Chicago Booth, Lesley is a co-chair of the African American MBA Association and pursuing Digital Product Management and Entrepreneurship.

## 1. *Introduction*

Student perspectives serve as a critical yet often untapped resource for institutional reform and interventions. When students partner with faculty and staff to examine educational challenges within their institution, students become empowered change agents equipped with a deeper sense of belonging and leadership. One way that students can demonstrate leadership as "pedagogical partners" with faculty is through student advisory boards. Student advisory boards like the Diversity in STEM Advisory Committee (DSAC) take understandings acquired from student experiences in STEM and engage campus leadership in challenging yet necessary conversations about student belonging, equity, and bias that challenge student belonging and learning. In this work, former DSAC students will reflect on their experiences in the sciences, and we will contextualize their voice within some of the themes of this volume, namely ideas around course changes, professional development, and classroom practice as well as examine how student perspectives are utilized in the case studies depicted in this volume

Many underrepresented minority students enter college with a commitment and interest in Science, Technology, Engineering, and Math (STEM) majors. Yet, many change their intended path of study following their initial experiences in the gateway STEM courses (Olson & Riordan, 2012). This is especially the case for Math (Olson & Riordan, 2012). The challenges that face retention in STEM are numerous (Sithole et al., 2017), yet these issues are compounded when students do not feel as though they belong within the scientific community (Wilson et al., 2015; Rainey et al., 2018). While addressing pedagogical challenges is a primary necessity, we must also aim to create an equitable and diverse community for science and mathematics learners to effectively transform STEM education. Equitable and diverse learning communities intend not only to impart knowledge and skills to the learner, but also to empower students to become active contributors to this environment (McGee Banks & Banks, 1995). To achieve equity requires two elements: (1) fairness, and (2) inclusion; often, students do not feel STEM classes possess these attributes. For learning environments to be equitable, there must be interventions that work in widely diverse settings and for diverse student populations. Often, methodologies to increase equity in STEM involve creation of programs or initiatives intended to promote student success and retention, however they infrequently involve the voices of students (Gibau, 2015). Underrepresented students with interest in the sciences endure various challenges despite not having interventions that align with their experiences. While there are interventions empowering students to use their voice to promote change (Nthontho, 2017), these approaches tend to be underutilized when developing and implementing many STEM-based interventions (Chase, 2020). There is an abundance of research on interventions employed to address the lack of underrepresented students in the sciences. However, educators should more routinely consider the voices of underrepresented students themselves in the design and implementation of various programs as an intervention (Gibau, 2015). There is a reemerging and steadily growing body of knowledge and literature related to the impact of student voices in higher education interventions, like instructional and curricular design (CS 6 , CS 29 , CS 12 ), and retention (Barker, 2018; Bron & Veugelers, 2014; Latham & Gross, 2013.; Richardson, 2019.), yet there is still a need to better understand how these student perspectives align specifically with institutional initiatives and interventions in STEM (Petersen et al., 2020).

To better understand the role of student voice, and why it should be an important aspect of any institutional conversation, it is important for us to first describe what we mean when we use the term “student voice.” Student voice is a concept in education research and reform that has been used to highlight the absence of student voices in decisions around teaching and learning. Since its emergence in the 1990s, the term “student voice” and our understanding of it has evolved and been debated over time (Cook-Sather, 2006). Seminal work by Cook-Sather explores the nuance behind student voice lexicon and states:

> *Currently, many people are using the term “student voice” to assert that young people have unique perspectives on learning, teaching, and schooling, that their insights warrant not only the attention but also the responses of adults, and that they* *should be afforded opportunities to actively shape their education. (Cook-Sather, 2006, p.28)*

Student voice implies students participate in deciding which structures and practices influence their educational experience ("Student Voice," 2019). This means the student voice is action oriented and should have an effect. As described by Lundy (2007, p.933), four elements are necessary for student voice, (1) space: students must be given the opportunity to share their experiences, (2) voice: it must be easy for students to share their experiences, (3) audience: student's point of view must be listened to, and (4) influence: student experiences must be acted upon, where deemed appropriate. Cook-Sather's work notes the presence of the terms "rights", "respect", and "listening" in literature where student voice is used yet caution that to move beyond merely evoking these terms requires a cultural shift in which administrators, faculty, and policy makers value the change promoted by students. For exploring this work, Lundy and Cook-Sather's characterizations of student voice will be the basis for the reflections provided as they capture insights into student learning, promote action, and offer poignant examples to educators as how to use this perspective to shape STEM education reform.

Student voice, having both positive and negative aspects, warrants care when used to describe student perspectives to ensure it does not reinforce unjust power dynamics and undermine institutional reform efforts (Cook-Sather, 2006). In higher education, student voice has often found a home as "student feedback" in the form of teaching evaluations, yet these can be biased against faculty, namely women and people of color (Cook-Sather, 2020; Hamermesh & Parker, 2005; Mengel et al., 2019). Consequently, educators can be apprehensive about faculty-student partnerships that give student voices excessive authority (Cook-Sather, 2020). However, "pedagogical partnerships" between faculty and students that equally involve both voices, allow students to make meaningful contributions to initiatives and establish some autonomy in the decision-making process (Cook-Sather, 2020). When considering equity-seeking groups, recent work by Cook-Sather examines how student voice can be utilized to better serve traditionally underrepresented students:

> *Because of the potential of participating in pedagogical partnership to increase students' sense of belonging, deepen their academic engagement and encourage their persistence, all of which are linked to success, student partners offer perspectives that are particularly important to defining and supporting success for those traditionally underrepresented in and underserved by higher education. As one student partner put it, the presence and participation of students from equity-seeking groups helps 'disrupt the exclusive nature of higher education and helps level the playing field to allow for students to achieve more than they thought possible before their voices were included' (Cook-Sather, 2018, p.4)*

Partnerships between faculty and students have the potential to contribute to student success by disrupting exclusionary practices and hierarchical institutional norms within the higher education sphere (Mercer-Mapstone & Bovill, 2020). Students can become empowered by partnerships with faculty/teachers (DeFur & Korinek, 2010) and exhibit more engagement in their education (Cook-Sather, 2007). Student voice has also proven effective in encouraging students to value long-term success and stress the significance of democratic values such as participation and leadership (Bron & Veugelers, 2014). Applying student voice as a tool to establish students as educational partners promotes inclusion (Yonezawa & Jones, 2009) and advances equity when diverse perspectives are incorporated. Student voice will be the basis for the reflections provided in this chapter and we acknowledge that the individual experiences and perspectives described herein should not be taken as broad generalizations. Instead, the goal of the student voices provided in this work are to demonstrate the unique perspective that students can provide when space, voice, audience, and influence are prioritized within a student advisory committee of an institution.

## 2. *The Student Advisory Committee*

Institutions have investigated various mechanisms to elicit student voice to better support student engagement and success (Booth et al., 2013; Jensen, 2020; King Miller, 2015; Yonezawa & Jones, 2009). Student advisory boards offer a formalized opportunity for students to become active partners in institutional reform efforts (Bacon & Bloom, 2000). This chapter will discuss the creation of one such faculty-student partnership to elicit the voice of students engaged in STEM-based areas of study (biology, chemistry, computer science, engineering, environmental studies, mathematics, physics, statistics, psychology, and pre-health) with backgrounds from underrepresented groups as defined by the National Institutes of Health (*Women, Minorities, and Persons with Disabilities in Science and Engineering: 2021*, 2021.). This student-based advisory board was created in 2016 at a highly selective, top-tiered private research-intensive institution we will refer to as Dunshire University. This committee was created in partnership with faculty and staff to bridge the hierarchical divide within the institution. Moreover, the Diversity in STEM Advisory Committee's (DSAC) goal was to take the insights gained from students who shared their academic experiences on campus, inside and outside of the classroom, to help promote broader diversity and inclusion for students of color in STEM pathways. Students often shared insight regarding mathematics courses required by their academic path and this information is subsequently shared with campus leadership. Using demographic data obtained by our university research engagement office, students of color with interest in STEM disciplines were identified and invited to join DSAC. Students self-selected and provided short purpose statements regarding why they desired participation in this committee to determine level of commitment. Purpose statement responses were overwhelmingly related to students' desires to help junior classmates avoid the same pitfalls and negative experiences in STEM courses. Meetings for DSAC took place in a campus identity center to ensure students felt the environment was safe to be vulnerable (Patton, 2010). Meals were provided to students during meetings to safeguard student's attendance against competing events with food being served. This point may seem trivial but served as an added logistical measure because meeting times were often during dinnertime and some participants expressed some level of food insecurity (Payne-Sturges et al., 2018).

Students on this advisory board shared their academic experiences on campus, inside and outside of the classroom, to help promote broader diversity and inclusion for students of color in STEM pathways. Contributors were asked questions akin to, "Was there ever a time in the classroom when you became more aware of your race?" or "What experiences or interactions with your STEM faculty or peers stuck out this week?". Students often voiced insight regarding mathematics courses required by their academic path because of its connection to many STEM disciplines. DSAC members also shared experiences that reflect biases they encountered that affect their learning environment, relationships with instructors, and their sense of belonging in the institution's scientific community. One such student account highlighted the biases faculty can possess about student aptitude based on race:

> *In my final project we were tasked with creating a business case and pitch for a very specific niche topic. I felt quite strongly about my topic and decided to do research on my own a couple weeks before that project was due to put together the best business case. I had rehearsed on this final project probably more than I had ever iterated on a final project for a class out of pure excitement and when I turned in the project two days later the professor emailed me back saying, "Hi Lesley, thanks for this. Are you sure this is your work" I felt sick to my stomach. I went up to my professor after class and said, "thanks for reading my final project. I understand you're confused if I had done this work. I cited all my sources at the end of the paper using the format that you asked for. I'd be happy to give you a little bit more context as to how I extract information from those sources and arrived at my own opinions, but I would love to make sure that you understand this is my work" and he said to me again, the only black person in class, "Oh yeah. OK. if it's in the end that's fine. I just didn't think someone like you could produce that work." (Lesley)*

Stories like these were shared to the group, so members would be required not to disclose topics discussed outside of sessions to maintain anonymity. Students were also encouraged not to divulge faculty names as a secondary means to protect their identity. Accounts of student experiences were de-identified and discussed in administrative meetings with deans, faculty, and staff. Their feedback would contribute to the design and broader programmatic implementation of campus wide interventions that promote equity and inclusion in STEM. As part of a research study, one DSAC student reflected on her experience in the STEM environment at her institution by asserting, "you don't want to come into a broken system" (Tremaine et al., 2020). Her perspective highlighted the tension students experience when reconciling how institutional efforts still miss the mark. An important aspect of DSAC was the staff-student partnership that helped validate students' voices while also assisting them with contextualizing their experiences. This ensured students remained introspective, critically thinking, and without victimization.

Challenges with such a program could be the lack of diverse faculty with whom students feel comfortable sharing their experiences (Guiffrida, 2006). As Golden et al. express in their case study (CS 12 ), sensitivity to facilitator identities is necessary due to the impact this can have on student willingness to share and influence the safety of the sharing space for BIPOC students. Facilitating the creation of faculty-student partnerships like DSAC could necessitate tapping professors within your institution from similar identity groups and would mean increased service-leadership commitments that can disproportionately impact faculty of color (Trejo, 2020). Identifying students to join "in the work" could also promote an elitist approach that could prioritize the voices of students who are already "privileged and engaged" (Mercer-Mapstone & Bovill, 2020). Lack of funding also could pose a challenge as we often discussed student experiences in STEM classes over a meal to promote participation and community. Lastly, having a safe space for students to share their stories was an important factor in implementing this committee. Often, the stories shared by students were painful to recount and the vulnerability that many students displayed made it important to prioritize having a space where students felt safe and comfortable. Institutional adoption would require faculty commitment and student buy-in to be effective. Students who participated in DSAC felt encouraged that their voices would reach administrators and help students coming after them. We also encouraged student participants to include their service as members of DSAC on their CV as a demonstration of leadership, which promotes equity and student agency ("Student Voice," 2019).

## 3. *Student Perspectives*

Student voices are a powerful tool for instructors and institutional leaders to gauge program buy-in and student engagement. Institutions have invested considerable time, research, and funding attempting to address attrition in STEM through curricular overhaul and program implementation (Granovskiy, 2018), yet student voices are consistently underutilized as interventions in STEM. The student perspectives featured herein reflect the experiences of former DSAC students who are now alumni. These students obtained bachelor's degrees in areas of biology, sociology, and environmental science and policy. Although each student began their academic journey with a STEM or health professions career in mind, only one student is currently pursuing a STEM/health professions career. Their perspectives will shed light on their experiences in STEM courses, what factors contributed to their career shift, whether their STEM environments were equitable and inclusive, and how instructors can consider student voices in the design and implementation of future STEM courses.

### 3.1. *Ricardo.*

Ricardo, a biology pre-med major, became involved in DSAC his sophomore year. As an ever-vocal contributor to our conversations, he often reflected on his experiences in the classroom noting the disparities that existed in the spaces he occupied. For greater context in this chapter, Ricardo recounts his experiences with the Math 105/106 sequence at Dunshire University. This two-semester course serves as a course equivalent to the institution's standard one semester calculus course. This two-semester option is largely for students without AP/IB credit earned, little to no prior calculus, and/or certain SAT scores. Here he provides educative introspection regarding his experience in STEM classes:

> *As a pre-med student, I was no stranger to the STEM classrooms at my Institution. During my first semester on campus, I took both Chemistry 101 (Core Concepts in Chemistry) and Math 105 (Laboratory Calculus and Functions I). The following semester, I would take Chemistry 201 (Organic Chemistry), Math 106 (Laboratory Calculus and Functions II), and Neuroscience 101 (Biological Basis of Behavior). I took note of the fact that a greater proportion of my peers were students of color, particularly Black and Latinx students, in my math courses, than in my Chemistry and Neuroscience courses.*
>
> *Math 105 / 106 is one of several math pathways offered at my Institution, particularly for students that either a) had limited calculus exposure in high school, or b) took the AP test but did not earn credit. However, matriculating into university without any AP credit in Math does not bode well for students interested in the Math major. In addition to Math 105 and 106 (which are not offered simultaneously),* *prospective Math majors are required to take 4 additional prerequisites (such as Calculus II, Multivariable Calculus, and Linear Algebra). This is in addition to the 7 courses required for the A.B. Math degree, or 8 courses required for the B.S. Math degree.*

*Any expectations or aspirations to graduate with such distinctions were leveled immediately. On our first day of class, our professor made it clear to us that, as students of Math 105, none of us would go on to become Math majors. This same professor, however, became an enormous champion and mentor to me, and built in me the confidence that could move mountains. They were incredibly patient, made themselves available at all hours of the day for individual support, and most importantly, met our academic and non-academic struggles with nothing but genuine kindness and understanding. At one point, falling in love with the subject and finally finding myself within a true support structure, I had even convinced myself that I too could major in Math. Throughout my first 2 years at university, this professor wrote several letters of recommendation for various internship and research opportunities. Although I did not end up pursuing the Math major, as I reflect on this experience post-graduation, I realize that I was so fortunate to have built such a powerful academic foundation and mentorship network early on.*

*Coming back in the spring to Math 106, I was excited to continue building on my positive experience from the previous semester. Unfortunately, circumstances did not favor my expectations. Instead of having the same professor from Math 105 as my instructor, I was instead in a section led by one of the Math's department's graduate TA's. Understandably, this TA, though unanimously a talented practitioner, was not able to extend the same quality of instruction or student support as my professor from the fall. Several of my classmates that made the same transition had similar observations. In fact, my success in Math 106 is in large part due to the individual appointments, office hours, and email exchanges with my professor from the prior semester.*

Ricardo's experience demonstrates the power of authenticity in mentorship and the importance of building an academic network. Interestingly, Ricardo was not deterred from continuing in Math 105 after his professor's assertion that it was likely students in the course would not be able to successfully pursue a math major. Luckily, this comment did not prevent Ricardo from fostering a positive mentorship relationship with this professor. However, it is likely that other students in the course could have felt discouraged, which is evidenced through the case study by McNeill and co-workers (CS 19 ) that suggested Black and Latinx students interpreted a course drop as racialized pressure and promoted doubt with regard to their sense of belonging. Despite this event, Ricardo felt the professor allowed students to show up as their whole, authentic selves and extended understanding and patience within the discipline. While math can be intimidating, this student was able to grow in confidence, in part, because of the trust and supportive relationship garnered from this professor. Ricardo's account also illuminates how pedagogical training for graduate student workers can aid student retention and boost success. In the work of Bennett et al. (CS 3 ), attention to graduate students as the primary workforce for introductory calculus classes highlights the need for better training practices. This case study discusses how they incorporate inclusive teaching practices into their training program for novice instructors while increasing support for instructors around DEI issues.

3.2. *Kristel.*

Kristel, a sociology pre-med major, became a member of DSAC through recruitment from a pre-medical scholars' program for minority students. Kristel was always candid about her vulnerabilities and how her identity as a black woman shaped her experiences in STEM. In this section, she provides deeper context regarding her path and how her background and the culture of STEM at her institution collided:

> *The most difficult mathematics course that I had taken in high school was AP Calculus A/B. My high school could not offer the second AP Calculus course, because we did not have any teachers who could teach it. Nonetheless, I hoped that my math background would be an adequate foundation for the math and science courses that I would be taking at my university. I soon came to realize that many of my classmates had taken higher mathematics courses. This difference in knowledge and practice meant that I had to take the slower-paced courses to fulfill my pre-medical requirements: Math 105 and 106. Most of our peers could start off at Math 111 and 112 or even skip to 212. All these adjustments are based on SAT score and AP exam scores. In my frustration with the system, I asked several students if I should take 111 and get my math requirement "over with" or strengthen my math foundation and take both Math 105 and 106. I was met with an influx of comments suggesting that my performance and the grade earned in one math class was better than subjecting myself to two math courses. However, I agreed with older peers that argued it was more important to have a stronger understanding of calculus, avoid only focusing on grades, and revisit math over a longer timeline, especially since my high school did not have a strong AP math department. I did not mind taking math at my university at a slower pace, since I value genuinely understanding topics to master them.*
>
> *I have spent a great amount of time reflecting and troubleshooting reasons why I did not thrive in STEM courses at my Institution. I only remember having one or two classmates that I could study with and prep for exams, since the collective class mentality was majorly competitive. Consequently, the math curve would cause everyone to hyper focus on the letter grade and beating the curve thus limiting class collectivism in the spirit of competitiveness and 'learning.' Except this wasn't learning. We did not have enough time to learn, while we had multiple assignments due. Although we had mathlab, where students could get additional assistance with math, many people did not have any interest in coming together and work through assignments as a mechanism to limit how many people CAN do well in the class. These are simply the behaviors of the STEM community at my university. Some people benefited from the curve, but I remember only ever being average or slightly below the average with each exam. My testing anxiety worsened significantly in my math and science classes, especially since I would leave questions blank, fill out as much as I knew, or my only goal was to complete as much of the test as I could. After a while, I think I expected myself to just do worse than my peers. I will say that I did not attend office hours religiously, since I would have conflict classes or extracurricular activities that I would prioritize. Sometimes professors would say, "Don't come to office hours without trying the homework or other problems," and it would deter me to attend if I* *genuinely did not know how to attempt certain problems. I wish the math professors could have centered the students' learning over meeting an arbitrary standard of understanding in a class.*

In Kristel's experience the classroom culture, influenced by course infrastructure, promoted inequity, and divide in her calculus classes. This environment led to marginalized students feeling isolated and further reinforcing power dynamics that limit access to success. Like the work of Akin and Viel (CS 1 ), Kristel experienced a "block grading system" that bred inequity and discouraged collaboration and cooperation due to the curved standard in the courses. This approach incentivized a competitive approach to calculus in an attempt for students to "get over the curve". As Kristel suggests, this case study addresses these issues through curricular overhaul and some experimental changes that are student success centered. In addition, Kristel's experience with considering a higher level math course to get it "over with" is supported by the work of Stacy (CS 23 ) where students would opt to take a differential calculus course that they were unprepared for. Stacy's study highlights how a "stretch differential calculus" course was created to respond to the DFW rate and placement challenges in differential calculus that deepened marginalization. Of note in this case study was once again the power of student voice, when a black student expressed feeling a lack of belonging due to the differential calculus course. In this case study, Stacy discusses hearing the department was not inclusive from this student, which motivated them to try this "stretch" course approach. This student's voice also served as evidence to "resistant" colleagues that it was no longer acceptable to take a position of inaction regarding BIPOC lack of sense of belonging in math classes at their institution.

### 3.3. *Lesley*

Lesley, an environmental science and policy major, began her academic journey with the hope that college would be a true meritocracy where her commitment, hard work, and stellar high school preparation would be critical nodes of her undergraduate career. Lesley's involvement in DSAC allowed her to vocalize the realizations of her identity and the perceptions others formed about her identity, regardless of the fruits of her labor. Here she will unpack those perceptions and how it shaped her trajectory.

*I loved environmental science since my high school teacher spurred a fire in me and I was always excited about learning about conservation and how we can better interact with the planet around us. With that lens, I started my educational career at my university and was really excited to try to pursue environmental science through my academic studies as well as my extracurricular activities. I took Chem101 in the spring semester of my freshman year and I luckily had quite a good foundation in chemistry from my high school experience. I loved all my science teachers and because of that, I felt very strong coming into the program, but Chem101 was notoriously one of those classes that was designed to "break you down". Now, I say that jokingly, but it truly was a sentiment among not only my black counterparts but also the class in general. This class was offered as a flipped classroom where our lectures were focused on how to solve problems together, so it was our responsibility to learn all the information prior. This was significantly challenging for me even though I loved chemistry and other STEM programs from my high school experience.*

*It was the first time I ever had a flipped classroom experience, but I prepared by spending the first couple of weeks before the semester starting to read through the textbook and familiarizing myself with the material. I was also able to get my hands on a couple of resources that others had provided from those that had taken the class before, so I felt like I had a pretty good foundation going into a challenging class. Adding to the complexity of this flipped class structure was that we were put into groups, and I was the only person of color in my small group of about five people. I remember trying to chime in sharing the answers to the questions I worked on the evening prior to class, and I distinctly remember my classmates saying, "well how did you arrive at that?" I went through the process of how I arrived at my answer, citing my sources, and explaining what I learned from the textbook, and I remember a white female classmate said "well, let's just double check your answers with the TA." This was frustrating to me, and I remember being annoyed and wondering why everybody else in this group was sharing answers with each other, were helpfully challenging where others answers came from, but nobody else solicited the response of, "let's make sure the TA checks your answers". After that encounter, I stopped connecting with the group and instead, I figured I'm just going to ace the class and show them who's really in charge and show them what they really messed up, so I did. It was unfortunate. Looking back, I wish that I would have challenged her question, but in the moment, I felt so frustrated and angry, and I did not have the agency or opportunity to advocate for myself. I do wish that in that situation there was more support from the professor, especially when it comes to creating diverse small groups. Also, I wish that there was more guidance from instructors in determining what is appropriate communication and interactions within those groups, as I certainly felt quite ostracized. Unfortunately, I did have the unfavorable experience of being discriminated against because I was the only non-white member of our group, but my calculus class in my sophomore fall was a much more enjoyable experience.*

*Though calculus was considered another difficult "weed out" course I decided that I was going to really take advantage of office hours and go there every single week. Jokingly, I guarantee you that the TA teaching my lecture was tired of me by the end of that semester. I spent every single week, every single office hour, going in asking questions and doing my homework in those office hours. I think that was truly the key to my success and that TA, while he occasionally displayed his visible frustration, was otherwise open to having this experience with me. I was often the only one in office hours and because of that he was able to really spend time to help me understand how to reframe my existing knowledge of calculus. Another element that really helped my calculus experience was having a buddy throughout the entire semester. This was extremely informal, but the camaraderie and environment that was created was nice and intimate. It was a totally different experience than the experience that we had with chemistry 101. I did not want to major in math, so beyond meeting my major requirements there was no need for me to take additional calculus, however, statistics was necessary.*

*My experience in statistics made me realize statistics is much more than understanding concepts, especially at my university, it involves coding, and I was not prepared to code. This disparity in the course highlighted how students with existing coding experience were immediately advantaged over students with little to no prior*

> *coding experience. I remember sitting in the library computer labs teaching myself R programming at all hours of the night and feeling like, "wow, I'm really not cut out for this". My experience in statistics led me away from the STEM and mathematics experiences that I initially sought out. I decided to withdraw from my second statistics class due to the R programming skill required in the class and how the instructor was quite removed from the experience. I knew quite quickly in the semester that I was not going to get the support I needed in the course from my professor, and as a result I wasn't going to succeed.*

Lesley's testimonial emphasized the importance of identity and how interactions around identity can influence the student experience. Her voice also suggests the need for faculty to prioritize vulnerability and create space for that in their classrooms. Lastly, Lesley's experience in statistics underscores the need for equity in curricular design. Her experience coming into statistics with no coding exposure and yet needing this background to do well in the course highlights a major disparity in statistics at her university.

## 4. *Voices in Context*

The premise of this volume is to highlight pedagogical approaches in the areas of mathematics education that address diversity and inclusion through several topics. Subsequently, this chapter would be remiss if it did not revisit the power of student voice within the context of the most relevant themes examined in this volume. Several case studies throughout this volume incorporate student perspectives because of the power it holds to help transform interventions. Golden and Hagman (CS 12 ) illustrate how designing an intervention can shift when student experiences are integrated into the course design process. Students in this study journaled about their course and provided feedback in a manner that necessitated restructuring of the course to better meet their needs.

### 4.1. *Classroom Practice*

As highlighted in earlier chapters of this volume, educators can explore classroom practices that examine frameworks for building more equitable and inclusive courses. When bearing in mind Kristel's voice, the "Rehumanizing Mathematics" framework holds deepest significance. As evidenced by the case study conducted by MacArthur (CS 17 ), their approach to rehumanizing exam structures resulted in lower DWF rates and increased final grades for the low and medium score groups. This action-oriented approach also utilized student voice through focus groups and surveys where students reported higher levels of confidence, lower levels of anxiety, and a stronger feeling of belongingness compared to traditional exams. Designing for collective responsibility rather than individualism and competition brings this classroom practice into alignment with student expectations. In addition, the case study of Wagon et al. (CS 29 ) examines how to promote equity in a STEM math pathway program and harnesses the power of student voice as one informal metric for program evaluation. Through the work of Davis (CS 8 ) we observe that student perspectives are incorporated via intake and exit surveys to evaluate the efficacy of a peer-tutoring program. With regard to classroom practice Kristel reflects:

*It should never be about the curve, the grade, the competition between peers, but instead fueling the curiosity for learning within every student. Learning should be fun. It should be enriching. It should make you want to continue building relationships with your peers and professors and to continue doing similar classes or opportunities that are related to that topic. It should not make you feel like you deserve to be below average from your peers because the class structure is not conducive to certain students' success. Professors should consider, "how was this a conducive environment for collaboration, mastery, or learning?" I truly believe that the structure of the calculus curriculum with homework heavy approaches, high-stakes summative assessments that count for a considerable percentage of our course grade, course curves, and the competitive atmosphere left many students behind and in the learning margins. Unfortunately, math at my university did not fulfill me, but it was an experience that illuminated how real imposter syndrome is and until the structure of these courses are changed to be inclusive and engaging for ALL students, students will continue feeling isolated, frustrated, discouraged, and on the margins. (Kristel)*

### 4.2. *Professional Development*

While diverse and insightful, the case studies on the development of change agents in mathematical instruction shared the common perspective that professional development is necessary for all instructors, at every level of training. Several case studies in this volume speak to the necessity for effective and appropriate professional development. Based on their experiences, our student contributors echo these sentiments:

*As I continued to pursue my studies at my university, the dichotomy between my Math 105 and Math 106 experience stuck with me. What interventions can universities pursue to better train graduate students for instruction? If there are departments that rely on graduate student labor, what are the learnings and best practices they can implement? Most importantly, how can universities ensure that first year students, particularly those that have been afforded less academically enriching opportunities in high school, receive access to a powerful foundation and support network early on? (Ricardo)*

*I would love to see the integration of inclusion policies not only among students but also for professors themselves. I'm sure if my professor senior year was coached in how to appropriately and equitably manage and teach a diverse classroom, a comment like, "I didn't think a person like you could produce that work" would not have existed. I also think a comment like, "well, we should have a TA check all your answers" might have existed, but there would have been the permission and agency to challenge that. Being empowered to take that experience to my professor or to my TA and not be fearful of retaliation would have put me in a better position. There is so much room for opportunity within the STEM ecosystem to improve the experiences of students of color and ultimately increase the pipelines of minority academics as well as STEM professionals. (Lesley)*

### 4.3. *Course Changes*

In this volume, course redesign to establish equity-directed practices is examined. Course transformations that involve major overhauls have the benefit of endless possibilities. Changes were anchored in the department's strengths and the needs of students. For example, the work by Burn and coworkers (CS 4 ) illustrates how developing practices that accelerate students to college level math broadens the options for paths through the discipline. By reflecting on her experience, Lesley shares how changing the size of her course transformed her experience:

> *Where possible, a huge part of having positive experiences for me was smaller class sizes and access to the instructor like my experience in calculus. Having that small class size made our environment safe so we could all be vulnerable about the things that made us nervous regardless of race or background. Also, having access to the instructor and their consistent encouragement to come to office hours and supporting us to ask questions in class meant that it was OK to engage with the professor. That was such an incredible experience in calculus and looking back again, I don't have a passion for math, but I wish I took more calculus classes. (Lesley)*

## 5. *Conclusion*

Student voice is a powerful instrument for creating broader pathways for equity and inclusion in STEM. While many options for partnering with students exist, student advisory boards show promise as inclusive ways to empower members and boost leadership capabilities. Student voices also promote equity when power dynamics are dismantled, and diverse groups are incorporated. By having space, voice, audience, and influence, student advisory boards like the Diversity in STEM Advisory Committee (DSAC) encourage students to use their voice to impart institutional change. Participants in DSAC shared their experiences and discussed how interventions have greater success when student expectations have alignment with institutional initiatives. We recognize and acknowledge that individual voices provided in this work do not support broad generalizations nor speak to all students impacted by the themes discussed herein. The voices in this chapter support the premise that student voice is an important, underutilized, and necessary component of the conversation surrounding STEM reform. We realize case studies around course redesign, professional development, and classroom practice have alignment with student voices featured herein and further support the use of student voice as an intervention in STEM.

### 5.1. *Guiding and Reflective Questions on the Role of Student Voice as Institutional Tools*

- In what ways, if any, does your institution currently lift up student voices?
- Are student voices valued elements of curricular and programmatic reform efforts?
- How do you incorporate student voice that fosters pedagogical partnerships?
- In what ways can your institution dismantle power dynamics that ensure students' voices are listened to and followed up with action?

# Part 3: Illustrative Case Studies

# 1 Equity in Grading Systems: Moving Away From “Fair” Towards Transparency and Inclusion in Coordinated Calculus Courses

**Victoria Akin**

*Duke University*

Victoria (Tori) Akin is an Assistant Professor of the Practice of Mathematics at Duke University. Tori has taught and coordinated courses in the Calculus sequence as well as Probability. Tori’s current research interests lie primarily in math education: specifically graduate student professional development, and the effect of intervention programs on girls' attitudes and beliefs about math.

**Shira Viel**

*Duke University*

Shira Viel is an Assistant Professor of the Practice of Mathematics at Duke University where she teaches and coordinates primarily in the Calculus sequence. Her current research interests focus on math education, including graduate student professional development, alternative grading methods, and students’ sense of community.

*Abstract: Coordinating grades across multiple sections and terms of large service courses is a known challenge, with multiple competing factors. We describe a historical Calculus-sequence grading system that prioritized “fairness” and “consistency” through norm-referenced grading to preset benchmark averages, and the system’s unintended inequitable consequences. We then detail the course coordinators’ varied efforts, under the auspices of the department and Calculus Committee, to reform their grading and assessment systems to make them more inclusive, equitable, and transparent. Experimental changes range from communication-focused (explicit and published formulae for grade computation) to full structural rehauls (mastery grading with mandatory non-letter grades). Design, implementation, and evaluation are ongoing and have been impacted by the continued global health crisis. Yet preliminary evidence is positive, with a decrease in Withdrawals, Ds, and Fs (DFW) across the Calculus sequence, and anecdotal testimony of improved student attitudes.*

## 1. *Introduction*

In the Fall of 2018, the Department of Mathematics at Duke University began implementing new grading and assessment systems in its coordinated Calculus courses. These changes were initiated by the course coordinators, including this article’s authors, in response to a concern that the original system had unintended inequitable consequences.

The coordinated Calculus courses at Duke, from single through multivariable calculus (see Table 1), are overseen by the Calculus Committee. For over 10 years, all such courses operated under the “block grading system”, described in detail in the next section. To keep grades “fair and consistent” across sections, semesters, and courses, grades were curved, or norm-referenced, with instructors expected to keep the average letter grade in each Calculus course roughly the same from year to year. Throughout this chapter, we refer to the departmental guideline for the average course letter grade as the course benchmark. In keeping with historical precedent, lower-numbered courses had lower benchmark average letter grades.

Equity issues with the block grading system were raised internally by faculty and graduate student instructors, identified in formal and informal student feedback, and discussed with an external reviewing committee. In the Summer of 2018, many factors encouraged change: new teaching faculty (the authors) were hired in 2017 and 2018, there was a department-wide culture of experimentation fostered by the chair, feedback from an external reviewing committee sparked a conversation within the Calculus Committee, and University administration was focused on enhancing large introductory courses. At this time, five members of the teaching faculty, including the authors of this article, drafted a letter to the Calculus Committee requesting the authority to make broad changes to the grading systems. The Calculus Committee approved the request, and the course coordinators began implementing changes immediately in the Fall of 2018. The Calculus Committee agreed to equalize benchmark grades across the Calculus sequence, raising many courses' previous benchmarks. The Committee also gave coordinators increased agency to adjust the structure of their grading and assessment, removing the requirement that all abide by the former block grading system. The global health crisis of 2020, which moved most instruction online and led to new grading policies across Duke's Trinity College of Arts and Sciences, has catalyzed further changes to the Calculus grading systems (Duke University, 2013). In the Fall of 2020, all Calculus I and Calculus II courses were graded mandatory Satisfactory/Unsatisfactory, providing an opportunity to try new forms of assessment and grading without concern for letter grade averages.

***Table 1***

*Overview of the seven Calculus courses at the University overseen by the Calculus Committee, listed in numerical order. During primary semesters, all but Calculus III-a are run as multi-section coordinated courses: the primary term and average course enrollment over all primary semesters between Fall 2008 and Spring 2020 are italicized. (The primary term for a particular course is the semester (Fall or Spring) in which that course has the larger enrollment and the greater number of sections. E.g., in the primary semester, Calculus II-a has 7-9 sections, in the off-semester Calculus II-a has 2 sections and is not coordinated in the same way). Each section is generally of size 20-30 students.*

| **Calculus I** | **Calculus II** | **Calculus III** |
|---|---|---|
| Calculus I-a (I-a.i and I-a.ii)<br>I-a.i: *Fall; average enrollment, 90*<br>I-a.ii: *Spring; average enrollment, 59*<br>Two semesters of differential Calculus and integration with integrated precalculus. Both semesters are required for Calculus I credit. | Calculus II-a<br>*Spring; average enrollment, 199*<br><br>Integration techniques and applications, and series. Intended for students who took Calculus I at Duke (either Calculus I-a.i/I-a.ii or Calculus I-b). | Calculus III-a<br>*Fall; average enrollment, 86*<br><br>Large single-section Multivariable Calculus and Linear Algebra for economics majors. |
| Calculus I-b<br>*Fall; average enrollment, 243*<br><br>One semester of differential Calculus and integration. | Calculus II-b<br>*Fall; average enrollment, 150*<br><br>Integration techniques and applications, and series. Intended for students with incoming Calculus I credits (5 on the AP Calculus AB exam or 4 on the Calculus BC exam). | Calculus III-b<br>*Fall; average enrollment, 377*<br><br>Multivariable Calculus. Taken primarily by students with planned majors in Physics and/or Engineering. |

*Students must have Calculus II transcript credit to enroll in the Calculus III courses; all other prerequisites are recommended, but not enforced, with guidelines based on prior coursework, standardized test scores, and a first-day placement assessment.*

The grading and assessment systems in Duke's coordinated Calculus courses are still evolving. The Department of Mathematics continues to seek student and instructor feedback as well as demographic, GPA, and retention data but has not yet undertaken any formal assessment. However, preliminary results indicate that moving away from "fair and consistent" towards inclusive and transparent is the right direction for equity.

## 2. *Background: The Block Grading System*

Calculus at Duke is primarily taught through multi-section, coordinated courses overseen by the Calculus Committee (see Table 1). Some course components are common to all sections. Each section is generally of size 20-30 students. The textbook, the daily syllabus with suggested homework problems, the final exam, and the approximate weighting system for course grades are set and overseen by the course coordinator, who typically teaches one or two sections of the course. Other aspects of coordinated Calculus courses are determined separately for each section. Lectures, class activities, and all assessments other than the final exam (e.g., midterm exams, homework assignments, and quizzes) are typically created by individual section instructors.

Historically, the Duke Mathematics Department and Calculus Committee prioritized "ensur[ing] that grades are assigned fairly and consistently" across sections, semesters, and courses in the Calculus sequence using a "block grading system" (Supervisor of First-year Instruction, 2018). Under this system, all students in a coordinated multi-section course took the same final exam. Then, a curve was set by the course coordinator and the Supervisor of First-year Instruction so that the average final exam letter grade remained approximately the same from semester to semester, adhering roughly to departmental benchmark averages. (Recall that a course benchmark is a departmental guideline for the average course letter grade.) Lower-numbered courses had lower benchmarks. Finally, each section instructor was expected to use the average final exam letter grade of their section as their section's average course letter grade, and, in consultation with the course coordinator, curve individual student letter grades to meet it.

A general sense of this grading system was conveyed to students through a common `Block Calculus Policies' document: "[I]nstructors use information from the block final exam to gauge their sections as a whole in comparison to other sections, creating a context in which to gauge individual students with respect to the grading standards of the Department of Mathematics" (Duke Block Calculus Course Policies, 2018). Hence instructors and students had to wait until the final exam to fully and accurately interpret performance on earlier assessments. The document explicitly stated that "Letter grades assigned on work graded before the final exam are assigned without the benefit of the context [of the common final], and thus should be viewed only as estimates" (Duke Block Calculus Course Policies, 2018).

Duke's Calculus Committee first implemented the block grading system of maintaining benchmark average grades many years before the authors joined the department. Anecdotally, these benchmarks were intended to prevent an upward drift in course grades and to keep the meaning of letter grades consistent so that grades could be used as indicators of "proficiency."

The desire to curb grade inflation and use grades as a mechanism to rank students for future employers is not uncommon in academia (Schneider & Hutt, 2014).

### 3. *Motivations for Change*

In the Summer of 2018, the Calculus Committee agreed to change the block grading system. The change was motivated by feedback and data from students, instructors, and an external reviewing committee, which highlighted the following issues of equity, transparency, and inclusion:

1. Equity: Lower-numbered courses had lower benchmark averages (Figure 1), potentially disproportionately disadvantaging members of minoritized communities in STEM (Figures 3-4) and perpetuating existing inequities among incoming students.
2. Equity: Students in different sections of the same course expressed having very different experiences with learning environments, instruction, and grades.
3. Transparency: Students could not know their precise standing in the class before the final exam, potentially adding stress to both students and instructors.
4. Inclusion: Individual students were evaluated in comparison to other students and other sections (a norm-referenced system) potentially increasing competition and decreasing the sense of community and belonging.
5. Inclusion: Students perceived average grades in Calculus courses as “too low,” impacting both their desire for and access to continued studies in STEM.

We delve into these issues and the Department’s ongoing efforts to address them in the next section.

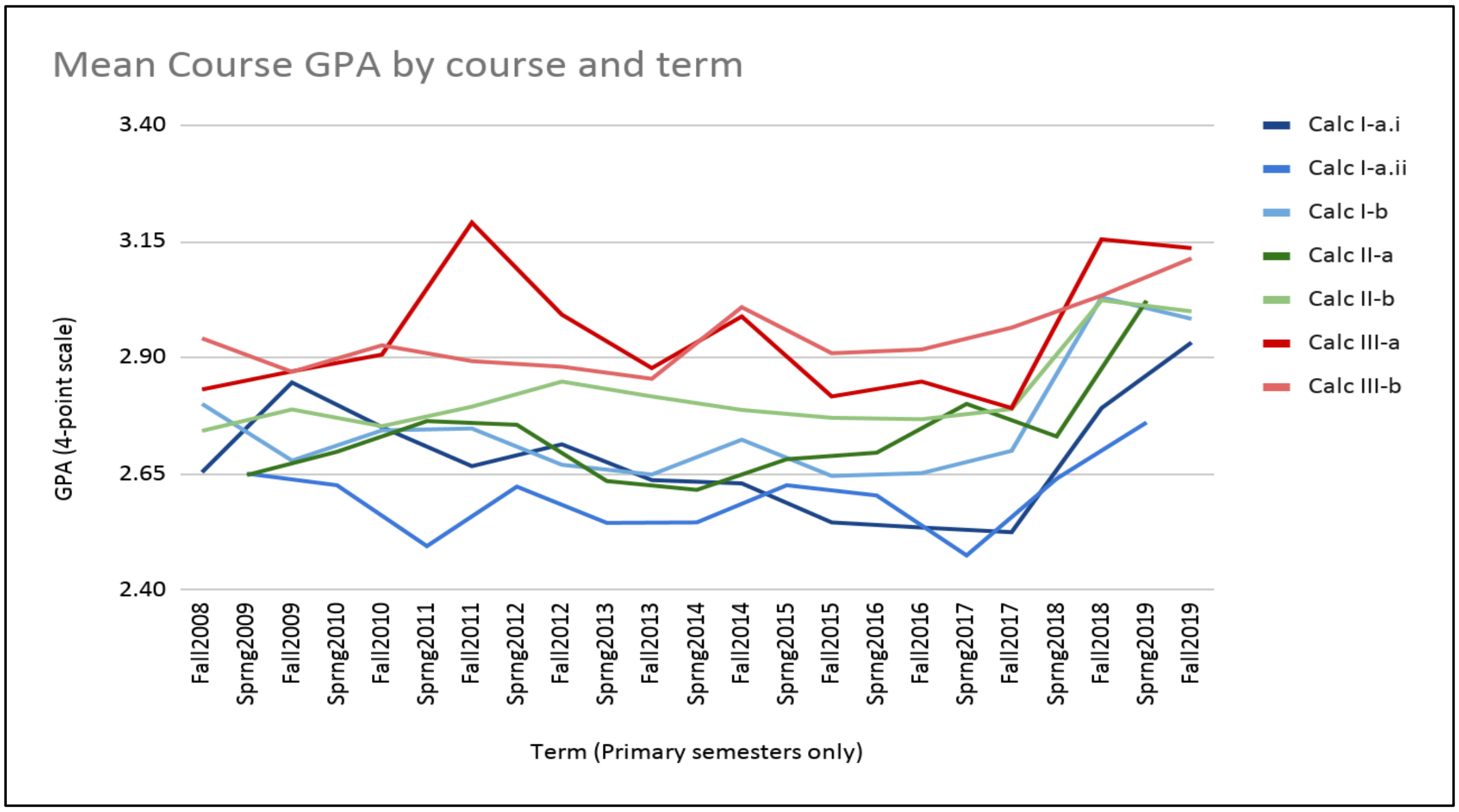

***Figure 1:*** *Mean course GPA in Calculus courses during primary semesters.*

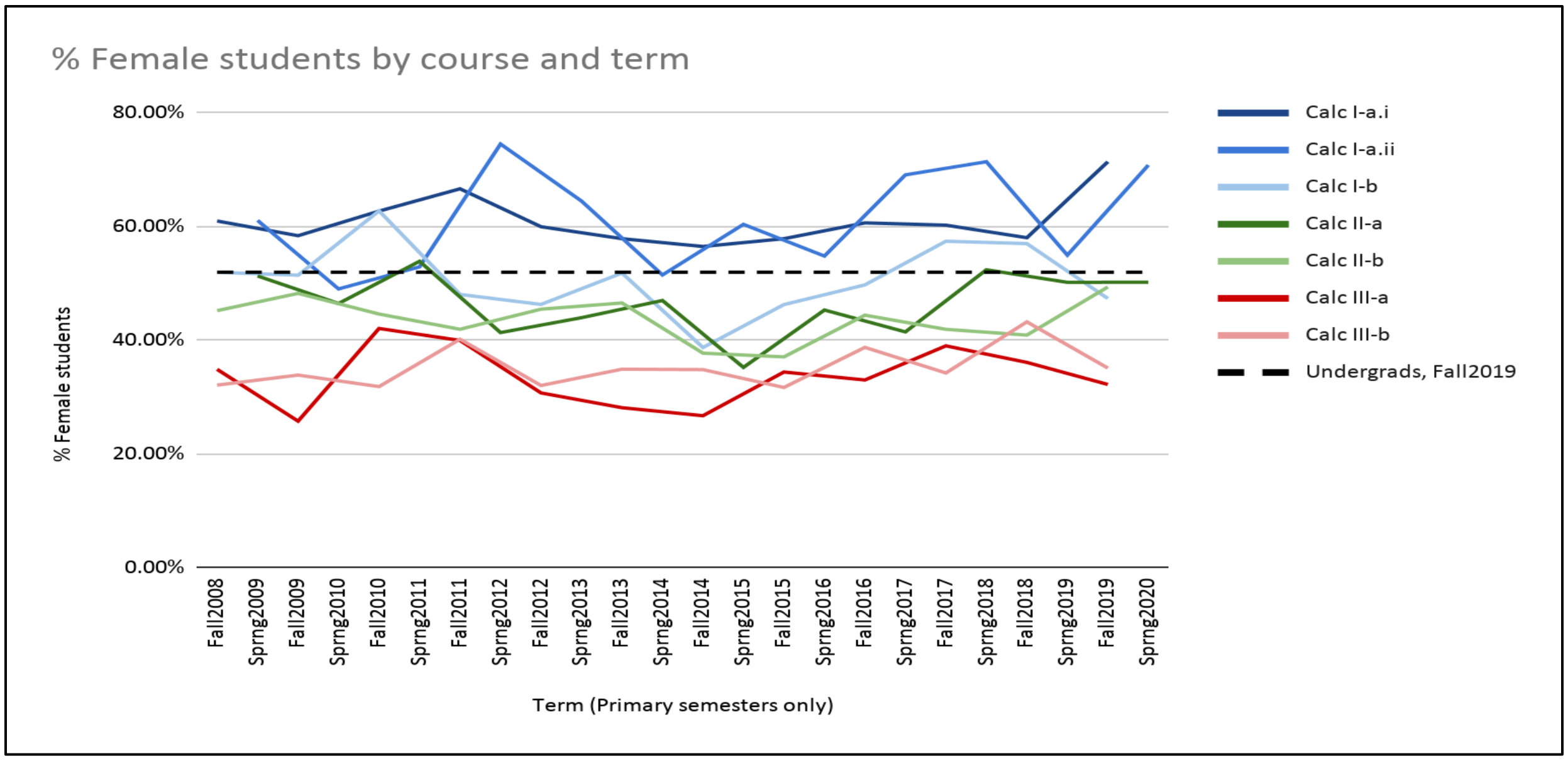

***Figure 2:*** *Percentage of self-identified female Calculus students during primary semesters*[12]

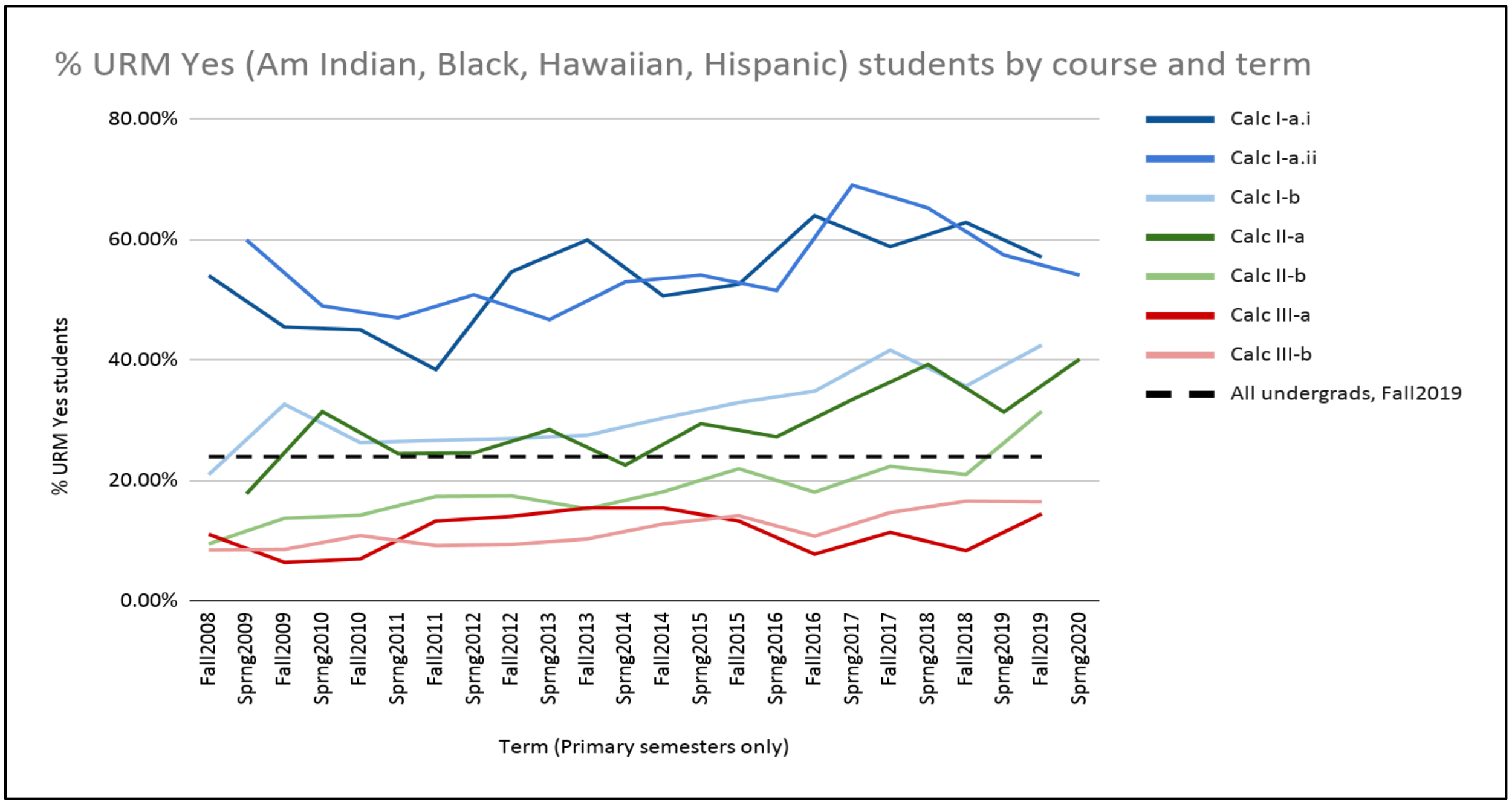

***Figure 3:*** *Percentage of self-identified "URM Yes" (Am. Indian, Black, Hawaiian, Hispanic) Calculus students during primary semesters*[13].

[12] Student sex is self-identified on applications to the University and can be updated through the Registrar. Current categories are `Female' and `Male', though choices will be more fluid in the future.

[13] Student ethnicity is self-identified on applications to the University and can be updated through the Registrar. The term URM, which stands for Under-Represented Minority, is that provided by the University, which categorizes students as `URM Yes', `URM No', or `URM Unknown.'

## 4. *Changes and Evidence*

Key changes to the systems of grading and assessment in coordinated Calculus courses since the Summer of 2018 are summarized in Table 2. These changes, conceived and implemented by course coordinators under the auspices of the Department of Mathematics and the Calculus Committee, are primarily aimed at increasing grading transparency and fostering a more inclusive and equitable learning environment by decreasing competition and emphasizing a growth mindset.

***Table 2***
*Overview of changes to Calculus grading and assessment systems from Fall 2018-Fall 2020. More recent changes are described in Section 5.*

| **Calculus I** | **Calculus II** | **Calculus III** |
|---|---|---|
| *Fall 2018:* Benchmark mean course grades equalized across ALL calculus courses:<br>In particular, benchmarks *raised* in Calculus I and Calculus II. | | |
| Calculus I-a (I-a.i and I-a.ii)<br><br>Fall 2018<br>Leeway for section instructors in assigning course grades. A section's final exam average still impacted the section's course grade average, but less rigidly than in years prior.<br><br>Spring 2020<br>Test correction cycle for all midterm exams.<br><br>Fall 2020<br>Mandatory no-letter grades: (Un)/Satisfactory only.<br>Mastery grading, no exams, resubmission cycle. | Calculus II-a<br><br>Spring 2019<br>Common midterm exams administered across many sections of the same course, as recommended by the external reviewing committee report to improve the grading system.<br><br>Spring 2020<br>Common midterm exams administered across many sections of the course.<br>Test correction cycle for all midterm exams.<br>Formative assessments with multiple attempts. | Calculus III-a<br><br>Fall 2019<br>Weekly quizzes as formative assessments: unit exam score replaced quiz scores if higher.<br>Non-norm-referenced grading explicit: "Grading is done on an absolute, but adjustable scale, given below. There is no curve. This scale may slide down, but will not go up. (Most likely, it will not slide.)" |
| Calculus I-b<br><br>Fall 2019<br>Explicit formula given for mean section letter grade: "0.65(B)+0.35(section average final exam grade)". Explicit statement that the average letter grade in the course is a B.<br><br>Fall 2020<br>Mandatory no-letter grades: (Un)/Satisfactory only.<br>Mastery grading, no exams, /retake cycle. | Calculus II-b<br><br>Fall 2018<br>Explicit weights given for all assessment categories, including 35% for final exam. All language from prior years describing the weighting system as non-binding or approximate removed.<br><br>Fall 2020<br>Mandatory no-letter grades: (Un)/Satisfactory only.<br>No exams. | Calculus III-b<br><br>Fall 2018<br>Each section instructor curved each midterm exam and announced the precise corresponding letter grades to be used in the computation of total course grades. |

The process of change is ongoing and incremental. For example, changing benchmark grades has improved student outcomes (see Figure 4 for a comparison of DFW rates), but equity and inclusion go beyond achievement, and students' sense of community, identity, and self-efficacy are all impacted by grading "on a curve" (Hughes et al., 2014, Seymour & Hewitt, 1997; Thiry et al., 2019). In striving for equity and inclusivity, the authors hope to eventually depart from norm-referenced grading entirely.

Certain structures in Duke's Mathematics Department carry inertia and encourage incremental change. For example, course coordination and Calculus Committee oversight support instructors and increase stability and coordination has been shown to enhance student success (Rasmussen et al., 2019), but also contribute to a department culture that can sometimes slow reforms.

In particular, the block grading system, with its goals of fairness and "consistency of meaning," was so entrenched that deviations had to be gradual (Supervisor of First-year Instruction, 2018). For example, multiple experienced graduate student instructors were immediately critical of the new systems. Throughout the Calculus II-a classes taught during the Spring of 2019, several graduate student instructors voiced concerns that grades might not be fair across sections, and/or the awarded grades were too high and included too many B's. Although the benchmarks in the block grading system were originally intended as flexible guidelines (see the variation in grades in Figure 1), they tended to be interpreted far more strictly by instructors with less experience and/or seniority. One graduate student commented:

> *As a first-time instructor, I appreciated the structure of the [block] grading system. In the end, I knew my grades would be exactly as generous or restrictive as the department wanted, so pleasing/disappointing the course coordinator or others above me was not something I had to worry about.*

While changes are sometimes slow, they are careful. The Calculus Committee provided advice and oversight for course coordinators who were trying new systems. The Department generally encourages innovation, and we anticipate continued reforms in the future.

To examine the impact of changes in grading and assessment in coordinated Calculus courses, we reviewed mid-semester surveys, end-of-course evaluations, and DFW rates. We would like to emphasize that our evidence of improvement is not the result of a carefully designed research study, and we do not attempt to identify which specific changes had an impact. Instead, we present post-hoc arguments that policy shifts largely improved student experiences by focusing on equity, transparency, and inclusion.

### 4.1. *Equity*

As a department, we are working towards building equitable grading and assessment systems in which every Calculus student has the support and opportunities they need to succeed. While equality is not equity, a first step toward this goal is addressing the inequalities inherent in the original block grading system. The changes described below, of raising benchmark average course grades in lower-numbered Calculus courses to match those in their higher-numbered counterparts, and of addressing wide variation in raw midterm exam grades, appear to be a step in the right direction.

Historically, lower-numbered Calculus courses had lower average grades (Figure 1). These same courses had higher proportions of students who self-identify as female and/or American Indian, Black, Hawaiian, or Hispanic (Figures 3-4). By using average grade benchmarks that were lower in lower-numbered courses, the Department disproportionately assigned lower grades to students from minoritized groups. As Calculus placement was based primarily on standardized test scores and course credit from AP and IB programs, this perpetuated existing inequities among incoming students at Duke by, on average, more severely impacting the GPAs of students entering the sequence at lower-numbered courses. Since course grade benchmarks were equalized in the Fall of 2018, the discrepancy in average course grades between lower- and higher-numbered Calculus courses appears to be decreasing (Figure 1). Further, both DFW and W rates dropped across the Calculus sequence, with the steepest drops in the lowest-numbered courses (Figure 4).

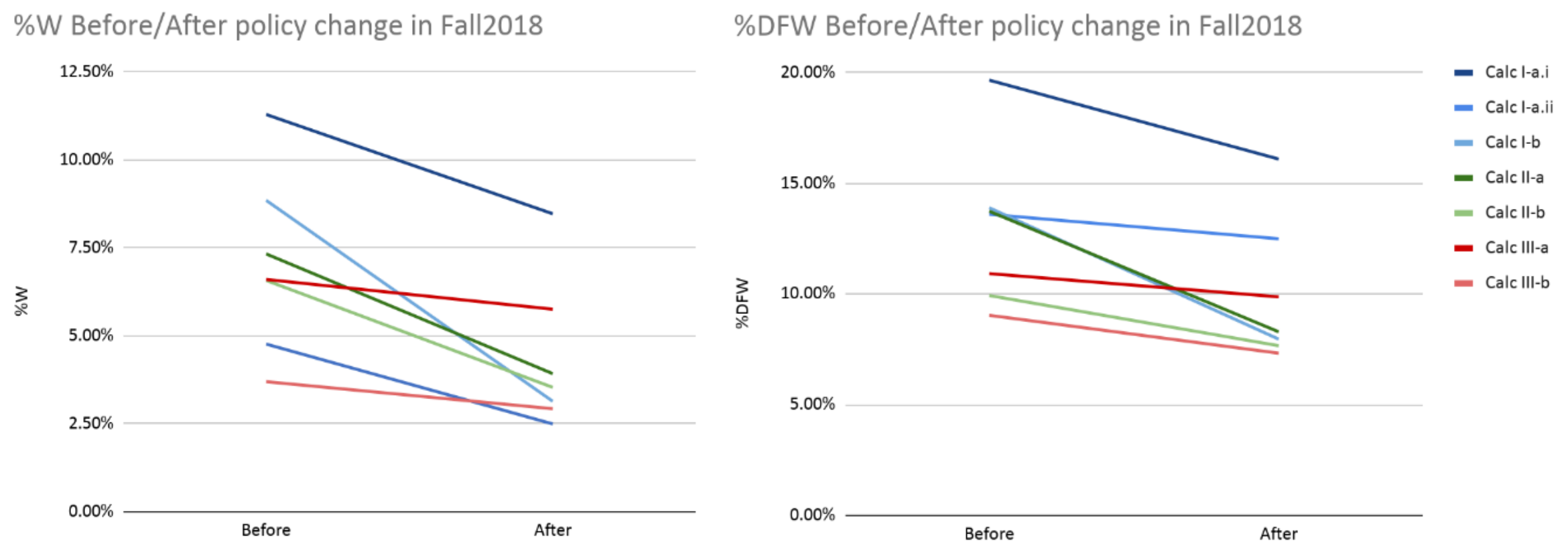

***Figure 4:*** *DFW and W rates in Calculus courses before and after Fall 2018. Displayed are `weighted averages': sum(enrollment*rate)/sum(enrollment), over Primary semesters only.*

Before the Fall of 2018, many student end-of-course evaluations of teaching referenced the “unfairness” of students’ varying experiences with grades within different sections of the same course, as well as with low raw grades in general. For example, raw section averages on the first midterm in Calc I-b ranged from 63.84 to 84.17. Again, equality is not equity, but perceived inequality can be demoralizing for students, as can low scores. Indeed, “Discouragement and loss of confidence due to low grades in early years [of STEM courses]” is an important factor in student decisions to leave a STEM major, particularly for students from historically minoritized groups in STEM (Thiry et al., 2019, p.105, 107). Two types of changes have been implemented since Fall 2018 to address these issues of raw grades and their variation between sections. First, some course grading systems have decreased the dependence of course grades on a section’s final exam performance (see details in Table 2). The previous heavy dependence incentivized many instructors, particularly inexperienced ones, to write challenging midterm exams with low raw scores. Second, some courses have tried to increase similarity in midterm assessments across sections. For example, common midterms were administered in the Spring 2019 Calculus II-b classes. End-of-course evaluations still contained negative comments about a system of comparison in which students with different instructors and classroom experiences were held to the same standard, so common assessments may not be enough to create an equitable experience.

## 4.2. *Transparency*

There is extensive literature on the value of transparency in instructor expectations, particularly on the benefits of using and sharing rubrics (Jonsson & Svingby, 2007; Reddy & Andrade, 2010). Before the grading system changed in the Fall of 2018, students repeatedly expressed frustration on end-of-course evaluations over their inability to know or predict their overall course grades before the final exam. The report by the external reviewing committee confirms that uncertainty about grading was a source of stress for students. Below, a first-time graduate student instructor reflects on teaching and grading Calculus I:

> *The way the grading was set up meant that we could not give our students definite answers about their grades until the end of the semester. Although I appreciate the need to be fair across sections, the lack of clarity translated to increased stress for the students. I felt dishonest when answering their questions vaguely, and felt some pressure to teach "to the test" since their grade distribution hinged on their final test performance. Students told me they are warned even before coming to campus by fellow undergrads that the final is one of the scariest tests they will ever take at Duke. This combination of mysterious grades and a scary final made it incredibly difficult to teach for the first time, help students relax, find joy in calculus and make the jump from high school math to college math. I am also worried that this set up was especially toxic for my students who were underrepresented in STEM.*

Beginning in the Fall of 2018, four different models were implemented to increase transparency in the computation of course grades (see Table 2 for details). Common features of these models included explicit and binding weighing systems and a decrease in the impact of the final exam. In the mid-semester Department survey administered in Fall 2019 and Spring 2020[14], students were asked, “How would you describe your understanding of the course grading system?” Across all Calculus courses the percentage of responding students who reported either “Somewhat Understand[ing]” or “Understand[ing]” was above 85%. Unfortunately, we have no data before grading reforms but believe this is a significant improvement. Partial evidence of this effect comes from analyzing end-of-course evaluations from the past five years. Before policy changes, there were multiple student comments on end-of-course evaluations requesting that the Department disclose grade distributions, but the number of this type of comment decreased dramatically after grading reforms began in the Fall of 2018. In the Spring of 2019, the number of complaints about uncertainty regarding course standing before the final dropped to zero in Calculus II-a.

## 4.3. *Inclusion*

The Calculus sequence at Duke has a reputation on campus as a ruiner of GPAs. As the external reviewing committee writes in their report, “Students refer to math at Duke as ‘Duke math’ because they see the math courses as unnecessarily difficult. Unfortunately, students are working really hard and still getting D’s, which ends up shaping many students’ career trajectories and resulting in students choosing paths based on the number of math requirements” (Bressoud et al., 2019). As mentioned above, low grades are a particularly important reason why men of color and women of all races leave STEM. "Men of color and women of all races and ethnicities were both more likely to report that grade-related problems had contributed to their switching decisions (69% and 67%, respectively) than were men overall (51%)" (Thiry et al., 2019, p. 286), and Kopparla (2019) suggests that DFW rates in introductory college math courses play a particularly important role in STEM retention (2019). Thus, the increase in average course grades (Figure 1) and decrease in DFW rates (Figure 4) after the policy changes in the Fall of 2018 indicate that the Math Department is making strides towards improved retention of these student groups.

---

[14] The Spring 2020 mid-semester survey was administered prior to Duke’s transition to remote learning.

Further evidence of improvement in Calculus students' sense of inclusion comes from the mid-semester Department survey. Since the Spring of 2019, students have been asked to rate the course climate on a 5-point scale from "Excluding and Hostile" to "Including and Friendly." Survey results show an increasing trend in the average rating across Calculus courses over the three semesters after Fall 2018[15], when grading policy changes began. Interestingly, among Calculus courses, the highest ratings tended to correspond to courses with a single section or two sections run by the same instructor; that is, the climate is viewed more positively in courses that do not compare grades across many sections.

The correspondence on midterm survey results between students' increased sense of inclusion and decreased sense of comparison corroborates other evidence that norm-referenced grading contributes to an exclusive climate. In a norm-referenced system, grades may not reflect a student's understanding of the course material or sense of their progress (Raymond, 2013; Sadler, 2009; Thiry et al., 2019). Norm-referenced grading also contributes to a competitive learning environment and culture, another important factor cited by students leaving STEM, particularly those who come from historically minoritized groups (Hughes et al., 2014; Seymour & Hewitt, 1997; Thiry et al., 2019). Strong peer support is particularly important for retaining historically minoritized communities in STEM, and reducing competition is an important step toward better inclusion (Palmer et al., 2011). The Department has begun using non-norm-based grading in some Calculus courses (Table 2), and the authors see great promise in continuing in this direction. Several faculty members, including the authors, have become particularly excited by what is broadly referred to as mastery grading, a criterion-, rather than norm-referenced system in which students are scored based on which and/or how many specific learning objectives they have mastered rather than on `total points earned', and are given multiple opportunities to demonstrate this mastery (Cilli-Turner et al., 2020; Iamarino, 2014; Owens, 2015; Tomlinson & McTighe, 2006). Calculus I-a.i and Calculus I-b piloted mastery-based systems in the Fall of 2020. By encouraging students to review mistakes and re-submit work, such systems promote the growth mindset that students can improve with time and effort, which contributes to resilience and success in STEM (Dweck, 2015). Even more, a growth mindset can protect underrepresented demographic groups against stereotype threat and can counterbalance the effect of lower socioeconomic status on achievement (Aronson et al., 2002; Claro et al., 2016).

[15] The Spring 2020 mid-semester survey was administered prior to Duke's transition to remote learning.

### 4.4. *Broad Evidence*

In a review of end-of-course evaluations for multi-section Calculus courses from Fall 2015 to Spring 2020, we saw a drop in the number of comments mentioning “grade”, “grading”, or “curve” after the implementation of grade policy changes in Fall 2018. The end-of-course evaluations painted a picture of generally more satisfied students, although students continued to note ways in which the grading system could be improved. As noted above, policy changes also resulted in an increase in average course grades (Figure 1) and a decrease in DFW rates (Figure 4). Finally, we observed early evidence of a culture shift around Calculus at Duke. Historically, students were cautioned by peers and advisors to take the fewest number of Calculus courses possible regardless of interest or placement, but the authors heard anecdotally from administrators and faculty in other departments that this narrative may be changing.

## 5. *Looking Forward*

We have taken steps toward identifying equity, transparency, and inclusion problems with the old block grading system and implementing new policies in its stead. These changes are a work in progress. Further potential improvements include eliminating curves, increasing transparency in expectations and grade computation, decreasing the disparity in experience between sections, and making effort a larger factor in course grades. We are particularly enthusiastic about moving away from high-stakes assessment and norm-referenced grading and towards formative assessment and mastery grading.

Across the sequence, we have been dealing with the same pandemic-related issues that affect all of academia including student mental health, differing incoming math experience from primarily remote high school learning, and changing policies for social distancing and in-person learning. As a result, all calculus courses were fully remote from the Fall of 2020 to the Spring of 2021. Calculus courses then primarily returned to in-person instruction in the Fall of 2021, with some remote/hybrid components to accommodate varying circumstances. All Calculus I and Calculus II courses in the Fall of 2020 and the Spring of 2021 were mandatorily graded as Satisfactory/Unsatisfactory, with no letter grades and therefore no curve. These courses all eliminated both midterm and final exams, and Calculus I-a.i and Calculus I-b experimented with a mastery-based assessment structure. Since the initial submission of this chapter, the Duke Mathematics Department has continued implementing changes to the grading systems in the Calculus sequence. Although these changes were not included in this chapter’s data, we detail the new changes here as part of Duke’s continued progress. After the return to in-person instruction, Calculus I-a grades have continued to be mastery-based and Satisfactory/Unsatisfactory grades are mandatory. In Fall 2021, Calculus I-b was structured as two large lectures taught by the same faculty member in an effort to make the student experience more consistent. In Fall 2022, Calculus I-b will be formatted as one large lecture section (with approximately 150 students) taught by the course coordinator and four small sections (with approximately 35 students) taught by graduate student instructors. The course description for Calculus I-b has been correspondingly updated to provide students guidance on choosing the class size that is best for them. Calculus II-a and Calculus II-b have continued efforts to make their grading systems transparent to students and to make in-class exams more similar across different sections. In Fall 2021, a non-primary semester, Calculus II-a experimented with offering students two grading system alternatives. Students could choose between a more traditional grading system in which high-stakes exams were heavily weighted or an alternate system in which students were asked to submit class notes and student exams were weighted less heavily.

There are multiple confounding factors, but the Department looks forward to investigating the impact of these changes on students' mindsets, sense of belonging, and identity. We will continue adapting the Calculus grading and assessment systems to make our learning environments more inclusive and equitable. We will also look for new ways to evaluate and improve our impact, particularly on the experiences of students from historically minoritized groups in STEM. The Department has joined a team of universities across the nation as part of an NSF-funded project to investigate and implement equity-focused changes in introductory math classes through networked improvement committees. Finally, the authors, as well as the Department at large, will continue active efforts outside the classroom to increase access to and inclusion in mathematics within Duke and the larger community. We hope that with this work, Calculus courses may act less as gatekeepers and more as launching pads.

# 2 Using a Multi-Tiered and Data-Driven Approach to Support Student Access and Success in Introductory Mathematics at a Large, Hispanic-Serving Institution

**Babette M. Benken**

*California State University – Long Beach*

Babette M. Benken, Ph.D., is the Richard D. Green Professor of Mathematics Education in the Department of Mathematics and Statistics and the Director of Graduate Studies in the College of Natural Sciences and Mathematics at the California State University, Long Beach. Her primary research areas include mathematics and STEM education, teacher learning, and graduate education.

**Krzysztof Slowinski**

*California State University – Long Beach*

Krzysztof Slowinski, Ph.D., is the Associate Dean for Academic Programs, Evaluation, and Advising in the College of Natural Sciences and Mathematics at the California State University, Long Beach.

*Abstract: This case study highlights details of, outcomes, barriers and recommendations stemming from multiple large-scale initiatives instituted by the Department of Mathematics and Statistics (in partnership with college leadership) at the California State University, Long Beach. These initiatives include multiple measures for course placement, intrusive advising with an early alert program, redesigned courses, and ongoing analysis of data. Taken together, the impact of these four initiatives over the past few years has expanded opportunities for students to successfully complete introductory calculus pathway courses for STEM majors, suggesting that these types of college-wide collaborative reforms supported by rigorous data analysis can lead to student achievement and access. Furthermore, it has done so without increasing time to degree. And, perhaps most notably, the increase in passing rates is beginning to support the decrease in achievement gaps.*

## 1. *Introduction*

California State University, Long Beach (CSULB) is a large, public Hispanic serving institution (~39,400 students, 46% of undergraduates are Hispanic/Latino). CSULB is part of the nation's largest public university system (CSU), serving more than 482,000 students across 23 campuses. Nearly half of the state's bachelor's degrees are awarded by the CSU. The CSU is also one of the most ethnically and racially diverse university systems in the U. S., with the student body of CSULB being no exception. For example, in 2021 CSULB was ranked 4th in the nation for awarding bachelor's degrees to underrepresented minorities. The campus' diverse student body has also experienced growing enrollments of both underrepresented and underserved students–over 84% identify as students of color, approximately 60% are female, over 55% of undergraduates identify as first generation, and approximately 75% of undergraduates receive some form of financial aid. Recently, CSULB was ranked No. 1 in the nation by Education Reform Now for its "social mobility impact"; this index highlights 614 four-year colleges and universities where students with verified financial need (measured by Pell Grant) are more likely to graduate than drop out and where federal loan repayment and default rates are better than the average for four-year institutions (Murphy, 2020).

Within the Department of Mathematics and Statistics, in spite of ongoing program improvement (e.g., Benken et al., 2015), analyses of enhanced institutional data revealed that too many students were still struggling in multiple introductory courses (e.g., precalculus, business calculus) and achievement gaps existed; this was particularly true for students deemed not yet ready for university-level mathematics. Our findings were similar to existing research that has shown that traditional postsecondary mathematics can be a primary barrier to degree completion and equitable outcomes for students (U.S. Department of Education, 2017). At CSULB this problem extended beyond mathematics and involved course completion in many STEM disciplines. Figure 1 shows overall DFW rates from 2009 to 2019 for all courses (total>15,000 students) offered in the College of Natural Sciences and Mathematics (CNSM) at CSULB. For this data, "underrepresented minority" (URM) includes "African American, American Indian and Hispanic/Latino/Latina" students; "NURM" signifies those students who are not identified as "underrepresented minority" (e.g., white, Asian) It is clear that while we have been making systematic progress in decreasing DWF rates across all courses in CNSM, opportunity gaps persist and have been a matter of major concern.

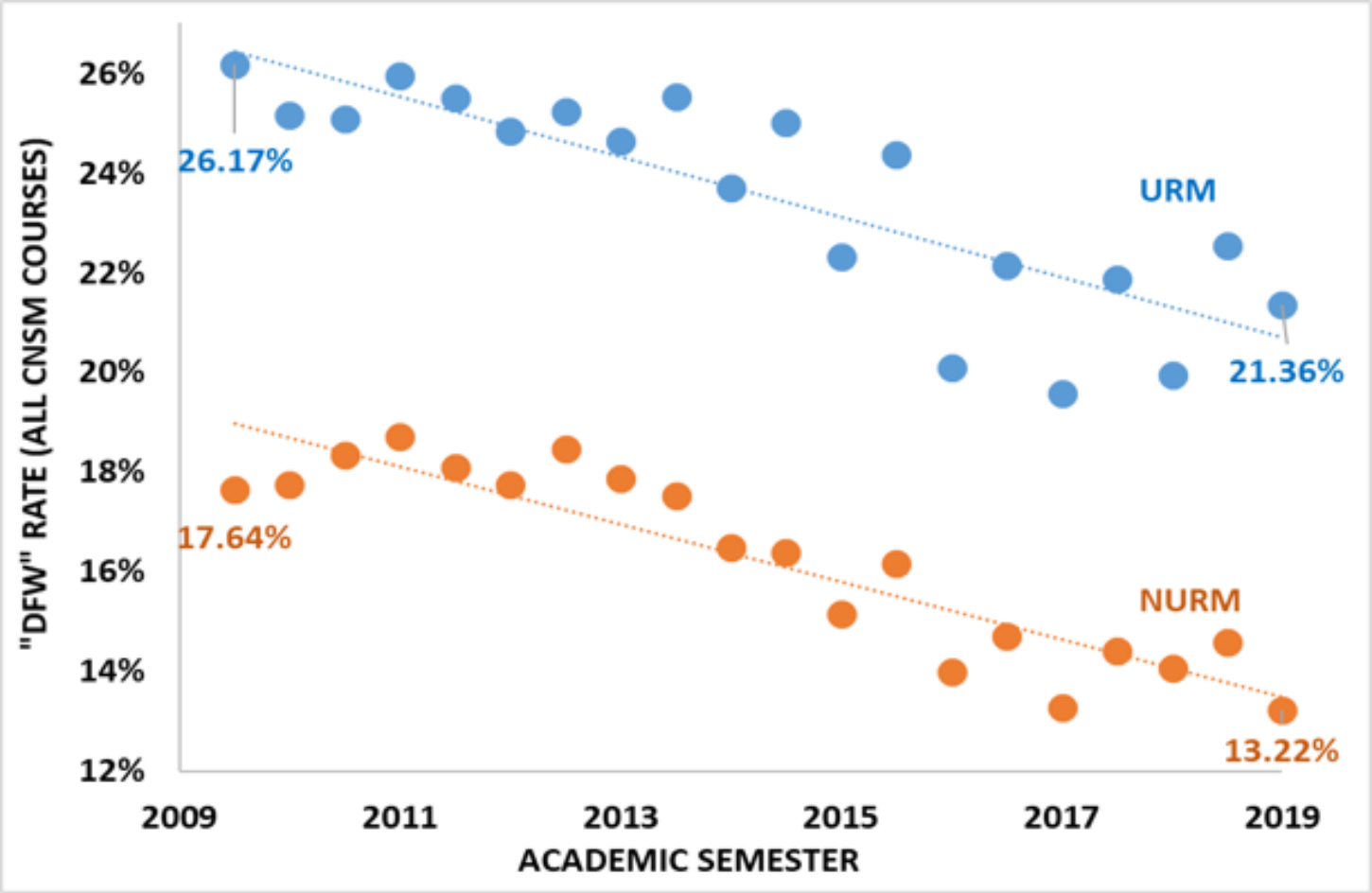

*Figure 1: DFW Rate (all CNSM Courses) by Year and Minority Status*

Furthermore, in 2017, the CSU system issued a mandated Executive Order that called for elimination of non-credit bearing developmental education beginning fall 2018 to improve student success, increase degree completion, and close persistent equity gaps (CSU System, 2017). In 2011, for example, over 30% of first-time freshmen in the CSU required remediation in mathematics (CSU System, 2012, 2019). Furthermore, the majority of these students were from underrepresented subgroups (e.g., 59% of African American first-time freshmen fell into this category). Although intended to help support student success in mathematics, remediation can have negative consequences for students. In some cases, it can also become a barrier for future academic achievement (Noel-Levitz & CAEL, 2006), thereby preventing access to higher education and career choice, including careers in STEM.

In response to these challenges, the Department, in partnership with college leadership, instituted multiple initiatives. These initiatives are outlined below. Outcomes, barriers, and recommendations are provided.

## 2. *Program Context*

The Department of Mathematics and Statistics is housed in CNSM and comprises 42 tenured/tenure-track faculty members across four disciplinary subgroups (pure, applied, statistics, math education). Of the faculty, approximately 29% are female, 52% are white, 17% are Asian American, 5% are Hispanic/Latina/o, and 2% are African American. It additionally employs approximately 35 lecturers and 30 teaching assistants (as instructor of record) each year. It offers a wide range of mathematics, statistics, and mathematics education courses at both the undergraduate and graduate levels (4 undergraduate majors, 3 minors, 4 master's programs). In fall 2019, the Department had approximately 455 majors (declared and pre) and 160 graduate students. Furthermore, it served numerous other majors through general education (GE) courses, including those for teacher content preparation (~6,700 students, 1,680 FTES in total). The Department plays a critical role in providing quantitative reasoning / mathematics courses for all freshmen. For example, in fall 2019, the Department offered GE foundation mathematics courses for 4,047 first time freshmen, including 679 freshmen in pre-calculus algebra, 778 freshmen in calculus 1-2 courses, and 349 freshmen in business calculus.

## 3. *Initiatives*

### 3.1. *Initiative #1–Multiple Measures for Course Placement*

To expand access to STEM pathways and increase achievement in both precalculus and calculus courses (and thus subsequent mathematics courses), CSULB instituted multiple placement measures. Historically, narrow placement practices (i.e., high stakes tests with standardized cutoff scores) at postsecondary institutions have shown to be ineffective (Liston & Getz, 2019). More recently, national data indicates that using multiple measures can be effective for both community colleges and universities (e.g., Barnett et al., 2018; U.S. Department of Education, 2017). At CSULB, an examination of data for entering freshmen revealed that SAT/ACT scores and high school GPA predicted success in Precalculus Algebra and students with low scores would need expanded support. We now use additional tools (e.g., ALEKS) and offer summer support (Early Start) for skill-building prior to beginning the course. Similarly, for calculus, we utilized findings from longitudinal data to identify a combination of SAT/ACT and AP Calculus scores for placement. We provide students with early information regarding placement schemes (http://web.csulb.edu/colleges/cnsm/advising/incoming-freshmen/math-placement.html )) and encourage them to brush up on their math skills during summer (https://web.csulb.edu/programs/hsi-stem/students/step-into-stem/). These programs intertwine math instruction with programming focused on growth mindset, metacognition, STEM identity and sense of belonging, and intensive peer mentoring (elaborated later).

Implementing these multiple measures with detailed pathways that outline recommendations for both course placement and parallel support (based on the combination of measures) has resulted in increased passing rates in all introductory courses. For example, from 2017 to 2019 there was a nearly 10% gain in passing rate (ABC) for Precalculus Algebra (Table 2).

*Table 2*

*Precalculus Algebra GPA & DFWU Rates by Semester and Minority Status (URM at CSULB = African American, American Indian and Hispanic/Latino/Latina)*

| MATH 113 Fall Semesters | NON-URM Avg. GPA | URM Avg. GPA | GPA Gap | DFWU Rate |
|---|---|---|---|---|
| **2019** | 2.56 | 2.28 | 0.28 | 18% |
| **2018** | 2.65 | 2.17 | 0.48 | 22% |
| **2017** | 2.38 | 1.96 | 0.42 | 27% |

Of particular importance is the performance of students in precalculus who, in the past, would have been deemed "remedial" or not fully algebra-ready based solely on MATH SAT score. These students (MATH SAT <570) are now placed in Precalculus Algebra via ALEKS PPL (with or without participation in Early Start or other summer bridge programs) and are provided additional fall semester support via dedicated supplemental instruction that is led by CSULB students and learning communities (including caseload-based peer mentoring), https://www.csulb.edu/college-of-natural-sciences-and-mathematics/freshmen-scholars-learning-community. Even this subset of students had a passing rate of 70% in 2019 and thus satisfied the campus' general education (GE) requirement and could also still pursue a STEM pathway. Furthermore, students with SAT scores within the bottom 10% of the class (SAT<510) still demonstrated a passing rate of 65%. In the past, these students would have required 1-2 semesters of non-credit bearing remedial math courses before even attempting Precalculus Algebra, which historically increased their time to degree and sometimes dissuaded them from continuing to pursue one.

This increase in passing rates has expanded access to STEM pathways and supported students' growth mindset and timely graduation. Since the majority of our students are underrepresented (majority are Hispanic/Latina/o), this implies that more of these students will have the opportunity to pursue STEM majors that require calculus and will not be deterred or delayed by retaking prerequisite courses. Moreover, these gains in Precalculus Algebra did not result from a reduction in rigor or under preparation for calculus courses. In fact, about 1/3 of the freshmen enrolled in Precalculus Algebra in their first semester attempted Calculus I the following semester and the completion rate for this group was 75% (vs 76% for all other students).

Utilizing multiple, and often complex, placement schema takes time for both design and implementation. Thus, we recommend training for both faculty and staff advisors, setting aside time in summer to analyze incoming students, and flexibility for last-minute modifications to course offerings (we usually have numerous course changes late in summer). Furthermore, the targeted support based on this analysis should include components that instill a growth mindset and sense of belonging and connection to the college/program (e.g., peer mentoring).

### 3.2. *Initiative #2–Intentional Advising*

In order to provide needed support for student success in all mathematics and science GE courses, we have instituted an intentional approach to advising, i.e., *intrusive advising* (Earl, 1988; Varney, 2007). In CNSM, this approach includes an Early Alert program (Bernacki et al., 2020) and additional mandatory academic advising. Early Alert is an academic advising and student support intervention plan designed to identify students experiencing academic difficulties mid-semester through faculty referral following college-level data analyses (https://web.csulb.edu/colleges/cnsm/advising/early-alert.html). Students are encouraged to meet with professors and advisors (department and/or college) to identify causes and strategize interventions to help them achieve satisfactory progress. Mandatory advising for students not meeting milestones includes assessment of learning skills and study habits by the college, and referrals to campus resources (e.g., tutoring, health services). Our college has since developed a *How to Succeed* website that provides an extensive collection of resources for students (e.g.., videos on growth mindset and how to overcome failure, links to campus resources, systems for locating and setting up appointments with advisors – https://web.csulb.edu/colleges/cnsm/advising/how-to-succeed.html, https://www2.calstate.edu/csu-system/news/Pages/growth-mindset-success-2019.aspx).

We have found accurate course placement combined with intentional advising and additional course-specific student support to be effective in improving overall retention and thus access to STEM majors. For example, from 2012 to 2016, the 3-year retention rates within the STEM disciplines for science and mathematics pre-major freshmen improved from 20% to 35% for URM students and from 35% to 50% for students who are not designated URM at CSULB. While these gains appear to be associated with changes in freshmen placement in calculus and advising approaches introduced during the 2013-2014 academic year, the large opportunity gap remains a major concern for access and inclusion. In this context, our most recent success with the Precalculus Algebra course seems promising in our efforts to eliminate opportunity gaps in STEM retention.

During interviews conducted in fall 2017 with students who switched from their initial STEM major, we found that students often feel isolated and do not perceive a sense of community and belonging in introductory courses. At a school that has many commuter students, it is sometimes more difficult for classmates to meet. This lack of a sense of community is especially detrimental for Latina/o and other underrepresented students, who are more likely to feel isolated (Nettles & Millet, 2006) and/or develop imposter syndrome (Ewing et al., 1996). Part of this advising effort has thus included developing a peer mentoring program and learning communities (e.g., precalculus and general chemistry students having the same peers and teachers across both courses). The learning communities target gate-keeper courses that have historically prevented students from pursuing STEM majors (Song et al., 2020). In these communities, the instructors follow the students, who register in reserved and paired course sections. Since implementing these initiatives students report on our annual survey that they are happier and feel a sense of belonging to the program. For example, in 2019, 88% of students completing the survey who face the greatest barriers to access (i.e., 1st generation, low income, and/or underrepresented) indicated that they anticipate their sense of belonging to grow over time between their first and second years based on their experiences in CNSM (number survey completers was 430).

Setting up systems for early and intentional advising that include ongoing analysis of student success data requires a team–one that is collaborative and includes faculty, advisors, and administrators. Furthermore, creating new programs and repositories of student resources takes time, particularly if this process requires generating new resources (e.g., processes, video) and/or securing university funding. We recommend recruiting a college-wide team and beginning with an analysis of current advising supports and passing rates for a subset of introductory courses with the lowest passing rates. We found this early work can lead organically to the development of essential resources (e.g., website, university funding) and programs (e.g., peer mentoring, learning communities) once key strategies are identified and have been shown to positively impact achievement.

### 3.3. *Initiative #3–Redesigned Courses*

CSULB has been actively engaged in the redesign of lower division GE mathematics courses since 2013. The Department (and campus community) had concerns about low passing rates, and thereby subsequent access to and retention in STEM, and time to degree. With recent CSU support following the aforementioned CSU Executive Order, in 2018 our faculty completed redesign of Precalculus Algebra, Business Calculus, Calculus I and II, and Survey of Calculus (life sciences). The redesign includes adaptive homework, intrusive advising, supplemental instruction, hybrid courses, learning communities, stretched courses, and growth mindset interventions.

For example, as discussed above, passing rates were low for Precalculus Algebra (Figure 2). Many of our students were taking 1-2 remedial courses, only to then subsequently struggle to pass or not be successful in Precalculus Algebra. This outcome resulted in many students needing to take yet another mathematics course to satisfy the GE requirement. Additionally, the overwhelming majority of those students who passed with a “C,” struggled to pass the subsequent course (e.g., Business Calculus, Calculus, Precalculus Trigonometry). In fact, for 2016 freshmen science/math cohort, only 14% of students who started in pre-baccalaureate mathematics courses persisted in STEM beyond 2 years. Likewise, only 35% of students receiving a “C” in their first semester of Precalculus Algebra persisted beyond 2 years (vs. 72% for those receiving an “A”).

Thus, for Precalculus Algebra (Math 113) we designed a stretched version that divides the content from the 1-semester course across two semesters (Math 112A & 112B). The first semester of the stretched version is also designated GE and thus if students decide not to continue in the STEM pathway, they will have at least satisfied the University’s quantitative reasoning requirement. For students deemed not yet ready for college-level mathematics, the Department also designed a co-requisite (1 unit, Math 92) course that provides additional time for students to collaboratively explore parent-course tasks and receive individualized instruction. This course is also open to, and often accessed by, students without this entry designation. What we did not anticipate is that many of the students who took the co-requisite course and intended to continue to the second half of Precalculus Algebra would request a similar co-requisite course for that second half. We designed the co-requisite for the second half of the stretched course and piloted one section of it in spring 2022. Advisors encouraged those students passing with "C" to pursue this support course, which provided an inclusive opportunity for a smaller learning community that provided differentiated 1-on-1 instruction. The course proved to support student success; 85% of students enrolled in the support course passed the parent course with C or better, while only 76% of students in the parent course passed overall.

Notable success of the complete redesign effort is in the calculus pathway with ABC rates for the first-time freshmen (2018) at 77% in Precalculus Algebra, 74% in Calculus I for Life Sciences, 77% in Calculus I for Physical Sciences, and 83% in Calculus II and III. These courses provide an effective and optimized calculus pathway. That is, students who enter Calculus I (based on ALEKS and/or SAT score) and Calculus II or III (based on appropriate AP Calculus AB or BC score of 4-5, respectively) stand the same chance of passing the course (above 80%) as students who enter the course by completing the CSULB prerequisite (i.e., Precalculus Algebra or Calculus I/II course). Furthermore, access has expanded and time to degree decreased. For example, in Fall 2017, 665 of CSULB first time freshmen required remediation, yet only 189 of these students earned GE math credit by the end of their freshmen year. By comparison, in Fall 2018, 704 CSULB freshmen were designated as "in need of additional support in math" (not identical but similar designation to previously "remedial" students) and 420 of these students completed GE math by the end of the first year. Hence, in the first year of implementation we have doubled GE credit attainment for this category of students. Central to this effort has been expanded coordination of courses (curriculum and structure) with targeted training for instructors and teaching assistants, as well as scheduling that allows teaching assistants to concentrate on only one course preparation.

### 3.4. *Initiative #4–Data Analysis*

Our significant ongoing effort in data analytics involves quantitative and qualitative methods and multiple sources of data to assess STEM students' experiences and inform course placement (e.g., class size and instructor data, major migration patterns, course completion statistics, passing rates in current and future courses, interviews with students who switched from initial STEM major, graduation rates). We are now expanding these efforts to include social-cognitive and affective traits in conjunction with admission, student success, and demographic data to identify optimal pathways to timely graduation while focusing on reduction of equity gaps and ensuring productive exploration. Examples of traits include: career commitment; student agency (e.g., grit); time & study environment; student/teacher responsibility; academic affect (e.g., test anxiety, perceptions of academic stress); academic opportunity (e.g., external commitments, finances); math confidence; sense of belonging; science identity; motivation). Based on the survey responses from 2015-18 science/math freshmen cohort students (n=350) we have assembled a data set that we used to develop a predictive STEM retention model for CSULB students.

This rigorous approach to analysis of STEM student success requires institutional support for time to assess appropriate data sources, collect and analyze data, and then determine relevant recommendations. At CSULB, we leveraged existing committees and research and scholarly funding opportunities to form college-wide teams (administrators, faculty/lecturers, advisors, students). These teams worked collaboratively to expand what department chairs and college staff could have accomplished alone. Given that we used university-level funding and initiatives it also provided a vehicle through which participants could garner recognition and an internal opportunity for scholarly productivity.

## 4. *Conclusion*

Taken together, the impact of these four initiatives over the past few years has expanded opportunities for students to successfully complete introductory calculus pathway courses and STEM majors, suggesting that these types of college-wide collaborative reforms supported by rigorous data analysis can lead to achievement and access. Furthermore, it has done so without increasing time to degree. And, perhaps most notably, the increase in passing rates is beginning to support the decrease in achievement gaps.

# 3 Increasing Inclusion in Large Enrollment, Uniform Math Courses: Instructor Training and Course Assessment

**Hanna Bennett**

*University of Michigan, Ann Arbor - Department of Mathematics*

Hanna Bennett is a Lecturer in the Math Department at the University of Michigan. As a Co-Director of the Introductory Program and course coordinator for the program, her duties include training and mentoring new instructors and finding creative ways to balance the needs of hundreds of students, dozens of instructors, and a very small number of course coordinators. She is particularly interested in Inquiry Based Learning and inclusive teaching.

**Susan J. Cheng**

*University of Michigan, Ann Arbor - Center for Research on Learning and Teaching*

Susan J. Cheng is an instructional consultant with the Foundational Course Initiative at the Center for Research on Learning and Teaching at the University of Michigan, where she specializes in assessment and analytics of large introductory courses. She is also an affiliate in the Department of Ecology and Evolutionary Biology. Her research examines student learning of environmental science, biology, and math, as well as the impacts of classroom climate on student learning.

**Paul Kessenich**

*University of Michigan, Ann Arbor - Department of Mathematics*

Paul Kessenich is a Lecturer in the Math Department at the University of Michigan, where he is a Co-Director of the Introductory Courses. His duties include coordinating large courses with many instructors, course design, and training and mentoring of instructors, including faculty and graduate students. Interests include instructor training (both in the classroom and more individualized tutoring), interactive and engaged teaching methods, and course design to promote equity and inclusion.

**Elaine Lande**

*University of Michigan, Ann Arbor - Comprehensive Studies Program*

Elaine Lande is a Lecturer in the Comprehensive Studies Program at the University of Michigan, where she teaches introductory mathematics courses. Her research has focused on teaching and learning mathematics at the tertiary level, including the difficulties and constraints instructors face, classroom interaction, and students' understanding of mathematical concepts. Her current scholarly interests focus on classroom and institutional factors that impact the inclusivity of mathematics.

**P. Gavin LaRose**

*University of Michigan, Ann Arbor - Department of Mathematics*

Gavin LaRose is the Karen Rhea Collegiate Lecturer in the University of Michigan Department of Mathematics and manager of instructional technology applications used in the Department. His research background is in applied mathematics, but his current scholarly focus is on disciplinary educational research and issues including instructor training, technology, and inclusive teaching. He has taught at both a small liberal arts university and at the University of Michigan, and his professional involvement has included work with the Michigan Section of the MAA, MAA Project NExT, the WeBWorK Project, and committees of the MAA, AMS, and AWM.

*Abstract: There are issues in building inclusion courses in large enrollment introductory courses taught at large universities that are unique to that venue: much of the instruction is done by generally unsupervised beginning teachers, assessment must be consistent between sections of*

*the course, and the course structure must be sustainable at scale. We describe work to address these issues, which focuses initially on our course preceding calculus, and which is aligned along two dimensions of change. The first dimension of change is improved support of instructors on issues of diversity, equity, and inclusion. We have integrated an inclusive teaching workshop into our existing training program for new instructors, and implemented complementary workshop sessions for instructors during the following semester. These efforts dovetail with more general efforts to increase the visibility and support of inclusive teaching in the department. The second dimension is the revision of the assessment in the course to include a significant mastery component that allows students greater ownership of their progress through the course and work in learning. This changes the assessment in the course from a model where 95% of students' grades were determined by three high-stakes exams, to a model where students' grades are determined by mastery assessments, comprehensive analysis of students' daily work, and two low-stakes exams. Preliminary evaluation of these changes indicates that they are a positive change to the course, producing effects that align with our objectives.*

---

> *"Our basic task [as teachers] is to determine what we mean by mastery of the subject and to search for the methods and materials which will enable the largest proportion of our students to attain such mastery." --Benjamin S. Bloom (1968)*

---

## 1. *Introduction*

In this case study, we describe our work on two tightly connected projects in our introductory mathematics courses that are designed to create a more inclusive learning environment for all students, including, and specifically motivated by, students from groups historically underrepresented in STEM. With these projects we seek to build our instructors' interest in and effective use of inclusive teaching practices; and to shift course assessment to emphasize mastery learning, helping promote a growth mindset in our students. This work is done in the context of a department in a large, public, historically white university[16] that supports effective instruction and inclusive teaching. Our current work is supported by departmental faculty and takes advantage of university-wide efforts on inclusive teaching; in the absence of this support, our work would be orders of magnitude more difficult.

## 2. *Motivation*

In the epigraph above, educational psychologist Benjamin Bloom brings out the importance of teaching so that all students can master a subject. This is at the heart of the work we describe in this case study, which is focused on our three introductory mathematics courses: our course preceding calculus and the first two semesters of calculus. In the regular academic year, these courses enroll about 4500 students total. Despite their large sizes, the courses are taught in small sections (capped at 18 students), in a student-centered manner that emphasizes student collaboration on in- and out-of-class work with support from the section instructor. This pedagogical model was implemented with calculus reform in the early 1990s, when the courses adopted a strong conceptual emphasis in which routine skills (e.g., differentiation in Calculus I) are assessed almost exclusively through mastery ("gateway") tests (LaRose, 2018; LaRose & Megginson, 2003) and other assessment is focused on deep conceptual understanding of the mathematics that students learn. Their structure is validated by recent research: small, collaborative classes with active student work are associated with better outcomes for women and groups underrepresented in STEM (Freeman et al., 2014; Laursen et al., 2014). Typically, 70% of the sections of our introductory math courses are taught by graduate students, 15% by post-docs, 13% by continuing full-time non-tenure track faculty[17], and 2% by tenured or tenure-track faculty. Importantly, each course is managed by a faculty course coordinator who designs the assessments and syllabi and supports all instructors and students in the course.

---

[16] As with many flagship public research universities, ours has never managed a large number of minority students, and its efforts in that direction have been hamstrung by a constitutional mandate eliminating direct consideration of race in admissions since the mid-2000s: Black enrollment was consistently about 7% before the change, and about 4% since [12,13].

Inclusive teaching and responsiveness to diversity is essential in these courses, especially in the course before calculus (where we first introduced mastery assessment, as described below). That course's enrollment has a higher fraction of first-generation students (15%), underrepresented minorities (35%; defined as self-reported race/ethnicity groups other than white and Asian), and women (58%) than the university as a whole (7%, 13%, and 50%, respectively). Moving from this course into calculus, we see drops in the fractions of students in all of these categories. Improving the learning experiences, outcomes, and retention of students from groups historically marginalized in mathematics motivates these changes in our courses, by which we strive to help redress the differential barriers they face in higher education.

To this end, the current work is focused along two dimensions of change. Because many of our instructors are at the start of their teaching careers, we have a week-long training program for them. The first dimension extends this training and provides ongoing support to facilitate their growth as inclusive teachers. The second dimension is in the revision of the assessment in these courses, starting with our course before calculus. For that course we have developed a consistent set of assessments that promote mastery learning and are sustainable at a very large scale. Below, we describe these two dimensions. Recognizing the iterative nature of course design, we also share lessons learned for how to improve our efforts in the future and conclude with a discussion of the evolution of the work we describe here.

## 3. *Dimension of Change 1: Training for New Instructors and Ongoing Teaching Support*

To support instructors who are new to teaching in our introductory math courses, we run a training program during the week before classes start in the fall term.[18] The overarching goal of this program is to motivate instructors to teach in an interactive, inquiry-focused manner, and to provide them with the experience and tools they need to do this effectively. The program format and topics have evolved over the years, and currently focus on experiential practice. The participants practice and receive feedback on lecturing, experience a model class as students to see how it can promote learning, and practice running this type of class themselves. There are also sessions on the logistical and administrative aspects of teaching. The program has also a strong emphasis on getting buy-in from the participants and on how they can get student buy-in when they go into their classes.

---

[17] We aim to hire instructors only into full-time positions, a goal which we meet in many but not all semesters.

[18] The instructional team in charge of this program includes coordinators and graduate student co-coordinators of the introductory courses, undergraduate directors, and some other associated faculty.

We have added an inclusive teaching workshop to the first day of this program and made complementary changes to the content of the programming in the remainder of the week. The two-hour workshop, run by our campus Center for Research on Teaching and Learning, follows a theatrical skit and discussion that bring out issues of inclusive teaching. The workshop starts by introducing instructors to the idea of teaching inclusively and why it matters, emphasizing that as inclusive teachers they must be deliberate and aware of the impact of students' identities and systemic inequities. It then provides concrete strategies that they can use to increase transparency about expectations and evaluation criteria, cultivate a sense of academic belonging in their classrooms, and increase the structure of classroom interactions to promote an inclusive classroom community. Changes in the programming in the remainder of the week include use of these strategies in sessions on group work, and a focus on norms and expectations for classrooms and student teams.

Recognizing that developing and refining one's teaching philosophy takes time, we also provide ongoing support to our instructors. In their first semester teaching, all new instructors receive one or two class visits to provide constructive feedback on their teaching and on teaching inclusively. Additional support is provided by weekly course meetings to discuss course material and the logistics of teaching it, and in our course preceding calculus[19] we have added to these two follow-up workshop sessions on inclusive teaching for instructors. We also solicit student feedback mid-semester for the instructors to use to improve their teaching.

There are several key lessons from this dimension of change. First, graduate students found activity-based sessions in the new instructor training program (especially those pertaining to group work and inquiry-based learning), and informal discussions with other instructors and with course coordinators, to be particularly valuable. Second, these informal supports continue to be rated as helpful for instructors with more experience. And third, additional mentorship and feedback on teaching, as well as more opportunities to build a peer support community, were the most frequent suggestions for further support. Finally, we note also an ongoing challenge for our program: while many of our instructors are fully invested in the work of teaching (and teaching inclusively), this competes with the academic norm of the preeminence of research. As such, a challenge we face is to shift the mindsets of (some of) our new instructors so that being an effective and inclusive instructor is considered a fundamental part of being a successful mathematician.

[19] Because this is the course on which we focused initially, these workshop sessions were done for those instructors. We expect that this will expand as we focus also on Calculus I and II.

This work to help new instructors teach inclusively exists in the context of a wider effort in our department to create a community of instructors who teach inclusively, and to build a consensus in the Department that this is necessary. We have an ongoing Learning Community on Inclusive Teaching (LCIT) that meets on a monthly basis to provide an informal space in which instructors are able to discuss selected readings on inclusive teaching, raise questions they have about teaching, and recommend things they have done in their classrooms that have worked well. In many respects, this provides an existent space in which some of the informal support that instructors have as they start to teach can continue, with a specific focus on inclusive teaching. The LCIT meets over lunch, and for its first two years of existence has been supported by small university grants that have allowed us to provide lunch for the participants.[20] While participation in the learning community is voluntary, it accomplishes—for its participants—the goals of providing a venue for continued exploration and support of inclusive teaching in the department. It has engaged a significant number of instructors who teach our introductory courses. Reaching all instructors in the Department and building a consensus that this work is essential presents a much bigger problem that requires systemic change at a larger scale than the courses in which we are currently working. This is work we must do but is beyond the scope of the project we describe here.

## 4. *Dimension of Change 2: Course Assessment*

Core to the philosophy and pedagogy of our introductory courses is an emphasis on conceptual understanding, learning as a growth process, and a small class, active-learning instructional model[21] that has been shown to be supportive of all students in mathematics. To better align the primary assessment in our course preceding calculus with this set of core principles, we have developed and pilot-tested a new assessment structure for the course that includes a significant mastery learning component. This replaces a model in which 95% of students' course averages was generated by three timed (and necessarily high-stakes) exams and a scaled grade distribution[22], which has been shown to put women and underrepresented minorities in STEM at a disadvantage (Madeaus & Clarke, 2001; Piontek, 2008). By contrast, a repeatable, low-stakes assessment that allows students to build understanding and ownership of the course material addresses many of the factors that have been shown to drive these students away from STEM fields (Seymour & Hewitt, 1997).

In the new assessment model, 50% of the points determining students' grades are earned on repeatable mastery assessments, and a further 10% from work done in class and class assignments. The remaining 40% is derived from two timed but much lower stakes exams that are shorter and significantly less formidable than the previous three. This model has two key features. First, it constitutes a comprehensive revision to the assessment in the course, in which the components are realigned for consistency with the course, other assessments, and course goals. And second, it aligns the grading structure of the course with the core goal that students learn in a growth mindset (Dweck, 2006) manner, which aligns with our core course goals.

---

[20] The department also has an Inquiry Based Learning Center which also holds informal lunches, the focus of which clearly complements the inclusive teaching lunches.

[21] As a general rule of thumb, we expect instructors to lecture in blocks of 20 minutes or less, and that at least half of the time in class will be spent on exploratory activities that students work on in groups.

[22] We note that while scaled based on a historical understanding of the meaning of students' performance, grades were not “curved” in a traditional sense: there was no fixed fraction of students who received any given grade, for example.

The development of the mastery assessments followed a model characteristic of such work: we determined first a comprehensive list of learning objectives for the course, and then created a set of repeatable assessment problems with randomized parameters that allow us to determine whether students have mastered each. These are grouped into eight assessments, each having five questions, which we administer using our on-line homework and testing system (WeBWorK, n.d.). Because our goal is for students to master the learning objectives, we require that students answer four or five of the five questions correctly to receive partial or full credit for the assessment. In total, students can earn 35% of the points in the course by completing these assessments. At the end of the semester, we have a final mastery assessment which includes topics from all preceding assessments and is ten questions long. These scores contribute an additional 15% of the points towards students' course grades. For all mastery assessments, students can practice as often as they like online, and receive credit by obtaining a high enough score when taking it in a proctored lab.[23] The lab is open 40 hours a week and students have approximately two weeks to complete each mastery assessment. We have also allowed students to choose two assessments to reopen at the end of the semester.

To complement the use of the mastery assessments, we made two additional changes to the course. First is an increase in the assessment weight of the work that students do in and for their class section. In the past, students' work on the web homework for the course was the last 5% of the credit in their course grade. Now, 10% is determined by their work on the web homework, quizzes in their section, and team homework done with other students. The team homework requires the solution of significant, conceptual problems, which students submit in carefully written solution papers describing their work and the mathematics they used. Both the problem work and the solution paper are produced collaboratively by the team. The final change to the course assessment is, of course, in the timed exams. We changed from two midterms and a final to two exams, one a midterm and the other a non-comprehensive exam toward the end of the semester. While the exams were in the past 8–12 problems long, with a 90-minute (for midterms) or 120-minute (for the final) time limit, the exams are now about 6 problems long, with a 90 minute time limit. The exams also focus specifically on the learning objectives that we are unable to evaluate easily on the mastery assessments: higher-order problem solving, graph sketching, and the description of mathematics in written form.

The grading scheme for the course then sets thresholds for each letter grade (80 points for an A, 65 for a B, etc.) and a minimum number of points obtained on the mastery assessments (30 for an A or B, 24 for a C, etc.) and the final mastery assessment (9/10 for an A or B, 8/10 for a C, etc.).

---

[23] ...except that in the winter of 2020 we abruptly taught the course online for the second half of the semester—an instructional mode that continued through the spring of 2021. In this case we were unable to proctor the assessments, and students self-proctored and submitted their hand-written work online. We also experimented with non-punitive online proctoring support for the credit-bearing assessments. Starting in the fall of 2021 we returned to in person instruction and assessment and discuss further changes to our assessment model in our conclusion.

Overall, this change in the course appears to have worked very well, even with the disruption caused by the onset of the COVID-19 pandemic midway through the first semester in which it was implemented. The pandemic necessitated an abrupt shift to online learning (with a corresponding and dramatic increase in the support needed by both instructors and students), which undermined our assessment plans for the project. The change to remote learning also makes it difficult to draw quantitative comparisons between semesters with and without mastery assessment: not only did the instructional mode for the course change, but the change and the pandemic itself have had a disproportionately negative impact on many of the populations that are underrepresented in STEM. That said, we have qualitative data that suggest that the changes have had a significant positive impact.

Key among these qualitative measures is that all students in the pilot who put in sufficient work to persist to the end of the semester passed the course. With the previous, exam-centered, assessment we have not been able to make this claim about any semester (even in the absence of a pandemic and its impact on the course). A related feature of the new assessment model is that now students have much better information, from a range of assessments making up almost half the credit in the course, to use when they are making the decision to continue in or to drop the course. In the previous format many students would either drop the course or give up after receiving a low grade on the first midterm, while others would stick with the course because the first midterm was only 25% of their grade—even though from the instructor's perspective it was unlikely that the student could succeed.

The new assessment structure also is better able to deal with exceptions necessitated by student circumstances or other factors external to the course: if a student is unable to complete a mastery assessment because of illness or other factors, it is straightforward to reopen it for them. This is not the case when the preponderance of the assessment is by synchronous, timed exams. The change has an indirect impact on how equitable the course assessment is, in that it makes it easier to create a course policy that is responsive to all students' needs and external constraints.

Anecdotally, instructors reported far fewer complaints from students about the exams, and students—while feeling at times overwhelmed by the work demanded by the mastery assessments—were strongly in favor of the change. Instructors also reported that students gained confidence as they passed mastery assessments. The feeling of accomplishment, that they could do the course work, and that their performance was in their control were all positive changes in student attitude. There is evidence supporting these anecdotal reports: students' persistence, as measured by the number of attempts that they made on each of the assessments, continued at a high level throughout the semester (the average number of attempts per student on the last four assessments was slightly higher than that for the first four). This change is important: when faced with three exams some students felt that their grade was not in their control but would instead depend on what the exam looked like or how the grading scale turned out. This premise carries with it significant negative implications for student learning: for example, rather than truly figuring out a part of a difficult topic, students would either hope it did not come up on the exam or hope that others would also struggle with it so that they could get by on the grading scale. Or they would guess or omit the solution to one part of the topic and bank on losing at most a few points. By contrast, with the mastery assessments students aiming for an A or B in the course can't just skip over a topic, because of the minimum mastery point requirements for the different letter grades. Further, students can explicitly identify the ideas they were struggling with and continue to work on them until they attain mastery. This overall evidence of a transition to a growth mindset is another, though indirect, measure of our potential at supporting underrepresented students, as adoption of growth mindsets in classrooms has been shown to have a positive impact on student success, particularly for groups that are traditionally underrepresented (Hill et al., 2010; Sisk et al., 2018).

## 5. *Conclusion*

In this case study, we have described a set of initiatives contributing to a fundamental change in an introductory mathematics course and its instruction. These changes produced a more inclusive and responsive experience for all students, especially for students from communities that have been historically excluded or marginalized in college math courses. In doing this, we recognize that this change requires work from individual class sections and, especially for a large course at a large university, changing the structure of the course. We further acknowledge that this type of change requires the support and investment of many stakeholders: our work on this course has significant and ongoing support from our Department, our Center for Research on Learning and Teaching (including its “Foundational Course Initiative” program), and our College and University.

Importantly, this work is ongoing. We describe here the core aspects of the changes made to the structure of our course before calculus, which we continue to revise in an iterative manner. We update our new instructor training program annually and are beginning to update Calculus I to include mastery assessment. We expect to then revise Calculus II. In our course before calculus we have made changes to the model described above to address logistical challenges, student workload, and test accessibility. We have returned to our proctored labs for credit-bearing assessments, increased the hours the labs are open to 55 hours per week, and increased the number of labs that are available to students. We have also condensed the mastery assessments to reduce the number and amount of time in the semester that students must be focused on their completion. The resulting assessment model has six mastery assessments of seven problems each (in total, 40% of students’ course averages), has increased the weight of the class work and web homework (to 15%), and has three exams (each worth only 15% of students’ course averages). We fully expect that this model will evolve as we continue our work.

Finally, we stress that the work we describe here is only a tiny start to a process that transcends any one course or institution. Our history has created a world in which inequity and differential access to learning are ingrained in the structure of our higher education, society, and culture, and these inequities must be addressed on that scale.

# 4 Promoting Racial Equity in the STEM Math Pathway in Community Colleges

**Helen Burn**

*Highline College*

Helen Burn, Ph.D., is an instructor in the Department of Mathematics and director of the Curriculum Research Group at Highline College. Her research focuses on community college mathematics curriculum including reform of pre-college mathematics and college algebra, and supporting adjunct faculty and the partner disciplines. Dr. Burn is chair of the Pathways Joint Subcommittee of the American Mathematical Association of Two-Year Colleges. Dr. Burn is a PI on the NSF Collaborative Research Grant: Transitioning Learners to Calculus in Community Colleges (TLC3).

**Chauntee Thrill**

*Appalachian State University*

Chauntee Thrill, Ph.D. is an Assistant Professor of Higher Education, in the Department of Leadership and Educational Studies at Appalachian State University. Dr. Chauntee Thrill's research interests center on postsecondary education, focusing specifically on college readiness/remedial education, affordability, as well as access and success for underserved student populations, especially within the community college sector. Dr. Thrill is a research associate on the NSF Collaborative Research Grant: Transitioning Learners to Calculus in Community Colleges (TLC3).

**Vilma Mesa**

*University of Michigan*

Vilma Mesa, Ph.D., is Professor of Education and Professor of Mathematics at the University of Michigan. Dr. Mesa investigates the role that resources play in developing teaching expertise in undergraduate mathematics, specifically at community colleges and in inquiry-based learning classrooms. She served as Associate Editor for the Journal for Research in Mathematics Education from 2000-2004 and is currently serving as associate editor for Educational Studies in Mathematics. Dr. Mesa is a PI on the NSF Collaborative Research Grant: Transitioning Learners to Calculus in Community Colleges (TLC3).

**Eboni Zamani-Gallaher**

*University of Illinois at Urbana-Champaign*

Eboni Zamani-Gallaher, Ph.D., is Professor of Higher Education/Community College Leadership and Director of the Office for Community College Research and Leadership (OCCRL). She previously served as Associate Head of the Department of Education Policy, Organization, and Leadership and Associate Dean of the Graduate College at the University of Illinois at Urbana-Champaign. Her teaching, research, and consulting activities largely include psychosocial adjustment and transition of marginalized collegians, transfer, access policies, student development and services at community colleges.

**J. Luke Wood**

*San Diego State University*

J. Luke Wood, Ph.D., is Vice President of Student Affairs & Campus Diversity and Dean's Distinguished Professor of Education at San Diego State University. Wood also serves as the Co-Director of the Community College Equity Assessment Lab (CCEAL), a national research and practice center that partners with community colleges to support their capacity in advancing outcomes for underserved students of color. Dr. Wood is a PI on the NSF Collaborative Research Grant: Transitioning Learners to Calculus in Community Colleges (TLC3).

*Abstract: Transitioning Learners to Calculus in Community Colleges (TLC3) is a mixed-methods research project focused on identifying practices that promote the success of African American, Latinx, Native American, and Southeast Asian students as they transition into and through mathematics courses required for degrees in science, technology, engineering, and mathematics (STEM). This paper focuses on accelerating students to college-level mathematics courses and supporting students' mathematical learning inside and outside of class, and explains the potential of these practices to enhance the success of racially minoritized students in the STEM math pathway. Accelerating students to college-level mathematics involves optimizing the number of courses in the STEM math pathway and ensuring accurate placement of students into their initial mathematics course. The findings reveal a need for equal attention to acceleration at the precalculus and calculus level in addition to the developmental level. Welcomeness to engage was used as a general instructional practice for all students, and there was limited use by instructors of other instructional strategies to be inclusive of racially minoritized students in their classrooms. Students access out-of-class support for mathematics learning mainly through faculty office hours, campus tutorial centers, and dedicated space. While we highlight practices and policies implemented at our case study colleges, which were identified as being successful in supporting African American, Latinx, Native American, and Southeast Asian students in moving through the STEM math pathway, the findings align with prior research that indicates there are practices in mathematics classrooms that uphold whiteness as normative and dominant.*

## 1. *Introduction*

Despite the framing of equity and opportunity gaps in mathematics around race/ethnicity, efforts to close gaps tend to take a race-neutral approach that "lifts all ships" but does little to advance racial equity in mathematics. Colorblind approaches reinforce whiteness as normative in mathematics, creating barriers to racialized student engagement, learning, and attainment in mathematics (Battey, 2013; Battey & Leyva, 2016; Davis & Martin, 2008; Martin, 2009; Gutiérrez, 2012, 2017; Palmer & Wood, 2013; Wood & Palmer, 2014). *Transitioning Learners to Calculus in Community Colleges* (TLC3) is a research project focused on identifying successful practices in mathematics programs in community colleges, through an equity-conscious lens (Burn, Mesa, Wood, & Zamani-Gallaher, 2016). We study practices that promote the success of African American, Latinx, Native American, and Southeast Asian students as they transition into and through mathematics courses required for degrees in science, technology, engineering, and mathematics (STEM).

The TLC3 study is based on survey data collected from 455 (44%) of the nation's 1,023 public, primarily associate degree-granting institutions (hereafter referred to as community colleges) and case studies of four Minority Serving Institutions (Table 1).[24] Conceptual support for the work includes Harper's (2010) anti-deficit framework, Wood et al. (2015) enhanced practices to support men of color in community colleges, and the Mathematical Association of America's National Study of Calculus (Bressoud et al., 2015; Burn et al., 2015). This paper focuses on two areas of practice: accelerating students to college-level mathematics courses and supporting students' mathematical learning inside and outside of class. We draw on the literature to explain the potential for each area of practice to contribute positively to the success of African American, Latinx, Native American, or Southeast Asian students in the STEM math pathway. Following this, we showcase practices at the PBI and Tribal College around accelerating students to college-level mathematics. In the area of supporting students' mathematical learning inside and outside of class, we provide positive examples gleaned from classroom observations and highlight practices at the AANAPISI and HSI colleges.

[24] Minority Serving Institution (MSI) designations include Asian American, Native American, and Pacific Islander Serving Institutions (AANAPISIs), Predominately Black Institutions (PBIs), Historically Black Colleges and Universities (HBCUs), Hispanic-Serving Institutions (HSIs), and Tribal Colleges and Universities (TCUs). With the exception of Tribal Colleges, MSI status is based on enrollment: HSIs enroll 25% or more Latinx, PBIs enroll 40% or more African American students, AANAPISI enroll 10% or more Asian American, Native American, or Pacific Islander Students.

## 2. *TLC3 Case Study Colleges*

We selected four Minority Serving Institutions (MSI) for case study (Table 1). The selection criteria included the college's responses to the TLC3 National Survey of Community College Mathematics Chairs (Burn, Mesa, Wood, & Zamani-Gallaher, 2018), available state-level mathematics attainment data, the college's expressed commitment to supporting racially minoritized students within their MSI designation, and their willingness to participate in the case study. Each case study college received local Institutional Review Board approval. Data collected during site visits conducted in 2018-19 included 31 interviews with faculty or staff, five student focus groups (45 students), and 26 classroom observations (715 students observed and surveyed). Table 1 shows the colleges selected for case study, their MSI designation, location, total enrollment, and number of full- and part-time mathematics faculty.

*Table 1*
*Transitioning Learners to Calculus in Community Colleges (TLC3) case study colleges[a]*

| College MSI[b] Designation | Location | Total enrollment (% in MSI designation) | Number full-time (part-time) mathematics faculty |
|---|---|---|---|
| AANAPISI College[b] | West, United States | 7,765 (13%) | 9 (14) |
| HSI College | West, United States | 12,521 (67%) | 14 (11) |
| PBI College | Midwest, United States | 4,232 (55%) | 14 (12) |
| Tribal College | West, United States | 367 (100%) | 2(6) |

[a]Data from the Integrated Postsecondary Education Data System (IPEDS) except full/part-time mathematics faculty was reported in the TLC3 national survey (https://files.eric.ed.gov/fulltext/ED592079.pdf)

[b]See Footnote 1 for MSI designations

## 3. *Accelerating to College-Level Mathematics*

The STEM math pathway in community colleges is lengthy and can include courses that range from developmental mathematics to Calculus II. Accelerating students to college-level mathematics involves optimizing the number of courses in the STEM math pathway and ensuring accurate placement of students into their initial mathematics course (Hodara, 2019). To date, many colleges have undertaken efforts to ensure the accurate placement into initial mathematics. However, once placed, students have access to a wider variety of course options, inclusive of acceleration, at the developmental level than at either the precalculus or calculus level (Burn et al., 2018).

Acceleration is an equity strategy. Indeed, 34% of historically underrepresented racially minoritized students (African American, Latinx, and Native American) who enter STEM fields require at least one remedial mathematics course compared to 18% of their Asian American and White peers (Herrara & Hurtado, 2014). The additional mathematics course load has been shown to slow academic progress (Herrara & Hurtado, 2014; Hodara, 2019). Racially minoritized students disproportionately begin at the developmental level for different reasons including unequal opportunities or inadequate preparation prior to enrolling in college, length of time since last math course, or because they are erroneously placed into developmental mathematics by an institution's placement process (Atuahene & Russell, 2016; Hodara, 2013, 2019; Melguizo et al., 2016). Being misplaced into developmental mathematics can signal to students that they do not belong in college (Cox, 2009), and the effect can be greater for historically underrepresented student subgroups (Wood et al., 2015). Alternatively, when students feel fairly placed into their initial mathematics course, even if they are placed at the developmental level, placement may be seen as having a positive influence on their success.

The PBI college in our study undertook efforts to accelerate students to college-level math by making changes to developmental mathematics and mathematics placement. First, the college began offering accelerated developmental courses by revising their existing classes from three- to two-hour courses offered in an 8-week format. This change reduced the overall number of developmental mathematics credits from nine to six hours and the amount of time students remained in the developmental sequence while covering the same content. At the time of study, this was a new initiative and data was not available about its effectiveness. Second, to ensure more accurate initial placement, the college began offering multiple options for placement including high school mathematics grades for Algebra I and Algebra II, ACT/SAT scores, or the ALEKS placement test. The college implemented a mandatory review enforced by the college placement testing center of no less than three hours prior to attempting the ALEKS test. At a cost of $25, students are provided an access code for online study materials that allows them to take assessments and identify mathematical topics they need to review, and the ability to retest up to five times but limited to retesting twice in one semester. In addition, the college adopted a policy whereby faculty members teaching the three lowest levels of developmental mathematics courses can advance students to the next course by offering students an assessment that affords students the opportunity to "test out" during the beginning of the term by demonstrating mastery of content. One student described the process as,

> *When I first took the placement test, they put me in a [developmental] class. But when we get to the class it was like okay, we're going to give you this test to see if you really belong here. When I took the test, I ended up being removed from that class and being placed into college level.*

The Tribal College in our study likewise adjusted both developmental mathematics and mathematics placement. The college modularized the two semester-long developmental courses with the highest failure rates into two 8-week segments. The new design enabled students who had not mastered the course content by the end of the first segment module to retake the material in the second part of the semester, leading to dramatic increases in completion rates. In addition, students reported that they preferred the modules because they did not feel stressed to learn all the material at once. The college also offered students multiple placement options that are presented to students during orientation. Options included prior college mathematics courses, ACT, PSAT, high school transcripts, the Partnership for Assessment of Readiness for College and Careers (PARCC) assessment, and the ACCUPLACER test. The Academic Coaching Center prepares students for taking the ACCUPLACER exam, a preparation that has been effective in increasing ACCUPLACER scores. The institutional policy is to make sure that students are placed adequately, and the college has two strategies to implement the policy: (1) challenge and retake the ACCUPLACER score and (2) allow the mathematics chair to recommend students for higher placement. A higher placement decision can also be influenced by student performance in the Summer Bridge program, a program for all incoming new students.

## 4. *Supporting Students' Mathematical Learning*

### 4.1. *In-Class Support*

The dominant instructional method observed in mathematics classrooms in community colleges is interactive lecture, a mix of lecture with opportunities for students to engage in the course content mainly through fielding questions or working problems in class (Blair et al., 2018; Burn & Mesa, 2017). Community college mathematics faculty using this instructional method described creating conditions for student engagement by building relationships with students, teaching for understanding, and leveraging personal attributes such as instructor approachability and charisma (Burn & Mesa, 2017). It is important to underscore that Wood et al.'s (2015) research establishes that factors under the direct control of faculty members--such as positive faculty-student interactions and classroom cultures that support, affirm, and validate students--contribute significantly more to the success of men of color in community colleges than do factors such as students' background characteristics or environmental pressures. However, African-American and Latinx students may be apprehensive to engage with white faculty[25] because of their lived experiences with racism in the form of ascriptions of low intelligence, assumptions of criminality, and racial microaggressions (Wood et al., 2015). In contrast, research also shows that Southeast Asian students (e.g., Cambodian, Hmong, Khmer, Laotian) may hesitate to engage with faculty due to cultural factors around respect for authority (Museus & Kiang, 2009).

Mathematics faculty can increase positive and productive faculty-student interactions through practicing welcomeness to engage (Wood et al., 2015). Welcomeness to engage entails faculty creating conditions inside and outside the classroom that communicate to students that their engagement is not just welcomed but also desired, with the latter being a crucial aspect. We observed features of welcomeness to engage in 21 of 26 mathematics lessons observed. The features were grouped into three categories: basics of welcomeness, physical aspects of welcomeness, and welcomeness to engage in mathematics (Burn, Gerhard, Blevins, & the TLC3 Research Team, 2019). The basics of welcomeness included arriving to class early or on time, interacting with students as they entered class, and learning and using students' names. Physical aspects of welcomeness included looking at and acknowledging students when they asked questions, gesturing to acknowledge students (i.e., pointing, nodding), circulating and engaging openly with students while they worked on problems, and leaning in when working with students.

[25] In two-year colleges, 77% of permanent full-time faculty and 78% of part-time faculty identify as non-Hispanic white (Blair, Kirkman, & Maxwell, 2018, p. 188).

We also observed ways that instructors welcomed students to engage that were specific to mathematics: acknowledging the difficulty of the material while assuring students that they could do it, encouraging students to "take a shot" at a problem, and finding ways to invite more students to participate to mitigate the tendency of one or two individuals to dominate in-class questions. We noted that instructors would keep students engaged in doing mathematics by reframing a mistake or by validating an approach a student had taken on a problem, for example by saying "factoring, I like it!" or "right" when a student did well. It is noteworthy that in the classes observed, welcomeness was used as a general instructional practice for all students rather than being a strategy used intentionally to bolster the engagement of racially minoritized students. We also noted limited use of other instructional strategies to be inclusive of URM students in their classroom. For example, we identified instances of using empowerment strategies for URM students (e.g., acknowledging a student's or group's contribution, assigning leadership tasks, or offering students praise) in only five of 26 classrooms observed. We identified instances of culturally relevant teaching (e.g., connecting the course content to contemporary topics or real-life experiences or using vocabulary used by specific student populations) in 10 of 26 classes observed.

### 4.2. *Out-of-Class Support*

Faculty office hours and campus tutorial centers are the two main ways students seek support outside of class. Students also benefit from having access to specially-designated space on campus, separate from a tutorial center, to gather to work on mathematics. Such space can be crucial on campuses where tutoring for upper-level STEM math pathway courses is less available (Burn et al., 2015). As with in-class support, cultural norms and expectations can influence students' help-seeking behaviors outside of class. Masculine identity may interfere with a students' willingness to seek help (Wood et al., 2015). The model minority stereotype can make Asian American and Pacific Islander (AAPI) students reluctant to seek out-of-class support because doing so violates the stereotype that Asians are academically gifted and do not need assistance (Museus & Kiang, 2009).

At the AANAPISI college, in addition to faculty office hours, students could access the campus tutorial center where they sign up and have access to trained tutors for all levels of mathematics. The center was overseen by a full-time mathematics faculty; at the college, more than half of the mathematics faculty devoted time to the tutoring center. In addition, the college had a large, open atrium in proximity to faculty offices that provided an optimal space for students to work and to access computers. Two of the five instructors interviewed at the AANAPISI college described the extra steps they take to build relationships with their Southeast Asian students to bolster their use of office hours and tutoring. One instructor at the college described,

> *At least for me, that's where that kind of personal dynamic and my engagement with them in the class really comes in [ . . . to] create a relationship with that Southeast Asian student, which is definitely a harder sell in many cases, but usually one that's accomplishable.*

At the HSI College, students in the STEM math pathway are encouraged to interact with each other and their instructors in a STEM-focused learning lab created by the college. Within the lab, students can work with preferred tutors and with supplemental instruction student leaders. The space also provides whiteboards and laptops to aid in students' learning. To encourage and promote the use of the learning lab, tutors conduct classroom visits to encourage students to seek out support through the lab. One instructor shared that their class takes a field trip to the STEM-focused learning lab on the first day of class to familiarize the students with the resource. Additional out-of-class support occurs through Saturday math classes and free boot camps that prepare students to succeed in the next course in their mathematics sequence. A participant in a student focus group described their experience with a Math Academy offered by the college and how it can enable students to accelerate in mathematics. The student described,

> *There's also something called Math Academy that you can do over the summer or the winter. It's like a 2-week math program that's a few hours every day. It's really big here, and the spots go out really quickly. I helped tutor for the Math Academy and I participated in it myself, and it really helped, and it counts. It allows you to skip a course up to precalculus if you do math academy every semester and winter.*

### 5. *Implications*

The possibility of racial equity in mathematics lies in creating new approaches and normative structures with intentionality towards who the beneficiaries might be (National Academies of Science, Engineering, and Medicine, 2019). Our data is reflective of national trends around redesign of mathematics placement and developmental mathematics, suggesting these twin efforts are becoming normative practice with positive benefit to URM students. However, achieving racial equity in the STEM math pathway requires equal attention to acceleration at the precalculus and calculus level, which pose significant structural inhibitors to URM student success (Palmer & Wood, 2013; Wood & Palmer, 2014). It was disheartening to see limited use by instructors of other instructional strategies to be inclusive of racially minoritized students in their classrooms. This finding speaks to the need for faculty professional development targeting strategies to work with students of different ethnic and cultural backgrounds based on an understanding of how the norms of whiteness create barriers to racialized student engagement and learning in mathematics.

## 6. *Conclusion*

Community colleges are seen as a vehicle towards addressing racial disparities in STEM because they enroll higher proportions of racially and ethnically diverse student populations. While we highlighted practices and policies implemented at our case study colleges, which were identified as being successful in supporting African American, Latinx, Native American, and Southeast Asian students in moving through the STEM math pathway, the findings align with prior research that indicates there are practices in mathematics classrooms that uphold whiteness as normative and dominant. Instructional practices determine which groups are recognized as contributing to the knowledge base, and thus become more prepared and successful. Such practices still ignore non-white students and continue to serve as a barrier to access and progression for racially minoritized students pursuing the STEM math pathway.

## 7. *Acknowledgement*

Support for this work is provided by the National Science Foundation's Improving Undergraduate STEM Education (IUSE) program under Awards 1625918, 1625387, 1625946, and 1625891. Any opinions, findings, and conclusions or recommendations expressed in this material are those of the author(s) and do not necessarily reflect the views of the National Science Foundation.

# 5 Stretching Calculus: A Yearlong Calculus Class that Stretches Minds

**Martha Byrne**

*Sonoma State University*

Martha Byrne is an Assistant Professor in the Department of Math and Stats at Sonoma State University where she works at the intersection of math and math education. She is very involved in the content mastery preparation of pre-service teachers at SSU.

**Natalie Hobson**

*Sonoma State University*

Natalie Hobson is an Assistant Professor in the Department of Math and Stats at Sonoma State University. She works closely with students on the teaching track and on the pure track of the major and teaches content courses for pre-service elementary teachers.

**Izabela Kanaana**

*Sonoma State University*

Izabela Kanaana is a Professor in the Department of Math and Stats at Sonoma State University. She is involved in curriculum development and teaching calculus, as well as discrete mathematics courses.

**Brigitte Lahme**

*Sonoma State University*

Brigitte Lahme is a professor in the Department of Math and Stats at Sonoma State University. She is very involved in K-16 curriculum development, with a focus on calculus and working with first year college students.

*Abstract: At Sonoma State University, as part of the transition away from remediation, we developed and implemented a yearlong calculus I course for students who would have previously been placed into remedial courses. The specific components of our stretch calculus program are intended to make the learning environment more equitable and inclusive in order to better support first-year students who were determined not college ready in mathematics. These components include mindset interventions, active learning, group work, and peer support. The active learning principles employed in our stretch program are motivated by the TRU framework, developed by The Mathematics Assessment Project. We discuss in this paper various aspects of success of the program. While we are proud of our program and our students' accomplishments, throughout the process of development and implementation, we confronted several obstacles and challenges, and from these, we made recommendations for other schools attempting such a model.*

## 1. *Program Setting*

Sonoma State University (SSU) is one of the 23 campuses of the California State University System (CSU). Of the 8,565 SSU undergraduates enrolled in the Fall 2018 semester, 51.6% identified as other than white only, 36% identified as Latinx, and 22% were first-generation college students. SSU generally accepts all qualified students who apply from high schools in its service area, and 99% of students come from California. SSU experienced an 80% increase in the percentage of Latinx undergraduates enrolled between Fall 2013 and Fall 2018 (1,590 to 2,871). In February 2017, SSU became a Hispanic Serving Institution (HSI) and is now one of 21 CSU campuses so designated. At the CSU level, students' status as an underrepresented minority (URM) is collected. At SSU, 95% of these URM students are Latinx. As such, though some of our discussions or data is on URM students, the majority of this population are Latinx. Despite comprising 36% of the university's student body in 2018-19, Latinx students only received 24% of undergraduate degrees in the School of Science and Technology (SST).

In fall 2017, aligned with a national trend, SSU began to transition away from remedial courses in order to better serve the students arriving on campus in need of more mathematical support. Prior to that fall students arriving at SSU without AP (or other college-level) credit all took SAT or ACT as placement exams to determine whether or not they were "ready" to take a college-level math course. Students who were not deemed ready were placed into one of two remedial algebra classes for which they were charged regular tuition rates, but for which they did not receive any credit toward the 120 unit minimum required for graduation. The students who were placed into these remedial courses came disproportionately from low-income households and poorly ranked school systems when compared to the rest of the student body, and they were also disproportionately URM students. Our department began our redesign project to provide a supported path to success for these students who had been badly served by educational systems previously.

What we have done is developed and implemented a yearlong calculus I course (stretch calculus) specifically for students who would have been placed into remedial courses. We were motivated by a desire to support students academically and improve STEM retention and 4-6 year graduation rates. In particular, as a department, we wanted to provide students more agency and ownership over their mathematics. Historically, calculus has been a course with large opportunity gaps (with respect to GPA and persistence in STEM). We expand on the impacts of the stretch program on these gaps in the results section below. We acknowledge that a focus on GPA gaps and retention rates is limited at best, but we do not have the data we would like to better address the complex issues in diversity, equity, and inclusion.

## 2. *Motivation and Design Process*

The students for whom we designed our stretch calculus program, historically would have been required to pass at least one semester of a non-credit bearing algebra course before registering for a college credit bearing math course. They also would have needed to take either two or three semesters of math before having the opportunity to take calculus, which made most STEM majors essentially inaccessible. While we were aware that this population disproportionately comprises URM students, we did not focus specifically on the racialized experiences of the students in the course design. We were motivated to provide better access and pathways to the entire group of poorly-served students believing that by implementing several key components, we would support all students.

To design our course, we used several of the themes of the *Teaching for Robust Understanding* (TRU) framework (https://truframework.org/) which was developed by the Math Assessment Project. Central to TRU are five dimensions of powerful classrooms. Of these dimensions, we focused primarily on (1) *Agency, Ownership, and Identity* and (2) *Equitable Access to Content.* That is, the classroom activities should provide students with an experience to take ownership over constructing robust mathematical ideas and invite all students to actively engage in meaningful ways. Our focus on these dimensions was largely due to our departmental belief that success in mathematics has more to do with access and mindset than it does on previously measured mathematical skills.

During the 16-17 academic year, a task force of lecturers and tenure-track professors met regularly to review research literature and discuss course design. Members of the task force as well as other members of the department attended a workshop led by Cinnamon Hillyard from UW Bothell about mindset interventions, opening students' minds to a different view of and identity within the world of mathematics, and supporting students' transition to college. During this workshop we planned a series of interventions for the first two weeks of the semester to set the tone for the course. These interventions included:

- "Group Noticing": group members check in with missing peers at the beginning of class
- "What am I good at? How did I get there?": Students silently answer both questions, then all answers are displayed on the board. This usually shows that any learned skill requires lots of practice
- "Mathematical Growth Mindset": described below
- "Community agreements and norms for discussion.": described below

We also planned additional interventions at strategic times throughout the year, to help students be successful in college and support their sense of belonging. Students attended brief one-on-one meetings with their instructor to check in about adjusting to life in college and balancing responsibilities. Instructors discussed with students opportunities like scholarship applications, work on campus, and student life events. Since most students are STEM majors, before the registration date for Spring semester, the science advisor visited the class to help students with the logistics of registration as well as identifying available major and career options.

## 3. ***Course Content and Pedagogy***

### 3.1. *Developing Mindset and Learning Community*

The specific components of our stretch calculus program — including mindset interventions, active learning through group work, and peer support — all align with the TRU framework to make the learning environment more equitable and inclusive in order to better support first-year students. The course starts with an introduction to the concept of growth mindset and collaborative group work via a group assignment which includes watching a video from Jo Boaler discussing mindsets. In small groups, students collaboratively create an illustration of one idea from the video that resonates with them and present their illustration in class (see illustrations in Figure 1). The illustrations are displayed in the classroom throughout the rest of the year for reference. The instructors then support their students' development of a growth mindset with careful language that values struggle, making mistakes, and students' unique thought processes.

We have found that at first our students are particularly reluctant to share their work for fear of being wrong and showing weakness in front of their peers. It is essential that our instructors work to normalize making mistakes, and to shift the understanding to one that views mistake making as an essential part of learning, beneficial to the entire classroom community. Facilitating this shift requires frequent, overt conversations about mistakes to challenge and change students' previously held beliefs about what mathematics is and what it means to be good at math. We begin this work by explicitly setting norms for their work as a learning community and create a large poster as a class for display on the generated norms. The instructor starts the discussion by offering a list of community agreements such as "Be a Learner," "Take Responsibility for Your Learning," "Appreciate the Contributions of Others," etc. Students review the list in small groups and discuss which of the agreements resonate with them, which they think will be easy to follow and which will be hardest. They also get a chance to suggest additional agreements they think are missing. The instructor then makes a revised community agreement poster including the class additions, which is displayed in the classroom and signed by students and the instructor (see center of Figure 1).

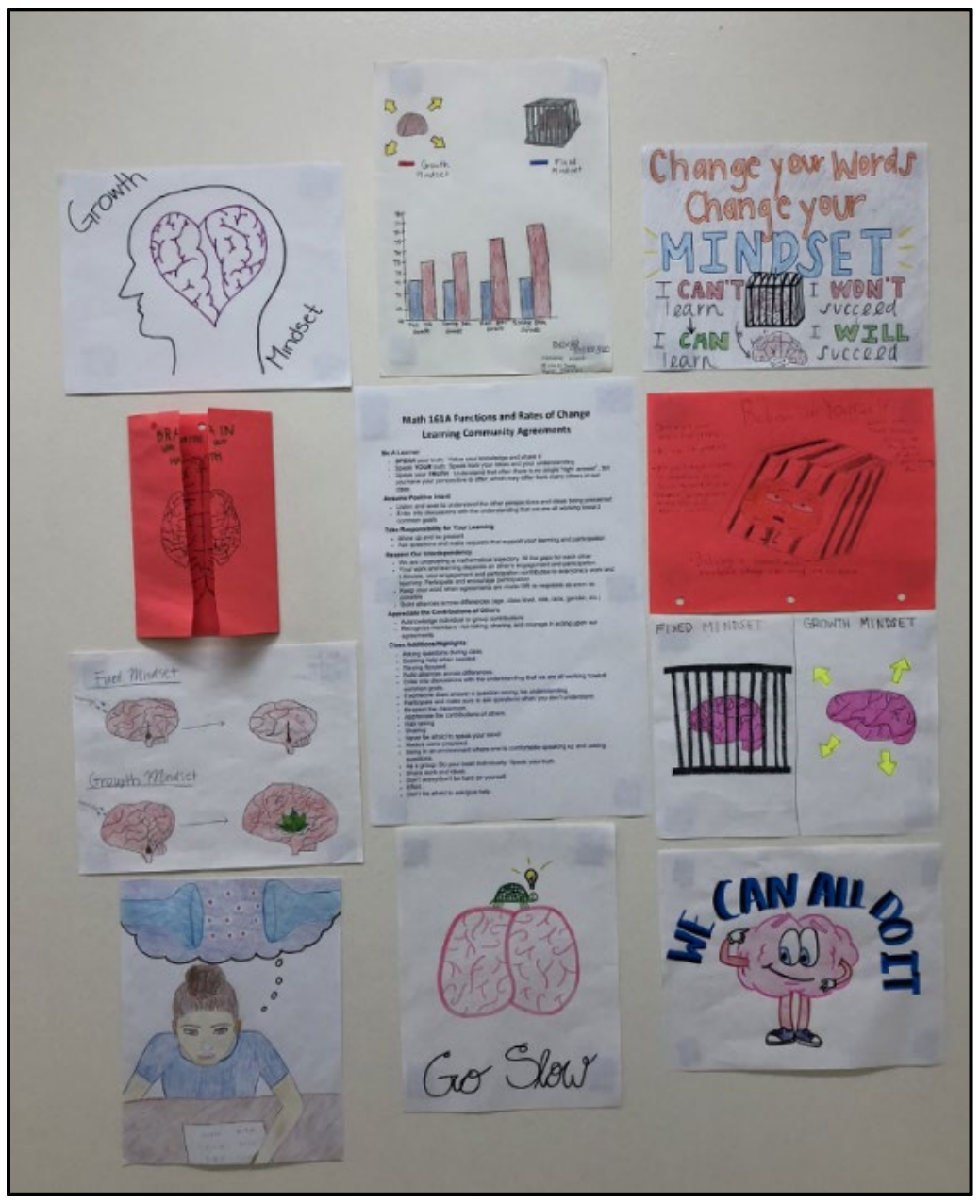

*Figure 1: Mindset illustrations and community agreements*

Active learning through group work dominates classroom activity throughout the year. We set an expectation of collaboration from the first day of class with group building exercises like building marshmallow towers and solving fermi problems. These activities break the usual pattern of problems done in a traditional math class where memorization and speed give some students an advantage while disadvantaging others. Instead, all students contribute to the group effort, there is no objectively right answer, and the activities value creativity and incorporate opportunities for revision. These principles are leveraged throughout calculus content as well.

Students are able to benefit from peer interactions through group work and their learning community led by an upper-division undergraduate student. This undergraduate student, the *Learning Community Mentor* (LCM) attends class with the students, setting an example and helping the instructor as needed, and leads a mandatory weekly learning session for students outside of the regular lecture time. These sessions are included as part of the course when students are registering, so students know before the semester starts what time commitment is required. The LCM is ideally a former student of the course and is selected by the course instructor for their ability and interest in fostering a collaborative learning environment. These students receive training and supervision from the university's tutoring program and are financially compensated for their time. Many of the LCMs are interested in careers in education.

### 3.2. *Calculus Content and Pedagogy*

Fostering a collaborative learning community supports students throughout the course, but this community would be less effective if we did not hold our students to the same rigorous expectations we have of all calculus students. Therefore, we start working with calculus topics from the start of the semester. We do not attempt to cover all the material students would see in algebra or pre-calculus; instead, we focus on the specific skills and concepts that support understanding of calculus. Beginning with discussions of rates of change, we review representations of functions using notation, graphs, tables, verbal descriptions, and equations. Topics are developed in the traditional order (following Active Calculus, by M. Boelkins available https://activecalculus.org/) using supplemental materials designed for group discussions. The first semester covers limits, the concept of the derivative, and differentiation rules up to the chain rule. The second semester continues with applications of the derivative, the concept of the integral, the fundamental theorem of calculus and integration by substitution. Right before moving from the conceptual understanding of differentiation to differentiation shortcuts, we review trigonometric functions. This unit starts with a field trip to the university's makerspace, where students make an electrically powered toy Ferris wheel and use it to graph the height of someone riding the Ferris wheel as a function of time. With this introduction to sinusoidal functions we reinforce the definitions of sine and cosine in terms of the unit circle. From this point on, trigonometric functions are part of the function catalog and are used throughout the rest of the year.

In the design of classroom activities, we kept in mind the principles of *Agency, Ownership, and Identity;* and *Equitable Access to Content* from the TRU mathematics framework. New course topics start with motivating examples that students can make sense of and elaborate on with their current understandings, even if they may not have the tools to fully investigate. The goal is to pique student interest and provide a need for new concepts and skills that will help to solve the problem. Throughout each unit, we want to give our students experiences where they can be successful problem solvers and where math makes sense to them, and provide them opportunities to make connections between ideas, learn from their mistakes, and revise their work to deepen their understanding. The following example is an outline of the story arc of one such unit.

### 3.3. *Example: Optimization - using derivatives to find maximum/minimum values*

At this point in the course, students have learned how to find derivatives of a variety of functions and they have used derivatives to investigate a variety of real world contexts. In this unit, students learn how to find solutions to real-world and mathematical optimization problems.

We start this unit with a hands-on activity. "What is the box with the largest volume that we can create by cutting out squares in each corner of a letter sized sheet of paper and folding up the sides?" Students construct a variety of open boxes using different colored pieces of paper for each size, visually inspect the classroom collection and make guesses as to which will have the largest volume. They then measure side lengths and compute the volume of the boxes. The class graphs the data of volume as a function of cut-out size. From the graph they estimate the size with the largest volume and then come up with a formula for the volume as a function of square length cutout of paper corners. Graphing the function shows that the largest volume coincides with a horizontal tangent line. The class determines that derivatives will be useful to find such a point.

From here we learn about methods of optimization including the first derivative test and the process of finding absolute maxima and minima of a function on a closed interval. At this point we revisit the open box example to find the absolute maximum of the volume function.

The culminating activity for this unit also serves as a major group project for the course: In groups, students design an “optimal” can of a given volume. The problem-solving process mimics the experience from the unit: They start by analyzing the dimensions and volume of a can and create other cans with the same volume, including making paper models. They then use calculus to find the can that minimizes the cost of producing a can of the same volume, given the material costs for the sides and the top/bottom of the can. They then present their work in a poster session (see Figures 2-4). The project includes three phases: (1) exploration, where students measure and create cans of a certain volume (2) analysis, where students use calculus to minimize the cost of the can, and (3) written report and presentation, where students create and showcase a poster to share their findings.

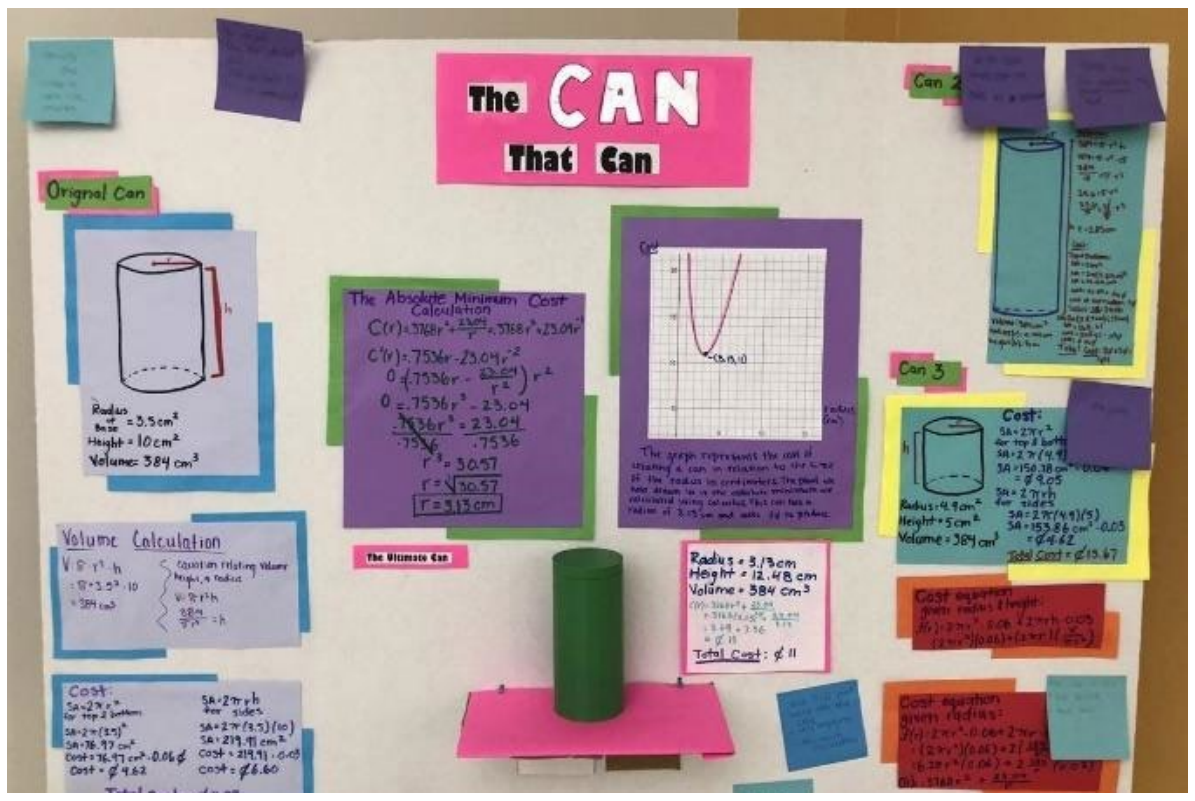

*Figure 2: Optimal Can Design Project*

While students often struggle through the process, in particular with the details of the calculus computations, their final results show creativity and initiative. Some students even used the 3D printer in the maker space (which the class visited earlier in the year) to make a model of their final can design.

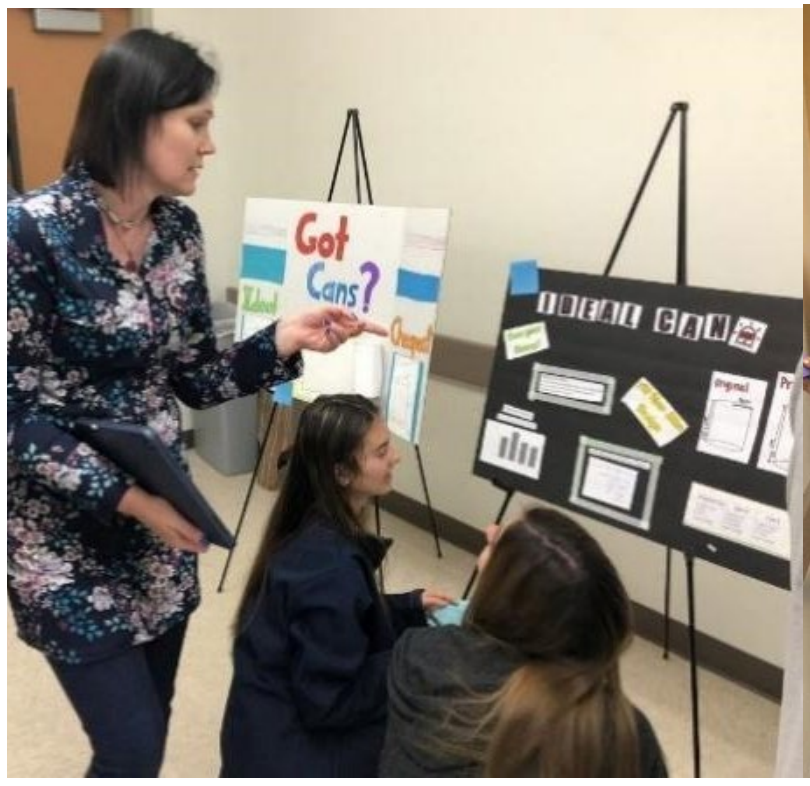

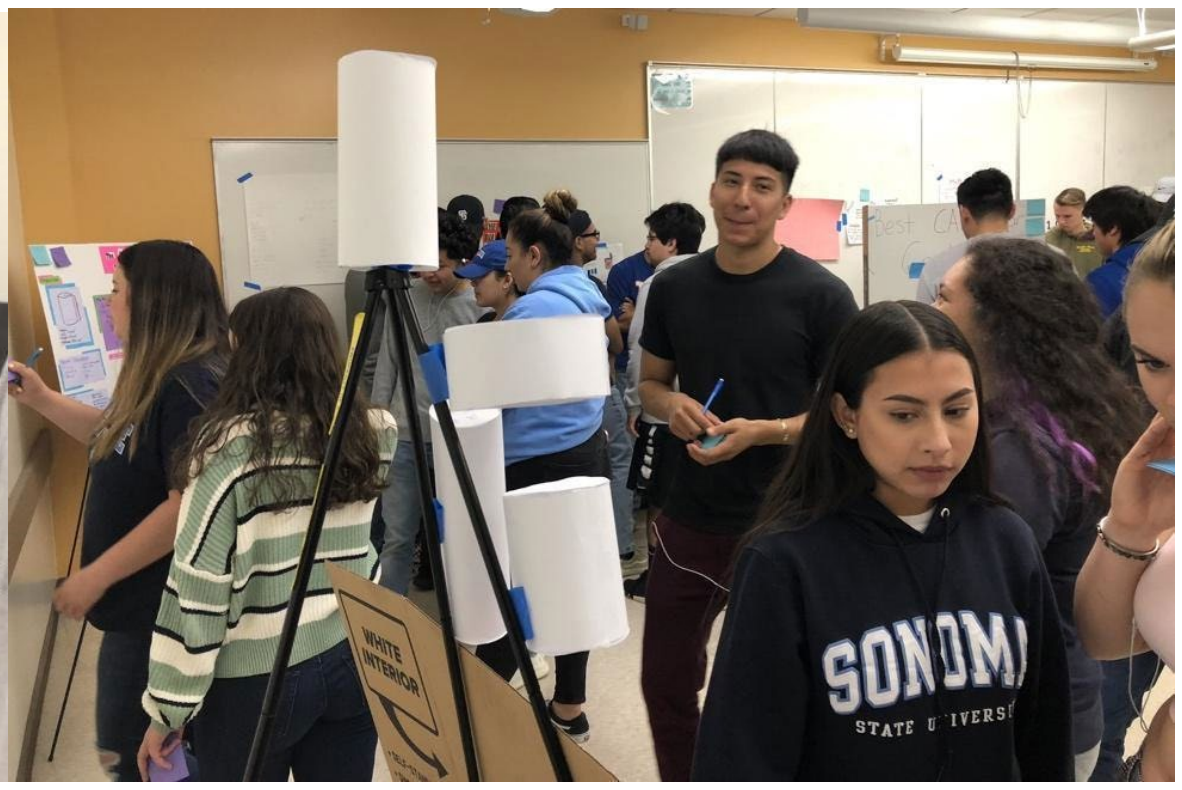

*Figure 3 (left): Instructor and students discussing a poster. Figure 4 (right): Poster presentation mini-conference with gallery walk*

All of the examples we describe, from developing mathematical growth mindset, to getting students comfortable with sharing their work, and to incorporating hands-on activities, illustrate how the relaxed pace of the course is particularly valuable. Instead of rushing through material for the sake of coverage, we can spend additional time discussing ideas in depth, increasing students' comfort level with the concepts, and developing students' computational skills.

### 4. *Administrative considerations*

For the first semester, students receive a grade of Credit/ No Credit, and assessments are a mix of group assignments, individual online homework using Webwork, and individual written homework assignments with opportunities for revision. The first semester includes two midterm exams, with opportunities and expectations to submit revisions. In place of a final exam, students complete a culminating group project and present their work during finals week. There are several benefits for the C/NC grades, some logistical and some pedagogical. Students who do not get credit often change their majors and do not wish to retake the class for a grade replacement. This grading scheme also removes some of the extrinsic grade motivation to allow students to focus on their learning, attitudes, and practices. Similarly, it also removes some of the pressure from the instructors and allows for more holistic grading.

The second semester includes a similar mix of assignments, with more emphasis on graded work that aims to prepare students for expectations in Calculus II, including a traditional final exam and letter grade. The second semester also includes the group project example described above on optimization that culminates in a poster presentation and written report.

### 5. *Impacts*

We designed this course through departmental conversations based on the belief that students' success in mathematics in general, and calculus in particular, is more closely linked to their mindset, confidence, and sense of belonging than to any particular skill set. This is not to say that particular mathematical skills and knowledge aren't important, but rather that those can be built over time and as needed if the student has the confidence and resilience necessary to do the work, and has genuine support.

Our program was specifically designed with the goal of closing GPA gaps, and increasing URM student participation in STEM courses and majors. Here we detail impacts of those measures for the students who participated in the stretch program in the 2018-2019 academic year (see Table 1). In particular, we highlight the success of students during their second semester in the stretch course, Spring 2019. This was the first year of full implementation of the calculus program, and due to severe academic disruption in 2020, we are not including reports of the stretch program from this most recent year.

*Table 1*

*Comparative DFW and GPA -*Students who did not complete either semester or did not enroll in the spring semester are not counted in the full year DFW rate for stretch calculus. Only the spring semester carried a letter grade.*

| Time Period | Course | DFW rate | Total Population | %URM | AVG URM GPA | AVG non- URM GPA |
|---|---|---|---|---|---|---|
| F12-F19[26] | calc I | 31% | 2789 | 32.6% | 1.89 | 2.17 |
| F18-S19 | calc I | 38% | 396 | 35.4% | 1.55 | 1.93 |
| F18-S19 | stretch calc, full year | 30%* | 58 | ≥47.5% | * | * |
| S19 | stretch calc, spring sem. only | 20% | 50 | 56% | 2.10 | 2.38 |

Of the 58 students enrolled in the stretch calculus course in the 18-19 academic year, 50 enrolled in the second semester, and 42 successfully completed the second semester. Within one semester, 30 had taken another math class, 21 of whom took calculus II. The DFW rate for the stretch students who took calculus II in the fall of 2019 was equivalent to the historical calculus II DFW rate at SSU. Additionally, the proportion of URM students in calc II increased substantially over historic and recent proportions (see Table 2).

*Table 2*

*Comparative URM representation in calc II*

| Time Period | Total Population in Calc II | %URM |
|---|---|---|
| F12-F19 | 1161 | 30.1% |
| F17 | 77 | 28.6% |
| S18 | 121 | 28.1% |
| F18 | 74 | 36.5% |
| S19 | 106 | 32.1% |
| F19 | 76 | 50% |

We do not have the data needed to determine if the increased representation of URM students in calculus II in Fall 2019 is specifically due to the 21 stretch students from the year before, but it is reasonable to assume that they contributed significantly due to the demographic breakdown of the stretch program and the abrupt change in URM representation in calc II. This

26 Our University System Data Dashboard allows us to look backwards from the most recent semester, so at the time of writing any historical data semesters through Fall 2019.

assumption that these students did have an impact in increasing representation, however, is further supported by noting that the only other semester presented here in which URM students made up more than a third of the calculus II student body was the first semester following the pilot of our stretch calculus program.

## 6. *Challenges and Recommendations*

Throughout our course development and implementation, we have confronted several obstacles and challenges. From these, we make recommendations for other schools attempting such a model.

The nature of the yearlong course can be challenging for schedules and student placement. Challenges also come from convincing others (e.g., faculty members, campus administrators) of the importance of the yearlong course particularly for students in majors with high unit counts such as engineering. We recommend creating advising materials and having ongoing conversations with STEM department advisors to get buy-in for their students' placement into stretch calculus, and to ensure that students do not encounter conflicts in the second semester with required major courses.

The peer support, offered through LCMs, is an essential component in these courses, and should be incorporated in any implementation of such a program. However, the peer mentors are also students and are sometimes not able to meet the needs of the course due to their own course obligations. We are working closely with the university's tutoring program to hire and support LCMs with ongoing training and supervision, and LCMs have priority registration to make their schedules work. In our pilot, LCM sessions were voluntary, and in the second year, we made LCM session attendance mandatory to provide more structured time for students to work on the course material. However, asking the undergraduate LCMs to productively hold the attention of up to thirty students for two hours was too high of an expectation, and for next year, we are planning to make one hour of attendance mandatory and the second hour voluntary.

Our program includes a co-teaching model that enables new instructors to the stretch course model to observe instructors and be trained to implement interventions in the year before they teach their own sections of stretch. This has been a key component of the success of the program, and we recommend implementing some professional development work for teachers of stretch courses.

At the same time, there are ongoing concerns about the financial sustainability of both the co-teaching model and the LCM model. So far we have been successful in working with the administration to support this professional development, and the work has also been supported by funding from the CSU system. Going forward, we are hoping that the data-supported success of the program will help to secure funding to continue our PD model.

Reflecting on the work of designing, piloting, and implementing the stretch calculus program, we found several additional experiences as essential. We developed this stretch calculus program as part of a wider change from remedial math to supported general education courses, and it was taken on by the entire department. It included a department-wide push to increase active learning and to create a more inclusive learning environment for all our students. Faculty collaboratively learned about supporting first-year students' transition to college and creating a culture of inclusion in two summer workshops facilitated by outside consultants which were not mentioned in the above. Many of the course interventions were designed as part of these workshops and ensuing collaborative course design. During the pilot and implementation of the courses, instructors frequently met (most at least every other week) to discuss successes and challenges and to support each other to grow and persevere as they changed their teaching practice. The work was, and continues to be, challenging but rewarding and it contributed greatly in moving us closer to a department culture of inclusiveness and belonging.

# 6 Redefining Success: A Multifaceted Approach to Assessing Departmental Change in First-year Mathematics/Statistics

**Judith E. Canner**

*California State University, Monterey Bay*

Dr. Judith E. Canner is a Professor of Statistics in the Department of Mathematics and Statistics at CSUMB. She is the Statistics Program Coordinator and past Quantitative Reasoning Assessment Coordinator for the campus. Her professional and research interests include expanding statistics and data science education, quantitative reasoning assessment, and data-driven pedagogical innovation.

**Jennifer Clinkenbeard**

*California State University, Monterey Bay*

Dr. Jennifer Clinkenbeard is an Assistant Professor Mathematics in the Department of Mathematics and Statistics at CSUMB. She is the current Quantitative Literacy and Finite Mathematics supervisor and coordinator, as well as incoming Quantitative Reasoning Assessment coordinator for the campus. Her research interests include mathematics education at the undergraduate level.

**Steven Kim**

*California State University, Monterey Bay*

Dr. Steven Kim is an Associate Professor of Statistics in the Department of Mathematics and Statistics at CSUMB. He serves students in introductory statistics, statistical theory, generalized linear model, simulations, and Bayesian inference. He studies statistical methods in dose-response modeling and collaborates with researchers in diverse scientific areas.

**Alison G. Lynch**

*California State University, Monterey Bay*

Dr. Alison G. Lynch is an Associate Professor of Mathematics in the Department of Mathematics and Statistics at CSUMB. She is the Calculus I course supervisor and the Math/Quantitative Reasoning representative on the campus's General Education committee. Her professional and research interests include student thinking in undergraduate mathematics and supporting students through mathematical transitions.

**Peri Shereen**

*California State University, Monterey Bay*

Dr. Peri Shereen is an Associate Professor of Mathematics in the Department of Mathematics and Statistics at CSUMB. She served as the Finite Mathematics coordinator for the 2018-2019 and 2019-2020 academic years. Her professional research interests include Lie algebras and representation theory, combinatorics, and voting theory.

**Alana Unfried**

*California State University, Monterey Bay*

Dr. Alana Unfried is an Associate Professor of Statistics in the Department of Mathematics and Statistics at CSUMB. She serves as the Introductory Statistics supervisor and coordinator. Her research interests focus on statistics education and consulting, including the development and validation of instruments assessing attitudes towards statistics.

**Jeffrey Wand**

*California State University, Monterey Bay*

Dr. Jeffrey Wand is an Assistant Professor of Mathematics in the Mathematics and Statistics Department at CSUMB. He is the Precalculus course supervisor and coordinator. His professional and research interests are Lie and representation theory, combinatorics, and he has even dabbled in some Bayesian statistics.

*Abstract: In 2018, the faculty in the Department of Mathematics and Statistics redesigned all general education courses and implemented support structures to increase student success. To assess whether these changes created equitable learning outcomes for our students, we worked with the Office of Institutional Assessment and Research to obtain the data on student grades in our general education and corequisite mathematics and statistics courses and associated demographic data. We also conducted surveys and focus groups to determine how both corequisite courses and our pedagogical practices impacted the attitudes of students towards mathematics/statistics, students' identity as mathematicians/statisticians, and their sense of belonging. While we observed no difference in GPA performance between first-generation students and others, we did observe ongoing equity gaps in student performance by gender, Pell eligibility, and race/ethnicity. The data collected on attitudes and experiences through surveys and focus groups allowed us to identify areas of improvement in our general education courses as we work to eliminate the observed gaps in equity.*

## 1. *Introduction*

In August 2017, the Chancellor's Office of the California State University (CSU) System issued an executive order that eliminated all remedial courses from the 23-campus system and required all first-year students to be placed in general education mathematics/statistics courses beginning in Fall 2018. The order was a major change for California State University, Monterey Bay (CSUMB), where 25-35% of first-year students had previously started in remedial mathematics courses. The Department of Mathematics and Statistics at CSUMB saw the mandate as an opportunity to redesign all of our general education courses (including Precalculus and Calculus I) and implement support structures to improve learning for all students. Redesign was particularly necessary as our general education courses would now include a highly diverse body of students in terms of their preparedness for college-level mathematics and statistics. We assumed that our existing general education courses were not serving our students as well as they could be, but we did not have processes in place for collecting, analyzing, and using data to make meaningful change. In the absence of systematic knowledge about our local data, we made curricular changes based on best practices in the research literature. We began offering corequisite courses in place of mathematics remediation to give additional support to any students requesting it. This case study discusses our approach to assessing this redesign during the 2018-19 academic year (the first year of following the executive order).

Prior to the curricular redesign, CSUMB required first-year students to take up to two semesters of non-credit bearing remedial algebra. Remedial courses increased students' time to graduation, and students who began in remedial mathematics were less likely to earn a bachelor's degree (Jaggers & Stacey 2014). When the CSU System eliminated all remedial courses, CSUMB created four corequisite support courses to pair with four general education mathematics or statistics courses. These support courses are taught by faculty members who also teach the corresponding general education course. Each support course meets once per week for two hours, focusing on developing mathematical skills and concepts needed for current course material as well as student success strategies. The support courses are open to any student in the corresponding general education course, and students self-select into the corequisite support course following a directed self-placement recommendation and/or a recommendation based on their high school courses, grades, and standardized test scores. Since the support courses are not required, monitoring students' course-taking behavior, as well as the academic outcomes of students who opt in and out of support, are a part of the Department's data collection and evaluation process as well.

In response to the existing equity gap shown in the systemwide and campus-wide data of graduation time and academic performance, we redesigned the general education courses to create more opportunities for equitable participation and student engagement through Reading Apprenticeship and Complex Instruction. Reading Apprenticeship is a pedagogical framework designed to apprentice students as readers in a specific discipline (Schoenbach et al., 2012). The framework attends to four dimensions: social, personal, cognitive, and knowledge-building. Complex instruction is a pedagogy initially designed for K–12 students that focuses on groupwork, defined as "students working together in a group small enough so that everyone can participate in a clearly assigned learning task" (Cohen & Lotan, 1997, 2). Groupwork in Complex Instruction has three defining characteristics to promote equitable participation by all members in a group: group-worthy tasks supporting multiple abilities; autonomy of the group through norms and roles; and attending to status through individual and group accountability (Cohen & Lotan 1997). We adapted both pedagogies for use in our college mathematics and statistics courses. CSUMB faculty members engaged in professional development around these pedagogies prior to and after implementation.

To adapt Reading Apprenticeship to the discipline, faculty focused on incorporating reading and reading strategies into class assignments and activities. For example, students were assigned a selection of text and in a five-column metacognitive reading log they were asked to record: a concept, their understanding before reading, new ideas and examples from the reading, and what they think these ideas and examples mean. After small-group and class discussion, they were asked to complete the final column by writing new understandings after the discussion.

To adapt the Complex Instruction framework, instructors first created group-worthy tasks for class time. Class time shifted to focus primarily on small-group discussion and student presentations of their findings. To facilitate these groups, instructors assigned roles that described a clear responsibility each student had to their group. (More details on Complex Instruction can be found in Cross-Cutting CHAPTER 3.)

Instructors also blended the two frameworks together; for example, using Complex Instruction group structures to facilitate small-group discussion of a reading log. Over time, instructors have revised their own teaching practices and understanding of how to best adapt each pedagogy in their mathematics and statistics courses.

To assess whether these changes created equitable learning outcomes for our students and to continue to improve our program, we identified a need to develop a systematic process for evaluation. The process we developed involves collecting and analyzing a wide variety of data (including grades, attitudes surveys, experience surveys, and focus groups) to create a substantive picture of what is working and where we can improve to better serve our students. By taking a multifaceted approach to assessment, we can understand which changes have been successful and for which subgroups. Our approach enables us to continually improve our first-year courses and progress toward our goal of equitable learning outcomes.

Our working definition of "equitable outcomes" is informed by the work of Rochelle Gutiérrez (2009), including both the dominant axis of access and achievement and the critical axis of power and identity. (This framework is described in detail in Cross-Cutting CHAPTER 3.) Our goal is to attend to the four dimensions of equity to shape our course design, implementation, and data collection/analysis. We address questions of achievement in our data collection by considering equity gaps in course GPA and pass rate by Pell eligibility, gender, race or ethnicity, and first-generation status. We examine how the creation of co-requisite courses and course redesign impacts access (to credit-bearing courses and opportunities to learn) through multiple measures including: time to successful completion of a course, course taking behavior, success in subsequent courses, and retention over time. We evaluate the identity component of the critical axis through the use of student questionnaires and focus group data to better understand students' sense of belonging, mathematical access, being "doers" of mathematics, and developing one's mathematical identity. Power is the one aspect of equity that we have not yet explored through data collection. However, pedagogical and design choices, including Complex Instruction, Reading Apprenticeship, and directed self-placement were implemented with student agency as a goal. While these considerations do not fully address issues of equity, using this framework as a guide for course design and data collection helped us to go beyond only questions of achievement.

## 2. *Using Institutional Data to Examine Achievement, Access, and Equity*

In order to evaluate the impact of our corequisite courses and our pedagogical practices on achievement and access, we worked with the Office of Institutional Assessment and Research (IAR) at CSUMB to obtain the data on student grades in our general education and corequisite mathematics and statistics courses and associated demographic data (e.g., ethnicity, Pell eligibility, SAT/ACT scores). Institutional data enabled us to not only evaluate our program, but also enabled us to compare our new program to the prior program of required remedial math courses before enrollment in general education courses. It also enabled us to track student success in subsequent courses in the Calculus sequence. In collaboration with IAR, we set up a data report that contains all the necessary information to support evaluation of achievement and access (as we have defined them), including grades in relevant courses and demographic data. For sequenced courses (e.g., Precalculus and Calculus I), we also used institutional data to evaluate student success in subsequent courses. We could then supplement the data collected at the institutional level with data collected within our classrooms.

## 3. *Using Surveys and Focus Groups to Examine Identity and Equity*

We approached the evaluation of identity through both surveys and focus groups to understand how both corequisite courses and our pedagogical practices impacted the attitude of students towards mathematics, students' identity as mathematicians, and their sense of belonging.

### 3.1. *Attitudes Survey*

At the beginning and end of each semester, we administered the Attitudes Toward Mathematics Inventory (Tapia & Marsh, 2004) to each mathematics GE course and the Survey of Attitudes Toward Statistics (Schau et al., 1995) in introductory statistics. The pre-surveys gave us a baseline picture of our student population and their attitudes toward the subject matter. Together with the post-survey, we analyzed shifts in student attitudes within each course. These shifts are especially relevant to the calculus sequence, as attitudes are known indicators for retention and success of students in college calculus (Pyzdrowski et al., 2013).

### 3.2. *Experiences Survey*

We developed an experience survey given at the end of each semester to measure student perceptions of the redesigned courses. Students were asked about interactions with peers, their instructor, class climate, and how beneficial course design elements such as reading assignments and group activities were for their learning. A companion experience survey was also given in the corequisite support course. The questions were adapted from the National Survey of Student Engagement (Kuh, 2001) and written by faculty.

### 3.3. *Focus Groups*

The final piece of data we collected to help us critically evaluate the changes we implemented to our curriculum were focus groups. Our campus Teaching, Learning, and Assessment staff conducted in-class focus groups in multiple sections of each course. In the focus groups, students responded to broad questions about the course, as well as targeted questions about in-class activities, reading assignments, support courses, and students' sense of belonging in class. The focus groups took place during class time so all students could participate. Students first discussed the questions in groups, then shared their responses with the class so others could agree or disagree. Since the focus groups were led by outside facilitators, students felt comfortable sharing their honest opinions. One challenge with the focus groups is that they require in-class time. Instructors were reluctant to give up 45 minutes of class time, especially for classes that only meet 160 minutes per week. As such, we only plan to conduct focus groups when we make major course changes.

## 4. *Preliminary Analysis*

We began to evaluate the equity in access and achievement in the 2018-19 academic year based on the narrow definitions of performance-based GPA and pass rate gaps between different groups in general education courses. For example, Precalculus students who completed remedial math courses had a pass rate of 50% compared to the pass rate of 82% for students who did not require remedial math prior to Precalculus. In contrast, under the new changes, the pass rate was 76% for Precalculus students who completed and passed the corequisite course, and it was 82% for students not enrolled in the corequisite support course. With the elimination of remedial courses, the results seemed remarkable because some students are enrolling in Precalculus up to a year earlier than they would have previously. We used a regression analysis to evaluate various forms of equity gaps (first-generation status, race/ethnicity, gender, and Pell-eligibility) while adjusting for the semesters (fall or spring) and the courses (the four GE courses). We observed no difference in GPA performance (4-point scale) between first generation students and others ($p = 0.653$). However, we observed ongoing gaps in student performance for other group comparisons. For example, male students showed a lower average GPA for their class compared to female students (difference = 0.17 points, $p = 0.018$). We also observed that Pell-eligible students had a lower average course GPA compared to students who were not Pell-eligible (difference = 0.24 points, $p = 0.001$). Finally, when we compared to white students, there was a lower average GPA among African Americans (0.56, more than one half of full letter grade), Hispanics (0.43), Asian Americans (0.15), and Others (0.26) ($p < 0.001$).

We also tracked the first cohort into Calculus I and found that students who enrolled in the support class for Precalculus performed lower in Calculus I than their peers who had not enrolled in the support class even when controlling for their Precalculus grade. The use of institutional data to observe the continued performance of students in the Calculus sequence allowed us to discover that, although the support class is getting students through Precalculus, students may need additional support to be successful in Calculus I. The benefit of the analysis of institutional data is that it can be replicated each semester quickly to identify equity gaps. The analysis of institutional data does not tell us *why* we observe the gaps in GPA performance and pass rates in subsequent courses. The data collected on attitudes and experiences through surveys and focus groups allow us to identify areas of improvement in our general education courses as we work to eliminate the observed gaps in achievement and access.

Based on the attitudes survey, we found that students who enrolled in Precalculus during the 2018-19 academic year reported higher self-confidence and enjoyment of mathematics after completing the course when compared to before completing the course. In contrast, the experiences survey showed that although more than three-quarters of Precalculus students indicated that they benefited from working with peers during class and that in-class activities helped them to learn the material, but only 57% of Precalculus students rated their overall course experience as “good” or “excellent.”

Informal student feedback and results from the experience survey indicated that the shift toward active learning made students feel as if they are not learning as much as they would in lecture-based formats. Recent studies show that students’ *feeling* of learning may actually be lower in an active environment, even though their actual learning is higher (Deslauriers et al., 2014). We were able to quantify the pushback using the survey, noting that by the end of the semester, pushback diminished and students appreciated the active learning approach. Regardless, we recommend transparency with your students to remind them regularly of the benefits of active learning and why we structure our courses in a particular way.

From the focus group summaries, we gained a clearer picture of how students perceived their course and what aspects of the course they felt supported (and hindered) their learning and sense of belonging. In Precalculus, students across sections remarked on the value of group work in supporting their sense of belonging. One student noted that, “I feel like I really did belong in this class because we were in our groups and I felt comfortable asking questions because there was like only two other people that you were asking. You did not have to ask in front of the entire class, you know. So it was really comfortable.” Another group expressed that, “We feel like every student in this class can help each other, feel like it's a positive learning environment.” However, students who did not feel a sense of belonging might not feel comfortable saying so publicly, so more focused data collection is needed in this area.

We are just beginning to gain a fuller picture of how our corequisite courses and pedagogical practices impact equity in achievement, access, and identity. We recognize the limitations of our measures, but we believe that even imperfect measures are valuable if they help us improve our approach and our courses.

## 5. *Final Reflections and Recommendations*

As we reflect on the last two years of our evaluation process, it has been critical to have support from all levels of our institution. Our data collection processes have involved students, instructors, course coordinators, and staff from other offices on campus. We also received funding from the department and university for data collection and analysis. To meet the challenges of ongoing data analysis, we collaborate with student researchers, who receive pay or credit, to help us clean and analyze the data, a great service and time-save for the faculty and a valuable experience for the students.

Our work also benefited from the collegial and collaborative nature of our department and its faculty. In each key area (e.g., curricular changes, leadership, data collection, data analysis), faculty with relevant expertise made recommendations, then the department made decisions by consensus. We recognize that decision-making is more contentious in some departments, which would cause this work to be more challenging.

One limitation of our analysis is that we cannot distinguish between the impact of the switch to a co-requisite model and our improved pedagogical methods since these treatments were applied at the same time and are confounded. For example, it is possible that the pedagogical changes alone led to the improvement in pass rates for Precalculus students needing additional support, and it is also possible that the corequisite course structure alone led to this improvement. However, the additional sources of assessment that we used, such as experiences surveys and focus groups, provide some evidence that both the corequisite model and pedagogical innovations have both led to improved outcomes, although we cannot tease out exactly how much credit to attribute to each feature.

We also note that our analyses thus far have been fairly segmented, looking at institutional data, attitudes data, and experiences survey data separately. Moving forward, we plan to look more holistically at these data. By analyzing these pieces together, we may be better able to understand how our courses are serving a diverse student body and how we can improve.

As we continue to improve our data collection and evaluation processes, we have discovered numerous areas for improvement that we would recommend departments attend to in their planning.

First, we neglected to gather data from faculty teaching the newly redesigned courses. Though we gathered qualitative data informally through faculty meetings, a more formal process would have provided more immediate realization of the need for ongoing professional development for current and new faculty instructors. We are now working to provide regular professional development within our department, especially for new instructors.

Second, our data collection does not currently capture the Power dimension of equity, despite our attention to this dimension in the curricular redesigns. Although we recognize this as a limitation, we have not yet found an effective strategy to assess this dimension. We continue to seek ways to measure power, especially in terms of student ownership of their learning and ownership of the curriculum.

Finally, we do not yet have efficient processes for sharing data with course instructors and coordinators for the purpose of making further curricular revisions. Our goal is to develop a set of metrics that support the evaluation of student success, using multiple metrics, and other information that will support our decision-making in course revision. We plan to move to automated reports, both during the semester to easily identify struggling students and at the end of the semester to facilitate continuous improvement efforts. As we move towards a balance between manual data analysis and automated data analysis for quantitative data collection (institutional data, the attitude surveys, and the experience survey) and ongoing qualitative data collection (free response from the experience survey and focus groups), we hope to more appropriately and quickly respond to any need for change within our courses to support student success and close equity gaps.

# 7 Addressing Diverse Equitable and Inclusive Issues in a Four-year Private Hispanic-Serving Institution

**Cheng Chang**
*Mercy College*

Dr. Cheng Chang is an Assistant Professor in mathematics at Mercy College. He is the coordinator for the Smart Math Initiative since 2016.

**Zhixiong Chen**
*Mercy College*

Dr. Zhixiong Chen is a professor in math and computer science at Mercy College. He is the director of the cyber education center since 2019 and was the Chair of the department of Math & CS at Mercy College from 2014-2019.

*Abstract: We introduce a comprehensive initiative at Mercy College, a four-year private Hispanic serving institution with a large percentage of commuter student population, to address diverse, equitable, and inclusive issues in first year College Algebra classes. The three main components of the initiative: classroom embedded peer tutors, course redesign, and course coordination are explained and shared with implementation ideas. We then present some assessments in the progress of student performance, including the reduced disparity between different ethnicity groups and genders.*

## 1. *Background*

Mercy College is a private four-year college in the New York metro area, with total student enrollment of 10,577 in 2019-2020. The College is a Hispanic Serving Institution with 70% women students. The enrollment of calculus classes has been consistently low despite a relatively large population of STEM majors, moreover, the women students' ratio in those classes is significantly less than 10%. To put it into context, there is a large percentage of first-generation college students, and students from low-income families in our school. Many students entering the College often have not had a sound background in foundational mathematics, essential to success in STEM. In fact, faculty report that many students, including those in the target population who plan to pursue STEM, drop out of their majors because of insufficient preparation as well as difficulty in absorbing concepts, especially when delivered in traditional lecture formats. Compounded by the fact that commuter students in the College made up 60% of the population, a concern shared by faculty is that students in general do not have a strong sense of community. This lack of inclusiveness is more likely being significant for the limited number of women students taking calculus classes.

We recently organized a math symposium[27] on redesigning math for student engagement and success, to share experiences among some regional institutions of similar type. In our own case, we identified three most challenging points related to diverse, equitable, and inclusive (DEI) issues. First, the majority of our students entering the College lack the math skills essential to success in Calculus, as basic as College Algebra. Based on our admission data, most students did not have any AP class experience in high schools, many of them did poorly in high school math classes as well as standardized tests like the New York State Regents Exams, some others had a multi-year gap after graduating high schools. Most STEM students at the College must start from College Algebra before they can move on to Precalculus and then Calculus. Secondly, we noted that many students who planned to pursue STEM dropped out of their majors because they could not obtain the required grade in College Algebra to move forward. This is one of the main reasons that the enrollment of calculus classes is consistently low. Lastly, a large portion of basic math classes are taught by adjunct instructors, the turnover ratio among them is not small. Due to limited resources at the College, most adjunct instructors do not have permanent office space or desks on campus. Most of them teach for multiple institutions and usually have a very heavy teaching load accompanied by busy commuting, so they usually come and go rapidly. One such consequence is that many students do not have equitable learning experiences and some students are still struggling with College Algebra even in the calculus class.

[27] https://mercycollege.prod.acquia-sites.com/liberal-arts/math-symposium

In order to have a robust and successful calculus sequence, we need to redesign our basic math curriculum and build a strong foundation for students, so that they can channel through to the calculus class. Especially for our underrepresented population like first-generation college students or women students in the STEM fields, they may not have received the attention, support, and resources in their high school education, so the academic preparedness gap for them is detrimental in college math courses and could well be the reason for eventual college drop off.

## 2. ***Smart Math Initiative***

The name of our program to address DEI issues is called Smart Math Initiative. It is a combination of math courses redesign, introduction of in-class peer student support, and departmental course coordination. The Initiative has two main goals. For one, it is designed to set an expectation in teaching pedagogy and provide a curriculum standardization across all basic math classes, to diminish the vastly different teaching styles among adjunct instructors, guarantee the rigor and consistency of the course contents, and accommodate each student with an individualized and self-paced study plan. For another, the Smart Math Initiative aims to provide students, especially underrepresented students like women, an abundance of support and encouragement inside and outside of classrooms. There are three main components: embedded peer tutors, controlled emporium, and coordination.

### 2.1. *Embedded Peer Tutor (EPT)*

The Embedded Peer Tutors (EPTs) are in-class student teaching assistants who have previously taken and received an A in a math class equivalent to or higher than College Algebra. They are currently enrolled undergraduate students of any major, as early as second semester in the freshman year. Unlike working in tutoring centers, EPTs work directly with students inside classrooms, in conjunction with and complement to the instructors.

We usually hire 10 to 20 EPTs each semester, depending on their schedules and availability. Each semester, all the EPTs must go through a training before they work inside the classrooms. By the nature of the emporium model which will be elaborated on the next component, a class is usually divided up into small groups or even individuals and have separate discussions or learning activities from time to time. Therefore, EPTs are supervised and directed by instructors to assist students, in the form of mini lecturing a learning topic, facilitating a group discussion, or helping students self-reflect learning strategies and form a good working habit.

We found, through the end of semester student evaluations, that the emotional support and encouragement provided by EPTs are indispensable to peer students, especially the students who have a lot of fear in learning mathematics. The role of EPTs in creating an inclusive learning environment is irreplaceable by instructors, as many students stated in the course evaluations that they felt less embarrassed and were more willing to ask questions to EPTs and they felt that the ways of explanations by EPTs are easier to understand and are more relatable to students.

Interestingly, the majority of our EPTs are women, most of whom are STEM major students. In Spring 2019, we hired 13 EPTs, out of which 10 are women. The female EPTs ratios in Fall 2019 and Spring 2020 are 8 out of 14 and 6 out 8, respectively. One potential explanation is that the women students' percentage in the College Algebra classes is larger, which is consistent with the overall College population.

In the post surveys conducted among the EPTs, the majority found the working experience as being enjoyable and beneficial to their college career development. Many EPTs eventually enrolled into a 5-year master program in teaching education degrees in STEM.

From a financial point of view, EPTs receive decent stipends for their work, which is an indication of their importance and the appreciation from the program. The close involvement of EPTs in the classrooms, as well as their achievements in helping peer students, have strengthened the inclusiveness and community bonds among students, and hence could help improve the retention rate.

### 2.2. *Controlled Emporium Model in Smart Math Lab*

The emporium model has been explored for many years across the states, we make some small tweaks in implementing to fit our unique needs. For instance, a significant portion of our student body is nontraditional college students, some have full-time jobs outside the school while maintaining a student status. We accommodate these students by using a combination of engaging lecturing (in contrast to no lecturing in a conventional emporium learning environment) and self-paced active learning inside the classroom, so that the students have the means to continue their learning outside the classroom on their own.

All the Smart Math Initiative classes take place in the designated computer labs, called Smart Math Labs (Figure 1), to create an equitable and regulated classroom environment. Actually, when we first piloted some classes in the traditional classrooms and asked students to bring their personal computers, a lot of unexpected issues occurred, for instance, some students did not have personal computers and used cellular phones instead, which became a huge distraction in learning. These labs have open floor concepts, with a maximum capacity of 30 students. The design of the labs is structured to welcome class discussions, nurture a support environment, and create bonds among students.

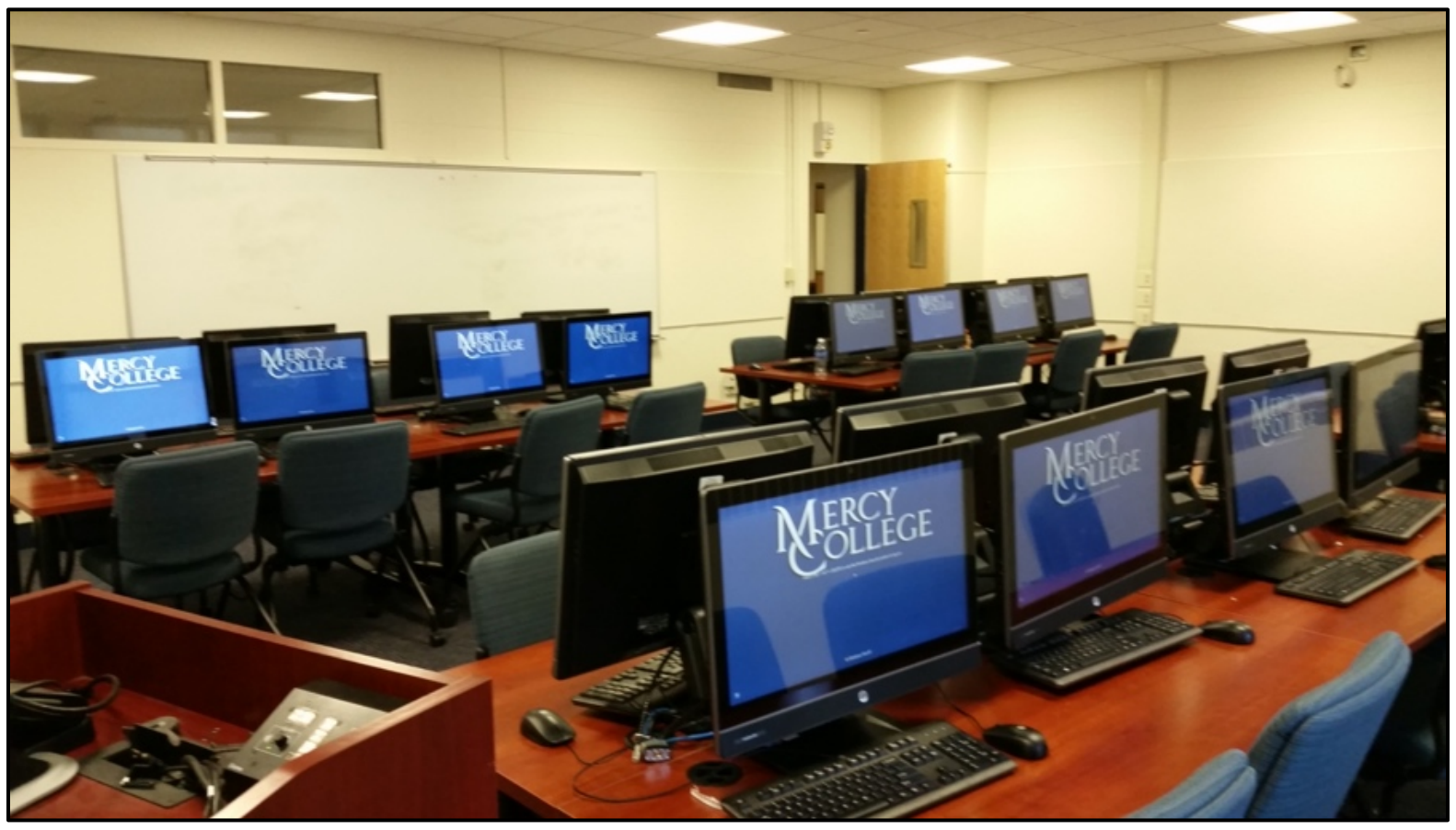

*Figure 1: Smart Math Lab*

In terms of the course contents, they are divided up into modules, delivered through a learning platform like MyOpenMath or ALEKS. Outside the classrooms, students' progress through the program at their own pace, by a standardized benchmark timeline. With self-paced learning, students tend to dive into different topics without realizing the connection and the big picture of them, and this will be addressed inside the classrooms.

In a typical 2-hour class meeting, an instructor will motivate and teach between half an hour and an hour, on a specific topic outlined in the benchmark schedule. The lecture will provide additional examples as well as connections between different topics. The rest of the class meeting will be individual interactions and/or small group discussions between the instructor and the students, assisted by EPTs.

In order for our controlled emporium model to be successful, we spent a lot of effort in workshops and training for instructors, especially part-time instructors. Following closely the Instructional Practices Guide (Abell et al., 2017), some high standard and consistent teaching materials including well designed group discussion worksheets are developed and updated regularly for all instructors. These efforts are made to guarantee that students receive maximum support, and instructors, full-time or part-time, have the same standard in delivering the courses.

## 2.3. *Coordination and Technology*

The College designated a full-time faculty within the math department as the Initiative Coordinator, since the college-wide initiative is recognized as one of the key factors in the student success and retention. Moreover, the Initiative in its infancy requires constant monitoring. The key responsibility of the coordinator is coordinating with different departments and there is one course release for the administrative duty each semester. The coordinator meets with the institutional research for data and performance analysis, works with the center for teaching and learning on workshops and class observations to support adjunct instructors, and also works with the tutoring center to recruit and supervise EPTs.

The textbooks, videos, and other readings for the courses are all available electronically online, conveniently accessible through mobile devices as well as computers. Instructors receive electronic reports on students' weekly performances, the coordinator monitors the learning statistics across sections. A typical weekly report includes the topics that the majority of the students are struggling with, the total number of hours that students spent on the learning platform, and the average grade of the weekly assignment. Through reading the reports prior to class meetings, instructors can identify the weaker students and help them more efficiently during class discussions.

Finally, the coordinator maintains an online forum, to create a virtual community among instructors. A shared site for instructors and EPTs is regularly updated with teaching notes, pedagogical articles, extra exercises and detailed solutions, etc.

### 3. *Assessment*

One of the key metrics in measuring the success of the Smart Math Initiative is by the completion rate. The two courses redesigned are Math115 (Math for Liberal Arts) and Math116 (College Algebra). Our focus in this article is on Math116, as it is required by STEM majors and its population is the potential pool for future calculus students.

We collected the data on the Math116 completion rates for the most recent five different cohorts. By a cohort, we mean all first-year students enrolled in fall semester, excluding continuing and transfer students. We track each cohort up to three years, since some students decided not to take Math116 in the first year or second year, some of them waited until the third year or beyond. The passed percentage is based on the number of students who took and passed Math116 with a grade of C or better, out of the total number of students in a cohort.

Figure 2 shows the trending of the first year passed percentage. We believe that the Initiative, which began in Spring 2016, could be one of the driving factors in the stead increase until Fall 2017. There are other supplementary actions (also implemented around Spring 2016) contributing to the increase in the first year passed percentage, along with the Initiative. To name a few examples, we eliminated pre-requisite requirement for Math116; emphasized the importance for students taking Math116 in their first semester; and we also adopted an OER textbook and homework system, to alleviate financial burden on students. On the other hand, there is a slight dip in the first year passed percentage in the Fall 2018 Cohort, it could be explained by another initiative we implemented in that year, to improve the course rigorous and knowledge retention.

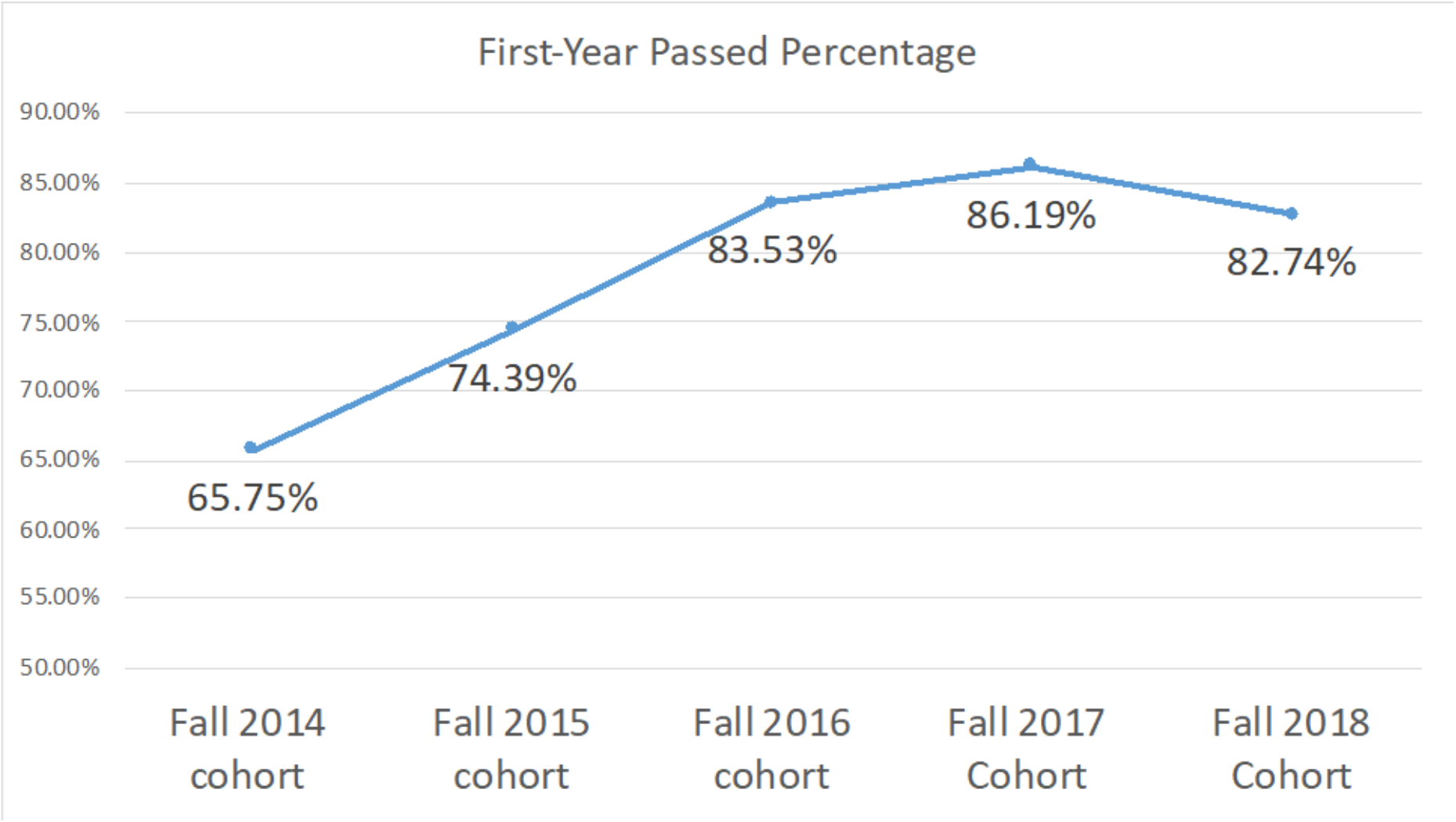

*Figure 2: MATH116 First-Year Passed Percentage*

Next, we want to examine how the EPTs help the classes in passing rates. Table 1 compares the pass rates of Math116 in Spring 2019, categorized by sections with EPTs and without. There is a clear contrast in the DFW rates, and the sections with EPTs significantly outperformed those without EPTs.

We must point out that none of the online (modality) sections had an EPT, and historically speaking, online sections usually have much higher DFW rates than in-person sections at the College. On the other hand, we received overwhelmingly positive feedback in the EPTs program, from students (through course evaluations), EPTs (through open response surveys), and instructors (through informal surveys and solicited feedback). These have strengthened our belief that EPTs can make significant impacts in student learning.

*Table 1*
*Spring 2019 Math116 EPTs Analysis*

| Math116 Numbers | Course Grade A-C | DFW | Math116 Percentages | Course Grade A-C | DFW |
|---|---|---|---|---|---|
| **Not-EPTs** | 36 | 58 | Not-EPTs | 38.30% | 61.70% |
| **EPTs** | 198 | 118 | EPTs | 62.66% | 37.34% |

Regarding the performance of different ethnicity groups, Table 2 provides the consolidated grade distributions within each group for the last six consecutive fall semesters. To make a direct comparison, we normalized the number of students in each group to 1 and focused on the DFW students. In the absolute ideal scenario, where all ethnicity groups perform equally well, the percentage of each group among the DFW students will be 25%. Figure 3 displays the distributions, and it clearly indicates that since the Initiative fully implemented in Spring 2016, the disparity among different ethnicity groups is reducing. In other words, the percentage of each ethnicity group among the DFW students is for the most part monotonically reaching closer to 25%.

Table 2 also shows a noticeable decrease in overall passing rates in Math116, but it does not necessarily indicate a negative outcome. For example, we eliminated pre-requisite requirement for Math116 since Fall 2016, thus, significantly more students were able to take Math116 in their first semester. As a result of that and the Initiative, both the absolute number of students passing Math116 with a C or better as well as the first year passed percentage in Math116 were increasing until Fall 2018 (Figure 2).

*Table 2*
*Math116 consolidated grade distribution within ethnicity groups*

| **Fall 14** | Asian | Black | Hispanic | White | Overall | **Fall 15** | Asian | Black | Hispanic | White | Overall |
|---|---|---|---|---|---|---|---|---|---|---|---|
| A-C | 96.77% | 70.75% | 77.58% | 90.45% | 80.68% | A-C | 90.32% | 66.67% | 73.50% | 91.37% | 78.84% |
| DFW | 3.23% | 29.25% | 22.42% | 9.55% | 19.32% | DFW | 9.68% | 33.33% | 26.50% | 8.63% | 21.16% |

| **Fall 16** | Asian | Black | Hispanic | White | Overall | **Fall 17** | Asian | Black | Hispanic | White | Overall |
|---|---|---|---|---|---|---|---|---|---|---|---|
| A-C | 82.76% | 59.54% | 69.09% | 83.77% | 70.41% | A-C | 87.50% | 71.67% | 76.20% | 83.10% | 75.83% |
| DFW | 17.24% | 40.46% | 30.91% | 16.23% | 29.59% | DFW | 12.50% | 28.33% | 23.80% | 16.90% | 24.17% |

| **Fall 18** | Asian | Black | Hispanic | White | Overall | **Fall 19** | Asian | Black | Hispanic | White | Overall |
|---|---|---|---|---|---|---|---|---|---|---|---|
| A-C | 75.00% | 66.67% | 69.70% | 77.56% | 70.46% | A-C | 83.33% | 64.68% | 67.84% | 75.71% | 69.41% |
| DFW | 25.00% | 33.33% | 30.30% | 22.44% | 29.54% | DFW | 16.67% | 35.32% | 32.16% | 24.29% | 30.59% |

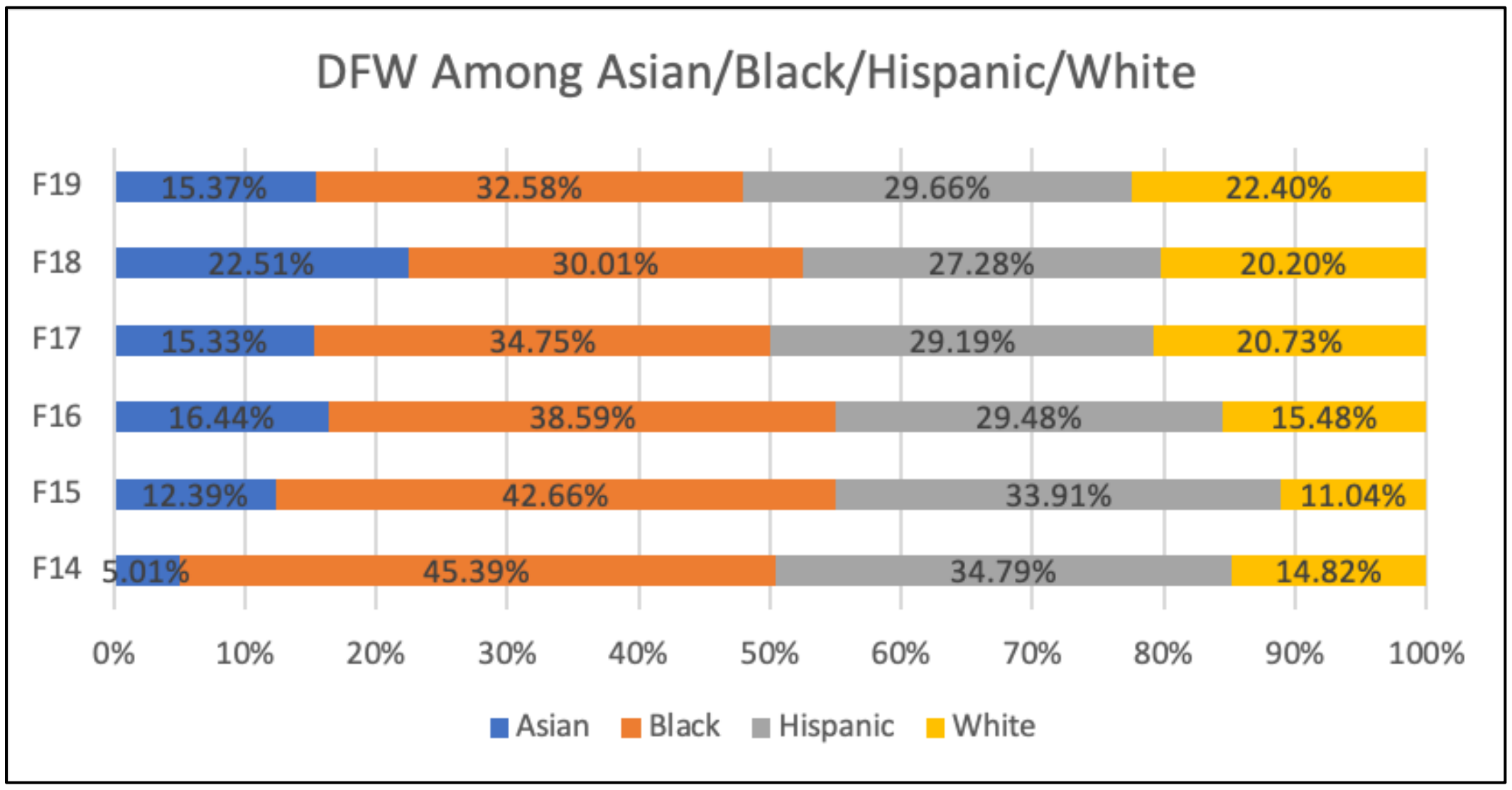

*Figure 3: Math116 DFW rates among Asian/Black/Hispanic/White*

Similarly, Table 3 lists the consolidated grade distributions within each gender over the past six fall semesters, and Figure 4 displays the direct comparison between genders. For the reference, in an ideal situation where male and female students perform equally well, the percentage of male or female among DFW students should be 50%. We observe in Figure 3 that the percentage of females among DFW students is gradually shrinking below 50% and stayed below 50% thereafter.

*Table 3*
*Math116 consolidated grade distribution within genders*

| **Fall 14** | Female | Male | **Fall 15** | Female | Male |
|---|---|---|---|---|---|
| A-C | 49.51% | 50.49% | A-C | 48.43% | 51.57% |
| DFW | 51.68% | 48.32% | DFW | 55.71% | 44.29% |
| **Fall 16** | Female | Male | **Fall 17** | Female | Male |
| A-C | 50.36% | 49.64% | A-C | 50.23% | 49.77% |
| DFW | 49.10% | 50.90% | DFW | 49.29% | 50.71% |
| **Fall 18** | Female | Male | **Fall 19** | Female | Male |
| A-C | 53.49% | 46.51% | A-C | 50.99% | 49.01% |
| DFW | 42.11% | 57.89% | DFW | 47.75% | 52.25% |

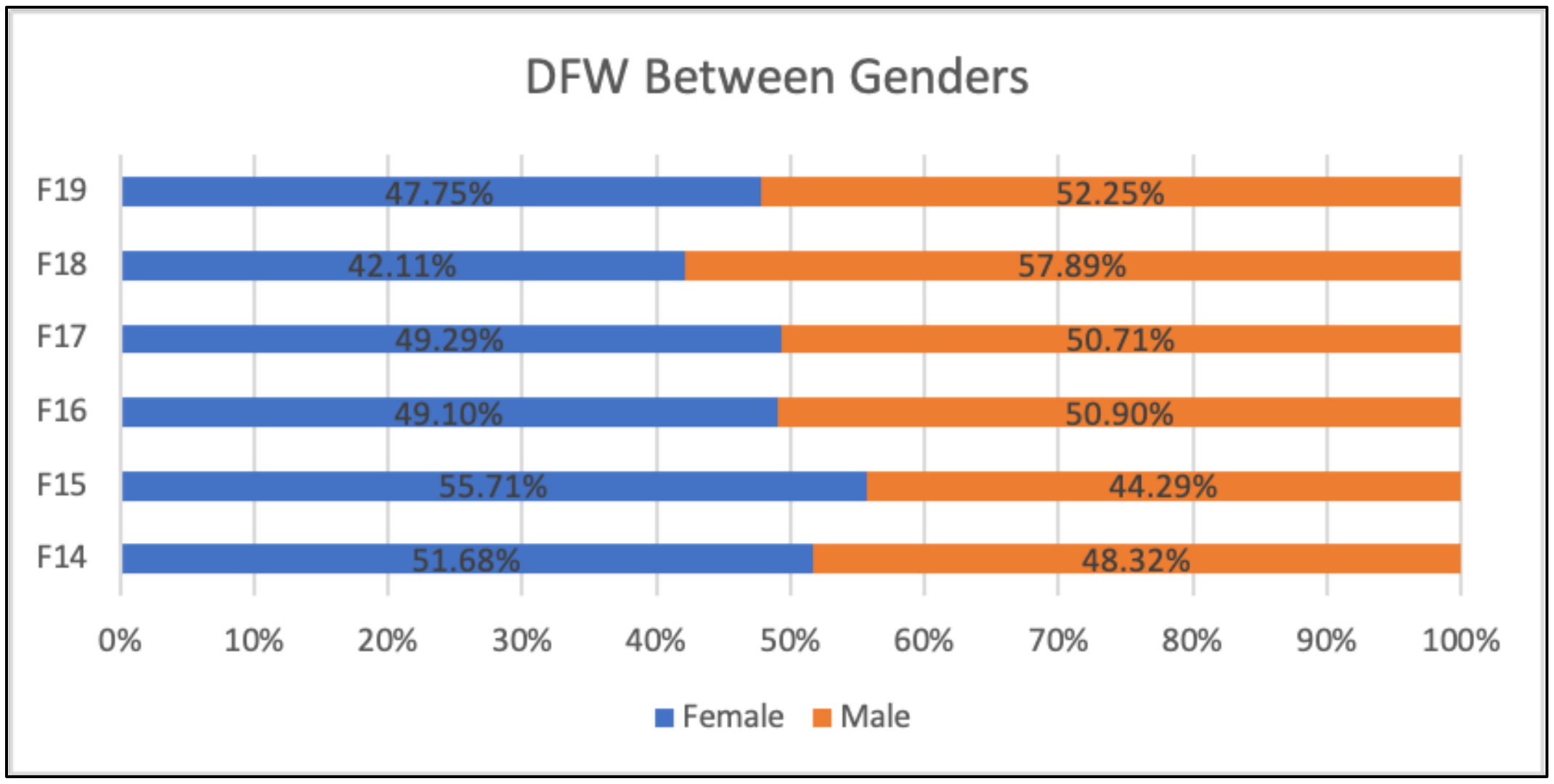

*Figure 4: Math116 DFW rates for male and female students*

## 4. *Implications and Recommendations*

Since we started the Smart Math Initiative, the EPTs program has always been a core component. For students, they not only look up to the EPTs as role models, but also, they see in themselves that they can be as successful as the EPTs. Another unexpected outcome in EPTs is that a large percentage of them are women students, around 75% in two of the most recent three semesters.

In the beginning of the Initiative, the expenses of EPTs were funded by an external grant through the Department of Education. After the grant ended, the College picked up the expenses from the college budget. Fortunately, some of our EPTs qualify for the Federal Work-Study Program, so a portion of the expenses are covered by the federal financial aid. Going forward, the College is going to charge a course fee for all students taking Smart Math Initiative classes for the sustainability of the program. The majority of the course fee will be used for the institutional license for the commercial learning platform (ALEKS), so that students will have access to it from day one at a highly discounted price. The rest of the course fee, which consists of approximately $15 per student, is still not fully enough to cover the expenses of EPTs, but the College will make up for the deficit.

The main capital investment for the Smart Math Initiative is converting conventional classrooms into computer labs. We have multiple campuses spanning three different locations, and the College converted six classrooms total across campuses averaging two per campus. Currently, the tests and exams component of the Smart Math Initiative courses is on paper, to avoid cheating on computers, however, we are considering asking the IT department to utilize the lockdown browser feature in the future, so that the tests and exams can move on to computers as well. Tests digitation not only provides students instant feedback on the test performance, but also helps the department monitor the progress across all the different sections.

Despite the initial success of the Smart Math Initiative, there are areas we are actively seeking to improve in the future. For one, we noticed that some students successfully passed a Smart Math course but struggled in subsequent major courses. They seem like not remembering how to apply things learned in the past or completely forgetting. For another, some students prefer conventional lecturing over the controlled emporium model, we want to serve this group of students while not going back to the old model. There are ideas about tackling these issues and the prospect of the Smart Math Initiative has become more exciting than it has ever been before.

# 8 Supporting Underrepresented Minority Students in STEM Through In-Class Peer Tutoring

**Tara C. Davis**
*Hawaii Pacific University*

Tara Davis is an Associate Professor of Mathematics at Hawaii Pacific University. She is interested in active and inquiry-based learning, undergraduate research, mathematics of peg solitaire on graphs, and geometric group theory.

*Abstract: In this article we describe our use of undergraduate students as in-class peer math tutors in inquiry-based classrooms. We aim to elevate students who are representative of our institution's primarily female, racially and ethnically diverse student body to positions of power by presenting them as role models and experts. We aim to break down barriers to students who need help, therefore promoting equitable access to tutoring for all math students. We will discuss impacts of this tutoring on student performance in their math classes, as well as provide student reactions to the intervention. We will also describe the practice in detail so that other institutions would understand how to implement a similar program.*

## 1. *The Context of the Program*

This article describes an intervention of supporting students in STEM through in-class peer tutoring created to make progress on several ongoing issues that we have perceived at our institution. The first issue is that our students are not utilizing the peer tutoring opportunities available to them. Our school, like most, has free peer tutoring available to students, but often students needing help do not actively seek out tutoring. By bringing peer tutors into an active learning or inquiry-based math classroom, we aim to break down this barrier. Another issue is that our students are not succeeding at desired pass rates in our math classes, ranging from developmental math through mid-level classes including Calculus (Hagman et al., 2015) and Linear Algebra. Therefore, one objective of this practice is to improve student pass rate and grades.

The third issue, which relates most directly to issues of student equity, inclusion and diversity, and will therefore be the focus of this article, is that of student belonging and availability of appropriate role models. Although the student body at our institution is primarily female and is racially diverse, (we generally rank in the top 10 schools in the Nation for student racial diversity) our math faculty is entirely white and North Asian, which is in line with the overall faculty demographics which are more than 80% white or Asian. Therefore, this program promotes diversity of role models that students have access to by elevating women, and well as Hispanic, Native Hawaiian, African American, Native American and Alaska Native students as peer tutors. These are all groups that are historically underrepresented in STEM. By providing role-model tutors, we aim to encourage the enrolled students to improve their mathematical skills and confidence in the class, but further, to encourage the enrolled students to persist. The practice also supported the peer tutors financially through an NSF grant as well as academically by providing course credit, and by elevating them to a public role of expertise in mathematics.

We hoped to encourage underrepresented STEM students to graduate and go on to work in the STEM fields by providing role-model tutors. The instructors were also supported by the teaching assistance, for example, instructors were able to spend more time helping each student since they did not need to attend to as many students at once, and the peer tutors helped foster and facilitate a positive active learning environment by engaging their peers with math activities during class time.

## 2. *Description of the Program*

### 2.1. *Development of the Program*

We will begin with the history of the program. Our institution has four concentrations within the math major, one of which is a Math Education concentration. This is more like a pure math major, and does not include an education dual degree or teaching certification, but it does include one required math tutoring class. This class, Math Education Practicum, is open to students from all concentrations, and it focuses on providing students with practical classroom tutoring experience, alongside the investigation of issues of teaching and learning mathematics. Originally, the student tutors worked in developmental math labs. Organizational restructuring and changes in leadership over the years lead to the eventual removal of these developmental math labs. Placements were required for the math tutors enrolled in the Practicum class. This led to the practice of utilizing the tutors during class time, rather than in an outside of class tutoring lab.

There is a second dimension to the development of this program. In addition to requiring placements for the Math Education majors enrolled in the Practicum class, we required employment for our students funded by the Louis Stokes Alliance for Minority Participation (LSAMP) [IOA-LSAMP NSF grant #HRD 1826864]. At the same time as the developmental math labs were being phased out, we were redefining the scope of student work under the LSAMP grant. Students meeting the criteria were eligible to be paid for tutoring work, and these LSAMP students constituted our second pool of potential in-class peer tutors. To qualify for LSAMP, students must be American citizens or permanent residents, maintain a 3.0 GPA, have declared a major in a STEM field, and be an underrepresented minority (African American, Alaskan Native, American Indian, Hispanic or Pacific Islander). The main objective of HPU's LSAMP is to support Native Hawaiian and other Pacific Islander students.

## 3. *Current Status of the Program*

In the current iteration of this program, students are assigned to work with a math faculty mentor as peer tutors in active learning, or inquiry-based math classrooms. It should be noted that on rare occasions there are some in-class peer tutors working at our institution who come from the Math Practicum pool, who are from groups over-represented in STEM (e.g., Japanese males. However, all LSAMP students, and the vast majority of the non-LSAMP students, are women, Hispanic, African American, Native American, or Native Hawaiian). For example, in 2018, the Math Practicum class had five students enrolled, three of whom were Pacific Islander LSAMP students, one was Hispanic, and the other was international. Another semester the class was over 65% Native Hawaiian. The program as it relates to this article will focus on the use of all student tutors, not exclusively on those from any racial or ethnic category, though as noted above, the majority of our peer tutors belong to groups that are historically underrepresented in STEM.

The author, who oversees the peer tutoring program and serves as the LSAMP Campus Coordinator, would send an annual solicitation the math faculty, looking for those who wish to support students, and who agree to utilize active learning pedagogy (so that the tutors have work to do during class). Those faculty wishing to supervise and mentor a tutor were provided with one. Student tutors were selected through a formal application process for the LSAMP program, and to be eligible to participate, the student must be from a racial or ethnic group the NSF identifies as being underrepresented in STEM. In addition to the requirements for the grant, the student should have earned a high grade in the class they would be tutoring. Students were trained either through the Practicum class, or by their faculty mentor; oftentimes the tutor already had a relationship with their mentor, and had taken the class they would be tutoring so they were familiar with the material and teaching methodology. The tutors were paid an hourly stipend. The author coordinated the employment with the HR and grants office, and placement with a faculty mentor, who managed the tutor's daily activities and assignments. There are approximately 5-10 student participants in the program per year, where each tutor is able to reach dozens more students through their daily tutoring work.

The tutors worked for several different math instructors, in classes from developmental math and pre-calculus, through mid-level calculus and linear algebra, to upper-level classes like proof writing and abstract algebra. The student tutors met with their faculty mentor prior to the start of the semester to have all expectations and duties laid out, and the tutors request any needed aides, such as advanced copies of the class handouts, access to the class Blackboard page, or solutions to the class problems. The tutors were told what to focus on during class time, and were encouraged to ask questions of the course instructor about their duties or the course material, including math questions. The tutors who had taken the Math Education Practicum class, which was the majority of tutors, had training on facilitating student problem solving, asking guiding questions, and focusing on keeping the work with the students.

The students enrolled in the math classes who were receiving in-class tutoring primarily worked in groups at the white boards (Davis, 2019) instead of listening to a lecture during class time. A typical class period would consist of a brief introduction and motivation by the instructor, which would include the posing of a problem or set of problems for the class to work on. Then time would be given for the students to work on the problem together in previously assigned small groups at their designated white board. During this time, the peer tutor(s) would be sent to listen to the discussions and interact with the small groups as needed. Oftentimes the enrolled students would raise their hands or call for the tutor by name. The instructor would also be available to help the small groups. Therefore, students were given individual attention each class period and would have the opportunity to ask questions in their small group to either the instructor or the peer tutor. The students often did ask questions of the tutors, who could then relay them to the instructor, who could reconvene a whole class discussion.

The program serves different math populations by focusing on the needs of the enrolled students and the faculty mentor. A developmental math class, for example, will focus more on students practicing routine problems, and the tutors may be required to help students check their answers, or answer questions in small groups. In a calculus or linear algebra class, the enrolled students may be engaged in a more complex task, which could take a large portion of a class period to complete, and in this scenario the tutor would be expected to facilitate small group discussion by listening to and giving advice on strategy, in addition to helping with basic computations and formulas.

## 4. *Evidence of impact towards diversity, equity, or inclusion (DEI)*

In this section, we will report on the progress that we were able to measure towards achieving the stated goals of this project, as well as additional evidence toward the impact on DEI. We emphasize that the data set collected is small scale, and more data collection would be needed to draw generalizable conclusions. Since the time that this project was implemented, our institution has continued to recruit students and professors to participate in our program. We have expanded to using students from outside of the math major, from computer science and engineering, as role model peer tutors in math classrooms, and these peer tutors provide the students with insights into the ways that the math will be applied in other majors. There was a significant struggle with continuing this program during the COVID-19 pandemic, although we attempted to use the tutors via breakout groups in zoom at the onset, we have since returned to normal learning conditions.

### 4.1. *Survey Data*

Informal intake and exit survey data in Calculus II, using 2 tutors for 21 students, and Linear Algebra, using 1 tutor for 20 students, was collected during Spring 2018. The tutors were one white female, one Hispanic female and one Native Hawaiian male. The surveys were created by the course instructor and not tested for validity or reliability. In terms of student population, in the Linear Algebra class,50% of the students were female, Alaskan Native, American Indian, Hispanic or Pacific Islander; in the Calculus class, 38% of the students were female, Alaskan Native, American Indian, Hispanic or Pacific Islander. Of the Linear Algebra students, 29% had been in a class previously that utilized peer tutors; 26% of the Calculus students had been in such a class before. We found that 85% of all students taking post-survey were encouraged by the peer tutors to continue in STEM; 96% were encouraged by the peer tutors to graduate; 88% viewed the peer tutors as role models. This was encouraging in light of our goal to increase availability of appropriate role models.

In terms of data about how much time students spent in tutoring, Linear Algebra students spend 2.94 hours on average during the semester outside of class with the tutor; Calculus students spend 0.8. This raises a question of why the students do not seek more out of class tutoring, when they indicate that they find in-class tutoring helpful.

On free response questions, the students responded in an overwhelmingly positive way to the tutors. Students appreciated the ability to have another perspective: “If I didn't understand the concept the way the teacher taught it, the tutors explained it a different way that helped me understand better.” The students working in small groups did not have to wait long to get attention when they were stuck, simply due to there being multiple helpers: “Having (the peer tutor) in class is great because I can ask questions when (the instructor) is working with other people". One student's response explicitly addressed gender: "They made it clear that people (ladies in particular for this class) can make it through STEM degree plans". There were a few quotes that hinted at the students comparing themselves to the tutor: "They inspired me to become better at math" and "It made me want to be as good at math as them!" as well as at persistence: "They really encouraged me to keep going." While it seemed clear from the free responses that the tutors were helping the students to ask questions and engage with the math, if we were to collect more data, we would have liked to have survey questions that specifically addressed belonging.

The peer tutors also completed intake and exit surveys, which reflected very positive attitudes. All of the tutors taking surveys had been enrolled in the faculty mentor's classes before as students when peer tutors were employed. In describing the impact that tutoring in this class had on them, one tutor stated, "I enjoyed helping people understand concepts and it has been interesting for me to see the students thought process." All tutors agreed that the tutoring improved their math skills, confidence, and impacted their decision to continue to major in math. As one tutor said, "Teacher roles can be led by students. Gives hope and shows that everyone is capable of doing math." Since 2018, one of the three peer tutors used during spring 2018 added a math major to her oceanography major and completed undergraduate math research, presenting at multiple national conferences and having a paper accepted to a peer-reviewed math journal, and the other two graduated and become secondary math teachers, one receiving a Master's degree. We have also had multiple students from classes who received peer tutoring go on to become peer tutors themselves.

### 4.2. *How peer tutoring can serve as an equity and inclusion practice*

We feel that this intervention is in line with the suggested equity practices given by Elllis Hagman (Ellis Hagman, 2019) because we are, using real-life applications, proofs and numerical methods of (Wawro et al., 2013) and (Greene, n.d.). In addition, we offer the in-class undergraduate peer tutoring as an equity and inclusion practice meant to assist student achievement. As indicated by Gutiérrez, (Gutiérrez, 2008), peer tutoring promotes equity and inclusion by showcasing the excellence already achieved by the peer tutors and providing enrolled students access to math experts who may resemble themselves more than their professors do. The tutors can impact the students' achievement by providing themselves as an example to look up to and to emulate in order to achieve academic success. It gives power to the women, Hispanic, and Native Hawaiian groups, by providing status to the tutors as experts and role models and hopefully may empower the students in those same groups to feel pride and hope for their future in math. Bringing tutoring into the classroom also creates an expectation that all students need help and should access tutoring on a daily basis; this is equitable because it makes explicit the belief that everyone can succeed in math with work and support, and that it is normal to struggle and need help as well.

## 5. *DFW Rates*

In the Section below we will go into more detail regarding data collected during Spring 2018 in Calculus II (21 students), and Linear Algebra (20 students) classes. The choice to use mid-level math classes for the intervention rather than developmental came in part from the fact that there is an overwhelming amount of data suggesting that underrepresented students are leaving the STEM fields not only in the developmental math classes but in those mid-level classes such as Calculus (Hagman et al. 2015). Moreover, Linear Algebra is a transitional course, being one of the first non-computational math classes that a student experiences. Students from a variety of STEM fields including Math, Engineering, Computer Science and Oceanography routinely take this class at our institution. Data from our institution’s Math Program Review shows that Linear Algebra at our institution historically has an almost 20% DFW rate on average over the last five years, and this rate was as high as 31.8% in AY 2013-2014, suggesting that this is indeed a high-risk class.

In the Calculus II class, the DFW rate was 33.3%, which is unfortunately higher than has been the average (21.3%) or the median (20.14%) over the last five years. We do not necessarily take this to mean that the tutors were not effective, as other factors may have come into play, such as a weak class, or the fact that this was the first time the instructor had taught this class in several years. The Linear Algebra class fared much better, with only a 10% DFW rate. Again, we cannot attribute this entirely to the peer tutoring, as the instructor was very experienced and the students were quite strong, but we take it as a positive sign. We are not attempting to make any overarching conclusions based on this small dataset, but we do provide it so as to have something quantitative in addition to the fact that the instructors and students have anecdotally been very satisfied with the intervention.

## 6. *Pre and Post Assessments*

### 6.1. *Calculus II*

We used the Learning Assistant Alliance, LASSO, (‘Learning Assistant Alliance’, n.d.) instrument to do a pre- and post- test to see how the Calculus student's math content knowledge (of pre-calculus) has changed, and we compared our data to national data. The LASSO pre-calculus inventory seemed to reflect improvements from the pre- to the post- assessment, suggesting that student’s pre-calculus skills, and perhaps by extension, their overall math skills improved upon taking the course. There were 9 students completing both surveys, and who took more than 300 seconds and answered more than 80% of the questions. The results of the pre- and post-assessments are shown in the table below and measured by Cohen's d of 1, which indicates a large measure of student improvement.

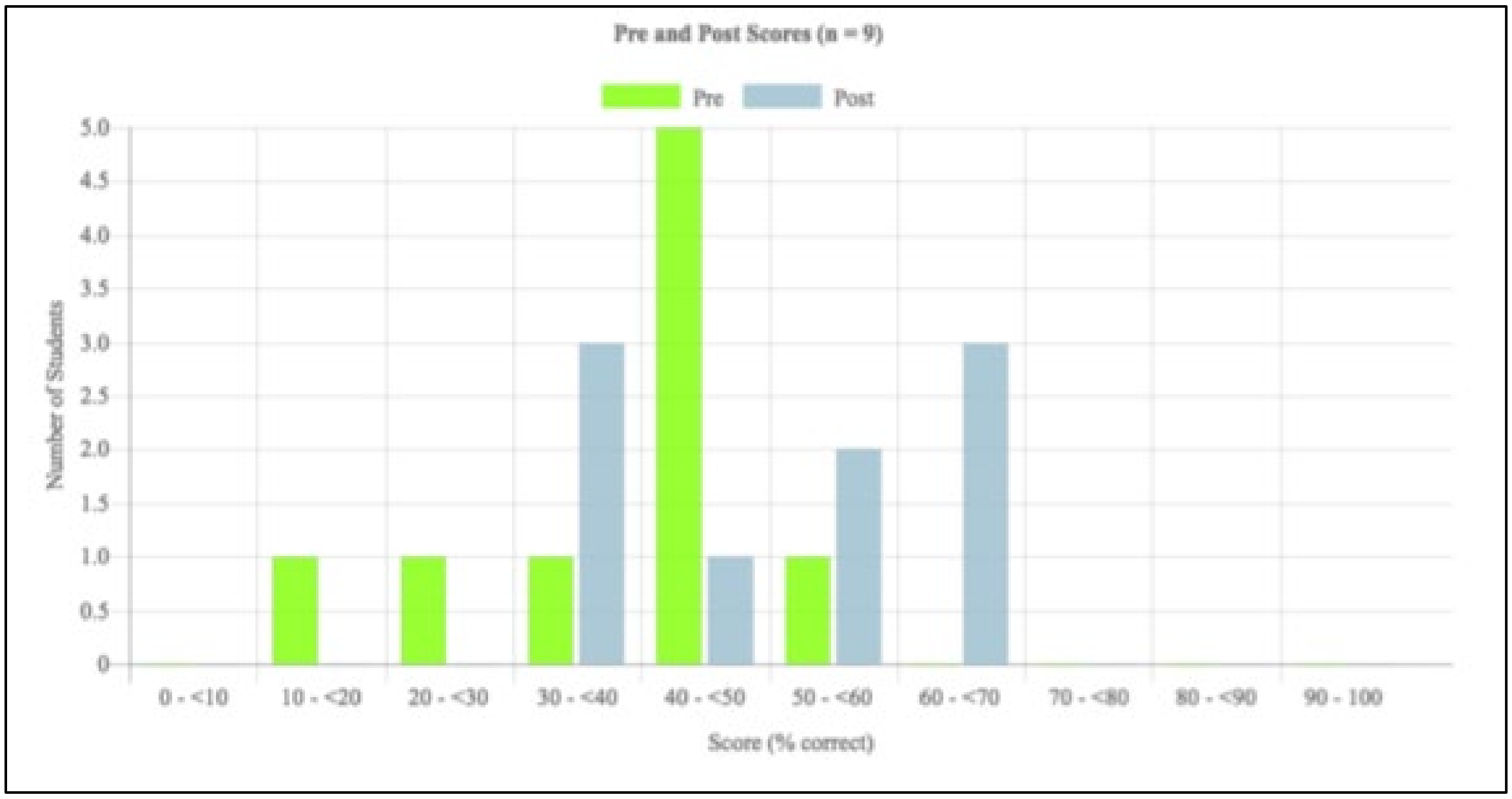

*Figure 1. Pre and post scores bar plot*

We recognize the limitations of this metric, given that only approximately half the class completed both surveys, and that it is a measurement of pre-calculus skills, not Calculus II skills. Moreover, it is a distant proxy to measure the effect of having peer tutors present in class; however, we use it as part of a triangulation of the overall impact.

## 6.2. *Linear Algebra*

We used as a post-test in the Linear Algebra class and can compare it with national data collected as part of the Teaching Inquiry-Oriented Materials: Establishing Supports (TIMES)[28] program as well (Haider, 2018) and (Wawro et al., 2013). The post-test was developed by researchers to compare linear algebra content knowledge between traditional and inquiry-oriented classrooms in schools across the country. The Linear Algebra students did not complete a pre-test, as they were not expected to have any previous knowledge; however, the role-model peer tutor did complete both the pre- and post- test. There were 12 of the 20 students who completed the post-test. These students did coincidentally have lower than average grades compared to others in the class on the final test and overall. The class average on the assessment was 63.4%. This is lower than the national averages for TIMES(Inquiry-Oriented Linear Algebra) students (67%), higher than average non-TIMES schools (57%) but not lower than all of the TIMES classes surveyed in (Haider, 2018). Therefore, some of the weaker students in the class still scored higher than those students across the nation who had traditional instruction. The tutor improved from pre- (55.5%) to post- (66.6%) assessment, though his score still seems low.

28 see https://times.math.vt.edu/ for details about Teaching Inquiry-Oriented Mathematics: Establishing Supports (TIMES) grant

## 7. *Implications and recommendations to other programs.*

Possible obstacles and barriers to implementing a similar program at other schools would include lack of instructor pedagogical freedom or large class size. Also, recruiting, managing, supporting and training tutors could pose a difficulty, depending on student population and institutional culture. The author checked in regularly with any tutors who were not working under her as well as with their faculty mentors. This program was implemented at a school with a very small math department (approximately six faculty members) and this management activity could become more onerous at larger schools. Another obstacle towards implementation would be finding compensation for all tutors.

Any institution that has math faculty applying active or inquiry-based learning in the classroom would be ripe for implementing this practice. The proposed duties of the tutors included an emphasis on peer tutoring in the classroom, but could be supplemented by peer tutoring in the tutoring center, holding office hours, grading, holding test review sessions, typesetting classroom materials in LaTex, or other teaching-related duties as assigned by their faculty mentor. As noted earlier, certain additional resources are available at our institution for the successful conduct of the project, namely, many students have already been trained to work as peer tutors by completing the Math Education Practicum class and therefore there is a pool of student workers already available. Any institution with a Math Education concentration would do well to recruit from this population.

# 9 Case Study of an Emerging Scholars Program at a Predominately White Institution: Analyzing the Success of Underrepresented Students in an Inquiry-Based Calculus Sequence

**Jessica Deshler**

*West Virginia University*

Jessica Deshler is Professor of Mathematics in the West Virginia University Department of Mathematics and Assistant Dean for Graduate Studies in the Eberly College of Arts and Sciences. Her research in mathematics education focuses on equity in mathematics and graduate student development.

**Edgar Fuller**

*Florida International University*

Edgar Fuller is Professor of Mathematics and Director of the Center for Transforming Teaching in Mathematics at Florida International University. He studies student learning and engagement in mathematics while trying to improve access to mathematics at all levels. He is currently co-PI on the NSF Catalyzing Change in Calculus project working to develop active and engaging coursework in calculus that leads to better access to STEM careers for all students.

**Tim McEldowney**

*West Virginia University*

Tim McEldowney is a postdoctoral researcher in the West Virginia University Center for Excellence in STEM Education. His work studies how critical transition points serve as barriers for minoritized students in STEM and how equity programs can support students through these transitions. He is currently PI of an NSF funded project to study how student knowledge about the graduate school application and admissions processes acts as a barrier to earning advanced degrees in mathematics.

**David Miller**

*West Virginia University*

David Miller is Professor of Mathematics and Assistant Director of Undergraduate Studies in the School of Mathematical and Data Sciences at West Virginia University. His research in mathematics education focuses on students' understanding of mathematical proof, professors' assessment of students' proof productions, and students' understanding of their Professors' proof feedback. He is currently co-PI on the NSF Louis Stokes Alliance for Minority Participation project for the Kentucky-West Virginia Alliance.

**Matthew Pascal**

*Point Park University*

Matt Pascal is a Professor of Mathematics at Point Park University in Pittsburgh, PA. He earned his Ph.D. from the Department of Mathematics & Statistics at American University in 2003 and has worked for equity in mathematics education for over two decades.

**John Stewart**

*West Virginia University*

John Stewart received his Ph.D. in physics from the University of Illinois, Urbana-Champaign in 1994 and is a Professor of Physics and Astronomy at West Virginia University. His work attempts both to further knowledge about how physics classes produce understanding and how physics departments and universities can work more effectively to serve and retain all students.

*Abstract: The Emerging Scholars Program (ESP), over the last three decades and across the country, has incorporated active learning strategies resulting in students (in particular, underrepresented students) earning at least a half of letter grade higher than other students in Calculus (Hsu et al., 2008). In 2009, West Virginia University (WVU) adapted the ESP in Calculus I in an effort to support the success and retention of underrepresented STEM students by embedding group and inquiry-based learning into a special section of Calculus that targeted underrepresented and first-generation students. We anticipated that supporting students in courses in the calculus sequence would build a community of underrepresented Calculus learners and serve as a support mechanism to provide a strong foundation for retaining underrepresented STEM students through graduation. In this chapter we analyze the success data of the students that have progressed through the WVU ESP Calculus course and compare them to their non-ESP counterparts. The results show that underrepresented ESP students succeed in the Calculus sequence at a substantially higher rate than underrepresented non-ESP students. In addition, underrepresented ESP students succeed at similar levels (slightly better but not significant) as students that are not underrepresented in the ESP classes. Finally, underrepresented ESP students graduate at similar rates to students not underrepresented in ESP and at significantly higher rates than underrepresented students in non-ESP. We believe that ESP is a viable option for institutions that want to make an impact for retaining and graduating underrepresented STEM majors, help them overcome obstacles and barriers in STEM, and increase diversity, equity, and inclusion in Calculus.*

## 1. *Introduction*

The Emerging Scholars Program (ESP) concept, originally developed by Uri Treisman in 1977 at the University of California at Berkeley, arose out of an intent to increase diversity within the student population (Treisman, 1990). Treisman observed that a core problem facing students in their academic work was a lack of community with shared experiences around developing and understanding mathematical ideas. Without this sense of community students can experience isolation that poses significant barriers to success. ESP implementations build communities of students around these shared experiences and identities in order to mitigate this sense of isolation and have been shown to increase success among underrepresented minority (URM) students (Asera, 2001; Hsu et al., 2008). We present here a case study of how we developed and implemented an embedded ESP into a calculus sequence at West Virginia University (WVU). As West Virginia's Land-Grant Institution, WVU is committed to producing the bachelor's, master's and doctoral level STEM graduates that meet West Virginia's workforce needs by educating nearly 30,000 students. The state of West Virginia is fully located within Appalachia; the student population of WVU is predominately white with 8% of students being URM and 13% being first-generation students (the first in their families to earn a four-year degree). For this study, students are classified as URM if they identify as Black, Hispanic, Native American, or Hawaiian/Pacific Islander. In this chapter, we describe our motivations for implementing a program for URM students at our institution, obstacles we faced and the modifications we made to the structure and targeted recruitment, as well as the successes we saw, with some suggestions for overcoming obstacles that may be encountered during implementation in a similar predominately white institutional setting.

## 2. *Motivations*

Two events facilitated the creation and implementation of an Emerging Scholars Program (ESP) at WVU. First, WVU received a National Science Foundation Louis Stokes Alliance for Minority Participation grant as part of a ten-institution alliance in Kentucky and West Virginia (KY-WV LSAMP). This grant provided a mechanism for the creation of additional programming to support URM students that would lead to an increase in their success and retention in mathematics classes and in turn support them through their STEM majors.

Second, two WVU mathematics faculty members (co-authors on this chapter) participated in an MAA PREP Workshop on Emerging Scholars Programs in 2009. During the workshop, participants were given time to brainstorm the design of the program that would be most beneficial to their students. The year before the workshop, we had developed out-of-class Calculus Seminars to increase engagement and provide extra support for students in calculus, but had difficulty getting students to regularly participate. This combined with our low number of URM students and a desire to implement an inquiry-based approach different from other calculus sections led us to decide that a specially designated section of calculus, instead of out-of-class activities, would be more conducive to promoting community building among our URM students.

## 3. *Overcoming Obstacles*

As we moved forward with planning the ESP at WVU, we faced several obstacles including determining the structure of the course and recruiting students to enroll in the designated section. At the time, our other calculus classes met for 50 minutes per day for five days a week. We did not believe this would be the most efficient structure for integrating the active and inquiry-based learning practices we wanted to be the foundation of this course, so we decided to hold class meetings fewer days each week for longer time periods. The department chair supported the decision to offer a special ESP section of Calculus I that met for 110 minutes per day for three days a week (a total of 330 minutes per week, more than the 250 minutes provided to students in regular sections of the course). We found that students were willing to spend more time in class if it meant fewer days of class meetings, as it helped free up alternate days in their schedules for other classes. In addition, given our student demographics, with a very small number of URM students previously spread across sections of our calculus classes, we hoped an embedded model would allow us to recruit more URM students into one class and would help with community building. We were not concerned with segregating the students into a separate class, as we believed Treisman's work showed the value of community, and that implementing ours in this non-standard way was appropriate given the institutional context. We also knew our URM population was non-homogeneous; that is, our students come from a variety of underrepresented and intersectional backgrounds, and therefore the class still had a mix of identities and hence diverse within itself. We also considered long-term sustainability of the program and knew that this would be a viable option for continuing to serve these populations of students in a meaningful way. The department secured dedicated space from the university and worked with the university to renovate and procure equipment to provide the learning environment we wanted in a departmentally controlled classroom in which we could schedule the classes using usual meeting patterns. By providing resources to remodel the classroom to be more conducive to group activities, and equipping the room with moveable furniture, large whiteboards, markers, and erasers for students to use in class, the repurposed space was transformed into a dedicated environment for the type of classroom practices being employed while also providing a common meeting and study space for students in the cohort. At the same time, we sought departmental support to teach the class in smaller sections while still counting the course as a regular part of an instructor's teaching load.

Once we had the course size and meeting pattern in place, we had to find a way to effectively recruit students. We worked with department and university administrators to obtain a list of calculus-ready URM students and sent each a personalized recruitment letter. The letter was addressed to the student, came from the professor teaching the course, and described their recruitment as a *'nomination'* to join ESP. We chose this wording to help students feel that participating in the course was an honor, not something they would view as extra work or extra commitment. The faculty member teaching the course and sending the letters had to be willing to work with students who accepted the nomination to rearrange their schedules and provide any help that they might need with the university registration system. This was particularly important because many of the students were in their first year of college and therefore new to the concept of signing up for courses at differing times and making changes to their schedule. When the recruitment efforts did not produce a sufficient number of students, the course designers increased the nomination pool by including first-generation college students, especially those from rural Appalachia. We anticipated that students from rural Appalachian backgrounds would benefit from the structure and community-building of the ESP courses in similar ways to URM students. We used the same recruitment technique with this population, and it allowed us to increase the classes to more sustainable sizes (though still smaller than our usual calculus courses) and worked especially well in later courses.

We also needed to consider developing ESP-style curricular materials that would incorporate the spirit of advanced problems into the structure of the embedded model of ESP calculus. Over the summer before the first implementation, the professor teaching the course used the database of Workshop Calculus problems (http://www.betterfilecabinet.com/) and specifically the Treisman resources (http://bfc.sfsu.edu/cgi-bin/prob.pl?Treisman_Workshop_Resources) to develop worksheets to include the following types of exercises: conceptual understanding (exercises where students developed concepts such as deriving underlying theory and formulas for the day's topic), traditional calculus problems on the topic (like those found in most textbooks) and non-routine, challenging problems that pushed the students beyond the traditional problems (see Deshler et al. (2016) for more information on the design of the course and curriculum).

Students were placed in groups of four at the beginning of the course and worked in their groups on a common whiteboard to encourage collaboration and promote mathematical communication. In addition to working together on problems in each class, groups had to present solutions to their peers periodically during whole class discussions and at the beginning of the next class to ensure they continued to work together outside of class, if they needed to complete the day's work. To build a community of learners who were invested in the success of their classmates we also rearranged groups regularly so that students got used to working with *all* their classmates, not just their original group. To assist in the daily activities of this dynamic course, we also had an undergraduate teaching assistant help facilitate the group work during class. In later years, these teaching assistants were former ESP students who had been through the courses themselves and understood the classroom environment we wanted to establish.

After the first semester of ESP Calculus I at our institution, ESP Calculus II was offered in the same format. The program continued to grow and now spans four semesters of mathematics – Calculus I through Differential Equations. One important community building feature of the course is that when an instructor agrees to teach ESP, they are agreeing to teach a cohort of students for the full four-semester sequence of courses. While students may be recruited at various points in the sequence to fill gaps created by attrition due to various reasons, many students progress through the full sequence with a single instructor.

## 4. *Effectiveness*

### 4.1. *Community-building*

As our primary vehicle to support student success was through building community in ESP, we include here two comments from students as they reflected on their experiences in the course and the ESP program. The first is from a student who was in ESP for three semesters of calculus, the second is from a student who enrolled in only one semester of ESP.

> *I would like to start by saying that this class was so much more than just math… One of the biggest things that made this class stand out was the diversity. We were mostly engineer[ing] [major]s … accustomed to be[ing] the only diverse students in most of our classes. This class felt like a safe space to be in with people who looked and acted like me for 6 hrs a week. Being able to be with the same people each semester really helps secure life-long bonds and allows you to really feel a part of something bigger than yourself. This class is so much more than a math class, because of all the different life skills and connections you make along the way.*

> *Students should be allowed to draw on their peer's strongest math subjects as to temporarily fill in for an area of which the others are lacking during learning as to refresh any past concepts that were previously not solidified. For example, this occurred when Paola could further elaborate on Euler's Method during class and review sessions so that the class can receive another person's perspective as to how [to] set up problems and solve them. This caused the class to recollect subjects like this and rarely require any further review until the end of the semester.*

## 4.2. *Success in Calculus I*

To understand the effectiveness of the WVU ESP in supporting underrepresented students, two measures were investigated: the rate at which WVU ESP students pass Calculus I and the rate at which ESP students graduate with STEM degrees. For this analysis, we focus on Calculus I which has the highest enrollment of the four ESP classes (and of calculus courses, in general, at our institution). During the program, over the course of 10 years, 172 students enrolled in ESP Calculus I. For the same period, 8325 students enrolled in non-ESP sections of Calculus I. Table 1 shows that URM students are substantially more successful in ESP Calculus I. Students are classified as passing the class if they earn a grade of *C*- or higher. Many of the non-URM students in the ESP section of Calculus I were first-generation college students.

*Table 1*
*ESP and Non-ESP Calculus I Outcomes from 2009-2019.*

| **ESP Calculus I** | | |
|---|---|---|
| | Pass | Fail |
| URM | 80% (*N*=82) | 20% (*N*=21) |
| Non-URM | 71% (*N*=49) | 29% (*N*=20) |
| **Non-ESP Calculus I** | | |
| | Pass | Fail |
| URM | 49% (*N*=298) | 51% (*N*=315) |
| Non-URM | 62% (*N*=4783) | 38% (*N*=2929) |

A chi-squared test showed the difference between URM and non-URM success in the ESP Calculus I course was not statistically significant, while the difference in success of between these populations in non-ESP Calculus I was statistically significant ($\chi^2(1) = 42.36$, $p < .001$). The difference in success of URM students between non-ESP and ESP Calculus I was also significant ($\chi^2(1) = 32.79$, $p < .001$). As such, ESP Calculus I was successful in supporting URM students to the extent that there was no statistically significant difference between the passing rate of URM and non-URM students in ESP Calculus 1; this was a significant improvement over Non-ESP Calculus 1.

## 4.3. *Graduation Rates*

The graduation rates of ESP and non-ESP students were also compared. Again, only data for students enrolled in Calculus I were examined. Students still active at the university were removed from the sample leaving 117 ESP students and 5,796 non-ESP students who were enrolled in Calculus I (since Fall 2009). Table 2 shows the percentage of students graduating with a STEM degree, graduating with a non-STEM degree, or failing to graduate.

For the ESP courses, the difference in graduation rates between URM and non-URM students was not significant. However, for the non-ESP courses, the difference in graduation rates between URM and non-URM students was statistically significant ($\chi^2(2) = 38.29$, $p < .001$). The difference in graduation rates for students who graduated with a STEM degree for URM students was significant (at the 10% level) between ESP and non-ESP groups ($\chi^2(2) = 5.85$, $p = 0.05$).

*Table 2*
*College Outcomes for Students Enrolled in Calculus I from 2009-2019.*

| **ESP Calculus I** | | | |
|---|---|---|---|
| | Graduate STEM | Graduate non-STEM | Did Not Graduate |
| URM | 53% (*N*=35) | 23% (*N*=15) | 24% (*N*=16) |
| Non-URM | 43% (*N*=22) | 20% (*N*=10) | 37% (*N*=19) |
| **Non-ESP Calculus I** | | | |
| | Graduate STEM | Graduate non-STEM | Did Not Graduate |
| URM | 39% (*N*=175) | 23% (*N*=105) | 38% (*N*=171) |
| Non-URM | 52% (*N*=2,762) | 23% (*N*=1,225) | 25% (*N*=1,358) |

We note the small number of ESP students limits the statistical power of our analyses but posit that the substantial difference in both the success rates and graduation rates strongly suggest the ESP at WVU is positively benefiting underrepresented students within our institutional context. For that reason, department leaders have consistently supported sustaining the ESP sections of calculus each semester.

## 5. *Conclusion*

The Emerging Scholars Program at WVU represents an effort to facilitate increased learning in mathematics for students known to not be well-served by traditional educational models. These increases are created by building dynamic learning communities where students can engage in the discovery of mathematical ideas, reinforcing their sense of identity as STEM learners. The structure of those learning communities directly emulates the ESP model developed by Uri Treisman (Asera, 2001; Hsu et al., 2008), which increased the success of students typically underrepresented in STEM by sustaining the connections they may not find in traditional calculus classes. Our model increased success in these groups, but it was necessary to build spaces for these new learning environments and to creatively build departmental support. These were substantial investments in the form of faculty time, student support and facility improvements, but these investments supported the success of dozens of students who might not have otherwise persisted in their programs of study.

# 10 Towards Student-ready Mathematics Departments Creating Mathematics Placement Messaging Within an Equity-minded Framework

**Gaye DiGregorio**

*Colorado State University*

DiGregorio earned a Ph.D. in Higher Education Leadership at Colorado State University. Her dissertation focused on the perceptions about mathematics placement and precalculus/calculus mathematics courses among college first-generation, low-income, students of color. As the Executive Director of the Collaborative for Student Achievement, she has administrative expertise in empowering students and staff within an inclusive community, developing and implementing collaborative campus wide initiatives, and overseeing organizational change in higher education.

**Jess Ellis Hagman**

*Colorado State University*

Hagman earned a Ph.D. in Mathematics and Science Education at San Diego State University & University of California, San Diego Hagman's work is focused on dramatically increasing the number and diversity of people who succeed in undergraduate mathematics – especially introductory mathematics courses that often function as a roadblock for STEM intending students, increasing the number of qualified and prepared undergraduate mathematics teachers, and increasing the number and diversity of people who use science to improve the teaching and learning of undergraduate mathematics.

**Heather Novak**

*Colorado State University*

Novak earned a Ph.D. and works as a director in the office of Institutional Research Planning and Effectiveness at Colorado State University. Her research focuses on exploring institutional policies and programming that facilitate undergraduate student success, especially for students with identities that are structurally underserved by higher education institutions.

**Steve Benoit**

*Colorado State University*

Benoit earned a Ph.D. in Mathematics from Colorado State University. He is the co-director of the Colorado State University Precalculus program and coordinates the math placement program and the distance instruction programs for the Mathematics department. Benoit was the principal software and web developer for the online systems used for precalculus instruction and placement programs.

*Abstract: Since the mathematics placement exam is often a student's first experience with college-level mathematics, it sets the tone for students and may significantly impact future mathematics success. In this case study, we describe one department's effort to change mathematics placement messaging with an emphasis on equity, utilizing the Five Principles for Creating Equity by Design from the Center for Urban Education (Bensimon et al., 2016). The Math Placement Exam terminology was changed to the Math Placement Tool, the website integrated language to support student success, and an interactive tool was created to assist students in learning what mathematics is required for their major. After the math placement messages were enhanced, increased rates for math placement completion and placement into college algebra or higher were observed for all students and those with structurally disadvantaged identities (First-Generation students, Pell-Eligible students, and Students of*

*Color). This is descriptive data and does not indicate a causal connection between messaging and these outcomes. For departments working within this equity framework, it is essential to have university-wide collaboration, disaggregate data based on student identities to assess student outcomes, and involve student feedback in the process. Creating asset framed mathematics placement messaging is an initial step to improve success in mathematics.*

## 1. *Introduction*

As mathematics departments and universities focus on increasing first-year mathematics success and access to STEM degrees, it is important to consider mathematics placement - the initial experience that impacts academic success - and center our efforts on enhancing the placement experiences of students with structurally disadvantaged identities. For instance, at Colorado State University (2016), approximately 79% of First-Generation, Pell-Eligible, and Students of Color place into college algebra; and 65% First-Generation, 64% Students of Color, and 63% Pell-Eligible students complete three credits of math during the first year. Although math completion is evident for a slight majority of these students, improvement is needed when considering the rates of all first-year students. Overall, there is sizeable numbers of First-Generation, Pell-Eligible, Students of Color enrolling in college, realizing the need to focus on improving their math experience and graduation rates. This data is shared to substantiate these known facts and acknowledge that attention to equity must go beyond access and achievement (Gutiérrez, 2007). Since the mathematics placement exam is often a student's first experience with college-level mathematics, it sets the tone for students and may significantly impact future mathematics success. In this case study, we describe one department's effort to change mathematics placement messaging with an emphasis on equity. This work was initiated and carried about by a team including an administrator in student affairs whose recent doctoral work studied college mathematics placement (Gaye DiGregorio), a mathematics education researcher who studies college precalculus and calculus programs nationally (Jess Ellis Hagman), and a mathematician who runs the placement program at Colorado State University (Steve Benoit). Our team came together to focus on placement in part because of chance and in part due to the right timing: Catherine Good, a psychology researcher, gave a colloquium presentation to our department discussing her and her teams work on growth mindset within college mathematics (e.g., Good et al., 2012), and she discussed how mathematics placement exams and other high stakes exams can emphasize a fixed mindset on mathematics intelligence. The mathematics department chair was struck by this, and the following day Jess and Gaye presented Gaye's dissertation findings to the chair and Steve.

Gaye's dissertation examined the mathematics placement experiences of eight First-Generation, Pell-Eligible , STEM Students of Color at a different university. One of her key findings was how the mathematics placement process reiterated a fixed mindset among these students (DiGregorio, 2018; Digregorio & Hagman, 2019). The colloquium and the dissertation findings, together with data from our department and from Institutional Research, convinced Steve and the department chair that we needed to make changes to our own placement process and messaging, while centering a focus on equity. We submitted a proposal for and were granted internal funding to support student success awarded by the Provost's office.

We approached these efforts from an asset, rather than deficit, oriented perspective. We make the assumption that these differences are not due to differences in interest (or ability) in college or STEM, but rather due to systemic barriers faced by some populations of students and not others. For this reason, we refer to First-Generation students, Pell-Eligible, and Students of Color as “structurally disadvantaged” to emphasize the role of the systems on these populations’ marginalization rather than on these populations’ underrepresentation. Consequently, we were attentive to how these programs are disproportionately negatively affecting students from structurally disadvantaged identities and worked to make changes to our program expressly to provide support for these students’ success (Hagman, 2019). Our work is done within an overarching *student-ready* framework of supporting students that attend college rather than requiring students to adapt to college, (Brown McNair et al., 2016).

## 2. *Mathematics Placement Messaging Initiative*

A first step to transform mathematics placement is to create messaging about the placement process within an equity-oriented framework. We wanted to move from viewing this messaging as mostly a student responsibility to a combined effort between students, mathematics departments, and colleges and universities (Bensimon, 2005). To do this work while centering equity, we used the Five Principles for Creating Equity by Design from the Center for Urban Education (shown in Table 1; Bensimon et al., 2016), in which equity-minded language focuses on student success goals with the institutional intent of supporting students that have been structurally disadvantaged in achieving these goals. While this work has the potential to benefit all students, this framework encourages those making changes to focus on the needs and experiences of students with marginalized identities rather than to approach student success equally across the student population.

*Table 1*
*Five Principles for Creating Equity by Design from the Center for Urban Education (Bensimon et al., 2016)*

| Principle # | Definition |
|---|---|
| **Principle 1** | Clarity in language, goals, and measures is vital to effective equitable practices. |
| **Principle 2** | ‘Equity-mindedness’ should be the guiding paradigm for language and action. |
| **Principle 3** | Equitable practice and policies are designed to accommodate differences in the contexts of student’s learning—not to treat all students the same. |
| **Principle 4** | Enacting equity requires a continual process of learning, disaggregating data, and questioning assumptions about relevance and effectiveness. |
| **Principle 5** | Equity must be enacted as a pervasive institution- and system-wide principle |

Beginning in spring 2019, with financial support from the university, our team worked together to transform Colorado State University’s mathematics placement program to be student-ready with an initial focus on the messaging. While we attended to all five principles for creating equity by design, the second principle (that equity-mindedness should guide the language and actions) most directly impacted this stage of the work. Bensimon et al. (2016) articulate two ways to have an equity-minded focus in language: adopt language that centers institutional/department responsibility rather than students' deficiencies, and to remove all language that blames students or uses deficit language (such as "at risk"). By focusing on this principle and using these actions to implement it, we narrowed our focus to attend to the messaging of the placement program, including the name of the program and how it is communicated across campus and on the program's website.

A fundamental change to illustrate the paradigm shift from a student responsibility to institutional support was to change the name of the mathematics placement *exam* (MPE). An "exam" suggests that students can fail, and the MPE - which was the commonly used acronym by faculty and staff across campus - was language not easily understood by incoming students and conveyed a message different from our intentions. Changing the name to the Math Placement Tool (without an acronym) beginning fall 2020, provides clarity for students and focuses on its role in helping students to identify the appropriate course placement for student success rather than an exam that can be passed or failed.

For fall 2019, website content was changed substantially to eliminate the long-standing focus on the logistics of taking the exam, emphasizing regulations and what was not allowed. This administrative approach was transformed to focus on supporting student success at Colorado State University and providing information from the perspective of an incoming student. For instance, the website now includes messaging such as "we want to support you in being successful in the placement process"; and "This site will help you … get placed into the highest course in which you can be successful". Figures 1 and 2 below show examples of the new supportive messaging.

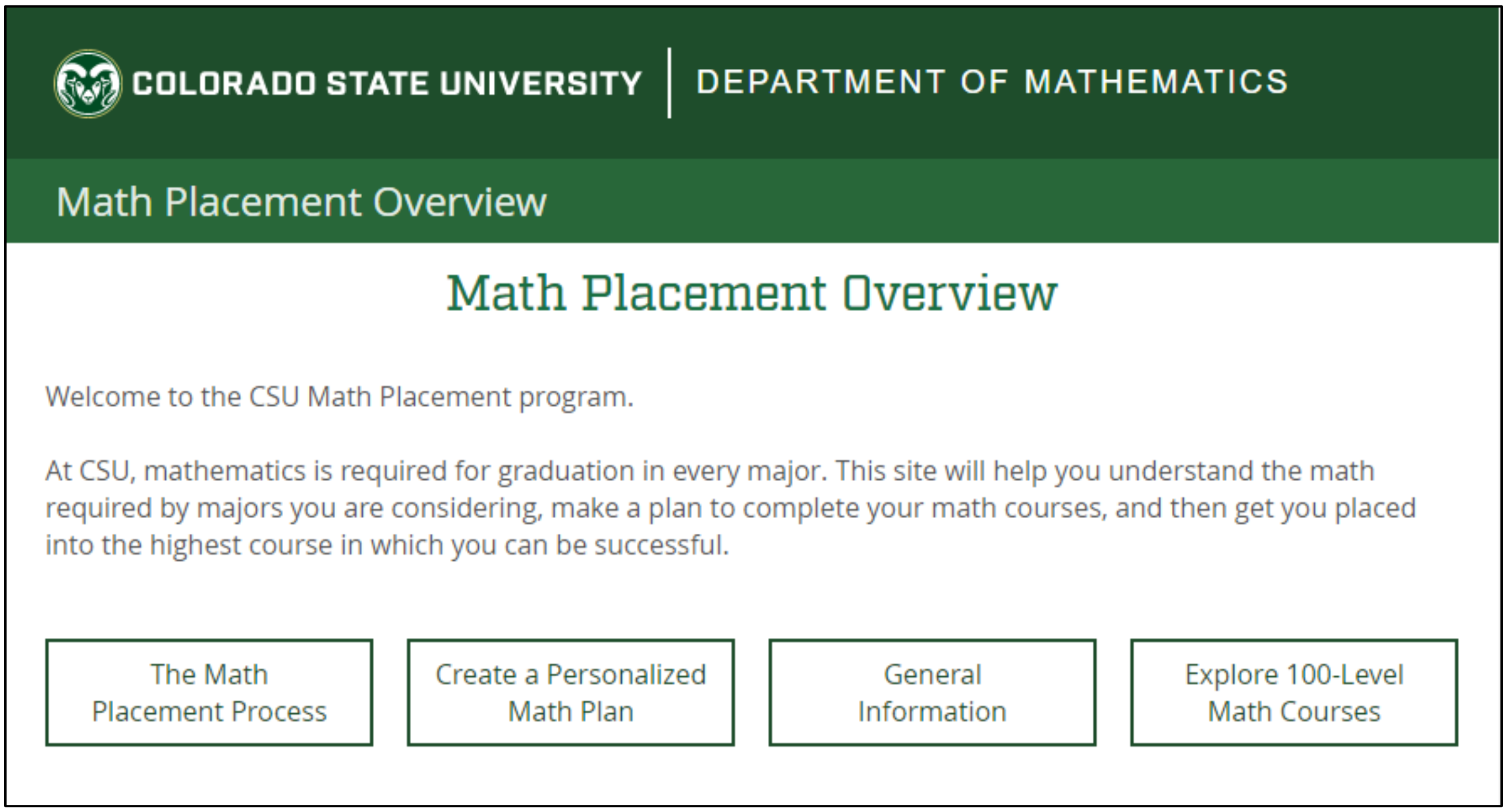

*Figure 1. The entry page for the department Math Placement website.*

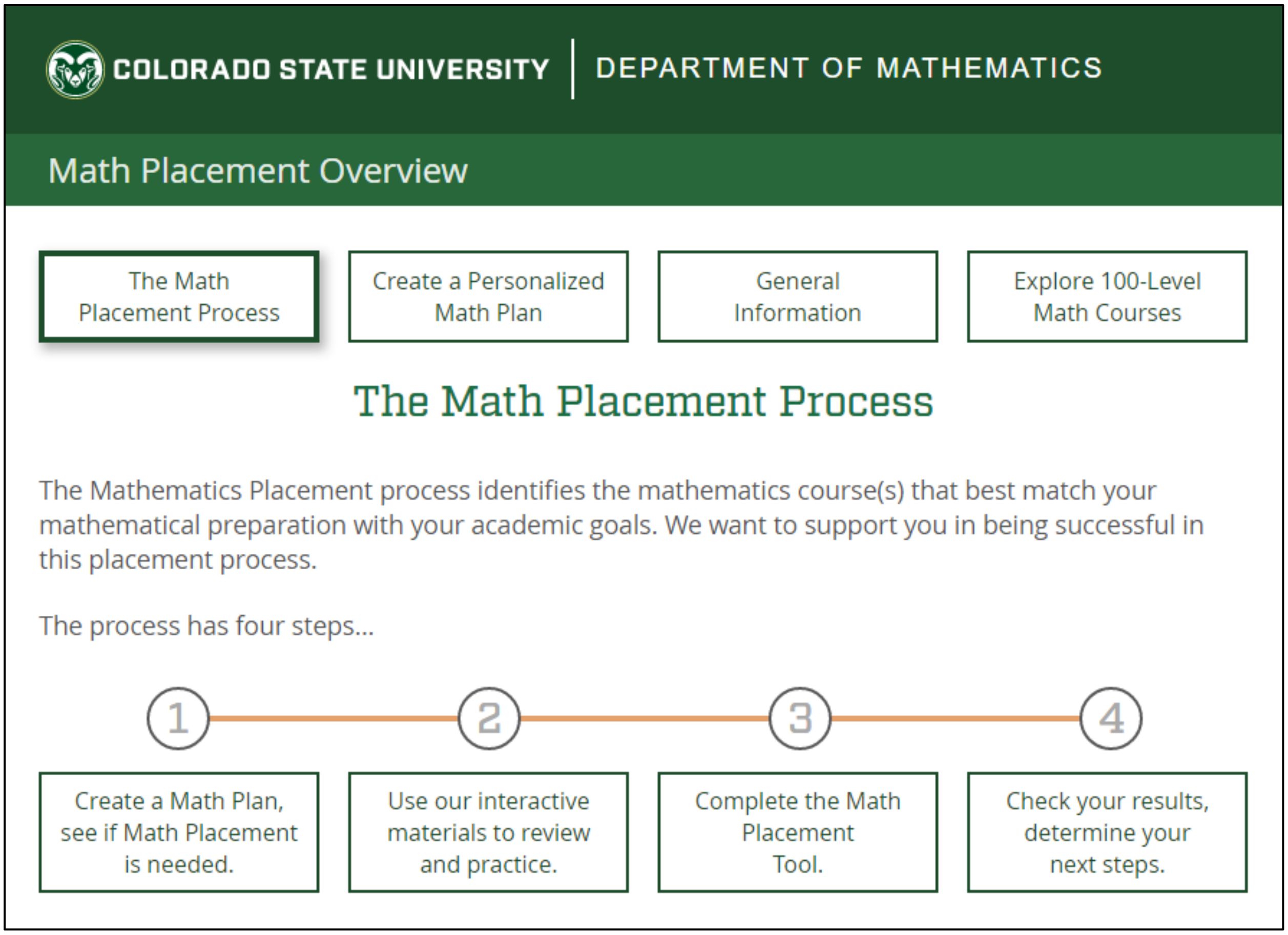

*Figure 2. Guiding students through the steps of the math placement process.*

In addition to the change from “exam” to “tool” and supportive website messaging, we changed several aspects of the content of the website to be aligned with students’ needs. First, a process was created for each student to select one or more majors of interest and discover what mathematics was required for them, providing a format for students to engage with mathematics placement information in an individualized manner. The web page where a student begins this process is shown in Figure 3. Additionally, new content was added to provide clear information about the content of the precalculus courses rather than listing the mathematics course numbers which an incoming student would not understand. Finally, information was laid out in a sequential manner with a step-by-step process for students, to more easily navigate through the placement process information, with clear indications in checklist form of the next step they should take, or that they have completed the process.

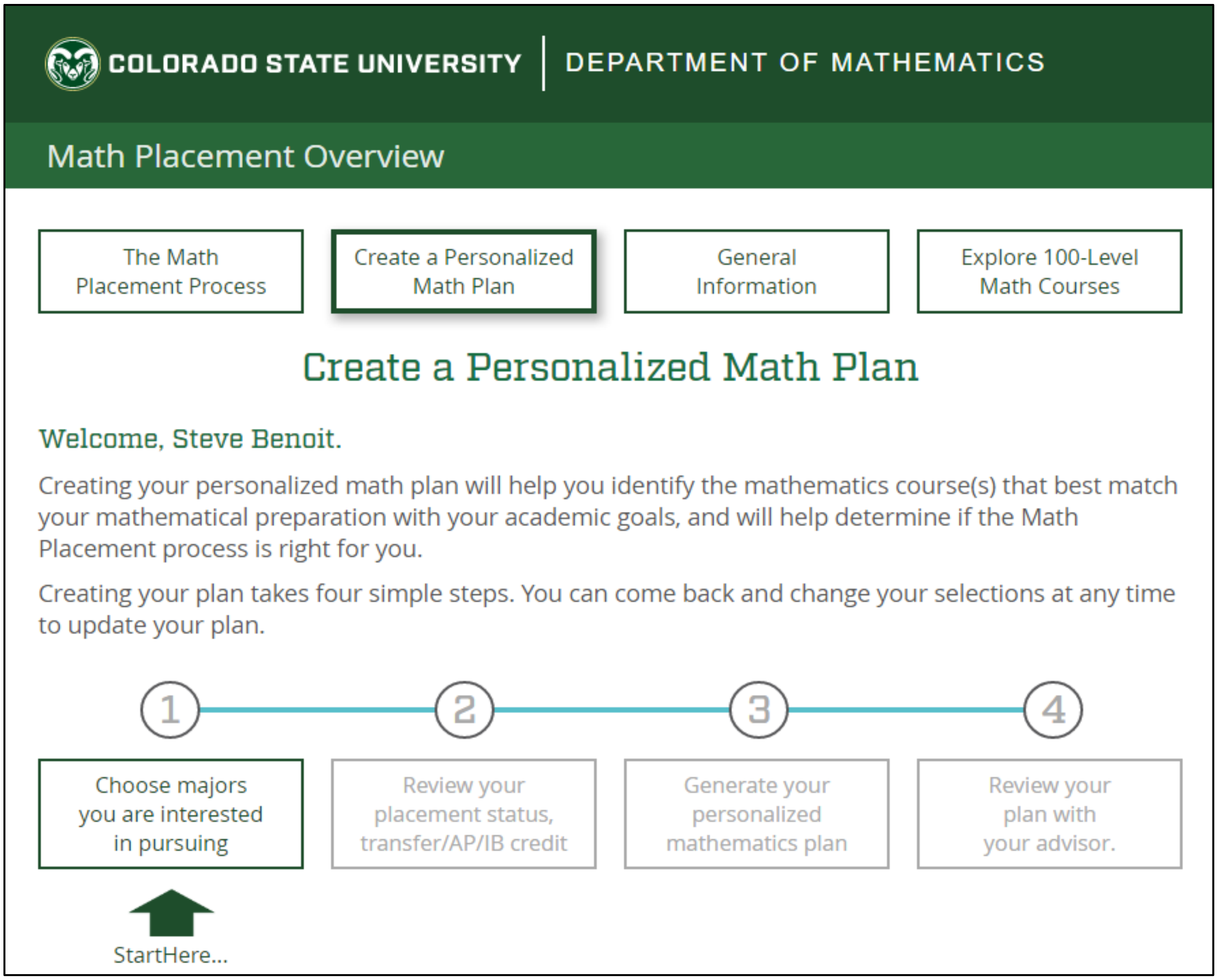

*Figure 3. Web tool to support students' creation of a personalized mathematics plan.*

The fifth principle of the Five Principles for Creating Equity by Design (Bensimon et al., 2016) emphasizes a need to focus on equity as an institution-wide effort. Traditionally, mathematics departments have taken sole responsibility for implementing the mathematics placement process. To endorse this university-wide, student success, equity-minded approach, enhancing mathematics placement messages was a collaboration between faculty and staff in the Department of Mathematics, staff in the Collaborative for Student Achievement (the central organization that supports advising and orientation of incoming students at CSU), staff in Admissions who work with prospective students, staff in the Registrar's Office supporting the university student portal, and staff in CSU's Creative Services group, with expertise on designing inclusive websites. This collaboration and institutional support to provide integrated communications about mathematics placement within overall university communication and registration processes for incoming students, was supported with institutional funds from the Provost's Office and was made possible by our multi-unit team.

While we do not want to diminish the huge importance of the financial support, we think the multi-unit team was more integral to our success. For other institutions seeking to make changes to their placement messaging, we believe starting with a team from within and outside the mathematics department is key: having team members from within the math department is necessary to know the details of the program and to impact the direction of the changes, while having team members outside the department create broader sustainable changes throughout the institution in admissions, orientation, and registration. University-wide efforts that enhance mathematics placement messaging within a student success equity framework are one step toward supporting students to be successful in mathematics.

Table 2 summarizes how the principles of equity design were applied to the math placement messaging.

*Table 2*
*Summary of how we enacted equity by design in the math placement messaging*

| Five Principles of Creating Equity by Design | Applying these principles with Math Placement Messaging |
|---|---|
| Principle 1: Clarity in language, goals, and measures is vital to effective equitable practices. | Removed unnecessarily confusing acronym (e.g., MPE). |
| Principle 2: 'Equity-mindedness' should be the guiding paradigm for language and action. | Focused messaging around support for students rather than logistics; changed language from "Exam" to "Tool." |
| Principle 3: Equitable practice and policies are designed to accommodate differences in the contexts of student's learning—not to treat all students the same. | Interviewed students with structurally disadvantaged identities to inform additional ways to support students in this process. Results can be shared upon request. |
| Principle 4: Enacting equity requires a continual process of learning, disaggregating data, and questioning assumptions about relevance and effectiveness. | Tracking math placement results longitudinally by various student demographics to inform future enhancement efforts. |
| Principle 5: Equity must be enacted as a pervasive institution- and system-wide principle. | Math Placement messaging was a collaboration between faculty and staff in the Department of Mathematics, staff in the Collaborative for Student Achievement (the central organization that supports advising and orientation of incoming students at CSU), staff in Admissions who work with prospective students, staff in the Registrar's Office supporting the university student portal, and staff in CSU's Creative Services group, with expertise on designing inclusive websites. |

### 3. *Tracking Mathematics Placement Completion and Results*

One way to assess the effect of changes to mathematics placement messaging is to track the percentages of students that completed mathematics placement prior to their first semester and placed into at least college algebra, disaggregated by student populations. Since research shows that completion of math in the first year is a predictor of college graduation, an increase in students placing into at least college algebra in their first term has the potential to affect overall college graduation (Colorado State University, 2015). At Colorado State University, all students are required to complete 3 credits of quantitative reasoning / mathematics, and, depending on a students' major, these credits can be earned in the precalculus or calculus sequence, or by taking courses for non-STEM majors (e.g., Mathematics in the Social Sciences or Patterns of Phenomena).

The tables below display the increased rates for math placement completion and placement into college algebra or higher for all students and those with structurally disadvantaged identities (First-Generation students, Pell-Eligible students, and Students of Color), after the math placement messages were enhanced.

*Table 3*
*First-Time, Full-Time Students that Complete the Math Placement Tool Prior to Semester Start (%)*

| | Overall % (N) | First-Generation % (N) | Pell-Eligible % (N) | Students of Color % (N) |
|---|---|---|---|---|
| Fall 18 | 62.2 (5280) | 57.7 (1187) | 57.8 (1128) | 59.9 (1482) |
| Fall 19 | 68.6 (5107) | 63.2 (1197) | 64 (1077) | 65.9 (1477) |
| Percentage Point Difference | 6.4 | 5.5 | 6.2 | 6 |

Rates and Population Sizes for Math Placement Completion and Placement in College Algebra or higher, and Enacted Equity by Design in the Math Placement Messaging

*Table 4*
*Tool Placement into College [SB1] Algebra or Higher among First-Time, Full-Time Students (%)*

| | Overall % (N) | First-Generation % (N) | Pell-Eligible % (N) | Students of Color % (N) |
|---|---|---|---|---|
| Fall 18 | 38 (5280) | 29.1 (1187) | 31.4 (1128) | 34.1 (1482) |
| Fall 19 | 42.5 (5107) | 32.7 (1197) | 32.9 (1077) | 37 (1477) |
| Percentage Point Difference | 4.5 | 3.6 | 1.5 | 2.9 |

It should be noted that not all students would be expected to complete the placement tool because some majors don't require college algebra. These percentages are descriptive and do not indicate a causal connection between messaging and mathematics placement outcomes, rather provide an initial step in tracking these outcomes. We recognize a tension that exists in using such data as an indicator of outcomes by bringing focus to comparing the outcomes of students from structurally disadvantaged identities to all students. Our goal here is not to gap gaze (i.e., emphasize the gaps in achievement between some students and others), but rather to acknowledge how the changes to the mathematics messaging affect different populations of students as one way to assess our changes.

As you can see from Figure 4 which is another illustration for this data, the rate of students completing the Math Placement Tool prior to the start of the semester increased in fall 2019 compared to fall 2018. The 6.4 percentage point increase was seen across the overall incoming student population (the blue dashed lines in the figure), and across each demographic group examined; however, students with structurally disadvantaged identities had smaller increases than the overall cohort. For example, the 6.2 percentage point gap in the figure represents the difference for Pell-Eligible students who completed the math placement in fall 2019 compared to fall 2018 (Pell-Eligible students in the fall 2019 cohort completed the math placement at a rate of 64% and that is 6.2 points higher than 57.8% in fall 2018)

Similarly, in fall 2019 First-Generation, Pell-Eligible, and Students of Color had higher percentages of students that placed into college algebra or higher, although these increases were less than the 4.5 PP increase overall. For the messaging to contribute to equitable outcomes we would expect the percentage point changes to be equal or greater for students with structurally disadvantaged identities.

These results could be associated with the removal of some structural barriers that have inhibited specific groups of students from receiving higher placements on the tool; however, the disproportionate gains indicate more work needs to be done. Such outcomes are consistent with the principles of Universal Design for Learning, which is an approach to curricular and program design that foregrounds the needs of a structurally disadvantaged population of students in the design process. These approaches focus on reducing barriers to success for a subpopulation of students, but often are helpful for all students (Rose et al., 2006).

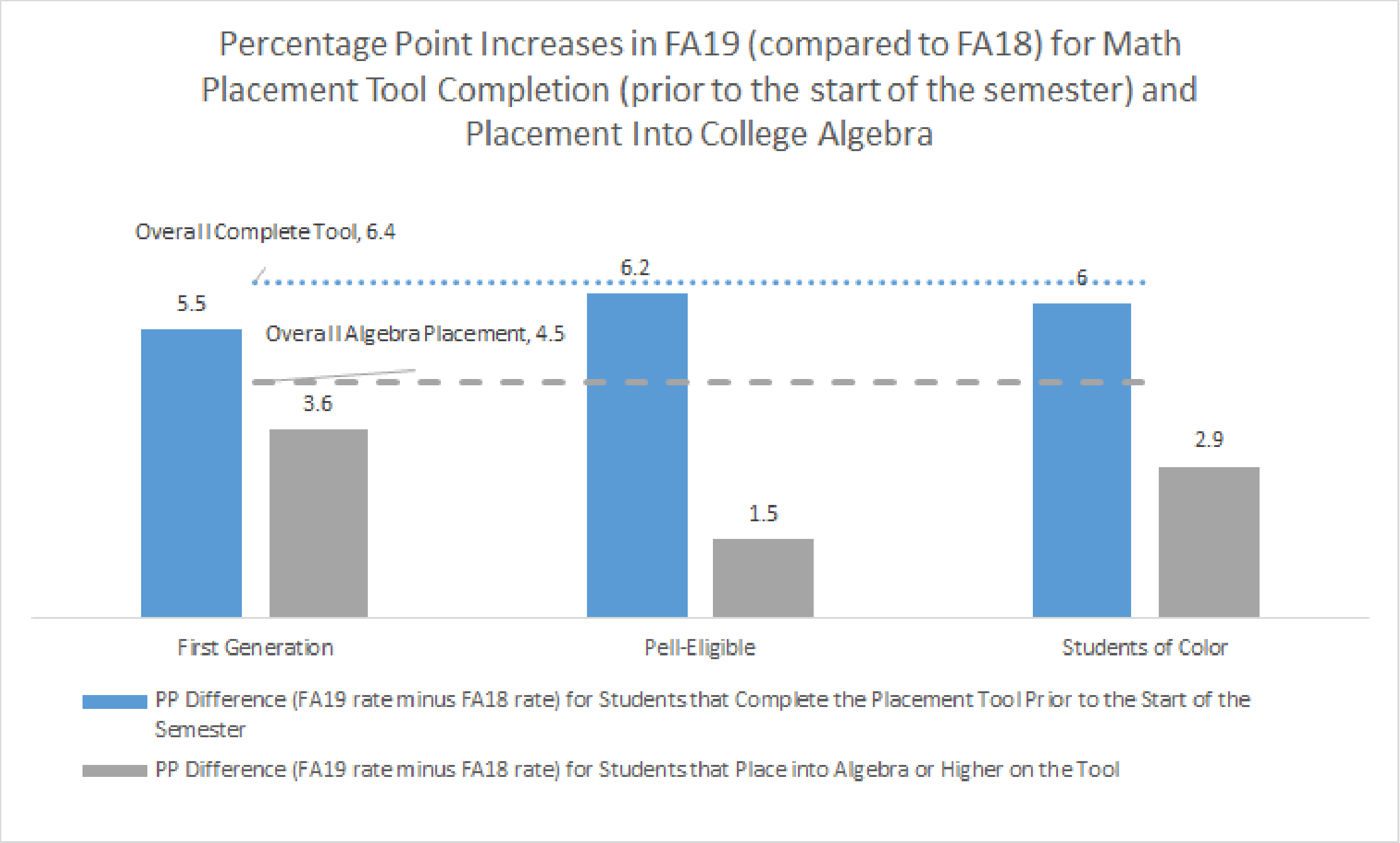

*Figure 4. Increases for Math Placement Tool Completion and Placement into College Algebra or higher after messaging changes.*

Once again, this assessment does not mean that enhanced mathematics placement messages resulted in these increased student outcomes. Continued longitudinal data of these outcomes is needed as the messaging continues to be enhanced, most recently with the change of terminology from the Math Placement Exam to the Math Placement Tool. Additionally, more qualitative research on the mathematics placement experiences for First-Generation students, Pell-Eligible students, and Students of Color at Colorado State University, is needed.

### 4. *Continued Efforts to Enhance Mathematics Placement Messaging*

Although these improved mathematics placement messages may be associated with the increase in the number of students completing mathematics placement and placing into at least college algebra, more work needs to be done. While our current efforts have been motivated by supporting some students but affect all students, we now need to specifically understand how these changes are perceived and experienced by structurally disadvantaged populations of students.

We have conducted individual interviews with students that are First-Generation, Pell-Eligible, and Students of Color who gave feedback about the mathematics placement messaging as well as the experiences with completing the Mathematics Placement tool. This feedback informs how to further enhance the mathematics placement messaging with a focus on students with structurally disadvantaged identities, and how-to better tailor the placement process to these students. Specifically, we will be exploring what barriers exist for each population of students. For example, we will explore how the messaging we use may leverage academic language in a way that is a barrier to First-Generation students less accustomed to this language. We also have developed a video included on the website from second year students on how they were successful in preparation and completion of the mathematics placement process. Using continuing student stories is a well-researched and recommended strategy to encourage student success. This feedback can also help determine the mathematics placement process for the future, to further encourage incoming students' success in mathematics.

The plans for future changes described here exist within a department that is continually attending more to issues of diversity, equity, and inclusion, but for which many in the department are focused outside of the first-year mathematics program (including mathematics placement). Because of this, these changes happened without much fanfare within the department, but have been very well received outside of the department (such as in Admissions).

### 5. *Recommendations to Enhance Mathematics Placement Messaging*

In their census survey of mathematics departments, Rasmussen et al. (2019) found that student placement was widely viewed as very important to precalculus and calculus programs. While student placement is important, it is also the aspect that departments are least successful at.

Placement is also an aspect of mathematics departments that affects students far beyond mathematics majors and STEM majors. Further, because of changes to what courses are allowed to be offered for credit in college and removing standardized tests as part of admissions (such as in California), many departments reevaluate their first-year mathematics offerings, and therefore must also revisit how to place students into these courses.

We view placement as a critical component of college mathematics programs that is in high need of attention, and one that is essential to address with a focus on equity. Our core recommendation for institutions to enhance the mathematics placement process is to constantly integrate equity within a focus on student success, such as by using the equity-minded framework (Bensimon et al., 2016) as a continual reference throughout the change process. For departments working within this framework, it is essential to:

- Have university-wide collaboration,
- Disaggregate data based on student identities to assess student outcomes, and
- Involve student feedback in the process.

We hope the equity-minded framework will assist structurally disadvantaged students in being more successful with mathematics placement and one step in changing mathematics departments to become *student ready* (Brown McNair et al., 2016). Improving student engagement and success in first year mathematics courses is imperative because of the correlation between early mathematics completion and graduation (Adelman, 2006; Colorado State University, 2015). It is time to support all students who are excited about STEM careers, especially those who have been structurally disadvantaged, and ensure that their mathematics placement does not discourage these interests, but rather supports them in thriving in their STEM pursuits.

# 11 Instructional Practices that Foster Success in Underrepresented Populations

**Edgar Fuller**

*Florida International University*

Edgar Fuller is Professor of Mathematics and Director of the Center for Transforming Teaching in Mathematics at Florida International University. He studies student learning and engagement in mathematics while trying to improve access to mathematics at all levels. He is currently co-PI on the NSF Catalyzing Change in Calculus project working to develop active and engaging coursework in calculus that leads to better access to STEM careers for all students.

**Adam Castillo**

*Florida International University*

Adam Castillo is a Visiting Research Assistant Professor affiliated with the STEM Transformation Institute and Center for Transforming Teaching in Mathematics at Florida International University. His research focuses on examining the impact of active learning approaches in mathematics, improving student success in postsecondary mathematics, and ensuring that equitable access to high-quality mathematics education is available to all students to help prepare them for success in college and their future careers.

**Charity Watson**

*Florida International University*

Charity Watson is a Clinical Assistant Professor in the Center for Transforming and Teaching in Mathematics and STEM Transformation Institute at Florida International University. Her research is focused on college-level students and their performance and persistence in STEM fields and in introductory-level mathematics courses. She is interested in factors that may specifically impact underrepresented minorities in mathematics and STEM.

**Pablo Duran**

*Florida International University*

Pablo Duran Oliva is currently a Visiting Research Assistant Professor at the STEM Transformation Institute at Florida International University. His research has focused on the impact of active learning approaches in mathematics, including classroom equity, changes in student attitudes and course performance. His research interests additionally include the development of mathematical thinking and curriculum alignment with student needs in STEM disciplines.

**Geoff Potvin**

*Florida International University*

Geoff Potvin is an Associate Professor in the Department of Physics and the STEM Transformation Institute at Florida International University. His research interests include: broadening participation in the physical sciences using an identity lens, with a particular focus on recruitment and retention across the secondary, post-secondary, and graduate boundaries. He has served as Chair of the Physics Education Research Leadership and Organizing Council and has been a member of the AAPT Committee on Diversity, AAPT Committee on Graduate Education, and the APS Forum on Education.

**Laird Kramer**

*Florida International University*

Laird Kramer is the Founding Director of the STEM Transformation Institute and a Professor of Physics at Florida International University. He focuses on facilitating institutional change through implementation of, and research on, evidence-based educational practices to establish equitable outcomes for all students.

*Abstract: Student success in calculus remains one of the major obstacles to students seeking to pursue careers in STEM nationally. In this work we present an active learning-based approach in calculus that leverages group work experiences, near-peer facilitation, argumentation and communication strategies to create a culturally relevant learning experience that allows students to learn and express mathematics naturally and based on their experience. By shifting instructional practices to more student-centered work, we find that students, especially those from traditionally underrepresented populations, are more successful and develop higher levels of confidence in mathematics.*

## 1. *Problem Statement*

Classroom approaches to instruction in mathematics have historically focused on the elucidation of a collection of ideas by an expert coupled with precise delineations of computational methods related to those ideas. Instructors describe concepts, illustrate these concepts with examples, and then ask students to replicate additional examples to confirm mastery of concepts. Such approaches may create difficulties for students whose background does not align with the traditional white, suburban experience (Russell & Atwater, 2005) since many of the norms observed in those instructional models create artificial boundaries in the classroom environment that discourages students outside those norms (Dika & D'Amico, 2016; Deshler et al., 2019; O'Hara, 2022). A number of studies have shown (Freeman et al., 2014; Ellis et al., 2016; Theobald et al., 2020; Miller et al., 2021) that classroom practices that focus on student engagement may lead to improved outcomes. At Florida International University (FIU), a large, urban research university with more than 58,000 students, of which 64.2% are Hispanic/Latine[29], this approach to instruction in its college algebra, precalculus and calculus courses resulted in success rates among its students that were consistently near 50%. In first semester calculus, success rates averaged 58% over the five-year period prior to the spring of 2018, at which point faculty in the Department of Mathematics and in the FIU STEM Transformation Institute collaborated to develop new approaches in this course that would foster student engagement, deeper learning and ultimately higher levels of success.

## 2. *Response to Problem*

First semester calculus at FIU historically enrolls more than 2000 students per year, and of these 74.1% identify as Hispanic/Latine. At the same time, 14.6% identify as Black/African American, and 9.2% report being White, non-Hispanic students. FIU has a strong commitment to continuous improvement in student equity and success in mathematics, and so the success rates in calculus represent pedagogical as well equity concerns. An effort was undertaken to leverage local experience with evidence based instructional practices into an approach to calculus instruction that would engage all students in learning experiences that resonated strongly with them. Recent projects such as the Characteristics of Successful Programs in College Calculus and Progress through Calculus, led by the Mathematical Association of America, have collected data that underline the importance of supported 'ambitious teaching practices' such as active learning approaches in mathematics. Guided in part by these results as well as successes observed in other disciplines, faculty working within the STEM Transformation Institute at FIU engaged in an effort to radically change the approach to calculus. The first semester of calculus was redesigned around a studio model that offers continuous opportunities for students to develop mathematical understanding that connects with their context and centers around their experiences.

29In this case study we use Latine to refer to students of Latin American descent from all genders and languages (Salinas & Lozano, 2021) in our population.

The curriculum is designed to foster student inquiry in calculus from notions of limits and rates of change to the fundamental theorem of calculus. This approach, the Modeling Practices in Calculus (MPC) curriculum, aims to promote inclusive teaching practices and support students as they develop in five main areas:

- Developing practices of mathematicians such as problem solving, reasoning, sense-making, and higher-order thinking inspired in part by the framework of Lin and Selden (Lin & Selden, 2009). We encourage students explicitly to notice patterns, make sense of the structures they observe, develop theories that explain them, and persevere through difficulties towards solutions. We model these as 'practices of mathematicians' and use active learning processes to help students see themselves as doing the work of mathematicians as they learn rather than memorizing facts. Engaging in the practices of mathematicians to learn calculus helps students understand how mathematical knowledge is created and develops the critical thinking, problem solving, and analytical skills students will need to be successful in mathematics. Instruction embedding practices of mathematicians helps to focus on the collective ways of problem solving, reasoning, and symbolizing established while discussing mathematical ideas (Cobb et al., 2001).
- Ubiquitous cooperative learning where students in small groups work to accomplish shared learning goals. Cooperative learning is one technique known for promoting equitable learning environments in mathematics (Gutierrez, 2002) as it fosters learning academic and social skills (Cohen, 1994; Johnson et al., 2007; Slavin, 1996) and creates positive interpersonal relationships (Johnson et al., 2014). Students work in groups, beginning each class with an overview of a new idea followed by a set of activities that build progressively from idea to idea. For instance, during the first week of class, students work together to develop an informal view of limits by understanding how a function changes over time. By working together to find to average rates of change over intervals that get shorter, approaching the point of interest, students begin to recognize that as the length of intervals decrease, the average rates of change start to approach a unique number (i.e., the instantaneous rate of change at the point of interest). After working for a variable period of time, students come back together with the larger class and present the connection between the average and instantaneous rates of change.
- Using argumentation/metacognition to engage in meaningful discourse, provide justifications and evaluate claims. The inclusion of instruction promoting mathematical argumentation can provide a deeper understanding of mathematics as students become generators of knowledge out of their reasoning and sense making (Yackell & Cobb, 1996).

- Establishing mathematical fluency through concept exploration. Fluency involves carrying out methods flexibly, accurately, efficiently and appropriately (Kilpatrick et al., 2001). Students build fluency by noticing mathematical relationships and using strategies through the study and discussion of various concepts in learning activities, as well as through tasks that promote reasoning and problem solving.
- A culturally responsive environment that allows students to construct their own understandings. Instruction that is culturally responsive makes learning more accessible, relevant, and meaningful for students (Ginsberg & Wlodkowski, 2009; Celedon-Pattichis et al., 2018). This includes the use of language that resonates with student experience as well as examples that connect well with their context, but it also means cultivating student expression of ideas from their own point of view so that their sense-making is developed internally instead of imposed externally. Our courses integrate Learning Assistants (LAs), trained 'near peer' undergraduate classroom facilitators prepared to support learning with groups (Otero et al., 2010). LAs are natural agents of culturally responsive learning, providing insights and connections from the point of view of a recent participant in the course. LAs provide an essential component that establishes a strong social connection surrounding the mathematical work being performed, and these connections enhance student ownership of the resultant mathematical knowledge. They provide a stronger connection to the ways in which students are constructing knowledge in the class and help mitigate the 'blind spots' that experienced mathematicians bring into dialogues.

For our active learning approach, we developed a set of activities that lead students in mathematical investigations. The investigations we chose develop core ideas such as limits, rates of change, related rates, and accumulation. The course is split into three units over which student learning is assessed using individualized exams followed ultimately by a comprehensive final. Class sizes range from 20 to 40 students typically, and groups are arranged to be three to four students. Group assignments are rotated throughout the semester.

In a typical class meeting, students are introduced to a new idea or topic at the beginning of class via a brief discussion. This time is used to set up a set of guided notes for students to work through in groups that cover introductory examples while leading them progressively to more and more formalized mathematical concepts, definitions, and theorems. After completing the notes, students are asked to summarize concepts and ideas in a whole class discussion briefly before beginning work on learning activities. These activities are intended to challenge these ideas and build skills and fluency by presenting problems that are similar to examples in the notes, but often require higher-order and critical thinking skills. While engaging in the series of learning activities, students proceed through a concept development with problems designed to lead them to essential insights and to challenge their mathematical toolset in ways that lead to increases in self-confidence, self-regulation, and effective meta-cognitive processes. One tool for helping students develop mathematical fluency is small whiteboards. These are used by student groups to share out problems at various points of class time called *board meetings that can also be used for summative assessment of understanding after the learning activities are* *completed*. Communicating solutions, even partial or incorrect ones, helps students test their reasoning and mathematical language.

Investigations are supported by both instructor guidance from an expert point of view and by LAs. LAs help students develop skills such as creating and defending ideas, making connections between concepts, and solving conceptual problems (Alzen et al., 2018). They also increase the flow of ideas from the students to faculty in general so that discussions are more strongly centered on student points of view. The classroom format and approach draw out student interaction in this way intentionally in order to enhance the connection of mathematical thought and concept development to student experience. LAs meet weekly in a group led by faculty from the course instructional team to discuss specific material and concepts for our calculus course as well as weekly in a seminar course to discuss pedagogy in the broader context of STEM courses. The seminar leader (or leaders if team taught) provides insight into the material being introduced as well as how students may interact with it. LAs also discuss techniques for getting students to engage with the development of their own mathematical ideas instead of simply showing students how to solve problems. Seminar leaders discuss this work with the entire course instructional team and provide that group with ways to work with LAs in the classroom that leverage the LA development in the seminar. The larger course instructional team meets every two weeks as a group to discuss LA participation in their classrooms as well as pacing, exam development, and other course logistics.

### 3. ***Outcomes and Results***

We have implemented our curriculum and classroom approach over the course of five semesters beginning with a two section, two instructor pilot in the spring of 2018. We increased the number of sections and instructors in subsequent semesters, reaching all but one section by the spring of 2020. During this implementation, the group of faculty working in the new curriculum would meet with the project leaders and prior instructors before the term to learn about the program and to plan for the upcoming term. Over the course of the semester, the team would meet weekly to reflect, discuss success and failures, revise material, and plan for upcoming classes. These ongoing meetings were an essential part of the project as they provided insight into the ways in which faculty could develop classroom practices that aligned with the model. Working with LAs and using a group work model on a daily basis requires sustained effort on the part of faculty, and the instructional group meetings provided sustained support in these areas so that faculty could manage any aspects of the new approach that became challenging. Discussions focusing on fostering group work, maintaining a pace in the course that was effective, and dealing with differences in student engagement levels were common and the group meetings provided a space for discussion and sharing of methods for responding to these challenges.

Students in this active learning environment experience significantly more successful student outcomes than those seen historically in other traditional, lecture-based classrooms. Overall, success rates, defined as an A, B or C letter grade assigned for the course, have averaged above 80% in sections using the MPC materials and practices while comparison sections, both those identified to be part of a research trial now being analyzed for publication as well as other sections not using the MPC approach, have averaged 63%. The increase in success occurred using assessments that had been historically used in prior semesters and, more importantly, that were used in identical format by both MPC and comparison sections indicating that the increases in learning outcomes are robust. Within that success, students from Hispanic/Latine as well as Black/African American groups have succeeded at rates similar to or higher than the total population. Hispanic/Latine students succeed at rates near 90%, and almost 80% of students identifying as Black/African American succeed. These outcomes represent a more than 30 percentage point increase from historical rates of success and are significantly higher than rates observed in paired comparison sections. Overall, more than 1000 students have participated in the MPC curriculum and these results have persisted across multiple terms and with growing numbers of instructors involved. Outcomes from this work also include higher levels of student satisfaction, value, motivation and self-confidence (Castillo et al., 2022), where a statistically significant increase was observed with small effect size (ANCOVA pre/post comparison, effect size partial $\eta^2 = 0.011$) for all students after experiencing mathematics investigations in a supportive environment as measured by the composite Likert scale based score computed from the Attitudes Towards Mathematics Inventory instrument. Within this scale of student attitudes, the main statistically significant increase in a sub-scale was observed in student self-confidence(ANCOVA pre/post comparison, effect size partial $\eta^2 = 0.02$). This suggests that MPC classrooms are more effective in helping students change their views on how capable they are in succeeding in other mathematics courses in their future coursework.

## 4. *Obstacles, Barriers, and Recommendations*

Students sometimes struggle in environments where they must engage in inquiry alone and so it is critical that faculty managing such classrooms develop the skills and understandings necessary to assist students effectively. A major component of this process must focus on working with Learning Assistants to facilitate interactions with student groups while they are working with new ideas. Professional development opportunities must be provided that support the use of these multi-level instructional strategies in order to foster best practices and to provide ample time for reflection in instructional teams. Faculty will also experience difficulties outside of the classroom and work must be done to facilitate conversations with other departmental faculty, deans and other administrators in order to ensure a successful and sustainable rollout.

Resources will need to be procured to support active learning classrooms and Learning Assistant support. The Learning Assistant Alliance (https://www.learningassistantalliance.org) is one avenue that institutions without an LA program can contact to support the implementation and assessment of the LA Model with goals of improving students' educational experiences and access across STEM disciplines. This model is one of social and structural organization that encourages and supports the adoption of research-based instructional strategies. Acting as a faculty development program, the LA program can offer resources and structures that lead to changes in values and practices among students, faculty, departments, and the institution (Close et al., 2018; Otero et al., 2010). To better support active learning, learning spaces should be re-envisioned to incorporate collaborative learning with whiteboards, tables and room for movement and interaction. This may result in decreases in student density and class size, but in many cases higher success rates will offset lower course capacity in the long-term. Additionally, it is not clear from our study how LAs and the active learning format impacted outcomes individually. The approach focuses on leveraging all aspects of the MPC curriculum into a powerful space to engage students in learning, and we sought to bring all the tools to bear synergistically.

Overall, students have consistently succeeded in second semester calculus at the same rate as their peers coming from traditional, non-MPC classes. There are many possible confounding influences on student success in the successor calculus II course and more work would be needed to understand how MPC students engage in the traditional lecture format after using active learning tools in calculus to better understand those outcomes. Indeed, it may be that success in calculus II would increase as well if an MPC approach were used there. This is the focus of future work.

## 5. *Summary*

We present here a summary of an effort to build more inclusive classrooms in a calculus setting at a large, urban, research-intensive university. Our goal was to provide a setting where students could engage in a collection of curated learning activities that would help them develop a mathematician's approach to understanding calculus concepts while improving their metacognitive skills in a way that facilitated growth and an improved approach to problem solving in new environments over time. Our outcomes indicate that the active learning foundation that we employed and the faculty support provided to bring about change in our classrooms led to improved student learning outcomes as measured by grades as well as other more focused metrics. A key component in implementing programs such as this will be consistently supporting faculty in their effort to lead their classrooms in a different way. This work includes empowering students and leading group work while mentoring Learning Assistants who help foster the active learning classroom environment that is the ongoing target of the process. This level of effort can be challenging to maintain but if supported in a sustained manner it leads to students who own their knowledge and are able to use it in ways that fit their individual vision of what it means to be successful in mathematics. Moreover, allowing them to grow and challenge themselves with these rigorous mathematical concepts will connect their vision with the discipline as a whole and lead them to long term success using mathematics.

# 12 Creating Space for Student Agency to Support Success in Online Pre-Calculus

**Shannon Golden**

*Colorado State University*

Shannon Golden is a third-year graduate student in the mathematics department at CSU. Her doctoral work will be in computational algebra, and her masters' work is an extension of the above work with an emphasis on analyzing how different influences affected students' ability to accept or reject their initial mathematical affect and their classroom behaviors, and how that impacted their success. Her goal is to enter academia to provide equitable instruction in mathematics at any level, using student agency in course design as a key tool towards achieving this goal.

**Jess Ellis Hagman**

*Colorado State University*

Jess Ellis Hagman is an Associate Professor in the Department of Mathematics at CSU. Her area of research is undergraduate mathematics education. Her work is driven by the goal of dramatically increasing the number and diversity of people who thrive and succeed in undergraduate mathematics-especially introductory mathematics courses that often function as a roadblock for STEM intending students. Her current research includes studying characteristics of successful precalculus and calculus programs, focusing on investigating ways departments can create diverse, equitable, and inclusive introductory mathematics programs.

*Abstract: As evidenced by this volume, many university introductory mathematics programs have recognized educational disparities and are currently exploring ways to better support students from structurally disadvantaged identities. Since many people in positions to make these important changes have themselves been successful within the existing systems, it is paramount that the students for whom these programs are intended to support are involved in identifying the issues they face and creating solutions. In a bottom-up approach, equity in mathematics is defined from the students' perspective, rather than based on outcome differences defined by the school administrators (a top-down approach). In this case study, we describe the creation of a supplemental precalculus course motivated by differences in outcomes and experiences between different populations of students, and implemented by giving students agency in its design. In this way, this course represents steps taken to shift to a more bottom-up approach to equity within a university introductory mathematics program.*

## 1. *Introduction*

At Colorado State University (CSU), precalculus is offered through five 1-credit courses in an online, mastery-based program currently called the Precalculus Center[30]. This atypical structure of precalculus was inspired by the idea that students should have the freedom to set their own pace, with the ability to choose how many and which courses they enroll in. While key aspects of how the system works will be highlighted here, more information can be found in the PAC*e* Student Guide.

[30] When completing this case study, the program had been named the Paced Algebra to Calculus electronically (PACe) Program.

Since these five precalculus courses are often used to satisfy the university's mathematics requirement and used as prerequisites to major classes, a large and diverse population of students take precalculus classes. Being fully online (with the exception of in-person proctored exams and optional in-person tutoring), the Precalculus Program allows for thousands of students to take precalculus classes; typically, a single semester has enrollment between 2000 to 3000 students. While the online learning environment is conducive to some, a lead adviser within the mathematics department identified an important population of students who were disproportionately not completing the first module in the precalculus sequence: women and students of color who expressed anxiety about math and/or the online format. These students were unable to take a subsequent precalculus course, and were prevented from enrolling in certain classes required for their major since they could not meet the mathematics prerequisite. In an effort to provide support to these students as well as to learn more about the challenges they were facing, a small, internal grant was awarded to create an in-person supplemental instruction where the students from this population would complete the first 1-credit precalculus course.

The population recruited for the course consisted of STEM-intending students of structurally disadvantaged identities. For this study, structurally disadvantaged identities include Students of Color, first-generation students, low-income students (as marked by Pell-eligibility), and women. We refer to these populations together as structurally disadvantaged to emphasize the effect that policies and systems in our educational programs have unintentionally and intentionally had on some populations (e.g., women, people of color, low-income people) and not on others (e.g., men, white students, upper-middle class people). It is well documented that women, students of color, low-income students, and first-generation students experience disadvantages in our universities and specifically in STEM, including microaggressions in their classes (Leyva et al., 2021; McNeill et all, 2020) and less access to high school math courses to prepare them for college coursework (Bressoud, 2021; Kolluri, 2018).

The rationale behind selecting our target audience was twofold. The focus on students in STEM arose because the mathematics prerequisites need to be completed in a certain timeline for students to be on the four-year track. The focus on students from structurally disadvantaged populations stems from the understanding that these students often come in with high math anxiety and lack a sense of belonging in STEM. Based on feedback from the advisor, research literature, and our own experiences, the research team conjectured that the Precalculus Program exacerbates the underlying anxieties experienced by these students. Moreover, the lack of support, connection, and communication in the online setting further isolates already disadvantaged students, making it extremely difficult for these students to reach out for help. The combination of these factors can easily lead to students disengaging and ending the semester with an "Unsatisfactory" on their transcript. To prevent these issues, the in-person course was purposefully designed to be responsive and flexible, allowing the instructors to adapt to the needs of the students as their needs were identified.

## 2. *Pilot Course Setup and Summary*

The first pilot of this course ran in Fall 2019 with 26 students pursuing STEM degrees, each of whom identified as a woman, a Student of Color, and/or first-generation college student. The instructional team consisted of 5 members: one tenure-track mathematics education researcher who identifies as a white woman (second author on this submission), one mathematics graduate student who worked as an Assistant Director of the PAC*e* center and identifies as a white woman (first author), one mathematics education undergraduate student who had completed Calculus I and identifies as Latino, one psychology undergraduate student who had previously failed the course herself and identifies as a mixed-race woman, and one undergraduate student who worked as a tutor for the main precalculus system and identifies as Latino. We were purposeful in creating a diverse instructional team--and sincerely acting as a team--due to the authors' understanding that our identities (as white women with/pursuing advanced degrees in math) may preclude us from understanding our students' needs in this environment.

Since the course was restricted to 1-credit, the course could only run for 50 minutes once a week. We strategically designed an additional 50-minute optional review session at the same time a different day, and made sure while advertising the course that students should keep that hour open. The course was designed so that 80% of students' grades would come from their performance in the Precalculus Program and the other 20% of their grade would come from our class' journaling requirement. Each week, the students were assigned a journal prompt, which we would then read, compile, and use to shape how the next class was run. Responses from these journal entries were also used to understand our students' needs and as data in this case study to share the students' perspectives on aspects of their experiences. Normally, the course would have deadlines split evenly throughout the semester; however, the in-person course had special permission to have a deadline schedule that started halfway through the semester to allow the instructional team more flexibility in creating a course that met our students' needs.

As an instructional team, we met once a week to discuss observations from both the classroom and the optional review session, important content for the week, and the weekly journal responses. In this discussion, we highlighted problems that had occurred or that we were anticipating and collaborated on solutions. This communication led to the class having a responsive structure, and we adapted both our beliefs about the goals of the course and the way the course was run. Our original teaching goals were to support the students in their development of mathematical knowledge, their change in mathematical affect, and their success in the course. In order to monitor our progress towards these goals, we planned on collecting feedback both formally with reflective journaling and informally with in class conversations. In the early weeks of the class, many students expressed their concern with the course; the main feedback, which will be discussed in further detail, that made us realize we needed to restructure our course design and implementation was that the students wanted a course to help them succeed in the Precalculus Program, not a course that focused on their mathematical development.

This student feedback helped us recognize that our implementation of these teaching goals did not align with the original goal of creating the course; in other words, to serve the students, the instructional team's job is to help the students pass the class so mathematics was not a barrier for them. This meant changing the teaching goals to be geared towards teaching test-taking strategies, different problem-solving methods, and how to use a graphing calculator. While these goals didn't align with supporting the development of mathematical knowledge as we conceive of it, they did support our students to be successful within the existing system.

In the beginning of the semester, we focused on explaining the intentions of why the course was created, delving into the students' mathematical confidence and anxiety, and working through prerequisite math skills needed to be successful in the first 1-credit precalculus course. We utilized a survey to find math concepts a majority of students seemed to struggle with and worked in small groups on these skills, each group being led by a different instructor. During these initial classes, we realized that some students were disengaged and that students were getting confused by our varying styles of how to solve a type of problem. In a journal entry, one student wrote "[t]here's just something about the class that isn't helpful to me at all, I get the concept and I think it's great but my classmates obviously are more ahead then I am and it really just makes me feel awful about myself," while another student said "I am happy there is a class about how to use [the Precalculus Program] and help us study for it, but as of right now we are not doing that."

After receiving this feedback, we realized that we had brainstormed how to help this population from different experiences and perspectives, but we had not let our specific students tell us what they were struggling with. With this in mind, we spent the fifth class period allowing the students to choose which instructor they wanted to have a conversation with about their sentiments on the course. Many mentioned they were afraid they were behind other students in the Precalculus Program that were not enrolled in our in-person section (since they had already started the content), they wanted to start the precalculus content and stop reviewing, and that they were interested in spending class time working in stations with other students working on the same section. Starting the sixth week and continuing for the rest of the semester, we spent our in-person classes having them grouped based on what required assignment(s) they were working on, so they could work with their peers and ask us questions as needed.

After the whole class conversation, there was never another journal response with a negative comment about the in-person class. In fact, most students had started benefitting from the new course structure as evidenced from data collected in student journals, in class conversations, and increased engagement with the course material, and were actively thankful for the opportunity. Some specific things they mentioned being thankful for include: the help provided from the instructional team, the relationships they formed with peers, and the accountability going to class provided them with. One student reflected "I'm being entirely genuine when I say I believe this was the most valuable class I had taken this semester. I have been taught patience, time management, self-responsibility, and, of course, math."

## 3. *Challenges Addressed in Second Iteration*

Based on the success of the first course in Fall 2019, the department allowed two more sections of these in-person precalculus courses to be run during Spring 2020; one course would be a combination of a prerequisite exam and the first 1-credit course, and the other would be a combination of the first two 1-credit courses. Because this was the second semester these courses ran, some challenges that were faced in the pilot course could be solved. For example, attendance became a major issue in the pilot class; the average attendance out of 26 students by the end of the semester was between 10-15 students. Another challenge was that students seemed to have a low level of engagement with email communication; often when a portion of students didn't show up, multiple targeted emails would be crafted to the different populations outlining strategies for how they could achieve success. Both of these problems were solved by making attendance a small portion of students' grades--this change made sense, since if a student is not benefitting from the in-person class, then they can simply take the online version. Now that most students were attending class on a weekly basis, communication could occur in-person, which has a stronger impact. Another thing the instructional team learned during the pilot was that a specific content unit in the course challenged many of the students more than the rest of the content. To combat this during the spring semester, students had an additional week to work on this material. One last problem was the disengagement of the students in the beginning of the term. Students were highly in favor of starting the content right away, most likely from a mix of anxiety and readiness to do mathematics, and because starting early means the deadlines can be spaced out more. So, in both spring classes, students began discussing mathematical material on the first day of classes to optimize time and to get the students in the mindset of doing mathematics.

## 4. *Implications and Recommendations*

Although in the spring we had information from the fall pilot about the population of students we were trying to serve, it would have been unrealistic to expect an identical copy of the course. The knowledge gained was extremely fruitful towards understanding the population we are continuing to target with these courses; however, each class has a unique dynamic that as an instructional team, we can only experience and adapt to versus try to prepare for. By giving the students agency in shaping the nature of the class, the instructional team can provide a course that is serving the students the best it can. To do this, it is important for instructors to express early on that they are not just willing, but expecting to change the structure of the course, so their feedback in this process is greatly appreciated. Instructors should give students this ability at any point, while also specifically prompting for their feedback through targeted journal entries, in-class conversations (both in one-on-one check-ins and small groups), and email conversations. Through this process, students can become agents of change: to take power in the classroom and create a space most conducive to them.

Regardless of the changes being made, if a program is geared to specifically support students from structurally disadvantaged populations, it must allow students to have agency over their learning. In this case study, we have described how we accomplished this agency at CSU within a supplemental instruction for our online precalculus course. We recognized that the mathematics department, the precalculus program, and the supplemental instruction were all run by white people with advanced mathematics degrees, which meant when we were creating how the supplemental instruction would run, we could not draw on our own educational experiences. Instead, to support students from marginalized populations who experienced mathematics anxiety, we needed to take a bottom-up approach and bring student voices and ideas into the course design. Otherwise, we risked creating a program, no matter how good our intentions were, that further perpetuated the inequities we sought to dismantle.

As a final note: due to the coronavirus, both classes in the spring shifted to fully online following spring break. Since the rapport between the students and the instructional team had already been developed, the use of online instruction was only needed to *maintain* these relationships. A future challenge that would be interesting to explore would be how to use online instruction to *develop* the relationships between the instructors and the students from the start of the semester.

# 13 Coordinating Courses Below Calculus Creating a Community of Practice

**Jacqueline Jensen-Vallin**
*Lamar University*

Jacqueline Jensen-Vallin is a Professor of Mathematics, interim department chair, and served for five years as the Director of the First-Year Mathematics Experience at Lamar University. She is also a former editor of MAA FOCUS, the newsmagazine of the MAA.

*Abstract: As a part of the State of Texas' requirement that state colleges and universities provide co-requisite options for most underprepared students (75% by Fall 2020), Lamar University has begun coordinating our first-year mathematics courses. For most of these first-year courses, college-ready and co-requisite students are comingled in the college-level courses, and co-requisite students are comingled from all college-level sections into co-requisite support courses. For the curriculum to run smoothly, we coordinate these courses on a very large scale. Weekly coordination meetings are held to discuss not only content, but also pedagogy in those courses. We work to create a community of practice, where we share information on the content and pedagogy, and search for ways to provide better support for all students in a first-year mathematics course. This article will discuss our model and results we've seen in creating our community of practice.*

## 1. *Background*

The State of Texas House Bill 2223 (HB 2223) made two sweeping changes to developmental mathematics (and English) and first-year mathematics (and English) courses. Both of these changes aim to decrease the time students spend in developmental courses, hoping to increase successful progression of students and shorten time to degree. At Lamar University, students who enter the university requiring foundational work are frequently Black/African-American and Latinx. While streamlining this path could affect many students, we expect significant positive effects for students at our minority-majority institution.

The first such change from HB 2223 is to require that public colleges and universities in Texas offer co-requisite courses for a subset of the students who would have previously taken stand-alone developmental coursework. In Texas, this requirement for developmental coursework is measured exclusively by the student's score on the TSI (Texas Success Initiative), a standardized test for students completing high school in Texas and entering college. In mathematics, students need a score of 350 on the TSI to be considered college ready. In 2018, Texas colleges and universities were required to have 25% of the non-college ready students in mathematics and English enrolled in co-requisite courses running parallel to their first-year college level course. In 2019, that increased to 50%; and in Fall 2020, we were expected to have 75% of those students enrolled in a co-requisite course.

The second large change is that students are eligible for state funding (and therefore the colleges and universities receive money based on the enrollment of these students) for no more than 9 hours of developmental coursework. This includes stand-alone and co-requisite coursework, as well as individually counting any retakes of those courses. For our department, this meant not only developing co-requisite courses to accompany our four first-year math pathways, but also shortening our developmental course sequence. Previously, we had offered three levels of developmental mathematics courses: Basic Math and Pre-algebra, Elementary Algebra, and Intermediate Algebra. Students might enter at any of these levels, but needed to complete the sequence before moving into a stand-alone college-level course.

These required changes lead to a complete redesign of all of our first-year courses, and a complete development of new curriculum for the co-requisite courses. In particular, we now have only one stand-alone developmental course, which is taken only by students entering one of our three algebraic pathways (college algebra, pre-calculus, and business math) and not needed for any students who are entering into the non-algebraic pathways (contemporary mathematics and introductory statistics).

While implementing these changes was stressful on a short timeline, we hoped that it would help students move more quickly through the first-year mathematics courses required for their degree and that it would greatly lower the number of students who get stuck in a cycle of remediation and never reach the first-year course.

**2. *Community of Practice***

Clearly there are very many challenges with such an initiative, especially one mandated by the state. One of our challenges was to create a position for a tenured faculty member (who would be named Director of First-Year Mathematics Experience (DFYME)) to lead this re-design. We wanted this position to be held by someone with some protection (from tenure) and seen by the administration to have some authority in the department. However, the DFYME designs and coordinates courses which are almost entirely taught by instructors. As a matter of fact, it's not common for tenure-stream faculty to teach first-year or developmental courses in our department. Also, previous to this, the instructors teaching college readiness were in a separate department and reported to their own chair, who was responsible for scheduling and evaluating those instructors.

As DFYME, it was important to me to create a community of practice – a group of people who care about their teaching and interact regularly to share ideas and improve the experience for everyone. We wanted to use this course design and redesign to integrate the new instructors into the department, take advantage of the strength and knowledge of those instructors (and those already in the department) who had been teaching college readiness and first-year mathematics courses, and create a community supporting continuous improvement of both the courses and our teaching techniques.

We had one year to design and create the updated courses. During that time, we took feedback from the development mathematics instructors. However, the extreme rewriting of the courses and pathways sometimes meant that we couldn't implement all suggestions. It's also important to note that we named each of the co-requisite courses "Fundamentals of…" to relate to its corresponding college level course. Both of these made forming the community of practice more challenging since some instructors felt this invalidated their expertise. Renaming the courses indicated significant change, and changing our approach to the courses sometimes meant abandoning previous methods. Since the Foundations courses had to coordinate with the college-level courses, the purpose of the course was changed from giving students a full and coherent story to supporting the work being done in the college-level course.

Course content was created by the DFYME, who built the homework assignments in the respective homework platforms, created due dates, and wrote a shared syllabus with synchronized pacing and similar grading expectations. While quizzes and exams were written individually by each instructor, the rest of the content was coordinated by the standardization of the commonly built homework assignments.

### 3. *Lead Instructors*

To combat this dichotomy between wanting a community of practice and having one person in charge of the curriculum, and to create an inclusive atmosphere, we initially tried having lead instructors for each of the pathway pairs. To that end, each group of college-level and complementary co-requisite instructors met three times per semester to discuss the pacing of the course, how the support courses were going, and ways to better support the students. It turns out that these meetings were not frequent enough, and there was often feedback that one instructor had skipped certain sections, while another instructor had covered everything and was three days behind schedule, while another was one week ahead. It was not a good community of practice and wasn't effective to coordinate the courses – we weren't doing enough to create regular interactions to support each other and to share ideas. Additionally, the lead instructor would report to the DFYME, who was not necessarily at the meetings. We hoped this would give people an opportunity to speak freely (which they did) but the meetings were not frequent enough to create a truly coordinated set of courses, nor a community around these courses, and the DFYME wasn't in the room to build trust with this community of instructors.

### 4. *Weekly Coordination Meetings with the DFYME*

The next iteration began in Fall 2019, and involved the DFYME meeting for 30 minutes each week with the groups of instructors for each of the following courses: Fundamentals of College Math (the only remaining stand-alone developmental course), College Algebra/Fundamentals of College Algebra (an algebra course for non-STEM majors), Precalculus/Fundamentals of Precalculus (a functions-first algebra course for STEM majors), Contemporary Math/Fundamentals of College Math (a liberal arts mathematics course), Introductory Statistics/Fundamentals of Statistics, and Precalculus II (a course for STEM majors focused on trigonometry). Set and standardized times were scheduled before the semester began, based on the teaching schedules of the instructors of each of the courses. Unfortunately, we were not able to compensate the instructors for their time either with a stipend or release time. We also didn't have a budget to provide food or snacks. The following set of expectations was circulated before the first meeting of the semester:

#### 4.1. *Norms for meetings*

- Recognize that everyone has expertise
- Be honest

- Share talk time
- Focus on what students can do, and how to help them reach the course goals.
- Look for solutions, not blame
- Focus on systems, not people

4.2. *Template for weekly meetings*

- What was your big success for the week?
- Did you meet the goals for the class this week (or will you by the end of the day)? If not, what topics do you have left to cover? What is your plan for "catching up"?
- What new technique or strategy did you try this week?
- What was your biggest challenge for the week?
- What are the big ideas for next week's content? How will you present those big ideas? What is the story we are telling with next week's content?
- Do you have a day/topic this week that you would like to have someone come observe?
- Other successes or challenges to discuss?

These guidelines were adapted from a workshop hosted by the Dana Center's Math Pathways (https://bit.ly/danactrpathways) consultants, who were most helpful in our design of appropriate norms and a focus for our meetings. It was very important to start from a place of success and then to move on to things that we needed to change or challenges we had faced. The goal here was to prevent the meetings from devolving into a complaint session, but instead for us to share good things as well.

These weekly meetings meant that every set of instructors had time with the DFYME every week, and since I was responsible for the pacing of the courses and the development of the curriculum, it meant that they had a forum for sharing concerns and suggestions. The norms meant that we recognize all of the expertise in the room, which substantially helped acknowledge the expertise of our instructors who were long-time teachers of the developmental curriculum and gave them a voice in a room with people who predominantly taught college-level courses.

These weekly meetings discuss not only content, but also pedagogy in order to create a community of practice among instructors of first-year courses. These meetings search for ways to provide better support for all students in a first-year mathematics course. This sharing frequently led to a deep conversation about how one instructor successfully covered a difficult topic in class that week. Since the implementation of these weekly meetings, we have seen changes in the department's ability to talk about teaching, gathered data about student success based on this synchronization, and created a forum in which we can discuss social justice issues in the classroom.

To supplement those course-specific meetings, we also implemented monthly teaching teas (now teaching luncheons), which also began in Fall 2018. These events allow faculty who are teaching across these first-year courses to gather simultaneously in a more convivial atmosphere. We had funding to provide snacks at these meetings. Topics for these gatherings include better use of the learning management system, homework systems, and the experiences of our incoming first-time in college students. Additional topics have included academic honesty, ways to support diversity in our classes, and new pedagogical methods that transcend individual courses. In one teaching tea, we specifically focused on ways we can make our classroom more inclusive, with some of us sharing specific techniques for inclusivity in class – this could be assigning groups for in-class work, making sure that students of color are not interrupted by their white peers, or taking extra time to learn names of all students. These larger meetings allow faculty who are not always in the same room for weekly course-specific meetings to discuss issues and common pedagogical concerns and innovations.

**5.** ***Observations and Results***

While lacking robust evidence, I can see in the department that the coordination and weekly meetings for the group of instructors involved in the first-year courses formed a community of practice, gave voice to many instructors, and transformed how we think about coordinated courses (people don't hate them as much as they originally did, I believe). Instructors attend every week, indicating that these sessions are valued, and when we have a new instructor join the fray, the experienced instructors talk to them about how helpful the weekly meetings are! These meetings also help instructors feel like their role in the department in teaching first-year students is important, when previously the focus in the department was on courses leading to STEM majors.

Additionally, we can investigate if this change affected student success – how are the pass rates in these redesigned courses for students taking the co-requisite courses and for those in the stand-alone courses? Overall, pass rates in the precalculus sequence declined during this implementation, but our university is lacking access to tracking grade data by ethnicity (other than to make periodic requests to institutional research for that data). We are continuing to track students who have taken the co-requisite and the stand-alone courses, and seek ways to better track our students of color as they progress, and to ask a number of follow-up questions:

- Did the passing rates decrease because of the increased number of co-requisite students in these courses?
- Did the decrease in passing rates indicate that we were creating more standardized expectations between sections? In other words, was the decrease because all instructors were now required to complete certain content affect the pass rates?
- Does the decrease in passing rates indicate that students who pass and continue to the next course are more likely to succeed? In other words, are we seeing an increase in the pass rates in the calculus sequence?
- What other supports can we build into these first-year courses that will increase the pass rates of the students as we have a larger percentage of underprepared students entering the co-requisite courses?

We are definitely moving towards a better community of practice, and being able to have difficult conversations with colleagues about pedagogy, equity, and course content will allow us to engage in continued improvements, leading to better experiences for our students and our instructors.

# 14 A College Algebra Intervention to Address Power Dynamics and Promote Students' Reasoning

**Heather Lynn Johnson**

*University of Colorado Denver*

Heather Lynn Johnson is a mathematics educator who investigates students' mathematical reasoning. She designs interactive, online tasks to provide students opportunities to expand their mathematical reasoning, and she teaches teachers to grow their students' mathematical reasoning.

**Gary Olson**

*University of Colorado Denver*

Gary Olson is a mathematician who teaches undergraduate mathematics and statistics courses and works to improve the professional development and teacher training for graduate teaching assistants. He is a passionate advocate for active learning and the use of TACTivities to promote reasoning and connections throughout the undergraduate curriculum.

**Amber Gardner**

*University of Colorado Denver*

Amber Gardner is a Pk-12 district math coordinator interested in critical and liberatory mathematics education. She works to design and evolve equity-centered math teaching systems that afford students' identity, power, access and achievement as empowered math learners.

**Xin Wang**

*RMC Research*

Xin Wang is a Senior Research Associate at RMC Research. Her research interest includes experimental and quasi-experimental designs, quantitative research and analysis, program evaluation, and educational technology. Xin currently serves as a Co-PI or lead evaluator for multiple research projects funded by NSF.

*Abstract: We describe a Professional Development (PD) intervention for College Algebra, to address classroom power dynamics and promote students' reasoning. We share the theory of change guiding the intervention, describe supporting materials for instructors, and discuss our efforts to promote sustainability.*

## 1. *Introduction*

The nature of students' learning experiences impacts their persistence in STEM majors (e.g., Seymour et al., 2019). College Algebra is a course that many students take on a path to STEM degrees, and it is fraught with challenges (Tunstall, 2018). Our focus is on a systemic instructional issue in early undergraduate courses such as College Algebra: An overemphasis on students' compliant answer finding, rather than their mathematical reasoning. We share a Professional Development (PD) intervention for College Algebra, to address classroom power dynamics and promote students' reasoning. We discuss how our theory of change informed our PD intervention, and we describe how artifacts, such as facilitation guides, can help to sustain practices beyond the intervention.

The year-long PD intervention consisted of face-to-face and online meetings, in conjunction with instructors' implementation of digital activities, occurred alongside broader department efforts to incorporate active learning and improve course passing rates in higher enrollment, lower division mathematics courses, including College Algebra. During our intervention, students of color were overrepresented in College Algebra. Per institutional data, 44% of undergraduate students were students of color; however, 63% of College Algebra students were students of color. Through the PD, we intended for instructors to learn to address classroom power dynamics via implementing digital activities, "Techtivities," designed to promote students' reasoning about relationships represented in graphs.

The editors for this volume have requested that authors address issues of diversity, equity, and inclusion in introductory math programs. We address a particular dimension of equity—power, which involves the voice and agency that students can have in a classroom, including who gets to do the talking and what gets to count as mathematics (Gutiérrez, 2012). For example, suppose students primarily encounter mathematical tasks that require them to find correct answers and fix mistakes. If students do have opportunities to share their reasoning, it is in service of explaining their answers. We argue that strategic implementation of mathematical tasks to promote reasoning, rather than answer finding, can serve as a catalyst to redistribute power in the classroom, to afford students more agency in their mathematics learning.

## 2. *Layering for Change: Positioning the PD Intervention*

To design and implement the PD intervention, we drew on a home-grown theory, developed by Heather Lynn Johnson, to instantiate a layering (Mahoney & Thelen, 2010) mode of institutional change. Operating from within, we worked to disrupt the status quo of College Algebra by integrating new elements, Techtivities, into the course. Because College Algebra was already packed with content, we anticipated that adding new elements might present challenges. Hence, we established relationships with stakeholders, including department chairpersons and course coordinators, to help the Techtivities to become part of the fabric of what counts in College Algebra.

Johnson's theory of change consists of four interconnected elements (DEEP): *Develop* innovative digital tasks (Techtivities) that privilege mathematical reasoning, rather than answer finding; *Embed* tasks within an existing course, connecting them to the high leverage content of function; *Extend* opportunities to instructors to examine their power to determine which student voices get counted (or marginalized); and *Provide* supports for instructors to implement tasks to promote students' reasoning. Figure 1 illustrates relationships between the four elements of our DEEP theory of change. Next, we explain how each element informed instructional practices in College Algebra.

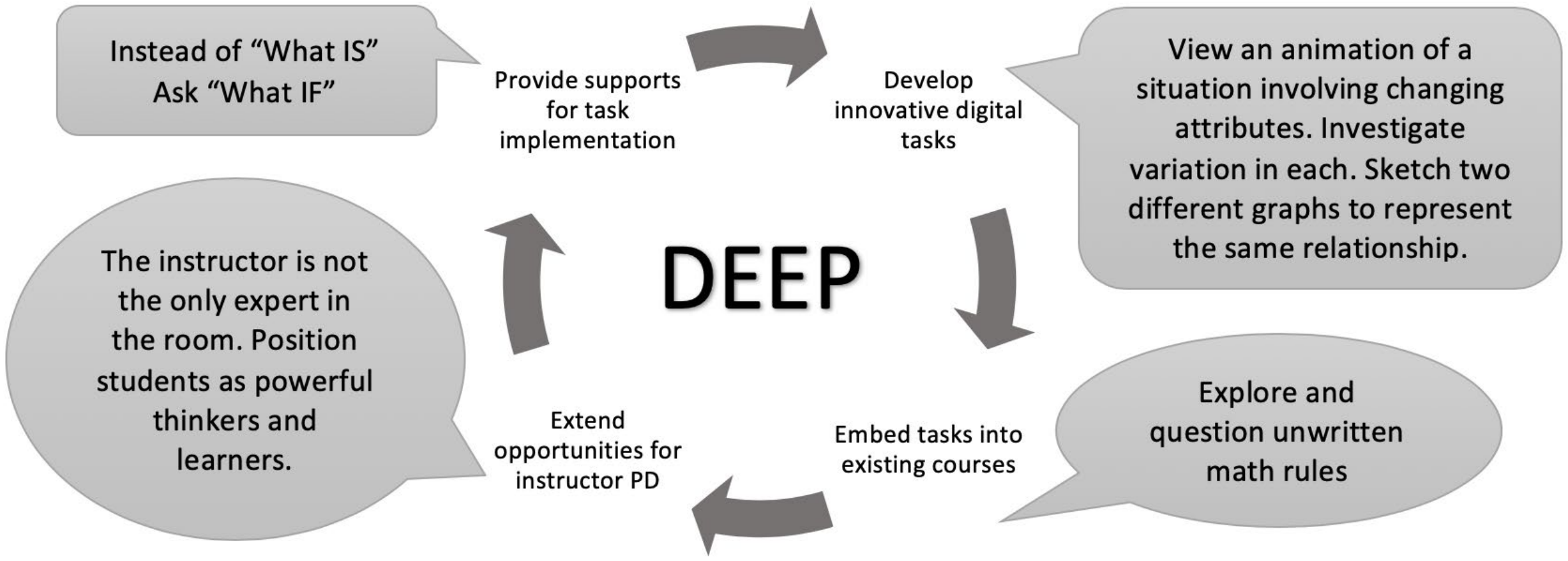

*Figure 1. The DEEP Theory of Change: Develop, Embed, Extend, Provide*

The Techtivities were a new element that we embedded into existing College Algebra courses. Developed in collaboration with Dan Meyer and the Desmos team, we designed a set of seven Techtivities (Johnson, 2018) to help students learn how graphs represent relationships between variables. One way we aimed to promote students' reasoning was through the kinds of questions we incorporated. In each Techtivity, students were to respond to a claim made by a hypothetical student. To address potential issues of bias, we selected student names to allow for gender ambiguity and phrased questions to center sensemaking, rather than judgments of correct or incorrect (Johnson et al., 2018). For example, one reflection question asked: "Chris says this graph shows height as a function of distance. Do you agree or disagree?" The graph in the activity represented height as a function of distance; however, distance was on the vertical axis and the graph did not pass the vertical line test. We wanted to make room for students to see sensibilities in another student's statement, even if it was not what they were expecting. Notably, students described connections between height and distance when responding to Chris's claim, which was useful not only for analyzing functions in the Techtivity, but also beyond (Olson & Johnson, 2022).

To position the Techtivities as something more than an "add-on" to an existing course, we linked them to key aspects of College Algebra content: Functions, Modeling, and Inverses. Yet, we recognized that innovative activities alone would be insufficient to change practice. If instructors implemented the Techtivities with a goal of fixing students' answers, rather than fostering students' reasoning, business would go on as usual, and the status quo would remain.

Through the PD intervention, we extended an opportunity for instructors to learn to distribute power in the classroom so that they would not be positioned as the sole experts in the room. In the face-to-face sessions, we engaged instructors in learning experiences to examine their own biases and discussed instructional routines that could make more room for student talk. In conjunction, we provided supports for instructors to implement the Techtivities. These supports included freely available facilitation guides (e.g., Olson et al., 2019), which situated the Techtivities in the larger scope of instruction, and included implementation suggestions and sample questions to promote students' active engagement.

**3.** ***The PD Intervention***

We enacted this intervention at a university located in the downtown of a large U.S. city in the Rocky Mountain west. If students intended to major in a STEM field at this university, and had not placed into a higher-level introductory course such as Calculus, they needed to pass through the College Algebra gate to pursue their goals. In this setting, College Algebra was divided into recitation and lecture components, with recitations facilitated exclusively by GTAs and lectures taught by either university faculty (typically nontenure line) or GTAs. Students attended recitations prior to lectures. Because recitations focused on active learning and student engagement, they were an ideal place to implement the Techtivities.

In Fall 2018 and Spring 2019, recitation facilitators implemented the Techtivities across all sections of College Algebra. In conjunction, we led a PD intervention to support implementation. In Fall 2018, only recitation instructors participated in the PD. In Spring 2019, both recitation instructors and lecture instructors participated. There were seven participants in Fall 2018, and six participants in Spring 2019 (three recitation instructors, three lecture instructors). Two recitation instructors and one lecture instructor in Spring 2019 were repeat participants. Each semester, the PD consisted of four face-to-face sessions and three online sessions. In the face-to-face PD, participants explored and/or debriefed the implementation of the Techtivities and engaged in learning experiences to investigate their power in the classroom. In the online PD, participants reflected on their implementation of the Techtivities via responding to Qualtrics surveys.

We designed the power-focused learning experiences to counter a systemic issue prevalent in introductory mathematics courses such as College Algebra: The notion that students just need to work harder to be successful (Gardner, 2020). We created space for participants to reflect on their own biases as instructors, because as instructors, they hold power to amplify (or silence) students' voices and to dismantle (or perpetuate) the status quo, which privileges some and marginalizes others. To do this, we identified aspects that comprise the status quo in College Algebra classrooms, such as which students get the most instructor attention and what counts as mathematics worth discussing, and we worked to make those aspects more visible. For example, in an exercise to uncover their own implicit biases, instructors reflected on which students they connected with, and whether those students were also most like themselves. One of our goals was for instructors to consider how their actions, intentional or unintentional, could impact students' opportunities to participate. For instance, if instructors tended to give more attention to certain students, those students who received less attention might feel as if they did not belong in mathematics, which could have ramifications for their achievement and persistence.

We offered the Techtivities as a tangible way for instructors to address power dynamics in the classroom. Because the Techtivities were designed to promote students' reasoning about graphs, we intended the implementation to make room for student talk that extended beyond what answers they found or if those answers were correct or incorrect. In the Techtivities, students received computer feedback on the graphs they sketch. The purpose was for students to explore and question, not to go back and fix their graphs so they look like the computer sketches. However, if instructors graded the Techtivities on accuracy, the exploration piece can get dampened. During a Spring 2019 PD, one recitation facilitator, Rachel, shared how she made efforts to maintain a focus on reasoning. She told students that they would be graded on effort, rather than accuracy, and she encouraged them to explore, without fear of being penalized for sketching an inaccurate graph.

Moves such as Rachel's support a redistribution of power, to shift the status quo in terms of what counts as mathematics. If students view answer finding to count more than other kinds of mathematical activity, they can miss opportunities to engage in mathematical thinking and reasoning, even when working on tasks designed to promote just that. Intentionally or unintentionally, instructors can position answer finding as the ultimate goal of mathematics, even though the work of mathematicians is far broader than that.

## 4. *Artifacts from the PD Intervention: Techtivity Facilitation Guides*

Artifacts from the PD Intervention include Techtivity facilitation guides (e.g., Olson et al., 2019). In the guides, we describe goals and moves to foster instructors' redistribution of power in the classroom, like what Rachel had discussed. For example, in each facilitation guide there are tips to amplify students' voice and position students as capable learners (see also Olson & Johnson, 2022). With these guides, we intend to provide supports that can extend beyond the PD intervention.

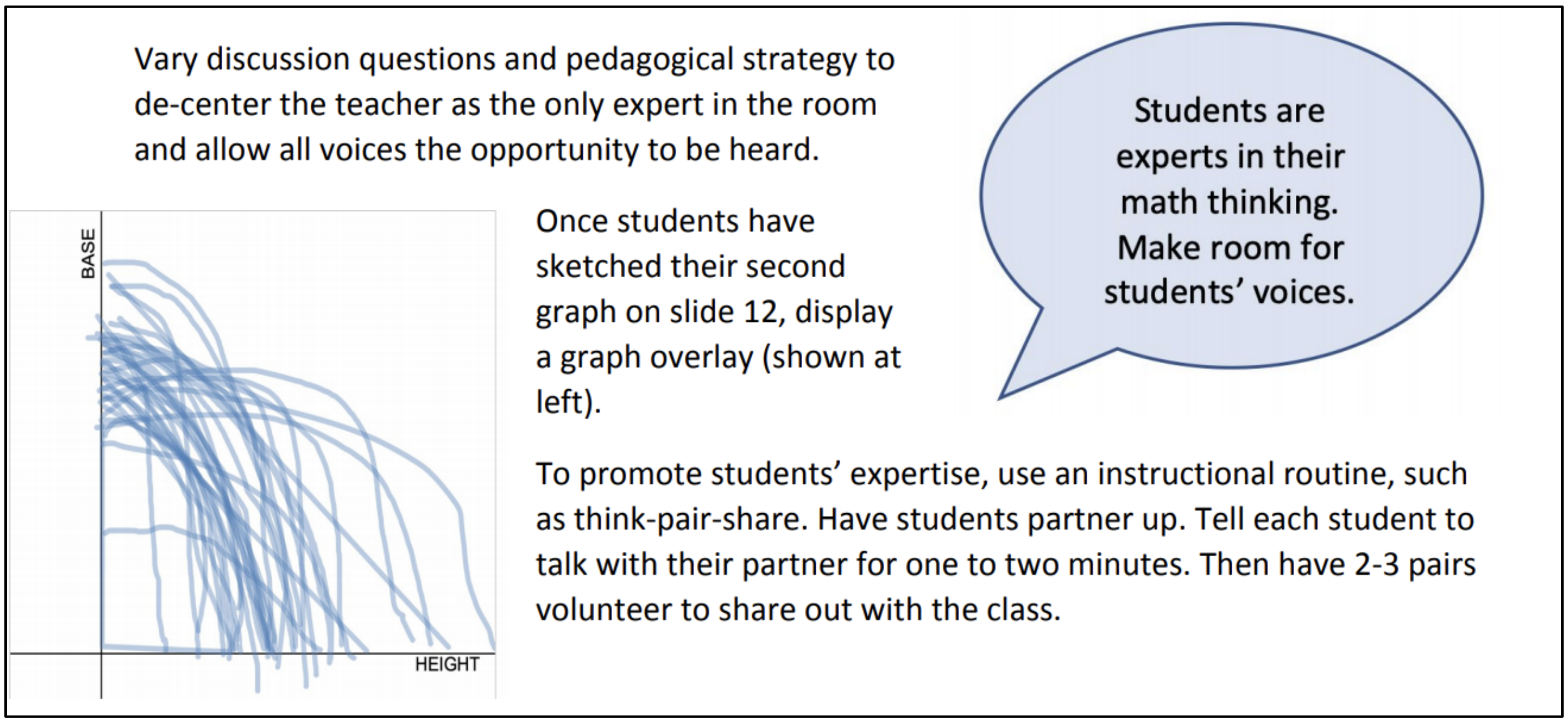

*Figure 2. Excerpt from the Changing kite and Dynamic tent facilitation guide (Olson et al., 2019)*

The facilitation guides have three main components. First, each guide begins with an overview of a key understanding; for example, points on Cartesian graphs represent a relationship between variables. The understanding is broad, to allow the Techtivities to connect to a range of content. Second, the guides provide specific ways for instructors to encourage talk and promote reasoning during students' work on the Techtivities. For example, there are discussion questions to promote students' small group conversations and instructional routines (e.g., think-pair-share) to foster students' mathematical thinking (See Figure 2). Third, the guides address ways to culminate the activity, in which students respond to a hypothetical student's claim. For example, instructors are directed to show different student responses and make room for students to ask questions and adjust their responses.

## 5. *Concluding Remarks*

All recitation instructors participating in the PD intervention agreed the Techtivities changed their own classroom instruction at least to some extent. They became more attentive to students' needs and feelings and cared more about students' deep thinking. Rachel's efforts provided evidence of how such changes can work to redistribute power in the classroom, to afford students more agency in their mathematics learning.

We address sustainability via operationalizing the DEEP theory of change (Develop, Embed, Extend, Provide) into existing university systems. Each year brings new instructors; artifacts such as the facilitation guides can support continued implementation. Furthermore, we have leveraged other structures, such as department orientations and seminars, to engage instructors in discussions around power dynamics in introductory university mathematics courses. Our next steps include developing professional learning communities, to promote a sense of belonging and ownership for instructors.

Instructional shifts to online and remote learning, brought about by the global pandemic, underscore the importance of mathematical reasoning. If online assessments and tasks focus only on answer finding, students can experience a mathematics of surveillance, monitored by proctors (virtual or human) looking for instances of academic dishonesty. Students in College Algebra need to experience mathematics that goes beyond just rote exercises to prepare them for some later math course. If opportunities to explore the breadth and depth of mathematics come only after students unlock gates that keep out far too many, everyone suffers.

## 6. *Acknowledgement*

This work is supported by U.S. National Science Foundation Grants DUE-1709903 and DUE-2013186. Opinions and conclusions are those of the authors.

# 15 Collaboratively Re-envisioning Calculus for the Modern Student

**Joel Kilty**

*Centre College*

Joel Kilty joined the Centre College faculty in 2009 and is currently an Associate Professor of Mathematics. He enjoys exploring the teaching and learning of Calculus, and developing innovative course materials to support the undergraduate curriculum. Joel loves spending time with his family, working in the garden, and spending time in the kitchen enjoying great food and wine.

**Alison Marr**

*Southwestern University*

Alison Marr is a Professor of Mathematics and holder of the Garey Chair in Mathematics at Southwestern University and is currently the co-director of the Enhancing Diversity in Graduate Education (EDGE) Summer Program. She enjoys doing research in graph theory and teaching mathematics at all levels including more interdisciplinary courses like her First-Year Seminar on television game shows. In her free time, Alison enjoys playing saxophone and traveling.

**Alex M. McAllister**

*Centre College*

Alex M. McAllister joined the Centre College faculty in 1999, where he has taught courses in mathematics, computer science, humanities, and philosophy. He currently holds the H.W. Stodghill Jr. and Adele H. Stodghill Professorship in Mathematics and he was awarded the MAA's Kentucky Section Teaching Award in 2015. Alex enjoys playing card games with his family and slowly renovating his 100+ year old home.

*Abstract: While traditional methods of teaching and learning Calculus have served some students well for many decades, national curriculum guides and problematic student outcomes highlight multiple changing realities and foster a greater awareness that asks the mathematics community to interrogate these traditional approaches. In response, the mathematics faculty at Centre College and Southwestern University (both small liberal arts colleges located in the southern United States) engaged in a thoughtful, collaborative re-envisioning project focused on evolving the entire calculus sequence to better meet the needs of the modern student. We decided to incorporate more modeling and realistic applications that utilize large data sets, include "new" ideas in each course for students who have studied calculus in high school, incrementally increase the challenge from one course to the next, and explicitly encourage persistence through the calculus sequence and beyond. In this paper, we discuss the development of this re-envisioned modern calculus sequence as an example of a successful curricular reform project, to include: (a) designing the courses with an eye toward inclusion and reducing barriers, (b) developing the courses independent of current textbooks, (c) making the hard decisions about what content stays, what content evolves, what content is added, and what content is let go, and (d) suggesting how to implement this process at other institutions.*

## 1. *Introduction*

While traditional methods of teaching and learning Calculus have served some students well for many decades, we have observed multiple changing realities that have forced us to interrogate these approaches. These changing realities include:

- evidence that standard approaches “filter” students of color, Pell-eligible students, and first-generation students (NSF, 2019),
- extremely different levels of student preparedness,
- the diverse career and continuing education paths of students,
- more sophisticated technologies, and
- access to large data sets that enable more realistic and more relevant applications.

In response, the mathematics faculty at Centre College and Southwestern University collaboratively engaged in a thoughtful re-envisioning project focused on evolving the entire Calculus sequence to better meet the needs of a broader swath of contemporary students, particularly Pell-eligible students, first-generation students, and students of color (for this paper that means American Indian, Alaska Native, Black, Hispanic, Native Hawaiian, and Pacific Islander students).

Centre College and Southwestern University are both private, undergraduate-only liberal arts colleges of around 1500 students, with six to nine tenured/tenure-track mathematics faculty. Both universities operate under a semester system and offer Calculus in a three-semester sequence. Most often, Business, Economics, and Computer Science majors, as well as some Psychology, Kinesiology, and Biology majors, exit the sequence after Calculus I; Chemistry and Biochemistry students typically continue through Calculus II; Physics and Mathematics majors and minors take the entire sequence. These courses are usually capped at 30 students.

## 2. *Motivation*

Our motivation to act was (and is) based on discussions at the national level, the cultures and data of our own institutions, and inspiration from other universities. On a national level, multiple studies have confirmed a multi-decade decline in the overall number of undergraduates choosing to pursue and complete science, technology, engineering, and mathematics (STEM) majors, particularly students of color and women students. Even more, just focusing on mathematics, the Mathematical Association of America’s (MAA’s) “A Common Vision for Undergraduate Mathematical Sciences Programs in 2025” offers a compelling call to action (Holm & Saxe, 2016). This 2015 report summarized the findings from studying seven curriculum guides from the leading mathematics organizations. The unanimous conclusion of all these guides is that “The status quo is unacceptable.” This call to action, in dialogue with our local context and our experiences as teachers, helped solidify our need to act and provided a framework for our approach to this re-envisioning project.

Over the last ten years the student populations at Centre and Southwestern have changed significantly, even though both are still predominantly white institutions. In Fall 2021, Southwestern enrolled 36.3% students of color (as previously defined), compared to 22.3% a decade before. Southwestern’s social and economic diversity has also increased over the past decade. In 2008, first-generation college students made up 7.6% of the entering class while in the last three academic years, first-generation students have comprised between 15 to 19% of the total student population. Finally, nearly one-third (31%) of current Southwestern students are eligible for Pell grants based on family income, a commonly used metric for identifying students with high financial need. Centre College has sought after and attained similar changes in its student demographics. In 2009, only 11% of Centre's students were students of color and, in 2014, only 9% were first-generation students. In contrast, the Class of 2023 consisted of 23% students of color, 20% first-generation students, and 23% Pell-eligible students.

Alongside these changes, we recognized that our introductory Calculus course was not serving these populations well. At Southwestern, the D/F/W rate for Calculus I was 22% over a recent five-year period, while that rate was 32% for students of color, 29% for Pell-eligible students, and 30% for first-generation students. This data was originally shared with the department as part of a grant writing process, and upon seeing this data we knew we had to act. At Centre, the D/F/W rate for Calculus I was 13% over a recent five-year period, and notably the rate was 10% for students of color. However, our D/F/W rate was 26% for Pell-eligible students and 18% for first-generation students.

In the light of this evidence both at our local institutions and on a national scale, we arrived at the discomforting conclusion that our standard approach to teaching Calculus "filters" students of color, Pell-eligible students, and first-generation students, rather than enabling their persistence and success. As such, we are obligated to respond. We must reexamine and revitalize our approach to teaching Calculus in concert with supporting the overall academic success of these untapped student populations.

On a local level, the cultures of our institutions include what might be described as a "constructive restlessness" with regard to curriculum design and effective teaching. Faculty are committed to a reflective thoughtfulness about the what and how of teaching with a particular eye toward enabling effective learning for all our students. And, more than just thinking about these questions, we are committed to action. It is also important to note that this action started at the department level. Before approaching partner disciplines, we wanted something concrete to compare and contrast to what was covered before. We give some specific examples from these conversations with partner disciplines in Section 3.

Finally, both of our institutions were inspired by the work done at Macalester College in converting their Calculus sequence to an Applied Multivariable Calculus sequence. Both Centre and Southwestern had connections with Macalester before this project. Centre had already created a version of the first course in the re-envisioned sequence that was modeled after the first course at Macalester and was taken by students who did not intend or need to take the whole Calculus sequence. Southwestern had recently undergone a 10-year department external review by a member of Macalester's faculty and part of the recommendations from that report were to reconsider our Calculus sequence. While we did not choose to fully implement the Calculus sequence developed at Macalester, there were elements of their sequence that we felt would help us improve the outcomes of our students of color, Pell-eligible students, and first-generation students.

## 3. *Development of Our Re-Envisioned Calculus Sequence*

Two of the primary goals of this Calculus re-envisioning project were:

- Reducing barriers to the persistence of students of color, Pell-eligible students, and first-generation students in STEM by creating improved pathways into and through the Calculus sequence that acknowledge the real and perceived advantages of some students, and mitigate the real and perceived disadvantages of others, and
- Improving student attitudes regarding both their confidence in their mathematical abilities and the perceived usefulness of mathematics in order to increase the number of students, particularly students of color, Pell-eligible students, and first-generation students, who persist through the Calculus sequence.

Toward achieving these goals, the mathematics faculty at both institutions actively collaborated to re-envision the Calculus sequence during a multi-day summer workshop funded by a grant from the Associated Colleges of the South, as well as through multiple conversations at national meetings and site visits to each college. During this workshop, involving almost all of the mathematics faculty from Southwestern and Centre, one of our main goals was to create a list of guiding principles for our approach to re-envisioning the Calculus sequence. Our process was highly collaborative, and no guiding principle was kept unless all faculty members from both institutions agreed. These guiding principles were intended to address the issues in the existing sequence which we felt contributed to educational inequity, build upon the framework developed in the MAA's "A Common Vision for Undergraduate Mathematical Sciences Programs in 2025" for curricular development, and provide a framework to help us make decisions on content as the process moved forward. We settled on the following list of guiding principles for our re-envisioned Calculus sequence:

1. The courses are organized in a way that promote educational equality and engage all students, with a focus on engaging students of color, Pell-eligible students, and first-generation students.
2. Based on a student's background, every course could be a good entry point into the sequence. Based on student interest, each course could be a good exit point from the sequence.
3. There is a smooth increase in difficulty across the sequence.
4. Some concepts and applications are intentionally layered in the sequence.
5. Every course incorporates activities where students utilize technology to further engage with the material.
6. Multiple perspectives are used to investigate topics. In the Calculus sequence, this includes differing perspectives such as numerical/analytical, discrete/continuous, and theory/application.
7. Every course has a persistent emphasis on multidisciplinary interactions.
8. Every course incorporates some theoretical disciplinary thinking.

After establishing these guiding principles, we began to explore various structures for a revised "modern" Calculus curriculum. Throughout this process we let our guiding principles lead the way and did not constrain ourselves to the typical presentation, topics, and organization found in standard Calculus textbooks. We split into three groups, one for each course, with each group including members from both schools. These groups developed proposals for what content to keep, what new content to add, and what content to remove based on our guiding principles. While each group focused on one course, we regularly shared out to the broader group to further develop ideas and determine how the different courses played together.

These conversations were not easy. Many hard questions were asked and discussed and, while we eventually all agreed to the course redesign, some faculty still had concerns. These larger group conversations and the acknowledgement of concerns were critical to landing on a final curriculum and to getting all faculty members on board and willing to try this new approach. Some of the conversations focused on our choice to move limits out of the first course. As mathematicians who learned and extensively taught Calculus with limits first, it was hard for us to imagine not starting with limits. In addition, we had many conversations about whether we were "watering down" the curriculum and about not forcing students to do hard algebra from day one in Calculus. We came to agree that significantly reordering the topics in the sequence would not result in a "watered down" curriculum, but instead would change the types of skills we develop in our students. While algebra skills are important, we believe those can be developed over time and throughout the sequence and should not be a barrier to completing the first course.

As mentioned before, the content for each course was also discussed with partner disciplines in parallel with the active decision-making conversations in our departments. Some content was added or moved to earlier in the course based on those conversations. In addition, Southwestern and Centre's topics list are slightly different based on these conversations and the need to serve our respective partner disciplines. For example, at Southwestern, computer science majors only take Calculus I and by removing limits from Calculus I the computer science faculty felt that students might be missing that important concept particularly as it relates to limits at infinity and their use in Big O calculations. Hence, we incorporated some discussion of horizontal asymptotes and what happens as our input values get large during our coverage of sigmoidal functions. The Chemistry Departments at Centre and Southwestern let us know that series were a vital part of their curriculum and since we were removing that from Calculus II (where their students stopped in the sequence) we might have a problem. After some discussion it was clear what they needed most was for their students to have seen Taylor polynomials. Hence, we included the topic of Taylor polynomials in both Modern Calculus I and II so students would have plenty of opportunity to understand the big ideas and could easily extend those polynomials to series in their Chemistry classes.

## 4. *Our Re-envisioned Calculus Sequence*

After multiple iterations, we settled on the following content for each course. This content list is still evolving, as is the order of topics, in response to our first experiences teaching these courses. A full day-by-day current list of topics for each course, day-by-day course learning objectives, and implementation details can be requested from the authors.

*Modern Calculus I*: The first course focuses on introducing students to modeling with various classes of functions using RStudio, dimensional analysis, derivatives, multivariable functions, partial derivatives, a graphical introduction to Lagrange Multipliers, and basic one- and two-dimensional integration.

*Modern Calculus II*: The second course begins with an introduction to differential equations, discusses approximation techniques with Euler's Method and Taylor Polynomials, introduces a few integration techniques, goes in depth with optimization of multivariable functions, covers double and triple integrals, integration over different coordinate systems, explores Monte Carlo methods of integration, and concludes looking at limits and their applications including continuity, L'Hospital's Rule, and improper integrals.

*Modern Calculus III*: The third course covers vector calculus, sequences and series, and looks at some of the theoretical ideas normally found in a first course including the formal definitions of derivative, limits, and integrals, introduces ideas related to the Mean Value Theorem and studies the proof behind the Fundamental Theorem of Calculus.

While the faculty at both Centre and Southwestern agreed on this basic outline of topics, the actual implementation details differ at our institutions as we seek to best serve our distinct populations of students. In addition, the course content and implementation details continue to evolve as we pilot-teach these courses and reflect on the experiences and insights of students and instructors to make improving adjustments. The primary goals of these adaptations is finding the best flow of ideas and the pedagogical approaches that will best enable student learning.

One common element to both institutions is the use of class worksheets to organize the class material and enable students to focus more on the course content and a bit less on taking notes during class. These handouts have been developed independently by the two institutions, and are shared among the faculty at each institution. Of course, instructors modify these handouts for their individual classes of students and to utilize the pedagogical approaches that they feel are most appropriate.

Having adopted a significantly re-envisioned calculus sequence, choosing an appropriate textbook to support student learning presents a challenge. At both Centre and Southwestern, we use *Mathematical Modeling and Applied Calculus* (Kilty & McAllister, 2018) for the first course in this sequence. For the second and third courses, Centre faculty use a traditional calculus text (hopping around sections in the order shared in the appendix) and class handouts; Southwestern faculty use OpenStax Calculus texts (Strang & Herman, 2016) and class handouts.

## 5. *Impact on Diversity, Equity, or Inclusion*

Ultimately, we hope our project helps create a more inclusive learning environment by providing a "more level playing field" for the students enrolled in Modern Calculus I. For example, many students from high-performing high schools, continuing-generation students, and white and Asian men have been exposed to the ideas of Calculus as part of their high school education before enrolling in a collegiate Calculus I course. This prior experience gives these students a distinct advantage as they work to solidify their understanding (in contrast with those trying to learn these ideas for the first time). In the classroom, the resulting atmosphere often makes those students without previous access to Calculus feel as if they do not belong or as if they are less able, because “all” the other students seem to “already know” Calculus (although, often, this “knowing” of Calculus can be more of a perception, than a reality).

Appropriate technology choices also help create a “more level playing field” for our students because the open-source software R and its integrated development environment RStudio are new to almost every student and free to all students, further reducing barriers for low-income students. With an eye towards inclusion, both Centre and Southwestern have set up a web server or virtual machine platforms that allow access to RStudio on any device. In a similar fashion, modifying the order and the list of topics studied in the Calculus sequence, particularly in Modern Calculus I, helps ensure that all students are learning new ideas. The result is a far more inclusive classroom atmosphere in which our students share a more uniform experience of engaging with and trying to learn the ideas and tools being studied. Even more, the inclusion of technology deeply enriches the experience of studying Calculus by allowing a substantive increase in both the quantity and the relevance of the “real-life” applications explored in class. We hope that these quantitative and qualitative differences create an increased desire for and greater subsequent follow-through on taking additional courses in the Calculus sequence.

Because our project is still in progress, we do not yet have data that supports the fulfillment of these aspirations. With the schedule of topics determined, we have now been focused on developing our approach to the content of each course. Specifically, we are working to address other important pieces of this re-envisioning project, such as including culturally relevant data sets throughout our unit on modeling and adopting culturally relevant pedagogies in the classroom.

In the direction of assessing the efficacy of this project, we are collecting attitudinal survey data from both students in the old (standard) Calculus sequence and students who complete the new (modern) Calculus sequence. At both the beginning and end of the term, students in each Calculus course are taking the “Attitudes Toward Mathematics Inventory” created by Martha Tapia (Tapia & Marsh II, 2004). In addition, we are monitoring student persistence data at both institutions in the old sequence and the new sequence. Our hope is to see improvement in these persistence rates for students of color, Pell-eligible students, and first-generation students. If we do not see the desired positive change for these groups, the faculty in our programs are committed to continue following the data and making further changes to our Calculus curriculum to achieve our guiding principles.

## 6. *Recommendations for Other Institutions*

First and foremost, we recommend that faculty at other institutions take a careful look at the data from their institution and work to develop a list of guiding principles for their curriculum re-envisioning. It is important to clearly identify the target audience of the revised courses and keep this audience in mind as decisions are made. Individual programs should then invest time reflecting on the best approach to meeting the challenges and capitalizing on the opportunities of their institution. However, cross-institutional collaboration can also be helpful in having discussions and helping all members of the department agree to this re-envisioning process. It is natural for faculty to learn about and consider models being used at other institutions, and we encourage this research and discovery process. However, we would also encourage institutions to not restrict themselves to existing models or textbooks currently available. Instead, envision new Calculus sequences that meet the needs of your students and provide more equal opportunities to all your students.

No matter how well designed, major curriculum changes to the Calculus sequence typically result in issues regarding transfer credit and student placement. Centre and Southwestern are both smaller, residential campuses with a relatively small population of students transferring into and out of our institution. Therefore, we choose to let our re-envisioning process be driven by a focus on what calculus sequence would best meet the needs of the overwhelming number of students taking these courses and who remain at our institutions. With that initial vision in mind, we then considered these other important issues. We acknowledge that, at many institutions, concerns about placement and transfer credit will be more pressing and must be considered at every step of the process.

At Centre, we place students based primarily on an in-house-developed placement test, and we re-interpreted the results of that test to place students in the most appropriate course in this new sequence. We continue to provide credit and placement for AP Calculus AB/BC scores in a manner similar to before implementing the revised sequence: - students with a 4 or 5 on the AP Calculus AB exam earn credit for the first course and are placed in the second course; students with a 4 or 5 on the AP Calculus BC exam earn credit for the first two courses and are placed into a traditional Calculus III course that we run alongside our re-envisioned Modern Calculus III course in the Fall Term. Transfer students are handled by the current department chair and the Registrar's Office on a case-by-case basis.

At Southwestern University, students without any Calculus credit are placed into the first course. Those with AP Calculus AB scores of 4 or 5 or transfer credit of Calculus I from another university are placed into the second course. The one space that causes some problems is those students with AP Calculus BC scores of 4 or 5 or credit for Calculus II from another institution. These students are given credit for a non-Calculus course in the math department (so still earning college credit) but are placed in the second course. Any student that takes all three classes at another university can transfer in the full sequence and be done with their calculus sequence at Southwestern.

## 7. *Acknowledgements*

This work was funded by a 2018 Associated Colleges of the South grant. We would also like to thank the members of both the Centre and Southwestern math faculty for all their contributions to this project and their willingness to try something new.

## 9. *Appendix: Detailed list of topics for each course*

| Topic | Modern Calculus I | Modern Calculus II | Modern Calculus III |
|---|---|---|---|
| 1 | Linear function modeling | Review of Calculus I | Vectors & dot product |
| 2 | Learning to work in R | Intro to differential equations | Cross product & gradients |
| 3 | Large data sets | Separation of variables | Planes & tangent planes |
| 4 | Transcendental functions | Slope fields | Space curves |
| 5 | Modeling with y=ex | Euler's method | Calculus of vector funcs. |
| 6 | Modeling with y=Cxn | Fund. Theorem of Calculus | Curvature & normal vec. |
| 7 | Modeling with sine | Substitution | Vector fields & flows |
| 8 | Sigmoidal modeling | Substitution & int. by parts | Line integrals |
| 9 | Sigmoidal modeling | Integration by parts | Conservative vector fields |
| 10 | Single-variable modeling | Taylor polynomials | Green's Theorem |
| 11 | Single-variable modeling | Lagrange polynomials | Parameterizing surfaces |
| 12 | Dimensional analysis | Numerical integration | Surface integrals |
| 13 | Dimensional analysis | Monte Carlo methods | Stokes Theorem |
| 14 | Intro to derivatives | Multivariable functions | Divergence Theorem |
| 15 | Interpreting the derivative | Partial derivatives | Review of Fund. Theorems |
| 16 | Basic derivatives | Gradients & Linear approx. | Sequences |
| 17 | Product & Quotient rules | Multivariable optimization | More on sequences |
| 18 | Chain rule | More optimization | Geometric series |
| 19 | Taylor polynomials | Closed set optimization | Divergence test |
| 20 | Global optimization | Lagrange multipliers | Comparison tests |
| 21 | Local optimization | Double integrals | Ratio & root tests |
| 22 | Concavity | More double integrals | Types of convergence |
| 23 | Multivariable functions | Triple integrals | Practice with series tests |
| 24 | Partial derivatives | Other coordinate systems | Power series |
| 25 | Vector operations | Change of variables | Taylor series |
| 26 | Gradients | Change of variables | New series from old |
| 27 | Lagrange Multipliers | Non rectangular triple integrals | New series from old |
| 28 | Accumulation | More multivar. integrals | Fourier series |
| 29 | Definite integrals | Spherical coordinates | Derivative proofs |
| 30 | Fundamental Theorem I | Monte Carlo methods | Definition of limit |
| 31 | More on antiderivatives | Limits (algebra,L'Hopital's) | Definition of integral |
| 32 | Fundamental Theorem II | Continuity/Squeeze Theorem | Mean Value Theorem |
| 33 | Integration in two-variables | Improper integrals | Fund. Theorem of Calculus |

# 16 Peer-Led Team Learning in Mathematics: An Effort to Address Diversity and Inclusion Through Learning and Leadership

**Janet Liou-Mark**

*NYC College of Technology*

Janet Liou-Mark was a Professor in the Department of Mathematics at New York City College of Technology for 22 years. She organized the Peer-Led Team Learning program at the college and served on the leadership team for many national grants. Her research focused on developing and evaluating programs that helped women and underrepresented minority and first-generation college students remain in school and successfully graduate with STEM degrees.

**Melanie Villatoro**

*NYC College of Technology*

Melanie Villatoro is an Associate Professor and Chair of the Department of Construction Management and Civil Engineering Technology at New York City College of Technology. Her outreach events target groups underrepresented in STEM and she seeks to increase diverse, qualified students entering STEM fields, particularly in engineering.

**Ariane M. Masuda**

*NYC College of Technology*

Ariane Masuda is a Professor in the Department of Mathematics at New York City College of Technology. She has extensive teaching experience, including teaching in Brazil and Canada. She is dedicated to improving student performance, and has been actively working on developing open educational resources. Her research interests lie in the area of Number Theory.

**Malika Ikramova**

*NYC College of Technology*

Malika Ikramova is a current graduate student in the Department of Nutrition and Dietetics at New York University. She graduated summa cum laude from the New York City College of Technology (BTech '19). Her research interests lie in STEM education and through her relevant experiences in research, she aspires to become a pediatric clinical nutritionist in an oncology setting.

**Farjana Shati**

*NYC College of Technology*

Farjana Shati is an Adjunct Professor in the Department of Mathematics at New York City College of Technology, where she graduated with a degree in Mathematics Education with academic honors. She has been part of the Peer-Led Team Learning Program since 2015 and aspires to work towards reducing world poverty and gender inequalities.

**Julia Rivera**

*NYC College of Technology*

Julia Rivera is an Adjunct Professor of Mathematics in the Department of Mathematics at New York City College of Technology and a continuing educator at Hostos Community College. She has been part of the Peer-Led Team Learning Program for many years and enjoys helping students. Her mentor, Dr. Janet Liou-Mark, inspired her to become a better educator for women in the STEM field.

**Victor Lee**

*NYC College of Technology*

Victor Lee is an Adjunct Professor of Mathematics in the Department of Mathematics at New York City College of Technology.  As a student, he was a Peer Leader and involved in the program for many years. His experience has led to his current placement as a graduate student in Mathematics Education at Brooklyn College.

*Abstract: The Peer-Led Team Learning (PLTL) model has shown to be an effective instructional method to support females, underrepresented minorities, and first-generation students in Science, Technology, Engineering, and Mathematics (STEM). The collaborative problem-solving setting, led by a peer leader, fosters learning that engages all the students. Six critical components are vital to the PLTL model: 1) the PLTL Workshop is integral to the course; 2)fFaculty is actively involved; 3) peer leaders are well trained; 4) the PLTL Workshop modules are challenging; 5) PLTL workshops are allocated time and space; and 6) there is institutional support. City Tech has implemented the PLTL workshops in selected foundation mathematics courses over the past five years because of the dismal pass and withdrawal rates. Overall results have shown that females, underrepresented minorities, and first-generation college students who actively participated in the PLTL workshops have higher course grades and lower withdrawal rates. Students are also afforded the opportunity to participate in the PLTL Leadership program. Through the PLTL Leadership program, females, underrepresented minorities, and first-generation college students (107 peer leaders in total) who have successfully completed their STEM degrees, are either in the STEM workforce or pursuing advanced STEM degrees. The PLTL model supports students who are academically disadvantaged and provides students with an opportunity to build their leadership skills and to create a pathway to graduate school.*

---

*Dedication: We would like to dedicate this manuscript in loving memory of Janet Liou-Mark, a role model and a champion for all who were fortunate enough to know her.*

---

## 1. *Introduction*

The Peer-Led Team Learning (PLTL) model is an effective instructional method to support females, underrepresented minorities, and first-generation students in Science, Technology, Engineering, and Mathematics (STEM). Each week, eight to ten students are grouped together in a collaborative problem-solving setting to work on an assigned module supported by a peer leader. The peer leader is a former student who has taken the course previously and is quite familiar with the course content. Six critical components are vital to the PLTL model: 1) the integration of PLTL workshop into the course; 2) Ffaculty is actively involved and supportive; 3) peer leaders are well trained in content and pedagogy skills; 4) the PLTL workshop modules are designed to be challenging; 5) PLTL workshops are allocated time and space; and 6) there is institutional support. The implementation of all the components provides a robust pathway toward undergraduates’ academic success in their STEM courses using this model.

New York City College of Technology (City Tech) is a Hispanic-serving institution where 34 percent of the students identified themselves as Hispanic and 29 percent as African Americans. Sixty-one percent report a household income less than $30,000, 62 percent are the first in their family to attend college, and 80 percent of incoming freshmen receive need-based aid. City Tech has implemented the PLTL workshops in selected foundation mathematics courses over the past seven years because of the dismal ABC pass rates. Each semester, one to three sections of each mathematics course has a designated PLTL section with an extra problem-solving hour per week with a peer leader. These mathematics courses are identified as College Algebra and Trigonometry, Precalculus, Calculus I, and Calculus II. The average ABC pass rates for these courses have been approximately 50 percent. The overall withdrawal rate was approximately 15 percent which is consistent across all the foundational mathematics courses.

### 2. *The Critical Components of the Peer-Led Team Learning Model*

The six critical components were used as the basis in adopting the PLTL model in the mathematics courses. Each component will be described in detail. The first critical component is that the PLTL workshop is integrated with the course. Undergraduates who registered for the PLTL sections have an extra hour attached to the course, similar to a zero-credit lab session. Peer leaders are assigned to a group of eight to ten students whom they work with until the end of the semester. The peer leaders are responsible for taking attendance and assigning a completion grade for the modules. A grade of 4 is given if the module is completely finished, a grade of 3 indicates a 75 percent competition and so forth. The goal of the PLTL workshop is for the students to develop meaningful discussions and to work collaboratively in solving the modules. The peer leaders are instrumental in helping the students find the solutions by providing hints, but the peer leaders are cautioned not to readily give out the answers which is different from tutoring. Moreover, the peer leaders do not grade each individual module, but they are cognizant that the group has the correct answer.

The second critical component is faculty involvement and support. Though the faculty may not be physically present during the workshop, they have an important role in making the workshops mandatory. They intentionally assign workshop participation as a percentage of the students' final grade and ensure that each weekly workshop aligns with the current lecture content. On the average, faculty would allocate 5-15 percent of the final grade for participation. Without the support of the faculty, students would start to wane from attending the workshop.

The third component is that peer leaders are well trained in content and pedagogy skills. There are two elements of peer leader training, a two-hour content orientation and a one-credit *Peer Leader Training in Mathematics* course. For content skill training, peer leaders are given the 12 modules before the semester begins and participate in a two-hour content orientation. Since they are not given an answer key, they must complete the modules before attending the content portion of the orientation that takes place a week before classes begin. Essentially by completing the modules ahead of time, they review the contents and anticipate the common problem areas their students may face. During the two-hour content orientation, their finalized answer keys are vetted and produced. For each module, meaningful discussions on the prior knowledge required, alternate solving techniques, and anticipated bear traps evolve during this time. For the pedagogy skill training, first-time peer leaders register for a one-credit *Peer*

*Leader Training in Mathematics* course offered by the mathematics department. The course content covers a brief overview of learning theories on behaviorism, cognitivism, and constructivism, facilitation strategies, and a session on stereotypes and diversity and inclusion. During the instruction portion of the orientation, peer leaders are coached on how to start the workshop with an icebreaker and how to set the guidelines and expectations so that the future workshops become seamless. Content and pedagogy skills are covered in these two elements; however, the practical application in the workshops is the best training the peer leaders can experience.

The fourth component is that the PLTL workshop modules are designed to be challenging. Since the modules are intended to be solved collaboratively, the problem may have different ways to be solved. This intentional design invites thought-provoking discussions, fosters healthy debates, and at the same time, it increases the students' problem-solving abilities. Depending on the mathematics course, there may be as few as two problems or as many as six to solve during the one-hour period. Two example problems have been included below.

> PEER-LED TEAM LEARNING MAT 1475 (Calculus I) Sample Problem – Foundations
>
> Sketch a graph of the distance traveled by a student over time from the first floor to the seventh floor of City Tech's elevator. Assume that the elevator stopped only on the $4^{th}$ and $6^{th}$ floors and spent approximately 3 minutes at each stop. Carefully state the assumptions you made. Does your graph represent a function?
>
> PEER-LED TEAM LEARNING MAT1575 (Calculus II) Sample Problem – Taylor Polynomials
>
> Compute $p_5$(x) at the given value for the following functions. (a) $f(x) = sin(x)$ at $a = 0$ (b) $f(x) = 1/(1-x)$ at $a = 0$ (Hint: Think about the geometric series from MAT1375.) (c) $f(x) = ln(1-x)$ at $a = 0$ (Hint: Compare the first derivative of this problem to the previous one.)

During the workshop the peer leaders encourage the students to work together to solve the problem and discuss their thought process as they work. Peer leaders are taught strategies for engaging different types of students in the workshop process. Their primary goal is to have the groups work together to develop a plan for solving the problem.

The fifth component states that PLTL workshops are allocated time and space. As discussed, undergraduates who registered for the PLTL sections have an extra hour attached to the course, similar to a zero-credit lab session. This allows for the student to include the workshop in their weekly schedule and report to the same classroom each week. Having a dedicated time and space for students to work is an ideal setting for students to freely speak and ask questions amongst their peers, including the peer leader. Ideally, a healthy community is formed where the advanced learners are instrumental in fostering understanding among their peers.

The last component urges institutional support. The institutional buy-in would assist in providing professional development for faculty on how best to mentor the four peer leaders assigned to their course. As a comprehensive institution, there are no graduate students to assist the faculty. This model uses undergraduates so the expectations are slightly different from what is required from graduate students. This extra mentorship from the faculty is valuable especially for peer leaders who are underrepresented minorities and the first in their families to go to college. In addition, the institution can assist in programming issues such as scheduling the additional one-hour workshops. More importantly, the institution will provide the stipends for the peer leaders. Typically, the peer leaders are given a $300 stipend for each course they peer led. This job may allow low-income and first-generation college students to stay on a commuter campus, so they can benefit from the academic support and extracurricular activities offered by the institution. Otherwise, these students will mostly seek off-campus jobs to fund their education and/or support their families.

Faculty support and proper allocation of time and space are two of the critical components that can contribute to reducing the impact of stereotype threat for female and underrepresented minority students. Stereotype threat is an individual's expectation that negative stereotypes about his or her member group will adversely influence others' judgments of his or her performance and that a poor performance will reflect badly on the member group. Stereotype threat can contribute to lower performance of minorities in STEM courses (Deemer et al., 2014; Meador, 2018; Rivardo et al., 2008; Starr et al. 2019). These two components ensure that all students attend. If the faculty allowed flexibility in workshop attendance, those who can truly benefit from the extra help might not attend because of stereotype threat. Moreover, the dedicated space and the small group settings provide a more personable environment in which stereotype threats are less likely to affect academic performance (Rivardo et al., 2008).

### 3. *Grade comparison*

Grades have been collected from the Fall of 2015 through the Spring of 2019 PLTL sections. Table 1 compares the PLTL sections using the Spring 2019 institutional data (all sections) as a comparison. Although the pass rates have been relatively the same, the failure and withdrawal rates of PLTL sections are smaller. Students can complete (ABCDF) their mathematics classes.

A small study was conducted using the Fall of 2015 through Spring of 2018 PLTL sections. The final course grades (out of 4.0) and the uniform departmental final grades (out of 100) were recorded for all the PLTL participants. Participants whose gender, first generation status, and ethnicity was not identified were excluded from the study. First generation college students are defined by being the first in their families to attend college. An analysis for African-American and Hispanic/Latino students determined no significant difference between the two populations, therefore the data is presented for underrepresented minorities, which includes African-Americans, Hispanics/Latinos, Native-Americans and Pacific-Islanders as one subgroup. Non-minorities include Caucasians and Asian-Americans. Results have shown that, among the participants of the PLTL workshops, females do overall statistically significantly better than their male counterparts especially in the lower-level mathematics. In addition, first-generation college students tend to do better than non-first-generation college students; however, the results were not statistically significant. Moreover, non-minority students did statistically significantly better than the underrepresented minority students. These positive results were attributed to the collaborations and community formed with their peers, being able to freely discuss solutions, and having a dedicated safe space to work out problem sets. Table 2 summarizes the results.

*Table 1*
*Grades of PLTL Mathematics Sections and 2019 Institutional Data*[31]

| | College Algebra and Trigonometry | | Precalculus | | Calculus I | | Calculus II | |
|---|---|---|---|---|---|---|---|---|
| | PLTL Sections | Spring 2019 Institutional Data | PLTL Sections | Spring 2019 Institutional Data | PLTL Sections | Spring 2019 Institutional Data | PLTL Sections | Spring 2019 Institutional Data |
| | (n=642) | (n=1592) | (n=662) | (n=926) | (n=349) | (n=622) | (n=778) | (n=385) |
| **Pass rates: ABC** | 54.0% | 39.9% | 57.3% | 47.9% | 57.3% | 53.7% | 53.2% | 43.9% |
| **Pass rates: ABCD** | 67.4% | 53.9% | 71.9% | 61.6% | 67.6% | 62.9% | 67.7% | 55.7% |
| **Withdrew** | 12.1% | 21.0% | 9.7% | 14.0% | 18.3% | 16.0% | 12.9% | 21.0% |
| **F** | 20.4% | 24.2% | 18.3% | 23.9% | 13.8% | 21.8% | 19.4% | 22.7% |
| **Others** | 0.0% | 0.0% | 0.2% | 0.0% | 0.0% | 0.0% | 0.0% | 1% |

## 4. *Final Considerations*

Workshop participants respond to student satisfaction surveys before and after their workshop experience. Survey questions ask about their satisfaction with the efficacy of the workshops including group dynamics, materials used, and time allotted for each activity. A Likert scale with one indicating “strongly disagree,” with three as “neutral,” and with five as “strongly agree” was used. Underrepresented minorities responded that the PLTL sessions were very helpful in understanding the materials, in providing a comfortable space for asking questions, and in succeeding in the course. The PLTL program can be categorized as a positive contributor to underrepresented minority retention in STEM courses. Determining the factors that positively contribute to minority STEM major recruitment and retention may allow for the elimination of stereotype threats that hinder academic success for minority students in the STEM fields (Meador, 2018).

The benefits of the PLTL instructional model are two-fold. First, it is designed to provide additional assistance for the students enrolled in fundamental mathematics courses, and secondly, the peer leaders are seamlessly enrolled into the PLTL Leadership Program. Through this leadership program, a community of practice consisting of new and experienced peer leaders is constantly forming and strengthening. This support network is where mathematical and general knowledge are shared, mathematical and leadership confidence are continuously building, the pursuit of STEM advanced degrees is encouraged and supported, and informal mentoring is practiced.

*Table 2*
*Means (Standard Deviations) and Z-test Results by Mathematics Courses*

[31] The “Withdrew” row includes W, WN and WU grades. W means Official Withdrawal (no academic penalty). WN means Unofficial Withdrawal (never attended; no academic penalty). WU means Withdrew Unofficially (participated at least in an academically related activity; no academic penalty). The “Others” row indicates course grade not given or incompletes.

| **Mean (SD)** Sample size | **Females** | **Males** | **Z-test** | **First Generation** | **Non-First Generation** | **Z-test** | **Underrepresented Minority** | **Non-Minority** | **Z-test** |
|---|---|---|---|---|---|---|---|---|---|
| **College Algebra and Trig. Course Final Grades** | 2.21 (1.29) *(n=132)* | 1.95 (1.36) *(n=249)* | z=0.03, p=0.06 | 2.25 (1.26) *(n=109)* | 1.97 (1.37) *(n=225)* | z=1.82, p=0.07 | 1.83 (1.33) *(n=223)* | 2.35 (1.32) *(n=146)* | z= -3.72, p<0.001* |
| **College Algebra and Trig. Final Exam Grades[32]** | 67.55 (22.11) *(n=132)* | 62.51 (23.98) *(n=249)* | z=2.06, p<0.05* | 67.44 (22.99) *(n=109)* | 63.56 (23.90) *(n=225)* | z=1.43, p=0.15 | 60.81 (23.94) *(n=223)* | 69.44 (21.64) *(n=146)* | z= -3.59, p<0.001* |
| **Precalculus Course Final Grades** | 2.652 (0.993) *(n=109)* | 2.01 (1.27 *(n=239)* | z=1.96, p<0.001* | 2.35 (1.2) *(n=124)* | 2.14 (1.21) *(n=224)* | z=1.96, p=0.12 | 1.98 (1.2) *(n=214)* | 2.59 (1.16) *(n=134)* | z=1.96, p<0.001* |
| **Precalculus Final Exam Grades** | 72.54 (18) *(n=109)* | 65.2 (21.05) *(n=239)* | z=1.96, p<0.001* | 68.9 (21) *(n=124)* | 66.6 (19.84) *(n=224)* | z=1.96, p=0.33 | 63.86 (20) *(n=214)* | 73.24 (18.9) *(n=134)* | z=1.96, p<0.001* |
| **Calculus I Course Final Grades** | 2.37 (1.38) *(n=103)* | 2.06 (1.20) *(n=238)* | z=1.96, p<0.05* | 2.24 (1.29) *(n=116)* | 2.15 (1.24) *(n=203)* | z=0.61, p=0.54 | 2.06 (1.22) *(n=163)* | 2.39 (1.24) *(n=131)* | z= -2.24, p<0.05* |
| **Calculus I Final Exam Grades** | 65.44 (22.06) *(n=103)* | 63.35 (21.11) *(n=238)* | Z=0.81, p=0.42 | 63.16 (22.10) *(n=116)* | 64.68 (20.84) *(n=203)* | z= -0.60, p=0.55 | 61.82 (21.42) *(n=163)* | 68.40 (20.10) *(n=131)* | z= -2.71, p<0.01* |
| **Calculus II Course Final Grades** | 2.36 (1.24) *(n=33)* | 2.27 (1.22) *(n=97)* | z = -0.36, p=0.72 | 2.36 (1.27) *(n=40)* | 2.24 (1.25) *(n=83)* | z=0.46, p = 0.65 | 2.01 (1.24) *(n=61)* | 2.14 (1.17) *(n=47)* | z= -2.27, p<0.05* |
| **Calculus II Final Exam Grades** | 70.15 (19.27) *(n=33)* | 70.28 (15.71) *(n=97)* | z= -0.03, p=0.97 | 71.70 (15.68) *(n=40)* | 68.47 (18.33) *(n=83)* | z=1.01, p=0.31 | 2.01 (1.24) *(n=61)* | 2.14 (1.17) *(n=47)* | z= -2.27, p<0.05* |

Peer leaders are trained in content and pedagogy skills. These pedagogy skills can be transferred to other subjects and disciplines. At City Tech, the PLTL program has expanded beyond math courses and has been implemented in science and engineering courses. Peer leaders have had the opportunity to lead workshops across disciplines. The benefits of peer leaders applying their skills to multiple disciplines include an increase in confidence, reinforcement of their existing content knowledge, and expansion of their influence circle for serving as a role model.

> *"Being a Peer Leader in different subjects has helped me grow not only as a facilitator but also as a leader. Always be prepared. The PLTL workshop format is the same for all these subjects. Although the concepts are different, the math does not change. Another important aspect to be a good facilitator is the willingness to do it. In my*

32 The final exams for the courses listed on Table 2 are uniform and prepared by a department committee.

*case, I wanted to be good at it, I wanted to share my ideas and how I get around the problems. It is a beautiful thing to see students understand a new concept."*

*-Construction Engineering Technology student and Peer Leader*

The PLTL program equips and empowers Peer Leaders with the leadership and confidence to work with students and faculty across the college.

Moreover, the students not only benefit from participating in the PLTL workshops, but they are also afforded the opportunity to participate in the PLTL Leadership program which formally trains students to become peer leaders. Studies have shown that students enhance their perceptions and confidence in STEM when they have role models who look like them (Hutton, 2019; Shin et al., 2016; Weber, 2011). Through the PLTL Leadership program, females, underrepresented minorities, and first-generation college students (107 peer leaders in total) who have successfully completed their STEM degrees, are either in the STEM workforce or pursuing advanced STEM degrees. Approximately 35% of the graduates have completed or are continuing with their masters or doctoral degrees in STEM. To address the issues of diversity, equity, and inclusion in the mathematical sciences, instructional models such as PLTL may support this effort. The PLTL model supports students who are academically disadvantaged, and it also provides students with an opportunity to build their leadership skills and to create a pathway to graduate school, because there are not many role models who look like them in the mathematical sciences.

# 17 Creating Humane Exam Structures in College Calculus Courses

**Kelly MacArthur**
*University of Arizona*

Kelly MacArthur (she/her) is currently a Postdoctoral Research Associate at the University of Arizona, after completing her Ph.D. in Mathematics Education from the Mathematical Sciences department at Montana State University. Her dissertation research centered the voices of STEM majors from historically marginalized groups to explore their views on humanizing undergraduate mathematics courses. The project reported on in this case study took place when she was an Associate Professor Lecturer at the University of Utah where she taught collegiate mathematics for over two decades.

*Abstract: Taking up a call to action given by Rochelle Gutiérrez (2018) to rehumanize math in order to move toward an equitable and inclusive learning environment for students, I changed assessments in my undergraduate Calculus II courses, to include small group discussions between students for a portion of each exam. This was additionally motivated by the MAA Instructional Practices Guide statement that "teaching and learning are forces for social change" (2018) and to align assessment strategies with active learning pedagogy. This assessment change, along with a consistently enacted classroom mission statement geared at both engaging students and creating a humane learning environment, produced higher exam scores, improved student confidence and increased students' sense of belonging. In particular, I explored how Black, Latinx and Indigenous students, who are the intended targets of the rehumanizing framework, and women experienced this exam change, given their historical underrepresentation in college calculus courses.*

## 1. *Why make a change?*

There is ample evidence to suggest that active learning classrooms, which include discussions and mathematical exploration among students, produce better learning results for our students than traditional lecture format (Freeman et al., 2014). As instructors, we want to create class structures that build rapport, increase confidence & sense of belonging and lessen anxiety among students through interactions with their peers and instructor. We want students to actively participate in building a more equitable and inclusive learning environment where people of diverse backgrounds are valued in their mathematical thinking and sharing. With these goals in mind, there has been much research, and subsequent impact, on inclusive teaching practices for what happens during class. However, we haven't yet seen widespread research or changes in teaching that put these same principles to play in assessments. This is a mismatch pedagogically, and an opportunity for growth. Exams in our college calculus courses typically comprise the highest percentage of a student's grade in the course and thus they are high-stakes and meaningful events for calculus students. Yet even as we have moved towards more active learning environments, our exam structure has mostly remained unchanged from the standard silent, solo testing model. It's time to extend our equity and inclusivity-focused efforts to high-stakes assessments.

Answering a call to change, given by Rochelle Gutiérrez, to rehumanize mathematics (Gutiérrez, 2018) and in line with the MAA Instructional Practices that “teaching and learning are forces for social change” (MAA, 2018), I changed assessments in my large-enrollment university Calculus II courses to include small group discussions between students for a portion of the exam. The goal of this structural change was to make the exam process itself more humane and inclusive, as well as to have the exam process more closely match what was already happening during regular class time, and then explore the impact on student learning. I adhere to the emphasis put forth by Gutiérrez on centering the work of rehumanizing classrooms, with the intention that this will move us towards equity as a peripheral goal. The central focus is creating a rehumanized experience for students, with a primary focus on Black, Latinx and Indigenous students.

> *”Unlike ‘equity,’ which can seem to represent a destination, ‘rehumanizing’ is a verb; it reflects an ongoing process and requires constant vigilance to maintain and to evolve with contexts. Moreover, rehumanizing is an ongoing performance and requires evidence from those from whom we seek to rehumanize our practices that, in fact, the practices are felt in that way.” (Gutiérrez, 2018, p. 4)*

## 2. *What was the change?*

After explaining to students that I felt strongly about aligning our exam structures with the consistently implemented discussion-oriented practices used in class, the exams were split into two pieces, namely a group portion and a solo portion, administered over two time periods. The group portion of the exam was administered first, with the intention to increase confidence and decrease anxiety for both portions of the exam. Additionally, students then had an opportunity to learn from each other during the group portion of the exam, which had the potential to improve their performance on the solo portion of the exam. Each student turned in their own group exam (earning their own individual grade on the group exam), and they are welcome to agree or disagree with their group members on the group exam questions. The group exam questions were more complicated or group-worthy types of questions, while the solo exam questions focused primarily on computational skills.

As the instructor, I created and assigned groups of 3-4 students for each exam. For the first midterm, I randomly created groups and then adjusted groupings to consider information I had about student needs, personalities, gender, anxieties, etc. being sure not to have any groups with only one woman in it, for example. For subsequent exams, I ordered the students' weighted average of the previous exam scores and divided the students into three groups according to scores. I then randomly selected one student's name from each group to form groups of three. Finally, I again attended to information I had about student needs, personalities, gender, anxieties, etc. to make small changes to those groups.

Over the last two years of administering exams in this way, there have been some slight updates made along the way to improve this structure from its original implementation, which have all been based on student feedback. These changes include (a) having students sit with their exam groups for two to three class periods preceding the group exam, (b) asking a few short, non-graded Likert scale questions (see Figure 1) on the solo portion of the exam, and (c) giving students 15 minutes of solo time at the beginning of the group exams, rather than

**Figure 1**

*Non-graded Likert Scale Exam Questions*

* Likert-scale questions asked at end of solo portion of exam:
(a) My group mates for the group exam seemed prepared for the exam and did their best to contribute.
(circle one)

| 1 | 2 | 3 | 4 | 5 |
|---|---|---|---|---|
| Strongly disagree | | | | Strongly Agree |

(b) I felt I was reasonably prepared for the group exam. I did my best and was able to contribute to the group exam.
(circle one)

| 1 | 2 | 3 | 4 | 5 |
|---|---|---|---|---|
| Strongly disagree | | | | Strongly Agree |

(c) I felt everyone in the group wrote their answers after understanding the question and did NOT just copy from someone else in the group.
(circle one)

| 1 | 2 | 3 | 4 | 5 |
|---|---|---|---|---|
| Strongly disagree | | | | Strongly Agree |

allowing the discussions to start immediately. Students were encouraged to sit with their exam groups two to three periods before the exam day after several students expressed that it helped alleviate social anxiety and allowed them time to strategize and practice discussing mathematics before a high-stress event, like the exam. Since class time was generally spent with time for students to discuss and respectfully debate mathematics, this class time with their exam groups was beneficial for students. It also allowed me and my TAs time to monitor the groups and adjust, if necessary. The non-graded, Likert scale questions were implemented to obtain immediate feedback about how the group exams were experienced by students and make changes going forward, if needed. This feedback also allowed me to make more informed grouping assignments for future exams. The decision to provide students with 15 minutes to work independently before embarking on group work was made very early on, based on survey data I gathered after the first exam given in this way. Most students reported that this exam structure worked better for them than traditional, silent, solo exams, but some students who reported leading their group discussions felt that their group members were not prepared for the exam. However, an equal number of students reported that they were prepared for the exam but didn't have enough time to process and think through the problems before group discussions began. As a result, they didn't feel like they were able to contribute as well as they knew they could. Thus, I made sure that students first had 15 minutes to work on the group portion of the exam silently and alone to ensure that everyone could contribute meaningfully to the group discussions.

## 3. *What were the results overall?*

### 3.1. *Student Population Data*

This project on assessment changes was done at a large R1 predominantly white institution (PWI). According to institutional data (at obia.utah.edu) for spring and fall semesters in 2018, domestic students of color made up about 24-28% of Calculus II students, and about 5% of Calculus II students were international students. The gender breakdown for the course sections for this study was about 28-30% women. (Note: There was no reported institutional data on non-binary or gender fluid students available.) I did not collect gender, race, or ethnicity information for the spring of 2018, but in the fall of 2018, out of 216 students who filled out at least one of the surveys, 61 students (28.7%) self-reported as women or gender fluid, and 23 reported as either Black, Indigenous or Latinx students (10.6%). The further breakdown of that group was zero Native American, two Black, 18 Latinx and three Pacific Islander or Native Hawaiian students.

### 3.2. *Grade Data Results*

These exam changes started in the spring of 2018, continuing in the fall of 2018. Additionally, I taught a large section of Calculus II in the fall of 2017 that served as a control group since there were no further structural or pedagogical changes, except for the exam changes, between the fall of 2017 and the spring/fall of 2018. Table 1 shows the aggregate grade point averages and the percentage of students earning D or F grades or the number of students who withdrew from the class (DFW) data for all Calculus II students across all sections taught from various instructors, for those three semesters. This data shows (1) my students from fall of 2017 were scoring about the same as other students with different instructors and (2) final grades increased with these changes to the exam structure starting in 2018, as compared to my previous courses and to the other instructors' courses.

*Table 1:*
*Comparison of Mean Grade Point Average and DFW Rates for Calculus II Students Across All Instructors*

| | Number of students | Mean GPA | DFW rate (as percentage) |
|---|---|---|---|
| **Fall 2017** | | | |
| My students before exam changes | 176 | 2.55 | 29.5% |
| All students except mine | 284 | 2.55 | 29.2% |
| **Spring 2018** | | | |
| My students after exam changes | 174 | 3.14 | 12.6% |
| All students except mine | 447 | 2.64 | 22.8% |
| **Fall 2018** | | | |
| My students after exam changes | 241 | 2.99 | 15.8% |
| All students except mine | 248 | 2.69 | 25% |

*Note:* The Mean Grade Point Average (Mean GPA) was calculated on a 4-point scale.

In fact, grades on exams and the final grades increased for all low and medium scoring groups, while the top scores remained high, with standard deviations becoming smaller than typical in these classes. Additionally, DFW rates significantly decreased.

### 3.3. *Spring 2018 Qualitative Results*

By looking at written survey responses and commentary from focus groups, collected from both 2018 semesters, we can learn more about why this change might have led to more positive outcomes. Students reported higher levels of confidence, lower levels of anxiety, and a stronger feeling of belongingness with this testing structure compared with traditional exams in prior math courses. The following quotes from students during the spring of 2018 are representative of larger themes repeatedly appearing in the survey data.

- Student 1: “It really helps relieve a lot of the anxiety. It also helps with the solo part because you know how the test is structured and you can learn new things working with other students in the class.”
- Student 2: “I really think the group structure helps both those who are doing great and otherwise. Explaining math and understanding it both require cooperation in my opinion.”

3.4. *Fall 2018 Survey Data*

A little more than halfway through the fall semester, my Calculus II students were asked to answer this question: “In your experience, describe how the group portion of the exams has made the classroom/learning environment more or less humane for students, compared to more traditionally structured exams.” Prior to this survey, I intentionally strayed away from using the words humane or rehumanizing in the context of teaching these calculus courses, to try not to influence their answers to this and other survey questions that had a rehumanizing lens. The survey answers were then coded for (A) positive, meaning they felt it made the environment more humane, (B) neutral or mixed, which included answers that were partly positive and partly negative, or (C) negative, meaning they felt it made the class less humane.

Out of the 216 students who answered the question, 84% of students overall answered that it made the class and their learning more humane compared to the traditional silent-solo model of testing. Compare this to a positivity rate, for humane-ness of the exams, of 79% for all women students, whereas 82% of white men answered positively. In each of these student categories, the percent of students who answered this question negatively ranged from 0-3%. The remaining percent of comments were neutral or mixed. Another insightful figure to report is that 90% of white women stated this exam structure was more humane than solo testing.

Considering Black, Latinx and Indigenous students more closely, the positivity rate overall for this group of students was 78%. Excitingly 93% of Black, Latinx and Indigenous men reported the exams as being a more humane experience for them, compared to solo exams. However, one concerning number to report is that Latinx and Indigenous women (there were no Black women students that semester) had a positivity rate of only 56%, with a neutral or mixed rate of 33%. This was admittedly from a small sample, only 9 such women. However, it’s a compelling number and thus is the target of an in-progress qualitative study (Dobie & MacArthur, 2021; MacArthur & Dobie, 2022) exploring why these women had contrasting experiences and what can be done to the assessment structure to indeed create a more humane experience for these students. In total, there is suggestive evidence that this group exam structure is considered more humane by most students, with particular emphasis on white women and Black, Latinx and Indigenous men, and overall, students felt that their learning was enhanced. With that said, there is clearly more work to be done here to increase humane-ness experiences for students who are at the intersection of two historically marginalized groups, namely Black, Latinx and Indigenous women.

## 4. *What were the results for Black, Latinx and Indigenous or women students?*

### 4.1. *Reported Benefits of Group Exams*

Repeated reasoning given by Black, Latinx and Indigenous or women students that supported their claim of the group exams being more humane included (1) learning from varied students' perspectives, (2) collaboration and cooperation, and (3) teamwork that mimics what they expect to see in their future careers. Table 2 highlights themes that were identified from students' feedback regarding their views of the increased humanity of group exams.

*Table 2*
*Illustration of Themes with Quotes*

| |
|---|
| Theme 1: "I enjoy it because you get to see how other people see the same concept or solve the same problem differently and makes the problem easier in a way." |
| Theme 2: "It allows for a more precise work. I believe three heads are better than one. It's also helpful because your partners can catch the silly mistakes you make under the testing pressure. I also believe having the group portion increases students desire to learn the concepts because they don't want to negatively impact another students grades." |
| Themes 1 & 2: "The group portion has made the classroom more humane by allowing the opportunity for students to learn from, help, and collaborate with their peers. I like that it emphasizes learning as a team and that having support from classmates is acceptable. It doesn't make it easy, but makes the hard questions more tolerable." |
| Theme 3: "Teamwork element gives is [sic] a more holistic belonging compared to being alone. In real life collaboration is how the world problems were solved." |

Embedded in these survey responses is an underlying theme of students embracing the group discussion during exams to express their rightful ownership of mathematical knowledge, independent of the instructor, thereby shifting the traditional power balance. Furthermore, the large number of positive responses revealed students felt a sense of belonging with each other both in the classroom and during high-stakes exam scenarios.

### 4.2. *Improvements Needed for Group Exams*

When looking through the comments from students who identify as women that were coded as neutral or mixed, there were some repeated themes that appeared. Out of 11 neutral or negative comments from women, two of them (both white women) mentioned being talked over in their group as the negative side effect of group exams. None of the men mentioned this impact as part of their negative responses. Five out of 11 neutral or negative coded answers from women also made comments in line with the hitchhiker effect, meaning the student was primarily concerned with their group mates not preparing responsibly for the exam, instead copying off the more studious students. These results suggest there is more work to do around developing group norms and paying closer attention to how students are grouped to better serve women.

## 5. *What are the challenges and implications for other institutions/instructors?*

This exam structure succeeded, in part, because students were used to discussing mathematics in class on a regular basis. For example, this class was taught as a flipped course where students were expected to watch a 20-minute lecture video before class and participate in both independent work and work with their elbow partner on calculus problems during class. The class mission statement[33] was referred to daily and enacted by students, TAs and the instructor. Additionally, these active-learning classes were taught by a long-time educator with mathematical content and pedagogical knowledge. It's prudent to recognize that simply changing exam structures, without ensuring a foundation of implementing equitable and inclusive pedagogy during class exists, will likely not have similar success for learning outcomes and belongingness of calculus students (e.g., Burdman et al., 2021; Reinholz et al., 2022; Reinholz & Shah, 2022). One of the main motivations and tenets of this research is specifically to align assessment structures with classroom pedagogies. As such, if a faculty member is teaching with a traditional lecture without consistently practiced discussions by students during class, this exam structure is not likely to engender a sense of belonging or necessarily improve confidence in students. Furthermore, even if an instructor is experienced with active learning techniques in class and attentive to inclusivity of students during class, the nuanced needs of group members and logistics of the group exams still need to be tended to. Lastly, the data highlighted in this paper is suggestive that students need to feel encouraged from the instructor to actively take ownership of their mathematical knowledge which means the instructor needs to be comfortable disrupting normative power dynamics in the classroom on a regular basis for group exams to work well.

As such, these assessment changes are scalable and replicable by educators and/or institutions who are already employing active learning strategies in their classes and have knowledge of inclusive teaching practices. Given such programs or faculty, this collaborative testing structure could be valuable in university calculus classes to improve student success, as measured by grades, and more importantly to improve student enjoyment, confidence, and their sense of belonging, particularly for historically underserved populations of students.

In closing, this is one attempt of many possible positive teaching moves that aims to increase students' sense of belonging in a way that helps them experience their college calculus class as an equitable and inclusive learning environment. Equitable assessment practices are always a work in progress as equity is not a destination, but a process, and this work is representative of that iterative, lifelong action.

[33] Class mission statement: "This is a kind, brave, inclusive, and failure-tolerant class."

# 18 Supporting the STEM Pathway at Appalachian

**Katherine Mawhinney**

*Appalachian State University*

Katherine Mawhinney is a Professor of Mathematics and the Assistant Chair of the Department of Mathematical Sciences at Appalachian State University. Her research interests include topology, undergraduate mathematics education, and K-12 mathematics teacher professional development.

**Katrina Palmer**

*Appalachian State University*

Katrina Palmer is a Professor of Mathematics in the Department of Mathematical Sciences at Appalachian State University. She is particularly interested in creating ways to actively engage students with mathematics both online and face-to-face. Her mathematics interests lie in numerical linear algebra applied to astronomical and medical imaging.

**Gregory Rhoads**

*Appalachian State University*

Gregory Rhoads is a Professor of Mathematics in the Department of Mathematical Sciences at Appalachian State University and teaches a wide range of undergraduate and graduate courses. His interests include applying complex function theory to minimal surfaces and investigating methods for improving the teaching of mathematics.

*Abstract: In this paper we describe an initiative aimed at providing additional support to students in calculus to improve progress through the course and keep students on the correct pathway to graduate on time. Most STEM disciplines require calculus but some students enter the course with inadequate preparation, leading to an unsuccessful attempt at completing the course. This puts the student behind in their program and reduces their confidence with mathematics. To address these issues, we have revamped our placement strategies, created a new support course for those students currently taking Calculus 1, and created a new fallback course to help those students who must drop Calculus 1 to prepare themselves for a successful retake the following semester. We will describe our efforts and provide suggestions from our experiences for those interested in implementing similar changes on their campuses.*

## 1. *Introduction*

The calculus-sequence is *the* mathematics pathway into STEM fields at most colleges and universities. Appalachian State University is committed to inclusive excellence: "... a process and outcome that promotes critical thinking about knowledge, seeks to understand and address historical and present-day inequities, and creates equitable classroom environments."[34] Rising to meet this challenge, the Department of Mathematical Sciences has been implementing a comprehensive effort to promote equitable access to and support through calculus.

[34] Appalachian State University - Center for Academic Excellence: https://cae.appstate.edu/inclusive-excellence

We started this initiative to support STEM students[35] by comparing success rates for our target populations: transfers, underrepresented, first generation, low income, and rural with success rates for all students.[36] We collected and compared student DFW rates (D, F, or Withdraw) in Precalculus and Calculus 1 for the general student population and for these target populations. These data are shown in Table 1.

As shown in the table, these target populations all have higher DFW rates than the entire population of potential STEM students. Our initial goal for this initiative was to reduce the DFW rates overall in Calculus 1 and to bring the DFW rates for these target populations more in line with the entire population.

*Table 1*
*DWF rates for STEM gateway courses from 2016-2019. Definitions of these groups can be found in the footnote at the bottom of the page.*

| | **Precalculus DFW** | | **Calculus 1 DFW** | | **Calculus 2 DFW** | |
|---|---|---|---|---|---|---|
| All Students | 27.7% | Difference | 32.0% | Difference | 33.3% | Difference |
| Transfers | 39.4% | 11.7% | 42.7% | 10.7% | 39.8% | 6.5% |
| Underrepresented | 33.2% | 5.5% | 34.4% | 2.4% | 36.3% | 3.0% |
| First Generation | 31.1% | 3.4% | 37.7% | 5.7% | 34.6% | 1.3% |
| Low Income | 33.4% | 5.7% | 35.8% | 3.8% | 37.5% | 4.2% |
| Rural | 29.9% | 2.2% | 37.6% | 5.6% | 37.7% | 4.4% |

## 2. *Focus on Inclusivity*

While lowering DFW rates is an easily measured goal, there are many different factors that affect student performance. The potential prior experiences of the students in our target populations will influence their interaction with their peers, instructors, and administrators, as well as mold their learning habits. These influences and their effects need to be acknowledged as well as addressed when considering inclusivity and equity for all of our students. Some challenges were anticipated at the beginning, and some were not identified until later in our process.

We decided to focus our efforts in two directions: placement and support. We wanted to ensure all of our students taking Calculus 1 had an academic preparation that will give them the best chance of success. Once they were in our Calculus 1 course, we wanted to create a support system that would address issues that challenge some of the students in our target populations.

---

[35] We focused on STEM students because we found the overall DFW rate of our STEM Pathway courses (Precalculus, Calculus 1 and Calculus 2) were significantly higher than the non-STEM 1000-level courses.
[36] Transfers are students who first enrolled at Appalachian as a transfer student.
Underrepresented describes any Appalachian student who does not self-identify as only white.
First Generation is self-reported.
Low Income is based on Pell Grant eligibility.
Rural students include tier 1 and 2 counties, as defined by the North Carolina Department of Commerce.

These program supports are designed to provide a more equitable mathematics experience, thereby allowing all students to proceed through their required mathematics pathway in a timely manner.

The following sections describe the program, evidence of support for inclusivity, and obstacles and recommendations for adoption by other institutions.

## 3. *Program Description*

Our initiative to support STEM students through calculus involves 3 primary steps: modifying the calculus placement process, creating calculus corequisite courses, and developing a fallback course.

### 3.1. *Placement*

The traditional view of placement procedures and testing is that these tools help ensure that students are "calculus-ready". We see these tools as helping to ensure that students are "calculus-eligible", recognizing that there is more to readiness than prerequisite skill and content knowledge. As Franco and Patel summarize, content knowledge is just one of the dimensions of college readiness indicating that students will be successful in college-level courses (Franco, 2013). Research has shown that academic behaviors and college knowledge are also important characteristics of successful students, so our supports aim to attend to these non-content needs as well.

Informally, faculty had been reporting the absence of strong algebraic skills in their calculus students, which led us to question the accuracy of our placement of students into Precalculus versus Calculus 1. Thus, in Fall 2017, we analyzed data from our department-developed Calculus Readiness Test and the subsequent student performance in Calculus 1. After evaluating this data, we saw that students with lower passing scores on the placement test struggled with Calculus 1 and had a higher DFW rate. We decided to require a higher score in order to place a student in Calculus 1. Among those students who were not successful on the Readiness Test, we noticed the differences in their precalculus preparation, some needing more and some needing less work than others. So, we simultaneously provided an alternate entry into Calculus 1 without the time cost of an entire semester for a precalculus course. If a student didn’t pass the Calculus Readiness Test, they could earn entrance into calculus by completing 85% of the online self-paced *Prep for Calculus* ALEKS[37] (Assessment and Learning in Knowledge Spaces, McGraw-Hill) course. This alternate entry allows each student to work on those topics that they specifically need, and gets them into Calculus 1 more quickly.

### 3.2. *Calculus 1 Corequisite*

We then created and implemented a 1-hr corequisite support course for Calculus 1. This support course was developed and first offered during Spring 2018, designed for any student who needed prerequisite content review or extra insight in the calculus content or related academic skills (for example, study skills, assessment preparation, and personal organization).

[37] www.aleks.com

Students could self-select into the course or be recommended by an advisor based on test scores, prior coursework, or the student's own apprehension towards calculus. Allowing students to self-select into the course resulted in an audience with a wide variety of content preparation and calculus readiness.

This course emphasizes active learning and reflection, attention to student self-regulation, academic skill development, and community building. Not only do students engage in tasks that increase their understanding of concepts, but they are actively doing mathematics in a safe space with classmates and with faculty who believe that they can be successful in mathematics. A typical weekly class meeting may include things like a quick review of current Calculus 1 concepts focusing on common misconceptions, active practice of Calculus 1 problems in small groups, reflection on study or test preparation skills, and/or a review of algebraic concepts required by the next Calculus 1 topic. Figure 1 shows two examples of sense-making activities that are used in the course.

**True or False? How do you know?**

$$(x^2 - y^2) = (x - y)(x - y)$$

$$\sqrt{x^2 - 1} = x - 1$$

$$\frac{\sqrt{x+h} - \sqrt{x}}{h} = \frac{\sqrt{x+h} - \sqrt{x}}{h} \cdot \frac{\sqrt{x+h} + \sqrt{x}}{\sqrt{x+h} + \sqrt{x}}$$

(a) True or False discussion items

| | |
|---|---|
| $(f \circ g)(x)$ | $f(x)g(x)$ |
| $f(x+h)$ | $f \circ g$ |

(b) Which One Doesn't Belong discussion item

*Figure 1: Two active-learning examples aimed at helping students make sense of mathematics and mathematical notation.*

Instead of traditional lecturing about definitions and properties, we often use examples like the ones in Figure 1 to have students question and explain their understanding of basic prerequisite topics, such as factoring, simplification, and function composition. Since these students who self-identify as not confident in calculus often do not want to offer answers to the whole class, students first consider the examples individually and then in small groups before we facilitate a discussion to uncover and resolve potential misconceptions about the topics.

In addition to these sense-making activities, we include several self-regulation tasks throughout the course. Self-regulated learning occurs when students monitor and regulate their own cognition and behaviors. These types of activities have students reflect on their awareness of support resources, engagement in the course, or their level of understanding of the material, all in a manner that is supportive and respectful. One example would be to have students rate themselves on prompts such as "I attend every class", "I ask questions during class", "I read the course textbook", and "I read my course notes", by holding up a certain number of fingers. The students are reminded of these aspects of learning over which they have control, and the instructor can use the results to gauge student involvement. Another example is a written assignment in which students describe their preparation for an assessment and how they felt once they completed the assessment. Students then reflect on the strengths and weaknesses of their preparation along with how their preparation will change for the next assessment.

Due to the positive responses from the students in the Calculus 1 corequisite, a 1-hour corequisite for Calculus 2 was introduced in Spring 2020, and a 1-hour Precalculus corequisite will start Fall 2020. All three of the corequisite courses will balance the active learning of course content with community building, self-regulation and development of general academic skills in the self-selected students who enroll in these courses.

### 3.3. *Fallback Course*

Many students who qualify for and take Calculus 1 realize early in the semester that they will not be successful in completing the course. Possible scenarios include a student with a low grade in Precalculus, a student transferring a precalculus course that had large content gaps, or someone with precalculus credit from many years ago.

Struggling Calculus 1 students had the following options:

1. stay in the course and try to earn a passing grade,
2. drop the course, with a penalty of using one of their 4 career drops[38], and enroll in the course again next semester, typically no better prepared
3. drop the course, with penalty, and enroll in Precalculus the next semester, which would not count for a student with credit for college-level Precalculus and would put them behind in their academic progress.

The fallback course was created for these students to provide them a better alternative while strengthening their preparation for calculus. Students can drop Calculus 1 and enroll in the fallback class without penalty (without using a career drop) anytime in the first 5 weeks of the 15-week semester. We chose to begin the fallback course five weeks into the semester in order to give students enough graded work in Calculus 1 for them to assess their chance of successfully completing the course, but also to allow 10 weeks for the students to complete the required material in the fallback course. We chose to have the students earn 3 credit hours for successfully completing the fallback course to offset the loss of 4 hours from dropping Calculus and 3 hours credit is an accurate reflection of the effort required.

The fallback course has 3 contact hours per week and focuses on algebraic techniques and misconceptions in an active and supportive environment. The class emphasizes active learning experiences that build conceptual understanding and support procedural fluency. Outside of class, students complete at least 85% of the online self-paced *Prep for Calculus* ALEKS course to document their progress and mastery. With successful completion of the fallback course, students are better able to enter their correct math pathway without penalty and are well prepared to succeed in Calculus 1 the next semester.

---

38 Courses dropped after the drop/add period but before the end of the ninth week of a fall or spring semester are referred to as 'career drops'. Students are limited to a total of four "career drops" during their undergraduate careers at Appalachian.

## 4. *Results*

Table 2 shows the DFW rates for Calculus 1 for all students during the first 4 semesters the corequisite support course was offered. The DFW rates are further broken down by those who enrolled in the support course and those who didn't. For example, 44.7% of the transfer students enrolled in Calculus 1 ended the course with a D, F or W, while 47.8% of transfers enrolled in both Calculus 1 and the support course ended Calculus 1 with a D, F or W. At first glance this data seems to suggest that the support course reduces the overall DFW rate but not for all target populations. However, the actual number of students in these targeted groups is low. For instance, there were a total of 17 underrepresented students enrolled in the corequisite course over the 4 semesters and of those 17, 7 earned a D, F, or W in Calculus 1 (for the 41.2% in the table). Also, since the course is intended to attract students who may not have performed well in the past, it would be better to compare the students in the corequisite course with students of a similar academic background (similar SAT, similar grade in Precalculus, etc.) who are not in the corequisite course.

*Table 2*
*Calculus 1 DFW Rates comparing those with the support course and those without.*

| | **Calculus 1 DFW Rates** | | |
|---|---|---|---|
| | All Calculus Students | With Support Course | Without Support Course |
| All Students | 33.6% | 29.9% | 33.8% |
| Transfers | 44.7% | 47.8% | 44.5% |
| Underrepresented | 33.3% | 41.2% | 32.8% |
| First Generation | 38.9% | 34.4% | 39.2% |
| Low Income | 36.7% | 30.8% | 37.0% |
| Rural | 37.0% | 38.9% | 36.9% |

While changes in the DFW rates were mixed, individual students commented on benefits of the support course beyond grade improvement. For example, in response to "What aspects of the course helped you learn?" students said, "... presenting questions in different ways, to make sure everyone could understand," and "Aspects such as more visual learning and deeper explanations." When asked about the big "takeaways" from the course one student responded, "... And to even go out of your comfort zone and ask some questions even if it might be awkward because it might feel like you're the only person asking questions. Like no, just go ahead and ask," and another said, "Probably just an understanding of how to access the tools like Math Lab, office hours, and a general do good at math system, or a secondary source of notes to look back on." A third student response to the prompt included, "Trying new ways of communication, trying new things, and always building that relationship with your teacher is helpful."

## 5. *Potential Obstacles and Recommendations*

The biggest obstacle we've found is the difficulty in getting the students who would most benefit from these support courses to register. Some institutions with support courses require students to enroll based on low test/placement scores, but we decided against the negatively-toned message this might send to some students, especially in our target populations. We chose to rely on general advisors and an advertising email to all students enrolled in Calculus 1 that touts the benefits of the support course. This approach empowers the student by giving them the choice, but some students may not have good indicators of their preparedness. Students who transfer in a college-level precalculus course currently do not take the placement test as they are deemed automatically prepared for Calculus 1. However, this preparation varies greatly depending on the rigor of their precalculus course and the time since taking precalculus. The creation of the fallback course as an option for students with precalculus credit but don't feel prepared for Calculus 1 may provide some incentive for students to voluntarily take the placement test to assess their preparation for Calculus 1. Also, email may not be the best way to reach some students. Is it reasonable for us to assume that first generation students would understand the benefits of a support course that doesn't count towards their degree program?

We also found it important for advisors to be fully cognizant of the goals of the courses so they can ask the appropriate screening questions to students and provide informed recommendations. In fact, implementing this effort has opened communication between department faculty and university advisors, which has yielded benefits beyond increasing enrollment in these support courses. In addition to advisors, calculus instructors also need to be aware of the goals and logistics of these courses so they can provide accurate information to students during the first week for the support courses and as the semester progresses for the fallback course. This also requires calculus instructors keep track of fallback course registration deadlines as well as communicate students' grades early in the semester.

Another potential obstacle is the non-standard structure of the support courses and communicating these differences to both faculty and administration. Faculty who wish to teach one of these courses need to be capable of and interested in providing an education experience that is not just "more of the same," engaging students in an active learning environment, attending to student self-regulation, offering formative feedback, and encouraging students to seek help during office hours or other times outside of class. Administration should be informed of the time commitments required for these courses by faculty and support staff and of the difficulties of measuring their impact, especially when making decisions about course assessment and future course offerings.

From an administrative standpoint there are two other aspects that need to be addressed. In scheduling, both support and fallback courses should be offered in a time period that would minimize potential student scheduling conflicts. Secondly, exceptions to university policies may be possible to help the students. At Appalachian, students were able to drop their Calculus 1 course and add the fallback course in the first 5 weeks of the semester without having to use a career drop. We also faced a university policy that if a student is enrolled in both Calculus 1 and the corequisite course and Calculus 1 is dropped, the corequisite is automatically dropped. We were allowed to petition for corequisite students to remain in the course if they dropped Calculus 1, but not register for the fallback course, so that only 1 career drop would be used.

## 6. *Conclusions*

Ultimately, supporting students in the STEM pathway is a collaborative effort. In order for it to be done well, there needs to be faculty, advisors, and administrators who are committed to helping students. Succeeding in calculus is essential for retaining STEM students and having them complete their degree.  Not all students enter the university with sufficient content preparation or have the academic skills needed to succeed in calculus, and special populations are especially vulnerable. Our initiative is our attempt to provide all Appalachian STEM students with the best possibility of completing Calculus 1 in a timely manner. This begins with students starting Calculus 1 with a sufficient content background and proper academic skills. By modifying the placement requirement and including the ALEKS option, we have given students a better idea of and a timely option for completing the prerequisite background. Our support course provides assistance on content, academic skills and self-regulation skills, again to level the field. The fallback course replaces a negative experience with a positive one that will build self-confidence which is important to our target groups. We will continue to improve our courses and hope to see greater student success. While we are still looking at more comprehensive measures of course assessment, student evaluation comments suggest these courses have had a positive impact so far. We hope that some of our techniques and experiences will be helpful for other institutions wishing to implement similar ideas.

# 19 “I Wish I Could Say ‘You Should Not be Here’”: An Analysis of Instructors and Students Contrasting Perceptions of a Racialized and Gendered Gatekeeping Practice in Calculus

**R. Taylor McNeill**

*Vanderbilt University - Peabody College of Education & Human Development*

Taylor McNeill <reagin.t.mcneill@vanderbilt.edu> is a doctoral student in Mathematics Education at Vanderbilt University’s Peabody College of Education. They previously received a Ph.D. in mathematics and worked as a mathematics professor, experiences which contributed to their interest in equity within higher education mathematics. Their research now focuses on the enactment of equitable mathematics instruction and minoritized faculty socialization in mathematics.

**Brittany L. Marshall**

*Rutgers University*

Brittany L. Marshall is a doctoral student in Mathematics Education at Rutgers University’s Graduate School of Education. Previously, she worked as a middle and high school mathematics teacher for almost a decade. Her research focuses on positive mathematics identity in Black girls and the teachers who help cultivate it.

**Luis A. Leyva**

*Vanderbilt University - Peabody College of Education & Human Development*

Luis Leyva, Ph.D. is an assistant professor in the Peabody College of Education & Human Development at Vanderbilt University. His research explores how interlocking systems of racism, white supremacy, and cis-heteropatriarchy shape historically marginalized students’ experiences of undergraduate STEM instructional and support practices. He has solo-authored two articles in the *Journal for Research in Mathematics Education* and has also published in *Cognition & Instruction*, *Harvard Educational Review*, and *Journal of Urban Mathematics Education*.

*Abstract: This case study discusses how conceptions of mathematics instruction as race- and gender-neutral impede inclusive teaching in undergraduate pre-calculus and calculus. Drawing on a larger study, we present a case study that exhibits the well-documented gatekeeping function of calculus. This event, titled course drop, features an instructor’s class-wide message that students who cannot complete steps for a problem quickly should consider dropping a course level or not taking Calculus 2. Instructors’ and students’ contrasting perceptions of this event reveal how instructors’ views of pre-calculus and calculus teaching as socially neutral can erase underrepresented students’ racialized and/or gendered experiences of instruction. In particular, we illustrate that while instructors framed the event as motivating to students as well as helping them in finding a supportive course placement, Black and Latin* students interpreted the event as increasing racialized pressure to succeed and causing them to question their belonging in STEM. These insights are leveraged to advocate for mathematics faculty professional development that promote race- and gender-conscious teaching practices rooted in critical self-awareness to foster more equitable learning opportunities in entry-level mathematics.*

## 1. *Introduction*

This case study discusses how conceptions of mathematics instruction as race- and gender-neutral impede inclusive teaching in undergraduate pre-calculus and calculus. We draw on data analysis from a larger study that explores instructors' and students' perceptions of instructional events, which Black and Latin*[39] students, as well as white women identified as discouraging. This study involved interviews with instructors and students that centered around written descriptions of instructional events. These events were developed to reflect themes in historically marginalized students' journaling of positive and negative moments experienced in calculus classrooms (see Battey et al., in press and Leyva et al., in press for more study details). Our case study focuses on data analysis for an event that exhibits the well-documented gatekeeping function of calculus (Ellis et al., 2016; Jett, 2013; Leyva, 2016). This event, titled *course drop*, features a professor's class-wide message that students should consider dropping a course level or not taking Calculus 2 if they could not complete steps for a problem quickly. Instructors' and students' contrasting perceptions of this event reveal how instructors' views of pre-calculus and calculus teaching as socially neutral can erase underrepresented students' racialized and/or gendered experiences of instruction.

Since mathematics is seen as free of culture (Battey & Leyva, 2016; Bishop, 1990), instructors may invoke reasoning based on logics of color-blindness and gender neutrality to shape their instruction. C*olor-blindness* is resistance to recognizing race as relevant in shaping social inequality (Bonilla-Silva, 2003). The counterpart to color-blindness, *gender-blindness* or *gender neutrality* (Acker, 1990), appeals to gender as irrelevant to social inequities, contributing to erasure of marginalization among women and gender nonconforming individuals in mathematics (Esmonde, 2011; Leyva, 2017). Colorblind, genderblind perspectives allow instructors and institutions to ignore inequities in mathematics instruction that disproportionately impact historically marginalized students, including Black students, Latin* students, and white women (McNeill et al., 2020).

The study's data sources include student journaling of instructional events as well as stimulus-response interviews with students and instructors around the same set of selected events. Participants included: (i) seven pre-calculus and calculus instructors at a large, public research university, including white and Asian men as well as Asian, Black, and Latin* women; and (ii) 20 students from underrepresented backgrounds across intersections of race and gender enrolled in precalculus and calculus. During interviews, participants were asked to read the following written description of the *course drop* instructional event, which was based on student journaling of events from actual pre-calculus and calculus classrooms.

> *During class, the professor said something along the lines of, 'If you do not know how to do these steps quickly you might want to consider dropping down to a lower class or consider not taking Calculus 2'. As a student who wants to major in a STEM field, it*

[39] Latin* is a term that encompasses fluidity of social identities. The asterisk considers variation in self-identification among people of the Latin American diaspora and origin (Salinas, 2020). Latin* responds to (mis)use of Latinx, a term reserved for gender-nonconforming peoples of Latin American origin and descent (Salinas & Lozano, 2019).

> *made me feel a bit uncomfortable because I did take a bit longer to do the steps. At the time of the event, I thought to myself if I was good enough to go on to Calculus 2.*

Participants were asked how they interpreted the event, including whether they had experienced something similar in the classes they taught or attended. Interviewers also asked participants whether they believed students would be discouraged by the event, probing their responses to see which student populations (if any) would experience it as discouraging and whether this would vary based on the race or gender of the instructor and/or student(s) involved.

**2. *Instructors' colorblind and genderblind event interpretations.***

Instructors largely perceived the *course drop* event as void of race or gender. Bhavik (Asian man) claimed any student, regardless of gender, would experience this event similarly due to shared learning goals among all mathematics students. He explained, "The gender should not affect this because as long as the person is a learner, the goal is to get through the course and learn the concepts. So yeah, for me, gender should not affect the interpretation." Joshua (white man) also perceived the event as race- and gender- neutral, as students universally struggle with calculus concepts.

> *No, this is literally a universal thing across the board that students in every... Even up to Calc 3... They still struggle with some concepts from Calc 1 and 2. This is universal, across the board, every student, every instructor.*

In this way, instructors perceived the event, and learning calculus in general, as a struggle for most students, regardless of their race and/or gender identities. Many students, however, interpreted the event as racialized and gendered. For example, seven students appealed to the underrepresentation of women, Black, and Latin* students to explain how this event may be experienced differently among students with different identities (Leyva et al., 2021).

**3. *Problematizing student motivation.***

Instructors also supported the professor's choice to deliver the class-wide *course drop* message, in many cases explaining that they perceived this message as motivating students to study and thus improve their ability to succeed in the course. Robert (white man) discussed how the point of this message is to make students realize "they have to take action."

> *I think it's [the instructor's message] intended to make them feel uncomfortable… You want to make the students realize that they have to do something, that they have to take action. They can't just ignore this problem. I think it's appropriate to make them feel uncomfortable… It's the whole purpose of saying this. 'Don't get too comfortable. You're not doing as well as you think you're doing… You can't just wish your way out of this.'*

Robert believes the discomfort that students feel from the *course drop* message will motivate them to work harder. In this way, Robert justifies the *course drop* message as being supportive of students.

While some historically marginalized students often also perceived the *course drop* message as motivating, they did not necessarily perceive such motivation resulting in improved course performance as Robert did. Many students perceived the motivational impact of the event increasing their level of stress when considering racialized-gendered rates of representation and achievement in mathematics, which can have deleterious effects (Battey et al., in press; Leyva et al., in press). To illustrate, Nadine (Black woman) described how this message would likely make it more challenging for her to learn, as the pressure to succeed as one of the only Black students in pre-calculus and calculus classrooms would make it difficult for her to focus on learning content.

> *I feel like minorities already... It's already difficult to get through certain classes because you don't have the support system. It's maybe just you and a couple of other kids. You're a little bit more isolated... You have a lot more pressure to succeed because it's not typical for... other people of your race to succeed. So, you have a lot more pressure on yourself to do well because... for me, it's like… people said I couldn't and therefore I'm going to and I'm going to do it well. When things like this happen, I'm like, 'Am I on the level the professor said? If I'm not, what am I doing wrong?' And again you're adding stress to your class and you're making it a lot more difficult for the student to focus on learning.*

Nadine perceived that she might feel more pressure to do well in the course to defy racialized expectations that Black students will be unsuccessful in mathematics. In this way, increased stress may not lead to productive motivation for underrepresented students of color like Nadine who face challenges of persisting in mathematics without a support network of same-race peers.

While Nadine and Robert similarly perceived the *course drop* event as increasing stress and pressure among students, Robert's perception of the event as race- and gender-neutral impairs his ability to recognize racial and gendered barriers in place for historically marginalized students. Nadine, for example, described how the event would result in her managing such stress on top of her academics and feelings of underrepresentation as a Black student. In this way, Robert, like many mathematics instructors, may inadvertently participate in racialized gatekeeping through their seemingly neutral instructional practices that they perceive as benefitting or supporting students.

### 4. *Problematizing student belonging in STEM.*

Instructors also supported the professor's choice to deliver the *course drop* message because it prompted students to switch classes or majors if they were not likely to be successful in the course or more advanced mathematics courses. Instructors largely viewed this as helpful as it supported students to find a course where they could be more successful or saved them from receiving a poor or failing grade. Stan (white man) discussed how students often persisted despite hearing messages like *course drop* and encountering "hard problems" designed to "wake up the students so that they study and work," though they were unlikely to succeed in the course. However, Stan spoke admirably of students who "face up" to realizing that STEM (science, technology, engineering, and mathematics) is not for them.

> *I had a student when I was teaching one time. The student came to me and said, 'I'm failing this course, I know I'm failing, I'm dropping this course, I'm getting out of engineering, this is not for me. And I want to thank you for helping realize that this is* *not the career I should be on.' This does happen occasionally where people face up to what the issue is... You know not everybody has to be an engineer, and so I told him it was okay.*

Similarly, Ester (Latina woman) supported the professor's *course drop* message because it helped students who she perceived to not belong in the class determine that they should drop the course.

*I have done this exactly… It's a good thing you're doing for this student. [If] the student is really having trouble, sometimes I wish I could say, 'You should not be here.' But I feel like I can't say that… [But also] I can't say, 'Oh it's going to be okay, you're going to pass.' It's not true.*

Stan and Ester's interpretations of the event suggest that certain students do not belong in calculus courses or STEM fields broadly. This perspective shapes their views of student attrition from mathematics courses as positive and supporting students to find the 'right' place for them.

Students acknowledged race and gender as factors affecting persistence following the *course drop* event. In many cases, racialized and gendered stereotypes of mathematical ability were raised to explain how students from historically marginalized populations may have their confidence in mathematics more readily eroded. Parker (Black man) interpreted the *course drop* message as disparaging toward students. He further explained that the *course drop* event would make him question whether he would be successful in the course, given the negative stereotypes around Black students' mathematical ability.

*Because for instance… there are these stereotypes that Black people aren't good in math. In this case, if it was especially towards a Black student, that would have discouraged that person really bad, because it's like, there is stereotype [sic] said about my culture or my race, and if a teacher is saying that stuff, that means that I'm just not good enough, I'm really not good enough to do this. I can't do it.*

Parker describes how racial stereotyping of mathematical ability would make him particularly vulnerable to feeling discouraged if he experienced the class-wide *course drop* event. Although Parker continued with mathematics, prior research has shown that racialized and gendered stereotypes of mathematics ability can contribute to student attrition from STEM (Beasley & Fischer, 2012).

Racialized and gendered rates of representation in STEM professional fields were also raised by students as affecting persistence following the *course drop* event. Victoria (Latina woman) discussed how gendered representation in computer science might cause her to think that she does not belong in STEM following the *course drop* event.

*If you're a female and you're not getting it, and then your professor, let's say he's a male, and he's telling you that you maybe shouldn't do this if you're not getting it, and you're like, 'Well, maybe that's why I don't see many female computer science majors. Maybe that's why I don't see many female surgeons that are known worldwide'... Maybe math isn't for me. Maybe STEM isn't for me. Maybe that's why they're all in public relations or something.*

The underrepresentation of women in STEM fields coupled with the *course drop* remark and her struggles with mathematical content signaled to Victoria that "math isn't for [her]." Therefore, situations like the *course drop* event can perpetuate exclusionary ideas of who belongs and succeeds in mathematical spaces, causing women like Victoria to reconsider pursuing STEM fields where they are underrepresented.

Stan's assessment that "not everybody has to be an engineer" brings him to see students dropping required mathematics courses or an engineering major as natural. Similarly, Ester perceives that students struggling with content should not be in her course. These orientations to ability and belonging in mathematical spaces frame student attrition as a positive consequence of instructional experiences. However, these orientations fall short in situating instructional remarks like in the *course drop* event in broader racialized and gendered realities of historically marginalized students. Parker and Victoria illustrate how an instructor's suggestion to drop a course collides with racialized and gendered representation in STEM as well as racialized stereotypes of mathematics ability to amplify their doubt about belonging in STEM. In this way, Stan and Ester's colorblind and genderblind framing of student attrition shapes seemingly neutral instructional practices like the *course drop* event that reinforce racialized and gendered inequities in the calculus sequence and STEM trajectories.

**5. *Implications for practice.***

Our case-study reveals that while faculty often hold perceptions of mathematics instruction as socially-neutral, systemic racism and sexism collide with seemingly neutral instructional practices to shape racialized and gendered experiences for Black students, Latin* students, and white women. In particular, students responding to the *course drop* event identified affective dimensions of learning, specifically stress, isolation, and belonging, as central to their interpretations of the event as racialized and gendered. Faculty invested in creating equitable learning opportunities, thus, can facilitate access to socially-affirming peer and instructional support among students from historically marginalized groups. Creating structured opportunities for students to collaborate and communicate both inside and outside of class (e.g., study groups across class sections) may support students' ability to form peer networks that strengthen performance and persistence. In addition, faculty should consider ways to reframe their support of students, both within and outside classroom teaching, to ally themselves with student wellbeing broadly, rather than limiting support to assistance with mathematical concepts and procedures.

This case study demonstrates that race and gender consciousness is a necessary, though insufficient, ingredient to equitable mathematics instruction. This is consistent with prior scholarship, which has shown that racial consciousness is necessary to disrupt norms and practices in STEM which perpetuate racial oppression (Haynes & Patton, 2019). However, relying on individual faculty to develop critical consciousness and reform their own teaching practice may be inadequate since it has been found that mathematics faculty, as a group, are poor advocates for diversity (Park & Denson, 2009) and often lack professional training and incentives to cultivate equitable teaching practices (Haynes & Patton, 2019; Walczyk et al., 2005). Thus, rather than placing responsibility on individual faculty to make mathematics instruction more equitable, departmental and institutional changes are needed to support faculty growth and prioritize inclusive teaching. For example, improved access to equity-oriented professional development may assist faculty in understanding racial and gendered impacts of their teaching practices, as well as direct them towards equitable pedagogies. This may encourage faculty to move beyond access-based conceptions of equity to instead center the impact of instruction on stress and sense of belonging among students who are minoritized by race and gender in the mathematics classroom. Departments may also consider redesigning evaluations of teaching to specifically invite student reflections regarding racialized and/or gendered experiences of classroom instruction. Prioritizing these and other metrics for equitable teaching in promotion and tenure criteria incentivize instructors to work at incorporating inclusive teaching practices while providing accountability for reducing racialized-gendered gatekeeping in introductory mathematics. In advocating for change, we encourage faculty to step outside their roles as instructors to consider how they, as colleagues, department chairs, and committee members, can help to institutionalize equity-oriented support for fellow members of the mathematics faculty.

# 20 Enhancing Academic Capital to Improve Equity and Inclusion in First Semester Calculus

**Tabitha T. Y. Mingus**

*Western Michigan University*

Tabitha Mingus is an associate professor of collegiate mathematics education at Western Michigan University. Her research interests are on student understanding of and reasoning in mathematics at the collegiate level and how pre-service mathematics teachers' experiences in learning mathematics impacts their teaching. She enjoys playing the violin, gardening, cooking, playing with her dogs, and hanging out with her ducks.

**Melinda E. Koelling**

*Western Michigan University*

Melinda Koelling is an associate professor of mathematics at Western Michigan University. Her research focuses on differential equations and computational neuroscience. In addition to euphonium playing, knitting, and ceramics, her interests include how mathematics is learned by her children and her students.

**Daniela Hernandez**

*Western Michigan University*

Daniela Hernandez is a faculty specialist in mathematics at Western Michigan University. She has taught and coordinated a variety of courses including developmental mathematics, precalculus, and calculus courses at WMU for over 20 years. In her spare time, she likes to travel, read, solve puzzles, garden, and experiences great joy from watching hummingbirds in her own yard.

*At Western Michigan University, the calculus sequence had unacceptable progression rates. The impact of these rates was not borne equally among all categories of students. Students of some races (e.g., black, multiracial) and family backgrounds (e.g., Pell-eligible, first-generation) either retook calculus and progressed slowly or switched from a STEM major. As a result of discussions in an interdisciplinary committee, calculus classes were coordinated to improve student outcomes and make the learning environment more equitable. In this paper, we describe our view of the problem at our institution, interventions taken, and results obtained. We argue that we improved student success while maintaining course integrity and standards. However, more work is needed to improve the progression of these groups.*

## 1. *Introduction*

University calculus provides gateway access to numerous STEM fields. Instructors have an obligation to foster the growth of all students to increase access to and inclusion in STEM careers. In the decade 2005-2015, the number of students failing calculus at Western Michigan University (WMU) rose to unacceptable levels. This attrition rate had several potential causes: a natural consequence of students' discovery of their true interests, student pragmatism about their commitment to a career path, and the educational environment in their courses.

We believed there was an opportunity to affect student attitudes and improve student success by improving the educational environment. This paper describes our work at our institution in Calculus 1 to increase student success for all groups. While all groups were challenged, certain populations (e.g., Pell Eligible, First Generation, African American, and Hispanic) faced greater challenges. Techniques implemented were chosen to help those that struggled most.

## 2. *Institutional and Program Diversity*

### 2.1. *Diversity at Western Michigan University*

Situated in southwest Michigan, WMU draws half of its students from the fifteen counties within a 50-mile radius of its main campus. Fall 2019 enrollment was 17,051 undergraduate and 4,419 graduate students. In Michigan, there are nearly three times more white-single race residents (7.48M people) than residents of color.[40] At WMU, there were two times more white-single race students (14,441) than students of color (7,029). The student body is more diverse than the surrounding community.

### 2.2. *Diversity in Calculus 1*

This case study focuses on an intervention in Calculus 1, a first course in a four-course sequence.[41] The typical class consists of students from various backgrounds. Racially underrepresented, women, students with financial need,[42] and those who are first generation tend to experience more challenges in the academic environment (Hurtado, Newman, Tran, & Chang, 2010). Consequently, we chose to examine the progress of these groups. In the Calculus 1 course, the two largest race/ethnic groups (after white-single race) are black or African American alone or Hispanic, each representing 7% of the course enrollments[43]. Nearly one-quarter of the students are female. For the 2018-2019 academic year, 21.7% of students were Pell eligible and 10.9% were classified as first generation.

## 3. *Student Progression in Calculus 1*

During 2005-2015, the failure to progress rate in Calculus 1 increased. Since DEWI[44] rates were steadily climbing while withdraw rates were declining, the increase in DEWI rates resulted from an increase in DC, D, and E grades. The disaggregated DEWI rates showed a high variability across sections.

The 2014-2015[45] DEWI rates were disaggregated by student characteristics. Gender did not appear to be a factor: 48.1% rate for females and 48.7% for males. Race/ethnicity did impact the rates (Figure 1). Compared to the overall rate: International, Asian, and Hispanic students' rates were lower; white students' were similar; and black and multiracial students' were higher.

---

[40] The second and third most common racial/ethnic groups are Black or African American Alone (1.36 M) and Hispanic or Latino (517 k). (Demographic data for the State of Michigan in 2018 taken from the website https://datausa.io/profile/geo/michigan#demographics).

[41] We refer to the courses as Calculus 1 (MATH1220), Calculus 2 (MATH1230), Calculus 3 (MATH2720), and Calculus 4 (MATH3740).

[42] Indicated by Pell-eligibility.

[43] International students were not identified by their racial/ethnic group.

[44] WMU assigns A, BA, B, CB, C, DC, D, E, X, NC, and I as grades. Students who receive a DC, D, E, X, NC, or I are included in the DEWI rate.

[45] The academic year prior to redesigning the course.

The race/ethnicity groups that merited special attention because of lack of representation in STEM are black, multiracial, and Hispanic students[46]. In addition to underrepresented race/ethnicity in STEM[47] fields, first generation status and Pell eligibility were associated with higher DEWI rates and also merit attention.

### 4. ***Leveraging academic capital formation to increase equity and inclusion in Calculus 1***

Having identified inequitable outcomes, we sought to identify corresponding programmatic and systemic barriers and planned interventions that would reduce them. Many barriers seemed to

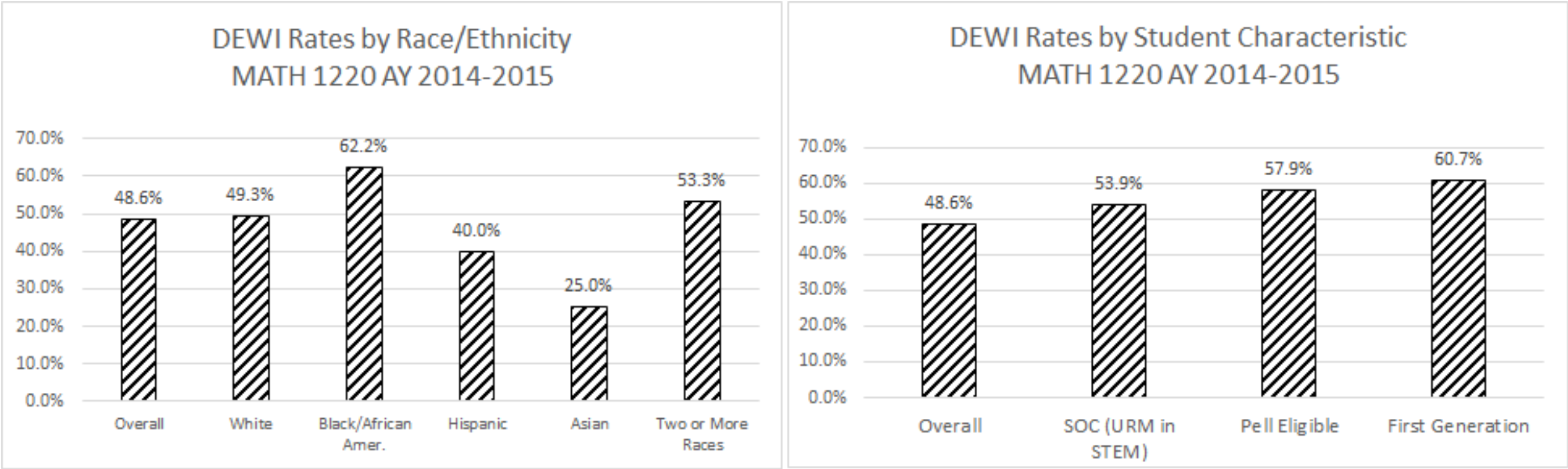

*Figure 1. Dashboard of DEWI rates in Calculus 1 (MATH1220) by Student Characteristics.*

*Note.* DEWI rates in Calculus 1 for student groups of concern: race/ethnicity (left) and students of color who are underrepresented in STEM, Pell eligible, and first generation (right) in 2014-2015 (year before start of redesign). These groups have DEWI rates significantly above the overall DEWI rate.

be related to academic capital. Academic capital is the "social processes that build family knowledge of educational and career options and support navigation through educational systems and professional organizations" (St. John, Hu, & Fisher, 2011, p. 1) and can also include knowledge and skills in an educational environment. Academic capital is an amalgam of tangible assets (money, books or technology) and intangible assets (access to networks, academic preparation and critical literacy, self-navigation/social agency, or college knowledge). Economic, cultural, and social capital affect the ability of an individual to acquire academic capital. Winkler and Sriram (2015) stated that "educational inequality is reproduced through the institutionalized process of rewarding the activation of social and cultural capital possessed primarily by privileged communities." Systemic inequities block some communities from developing academic capital which in turn makes these communities more vulnerable to existing systemic inequities.

---

[46] Underrepresented groups (URM): WMU uses the Higher Education Act definition of URM. For STEM fields, URM includes females. Students of Color (SOC): In this case study, students of color refer to any student that does not self-identify as white-single race. These include students that identify as American Indian, Alaskan Native, Asian, Black (not of Hispanic origin), Hispanic (including persons of Mexican, Puerto Rican, Cuban, and Central or South American origin), Pacific Islander, or other ethnic groups.

[47] This includes Black/African American and Hispanic/Latinx students, but does not include Asian, international, or female students who are white-single race.

Despite their acumen, grit, and interest in STEM occupations, students do not have equal access to opportunities and experiences that improve their chance for success. Familial role models that help negotiate classroom norms and institutional procedures are scarce for first generation students. Funding needed to support academic growth may not be available to students from lower socio-economic communities.   When these inequities are added to a higher education system that is more accessible to the affluent/privileged communities, the value of first generation and lower socio-economic students is unrealized.

In 2015, an intercollege calculus task force studied the issues surrounding the high DEWI rates. A team (a mathematician, a collegiate mathematics educator, and an experienced instructor) were tasked with improving DEWI rates and accomplishing the content goals of the course. Although close coordination was not a departmental norm, the team was encouraged to coordinate across sections.

The assigned team believed that thoughtful course redesign and coordinated instruction could increase some aspects of the students' academic capital. The redesign focused on the structure of the course (college knowledge), curricular (critical literacy and academic preparation) and instructional elements (engagement in learning), and support structures for students (self-navigation/social agency, college knowledge). The team also recognized that students at this age have not fully developed executive cognitive function. They are in the process of developing self-regulation skills and the ability to accurately assess the quality and depth of their own learning. McGuire (2015) argues that one way to help build students' academic capital is to provide them with specific training on metacognitive strategies. Cross-cutting the elements included in the course redesign were specific training in and development of metacognitive tools.

### 4.1. *Equitable Learning Environment: Coordination of Calculus 1*

Since the wide range of DEWI rates across sections pointed to a possible source of inequity, a more uniform approach across all sections was undertaken. In the Spring 2016, the aforementioned team began teaching the course. The team implemented several common elements such as course syllabus, class schedule (order and pacing), topics, student learning outcomes, grading scale and scheme[48]. An aspect of college knowledge is students' abilities to use information to determine pathways through educational systems. Knowledge of the content of courses, course sequences and the expectations of professors are information resources to which students may have minimal access or access to that knowledge may have limited impact on their choices because of other circumstances.  The team believed commonality was critical to "smooth out" some of the course variability that arises naturally from instructor choice and could lead to divergent outcomes for the student.

Throughout the semester, the instructional team met weekly, developed common formative and summative assessments, and established common mastery tests and grading rubrics. Weekly meetings allowed the instructors to pace one another, to check depth and breadth of coverage of topics, and to discuss difficulties students were having with the material across sections. It was the instructors' hope that this would facilitate students networking across sections. The instructors supported each other in choosing interventions to reach out to students before they began to struggle and shared knowledge of resources that could help students navigate through their programs.

[48] Grading scale is the percentage cutoffs for grade. Grading scheme is the assessment composition used to determine the grade (i.e., 30% exams, 20% homework…).

The team committed to increasing active participation of students during and outside of class. Group activities were used to introduce or reinforce material. The instructional team encouraged students to work collaboratively to improve their understanding. These changes targeted the development of students' academic engagement, social agency, and sense of belonging.

### 4.2. *Equitable Learning Environment: Student Preparation for Calculus 1*

Despite having satisfied prerequisites, anecdotal evidence suggested many lacked fundamental precalculus skills. The self-remediation tool ALEKS[49] was implemented to measure prerequisite knowledge early in the semester. The scores on this assessment showed student performance ranged from 9.4% to 84.2% on prerequisite material. Uneven preparation is a source of inequity. Unaddressed gaps in algebra and precalculus content could impact student success. The team added a mastery test that addressed prerequisite content. Students were given multiple attempts to achieve a high standard. ALEKS provided support for learning missing material, and students were encouraged to continue filling gaps.

### 4.3. *Equitable and Inclusive Learning Environment: Developing and Communicating Expectations*

Instructor expectations vary as a result of instructor topic strengths and preferences, as well as their experiences as students and teachers. Initially, the disparity of student expectations among the team members was significant, and frequent communication was essential to establish uniform expectations before the material was presented to the students.

With the addition of new members to the team, the revision of common expectations for students was necessary. To improve student and instructor experience, it was found that transparency in student and team expectations along with timely effective communications of these expectations was imperative. The team developed mechanisms designed to improve communication to students and to set a standard across sections.

Instructor expectations of student's level of justification and nature of evidence required in their solutions were communicated through writing assignments. These assignments were designed using the Transparency in Learning and Teaching Framework (Winkelmes, 2019; Mingus & Koelling, 2020) and used to guide students on how to write solutions. A feedback cycle helped the instructors to clarify and students to understand instructor expectations. While the incorporation of transparent assessments and increased student-instructor communication has the potential to improve all students' performance, these types of interventions are particularly helpful for students with less academic capital (Winkelmes, 2013; Winkelmes et al, 2015).

---

[49] ALEKS, from McGraw Hill, is an online adaptive learning assessment tool that provides feedback on the student's knowledge readiness levels and presents a path and an option for remediation.

Expectations of proficiency with skills and notation were addressed by many instructors using "rounds." Some students experience barriers because they have not had the opportunity to learn how to use appropriate notation to express their ideas. In rounds, the class is divided into groups, and each group is assigned a different problem to work on the board. After discussing all examples and analyzing and correcting mistakes, another round of problems is assigned. In our experience, rounds significantly improve student facility with skills and use of notation to express their ideas.

Expectations of work outside class were communicated by weekly checklists. In addition to listing required tasks, a checklist included suggestions of other tasks students could do to aid their learning: prompts to visit learning assistants or tutors, videos to watch for extra help, just-in-time assignments in ALEKS to remediate key skills, and suggestions of extra problems.

### 4.4. *Equitable and Inclusive Learning Environment: Supporting Students with Learning Assistants*

An instructor is often very different from the students, and there is a need for bridges between instructor and student to improve student involvement and increase inclusion of all learners.

Beginning in Spring 2016 semester, undergraduate learning assistants (LA) have been assigned to most sections of Calculus 1. LAs were provided with pedagogical training. LAs had excelled in the course and maintained a high overall GPA. While calculus content knowledge was important, the ability to be a co-learner was paramount. LAs supported student learning and facilitated discussions inside and outside the classroom. They met weekly with their instructor. They co-facilitated classroom activities, sent weekly emails/reminders, and conducted scheduled help sessions.

Our LAs helped create more inclusive educational experiences. First, LAs helped to include the mathematical reasoning of all students. In conversation with groups of students, LAs would work to understand student reasoning and support students in extending their thinking. Second, LAs helped to include student experiences in the management of the course. When meeting with the instructor, LAs would share student experiences. They gave the instructor insight into the student experience and helped the instructor understand how to improve those experiences.

Selection criteria attended to the ethnicity/race and gender representation of the LAs and resulted in a distribution similar to the population of students in the course. This provided successful role models of underrepresented groups in STEM fields that helped students visualize themselves as successful college students. LAs also acted as mentors and cheerleaders. LAs provided: study strategies and insights into the motivation and expectations of assignments; advice on how to make good choices in student-instructor relations, sources of information (advising, instructor, etc.), class selection, and time management; and support for academic work and taking action when needed (seeing instructor or advisor, finishing work).

LAs helped create a more equitable learning experience. Teams of LAs cooperatively conducted reviews for examinations. This helped compensate for instructional differences among sections. LAs scheduled additional student hours at more convenient (i.e., late and weekend) times and locations. LAs helped students determine what resources were needed and how to access them.

LAs helped students develop their academic capital, especially with respect to their “college knowledge” and content/academic preparation. LAs had a student level understanding of university processes such as registration, using the e-learning system, advising and programmatic requirements. LAs provided a peer resource to help students negotiate their pathway through the educational system. LAs held informal and impromptu advising sessions for students and connected them with advisors. Switching the course to a more student-directed, activity based format caused a level of discomfort among the students. LAs were able to help instructors “renegotiate” the classroom norms with the students. The shift in assignments that required a higher level of reasoning and writing was another area where LAs were able to facilitate communication between the instructor and the students regarding expectations and goals of the assignments.

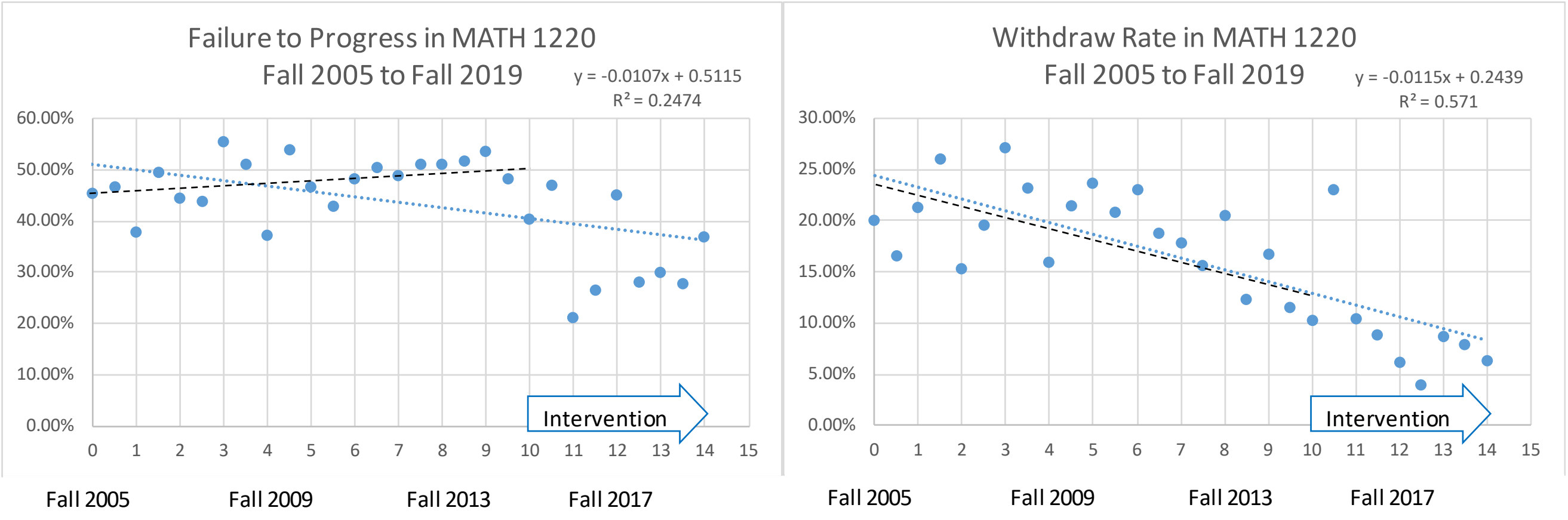

*Figure 2. Dashboard of failure to progress and withdraw rates in Calculus 1 (MATH1220).*

*Note*. Failure to progress (DEWI, including withdrawals) and withdraw rates from Calculus 1 from Fall 2005 to Fall 2019. Dashed line: least squares fit line (Fall 2005 to Fall 2015). Dotted line: least squares fit line (Fall 2005 to Fall 2015). Failure rates changed to decreasing while withdraw rates remained decreasing.

An institutional research study (Fitzpatrick, Hearit, and Easley, 2019) showed that students who met with their learning assistants at review sessions performed significantly better in the course. The effect of meeting with a learning assistant increased with each additional meeting. The study found that meeting with an LA at least once, rather than not at all, was linked with an average grade higher by more than half a letter grade (0.54 with $p<0.01$).

## 5. *Progress on Enhancing Diversity, Inclusion and Equity*

### 5.1. *DEWI Rates over Time*

Since the redesign and collaborative coordination of Calculus 1, the overall DEWI and withdraw rates in the course have declined. In spring 2016, a teaching team began collaborating, and all team members remained on the team through spring 2017. This team was particularly cohesive and collaborative; the DEWI rates dropped significantly. As new members joined the team, the process of integrating new ideas and coming to consensus created challenges for the team members and for timely communication to students; the DEWI rates were higher in those semesters. There is an overall downward trend in the DEWI rates (Figure 2).

Despite the effort to build a more uniform classroom experience, the spread of DEWI rates across sections has not decreased as the team had hypothesized (Figure 3). Some of this variability is to be expected (e.g., honors sections in fall semesters, cohorts of engineering

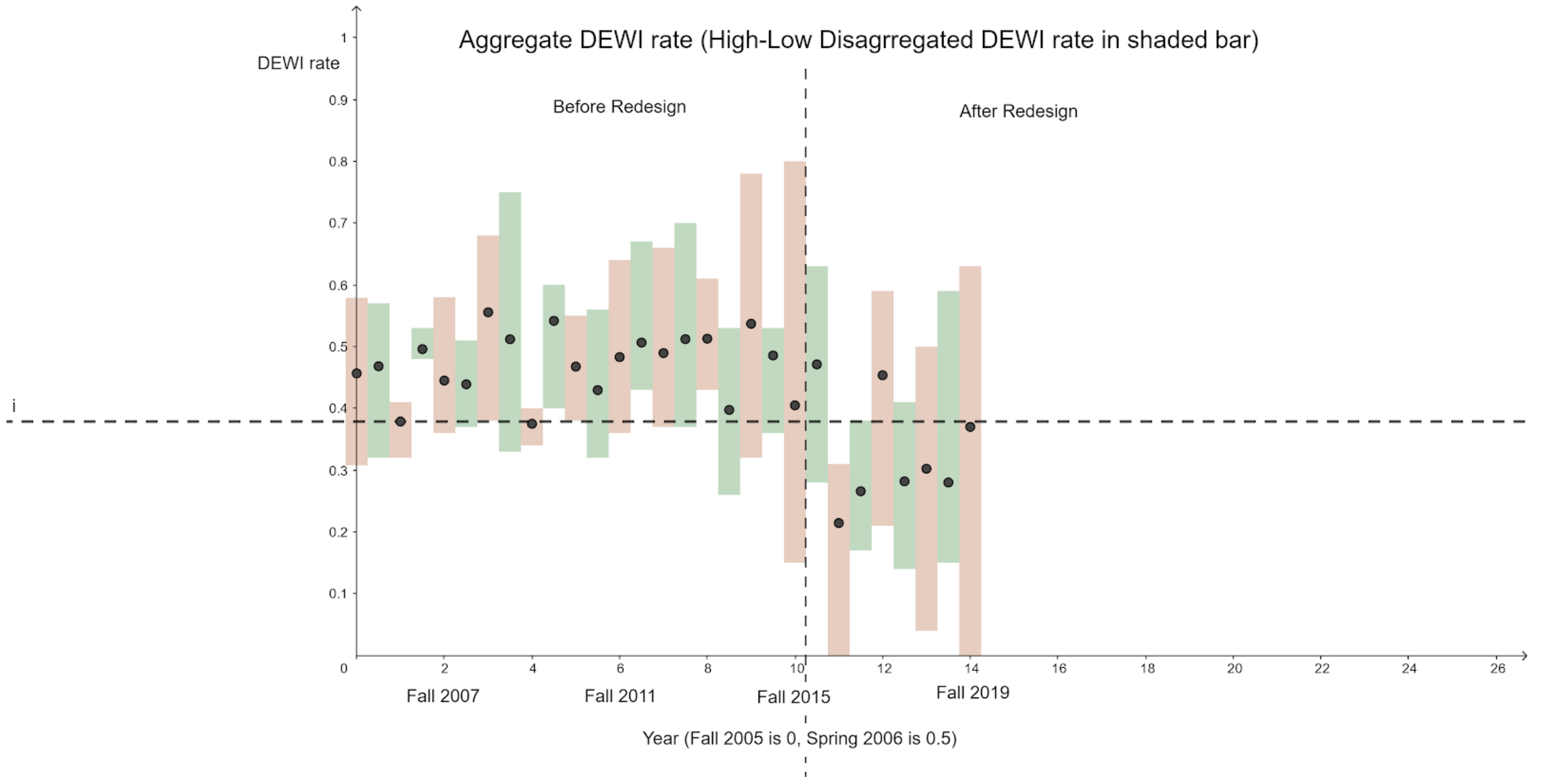

*Figure 3. Failure to progress rates 2005-2019 with disaggregated ranges.*

*Note*. Failure to progress (DEWI) rates for fall (brown) and spring (green) semesters from Fall 2005 to Fall 2019 (includes post-intervention). Dot shows semester average and bar extends from low to high DEWI rates across sections. Horizontal line shows lowest DEWI rates pre-intervention. Post-intervention, most semesters have average DEWI rates below the line. There remains a large variability in DEWI rates across sections and semesters. Adapted from figure in Mingus & Koelling (2020).

students placed into sections). Some of this variability is also due to the degree to which instructors understand, believe in, and are able to implement common priorities.

## 5.2. *DEWI Rates for Special Groups*

The DEWI rates are also lower for special populations (Figure 4). Female students' DEWI rates improved. There is improvement in the DEWI rates of students who are classified as racially minoritized and underrepresented students, Pell Eligible, or First Generation, but they continue to bear the burden of the DEWI rates.

## 5.3. *Progression Rates over Time*

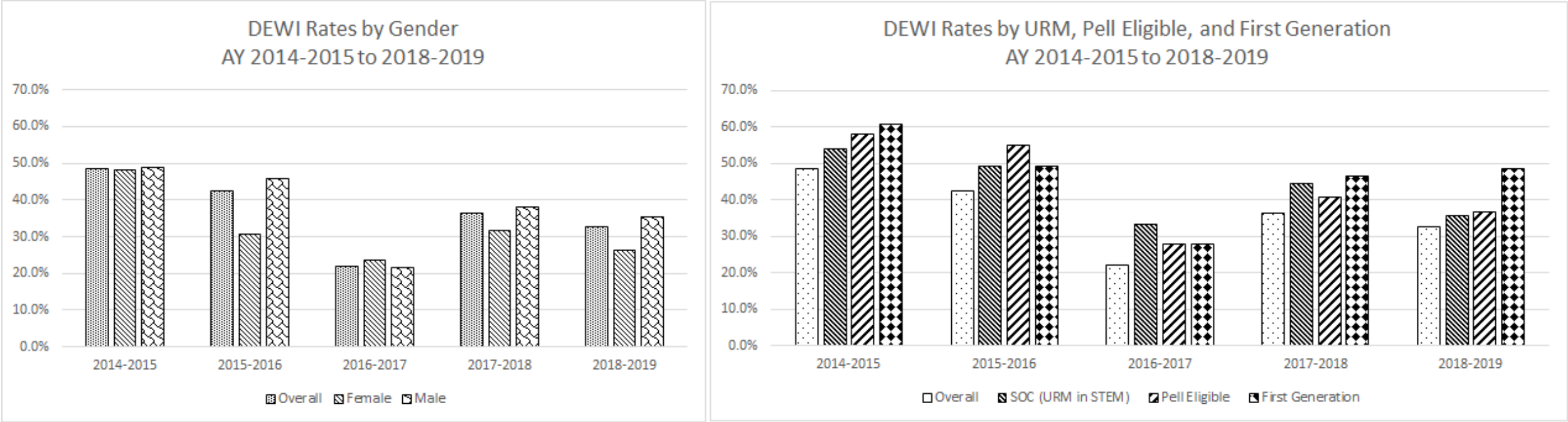

*Figure 4. Dashboard of DEWI rates for special groups of students.*

*Note*. DEWI rates for different student characteristics from 2014-2015 to 2018-2019. See footnotes 7 and 8 for definitions of URM and SOC. DEWI rates are lower post-intervention (Spring 2016-present) for all groups. Data for Pell Eligibility and First-Generation status were unavailable for AY 2013-2014.

In addition to improving success rates in Calculus 1, the team maintained high academic standards. To study our standards, we examined student progression[50]. The number of students progressing naturally declines because some majors require fewer semesters of calculus. To track student progression rates over the calculus sequence, the team examined the cumulative percent eligibility (CPE). CPE is the percent of students starting Calculus 1 who earn a grade in Calculus 2 that makes them eligible to start Calculus 3.[51]

Over the years 2013 to 2019, there is no decrease in the percentage of students achieving eligibility for Calculus 3. This suggests the improvement of the DEWI rates in Calculus 1 is not an artifact of lowering the academic standards in the course.

## 5.4. *Progression Rate for Special Groups over Time*

Between 2013 and 2017, CPE was increasing for all groups except international students (Figure 5). There was a dip in CPE in 2017-2018 for all groups. New faculty were incorporated into the team for the first time in 2017-2018. The original team developed and implemented changes for three semesters. New team members naturally had their own ideas of how Calculus 1 should be taught, and common expectations had to be adjusted to reflect the new makeup of the team. We did not appreciate the extent of the disparity of expectations. Since expectations were discussed as the semester progressed and material was taught, it was difficult to communicate them to students in time for the students to meet these expectations.

---

[50] For the data shown, students took Calculus 1 in a fall semester and progressed into Calculus 2 in the following spring semester for the academic years 2013-2014 through 2018-2019: redesign of the course began in spring 2016, thus three of these academic years were before the redesign and three were after. Students eligible for a class earned a grade of C or higher in the prerequisite. Pell-eligibility and first-generation status were unavailable for this data analysis.

[51] CPE is the percent corresponding to the fraction of students who start Calculus 1 which is the product of (1) the fraction of students who earned a grade in Calculus 1 out of those who start it, (2) the fraction of students who earned a grade that makes them eligible for Calculus 2 out of those who earned a grade in Calculus 1, (3) the fraction of students who earned a grade in Calculus 2 out of those who were eligible for it , and (4) the fraction of students who were eligible for Calculus 3 out of those who earned a grade in Calculus 2.

As a team, we are learning how to incorporate new members, build common expectations, and support students as they try to meet those expectations. The CPE recovered in 2018-2019, and we hope that what we have learned about working with new team members will allow us to maintain higher levels of CPE. Further, if the academic year 2017-2018 were removed from

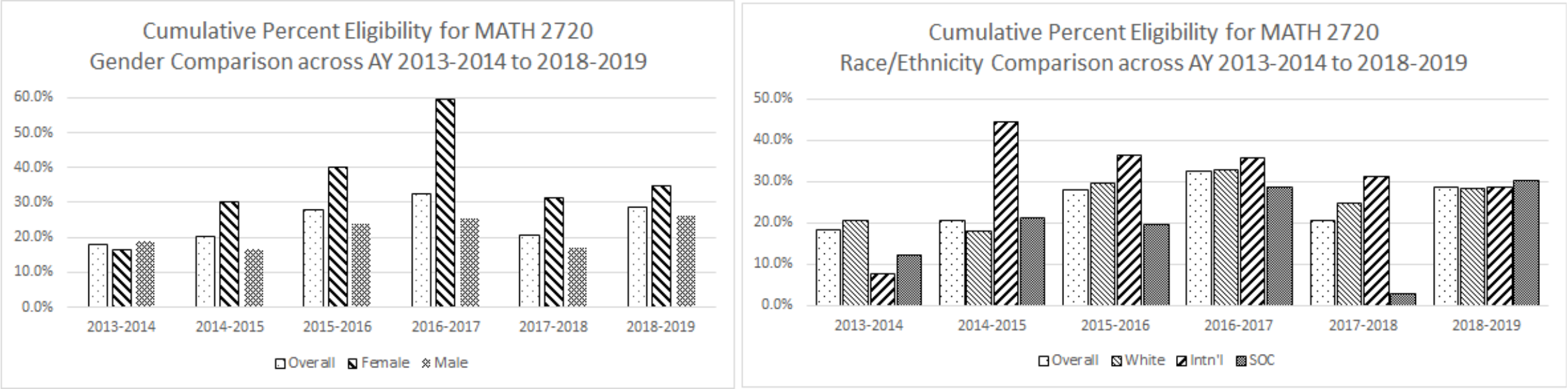

*Figure 5. Dashboard of cumulative percent eligibility for special groups: Calculus 3 (MATH2720).*

*Note*. Cumulative Percent Eligibility for special groups between 2013-2014 and 2018-2019. See footnotes 7 and 8 for definitions of SOC. For two of the three years after the redesign, the cumulative percent eligibility is more equitable. We hope that the lessons learned in 2017-2018 will help us keep it more equitable.

consideration, the CPE has become more equitable. Academic year 2018-2019 was particularly equitable. In most semesters, Calculus 2 instructors did not coordinate. However, 80% of Calculus 2 instructors in 2019 coordinated together and had been members of the Calculus 1 teaching team. This suggests that increased coordination across the entire calculus sequence would be beneficial for student success.

### 5.5. *Academic Capital*

There has likely been an increase in the critical literacy, academic preparation, and subject knowledge facets of academic capital. The decrease in DEWI rate and the maintenance (or increase) of CPE may indicate an increase in academic preparation and subject knowledge. In addition, there is evidence of increased critical literacy: ability to use mathematical reasoning to draw conclusions and investigate problems. A department assessment tool showed evidence of improved ability to communicate mathematical reasoning in writing from before to after the redesign (Levin, Mingus, Strom, & Wood, 2019). Instructors have also noticed increased student performance on problems that involve critical literacy in high stakes assessments.

Over multiple semesters, we have anecdotal evidence that some students are better able to navigate the college experience. Students have been more proactive in seeking information about programs, courses, instructors, and employment opportunities available in the mathematics department.

There is also anecdotal evidence of increased sense of belonging. Students have maintained working groups established in Calculus 1 in the next two courses in the sequence. Students have purposely scheduled classes to stay together, and students have continued to reach out to LAs and instructors after Calculus 1 was completed. Students have utilized the student success center and tutors rooms as a place to gather and collaborate. Students are more willing to have conversations about mathematics, their academic progress, their student experience, and possible career paths with their instructors and LAs. Students are showing increased persistence through the calculus sequence. Students are actively seeking to serve as learning assistants for lower-level mathematics courses with the hopes of maturing into LAs for calculus.

## 6. *Conclusion*

Gateway courses have a crucial role in developing students' understanding of the structures and academic expectations of higher education institutions. As we developed our understanding of our students' backgrounds and our teaching, we aimed to provide an equitable experience through standardizing course structure and coordinating instruction. We incorporated curricular, pedagogical, and holistic support systems to improve the experience and increase the academic capital of all students. The incorporation of learning assistants, active learning, and low-stakes assessments with high transparency of expectations effectively built academic capital. Students had a greater sense of belonging and felt more comfortable communicating their ideas.

Collaborative coordination is a foundation from which we can further decrease course and programmatic inequities and barriers. We have created a group of instructors who implement common features of the course built to address some known causes of inequity. This foundation prepares us to look for and to address other causes of inequity.

Coordination can also support instructors in the critical self-examination required to foster the acquisition of academic capital. Instructor attitudes and behaviors influence the development of student academic capital. However, our attitudes and behaviors are not always in concert with our beliefs. This inconsistency can be hard to recognize and reconcile. The teaching team becomes a support group to help instructors negotiate the changes they need to make in their teaching to better align their actions with their beliefs.

Gains in content knowledge, academic capital, self-efficacy, and sense of belonging can be ephemeral if they are not reinforced by faculty and university. Incorporating support structures into gateway courses can solidify these gains. Although faculty want to maintain high standards in their courses and ensure that the transition between courses in a sequence is smooth, individual faculty working alone is insufficient. More faculty and student support is needed to foster productive transitions.

With such student support, more minoritized and underserved students can succeed in their courses, thus building academic capital for their families and communities.

# 21 A Multi-pronged Approach for Closing Opportunity Gaps in Calculus

**Jesús Oliver**

*California State University, East Bay (CSUEB)*

Jesús Oliver received his Ph.D. from the University of California, San Diego. Prior to coming to CSUEB, he was a Lecturer at the University of California, San Diego and a Postdoctoral Fellow at the Institut Henri Poincare in Paris, France. Prof. Oliver's primary teaching interest is the development and implementation of active learning techniques in the Mathematics classroom.

**Julia Olkin**

*California State University, East Bay (CSUEB)*

Julia Olkin received her Ph.D. in Mathematical Sciences from Rice University in Houston, TX. She has been a Professor at CSU East Bay for 14 years. She has extensive experience running professional development workshops for K–12 teachers and faculty learning programs for college instructors, emphasizing incorporating an active learning approach.

**Andrew Stanciulescu**

*California State University, East Bay (CSUEB)*

Andrew Stanciulescu - first generation American and student researcher at California State University East Bay. I am majoring in Computer Science and Statistics with a Data Science concentration.

*Abstract: The main goal of the SEMINAL project at California State University, East Bay (CSUEB) is to create, support and disseminate active learning pedagogy in Calculus I. The SEMINAL CSU East Bay project seeks to fulfill these goals via a three-pronged approach: course redesign and curation, coordination with a gentle touch, and the cultivation of a robust community of practice. This multi-pronged approach has been implemented, via Dynamic Calendar, while closely attending to issues of diversity, equity, and inclusivity throughout. The mathematics department has recently seen a sizable reduction in the gap between passing rates for Black, Latinx, and Native American (BLN)[52] students versus non-BLN students over a three-semester period. Further data disaggregation also suggests higher rates of achievement at the top end of the grade spectrum in Calculus I for female students versus male students.*

## 1. *Introduction*

The Department of Mathematics at California State University East Bay (CSUEB) is a Phase 2 partner in the Student Engagement in Mathematics through an Institutional Network for Active Learning (SEMINAL) project funded by the Association of Public and Land-grant Universities and the National Science Foundation. The SEMINAL national effort seeks to study how to enact and support institutional change fostering the use of active learning mathematics. One of the main areas of emphasis of the (local) SEMINAL project at CSUEB has been to create, support and disseminate active learning pedagogy in all sections of Calculus I. The implementation of this framework also aims to promote inclusive teaching by attending to issues of diversity, equity, and inclusivity.

[52] CSU East Bay currently refers to this group as historically underserved students (HUS).

The SEMINAL CSU East Bay project seeks to fulfill these goals via a three-pronged approach: (1) course redesign and curation, (2) coordination with a gentle touch, and (3) the cultivation of a robust community of practice. Specifically, the following suite of changes has been implemented as part of the SEMINAL grant: a Mathematical Community of Practice (MCoP) for Calculus I that includes orientation for new instructors and monthly meetings; a Dynamic Calendar which acts as a pacing guide for the courses, as well as a repository of all classroom activities; and a redesigned curriculum to incorporate student engagement in the classroom (Oliver & Olkin 2020). The SEMINAL CSUEB leadership team and the MCoP also work towards building an inclusive, student-centered classroom community. Examples of how our community seeks to work towards this goal include: encouraging all instructors to learn students' preferred names through the use of nametags on desks; in the online environment we also make extensive use of the chat feature in Zoom as we notice that students frequently feel more comfortable writing their thoughts and questions rather than turning on videos and/or speaking out loud. Asking everyone to put an answer into chat elicits a much better response and allows the instructors to encourage some of the shyer students to participate. In a face-to-face environment, many instructors from the MCoP mix up the student groups during activities so that students gain experience working with a variety of different people, to build classroom community, and to ensure that all voices are heard.

The glue that ties together the three-pronged approach mentioned above is the Dynamic Calendar. This electronic calendar serves as a repository for essential instructional resources. One of the features that makes it "dynamic" is that it includes embedded live links to a set of in-class activities, videos, and worksheets tied to a weekly topic list. For example, the calendar includes links to digital worksheets used to facilitate scaffolded problem-solving, templates for formative assessment activities, and video assignments for flipped lectures. Because of the opt-in nature of the calendar structure, instructors can implement active learning strategies based on their comfort level and experience, reducing the barrier of entry for those new to exploring the use of active learning in their courses (Auster & Wylie 2006). When used in this manner, the Dynamic Calendar helps support course coordination with a gentle touch. The calendar is also dynamic in the sense that it can be edited based on instructor and student feedback, availability of new resources, and changes in topic selection. This allows for a community-wide adaptive instructional approach in which broader pedagogical recommendations can be made on a large scale to best meet students' needs. For example, in order to support students coming in with different levels of preparation, members of the MCoP are building a library of short videos for just-in-time learning, with the goal of covering review material that some students may feel weak in. Instructors also post suggestions in CSUEB's Learning Management System about practice websites such as the free problem generator in https://www.wolframalpha.com/pro/problem-generator/ to aid in just-in-time instruction. As for pacing of the course, during one semester, when all the Calculus I instructors noted that they pushed the first exam to a later week, that new pacing was reflected in an updated calendar. Two screen shots of the Dynamic Calendar are shown in Figure 1. The left screen shot shows details of Week 1, plus the color-coding in the calendar, and the right screen shot shows details further along in the course. All underlined text are live links to documents, activities, videos. For all exams there is a link to a repository of sample exams and practice reviews, especially useful for a first-time instructor of Calculus I.

California State University, East Bay
Department of Mathematics
Math 130: Calculus I
Calculus, J. Stewart, 8th Edition

Group Activity ABCD Activity Comment Video Diversity, Equity, Inclusion Activity
Mini Lecture includes group work

| Week | Section | Comments |
|---|---|---|
| 1 | <br>1.1<br><br><br><br>1.2<br>1.3<br>1.3<br>1.4<br>1.5 | Diversity/Inclusion: Personal Math History - Assign as homework.<br>Review Activity - review fundamental concepts<br>Group Activity: Diagnostic Activity + Diagnostic Activity Solutions<br>ABCD Activity Review Domain and Functions<br>Suggestion: review horizontal, vertical translations and shifts of function graphs during the semester as needed.<br>Mini lecture: Linear functions<br>HOMEWORK: Composition of functions: assign two videos for students to watch for 2nd class period, instead of doing a full lecture<br>• Intro to function composition(Video 1)<br>• Function composition (Video 2) |
| | | ABCD Activity: Function Composition (students have watched videos)<br>Mini lecture: exponential functions<br>Group Activity: Exponential vs Linear Functions<br>Sect1.4-ExponentialFunctions (CarGas) |
| 2 | <br>1.5<br><br>2.1<br><br><br><br>2.2 | Growth Mindset HW (Video + Google survey) (Word doc).<br>Mini Lecture: Inverses<br>Mini Lecture: Logs<br>Group Activity: Logarithms Sect1.5-Logarithmic Functions (CoffeeMug -- 2 questions).<br>Comment: Tangent line concepts are covered again in Section 2.7. It's Ok to go quickly. Focus should be on conceptual understanding of describing the position of a particle as a function of time. Emphasize average velocity here.<br>Mini Lecture: Definition of a limit<br>HOMEWORK: Reminder- tell students to bring electronic devices to next class. |
| 3 | <br><br>2.2<br>2.3 | Diversity: Living Proof Essays - could start assigning first one (4 essays throughout semester)<br>Mini Lecture on definition of limit (application, slope of tangent line);<br>Group Activity: Use electronic devices to run Desmos activity on Limits (teacher.desmos.com -- limit activity) (Winner!) |

*Figure 1: Two screen shots of the Dynamic Calendar*

In addition to serving as a repository for active learning activities and as a pacing guide, the Dynamic Calendar is an important tool for supporting the implementation of Diversity, Equity, and Inclusivity (DEI) approaches in the classroom. In particular, a sampling of DEI activities developed by our department are included throughout the Calendar, color-coded in turquoise. These activities seek to strengthen a sense of belonging, reinforce the growth mindset idea, and ask students to reflect on their personal background and prior experience with mathematics. To further inform discussions of DEI pedagogy in the classroom, our department worked on a definition of "teaching math equitably." This definition is codified as a guideline at the bottom of the Dynamic Calendar, shown in Figure 2. Here, we reinforce our role as instructors in the classroom, and remind everyone ways in which we can avoid propagating inequities in the classroom.

**Teaching Math Equitably** begins with the recognition that inequities due, among other things, to an individual's race and gender exist in many facets of our society, including within the teaching and learning of mathematics. As instructors we have control over the policies and norms we create in our classrooms. It is our duty to reflect upon our practices and be intentional in the creation of policies and norms that (1) halt the propagation of inequities and (2) cultivate equity. Examples of policies that may propagate inequities are:

- Zero tolerance policy for late homework
- Heavy emphasis on timed, summative exams.
- Heavy emphasis on traditional lecturing.
- Intimidation of students and/or establishing destructive hierarchies in the classroom.

Mathematics Communities of Practice

These documents were developed by the 2018-2020 Mathematics Community of Practice (MCoP) with support from the National Science Foundation and the CSUEB College of Science.

*Figure 2: DEI Statement in the Dynamic Calendar*

The Dynamic Calendar grew out of a need from the redesigning of Calculus I, and then later Calculus II and Precalculus. As we created a suite of activities, group work, and other tasks, we needed a convenient way to communicate and share this curriculum. Additionally, it serves as a useful orientation tool for new instructors. An example of a suite of in-class tasks, we created a packet of ABCD activities specifically for Calculus I. The ABCD cards are color-coded cards with the letters A, B, C, D, shown in Figure 3, that can be folded in such a way to show just one letter at a time, and can be used by students to communicate a response to a multiple choice question. Additionally, the entire card can be held out to signify "I don't know," which is useful information to gather from the class. In these activities, the multiple choice questions are worded to probe for conceptual understanding. In particular, knowing the "right" answer is not enough. This is due to the fact some of the answers are along the lines of, "yes, because …" or "x=2, because…" so that students are forced to be sure they know why they are choosing a particular answer.

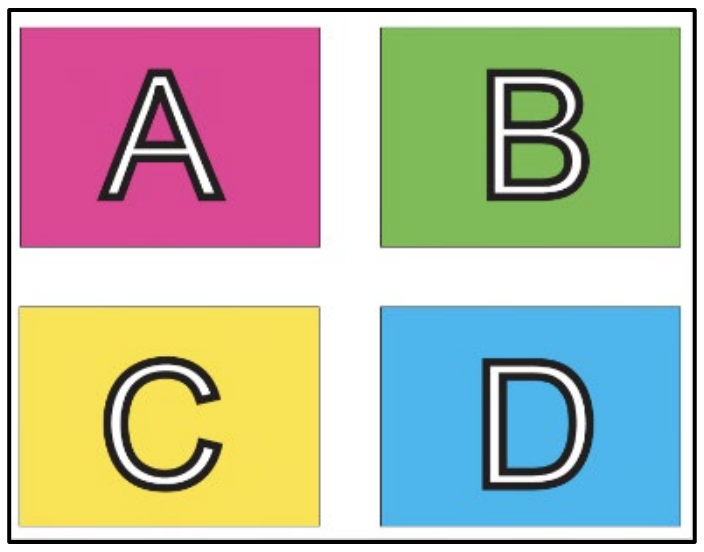

*Figure 3: ABCD Cards*

These tools do not exist in isolation within the CSUEB math department. CSUEB has a robust Mathematics Community of Practice that meets regularly, and supports each member using the resources outlined above. The meetings are generally grounded around issues of pedagogy, assessment questions, pacing, design and use of activities, as well as affective support for instructors. In this manner, the MCoP provides multiple points of contact for instructors seeking to incorporate active learning in their classrooms. The MCoP meetings are well-attended: at CSUEB, all but one instructor from all sections of Precalculus, Calculus I and II instructors attend fairly regularly the monthly meetings. MCoP meetings also provide an opportunity for instructors to have conversations about issues of equity in the classroom such as opportunity gaps in multi-section courses as well as differences in instructor grading criteria. The MCoP has had numerous discussions around creating more equitable assessments. In particular, members discussed various approaches for implementing lower stakes assessments, allowing notes during exams, and writing conceptual questions that are less computationally taxing. In particular, during these discussions we discovered that there was a very wide range of what instructors were allowing in terms of notes during exams. We felt this was inequitable and discussed having a minimum policy of a 3x5 index card for all assessments in all classes. This policy change is now officially implemented by the department. Instructors are allowed more, but this provides a threshold and encourages questions that move away from memorization and more towards conceptual understanding.

Although the rest of this Case Study seeks to document how this multi-pronged approach has helped to close opportunity gaps in Calculus I (a measure of student achievement), the authors would like to highlight how the MCoP seeks to attend to other dimensions of equity. A useful way to think of the different aspects of equity is provided by Gutierrez (2009). According to Gutierrez, there are four dimensions of equity in a classroom: Access, Achievement, Power, and Identity. All four dimensions need to be attended to in order to move towards equity in the classroom (see also Section 2 of this volume for further discussion of Gutierrez's work). For example, the codification of "teaching math equitably" in the Dynamic Calendar addresses the issue of Access, which is to ensure that all students, ideally, have equal access to course resources, whether that be materials, good instruction, and robust curriculum. Providing equal opportunities for success, such as allowing notes on assessments, displays elements of both Access and Achievement, as described by Gutierrez. MCoP members also help each other with inclusive teaching practices as they pertain to Power in the classroom. For example, we discuss ways to encourage student engagement during whole-class and small group discussions. If an instructor keeps calling on the same students, that instructor is denying others of opportunity. In particular we highlight approaches that do not resort to cold calling, and maintain a safe space within the classroom. If instructors have issues with student group dynamics, we also brainstorm ways to improve the classroom environment.

## 2. *DFW Rates and Data Disaggregation*

Up to this point, the CSUEB SEMINAL project has seen some early success in lowering the departmental DFW rates for Calculus I. At CSUEB, we consider the DFW category to include grades of D, F, and WU (Withdrawal Unauthorized). For context, the math courses at CSUEB are capped at 35 students, and there are between six and eight sections of Calculus I taught each term. These courses are taught by a mixture of tenure-track faculty and lecturers. Figure 4 shows a term-by-term comparison of the average DFW rates and Table 1 gives the associated DFW numbers and percentages. Historically, the departmental aggregate DFW rate for Calculus I hovered around 36.4%, represented by the gold bar in Figure 4. A transition year (2017-2018) occurred before SEMINAL where the department received a Course Redesign Grant from the CSU Chancellor's office, and the work on redesigning Calculus I and incorporating more active learning began. This academic year is represented by the three green bars, as this was the last year CSUEB was on the quarter system before converting to semesters in Fall 2018. By contrast, during the SEMINAL grant funded period, the DFW rate has consistently stayed near or below 20%.

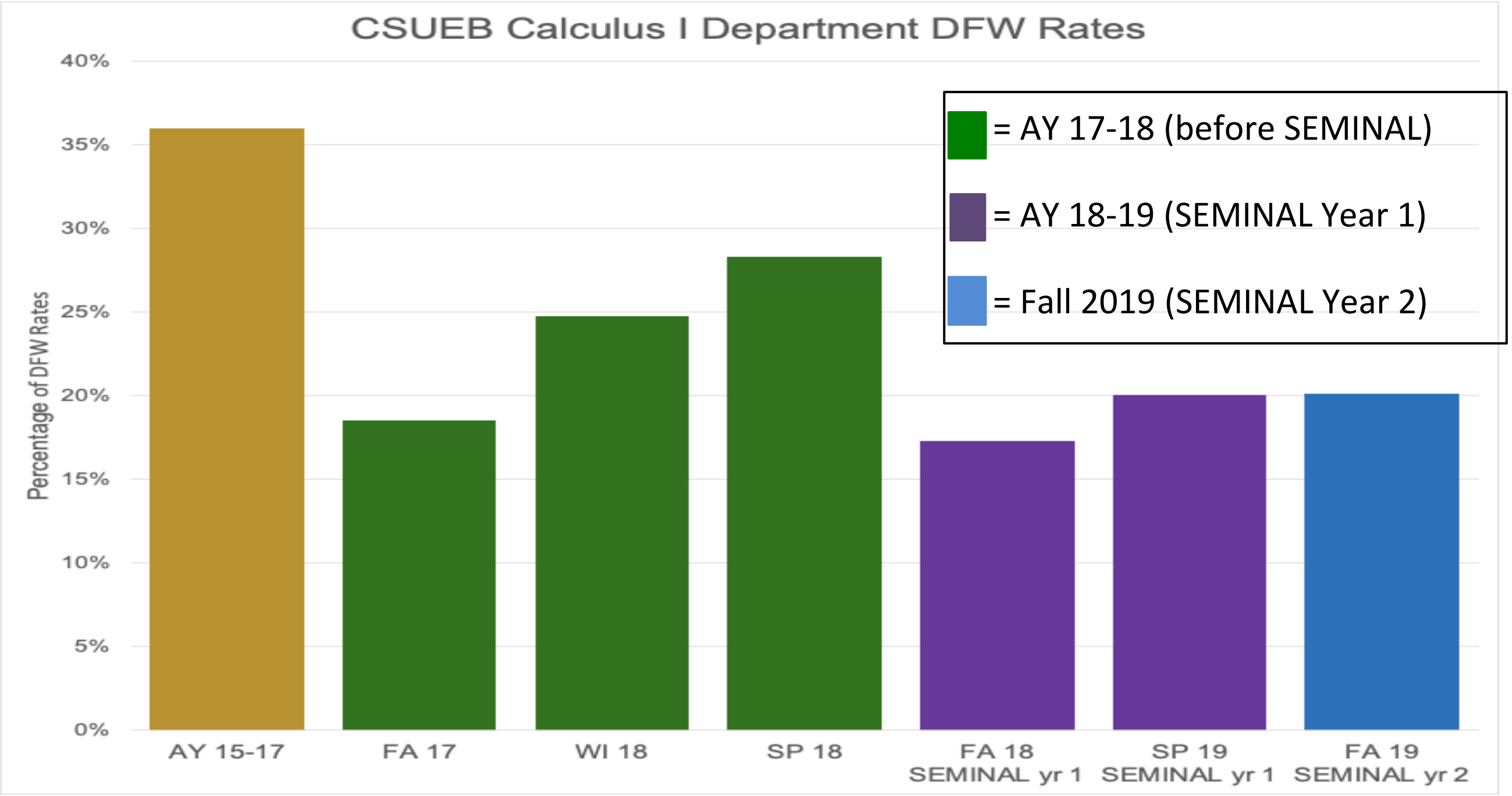

*Figure 4: Math Department DFW Rates*

*Table 1*
*Math Department DFW Rates*

| | N | #DFW | % |
|---|---|---|---|
| AY 15-17 | 922 | 336 | 36.4% |
| Fall 17 | 184 | 34 | 18.5% |
| Winter 18 | 178 | 44 | 24.7% |
| Spring 18 | 106 | 30 | 28.3% |
| Fall 18 | 231 | 40 | 17.3% |
| Spring 19 | 165 | 33 | 20.0% |
| Fall 19 | 189 | 38 | 20.1% |

While this is an improvement across the board, this data analysis raises the question of what effect, if any, the approaches implemented by the SEMINAL team at CSUEB are having on student groups that traditionally face the largest barriers to success. For example, are there student populations that are being negatively affected?

## 3. *Data Analysis: DFW Rates and the Impact on the Opportunity Gap*

Before diving deeper into the data and its interpretation, the authors would like to provide some important context: The SEMINAL effort at CSUEB focuses on deploying active learning approaches to teaching with the goal of helping to counter the inequities experienced by students due to, among other factors, race, and gender. Therefore, the leadership team believes that using the term "achievement gap" to describe the disparity between traditional and historically underrepresented students is ineffective, especially as it is framed as an implicit endorsement of the currently existing, inequitable educational system. Since the SEMINAL CSUEB effort seeks to increase opportunities for student success, we believe that the term "opportunity gap" more accurately reflects the leadership team's belief that it is the curriculum and the teaching that needs adjustment to better serve our students -- not the other way around.

To gain a sense for the impact of the project on students coming from BLN populations, the CSUEB SEMINAL team requested grade distribution data from the Office of Institutional Effectiveness & Research at our campus. In particular, we asked that the data be disaggregated by BLN versus non-BLN in Calculus I classes. The results from disaggregating data in this way are as follows: historically, the DFW rate for BLN students was 28.2% while for non-BLN students it was 18.1%. The opportunity gap, as measured by the difference of the non-BLN DFW rate from the BLN DFW rate, was 10.1%. A quick calculation shows that this difference was statistically significant (P= 0.0045). By contrast, during academic year 2018-2019, the DFW rate for BLN students was 17.9%, lower than the rate of 18.9% for non-BLN students. The opportunity gap thus improved to -1.0% during this academic year. This number suggests that BLN students achieved a passing grade in Calculus I at slightly higher rates than non-BLN students during academic year 2018-2019 and the opportunity gap was closed. Further calculation also shows that this difference was not statistically significant (P=0.4). Did these trends persist? In Fall 2019 the results were encouraging, but slightly less favorable: the DFW rate for BLN students was 20.7% while for non-BLN students it was 19.6%. The opportunity gap for this semester was 1.1% and that difference was not statistically significant (0.425). While this indicates that an opportunity gap reappeared during this semester, these results point to a sustained positive trend in closing opportunity gaps in passing rates across the board. This result is also in close alignment with recent findings from the STEM education literature that report that "active learning benefits all students, but offer disproportionate benefits for individuals from underrepresented groups" (Theobald et al., 2020). Lastly, we note that since non-BLN student DFW rates have stayed roughly around the 19%-20% range from Fall 2017 through Fall 2019, this data shows that *essentially all the gains* in lowering the DFW rate across Calculus I sections during SEMINAL grant funded period can be traced to increased academic success achieved by BLN students in Calculus I. Figure 5 shows this data as a bar graph, and Table 2 gives the associated numbers and percentages.

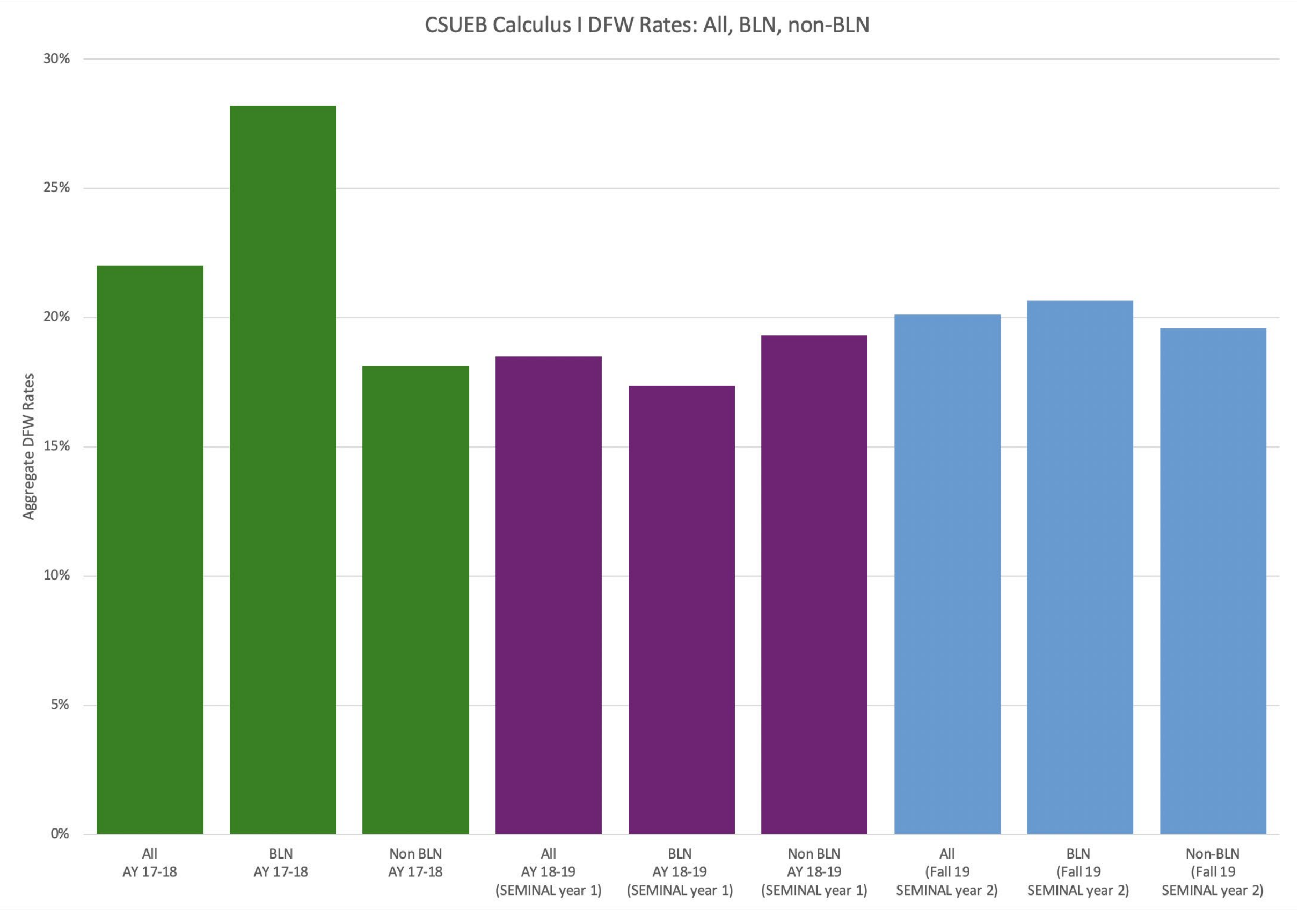

*Figure 5: Comparison of DFW Rates Among BLN and Non-BLN Students*

*Table 2*
*Comparison of DFW Rates Among BLN and Non-BLN Students*

| | N | #DFW | % | Opportunity gap |
|---|---|---|---|---|
| AY 17-18 All | 486 | 107 | 22.0% | 10.1% |
| AY 17-18 BLN | 188 | 53 | 28.2% | |
| AY 17-18 Non | 298 | 54 | 18.1% | |
| AY 18-19 All | 396 | 73 | 18.4% | -1.0% |
| AY 18-19 BLN | 168 | 30 | 17.9% | |
| AY 18-19 Non | 228 | 43 | 18.9% | |
| Fall 19 All | 189 | 38 | 20.1 | 1.1% |
| Fall 19 BLN | 92 | 19 | 20.7 | |
| Fall 19 Non | 97 | 19 | 19.6 | |

## 4. *Investigation: Master, Pass, and No Pass Rates*

Although the story so far is positive and the opportunity gap (measured by DFW rates) has decreased significantly, this still leaves open the question of what opportunity gaps, if any, we might have for students who pass Calculus I. Thus, we wanted to dig deeper and probe for potential equity gaps in grade data for students who passed this course. Based on examples from our SEMINAL partners at San Diego State University, we defined three categories for student grades: **Master, Pass and No Pass.** The Master category includes students who receive grades of A, A-, B+, B, or B-. The Pass category includes students who receive grades of C+, C, C-. The No Pass category includes students who receive grades of D, F or WU. Note that the definition of No Pass is the same as the definition of DFW grades -- this is done for consistency. These three categories allow for the study of potential opportunity gaps for students who earn a passing grade, while providing flexibility to further disaggregate the data by BLN or by gender.

In this subsequent analysis, we restricted ourselves to academic year 2018-19 Calculus I summary grade data. Table 3 shows the resulting breakdown. The No-Pass row reflects the -1% opportunity gap discussed above. Upon disaggregating this data by Master, Pass, No Pass categories, we see that the percentage of all students who master the course is 62.88%. Breaking this number down, for BLN students the mastery rate was 60.12% while for non-BLN the mastery rate was 64.91%. Computing the difference of these two numbers yields a mastery rate opportunity gap of 4.79%. We would hope that there would be no gap. As we continue to track these numbers in future academic years, if the gap persists, we aim to dig deeper into the potential reasons for this gap in the numbers.

*Table 3*
*Comparison of DFW Rates Among BLN and Non-BLN Students*

| AY 18-19 | n-ALL | %-ALL | n-BLN | %-BLN | n-nonBLN | %-nonBLN |
|---|---|---|---|---|---|---|
| **MASTER** | 249 | 62.88% | 101 | 60.12% | 148 | 64.91% |
| **PASS** | 74 | 18.69% | 37 | 22.02% | 37 | 16.23% |
| **NO-PASS** | 73 | 18.43% | 30 | 17.86% | 43 | 18.86% |
| **n-SUM** | 396 | | 168 | | 228 | |

*Master = grades of A, A-, B+, B,B-* *Pass = C+, C, C-* *No Pass = D, F, WU*

The SEMINAL CSUEB team performed a deeper analysis in order to further probe for potential opportunity gaps. Starting with the same summary grade data from academic year 2018-2019, the data was disaggregated by Master, Pass, No Pass categories as well as by student self-reported sex. Students self-selected from a choice of "female", "male" or "other" on their CSUEB profile. For Calculus I data, no student reported as "other" and thus, Table 4 lists two sexes. Table 5 shows that the mastery rate for all students in Calculus I was 62.66% whereas female students achieved a mastery rate of 69.88%. On the other hand, male students achieved a mastery rate of 57.83%. Thus, female students earned master-level grades at a rate of 12.05% higher compared to male students. When it came to Pass-level grades, the rates were 18.69% for all students, 12.65% for female students and 23.04% for male students, respectively. The No Pass rate (i.e., DFW rate) was 18.43% for all students, 17.47% for female students and 19.13% for male students.

*Table 4*
*Comparison Calc I AY 2018-19 based on Sex*

| AY 18-19 | n-ALL | %-ALL | n-Female | %-Female | n-Male | %-Male |
|---|---|---|---|---|---|---|
| **MASTER** | 249 | 62.88% | 116 | 69.88% | 133 | 57.83% |
| **PASS** | 74 | 18.69% | 21 | 12.65% | 53 | 23.04% |
| **NO-PASS** | 73 | 18.43% | 29 | 17.47% | 44 | 19.13% |
| **n-SUM** | 396 | | 166 | | 230 | |

The takeaway message from this academic year is that the female students are outperforming male students in Calculus I. Notice that if we had just concentrated on passing grades, we would have lost the story of variations in student grades. For this year, the male students were in the Pass group (C+/C/C-) at nearly twice the rate as the female students.

It might seem like the next obvious step is to further break down the data into female/male/BLN/Non-BLN. However, the more we break down the data, the higher the probability that we could encounter a small sample size. We opted to create a minimum cut-off for the number of students in each cell at 20. Our rationale was that below 20 students in a cell would not gain meaningful analyses and would run the risk of drawing incorrect conclusions. Finally, we thought having a minimum of 20 participants in each cell would be adequate as it represents 5% of the total number of students in the sample.

## 5. *Conclusion*

The SEMINAL CSUEB project has seen success in lowering the overall DFW rates via three major components: incorporating a robust community of practice, encouraging coordination via a dynamic calendar, and nudging all lecturers to attend to equity when incorporating active learning in the classroom. This success is also reflected when looking at the DFW rate opportunity gaps for BLN versus non-BLN students. Further disaggregation of these numbers to investigate the BLN versus non-BLN populations for Mastery-level grades shows a more nuanced picture of success in our Calculus I courses. This suggests that more work needs to be pursued in order to achieve more equitable outcomes at the A and B grade level. Disaggregation of this data by gender and Master-Pass-No Pass categories also indicates that female students are achieving grades of A and B at noticeably higher rates than male students. This opportunity gap has brought to light the need to keep track of whether this opportunity gap is a trend or just a one-off. Averaging over all sections of a course can hide information about equity gaps within the population of students who pass the course. Potentially, it may also hide the true variance in student grades each semester.

Looking to the future, we plan to track these trends in subsequent years. Also, we plan to investigate the retention of students in STEM majors. The SEMINAL CSUEB team is collecting and evaluating data on the success of students who pass Calculus I and continue on in a STEM class. We will investigate the grade distributions of BLN and non-BLN students who passed Calculus I at CSUEB to see if the gains in closing the opportunity gap persist for this student population in subsequent semesters. Additionally, we will parse the data to investigate if BLN students are persisting as STEM majors at the same rate as non-BLN students.

# 22 A Focus on DEI for Students Yields DEI for Instructors

**Hortensia Soto**

*Colorado State University*

Hortensia Soto is a professor in the Department of Mathematics at Colorado State University. She brings expertise in facilitating and researching professional development for K-16 mathematics instructors. She also conducts research on the teaching and learning of undergraduate mathematics using an embodied cognition lens.

**Cynthia Wyels**

*California State University Channel Islands*

Cindy Wyels is a professor in the Mathematics Program at California State University Channel Islands. She is a co-author of two $6m HSI-STEM grants aimed at increasing the number of Hispanic and low-income students majoring in STEM fields; her research interests include using institutional data to analyze educational interventions.

**Andrea Alt**

*University of Northern Colorado*

Andrea Alt is a graduate student at the University of Northern Colorado pursuing a degree in Educational Mathematics. Her research interests include leveraging embodied cognition, the use of movement and action within an environment, in the mathematics classroom to enhance learning and capturing the gestures of online mathematics instructors.

**Ciera Street**

*Colorado State University*

Ciera Street is a graduate student in the Department of Mathematics at Colorado State University, where is working on her Ph.D. in Mathematics with a research focus in Mathematics Education. Her research interests include examining student and faculty experiences in post-secondary mathematics courses to investigate supportive and equitable structures within collegiate mathematics classrooms and departments.

**Selenne Bañuelos**

*California State University Channel Islands*

Selenne Bañuelos is an Associate Professor in the Mathematics Program at California State University Channel Islands (CI). She is the CI Campus Lead for the STEM-Service Courses (SSC) Initiative. She has experience in co-creating and facilitating professional development sessions tied to the themes of the SSC initiative.

*Abstract: In this case study we describe a professional development program for collegiate instructors from Hispanic-Serving institutions aimed to help them develop equitable teaching methods that support Hispanic and low-income students. In an effort to achieve these goals we created a professional development program focused on building classroom community, teaching with a student-centered lens, and creating and implementing rich mathematical tasks, all aimed at promoting equity in mathematics courses. Here we share the pleasant surprise of how the program impacted the participants, specifically how they experienced higher levels of inclusion as instructors at their school. We provide evidence of how some participants felt more knowledgeable about integrating evidence-based teaching methods, experienced a shift in their identity due to their new-found community, and felt empowered to make informed decisions*

*about how to teach in their classrooms. In summary, we hypothesize a model of how attending to equity and inclusivity of instructors can support DEI of students and vice versa.*

## 1. *Introduction*

California State University Channel Islands (CI) is a regional public university that embraces its Hispanic-Serving Institution (HSI) status. Nearly 67% of CI science, technology, engineering, and mathematics (STEM) undergraduates are Hispanic and/or low-income. Analysis of regional data uncovered race/ethnicity-based inequitable outcomes in graduation rates, retention in STEM majors, and success in the mathematics courses necessary for STEM majors at CI and its three primary feeder community colleges. The *STEM Service Courses (SSC)* Initiative[53] was created to address these inequitable student outcomes. This large professional development (PD) program, entitled Pathways with Regional Outreach and Mathematics Excellence for Student Achievement in STEM (PROMESAS), focused on building classroom community, teaching with a student-centered lens, and creating and implementing rich mathematical tasks, all aimed at promoting equitable teaching practices. We interpret equitable teaching practices as techniques that "mitigate systemic differences in ways that people experience and are afforded educational opportunities, particularly differences that privilege one group over another" (Adiredja & Andrews-Larson, 2017; Gutiérrez, 2009). The choice to focus on these themes stems from the work of Bressoud and Rasmussen (2015) who describe seven characteristics of successful calculus programs. We also leveraged the work of Estrada et al. (2018) who claim that classroom social inclusion such as integrating group work, engaging students in meaningful application tasks, and bolstering community in the classroom can promote equity. As such equitable teaching was a prescribed goal of the PD and the three themes were mechanisms for achieving this goal.

In an effort to achieve this goal we sought to improve Calculus success rates, STEM majors' retention rates, and students' perceptions of their experiences via a grant-funded 5-year initiative for mathematics faculty from CI and the three regional HSI community colleges through PD emphasizing equitable teaching techniques. We particularly targeted closing the "achievement gaps" -- the persistent 7 - 13 percentage point differences in the passing rates for Hispanic or low-income students in Calculus preparation, Calculus I, and Calculus II courses. Here we focus on how this PD actually improved equity for another structurally disadvantaged population: faculty with high teaching loads, specifically the adjunct instructors who participated in the PD. In this report we describe the program and the impact that it had on the PD participants, referred to as fellows, specifically how they experienced higher levels of inclusion as instructors at their schools.

## 2. *Description of Program*

Each of three cohorts of fellows began by reading and reflecting on articles about equity and stereotype threat (Aronson et al., 2002; Herzig, 2005), high-cognitive-demand tasks (i.e., rich

[53] Part of Project PROMESAS - Pathways with Regional Outreach and Mathematics Excellence for Student Achievement in STEM (US Department of Education, HSI-STEM, $6m, #P031C160017

tasks) (Smith & Stein, 1998; White & Mesa, 2014), sense of community (McMillan & Chavis, 1986), and student-centered learning (Laursen et al., 2014). These articles were selected because they offered glimpses into our themes and overarching goal. We asked the fellows to reflect on how they integrated or addressed these ideas in their own classrooms. Some articles such as the ones written by Herzig et al. and Laursen et al. were well received by the participants, but they struggled with other articles. This struggle was due to the technical and theoretical aspects of the articles. For example, the Aronson, et al. article is a quantitative research article, thus, we asked the participants to focus on the implications presented by the authors, but this did not seem to be enough for the participants to understand stereotype threat. The sense of community article was also challenging for some fellows because they viewed it as "fluff" and not relevant to teaching. On the other hand, some fellows seemed to identify with the *membership* component of a sense of community wherein one feels a sense of belonging. Although the fellows had two and a half months to read and reflect on the articles, several fellows waited until the last week to complete the assignment, which might have also added to the burden.

These readings, the three themes of community, rich tasks, and student-centered learning -- and their criticality for creating equity -- were then emphasized during an initial 5-day summer institute. During the academic year, the fellows engaged in monthly day-long workshops in which they continued exploring the three themes and equity, shared teaching materials, and discussed successes and challenges encountered as they worked to transform their beliefs about teaching. Each cohort culminated their first year with a 2-day institute wherein they developed action plans for intentionally transforming their teaching during the following year. The cohort then participated in another summer institute and a second year of monthly workshops, overlapping with a new cohort.

During the institutes and follow-up workshops, we unpacked how the research contributed to equitable teaching practices that mitigated systemic differences as defined by the theoretical work of Adiredja and Andrews-Larson (2017) and Gutiérrez (2009). We also shared practical techniques as described by Estrada, et al. (2018) who claim that *kindness cues* of social inclusion can promote equity by minimizing macro/micro aggressions which include blatant or subtle acts of discrimination and maximizing macro/micro affirmations which include blatant or subtle acts of social inclusion. Macro affirmations include positive facial expressions (e.g., smiling), caring, physical touch, and prosocial actions such as helping or politeness. Micro affirmations include subtler actions such as positive body language, tone of voice, and actions of vulnerability. These cues of social inclusion promote shared values within the classroom that "can disproportionately attract students from underrepresented groups without direct targeting or categorization" (Diekman et al., 2019, pp. 202-203) and can help students better connect to their academic community.

In all of our activities we modeled practical approaches that help mitigate systemic differences, e.g., using kindness cues of social inclusion to promote equity by maximizing macro/micro affirmations. For example, when someone presented their work to a task, we clapped. At first, the fellows thought this was strange, but some fellows soon recognized how such an affirmation of celebration can build community in the classroom. We also celebrated incorrect solutions to tasks because these were opportunities for learning. Furthermore, on each day of the institute and during monthly workshops we ensured that each fellow had an opportunity to present or speak on behalf of their group so that they could experience having a voice. We also shared our own stories of struggle when learning mathematics, an action of vulnerability, and encouraged the fellows to share their stories with their own students. Besides modeling such practices, we discussed why we did this and the possible implications for promoting equity in the classroom.

During discussions some fellows expressed concern about the amount of time that it would take to implement rich tasks and student-centered teaching. Some fellows also questioned whether students really learn concepts using these teaching techniques compared to lecture style teaching. Thus, we shared Freeman et al.'s (2014) findings indicating that students in a lecture class are 1.5 times more likely to fail than students enrolled in an active learning classroom, and how active learning benefits all students. Some fellows did not understand the significance of this finding until we mentioned Freeman et al.'s comment that if this had been a medical study, the study would have been halted. In an effort to provide the fellows with other strategies of active learning we presented them with the MAA *Instructional Practices Guide* (2018), which we integrated in several of our activities. The fellows also completed the mindset quiz, which assessed whether the fellows believed that intelligence was an inherent trait or something that can be developed. The fellows were not required to share their score – we simply wanted them to be aware of their own belief system and how it might impact their students.

Given that DEI conversations are often challenging, we routinely facilitated conversations exploring how rich tasks, student-centered learning, and a sense of community can promote DEI. For example, we discussed how student-centered learning could invoke enthusiasm and curiosity, address various learning preferences, and provide students with choices regarding the media forms for certain assignments. We also provided examples and asked the fellows to create rich tasks with a low threshold and a high ceiling, where students have opportunities to explore, communicate, argue, struggle, and tinker. Finally, we discussed how the culmination of implementing rich tasks via student-centered learning could create a sense of community in the classroom where students feel nurtured, a sense of belonging, influential, and a shared emotional connection (McMillan & Chavis, 1986). While the SSC initiative was designed to promote DEI for students, we were pleasantly surprised to find it did so for the fellows as well.

### 3. *Evidence of Impact towards DEI*

This summary is based on the first cohort's journal prompts and audio-taped interviews. The fellows, whose teaching experience ranged from 1-30 years, wrote journal reflections after each of the monthly PD workshops; they reflected on new readings, how to integrate what they learned, challenges with integrating new teaching strategies into their classroom, and their goals for transforming their teaching. Evidence also stems from 60 to 90-minute audio-recorded interviews with eight (four adjunct and four tenure-track) of the 14 fellows after they completed the first year of the program. During the interview the fellows shared their experience with the PD, in particular with each of the themes and their relationship to enhancing equity.

Adiredja and Andrews-Larson’s (2017) model for addressing equity in the mathematics classroom encompasses three components: *knowledge, identity, and power*. Our work suggests that this model extends to faculty, not just students. For example, the fellows described how the SSC initiative provided them knowledge on how to facilitate student-centered classrooms and to create rich tasks. It appeared that having access to evidence-based teaching techniques and assistance from a mathematics education expert helped the fellows gain pedagogical knowledge. For example, in his interview Saul, who had eight years of teaching experience, reported that, “Through active learning at *PROMESAS*. . . my teaching IQ [has] gone up.” In turn, the fellows were able to provide their students with opportunities to gain knowledge through a student-centered environment with an increased level of rigor using rich tasks. For example, Braden wrote, “enabling [the students] to discuss problems among themselves . . . is a powerful tool that allows them to learn [and] discover on their own.” The fellows came to believe that implementing rich tasks via student discovery gave students access to deeper and cognitively challenging mathematics. Estrada et al. (2018) claim that such a classroom structure can foster positive classroom relationships, which potentially communicates kindness cues and affirms social inclusion especially for Hispanic and low-income students who may feel inferior. Therefore, not only can students gain more mathematical knowledge, but they might also begin to identify as mathematics learners.

The second component of Adiredja and Andrews-Larson’s (2017) model for attending to equity in the mathematics classroom is the dynamic and situated nature of identity. Again, the researchers discuss this as important for equity within the classroom, but we saw evidence of this during our PD for the fellows. In the interviews and journals, the fellows described a shift in their own sense of identity and belonging in the *SSC* community. Before engaging in this project, some of the fellows described feeling isolated and without a voice. For example, in his interview, Adam shyly admitted feeling overwhelmed, underappreciated, and invisible prior to joining the PD. He said, “if you’re just an adjunct . . . driving between schools, you are like a ghost going into different schools, teaching, and leaving.” Adam added that he felt as though he could not vocalize complaints. It was his belief that “most of the teachers that are teaching . . . these days are not complaining . . . [because] no one’s listening.” Other fellows also shared their frustrations with their part-time instructor roles and limited camaraderie with colleagues. In her interview Megan said, “as a part-time lecturer, you don't get that sense of camaraderie or community with your colleagues as much as you would if you were a full-time person. And I feel like *PROMESAS* gives [me] that community.” Matthias talked about immediately feeling welcomed into the community of fellows. He said, “[The *PROMESAS* team] created this environment that made it safe for all of us . . . I was able to talk to [the other fellows] as if I knew them for a long time.” Max, too, quickly felt integrated into the community of fellows, and expressed feeling more confident in himself when he collaborated with other educators.

Comments such as this align with Estrada et al.’s (2018) assertion that research training programs can become a gathering place where teachers work with and learn from each other. Battey et al. (2018) refer to such an environment as one that creates emotional space, such as when the fellows felt validated when they discussed emotions and experiences surrounding their teaching transformations. By acknowledging their teaching struggles and supporting others’ teaching efforts, the fellows challenged the notion that teaching intelligence is innate.

Instead, the fellows highlighted how perseverance through adverse teaching experiences can be a rewarding experience that builds teaching confidence, as suggested by Battey and colleagues (2018). The last component of Adiredja and Andrews-Larson's (2017) model for addressing equity in the mathematics classroom is the co-constructed balance of classroom power. The fellows credited the PD for feeling empowered and having a voice in the decision-making of their own teaching. Megan shared her new-found power when she remarked, "I'm much more okay with just trying something out and seeing if it works, and then, sometimes it doesn't, and I get much less upset about it than I used to." Similarly, Max expressed that as a result of the PD, he was able to align his teaching beliefs with his teaching practices. He indicated that before the PD and as an adjunct he felt pressure from administrators to teach the students a certain amount of material, which required teacher-centeredness. After the PD, the fellows appeared less threatened by the formal power structure of the educational institution because they could argue that the teaching techniques they were adopting are evidence-based. Furthermore, they were part of a community in which others supported their teaching choices.

The fellows reported that their own empowerment allowed them to empower their students by giving students more of a voice in the classroom where each student's method of completing a problem was viewed as valid and respected. For example, Miguel said, "these [rich tasks] exploit the students' familiarity with these past concepts and expand on these ideas. This enables them to take ownership over these [mathematical] ideas." The notion of student empowerment in a student-centered environment parallels Estrada et al.'s (2018) work with kindness cues affirming social inclusion in the mathematics classroom. Estrada et al. suggest that, in addition to positively impacting learning outcomes, problem-based learning "can foster positive classroom relationships, which potentially communicates kindness cues affirming social inclusion for all students" (p. 274). Similarly, Battey et al. (2018) attribute students' mathematical empowerment to mathematical ability and acknowledgement of student contributions, such as when a fellow named a concept after the student who articulated it.

## 4. *Implications & Recommendations*

Estrada et al. (2018) suggest that confidence in the social context of the mathematics classroom can increase student engagement, can promote a greater connection to the classroom community, and thus, can influence students' identities. This can emerge by developing or maintaining students' shared values with their institution, discipline, or field area, which the fellows attempted to do. For example, as part of his interview, Kyle shared, "I am bringing [examples from] . . . different fields just so there are things that kind of pique the interest of different students." Estrada et al. also argue that affirming cues of affiliation can increase equity and connect students to their academic community, which reduces social ambiguity and relates to greater persistence. Max's remark of, "[the students are] talking with each other . . . about problems, and that doesn't happen if you don't create that kind of sense of community" supports this claim. The fellows believed that their students' sense of identity or belonging in the mathematics classroom increased as a result of their own shift in identity combined with their new-found pedagogical knowledge.

We argue that successful PD programs designed to address equitable teaching practices need to integrate both practical and theoretical perspectives of equity (e.g., Adiredja & Andrews-Larson, 2017; Estrada et al., 2018) and that a sense of community is central to this mission as depicted in the figure below. It appears that as instructors become more knowledgeable about and empowered to integrate student-centered learning and rich tasks, they feel safe to integrate kindness cues which strengthens the community of their classroom. Similarly, as students experience kindness cues, the classroom community strengthens, creating a safe space to learn and to make mistakes, which can be empowering and transform students' mathematical identities. We believe that increasing DEI for our students requires first addressing DEI for our faculty, especially for non-tenure track instructors who often teach STEM-pipeline courses such as Calculus. Our work suggests that becoming part of a community can empower such instructors to transform their teaching.

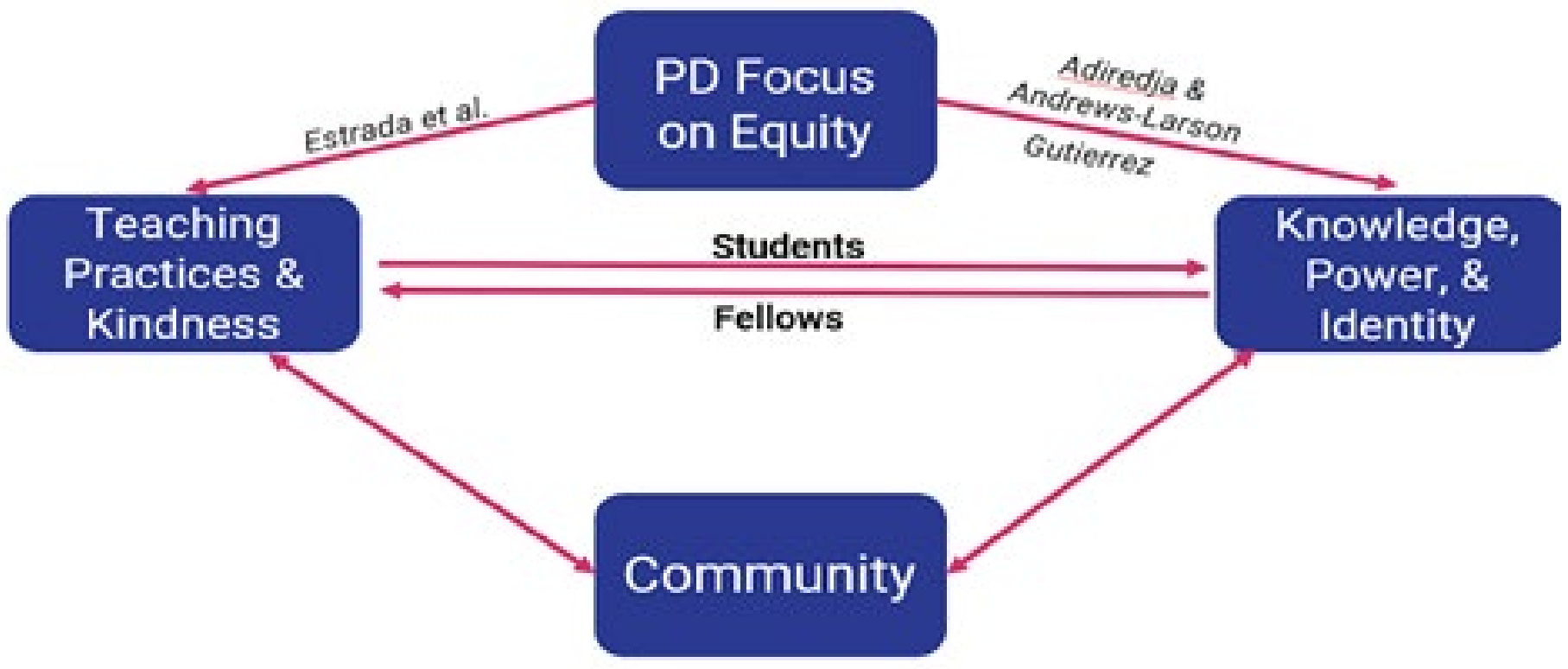

*Fig. 1 Interplay between Practical and Theoretical Perspectives on Equity*

Obstacles that our fellows encountered relative to better addressing DEI in their own classes included accepting that DEI is everyone's responsibility, insecurity in integrating new teaching techniques, push-back from students and colleagues, and a misconception that addressing DEI means not challenging students cognitively. These obstacles can be mitigated by creating a sense of community within departments and in classrooms. As of this writing, three cohorts have completed at least a year of the PD. Indeed, at each of our three community colleges we see that a critical mass of fellows have successfully shifted department cultures to embrace a sense of DEI as a collective responsibility. Support from chairs and deans and a concurrent initiative of the California Community Colleges to close gaps in retention and degree completion between students from marginalized populations (the largest contingent of which are Hispanic), low-income students, and white, non-low-income students have added valuable impetus. Since collecting these data, we have seen a shift in how administrators at these institutions embrace our work. We are currently collecting data from those administrators to better explain any systemic change that has occurred.

# 23 Undermining the "Bad at Math" Identity: Overhauling Differential Calculus to Empower Marginalized Students

**Emerald Stacy**
*Washington College*

Emerald Stacy is an Assistant Professor of Mathematics at Washington College, where she is also the faculty chair of the Diversity Committee, and the Founder and Director of the Justice, Equity, Diversity, and Inclusion (JEDI) Initiative. Founded in the summer of 2020, JEDI brings together faculty, staff, and students to create collaborative, inclusive social justice workshops that evolve to meet the needs of our community.

*Abstract: In the 2019-2020 school year, Washington College changed from a Pre-Calculus and Differential Calculus model to a year-long Differential Calculus (called Stretch Differential Calculus) with all as-needed algebra and trigonometry included. This paper will outline this new course and student responses.*

## 1. *Introduction*

Historically, Differential Calculus has had one of the highest D/F/W rates of all classes offered at Washington College. If a student wasn't ready for Differential Calculus, their only options were to take a one semester course in Pre-Calculus, attend a community college, or register for Differential Calculus anyway. Coming into the 2019-2020 school year, we made some changes. This article outlines the changes we made to our curriculum. Most notably, we stopped teaching Pre-Calculus, and now offer a year-long course in Differential Calculus with all the needed pre-requisite material included.

## 2. *The Problem*

Almost all students who take Differential Calculus at Washington College intend to major in STEM. Although students took a math placement test, the recommendations were not enforced so many students went straight into Differential Calculus unprepared.

Our Differential Calculus sections like many across the US were a blend of students who had passed the Advanced Placement exam, students that were prepared but had not yet tackled Calculus, and students that were still working on developing problem-solving skills and mathematical confidence. This mix created an atmosphere where many students, regardless of preparation, felt marginalized. Students of color were reporting that they did not feel that they belonged in STEM. Additionally, our institutional data shows that successfully passing Pre-Calculus did not indicate success in Differential Calculus.

In 2017-2018, the Provost requested ideas to help lower our D/F/W rate for Differential Calculus. In addition to a year-long sequence, we considered adding lab sections to Calculus. At the time, none of these ideas had unilateral support, or a champion to help push through any changes.

### 3. *Our Solution*

In April 2019, Dr. Eric Anthony Grollman, the founder and outgoing editor of *Conditionally Accepted*, a career advice column on Inside Higher Ed, came to speak at our school and led a workshop entitled *Becoming an Intellectual Activist*. Dr. Grollman told their own origin story, and discussed the barriers to activism within academia. The workshop was attended by roughly 35 staff, faculty, and students. I attended this workshop in my second semester at Washington College with the hope that I would leave with a clearer idea of how to be an activist at my own school while not endangering my career.

The second half of the workshop was devoted to time for open sharing, with crowd-sourced brainstorming for solutions. At this time, a black female-identified student reported that students of color did not feel at home in the science building or in math classes, and that they did not belong in STEM at Washington College. Hearing this directly from a student motivated me to try something – anything – to make our division more inclusive. I took this information back to my department and leveraged the student's words to convince resistant faculty that it was no longer acceptable to do nothing.

On an individual level, some of the math faculty were already doing standards-based grading or project-based courses, and because of the workshop, we wanted to make some program-wide changes to our department. Beginning in Fall 2019, we eliminated Pre-Calculus, began enforcing placement tests, and I designed a year-long course called Stretch Differential Calculus which was offered for any student who felt either academically or emotionally unprepared for Differential Calculus. A list of the topics covered in Stretch Differential Calculus are shown in Table 1. During the 2019-2020 academic year, I ran the course as a Special Topics course, and the course is now included in our course catalog as a standard offering.

**Table 1**

*Stretch Differential Calculus Course Topics*

| Fall Semester Topics | Spring Semester Topics |
|---|---|
| • Functions, polynomials, and lines<br>• The Derivative (interpretation, power rule, and tangent lines)<br>• Kinematics<br>• Exponential Functions and their Derivatives with Applications<br>• The Product Rule, Quotient Rule, and Chain Rule<br>• Characteristics of Graphs (increasing, decreasing, concavity)<br>• Extreme Values (local and absolute)<br>• Optimization | • Triangles, Trigonometric functions, and their inverses<br>• Derivatives of trigonometric functions & Optimization w/ Trig<br>• Circles<br>• Implicit Differentiation<br>• Applications of implicit differentiation (derivatives of logs, inverse functions, and related rates)<br>• 2-Dimentional Kinematics<br>• Infinite Limits and Limits at Infinity<br>• Summation Notation<br>• Definite Integrals<br>• Fundamental Theorem of Calculus |

## 4. *The Barriers*

One of the barriers to implementing this initiative was getting senior faculty buy-in. There was a pervasive belief that if a student wasn't ready for Differential Calculus, that they did not belong at Washington College. Also, some faculty were concerned that if a student didn't take Differential Calculus in their first semester that they wouldn't have time to finish a STEM degree. In practice, many students were failing Differential Calculus in their first semester of college, and had to retake it in Spring, or wait until the following Fall to take Pre-Calculus.

## 5. *The Course*

The teaching philosophy for this course is rooted in Feminist Pedagogy and Standpoint Epistemology. We are transparent with students about using evidence-based teaching practices that help minimize the opportunity gap, such as using an online homework system so that students get instant feedback (Serhan, 2019). We look at the data and outcomes of various educational research so that we will understand why the course is developed in the way that it is (Riegle-Crumb & Humpheries, 2012).

The class takes place in the beautiful and historical Rose O'Neill Literary House. There are couches, lots of places to work in small groups, a tea station, and a climate-controlled porch with a podium and projector. Since there is an entire year to stretch the course out for, there is a lot of time to adapt to the students in the room, slow down when we need to, and respond to how learning is happening. Additionally, we dedicate class time to discuss "how to college," and discuss the invisible syllabus of academia. Table 2 provides an outline of the course philosophy and other structures used to evaluate student understanding.

**Table 2**

*Course Philosophy and Formative Assessments*

| Course Philosophy | Online Homework | Written Homework | Student-Led Review Sessions | Synthesis Project |
|---|---|---|---|---|
| Emphasis on problem solving, communication, student agency, and creativity. | Short assignments that offer instant feedback and iteration. | Weekly written prompts for students to explore a topic. | Students prepare and lead a review session before an exam. | Students choose a topic from the semester to introduce with examples and tips for future students. |

### 5.1. *Day One*

To set the tone of the course and to promote group collaboration instead of diving into the syllabus, we spend the first day completing the Marshmallow Challenge. Students break up into teams of three, and each team is given twenty raw spaghetti noodles, a yard of string, a yard of tape, a pair of scissors, and a large marshmallow. Students are given 10-15 minutes to construct a freestanding structure from the items provided that can support the marshmallow. The winning team is the team that can construct the tallest structure to support the marshmallow. After the competition, we watch Tom Wujac's TED talk "Build a tower, build a team." The goal of the activity is to get students to work together on an activity that they don't think of as "doing math." We discuss the team dynamics, and each student writes a reflection after class to the following prompt:

> *Set a timer for 20 minutes, and write (in an informal style) about the problem-solving process from the Marshmallow Challenge. You might include: how you managed the power dynamic in the group, what you tried, what worked, what didn't work, or what you would do differently if we came back to class on Friday and did it again.*

I have found that students are able to dive into this activity and think critically about their group dynamics. This gives us a language to discuss group work through the semester, and what it means to be a problem-solver. The slower pace of the course allows us to spend time on activities that build connections between students and the mathematical concepts.

### 5.2. *A Typical Day*

On a typical day in a 50-minute class, we begin with a lecture, followed by an activity, and end with students presenting solutions to the activity.

*Lecture:* Class begins with us all sitting in a circle. I lecture for about ten minutes on an iPad with Notability that is projected on a screen. Before moving into the activity for the day, I ask if students are ready, or if they would like one more example. I try very hard to make sure students are driving the day, and that everything we do is to help them feel comfortable and ready to tackle a challenge.

*Activity:* Once they are ready, we break up into groups and they work on a worksheet to practice what I've lectured on for 20-30 minutes.

*Presentations:* Before class ends, we gather back together, and I put a pdf of the activity on the iPad and I sit down. Students take turns presenting their solutions on the iPad. When students make mistakes, other students offer support by way of cheering, asking leading questions, and reminding the presenter that "they can do it".

## 6. ***Math is Creative***

Near the beginning of the semester, I give a presentation arguing that math is a creative endeavor, and that we will approach problem solving as a creative practice. In the presentation, I share my own journey through college, how I made my living as an artist for ten years, and then ventured back into academia to go to graduate school. I argue that every time a problem is solved, the solution has been created – even if that solution has been previously created by someone else. Thus, a student could succeed in math based on their creative abilities, as opposed to "being a math person." My hope was that if I convinced students that this is a class that they can succeed in – even if they identify as "not good at math" – that they would then be willing to fight for success and support each other in the effort.

## 7. ***Written Homework***

In addition to online homework, students complete weekly written homework that is more holistic and conceptual, with prompts such as:

*Suppose you are given a polynomial function modeling the height of an object in terms of time. Describe how to determine how fast the object hits the ground. Include an example in your explanation.*

The written homework is to be typed to encourage them to treat the writing component as they would an essay. They are welcome to carefully write our equations by hand, take a picture and insert them into a word document, but their explanations should be typed. One of the learning outcomes for all Washington College students is that they develop verbal and written communication skills. In the math department, we encourage assignments and activities that help build those skills.

**8. *Review Sessions***

The class period before an exam is dedicated to student-led review. Each student is assigned a learning objective for the semester. When that objective is about to be assessed, they are responsible for leading the class through a review. How they choose to lead the review is completely up to them. Students might design a PowerPoint, offer worksheets with solutions, or develop creative games to play to practice. The goal of this activity is to give students an opportunity to practice synthesizing and communicating mathematical concepts.

**9. *Exams***

Each exam has four questions. Each question has multiple steps, but there is no scaffolding provided. For example:

*Suppose you are competing in a paper airplane competition. Now, these are no ordinary paper airplanes - they are trick airplanes and seem to defy gravity. The following function describes the height of your airplane in meters after t seconds:*

$$h(t) = -(t+3)(t+6)(t-4) = -t^3 - 5t^2 + 18t + 72.$$

*How fast is the airplane moving when it lands?*

Students are graded qualitatively, with each problem earning A, B, C, D, F, or 0. After the exam, students have a chance to correct their work to raise their grade up to one letter grade per problem. To have their grade raised, they must correctly redo the problem, and explain any errors they made originally. For example, if a student earns an F on each problem, but works very hard on their corrections, they could earn a D on the exam. Thus, a student would hypothetically earn an F on every exam problem and still pass the class through hard work. This is deliberate. We want to lower the stress of exams and encourage students to dig deep into the material.

**10. *The Synthesis Project***

Near the end of each semester, students choose a topic to create a guide for future students. Each guide should contain a complete overview of the topic, 2-5 examples, and tips and tricks. The guides are done in two drafts: the first gets peer reviewed as well as feedback from me, and then students submit their final draft. Although the product has varied significantly from student to student, students engaged with the project and created some quite excellent guides.

Students reported that they found the synthesis project to be a helpful assignment. The variety of topics and flexibility to work on the project independently or in groups was appreciated. One student commented

> *"The synthesis activity helped reduce stress for me as it did not have the formal label of an exam. The activity allowed me to choose a topic to teach other students which was extremely helpful when I was listening to other peer's projects. Having the concepts explained from a different perspective, in a different way definitely gave me a better grasp of the topics we needed to master."*

### 11. *Impact*

Although Washington College is a predominantly white institution, roughly half of the students that enrolled in this course in 2019-2020 are black women, and most of the students are first-generation college students. As seen in Table 3, the D/F/W rates for Differential Calculus drastically dropped with the new model. For context, I taught two of five sections of Differential Calculus in Fall 2018, and all sections of Stretch Calculus in Fall 2019.

**Table 3**

*D/F/W Rates for Differential Calculus Courses at Washington College, 2017 – 2021*

| | Differential Calculus D/F/W Rate | Stretch Calculus D/F/W Rate |
|---|---|---|
| Fall 2017 (n = 89) | 21.3% | N/A |
| Fall 2018 (n = 63) | 22.2% | N/A |
| Fall 2019 (n = 59) | 9.3% | 7.2% |
| Fall 2020 (n=69) | 7.3% | 7.1% |
| Fall 2021 (n=61) | 10.3% | 0% |

Not only was the D/F/W rate for Stretch Calculus 7.2% in the Fall of 2019, the D/F/W rate for Differential Calculus dropped as well. Using $p_0 = .213$, and $n = 59$, we have $p = 0.008061$ and thus the difference in D/F/W rates with the introduction of Stretch Differential Calculus is statistically significant.

Grades certainly don't tell the whole story, but the narrative comments in the Student Evaluations of Teaching have been supportive of the course. Students claimed that the class format was welcoming, and that even though the material was hard, the course was *"taught in a very friendly and open manner, which made it a lot easier to make mistakes and ask questions without feeling insecure."*

As of May 2021, we have run the sequence twice, and there is unilateral departmental and divisional support for the change. The course has now been added to the Course Catalog and approved by a campus-wide faculty vote.

# 24 The Effective Thinking Calculus Project at UT Austin

**Michael Starbird**

*The University of Texas at Austin*

Michael Starbird is a Distinguished Teaching Professor of Mathematics at The University of Texas at Austin with many teaching awards including the MAA's Haimo Award. He has written nine books and produced DVD courses in the *Great Courses Series* on calculus, statistics, probability, geometry, and the joy of thinking, and a MOOC entitled Effective Thinking Through Mathematics. Two books co-authored with Edward Burger are *The Heart of Mathematics: An invitation to effective thinking* and *The 5 Elements of Effective Thinking*, which has appeared in 17 foreign language translations.

**William R. Wolesensky**

*The University of Texas at Austin*

William Wolesensky started at The University of Texas at Austin in the fall of 2012. He received his Ph.D. in mathematics in 2002 from the University of Nebraska in Lincoln. While at the University of Texas, he has taught a wide range of undergraduate math courses and helped develop innovative first-year calculus curriculum and reform. His teaching excellence has been recognized most recently by a UT System 2020 Regents Outstanding Teaching Award.

**toufic chahin**

*The University of Texas at Austin*

toufic chahin is a second-year doctoral student in Education Policy and Planning at The University of Texas at Austin, having recently earned a Masters in STEM Education. He is researching the discourse of education policies that impact underrepresented minority students. Seven years as a high school math teacher informs his motivation to work towards equity and inclusion for all students. A Texas native, his family history traces its roots to the Mexican border and includes educators going back four generations.

**Bethany Miller**

*The University of Texas at Austin*

Bethany Miller is a graduate student at The University of Texas at Austin pursuing a doctorate in STEM Education. Her experiences working with K-12 students within minoritized communities has contributed to her interest in creating equitable spaces within the mathematics classroom.

**Joel Walsh**

*The University of Texas at Austin*

Joel Walsh is a fifth year Ph.D. student in STEM Education at The University of Texas at Austin and first-generation college student. He grew up attending school all over the U.S. and in his mother's home country of Korea, a study in contrast that led him to eventually teach high school mathematics in public schools in Los Angeles, California.

*Abstract: The basic premise of the Effective Thinking Calculus (ETC) Project is to view the mathematical challenges all students face with calculus as reflections of more general issues—essentially, issues about thinking strategies and issues around a sense of belonging. The Project consists of two bundled courses populated with exactly the same diverse set of Discovery Scholars students and having the same instructional team. The two bundled classes are a Signature Course titled "Elements of Effective Thinking" and a first semester calculus course. The Effective Thinking course (as well as the calculus course) encourages students to develop*

*powerful practices of mind through their own agency including: understanding fundamentals deeply, striving for meaningful understanding, learning from mistakes, raising questions, and being comfortable with change. The Effective Thinking course activities and assignments help students to apply those practices of mind to learning calculus, learning other subjects, and succeeding in life. The associated calculus course employs inquiry-based learning methods of instruction that directly connect the practices of effective thinking to the learning of calculus. Students discover mathematical concepts by working on guided challenges, often in groups. One of the goals of this method of instruction is to transform students from being purely consumers of knowledge to seeing themselves as producers of knowledge. The experiences of both courses involve productive struggle, learning from mistakes, and explaining ideas to classmates. The ETC Project is certainly still a work in progress, but the fundamental concept of overtly teaching effective thinking practices of mind and applying those systematically to a specific course (in this case, calculus) definitely shows promise.*

## 1. *Introduction*

A Texas state law requires The University of Texas at Austin (UT) to automatically accept any Texas high school graduates who were in the top 6% of their graduating class. This legislation was an attempt to have the University of Texas student body reflect the racial and economic diversity of the state. Texas high schools vary in quality and opportunities they offer for learning mathematics, so students arrive with a wide range in the quality of their preparation. To respond to this reality, UT has created several success programs designed to assist students identified as having a low probability of graduating in four years. One such program is called the Discovery Scholars Program. Historically, Discovery Scholars students' success rate in calculus courses was very poor, thereby contributing significantly to their lack of success in college in general. So, we created a new approach for these students called the Effective Thinking Calculus (ETC) Project.

The basic premise of the ETC Project is to view the mathematical challenges all students face with calculus as reflections of more general issues—essentially, issues about thinking strategies and issues around a sense of belonging. The Discovery Scholars Program has a purposefully diverse enrollment including many students from groups that have historically been excluded from the opportunity to produce mathematics. We mindfully and ardently take the position that all our students can understand mathematics meaningfully, can produce mathematics on their own, and can view mathematics as a part of who they are. The students discover what the instructional team all know, namely, that our students are strongly capable of doing and understanding mathematics meaningfully. In fact, we repeatedly, overtly encourage them to consider becoming mathematics majors--of course, acknowledging that most will choose other majors. Creating shared experiences of grappling with intellectual challenges while being supported and encouraged has allowed our students to succeed at a much higher rate than historically similar cohorts have done.

Some background: UT requires every first-year student to take a course designated as a Signature Course. All Signature Courses are taught by senior faculty members and are designed to give UT students a shared intellectual experience that includes college-level writing, speaking, reasoning, and problem-solving—typically using interdisciplinary subject matter. We took advantage of UT's Signature Course program to create the ETC Project. However, there is no necessity for having a university-wide requirement to replicate the structure of the ETC Project at other colleges and universities.

### 2. *Structure of the Effective Thinking Calculus Project*

The Effective Thinking Calculus Project consists of two bundled courses populated with the same set of Discovery Scholars students and having the same instructional team. The two bundled classes are a Signature Course titled "Elements of Effective Thinking" and a first semester calculus course.

The Effective Thinking course (as well as the calculus course) encourages students to develop powerful practices of mind through their own agency including: understanding fundamentals deeply, striving for meaningful understanding, learning from mistakes, raising questions, and being comfortable with change. The Effective Thinking course activities and assignments help students to apply those practices of mind to learning calculus, learning other subjects, and succeeding in life.

The associated calculus course employs inquiry-based learning methods of instruction that directly connect the practices of effective thinking to the learning of calculus. Students discover mathematical concepts by working on guided challenges, often in groups. One of the goals of this method of instruction is to transform students from being purely consumers of knowledge to seeing themselves as producers of knowledge. The experiences of both courses involve productive struggle, learning from mistakes, and explaining ideas to classmates.

### 3. *The 'Effective Thinking' Signature Course*

Some of the Effective Thinking Course assignments and activities overtly concern mathematics and calculus, but the bulk of the activities and assignments only tangentially touch upon calculus. *The 5 Elements of Effective Thinking* book (Burger & Starbird, 2012) provides a perspective on active understanding that underlies all the experiences and assignments throughout the semester. The structure and assignments are designed to support the ETC Project's ambition to create a community of students who are improving their thinking and learning in an uplifting atmosphere.

The overtly mathematical features of the Effective Thinking course include proofs of the Pythagorean Theorem, proofs of formulas for the areas of circles and other geometric figures, a geometric representation of the difference of squares formula, and other mathematical gems. These help to give students the idea that mathematics is meaningful rather than being a set of facts and procedures to memorize. Students are often amazed that there are reasons behind formulas.

During the semester, several 'Tell the Story of Concept Creation' exercises lead students to write papers explaining why the core ideas of calculus are natural and meaningful. For example, we ask students to use a car moving on a road to discover why the definition of the derivative is the natural consequence of this analysis. These assignments help students explore the meaning and definitions of the derivative, the definite integral, and the Fundamental Theorem of Calculus.

Other writing and speaking assignments are not overtly mathematical. One recurring theme involves puzzles. Students are asked to create strategies for approaching puzzles rather than focusing on the solutions. These specially designed puzzles do not refer to background knowledge, so students work together to create methods of problem-solving. The puzzles help create a general framework for understanding simple things deeply—a practice we regularly encourage students to develop throughout the course.

Other themes of the course include the scientific method, evaluating and using good evidence, data-based descriptions of social change, experiencing and describing art, and personal reflections on their own study practices and learning. In addition, the course participates in a UT wellness initiative called SHIFT, which strives "to change the culture around substance use from one of misuse to one of well-being." All these experiences are designed to create a mindset that embraces personal transformation, overtly manifested in their final project, which involves both a written component and an oral presentation. In that 'Change Your Mind' project, our students select a topic of personal significance to them and explore it by overtly applying practices of mind they have acquired to develop deeper and more nuanced insights into their topic. Students are encouraged to recognize how their relationship to their topic has shifted through deeper thinking about it.

Some of the writing and speaking assignments require students to use UT resources, specifically the University's Writing Center and the Public Speaking Center, to get third party constructive criticism. This requirement has the hidden agenda of encouraging students to realize that UT and the broader society have available resources specifically designed to help them succeed and that taking advantage of opportunities is part of successful living.

In addition to the academic content students engage with, the Effective Thinking course strives to create a lasting, inclusive, supportive community. Small weekly discussion sections led by TAs help to form cohesive relationships while helping students grapple with the academic content of the Effective Thinking course. In addition, ETC students have a support structure in the form of First-Year Interest (FIG) groups, designed to help students adjust to college life and expectations. Students in a FIG take many of their introductory courses together (our two classes and a chemistry course or an economics course, for example) and attend weekly meetings that are run by a peer mentor and a university staff member. In those meetings, students discuss how things are going, hear from visitors about life experience topics, and generally have someone attentive to the students' condition.

To encourage students to view their education as extending beyond the classroom, every year the professor invites all our students to his house for one large celebration where students can informally interact with their peers as well their professor, teaching assistants, and the professor's dogs. Here, the students share their interests, backgrounds, and experiences while connecting with others.

## 4. *The Calculus Course*

The Calculus course part of the ETC Project is a standard first semester course in differential calculus: limits, derivatives and applications of derivatives, early bits of integration, and the fundamental theorem of calculus. At UT this class is taught via lectures and biweekly TA sessions. In our Effective Thinking Calculus section, group work, inquiry, teaching to learn, and a reverse-flipped component are central features of instruction.

A typical class session consists of a topic being briefly introduced by an instructor followed by challenge questions posed to reinforce and augment the concept or technique. Students work in small groups with those seated near them while the instructional team and undergraduate learning assistants (who are previous year ETC students) go around the room helping. The class regularly reconvenes to see resolutions—often student-generated answers. This basic format is augmented with several wrinkles.

One wrinkle is the 'reverse-flipped' assignment. After almost every class meeting and due before midnight of that same day, the students are required to watch short online videos and take the corresponding short quizzes on the topic that was explored that day in class. The videos were created by various professors at UT to be used in flipped classes—that is, the videos were designed to be watched before class; however, we prefer to introduce topics in an interactive, inquiry format during class and then later have the students watch these direct instruction videos to reinforce and solidify their understanding.

Another wrinkle is the "Problem of the week." With perhaps a half hour left in a class period, the whole class is given a challenging problem, generally focused on an exploration or application. As usual, the students work on it in small groups. The twist? At the end of the hour, every student in the entire class receives the same grade on the problem—namely, the score equal to the percentage of students who submit the correct answer. Thus, students are encouraged and motivated to share ideas and convince each other of the correct approach. Periodic polls are taken and shared during the class to give a pulse on their classmates' thoughts. Certain weeks, such as the surprisingly contentious optimization week and limits week, were models of engagement and collaboration. The instructional team were often relegated to the sidelines as students covered all available boards while attempting to prove their arguments and persuade their peers. Of course, it is essential to have established a spirit of mutual support and cheerfulness in the class in order to have the students joyfully engage with one another. By the end of each of these hours, the students had generally reached a near consensus and found the experience uplifting.

"Calcsplanation" is our 'teach-to-learn' component of instruction. Each of the weekly Calcsplanation assignments consists of a prompt and outline of the calculus topic being treated that week. Students then create short videos that explain the topic and share the videos with their fellow students in small groups with facilitators managing the discussion by asking probing questions as if they themselves were students. The facilitators are ETC Project students from the previous year.

## 5. *Testing*

The instructional team share a common philosophy towards testing; namely that the principal aims of assessment are to let students demonstrate their understanding and to let students identify any misunderstandings so they can subsequently clarify their knowledge. One practice that embodies these beliefs is our approach to our students' common Pre-Calculus gaps. Within the first week of class, students take an anonymous diagnostic quiz covering a wide range of foundational Pre-Calculus knowledge and skills. Students then complete practice exercises corresponding to each missed question, after which they retake a different version of the test. Earning a 100% on the test is a mandatory assignment and often requires several attempts, but it is completed during the first few weeks of class.

The second somewhat novel assessment practice involves group tests. We have the advantage of teaching the two courses in back-to-back hours in the same room. So, for each midterm test, the first hour is used for an individual test. Then during the second hour, the students are placed randomly in groups of three to take another test. The group portion of the test contains questions that are either identical to or more rigorous variations of those found on the individual portion of the test. All members of each group receive the same grade, which is worth about 30% - 40% of a student's overall midterm exam grade. Our goal is to provide students with an opportunity to explain and discuss concepts immediately after they have worked as hard as they can on questions about those concepts. This group portion lets students gain confidence in their understanding or remediate misunderstandings immediately. We view those concentrated discussions during the group portion of the tests as among the most effective learning moments that the students experience.

## 6. *Impact of the ETC Project*

The ETC Project has run for five fall semesters with cohorts of 46, 81, 56, 37, and 46 respectively. The ETC Project is certainly still a work in progress; however, the fundamental structure seems promising. Overtly teaching effective thinking practices of mind and applying those systematically to a specific course (in this case calculus) seems helpful. At the same time, building an inclusive community in the context of meaningful learning creates a mutually supportive environment in which all students experience themselves as participants in the common journey of self-creation well beyond mathematics.

Assessment experts from the Dean's Office in the College of Natural Sciences at UT compared our Effective Thinking Calculus Project students with students of similar backgrounds from previous years. The data indicate that the students who participated in the ETC Project (the Treatment group) fared better than Discovery Scholars students who had attempted calculus in previous years (the Control group). The differences in some metrics were strikingly beneficial for students in the ETC Project—for example, the D/F rate, drop rate, and percentage who took a second semester of calculus.

A second study compared ETC Project students in a follow-on second semester calculus section with students in the same section who had taken other sections of first semester calculus. By the end of the semester, the ETC students' grade distribution was indistinguishable from that of the non-ETC students, suggesting that the treatment had been successful at leveling the playing field between the intervention population and the general student population. This equality occurred even though the UT math placement exam scores for some of the ETC students had been below the cutoff for placement into a standard differential calculus course.

We attempted one intervention that completely failed. In the second and third years of the project, we offered students the option of dropping back into a pre-calculus course that started one month into the semester for students who were struggling. At that point one of the midterm tests had been given and seven students one year and five the other year had done very badly on it—scoring from 25-55 points on a standard grading scale. We individually advised them to consider switching into the pre-calculus course. One year, not a single student chose that option, and the other year one student did switch. The students were simply not willing to leave our class. All those students succeeded in passing with a C or better.

## 7. *Conclusion*

The success of the ETC Project has been impressive; however, the students in the ETC Project had several advantages over the Control group beyond the structure of the classes. A significant factor was that every member of our instructional team is dedicated to helping every single student to succeed. None of us believe that a low D/F/drop rate is good enough—we believe every one of our students can succeed at first semester calculus. So, when a student was struggling, we knew about it and acted. The active classrooms with frequent small group discussions, the weekly TA-led sessions with student presentations and interactions, the common challenge of learning calculus, the First-Year Interest Groups, and the informal gatherings all contribute to creating a mutually supportive community and a strong sense of belonging to the class and to the university.

The Effective Thinking Calculus Project continues. Our immediate goal is to have students with varied previous mathematical opportunities learn first semester calculus meaningfully. But in addition, we want to raise students' standards of understanding and apply that higher level of understanding to everything they do—in school and beyond. We believe we have improved the experience each year, and we look forward to further improvement in the years ahead. Our aspiration is to transform students to adopt practices of mind that allow them to take joy in lifelong learning, to learn meaningfully in mathematics and in every other subject, and to become constructive, thoughtful, and reasoning citizens of the world.

# 25 Leveraging Classroom Data to Promote Equitable Instruction

**Amelia Stone-Johnstone**

*California State University, Fullerton*

Dr. Amelia Stone-Johnstone is an Assistant Professor in the Department of Mathematics at California State University, Fullerton. Her research interests include academic support structures for students in gateway mathematics courses, faculty and instructor professional development, and equity-minded and active instructional practices in undergraduate mathematics education.

**Robin Wilson**

*Loyola Marymount University*

Dr. Robin Wilson is a Professor of Mathematics at Loyola Marymount University. His current research interests include both low-dimensional topology and issues of equity and access in mathematics education.

**Daniel Reinholz**

*San Diego State University*

Dr. Daniel L. Reinholz is an Associate Professor in the Department of Mathematics and Statistics at San Diego State University and is also a Provost's Professor of Equity in Education. Their research focuses on disability justice, faculty development, and systemic change in STEM. Dr. Reinholz is a co-developer of the EQUIP tool (https://equip.ninja) for capturing and disrupting patterns of racial and gender inequity in classroom practice. Dr. Reinholz is also a co-founder of Sines of Disability (www.sinesofdisability.com) which aims to disrupt ableism in mathematics and beyond.

*Abstract: This chapter describes how classroom data can be leveraged to address inequities within calculus classrooms. Calculus classrooms are an especially important site for understanding equity because they often serve as a gatekeeper to higher-level mathematics and STEM careers. We discuss how the authors worked collaboratively to support equitable teaching in one of the author's classrooms (Wilson). We used the classroom observation tool EQUIP (Reinholz & Shah, 2018), to provide Wilson with data on patterns of student participation. We describe how these data were used to help create a better learning experience for students in the class who were historically marginalized in mathematics. We provide both classroom data and our own reflection on the process to help generalize our methodology for others who seek to address inequities in introductory college-level mathematics courses.*

## 1. *Background*

Participation in classroom discourse is an essential part of learning (e.g., Banes et al., 2019), and it can help shape students' attitudes and beliefs about mathematics (Herbel-Eisenmann et al., 2011). Research shows that opportunities to engage discursively are not widely distributed to all students. There are race-based and gender-based inequities in *who* gets opportunities to participate and *how* opportunities are distributed within the classroom (e.g., Ernest et al., 2019; McAfee, 2014; Sadker et al., 2009; Smith et al., 2019). These inequities can be attributed, at least in part, to an instructor's implicit biases (Ball, 2018). These implicit biases are expectations and evaluations about particular groups of people that operate largely outside of conscious control (Staats et al., 2017). In mathematics, problematic racial and gender stereotypes create a status quo that portray mathematics as largely masculine and White/Asian (Shah, 2017). Given these stereotypes, instructors are likely to inadvertently favor these groups of students, for instance, in who they call on to participate in a discussion. These inequities can arise even when instructors have explicit equity goals.

## 2. *Method*

### 2.1. *Positionality*

The lead author identifies as a Black woman. She co-facilitated the professional learning community (PLC) that is described in this chapter and has extensive experience with using data to support faculty professional development towards equitable teaching. The second author identifies as a Black man. He was a participant in the PLC and has an expressed interest and commitment to fostering equitable learning environments in his classroom. The final author identifies as a white man with multiple disabilities. He is one of the developers of the EQUIP observational tool, and was the lead facilitator for the PLC. He takes a disability justice lens to issues of racial and gender justice.

### 2.2. *Data Source and Analytic Method*

The current data come from a larger study in which three mathematics instructors volunteered and participated in an unpaid PLC that used data analytics to support instructional change. At the start of the PLC, each instructor met with the lead facilitator to discuss their equity-related instructional goals for their Calculus classroom. The PLC was organized around a series of reflection cycles where participants had their classrooms recorded and coded by the facilitators, and then participated in a debrief meeting about data analytics from their classroom observation. As a member of the community, Wilson was observed teaching calculus four times during a single semester and received data analytics from EQUIP after each observation. Within the context of the learning community, Wilson reflected on the data he received in collaboration with his peers and the facilitators (Reinholz and Stone-Johnstone), and set action plans to change his instruction.

The data analytics were generated from EQUIP (Reinholz & Shah, 2018). EQUIP is a free, customizable web app (https://www.equip.ninja) that supports the coding of classroom participation and automatically generates analytics to support instructor reflection. EQUIP focuses on participation at the level of individual students, and when combined with demographic information about students, it can generate analytics about groups. Thus, EQUIP analytics can answer questions such as: What percentage of the questions were asked to women in the class? Or what proportion of high-level questions went to Black students? The answers to these questions provide a starting point for deeper conversations about racial and gender equity in the classroom.

For this study, student demographics were collected using an open-ended student survey where students were able to self-identify their racial and gender identity. There were 36 of students enrolled in Wilson's Calculus class where in terms of gender, they self-identified as either man, woman, or non-binary. Students identified their race as either Asian, Black, Filipino, Latinx, Middle Eastern/North African (MENA), Mixed-White, or White. In this chapter, we reflect on the learning process that happened during the semester in which the analytics were available and also, we reflect from the perspective of one year later, as Wilson continued to change his practice even without access to further analytics. Here we provide Wilson's reflections, told in the first person, to tell the story from his perspective.

## 3. *Wilson's Reflections*

When I was first contacted about the opportunity to use participation data from my classroom to help make informed pedagogical choices I was intrigued. As an educator, I am constantly working and struggling to find ways to keep the students at the center of my classroom, to provide them opportunities for deep inquiry around the mathematics we are learning, and to "include all students in rigorous mathematical learning and mathematical identity building" (Laursen & Rasmussen, 2019). In many ways, I feel well prepared to effectively teach students from a wide variety of backgrounds and varying ability levels. Each day I set out with the goal of including "all" students in meaningful yet rigorous mathematics activities, however without data on who is participating, student participation is no doubt biased by my own blind spots which can create false perceptions about what is happening in the classroom. The data that was collected and shared with me throughout the PLC helped me see my blind spots in terms of both *who* was participating and *how* students were participating from class to class. Further, having the data broken down by gender and ethnicity provided an even clearer focus on what students I was and was not serving.

### 3.1. *Learning About Who was Participating*

Upon entering into the EQUIP study as a participant I expected to see participation relatively evenly distributed throughout the class. On the one hand, I thought the data we were going to look at would say more about my students than about myself as a teacher. Like maybe they were just shy, or maybe they didn't speak English as their first language, or maybe they work a lot outside of class and have less time to put in therefore had less to contribute. On the other hand, I was also nervous to see what biases of my own might show up in the data. These were the things on my mind upon entering into the EQUIP study.

What I was surprised to find, upon reflecting on the data with my colleagues in the PLC, was that my classroom data said much more about me and my own teaching than it did about any of my students and what they themselves were bringing into the classroom. The project helped me recognize the power that I had around the nature of students' participation in my class. In addition, what I saw in the data about the classroom interactions after each recorded class session did not always align with my own memory of what happened in the classroom. For instance, I was very surprised to learn that only 17 out of 36 students participated during the first round of classroom observations (see Figure 1). I was convinced that I had engaged every student in the class, and that more students participated than the data revealed.

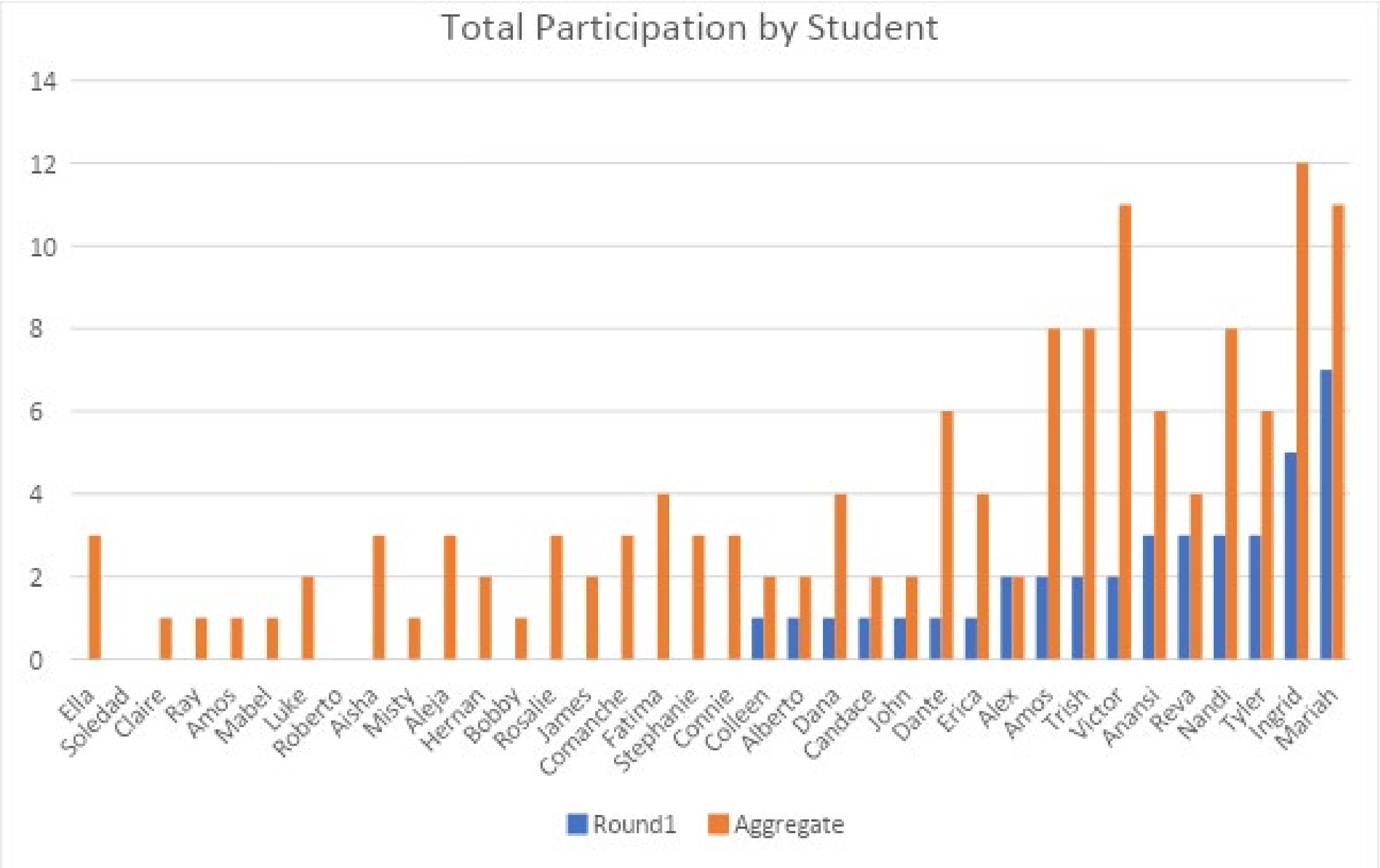

*Figure 1. Distribution of Student Participation.*

The data revealed to me that the three Filipino men, that I was fairly certain participated, actually did not contribute to the classroom discussion. Also, there were students (like one Latinx woman) who I did not think were contributing more than anyone else, who turned out to be participating on a regular basis in an unbalanced way. I'm not sure what these anecdotes say about my own biases, but my takeaway was that bias is complex, and unpacking one's own bias is a very non-trivial task. In fact, the task is so complex that I don't believe it's possible for me to ever be completely bias-free, and I've turned to other strategies to facilitate classroom participation more equitably. While there are lots of ideas out there about how to do this, the one that I have settled on as most effective for me is the use of index cards to call on students in a more fairly distributed way. In fact, by the end of the study I was able to engage about 94% of the students in whole class discussions over the entire observation period (see Aggregate in Figure 1).

### 3.2. *Learning About How Students were Participating*

In addition to the data about who participated in the course, another type of data that EQUIP provided was related to how students participated. There were times when I was sure that my students engaged in class discussions at a high-level (e.g., providing long contributions, and explaining their reasoning), but the data (see Figure 2) revealed that while a large percentage of students participated, they mostly contributed short answers (e.g., a few words, or a single sentence). This data pushed me to adopt strategies to encourage deeper engagement, such as providing think-time for students through think-pair-share activities, increasing wait time, and carefully monitoring student work during collaborative activities; and selecting students to share their multiple solution paths, including false-starts and small mistakes.

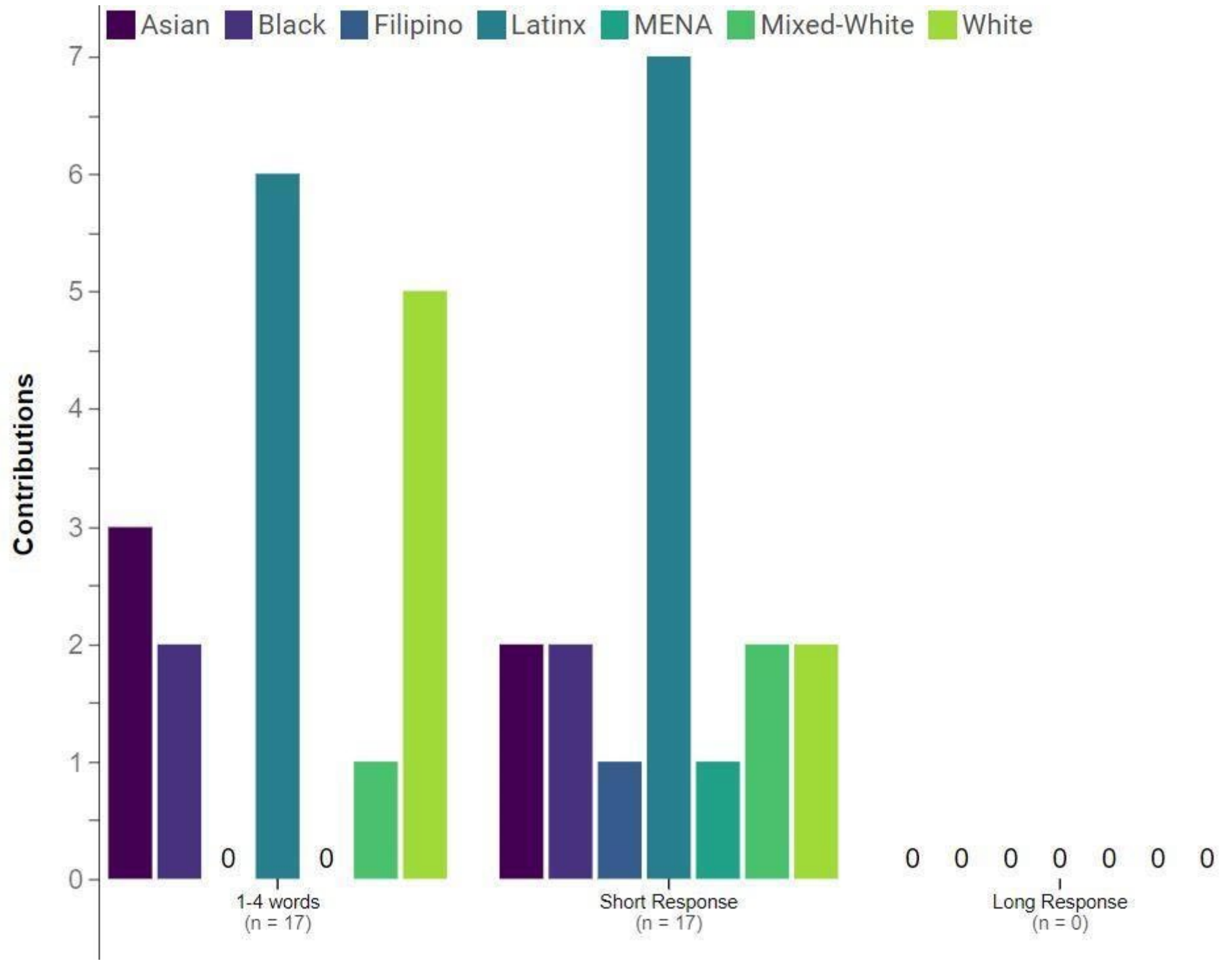

*Figure 2. Distribution of Talk Length by Race.*

As a result of this experience, I am no longer satisfied with short, one-word answers and it's become a habit to probe deeper into students' understanding. I do this by creating a classroom community where sociomathematical norms for inquiry are established and reinforced by me asking follow-up questions, choosing questions that provoke inquiry, and by promoting a culture where all contributions are valued. This can be seen in the following example:

| | |
|---|---|
| Wilson: | *Does anyone … someone I haven't heard from today that wants to explain how they came up with C? Anybody? Go ahead Tyler.* |
| Tyler: | *Umm so basically … the one thing that told me it was C of all else was the slope at x=0.* |
| Wilson: | *Ahhh!* |
| Tyler: | *It was negative which means the derivative for the 0 or the y-intercept would need to be a negative …* |
| Wilson: | *Hold up one second guys, because Tyler is giving us some real gems. So, he said that he was struggling at first, and then he realized that the slope at 0 … is it positive or negative?* |
| Various students: | *Negative.* |
| Wilson: | *That means that at 0 the derivative graph should be negative. And only one of them fits that, is that how you nailed it down?* |
| Tyler: | *Yep.* |

In this example, we were looking at various potential graphs of the derivative of a given function. It was important for me to have students explain why they chose the answer that they chose. I also found it important for the other students to recognize that Tyler did not immediately know the answer. The steps he took to select C were just as important as the answer itself.

An important lesson for me has been that getting all students to contribute in a meaningful way does not require all students to have the right answers. Meaningful student contributions take the form of “correct” answers, “incorrect” answers, and their own questions. And each of these types of contributions helps lead to a stronger understanding of the mathematics, not just for the students contributing, but for their classmates as well. As Rochelle Gutiérrez (2018) suggested, students do not have misconceptions, they have conceptions. It is up to us as their teachers to guide their conceptions toward a deeper understanding of the mathematics at hand. This experience has pushed me to think much harder about every student contribution, especially the “wrong answers”, and try to understand how the student is thinking about the mathematics so that we can all learn from their conceptions. In fact, the incorrect answers and clarifying questions often do more to push the mathematical agenda of the class forward than I or any short (but correct) answer can. For example, I find that the questions that my students have that they will often phrase as “stupid” questions, are likely representative of the understanding of a larger sample of students in the room. By making their thinking public I am given insight that can help not only drive that moment of instruction, but can also inform decision-making in the future when that particular topic arises. Therefore, it becomes critical to hear as many voices as possible in the classroom to increase the opportunities to highlight student thinking and to maximize opportunities for mathematical learning and identity building. At the same time this formative assessment allows teachers to learn more about their students as mathematical doers and thinkers.

### 4. *In summary*

My experience with EQUIP helped change the focus of my teaching from thinking about what I was going to do next mathematically, to focusing on who was going to participate next in the classroom and how to build from their contribution. The power dynamics around student participation in whole class discussions is very complex, and student identities influence their participation in this setting. At the same time these interactions (or lack thereof) are also critical in the development of students’ evolving math identities (Oppland-Cordell & Martin, 2015). The EQUIP data and dialogues surrounding it have supported my efforts to foster student learning and build positive mathematics identities. My participation in this project has helped me rethink what it means for students to participate in the classroom in a meaningful way, and it has raised my awareness about the fact that not all participation is equal. In addition, this experience has changed the way I teach and how I interact with students in ways that are proving to be long-lasting. Even now, more than a year after the original study took place, I’m finding that I’m still changing my teaching practices, because I now have a new lens through which I can observe my practice. This was afforded by having access to tangible data about the nature of student participation in my classroom.

## 5. *Facilitators' Reflections*

Confronting classroom inequities is difficult, requiring vulnerability and humility. We found Wilson's thoughtful reflections and openness to be integral to the progress we made as a group. More than mere observers, we found ourselves in true partnership with Wilson to improve student experiences.

As Wilson explained in a post-PLC interview, having data around what was happening in his classroom pushed him to see beyond his blind spots and to adjust his practice towards achieving his equity goals. This was evident in the episode Wilson described earlier, where a group of Filipinos participated way less than he expected. Because many instructors do try to build community in their classroom through informal conversations at the start and end of class, it is easy to conflate non-content-related student engagement with content-related engagement when thinking about who participates – EQUIP is regularly used to capture student engagement around course materials in particular. In recognizing this lack of student participation from this group of students, Wilson found ways to engage them to ensure that they were getting what they needed out of the learning environment. This included more individualized supports after class and more interaction with these students through group work.

In the biweekly PLC meetings, the instructors discussed emergent topics stemming from the data analytics they had received such as, effective ways of engaging students in lecture-based courses, and whether rotating groups or seating would impact participation of marginalized students (e.g., students identifying as women, Black, Filipino, Latinx, LGBTQIA+, etc.). It was important for the conversations to be driven by the instructors, instead of the facilitators, since instructors have valuable lived experiences which they can capitalize on in providing each other constructive feedback. As we (Stone-Johnstone and Reinholz) analyzed data from each cycle of data collection, we could see in Wilson's teaching practices how he consciously transformed aspects of his pedagogy to confront areas for growth that emerged from previous cycles. One big transformation that we witnessed was the huge contrast from the beginning of the PLC where a handful of students in Wilson's class would regularly dominate whole discussions to the end where he was able to engage 27 out 36 students in course content throughout the one hour and fifty-minute class. There is only so much one can do in a single semester, but for Wilson, elevating various diverse voices in his classroom was a goal that he was able to accomplish. And as he articulated, he is continuing to think about ways to create equitable spaces in his classroom, as well as ways of engaging students more deeply in course content.

## 6. *Discussion*

We write this chapter as mathematics educators who have deep commitments to equity. At the same time, we recognize that teaching equitably is an ongoing, and lifelong journey. Very often teaching is framed as an individual act, something that we do by ourselves when we walk into our classroom. In reality, our classrooms are reflections of the inequity in society writ large, and they are also a reflection of our own biases and ways of being. Rather than judging individual instructors based on their teaching, we feel it is imperative to work together, in community, to improve equity. Having access to real, meaningful data, and thought partners through which to digest the data is one key ingredient to this. The entire authorship team realizes that we all learned so much about equitable teaching during our time working together, even though only one of us was teaching at that time. There is great potential when we can work together like this, and we would love to see institutional sanction and support for more efforts like this.

# 26 CALCULUS: Crossing the Bridge to Success in STEM

**Elaine A. Terry**
*Saint Joseph's University*

Elaine is a tenured member of the Department of Mathematics at Saint Joseph's University. With extensive experience in teaching precalculus and calculus, she is especially interested in ensuring the success of underrepresented students in STEM. She enjoys living in the Mount Airy section of Philadelphia with her two cats, Tiger and Boom-Boom.

*Abstract: Calculus is often described as the gateway course to Science, Technology, Engineering, and Mathematics (STEM) education. However, it can be the mathematics course that prevents students from obtaining a STEM degree. This is especially the case for students from underrepresented groups. CALCULUS: Crossing the Bridge to Success in STEM is a four-week summer bridge program at Saint Joseph's University that is designed to equip first-year students with the tools necessary for success in a STEM major. Students attend workshop classes in three STEM subjects with the mathematics class serving as the foundational course of the program. Using a workshop style approach, students are introduced to calculus early concepts such as the limit and the derivative. They are encouraged to read, write, and discuss the underlying concepts of a first-year calculus course by working on worksheets entitled Previews to Calculus.*

## 1. *Introduction*

A private Jesuit Catholic university, Saint Joseph's is located in Lower Merion and Philadelphia counties. Founded in 1851 as Saint Joseph's College for men today the university is a coeducational institution with a student population of approximately 8,300 including undergraduate (day and evening), graduate and doctoral degree students. The undergraduate day program is approximately 77% white and 14% from an underrepresented group which includes African-American and Hispanic (non-Black). In any one academic year more than one hundred first-year students declare a STEM (natural sciences, mathematics, computer science) major. Many first-year STEM students experience difficulty handling the challenges of taking a college mathematics course (pre-calculus or calculus) and one or two lab science courses simultaneously. By the end of their freshman year some of them have made the decision to switch to a non-STEM major. Those that remain in STEM have the false belief that low and even failing grades will not prevent them from entering medical, professional or graduate school. While these are issues that may affect all STEM students, it is especially a difficulty for underrepresented students, who represent a small but significant number of STEM majors. Studies suggest that with early intervention underrepresented students can successfully complete an undergraduate STEM program. *CALCULUS: Crossing the Bridge to Success in STEM (CB-CTEM)* was developed in order to equip first-year underrepresented students with the tools necessary for success in a STEM major. CB-STEM, is a four-week summer pre-college non-residential program. It is the first intervention program at the University to address the gap in STEM education specifically for underrepresented (African-American, Hispanic (Non-Black) and First-Generation) students at the pre-freshman level. As a multi-faceted program, the primary goal is to provide incoming underrepresented students with the tools and resources that are beneficial for their first year in college as a STEM major.

To be considered for the program, students must meet the following criteria:

- Be admitted as a full-time student to the University.
- Be a first-time incoming college freshman with a declared STEM major.
- Be classified as a member of an underrepresented group which includes African-American, Hispanic, and First-Generation college students.
- Have completed pre-calculus or calculus in high school.

Students that are accepted into the program are expected to participate in academic and informational workshop classes. There is no cost to students to participate in the program.

The primary objectives of CB-STEM are:

- To increase participants' chances of passing their first college calculus course.
- To introduce students to: lecture format, classroom technology, laboratory class, faculty expectations, and college-level exams
- To help students gain an understanding of good study habits, techniques, and skills.
- To introduce students to STEM faculty.
- To expose students to available academic resources at the University.
- To help students to connect with and build community with other first-year STEM students.

There are three academic workshops that are designed to make students aware of the rigorous requirements of the Saint Joseph's STEM major. The foundation of the program is the mathematics workshop, which uses *Previews to Calculus* worksheets to help teach and reinforce early calculus concepts including the limit and the derivative. More information about the mathematics workshop as well as a brief overview of the biology and chemistry workshops follow. A fourth seminar is designed to give students information that will be beneficial to them socially, academically, and professionally.

## 2. ***Natural Science and Informational Workshops***

The three one-week academic workshops introduce students to the rigorous course expectations of the Saint Joseph's STEM curriculum. The workshops were organized in conjunction with two University professors from the departments of biology and chemistry. The three of us met to discuss the program and agree upon its structure. Each workshop met for sixty-minutes for four consecutive days. We discussed the difficulties that many students have with these courses as first-year STEM students. It was agreed that students would benefit by learning about course expectations, how to take notes, homework, and write lab reports. Students were also given information about studying for tests and test taking strategies. We all cited the lack of these skills as the primary reasons for lack of success as a first-year STEM student.

A fourth workshop is conducted by a University administrator who has experience working with underrepresented student populations. These are fifty-minute workshops that meet at least four times throughout the duration of the program. The aim of the workshop is to empower students with knowledge, skills and self-awareness that is necessary for success as a first year undergraduate STEM student. Topics of discussion are geared towards achieving success in the classroom and the importance of finding opportunities at the university that are relevant to STEM education. In addition, students are given information that will aid in the transition to college by helping them to find opportunities to become immersed in the university community as well as preparation for professional and graduate school and career opportunities.

### 3. *The Mathematics Workshop*

By the first day of the CB-STEM program all of the students have taken the University mathematics placement test. Because they have declared a STEM major, they have either placed into pre-calculus or calculus. In most cases, the students who placed into pre-calculus have not had calculus in high school. Those that place into calculus have a choice of two different University calculus courses based upon their declared major. The CB-STEM mathematics seminar exposes the pre-calculus students to the early concepts and language of calculus including limit, tangent line, derivative, and area below a curve. For the calculus students, the workshop serves as a review of calculus that gives them a deeper understanding of calculus concepts besides the usual drill and practice problems that many of them are accustomed to from high school calculus. The math workshop meets every day for four weeks for one hour fifteen minutes. Lectures are limited to thirty and no more than forty minutes. Students are required to take notes and encouraged to ask questions.

In order to reinforce the mathematical concepts from the brief lectures students are given *Previews to Calculus* worksheets. The worksheets are intended to challenge students to rethink how mathematics problems are solved. The exercises are written so as to reinforce conceptual understanding from the lectures and to stimulate interest in learning calculus in a less routine manner. There are worksheets that cover topics including the infinitesimal, slopes of graphs, rate of change, and area below a graph. The applications in the worksheets illustrate how calculus can be applied to other disciplines including, the natural sciences, business, and economics. Teaching assistants, current undergraduate STEM students, serve as tutors to help students with any problems they may experience while working on the worksheet.

The first worksheet, Preview #1, is presented below. It opens with a brief description of calculus; discussing the meaning of calculus as the study of *change* and how the tools in calculus are used to describe numerically how something is changing at a given instance. To see this the worksheet begins with a problem that students should be able to complete. They are encouraged to work the problem and write a brief explanation of what the value they obtain means. Average rate of change (ARC) is defined and discussed. Two problems are given that encourage students to use an average rate of change formula that is appropriate for the given problem and to discuss the values obtained. Units must be included with the answer to in order to get students to understand the idea of time dependent rate of change.

**Preview #1: What is Calculus**

In short, calculus is the study of change. You might ask: What is meant by change in mathematics? Change of what? The most important concept in calculus involving change is known as the derivative. Derivatives help to compute the rate of change of something; they help to answer the questions: How fast or slow is something changing at a certain instant? What is the rate at which something is increasing or decreasing at a given instant?

In this first preview an example is given that you should be able to work through. Following this example, two questions are asked that will illustrate the need for a mathematical tool that models movement, change.

**Example 1.1**

The number *N* of rabbits in a colony can be modeled by the polynomial function below:

$$N = f(t) = 120t - 0.4t^4 + 1000$$

where *t* is the time in months since observing began.

(a) What is the number *N* of rabbits in the colony when observing began? Briefly explain what your answer means.

(b) What is the number *N* of rabbits in the colony 4 months from when observing began? Does this represent an increase or decrease in the number of rabbits in the colony? Briefly explain what your answer means.

(c) What is the number *N* of rabbits in the colony 5 months from when observing began? Does this represent an increase or decrease in the number of rabbits in the colony? Briefly explain what your answer means.

**Definition**: Recall that for a given function *y=f(t)*, that the average rate of change, ARC, of *f* over an interval $[t_1, t_2]$ can be found using a difference quotient as follows:

$$ARC = \frac{f(t_2) \quad - f(t_1)}{t_2 \quad - t_1}$$

Use the ARC difference quotient for the following problems.

(d) What is the average rate of change of rabbits from month 3 to month 4? Include units with your answer and briefly explain what your answer means.

(e) What is the average rate of change of rabbits from month 5 to month 6? Include units with your answer and briefly explain what your answer means.

In calculus we would ask the following questions as it relates to the increase/decrease of rabbits in the colony:

Q1: What is the RATE OF CHANGE of the number of rabbits increasing/decreasing in the colony AT 3 months?

Q2: What is the RATE OF CHANGE of the number of rabbits increasing/decreasing in the colony AT 5 months?

The two questions at the end of the worksheet are discussed. What do they mean? How do we find the answer? Students are then asked to complete another worksheet that contains tables with the values that will aid in getting a better idea of rate of change at a given instance. They are then encouraged to review the results that they obtained and to use them to try answering Q1 and Q2. Discussion is held about the difference between average rate of change (ARC) and (instantaneous) rate of change. Again, units for instantaneous rate of change are discussed to reinforce that we are not just looking for the *number* of rabbits but the *rate* at which the number of rabbits is changing with respect to a certain time.

Use the rabbit colony function $N = f(t) = 120t - 0.4t^4 + 1000$ to complete the following tables with the appropriate values. Use the tables to help answer Q1 and Q2. Approximate to four decimal places.

**Table 1A**

| t, months | N, # of Rabbits | Time Interval | Average Rate of Change |
|---|---|---|---|
| 3.9 | | [3.9, 4] | |
| 3.99 | | [3.99, 4] | |
| 3.999 | | [3.999, 4] | |
| 3.9999 | | [3.9999,4] | |

**Table 1B**

| t, months | N, # of Rabbits | Time Interval | Average Rate of Change |
|---|---|---|---|
| 4.9 | | [4.9, 5] | |
| 4.99 | | [4.99, 5] | |
| 4.999 | | [4.999, 5] | |
| 4.9999 | | [4.9999, 5] | |

### 4. *2019 CALCULUS: Crossing the Bridge to Calculus COHORT*

While we were unable to recruit the 12-14 students that we proposed, we had four very astute students in the first cohort, one male and three females. Of the four, two had taken calculus in high school with the other two having taken pre-calculus. Three of them enrolled in a calculus course as they were declared STEM majors. The psychology major had the option of following the calculus path or taking a different general education math course. She chose to enroll in *The Whole Truth About Whole Numbers,* which is described as a number theory course for non-math majors.

As the table below indicates, from an academic standpoint, the students had a successful first year. However, there is some concern about the physics major due to the drop in GPA.

| Student | Major | Fall GPA | Fall Math | Sp GPA | Sp Math |
|---|---|---|---|---|---|
| #1 | Computer Science | 3.81 | Pre-Cal<br>B | 3.71 | Calc. I<br>A- |
| #2 | Psychology | 3.80 | Whole Truth<br>A- | 3.90 | NA |
| #3 | Physics | 3.14 | Calc. I<br>B | 2.93 | Calc. II<br>B- |
| #4 | Biology | 3.93 | Calc. I<br>A | 3.88 | NA |

A survey was administered via Survey Monkey. Following are results from some of the questions that were asked.

1. How helpful was the CB-STEM program to you adjusting to the academic requirements of SJU? Results: 1 responded a great deal; 2 responded a lot; 1 responded a little
2. How likely are you to recommend the CB-STEM program to other Freshmen STEM students? Results: 2 responded very likely and 2 responded likely
3. How satisfied are you with the STEM major you have chosen and the department in which your major is housed at SJU? Results: 3 responded satisfied and 1 responded neither satisfied nor dissatisfied
4. How likely are you to recommend SJU to someone seeking to study a STEM subject? Results: 3 responded likely and 1 responded neither likely nor unlikely

There are various reasons for the small number of participants in the program. Including the fact that we are unable to provide housing, students tend to work during the summer and/or may have other academic opportunities that are preferred to our program. Also, the program ran during the week of orientation for first-year students. This is a time when students are on campus for registration and other activities. We plan to be more cognizant of these issues in the future when planning. Two students that had been accepted into the program but chose not to participate confided in me that they regret that they did not participate in the program given the difficulties they were experiencing as biology majors. Both students have changed their major.

## 5. *Conclusion*

The primary aim of CB-STEM is to ensure that underrepresented students have the opportunity to successfully complete their degree in a STEM field at Saint Joseph's University. Beyond helping students feel academically prepared for calculus and a STEM major, CB-STEM also helps students acquire social-emotional skills, identify personal behaviors that may affect their success at the University, and identify and reflect on their goals for the future. Overall CB-STEM is important to creating a sense of belonging, a shared identity and stressing the importance of academic success in STEM.

# 27 Shifting Pre-Calculus from Gatekeeper to Gateway at a Hispanic-Serving Institution

**Elsa Q. Villa**

*The University of Texas at El Paso*

Elsa Q. Villa (evilla@utep.edu) is a research assistant professor in the College of Education where she serves as director of the Center for Education Research & Policy Studies. She has led and co-led numerous federally funded grants to increase student success. Her research interests include STEM identity, learning communities, gender, communities of practice, and transformative learning.

**Lawrence M. Lesser**

*The University of Texas at El Paso*

Lawrence M. Lesser (Lesser@utep.edu) is a Professor in the UTEP Department of Mathematical Sciences. His statistics/mathematics education scholarly activities in engagement and DEI are informed by diverse experiences, including high school teaching, directing UTEP's teaching center, and co-founding the TODOS journal *Teaching for Excellence and Equity in Mathematics*. He has published 120+ papers, been PI of recent NSF grants, and won several state/national awards for his teaching and scholarship.

**Erika Mein**

*The University of Texas at El Paso*

Erika Mein (elmein2@utep.edu) is Associate Dean of Undergraduate Studies and Educator Preparation in the College of Education at the University of Texas at El Paso, and Associate Professor in the Department of Teacher Education. Since starting at UTEP in 2008, she has been involved as PI or co-PI on a series of federal grants focused on Latinx student success, particularly in STEM fields. Her research from these projects has been published in journals such as *Action in Teacher Education*, *Theory into Practice*, and *Journal of Hispanic Higher Education*.

*Abstract: An intervention was funded by the U.S. Department of Education to improve student success in the 5-credit hour precalculus course at The University of Texas at El Paso, a research university at the US-México border with an 80% Latinx student body. The intervention was a collaboration between UTEP's College of Education and the Mathematical Sciences Department in the College of Science. The four-year intervention developed and involved teaching assistants in structured and sustained professional development grounded in the equitable and inclusive pedagogical practice of cooperative learning in order to promote the creation of math-focused learning communities in the precalculus course. Findings indicate the intervention increased student pass rates and their success in the subsequent Calculus 1 course.*

## 1. *Introduction*

Despite more minoritized undergraduates receiving STEM degrees, a systemic gap in degree attainment rates between Whites and students of color remains (NSF, 2019). This is due in part to passing rates in precalculus, the gateway to first-semester calculus required by all STEM majors (Cooper & Woodward, 2011). Various reasons abound for noncompletion of precalculus including a lack of sense of belonging which contributes to engagement and well-being in postsecondary education, especially for minoritized students, such as Latinx students who might be first-generation college students (Hurtado & Carter, 1997; NASEM, 2017; Gopalan & Brady, 2020). Equitable and inclusive college environments contribute to creating a sense of belonging, which has promise of improving retention and persistence, especially in precalculus, and thereby substantively impacting the path to graduation.

Thus, with funding from the U.S. Department of Education under its Minority Science and Engineering Improvement Program, The University of Texas at El Paso (UTEP) implemented an intervention in the 5-credit hour precalculus course MATH 1508 with an aim of improving pass rates and ensuring success for students. The intervention was grounded in the equitable and inclusive pedagogical practice of cooperative learning (Johnson et al., 1998) to promote the creation of learning communities, which are small organically-grown groups serving as academic and social support structures (Astin, 1999; Dagley et al., 2016; Shapiro & Levine, 1999). Such communities are particularly critical for UTEP, as it is a Minority- and Hispanic-Serving Institution (MSI/HSI) with a majority Latinx population (>80%), a population specifically identified as needing more attention by researchers (e.g., Hagman, 2021). Findings indicate the intervention, which was implemented over a four-year period, has had a positive effect on student pass rates and their success in the subsequent Calculus 1 course.

## 2. *Motivation*

The UTEP Department of Mathematical Sciences (hereafter referred to as "math department") has a longstanding history of revisiting and reforming its precalculus course. For example, over 25 years ago, the department divided both semesters of its precalculus sequence into four or five modules each to ensure a more fine-grained sequence of proficiency attainment. While this intervention improved success in the subsequent calculus course, it did not sustain for complex reasons.

The math department then merged the two-semester precalculus course sequence into a one-semester course so STEM majors could start their required calculus sequence one semester sooner. Contemporaneously, the department received a grant from the U.S. Department of Education to establish Peer Led Team Learning (PLTL), a nationally established model that has proven effective (e.g., Hockings et al., 2008). PLTL supplements a course with formalized study groups that meet for two hours each week and incorporate active learning approaches led by near peers who facilitate group discussions around key concepts. What has sustained after that grant is the 5-credit hour precalculus course's structure of four hours of lecture and two hours of workshop each week with class size limited to 45 students.

Then, several faculty members from the UTEP College of Education approached the math department chair in 2015 and expressed interest in submitting a U.S. Department of Education grant proposal to strengthen the 2-hour workshop, which had reverted to being run by students with only minimal guidance. For example, those teaching assistants (TAs) were instructed to let students form groups to work on weekly problem sets, but without intentional structure to the groups, students tended to work alone rather than with unfamiliar classmates, and the precalculus failure rate then was 35-45%.

These particular education faculty, including author Villa (2013), had been successful in an intervention with several cohorts of first-semester engineering students between 2013 and 2015. The intervention embedded team-based projects using cooperative learning principles with one team that was fully immersed in precalculus concepts. Students were placed into teams of six with each student assigned a different particular precalculus mathematical concept (such as functions, exponents and logarithms, coordinates of the unit circle, and trigonometric functions and identities) to learn well enough to teach the rest of the team. To prepare for the teaching portion, they would meet weekly with someone else from a different team who was assigned the same concept using available resources such as Khan Academy videos and precalculus books. Then at various times during the semester, they would teach that particular concept; these times aligned with departmental precalculus exams in order for students to prepare for the exam. Of course, these cohorts had leaders who (as upper-division STEM majors) were close peers and were available to answer questions and guide the students as needed.

Grant external evaluators found that 85.2% (N = 107) of treatment students passed calculus compared to 75.8% (N = 872) of non-treatment students ($p = 0.013$). The treatment students progressed toward degree completion faster than non-treatment students. Thus, these education faculty felt they could scale their effort by designing a similar intervention with the precalculus workshops. In 2015, they submitted a grant proposal to the U.S. Department of Education under the Minority Science and Engineering Improvement Program wherein they proposed team-based approaches similar to the engineering grant where students generated self-directed learning of mathematics.

## 3. *Description and Population*

UTEP is a majority Hispanic institution located on the U.S. border with México and is the only Mexican American-majority university classified by Carnegie as a very high research activity institution (Carnegie Classification of Institutions of Higher Education, 2018). Of UTEP's roughly 25,000 undergraduates, 81% identify as Hispanic-American, 32% are from families with annual household income below $20,000, and about 50% are first-generation college students.

### 3.1. *Student Participants*

Each semester, over 1,000 students enroll in the precalculus course with approximately 85% identifying as Hispanic, of which 45% are female. With the failure rate (i.e., DFW rate) between 35% and 45% prior to the intervention, a significant fraction of students would either repeat the course or drop out.

While we cannot conclusively state the reasons for that attrition, we can try addressing it with evidence-based high-impact practices in classroom settings that create more inclusive environments by cultivating equity and a sense of belonging among students. The core of our intervention is the use of cooperative learning, which NSSE (National Survey of Student Engagement) data show to be a particularly important "high-impact practice" for underserved students such as Hispanics and first-generation students, actively fostering their learning and improved completion rates (Finley & McNair, 2013).

## 4. *Teaching Assistants*

Prior to the intervention, teaching assistants were mainly graduate students who were international students along with a few undergraduate STEM majors. These undergraduate students were of Hispanic origin. Over the three years of our intervention implementation, the number of graduate students serving as precalculus TAs has decreased, especially those who are international; and the number of undergraduate students increased based on recommendations from the precalculus coordinator. From its upper-division undergraduates, the math department leadership agreed to fill TA positions primarily with undergraduate students who happened to be Hispanic given that over 85% of the undergraduate population is Hispanic. By choosing TAs from the undergraduate population instead of from the graduate population, the TAs were more likely to be closer in age, ethnicity, and bilingual fluency in Spanish and English to the students they were teaching. Indeed, it was observed that there was translanguaging occurring in the workshop among the precalculus students and the TAs. Translanguaging refers to the use of one's full language system, which can include multiple languages and varieties (García & Wei, 2014). In the precalculus workshops, students and TAs seamlessly moved between Spanish and English to discuss problem solutions and explanations.

## 5. *Intervention*

The education faculty designing the intervention drew upon a sociocultural theory of learning that views learning as a social and cultural phenomenon rather than as a phenomenon taking place only within a learner's mind. Grounded in the work of Russian psychologists and philosophers Lev Vygotsky, Alexander Luria, and Alexie Leontiev, this sociocultural approach focuses on social interaction, authentic engagement, and the learning environment. Its characteristics include authentic activities to promote learning, collaborative construction of knowledge, coaching and scaffolding at critical moments, and reflection to encourage metacognition (Cole & Engeström, 1993; Lave, 1991; Vygotsky, 1978).

The intervention aimed to provide the TAs with professional development prior to each semester and to comprise cooperative learning, inclusion, and an overview of sociocultural theory. The goal was to develop learning communities through the practice of cooperative learning to promote social interaction and authentic engagement, and thus foster student learning.

### 5.1. *Cooperative Learning*

Cooperative learning is the instructional, intentional use of small groups of three to four students who work together to maximize their own and each other's learning, creating student connectedness as problems are worked through (Villa et al., 2013). A meta-analysis found cooperative learning increases student achievement, self-esteem, and the formation of positive relationships (Johnson & Johnson, 1989). Moreover, the use of small groups creates a space for "extended social interaction…where students with different mathematics backgrounds [a]re able to learn from one another" (Mein & Esquinca, 2017, p. 287). As Mein and Esquinca (2017) note, the use of small groups provides a safe learning environment where "students' knowledge, identities and practices are valued" (p. 288), thus creating an inclusive environment.

In the first session each semester, students are asked to work as individuals through a set of review problems (i.e., intermediate and advanced algebra) for 15-20 minutes. They are then asked to turn to the nearest person and spend about 20 minutes reaching consensus on solutions. They are then asked to turn to the pair next to them and reach consensus as a group of four on a set of answers. They are told that TAs will randomly pick a group and an individual from the group to share a particular problem's solution. Notice that this announcement ensures that teams work together with everyone engaged, no one left behind, and no one wanting to let the team down since no one knows who will be the one selected.

The following week, the students are put into groupings different from the prior week's groupings to allow them to work with new students. A problem is identified for all student groups to work on; and, similar to the first day of class, they are told that TAs will randomly pick a group and a member of that group to present a problem solution to the class. These problems are solved on a digital tablet using an app (e.g., Notability) that records the explanation and the steps. This recording is presented to the class and then archived online in the learning management system Blackboard for future reference. The course shell then becomes a repository of problem solutions for students. Occasionally, the TA will pause the recording before the solution is completed and ask the class to turn to their teammates and discuss whether they agree with the steps taken.

### 5.2. *Teaching Assistant Professional Development*

In order to create precalculus workshops that foster a sense of belonging, it is important to prepare the TAs to oversee the workshops and create the cooperative groups. Thus, prior to each semester, a 3-hour training is held to provide professional development (PD). The PD agenda includes inclusion and cooperative learning and their theoretical underpinnings. Cooperative learning is taught in part by modeling it, with the various PD elements including small group work that allowed learning and experiencing the key elements of cooperative learning: positive interdependence, individual accountability, group processing (reflection), social skills, and face-to-face promotive interaction (Johnson & Johnson, 1989). For example, the technique noted in the previous section of the TA randomly choosing a team's student to explain a problem allows TAs to implement the elements of individual accountability and positive interdependence since it is in each student's interest to be prepared and to help ensure all teammates are prepared as well.

To internalize these key cooperative learning elements, the TAs are put into groups of three with each TA within a group responsible for learning one or two of the key elements. Each member is paired with another group's member who has the same assigned element. Together, they read and learn about the element(s), knowing that they need to return to their original group of three and "teach" the element(s) to the others. This method of dividing up material is called a cooperative learning "jigsaw" and gives each team member something to have responsibility for teaching to others as well as things to have to learn from others. After the activity, the lead PD instructor (one of the education faculty) has the group reflect on how the jigsaw's structure incorporated the key elements of cooperative learning. When the agenda for the first day of precalculus workshop is presented, the TAs select the activities that apply the key elements of cooperative learning, such as randomly selecting a student to explain a problem, which is an example of incorporating individual accountability. During the semester, TAs, education faculty, and the precalculus coordinator meet weekly to reflect on what worked well and what needed improvement. Together, they prepare for the following week with each TA taking turns on the workshop design.

Several of the undergraduate TAs worked one summer with one of the education faculty, whose specialty is sociocultural foundations in education, to create three short (under three minutes) videos: describing the workshop's underlying sociocultural theories of learning, how to build learning communities through teamwork, and activity footage that shows good and bad examples of workshop situations, such as groups that contain too many students or groups with someone who is disengaged. For example, the video might show one of the students too far away from the group and not interacting. The question posed is: What would you do? Such examples enable each TA to be able to more readily identify these situations and intervene in real time. Or they may see a group that is too large to ensure that everyone participates. TAs are reminded that in such a situation, they should intervene to break up the group and assign one of the group members to a group of three or create a new group of three. See Figure 1 for a screenshot of one of the videos created for TAs and see that video (with students' voices and avatars) at https://www.youtube.com/watch?v=XJRAyW5QHdc.

Last, the PD features a session on inclusion. To begin this session, TAs rate their overall high school experience from 1 to 10 where 1 is the worst experience they have had and 10 is the best. Then they are asked to rate their first semester in college. They have a whole group discussion on these experiences and generalize that the high ratings are where they felt a strong sense of belonging (e.g., being in the marching band or on an athletic/academic team, etc.) and a weak sense of belonging for the low ratings. This provides an experience for the TAs to draw on when they are introduced to the theory of communities of practice, which the PD faculty lead would describe as such: learners are invited into a community and, in the beginning, are situated on the margins, moving closer to the center as they are embraced by the experts and gain knowledge. If they are not embraced by those more expert, they might remain on the margin and never fully gain entry into the community (Lave, 1991; Wenger, 1998). Or, if they are not invited into the community, they assume an outsider position. The PD asks the TAs to revisit their high school and college ratings and reflect on whether they were insiders, on the margin, or outside the community. These reflections are shared in small groups and then discussed as a whole group.

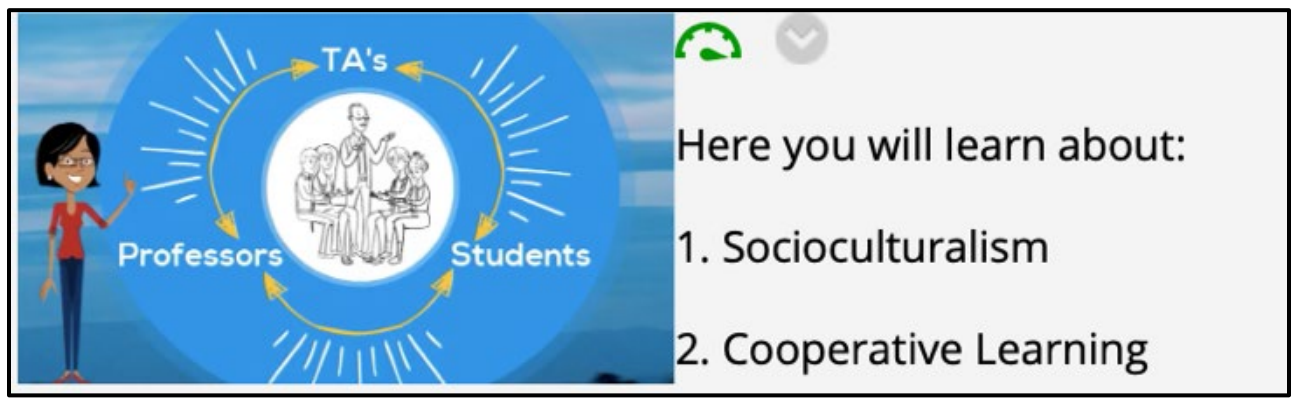

*Figure 1. Opening frame of TA video*

To promote inclusionary practices among the TAs, the PD concludes with TAs sharing what they should be doing as precalculus workshop leaders and these ideas are placed on a flip chart sheet that is posted on the wall. Each TA is given three sticker dots to "vote" on which ones deserve highest priority. Once this activity is completed, a discussion on the items with the most dots is held to reinforce how the TAs can facilitate a community of practice where all participants can feel as insiders, reinforcing the notion that students will feel a sense of belonging. This is important to all students, regardless of gender, sexual orientation, race, or ethnicity.

## 6. *Evidence of Impact*

The grant external evaluators analyzed data for five academic years starting with 2014-15, and Table 1 shows a general decline in failure rates (i.e., students withdrawing or earning D or F) from 38% to around 25%. While this decrease is good, it is important to ensure overall progress does not leave underrepresented groups behind. As an aside, evaluators found significant evidence that women outperformed men for precalculus, Calculus 1, and Calculus 2 ($p$ = 0.000, 0.001, and 0.008, respectively). While we cannot conclusively state why this occurred, that trend is not inconsistent with results of cooperative learning studies that go back a few decades.

*Table 1*
*Precalculus failure rates.*

| Academic Year | # enrolled | DFW rate = #DFW / #total enrolled |
|---|---|---|
| 2014-2015 | 1022 | 38.06% |
| 2015-2016 | 1094 | 23.86% |
| 2016-2017 | 1112 | 28.60% |
| 2017-2018 | 1167 | 25.62% |
| 2018-2019 | 1022 | 24.36% |

Another success indicator is that recently, including the grant window, *Diverse: Issues in Higher Education* has ranked UTEP in the top ten nationally for Mathematics and Statistics degrees conferred to Hispanic students at both the undergraduate and master's levels. For example, in the 2017-18 year, Hispanic students comprised 91% of its undergraduates graduating with these degrees, a percentage exceeding the Hispanic percentage of UTEP's student body.

## 7. *Barriers and Recommendations*

COVID-19 has had a major impact on the precalculus workshop. In the middle of the spring 2020 semester, the workshop (as with almost all instruction nationwide) moved from face-to-face to online. Because the TAs had already met the students for half the semester, it was a straightforward transition to maintain connections and synchronous interactions via the learning management system Blackboard.

However, the fall 2020 semester was challenging. The TA professional development was delivered virtually using Zoom with breakout rooms for small group activities. Having TAs experience their first semester as a TA in an all-virtual environment was a particular challenge, likely due to their lack of experience in these workshops when it was face-to-face and/or their own experiences in small group work that is deliberate and intentional with individual accountability built into its design.

When face-to-face sections resumed in fall 2021 (i.e., with the availability of COVID vaccines), the precalculus workshop meeting time remained online (but still synchronous) because it was recognized that some students might not feel safe working closely in groups if the workshop were conducted in person without assurance of any needed social distancing, masking, vaccination, or ventilation. Students in online asynchronous sections, however, often found it hard to attend a weekly synchronous workshop meeting. Therefore, beginning with the fall 2022 semester, those students were allowed to fulfill the workshop part of their grade by completing the Personal Study Plan in the WebAssign platform at their own pace. The math department and its chair continue to support the intervention past the end of the grant. In particular, (except for asynchronous sections) the cooperative learning strategy continues for the precalculus course. TAs continue having PD training (online) before each semester starts and then having weekly meetings throughout the semester. We continue exploring other online tools and approaches to facilitate the workshops as well as the PD training in an even more meaningful way. Faculty interested in implementing such an intervention at their institutions may want to find faculty in their education department willing to consult on how best to create a supportive workshop environment where all students develop a sense of belonging.

## 8. *Acknowledgements*

This material is based upon work supported by the U.S. Department of Education (USDoE) under #P120A160056. Any opinions, findings, and conclusions or recommendations expressed in this material are those of the authors and do not necessarily reflect the views of USDoE. Much appreciated are the UTEP doctoral students and (College of Education and Mathematical Sciences Department) faculty integral to the intervention's deployment and associated research.

# 28 Nivelando el Campo Educativo / Leveling the Educational Field

**Cristina Villalobos**

*The University of Texas Rio Grande Valley*

Dr. Cristina Villalobos holds the Myles and Sylvia Aaronson Endowed Professorship in the School of Mathematical and Statistical Sciences at the University of Texas Rio Grande Valley. Her research is in optimization, optimal control, and STEM education. She is Founding Director of the Center of Excellence in STEM Education, and Associate Dean for Strategic Initiatives and Institutional Effectiveness in the College of Sciences. Her work in STEM curricula, leadership and student mentorship especially for underrepresented groups can be summarized with the 2020 Presidential Award for Excellence in Science, Mathematics, and Engineering Mentoring (PAESMEM).

**Timothy Huber**

*The University of Texas Rio Grande Valley*

Dr. Timothy Huber is a Professor and Director of the School of Mathematical and Statistical Sciences at the University of Texas Rio Grande Valley. His research areas include analytic number theory, combinatorics, and special functions. He leads several grants and other initiatives intended to level the playing field for students in STEM. He is a recipient of the University of Texas System Regents' Outstanding Teaching Award.

**Jose Ponce**

*The University of Texas Rio Grande Valley*

Mr. José Ponce is a Lecturer in the School of Mathematical and Statistical Sciences at the University of Texas Rio Grande Valley whose research area is in bilingual instruction. Mr. Ponce is the coordinator of the Dual Language Certificate program in Mathematics and course coordinator of Elementary Statistics. Mr. Ponce has been recognized as a Community Engaged Scholar for his participation in workshops and discussions with community members and leaders of the Rio Grande Valley and is a member of PUEDE (Parents United for Excellent Dual Education), a group that advocates for dual language programs at local K-12 schools in the Rio Grande Valley.

**Josef Sifuentes**

*The University of Texas Rio Grande Valley*

Dr. Josef Sifuentes is an Associate Professor in the School of Mathematical and Statistical Sciences at the University of Texas Rio Grande Valley whose research area is in iterative methods in applied mathematics and numerical linear algebra. Dr. Sifuentes has been a co-PI on two NSF grants aimed at widening the pipeline to graduate programs to underrepresented minorities. He is also the co-director of the UTRGV LSAMP Summer Research Academy.

*Abstract: The University of Texas Rio Grande Valley is located in south Texas on the border with Mexico. UTRGV's goals include broadening student success, building students' self-confidence, and cultivating a sense of pride in the linguistic and cultural heritage of the Rio Grande Valley. This paper describes several programs that seek to level the educational field by building inclusive and supportive academic environments. The narrative provides a template for institutions developing similar inclusive academic initiatives. We describe the role of Spanish and English in the Valley and offer an overview of introductory dual language courses in the department and their impact on students. These courses are part of a Dual Language Certificate currently under development. The peer groups formed in these courses are powerful. They flip the script of English language dominance in the classroom to an environment where both languages are valued. The focus on cultivating equity and inclusivity through multiple modes of*

*communication and interaction continues in peer groups in the Calculus sequence. These groups are part of an effort to build an academic community to support learning and encourage collaborative problem solving. Promoting equity through multi-section course coordination is also discussed. Finally, broader components "para nivelar el campo educativo" continue beyond core courses into students' professional development through the Center of Excellence in STEM Education. Its promotion of pathways through the university and to broader opportunities is designed to increase the number of Latino students attaining STEM degrees and leadership positions across the Nation.*

## 1. *Introduction*

Given the diversity of the Rio Grande Valley and its proximity to Mexico, many students develop as dual-language learners – initially learning to communicate in Spanish at home and in English at school. However, students often do not develop both languages simultaneously in school, since English takes precedence over Spanish, and all academic work and teaching is conducted in English. Several efforts over the years have led some school districts to develop dual-language programs to empower and build self-confidence and cultural awareness in students. Bilingual education has been an integral part of the University of Texas Rio Grande Valley (UTRGV) since its creation in fall 2015. The university's commitment has included founding an institute that promotes activities and programs to create a bilingual, bicultural, and biliterate institution, referred to as the B3 Institute. As one of the largest Hispanic-Serving Institutions in the nation with 90% of its 30,000+ student body being Latino, primarily Mexican-American, UTRGV is building on its cultural and linguistic assets to promote diversity, equity, and inclusion through dual-language instruction, peer community building through active learning in Calculus classes, and holistic activities that increase the number of underrepresented groups in STEM.

## 2. *Dual Language Certificate in Mathematics and Student Experience*

In collaboration with the B3 Institute and as part of UTRGV's strategic plan, the Mathematics department began developing a dual language certificate (DLC) in Spring 2018 with the goal of providing an innovative and accessible educational environment that promotes success for its native Spanish-speaking students and those seeking to demonstrate proficiency in dual-language settings. The DLC proposal consists of 9 credit hours of math courses offered in Spanish and 6 credit hours of Spanish language proficiency. Since fall 2018, sections of College Algebra and Elementary Statistics have been offered in Spanish. Additional courses that will be offered in Spanish include Calculus 1, Calculus 2, Calculus 3, and Differential Equations.

Most language programs lack intersection with courses that provide content knowledge. The planned DLC program targets students who want to develop language proficiency directly in mathematics content areas. The certificate gives students credit for multilingual skills and leverages existing strengths. This is important since historically, educators have viewed Spanish skills as deficits. In the context of many years of persistent ethnic and racial discrimination, the addition of the certificate sends a message to those who grew up in Spanish-speaking households that being bilingual is a valuable skill.

In the early 21st Century, college-level dual language programs are still somewhat novel. It is important that such programs adapt organically to the needs of students and remain inclusive. Although initial DLC courses at UTRGV were offered exclusively in Spanish, they have evolved into courses that are truly bilingual and inclusive. Instruction is primarily in Spanish, but student responses and teacher feedback are multilingual. The goal is to meet students where they are at linguistically. The approach allows those who are not fluent in Spanish to increase their mathematics language skills through participation in a Spanish immersive environment.

Student recruitment for the bilingual courses is led by a DLC committee of faculty and includes on-campus advertising and online. Since the majority of students have a great experience in bilingual courses, a significant amount of advertisement occurs through word of mouth. Some students have expressed apprehension about the dual challenge of learning mathematics in Spanish, even if they spoke Spanish well. Students had concerns that the Spanish modality may compound existing challenges in the course. Student evaluations and pass rates comparable to other sections indicate that most students found the bilingual courses provided additional support. One student said "... This would be my 3rd [time] taking College Algebra and the way the lessons were structured helped me a lot." Some students reported feeling more comfortable in the bilingual setting and attributed their success to the course.

In the classroom, students are given the liberty to express themselves in any of the two languages or bilingually. Often, students become aware that they have skills necessary to understand academic material presented in a different language, skills they did not realize they had or could master. For example, students commented that "... even though the class was in Spanish, [the instructor] still made it very clear what we needed to learn. I was able to grasp all the material he explained...". Although instruction is in Spanish, instructors do not restrict students to use a particular language. Any means of communication is acceptable if it is effective and enhances learning. Students are encouraged to explore the Spanish language as a tool to increase the effectiveness of content communication. Some students communicate mostly in English. These students are not native speakers in Spanish, but since they were exposed to Spanish while growing up, they are able to understand it and write it. The instructor gives supportive feedback on the Spanish component of the class by correcting grammar or rewording sentences, for example. Feedback on grammar is constructive and generally not attached to students' numeric grade. The result is that students feel more comfortable with their level of Spanish and their ability to use it in an academic setting.

All assignments in the course are provided in Spanish; see Figure 1 for an example. This gives students the opportunity to learn new vocabulary and read Spanish in context. Students are placed in groups of three to five on alternate lecture days to solve problems together. The interaction is not only a rich exchange of ideas for mathematical problem solving, but also a diverse combination of English and Spanish. While doing group work, students interact and seek to integrate Spanish in their conversations. Subsequently, these same students are seen interacting outside the group setting in Spanish. A bond is developed between students that persists throughout their academic career. The experience in the course enhances students' communication abilities and reinforces the importance of their Spanish skills. The class provides a space where bilingualism is the norm. Students finish the course with a first-hand understanding that people from diverse backgrounds with different language abilities can come together and solve problems.

Tarea #3 – Álgebra - MATH 1314.81X
Nombre:____________________

Fecha de Entrega: Viernes, 21 de Febrero del 2020. **Deben mostrar su trabajo completo en todos los problemas, no solo las respuestas.**

1. Verifica que las funciones $f(x) = 3x + 4$ y $g(x) = \frac{1}{3}(x - 4)$ son inversas.

2. Encuentra la inversa de $f(x) = \frac{3x+4}{2x-3}$.

3. Supongamos que un fabricante de secadoras de ropa que trabaja con gas ha encontrado que, cuando el precio unitario es $p$ dólares, el ingreso $R$ (en dólares) es

$$R(p) = -4p^2 + 4000p$$

¿Cuál es el precio unitario $p$ que debe cobrarse para maximizar el ingreso? ¿Cuál es el ingreso máximo?

4. Grafica $f(x) = 4x + 2$ y su inversa en el plano de abajo. (No olvides la recta $y = x$).

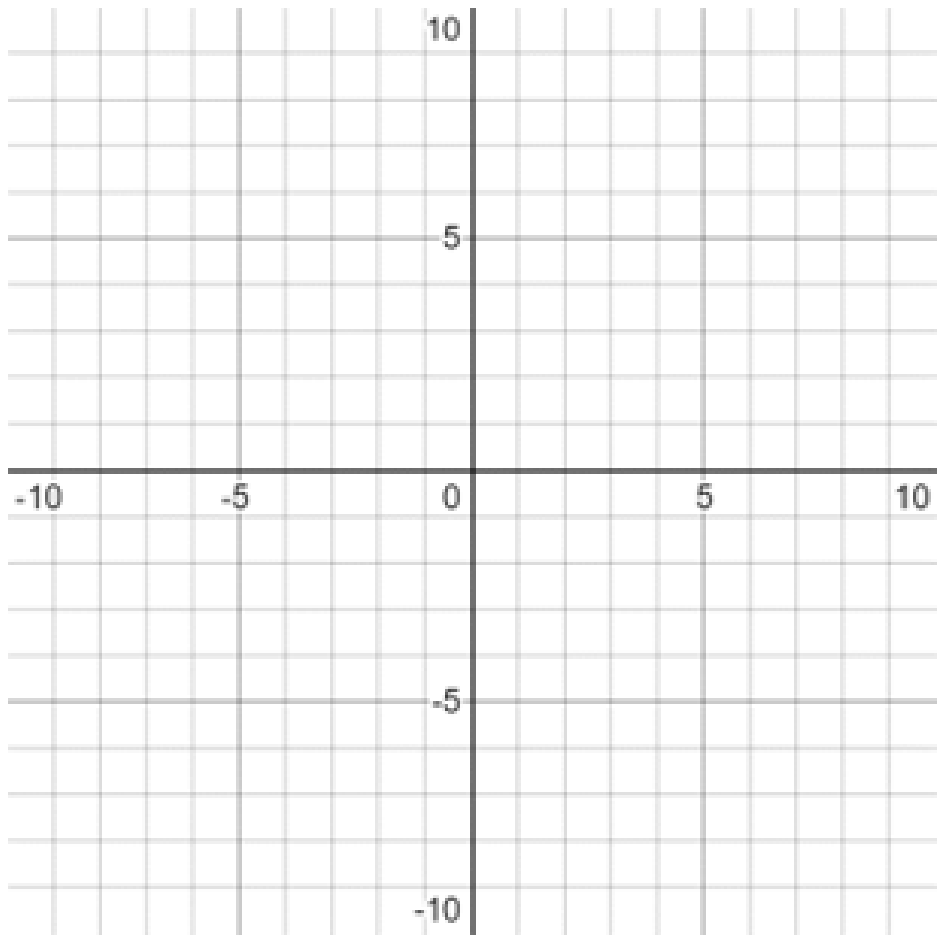

*Figure 1: Example of a College Algebra assignment written in Spanish*

### 3. *Continuing Dual Language Settings in Calculus 1 and Calculus 2*

To broaden student success, since fall 2016 the department has implemented collaborative problem-solving sessions in subsets of Calculus 1 and Calculus 2 classes. In these sessions, students discuss mathematics and solve problems together, often in dual languages. The idea behind these sessions is to incorporate active learning engagement that prioritizes peer-to-peer interactions to increase student success. The peer-to-peer interactions occur both within the group and with undergraduate Learning Assistants (LAs) who are available to help discuss the problem sets in multiple languages. Literature points to the benefits derived from engaging students in active learning sessions and providing a supportive, equitable, and inclusive academic environment (Theobald, 2020).

During the sessions, students work on problem sets that are more open-ended and conceptual than typical homework sets. Often, the links between concepts and applications are emphasized. In one set, for example, the construction of the derivative of the inverse tangent function is scaffolded following the formulation of the derivative of sine and cosine in class. Subsequently students are asked to compute the velocity of a particle whose position is described by the inverse tangent function. In working out the derivation in small groups, each student has the opportunity to actively engage in the lecture material and ask questions amongst their peers in whatever language they are most comfortable. At an HSI institution on the Texas-Mexico border, this language is frequently Spanish, or a combination of Spanish and English, reflecting the languages of their homes.

During one group session on the inverse tangent and velocity, a particular student could not contain his excitement when he realized that the calculation of the derivative of the inverse tangent, which required applying the chain rule, could then be applied to compute the velocity of the particle. He announced his discovery to the entire class in English. Then he explained it to his group in Spanish, his native language. His excitement, along with his eureka moment, spread to his group of peers. This is an example of the benefits described in the literature of small group learning environments. Students who are comfortable in the use of their language tend to feel less threatened about what others think about their language skills and are more likely to engage in classroom discussions that build positive STEM identities.

The LAs put peer-to-peer interactions at the forefront. While not every instructor can join discussions in Spanish and English, many of the peer LAs can. The LAs are undergraduate peer mentors who were successful in Calculus and have been hired to facilitate the group collaborative problem-solving sessions. They typically reflect the cultural make-up of the classes they serve. In addition, by often speaking a language many students feel more at ease in, they serve as a reminder that success in Calculus is not reserved only for those who look like their instructor. Having an instructor or LA who is supportive of bilingual interactions allows students to develop a deeper trust in those individuals.

Designing community-building collaborative problem-solving assignments requires a substantial investment of time. This can be alleviated by centralizing the effort, with input from instructors utilizing the materials. Course coordinators and participating instructors have increased coordination in core mathematics classes. Common group activities are designed to fit into each instructor's lecture schedule. Such coordination promotes equity and encourages instructors to further coordinate common assessments, textbooks, and syllabi, providing greater consistency among all sections of these courses.

### 4. *Continuing to Level the Educational Field Beyond Core Math Classes*

Another example of UTRGV's inclusion and diversity in leveling the educational field is through the work of the Center of Excellence in STEM Education which creates pathways from K-12 to the university and to postgraduate programs. The overall aim is to increase the numbers of Latino students attaining STEM degrees. The Center was established through a federal grant and was one of three Minority Institutions to receive the award. After the grant period, the Center continued to provide activities to the community, K-12 schools, and the UTRGV student population in an effort to bring awareness to STEM fields and retain students in STEM programs. As UTRGV serves a large Latino population in the Rio Grande Valley, it is our responsibility to prepare students for STEM careers and thus contribute to increasing the numbers of underrepresented minorities achieving STEM degrees. At the K-12 level, the Center provides hands-on STEM activities for local K-12 schools visiting the college campus. Activities are developed to have a “wow” factor that engages students and increases their interest in STEM. Student staff, who are generally Latino, develop and deliver the activities and thus serve as role models for the K-12 students. The activities are conducted in English, but student staff also guide K-12 students in Spanish whenever needed, especially when they notice K-12 students speaking primarily in Spanish. The idea is to make K-12 students comfortable, encourage participation in the activity, and acknowledge that Spanish is welcomed in STEM activities.

The Center also serves as a hub for college student academic, career, and professional development with goals that include increasing the pipeline of underrepresented minorities earning STEM undergraduate and graduate degrees and diversifying the STEM workforce. To achieve these goals, workshops are provided throughout the academic year on topics such as applying to graduate school, fellowships, external summer research programs, and tips for designing posters and research presentations at conferences. For example, panels of UTRGV students discuss their experience obtaining an internship and inform their classmates about the application process. Student panelists provide their perspective on the perseverance that is needed to apply and how to approach recruiters. Similarly, the Center organizes diverse panels of scientists to discuss their research career tracks in the three career sectors of academia, government, and industry. Their personal stories demonstrate to students that there is a place in research fields for them. Panelists are from diverse backgrounds, and perhaps more importantly, from underrepresented groups in STEM, such as Latinos and females, allowing students to identify with the panelists.

For three years, the Center has partnered with the UTRGV College Assistance Migrant Program (CAMP), which assists first year freshmen who are migratory or seasonal farmworkers (or children of such workers); historically these students are Latino and first-generation college students. To increase CAMP retention in STEM degree programs and to assist students to learn about STEM careers, the students attend 10 hours of STEM professional development workshops provided by the Center. The Center also provides students with experiential learning opportunities by pairing them with faculty research mentors who have established successful research programs with students. The objective is to introduce CAMP students to a research environment where they are guided by peers from the research group and where they become aware of STEM applications and discoveries. The program has had its challenges placing first-year students in a research environment, but the faculty research mentors have been extraordinarily welcoming to students in providing them with academic guidance and introducing them to their research groups. The hope is that these students will be involved in research activities and as a result will be retained in STEM fields.

## 5. *Conclusion*

UTRGV's Math department has championed inclusivity and equity through a Dual Language Certificate initiative and a linguistically inclusive active learning model of collaborative problem-solving sessions in lower-level mathematics classes. UTRGV's academic efforts to increase recruitment and retention among Latinos in STEM have been augmented by the Center of Excellence in STEM Education. These initiatives contribute to the development of a welcoming community of Latino scholars and increase the pipeline of students attaining STEM degrees at the undergraduate and graduate levels. These initiatives result in a strong group of leaders across the nation who can provide a voice on Latino issues. All three initiatives address the need to "nivelar el campo educativo" or "level the educational field" for Latino students and prepare them for success in their mathematics courses, their undergraduate STEM career, and beyond.

# 29 First-generation students' experience of entering a STEM major

**Amber Wagnon**

*Stephen F. Austin State University*

Dr. Amber Wagnon is an assistant professor of education studies at Stephen F. Austin State University, specializing in student empowerment and student writing as a mechanism for change. Her work also focuses on utilizing technology to engage learners beyond the classroom.

**Keith Hubbard**

*Stephen F. Austin State University*

Dr. Keith Hubbard is a professor of mathematics at Stephen F. Austin State University, focused primarily on studying and fostering STEM student success. He has served as principal or coprincipal investigator on six National Science Foundation projects focused on STEM student or prospective STEM teacher training and support. He has served as a GenJacks instructor and a GenJacks mentor.

*Abstract: In 2014, Stephen F. Austin State University initiated an admission program, called GenJacks, to support and facilitate success of their first-generation college students. GenJacks utilizes a cohort system designed to mentor and support students through three courses taken during their first semester of college. For STEM majors, this typically included College Algebra. We utilize embedded case studies of two focal students to explore the ways in which the GenJacks program contributes to the success and retention of first-generation students. Specifically, we highlight program methods which assisted one student to persist in an Engineering Physics major, but failed to assist another student in her STEM aspiration. We will identify the pedagogical practices employed by instructors which relate to the success of the focal students. This is important as one of the obstacles that many first-generation students face is guidance and assistance once enrolled, as they are encountering territory that is unfamiliar to both themselves and their family (Ward et al., 2012). Our chapter will explore ways in which guidance and mentoring can contribute to identity formation and the retention of first-generation students. Highlighting both the successes and obstacles of these students is vital as first-generation students are a group that remains inadequately understood. These findings can assist higher education institutions and policy makers in formulating plans and creating programs which increase the success rate of first-generation students.*

## 1. *Introduction*

Situated a few hours between Dallas and Houston, surrounded by giant pine trees, is the campus of Stephen F. Austin State University (SFA). In the spring of 2020, 5,798 undergraduate students enrolled at the university were classified as first-generation college students; defined as a student whose parents did not graduate from a four-year institution. This is a large percentage of the entire undergraduate population, making them a critical element for policy and program considerations. The institutional make-up of first-generation students mirrors national trends – 56% of undergraduates nationally are first-generation students (Center for First Generation Student Success, 2019).

Retaining these students in STEM programs is particularly critical, as recent research has revealed that fewer than 50% of students who declare a STEM major persisted to graduation (Higher Education Research Institute, 2010; National Academy of Sciences, 2005). Unfortunately, there are even larger disparities for first-generation students to persist and successfully complete STEM courses at institutions of higher education (McCarron & Inkelas, 2006). Studies have noted that some of the most profound attrition occurs during first year introductory courses (Dyrberg, & Holmegaard, 2019). In response to this information, a new program was created at Stephen F. Austin State University.

## 2. *GenJacks*

In 2014, Stephen F. Austin State University initiated an admission program designed to support and facilitate success of first-generation college students. The program, GenJacks, was designed to support first-generation students, particularly those from underrepresented groups, by providing them an academic community. Initially, GenJacks leadership sought a mathematics instructor to teach a terminal GenJacks mathematics course. However, as the GenJacks program grew, the need for serving first-generation STEM students, all of whom would need to successfully progress through calculus, became a priority. Prior to the creation and implementation of GenJacks, there were no resources specifically designed for first generation students of any major. The program utilizes a cohort system designed to mentor and support students.

GenJacks seeks to support students socially with workshops, community building events, and service opportunities. The program also supports students academically through enrollment in three linked courses taken during their first semester of college. Both the courses and the faculty are carefully chosen to maximize student support. For STEM majors, one of the GenJacks courses was College Algebra.

Many of the case studies in this collection center on the programs and policies. Such reflections have great value, but one of the key attributes of a case study is acknowledging the individual's experiences within a larger context. To that end, we have chosen to highlight two individual students, David and Irene, within the GenJacks program which serves a wide variety of students and should not be homogenized. We believe this is particularly fitting in this manuscript, as it seeks to honor the diversity of human experience and learning.

David and Irene were members of the 2018 cohort, which included roughly 100 students. Several courses were open only to students within the GenJacks program. Among them was a College Algebra course designed to prepare predominantly STEM majors for calculus. Thirteen GenJacks chose to take that course. (They have taken to calling themselves 'GenJacks.') This group came into college with slightly stronger mean SAT scores (1114 vs. 1083) than all first-year students, and stronger mean ACT scores (23.5 vs. 21.1). However, these strong scores did not uniformly translate into success in the college classroom. In fact, the pass rate across all sections of College Algebra that semester was 61% while the GenJack students had a 46% pass rate. Similarly, GenJacks had a first semester mean GPA of 2.57 while all first-year students had an initial semester GPA of 2.83. This gap persisted, with only 54% of the GenJacks persisting a full year. By comparison, first-year students in general persisted at a 66% rate. These differences in success highlight the continued need for programs like GenJacks in which first generation students are supported in their transition to higher education.

In the summer after graduating high school, David, an Engineering Physics major, and Irene, a Physics major, both enrolled in the GenJacks program through Stephen F. Austin. Both David and Irene are first-generation, Hispanic students. When asked to describe any additional demographics that they felt were important to their identity David explained that he was raised by a single mother who was disabled, and Irene described herself as being an only child. Despite similar backgrounds and majors, only David has persisted through the STEM program to graduation in the spring of 2021. David and Irene depict two different “middle outcomes” within the GenJacks cohort. David finished his first semester with a 3.36 GPA, somewhat above average, while Irene finished with a 2.29 GPA, somewhat below average.

Many circumstances contributed to these outcomes, including the climate and the connections made within the university and the GenJacks program. We constructed a list of potential interviewees that had experienced GenJacks College Algebra roughly three years ago, believing this would be an appropriate amount of time for them to holistically gauge their experience without having forgotten too many of the specifics. The semester GPAs and major changes of these students were examined to find representatives of different experiences and outcomes in STEM persistence. David and Irene were selected for this study as representatives of two outcomes of particular significance. They were both contacted personally and interviewed by phone. After interview analysis, they were asked follow-up questions via Zoom.

### 3. *GenJacks Participants*

As a junior high school student David never considered going to college an option; in fact, he revealed that his plan was to quit high school and “just work and try to help out my family.” But a junior high school teacher approached him about a college readiness program, asking him to join because many of his teachers saw “potential” in him. This event, this vote of confidence coming from an authority figure began to transform David’s self-image. This program also worked directly with the GenJacks staff at SFA.

On the other hand, Irene always “intended on attending a four year” institution but she was “hesitant as to how financially probable” it would be for her family. However, through direct assistance from her high school counselor she was awarded a local scholarship to attend Stephen F. Austin State University. Both David and Irene’s decision to attend SFA was influenced by the connections between the faculty and staff at their respective schools and SFA. It is important to remember that research has long claimed that one way to recruit and retain first-generation students is through partnerships between colleges and K-12 schools (Neislet, 1992; Ward et al., 2012). Despite the shared pathways to SFA, the way in which Irene and David interacted within the GenJacks program were different.

### 4. *Summer Leadership Academy*

Each student enrolled in the program and attended the Summer Leadership Academy (SLA). This event is a weeklong summer orientation for those enrolled in the GenJacks program, occurring on campus before the other students have returned for the fall semester. This orientation provides an opportunity for the GenJack students to meet their cohort members, to connect with peer mentors who have previously completed the program, and to work with the faculty and staff who operate the program.

While both Irene and David attended the SLA, David was more engaged in the in the weeklong program. Irene only attended the first day and then sporadically during the rest of the week. Irene explained that the first day “was very educational” and explained the differences between “things like subsidized and unsubsidized loans” and “a lot of stuff that I had not known prior.” This knowledge is important, as we know that many students of college educated parents have “help navigating, transitioning to college, and succeeding once there” but when institutions fail to meet these needs of first-generation students they often struggle (Horowitz, 2019, p. 11). But outside of this information, Irene admitted that she felt disconnected from the rest of her cohort. When she would periodically attend events during the summer orientation, she could see how connected her peers were, revealing that when she “showed up to the events” she “could see that the bond wasn’t there for me.”

Essentially, what was missing for Irene was a sense of belonging. A sense of belonging, when applied to a collegiate setting, is “students’ perceived social support on campus, a feeling or sensation of connectedness, and the experience of mattering or feeling cared about, accepted, respected, valued by, and important to the campus community” (Strayhorn, 2012, p.4). This connectedness is important. It was important to David, who reported feeling a sense of isolation during his initial time on campus. He revealed that “when I first got here, I was automatically feeling like people were staring at me and like they know that I'm not supposed to be here” and that he felt others had grown up with a “different standard of living.”

This feeling of ‘other’ complicates the often-difficult transition to college that many students face. However, unlike Irene, David reported an intense connection to his peers and the GenJacks staff. David explained that when he first began the SLA he was hesitant to get involved. Instead during the first two days he often slept in, opting out of events. He explains that “I didn't really come in with the mentality that I'm going to lean on them. So, I was trying to like to strengthen myself, and like become independent. But like, towards the end of it I started really making friends.”

These bonds were initiated and sustained by both his GenJacks group leader, and a peer who had already successfully completed two years of a STEM program. David referred to this peer as playing the role of “big brother,” sharing how this GenJack’s openness about his own family struggles encouraged Daniel that he could be successful despite his past.

These experiences, combined with the fact that David persisted toward graduation with a STEM degree, reiterate the importance of a sense of belonging. When students “feel securely connected with others in an environment or a group” then they “tend to have higher self-concepts, greater confidence in their academic skills, and rate their college experience as satisfactory” (Strayhorn, 2012, p.91). In addition to social programs, the designated courses for GenJacks were important in creating this inclusivity. To increase academic success and retention, higher education institutions must work diligently to create a climate that fosters a sense of belonging for first-generation students.

**5. *Designated Courses***

Such a climate also contributes to a student's academic progress. Through the GenJacks program, both David and Irene were enrolled in courses with their GenJacks peers. Enrolling students in courses as a group can assist students in navigating "the academic environment and reconciling the gap between expectations and realities" in the college classroom (Ward et al., 2012, p. 62). This was the case for David who explained that due to the close relationships he had formed over the summer they "could like rely on each other and work with each other when we got stuck." He attributed much of his initial college success to the strong relationships created by having the "same professor together" and always "studying together." He explained, that for him, these relationships felt like "a circle of like people you can depend on and rely on."

However, Irene revealed that she felt she "was like, the odd man out when it came to taking classes with everyone. I kind of just like did the classes by myself like I would any other class just because I hadn't built those relationships." But she explained that she did belong to a GenJacks group text messaging app and she found that "they were always asking [if anyone needed help] in that GroupMe. So, if I was ever confused, I could at least like go into the GroupMe and ask them there, and they would respond." This interaction enabled Irene to benefit from the group knowledge, though in a less critical way, as she viewed herself as being on the fringe of the group.

Faculty can utilize classroom interactions to further ensure a sense of belonging. David explained his college algebra professor "made the strongest effort to remember our names." He explained that the desire to learn each students' name seemed genuine and made him feel that his professor was "really interested in who I am, that I'm not just a student." While this may seem trivial, students, particularly first-generation students, view these sincere actions to mean "that they mattered to someone on campus," thus creating a sense of both acceptance and belonging (Strayhorn, 2012, p. 81). This feeling of mattering is critical for STEM programs to retain students, especially those who may feel "marginalized, unprepared for, or 'out of place' in STEM fields" (Strayhorn, 2012). Irene, enrolled in the same college algebra course, also praised the professor for creating a welcoming environment that she described as "energetic and welcoming." When initially created, leadership sought out instructors with a history of actively engaging students in the learning environment. As the program has continued, a new position was created: GenJacks mentors. These mentors are faculty or staff willing to invest in one to three students in a mentoring capacity from across campus. As other institutions seek to improve support for underserved student groups, it is important to intentionally foster venues for student focused faculty and staff to meet, exchange ideas, and grow in intentionality about serving students.

In addition to the social elements, the classroom instruction they encountered and the ability to engage in self-regulated learning were both key factors of success. Both Irene and David were required to create a study plan prior to each course assessment. Irene explained that "he'd have you like write out how your plan on studying for the test throughout the week before the actual day of the test." She explained that these study plans provided an "idea or concept of the amount of time" she needed to study. Through these study plans Irene and David were learning to engage in self-regulated learning. Self-regulated learning is defined as the "specific processes, strategies, or responses by students to improve their academic achievement" (Zimmerman & Schunk, 2001, p.5). The first steps in self-regulated learning are to understand how to set learning goals and then plan how to effectively achieve these goals (Nilson, 2013). Both of these steps were required when the students engaged in creating their study plans. We know that many first-generation students arrive at college without adequate preparation for the academic rigor they will experience at the university level (Atherton, 2014; Tobias, 1994; Ward et al., 2012). However, when faculty can embed opportunities for students to practice self-regulated learning within the course the opportunity for success for first-generation students can increase (Antonelli et al., 2020) While self-regulated learning was not a designed learning feature of the program it was an important element discovered when program reviewed occurred. This discovery also highlights the importance of engaging in intentional program reviews regularly.

Despite Irene's praise, she did not successfully pass her college algebra course. She interpreted this not as a reflection of the professor, but a reflection of her abilities, describing herself as "truly atrocious at math and science." Irene's academic struggles could be attributed to her self-concept, a person's evaluation of their abilities based on experiences and interpretations (Schunk & Pajares, 2005). Researchers have found that a person's belief in their academic abilities impacts their academic performance. This impact is more profound for first-generation college students (DeFreitas & Rinn, 2013; Hellman & Harbeck, 1997; Vuong et al., 2010). Perhaps some of Irene's academic struggles were a result of her diminished self-concept regarding mathematics.

Irene went on to explain that in high school she was enrolled in science courses "that related to the human body" and was successful in those. However, she was not enrolled in such courses that related to her scientific interests and perceived skill set, in her first collegiate academic year. This is an important consideration, as we know that first-generation students are less likely to choose a STEM major. Some institutions, like Southern New Hampshire University, are working to combat this by creating STEM majors that appeal to first-generation students, making them aware of a myriad of careers and graduate-school options. Perhaps such an option would have contributed to Irene's success and retention in the STEM program (D'Agostino, 2013).

## 6. *Conclusion*

Highlighting both the successes and obstacles of David and Irene is vital as first-generation students are a group "that remains largely unnoticed and poorly understood despite all of the research on students that has emerged in the past decades" (Ward et al., 2012, p. xiii). Unfortunately, limitations impeded some areas of program evaluation. The overall program was more driven by concern for students than concern for conducting research, so there was limited documentation of real-time student dispositions, student obstacles, and the effects of interventions. As other institutions seek to implement similar programs, it is important to consider appointing a math faculty member to serve as a program liaison. It is also imperative to set goals for assisting students and compiling data to inform research findings, both elements which are being added to the GenJacks program.

Such policies and partnerships can work to ensure programs are best meeting the needs of students. As we know, higher education can enable first-generation students to gain upward social mobility making it imperative that we examine the experiences of these students and develop policies and programs reflective of their needs.

## 30 Inclusive Calculus: Active, Applied, and Inquiry-based Learning with Linked Workshop and First-Year Seminar

**John Zobitz**

*Augsburg University - Department of Mathematics*

John Zobitz is a Professor of Mathematics at Augsburg University and is Augsburg's Data Science Program Director. His research areas include mathematical biology and environmental data science, specifically focused on integration of nearly-continuous ecological data with mathematical models. He is a member of the Ecological Forecasting Initiative, which is a grassroots organization aimed at developing real-time ecological forecasts.

**Pavel Bělík**

*Augsburg University - Department of Mathematics*

Pavel Bělík is Chair of the Department of Mathematics, Statistics, and Computer Science and a Professor of Mathematics and Data Science at Augsburg University. His research interests include the mathematical modeling of phenomena in nature and in engineering, and the computational aspects and the numerical analysis associated with the modeling process. In his teaching, he enjoys exploring ways to incorporate active and technology-supported learning into his classes.

**Suzanne Dorée**

*Augsburg University - Department of Mathematics*

Suzanne Dorée is Chair of the Department of Mathematics, Statistics, and Computer Science and a Professor of Mathematics at Augsburg University. She has been involved in the national conversation about the undergraduate mathematics curriculum including through the Mathematical Association of America, the Common Vision project, Transforming Post-Secondary Education in Mathematics, and the SUMMIT-P project. Dr. Dorée enjoys teaching mathematics at all levels using active and inquiry- based learning.

**Rebekah Dupont**

*Augsburg University - Department of Mathematics*

Rebekah Dupont is Director of STEM Programs and an Associate Professor of Mathematics at Augsburg University. She leads the collaborative S-STEM PRISM project and Augsburg's involvement in the LSAMP North Star STEM Alliance. She is interested in strengthening partnerships across organizations to increase the quality and equity of high impact practices in STEM.

**Jody Sorensen**

*Augsburg University - Department of Mathematics*

Jody Sorensen is a Professor of Mathematics at Augsburg University. She is interested in dynamical systems, history of mathematics, and curricular development. Dr. Sorensen won the Pólya Award for an article in the MAA's College Mathematics journal. As part of the NSF-funded SUMMIT-P project, she helped create active and applied activities for every day of Calculus I and II.

*Abstract: Augsburg University is a small, private university in Minneapolis, MN. Over the past decade we have intentionally substantially increased the diversity of our student body, especially students from minoritized and marginalized populations, first-generation students, and students with disabilities. For example, 56% of our students identify as people of color in 2020-21-20 compared to 34% in 2016-17. We took a comprehensive approach of rethinking calculus to create a student-ready program that welcomes and values all our students. In this*

*paper we discuss the rationale, implementation, and impact of five key elements of our program. 1) We partnered with science and Economics faculty, through the NSF-funded SUMMIT-P Project, to refocus calculus on authentic, applied, and contextualized problems relevant to students' lives and future studies. 2) We redesigned each class period to include small group, inquiry-based exploration and collaborative practice work at the boards; these activities now constitute the majority of class time. 3) In addition to standard tutoring, we developed a weekly Calculus Workshop where students work collaboratively on challenging calculus problems. The workshops are open to all, with an emphasis on building community among students from groups traditionally underrepresented in STEM. 4) We connected a section of first year seminar to Calculus I to provide students an opportunity to connect with alumni working in quantitative fields, learn about local STEM labs and companies, and to develop college navigation and learning skills. 5) We built robust transfer pathways for students who start at a community college (or other university), including securing scholarships for transfer students in STEM with support from the NSF-funded AugSTEM program. The collective effect of these changes – to the course content, pedagogy, and student support – is a lively, engaging, and welcoming calculus classroom.*

## 1. *Introduction: Motivation*

Uri Treisman, founder of the Emerging Scholars Program and director of the Charles A. Dana Center at University of Texas Austin has long highlighted the role mathematics plays in equity. As quoted in a recent *Education Week* article,

> *'Math[ematics] is even more important to upward mobility now than it was 20 or 30 years ago, because ... it's seen as related to your general ability to solve problems quickly,' Treisman said, adding that as a result, 'there's general anxiety and panic about equity issues for anything new, even though the current [calculus] pathway is a burial ground for students of color.'*

We teach at Augsburg University in Minneapolis, MN. Originally established as a Lutheran Seminary for Norwegian immigrants, our programs and students have changed dramatically over the past 150 years. In the past decade Augsburg has intentionally substantially increased the diversity of our student body. In fall 2019 53% of our undergraduates were students of color (1141 students) and 37% were first-generation (807 students) (Augsburg University, 2020). During this time the population of students who identify as (a) Black or African American or (b) Hispanic or Latino grew 10% and 16% respectively each year. In addition, around one third of our new students each fall are transfer students, making Augsburg University a popular transfer destination for students at local community colleges. See Figures 1a and 1b which illustrate how rapidly we made this change as an institution.

As elsewhere, our Calculus I students at Augsburg arrive with aspirations for their future and dreams of making an impact on the world. They are often nervous about the course because of its reputation for being difficult – especially so if they needed to complete Precalculus or other mathematics courses first. Societal, economic, and racial systemic barriers amplify the challenges of the course itself. As well-intended faculty members who care about our students' successes, we are disheartened to realize the extent to which we built barriers to student learning that had nothing to do with calculus itself. We imagined our classrooms as a place that supported everyone's learning and where all students felt comfortable (and safe) as they struggled to learn. We were already using active learning and tried to use relevant examples in class. We hoped that made us different, but even seemingly innocuous practices like allowing students to call out answers to questions we posed, frequently using examples from Physics, or expecting students to be quiet to learn reinforced the historically white, male culture of Mathematics. We discounted the impact of racial and financial stress on student learning. The department wanted to face the problem head-on, and embraced changes in our calculus courses that would open them up to students from all backgrounds.

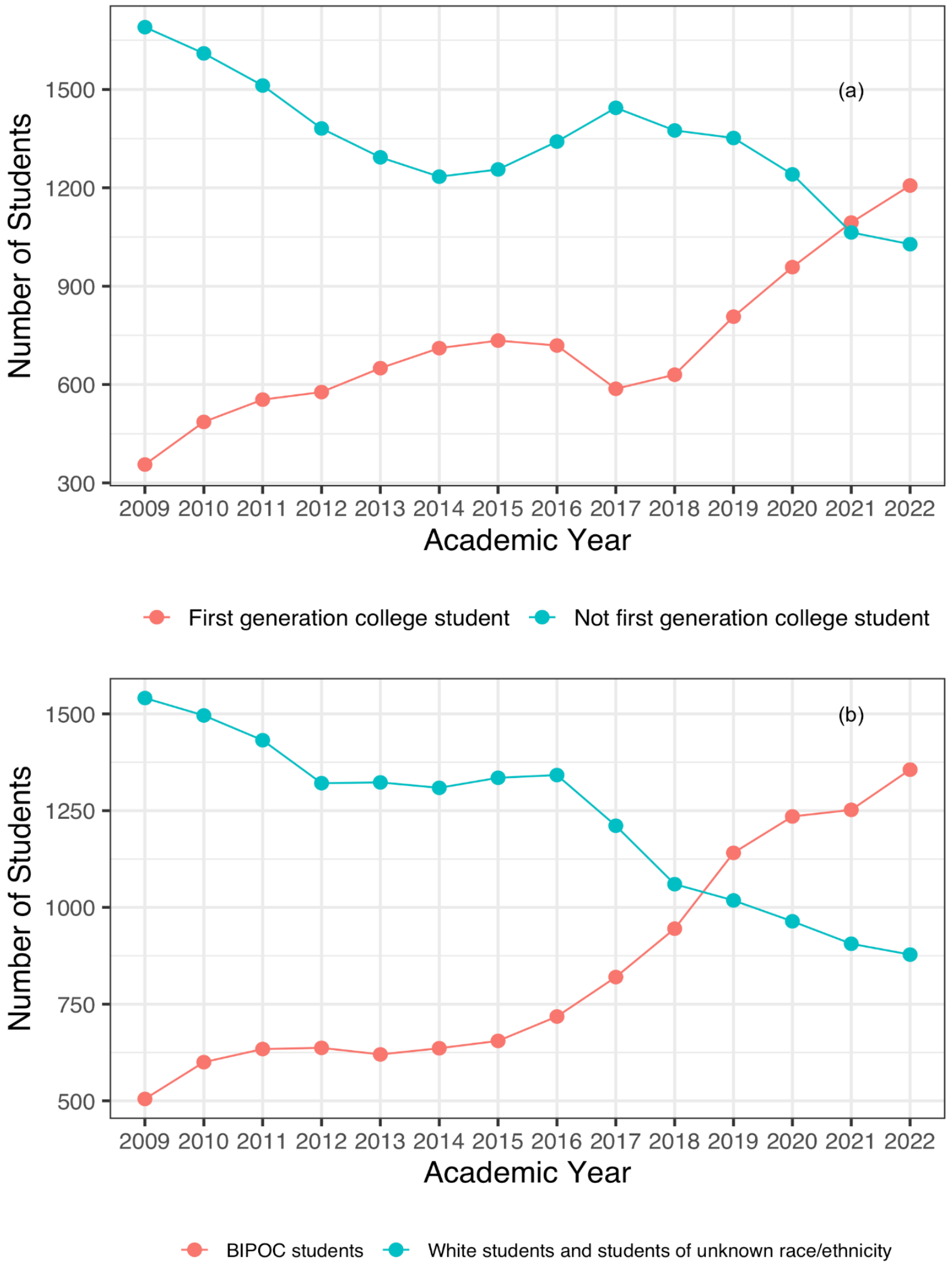

Source: Augsburg University Factbook

*Figure 1: Panel (a): Enrollment trends of traditional undergraduate Augsburg University students, comparing between first-generation college students and not first-generation college students. Panel (b): Enrollment trends of traditional undergraduate Augsburg University students comparing between BIPOC students and white students or students of unknown race or ethnicity.*

## 2. *A Multi-faceted Approach*

*Description & Population targeted:* We each began approaching our calculus sequence from different angles – active learning structures, rethinking the core topics, technology and labs, data-driven relevant examples, classroom practices attuned to different backgrounds and levels of preparation, and outside of classroom supports. Over several years of sharing and collaborating with faculty members across the department, these efforts coalesced.

We targeted our intervention to focus on students of color, first-generation college students, and transfer students. We define students of color as the aggregate of multiple ethnicity categories: American Indian or Alaska Native; Asian; Black or African American; Hispanic or Latino; Native Hawaiian or other Pacific Islander; or two or more races, as reported from the Integrated Postsecondary Education Data System (IPEDS) database.

***Rationale:*** We are experienced mathematics faculty members accustomed to improving the curriculum. Change was needed; we took chances and wanted to be creative. In addition, we hold scholarship of teaching and learning in high regard, so curricular work is part of, rather than tangential to, our work as faculty- scholars. Many of us were engaged in on-going professional development around diversity, inclusion, and equity. Our group included the Director of STEM programs (Dupont) whose job centers on improving student success in STEM, particularly for students of color and first-generation students.

We also had the context to make the changes. Augsburg University has a long-standing commitment to social justice and equity, which is a lens through which we view curricular change. Our rapidly changing demographic has been a center of campus-wide discussions and calls to action for a more inclusive campus (Augsburg University, 2020).

The urgency of our changing context encouraged us to take a critical look at calculus with attention to the course's role as the entry gate to STEM majors where the representation of students of color anecdotally did not match the general student population at Augsburg. Our goal was, and remains, to create student-ready calculus courses that welcome and value our students.

In this paper we discuss the rationale, implementation, and impact of five key facets of our renovation of calculus: Curriculum, Pedagogy, Linked Workshop, Linked First-Year Seminar, and Transfer Path- ways/Scholarships.

## 3. *Five Facets*

We implemented five key changes to calculus that we believe increased inclusion and equity in *Calculus I* and in the Mathematics major.

### 3.1. *Collaborating with the Partner Disciplines to Rethink Calculus Content: The SUMMIT- P Project*

The NSF-sponsored Curriculum Foundations Project (CFP) brought together partner discipline faculty to discern what mathematics their majors needed to know. Their findings are summarized in Ganter & Barker (2004) and Ganter & Haver (2011). Augsburg is a member of subsequent collaborative NSF-funded research project SUMMIT-P: A National Consortium for Synergistic Undergraduate Mathematics via Multi-institutional Interdisciplinary Teaching Partnerships (NSF Award #1625557; https://www.summit-p.com/; SUMMIT-P, 2020) that studies how mathematics and partner discipline faculty members can collaborate to improve mathematics courses in the first two years based on CFP findings.

At Augsburg, our team of three mathematicians, one chemist, and one economist are renovating our calculus curriculum. To start, we held listening sessions with additional faculty members from Biology, Physics, Chemistry, Environmental Studies, Economics, Finance, and Business to refine CFP recommendations for our local courses. We also reviewed textbooks from classes requiring calculus.

Based on this work we narrowed the content of our calculus courses to essential concepts, skills, and habits of mind. We wanted to be sure that everything we were teaching students would be useful and relevant. We also brought in a wider array of applications of calculus, especially from Biology and Economics, and committed to teaching calculus concepts in context rather than in the abstract. We switched textbooks to align with these changes (Briggs et al., 2018) and wrote daily exploratory activities and weekly labs which often start with data (not equations). These changes made the course more accessible to students from varied backgrounds, and increased interest in the course.

**Results**: Paring down the course to essential concepts and skills considerably reduced the need for rote algebraic manipulation which gave students with less algebraic facility greater access to the big ideas. This refocusing is especially important for students without extensive practice in rote algebraic manipulation from high school.

Using examples that are more familiar to students (e.g., Biology vs. Physics) made the courses more relevant and useful to students. In the 2018-2019 school year students completed an attitudinal survey. Of the 32 respondents, two-thirds reported they had already applied what they had learned in *Calculus I* in their non-mathematics courses. In addition, 80% of students believe that the things they learned in *Calculus I* will be useful to them after college. A similar percentage of these students stated what they learned in the course will be applicable to their futures.

### 3.2. *Consistent use of Active and Inquiry-based Learning: 20 + 20 + 20 = 70*

During class meetings, we spend the majority of time on active and inquiry-based learning (IBL) pedagogies that have positive effects on student learning (Conference Board of the Mathematical Sciences, 2016). Our calculus course meets four times per week: three 70-minute class periods and one 100-minute laboratory. We split the class into three 20-minute segments (with 10 minutes to use as needed). Students begin with an exploratory activity and work with classmates of their choosing. Next, the instructor recaps key ideas. Lastly, students work on practice problems with randomly assigned partners at the white boards.

We deliberately made the in-class activities carefully scaffolded - starting in a familiar place and then increasing in difficulty and sophistication. Students of all backgrounds can make a start and engage with the material - everyone has thoughts to share and questions to ask.

During these times, students support one another in learning and the instructor circulates to help. Studies by Freeman et al. (2014) and Laursen et al. (2014) have shown active learning has disproportionately positive effects for students from minoritized groups (Conference Board of the Mathematical Sciences, 2016). In Laursen et al. (2014) studying the impact of IBL by gender, women performed as well in IBL or non-IBL courses, but women's reported level of mastery was higher in IBL courses. According to Laursen, "IBL methods do not 'fix' women but fix an inequitable course."

**Results:** Using active and inquiry based learning pedagogies **every** student is talking in class, with each other and with the instructor. We see increased ownership ("When we did part (c),. . . " vs. "When you did Example 3, . . . ") and sense of belonging, particularly in the formation of outside-of-class study partners. Instructors are able to see and hear what students are doing and able to quickly identify where students might be struggling and support them.

### 3.3. *Developing Skills in Community: The Calculus Workshop*

For some students of color and first-generation students finding a home with a community of learners took additional support. This was especially true early on when in many cases there would only be one student of color in a given calculus section. In 2009 we (Dupont) began a Calculus Workshop where students work together to strengthen their mathematical and study skills. The workshop is open to all students, with targeted recruitment to first-year students of color or first-generation students.

The workshop sets high expectations for academic excellence. The goal of the collaborative activities is to develop student confidence in their ability to tackle difficult mathematical problems. The workshop problems contain a mix of previewing material to be presented in future classes, reviewing key concepts, or tackling challenging multi-step calculus problems.

A peer mentor assists the professor with student questions. The peer mentor also serves as a role model, helps students navigate the "hidden curriculum," and connect students with leadership, academic, and research opportunities such as study abroad, Fulbright awards, the McNair Scholars program.

**Results:** After the first year, workshop participants achieved a 0.6 higher grade point average (GPA) than the *Calculus I - II* students not in the workshop, despite having entered college with an overall lower average ACT score and GPA (Augsburg Now, 2011). Examining more recent data from 2016-2020 in the aggregate (25 different *Calculus I* sections) a higher percentage of workshop participants earned a grade of B- or better in *Calculus I* compared to students not in Workshop (Table 1). This pattern is consistent and higher for students of color or first-generation students.

*Table 1*
*Comparison between students who receive a B- or better in Calculus I and participation in Work- shop, also separated between students of color and first-generation students. We report aggregated numbers from 2016-2020.*

| Population | B- or better | Less than B- |
|---|---|---|
| All students in Calculus I | 209 (53%) | 184 (47%) |
| Workshop students in Calculus I | 65 (64%) | 37 (36%) |
| Students of color in Calculus I | 101 (50%) | 103 (50%) |
| Students of color in Workshop & Calculus I | 37 (65%) | 20 (35%) |
| First-generation students in Calculus I | 83 (54%) | 71 (46%) |
| First-generation students in Workshop & Calculus I | 36 (75%) | 12 (25%) |

The results in Table 1 indicate that for the intended Workshop population (students of color and first-generation students), participation in Workshop has a positive association in earning a B- or better in *Calculus I*. In addition, students sometimes enroll in the same section of subsequent courses with classmates from Workshop or continue study partnerships, strengthening community building.

### 3.4. *Envisioning a Career in STEM: The Augsburg Seminar*

Our students want to see a connection between their major and meaningful, financially-stable careers. Surprisingly, many of our students did not know that many of the top-ranked careers start with a major in mathematics, statistics, or computer science (U.S. News, 2020; Indeed, 2020). We wanted to help students of color and first-generation students see those connections to STEM and Business, but especially to Mathematics, Statistics, or Data Science.

We linked a section of Augsburg Seminar, our 1-credit first-year seminar, to *Calculus I* with a focus on career pathways in Mathematics, Statistics, Computer Science, or Data Science. Students register for Augsburg Seminar during summer orientation, and the course is open to all students.

We revised the Augsburg Seminar linked to *Calculus I* so students engage in a series of career planning and visioning exercises; explore profiles of BIPOC, LGBTQIA, and female mathematicians and computer scientists (Mathematically Gifted & Black, 2020; Lathisms, 2020; Project 5050, 2020; Sigma Pi Sigma, 2020; Henrich et al., 2019); and interview recent mathematics and computer science alumni for career advice. We also help students create an academic plan and mentor them through their coursework (often for all four years). The course also addresses inclusion and equity by discussing implicit bias (both in and outside mathematics) and intercultural communication. In addition, the seminar includes sessions on financial aid, study abroad, advising, and mental health to support student success.

**Results:** Students enjoy learning more about careers and, especially, interviewing alumni. Students appreciate that they can explore profiles of a variety of mathematicians, expanding their view that mathematics is an exclusively white, male field (Stinson, 2013) . As a predominantly white faculty, we believe these connections help address the lament that “they can’t be who they can’t see.”

Because the Augsburg Seminar is both major specific and linked to *Calculus I* we hoped the changes to the Seminar were associated with a positive outcome of B- or better in *Calculus I* for Mathematics and Computer Science (MSCS) majors. We examined data from 2016-2020 comparing the grades between students in Augsburg Seminar to students in *Calculus I* in the aggregate as well as based on populations of interest (students of color and first-generation students). We do not report percentages as in Table 1 because the number of students in this population is smaller.

*Table 2*

*Comparison between Mathematics and Computer Science (MSCS) majors who receive a B- or better in Calculus I and participation in Augsburg Seminar, also separated between students of color and first-generation students. We report aggregated numbers from 2016-2020.*

| **Population** | **B- or better** | **Less than B-** |
|---|---|---|
| MSCS majors in *Calculus I* | 83 | 63 |
| MSCS majors and Augsburg Seminar students in *Calculus I* | 23 | 12 |
| MSCS majors and students of color in *Calculus I* | 39 | 36 |
| MSCS majors and students of color in Augsburg Seminar & *Calculus I* | 13 | 5 |
| MSCS majors and first-generation students in *Calculus I* | 39 | 22 |
| MSCS majors and first-generation students in Augsburg Seminar & *Calculus I* | 14 | 2 |

The majority of MSCS majors who participated in Augsburg Seminar received a B- or better in *Calculus I*, and this pattern was consistent - and in some cases disproportionately larger - when controlling for a particular subset of the population. Anecdotally we know that students in the linked Augsburg Seminar have a higher-degree of self-selection to study mathematics than their peers, which may be a confounding variable in this association.

### 3.5. *Supporting Transfer Students: The AugSTEM Program*

Approximately one-third of our incoming class each year consists of transfer students. These students have additional needs, including building a connection to their new classmates and professors and learning to navigate a new institution with different expectations. Some transfer students also carry family and work obligations outside of school, which adds to their pressure to succeed. Many started at a community college or state university because it was more affordable than Augsburg.

We designed the AugSTEM Scholars Program to relieve some of their financial burden and to facilitate an easier transition for transfer students. Augsburg has received two awards under the NSF Scholarships for Science, Technology, Engineering and Mathematics (NSF S-STEM) program, Award #1154096 (2012-2016) and Award #1565060 (2016-2020).

Students in the AugSTEM Scholars Program receive personalized support for their needs, starting with an Individual Development Plan (IDP) as well as the CliftonStrengths assessment [20]. IDPs can take different forms and appear in many workplace and educational settings (CIMER, 2020; My IDP, 2020; Bosch, 2017). We directly connect AugSTEM Scholars with recent alumni about navigating the transition to post-bachelor's opportunities.

To support students and advisors at Community Colleges, one of our staff members facilitated communications between department chairs, faculty members at Augsburg, and area Community Colleges to develop 32 STEM Transfer Guides at eight community colleges. The transfer guides are shared with community college advisors and to students at events targeting STEM transfers.

**Results:** Over the period 2012-2019, the two S-STEM grants supported 111 scholars of whom 92 graduated in STEM and 18 are continuing in STEM (99% retention). Of these scholars, 41% were transfers. Five students received Graduate Fellowships for their academic achievement – two of these students began at community colleges and two were NSF-identified underrepresented minorities. Table 3 shows a comparison of the number of STEM graduates overall and the number who started as transfers from prior to getting the first grant in 2012 to 2018. Three-year averages were used to compare the growth over this time period of all graduates (11%), STEM graduates overall (56%) and STEM graduates who began as transfers (115%).

*Table 3*
*Three-year averages of Bachelor's Graduates by STEM and Transfer Status since AugSTEM Scholars began in 2012. Source: Augsburg Office of Planning and Effectiveness (Day program students only)*

| Year | 2009 | 2010 | 2011 | 2016 | 2017 | 2018 |
|---|---|---|---|---|---|---|
| All Graduates | 356 | 382 | 426 | 456 | 385 | 448 |
| All STEM Graduates | 64 | 43 | 76 | 107 | 87 | 91 |
| Transfer STEM Graduates | 18 | 15 | 16 | 36 | 35 | 34 |

Even with the AugSTEM scholarships of up to $10,000 per year, many students still experience financial need. Mentors help students apply for additional scholarships or find employment with a higher hourly wage. We also established a new donor-funded scholarship specifically for 5th-year STEM transfers.

## 4. *Conclusions*

Our systematic changes addressed all dimensions of how students engage with calculus, inside and outside of class. Perhaps you may be wondering how you can make change at your institution. As a place to start we encourage you and your colleagues to reflect on the following questions:

- What motivates you to change your curriculum?
- What is the department experience in changing curriculum and pedagogy?
- Where is your department in terms of training for diversity, equity, inclusion, and anti-racism?
- What tools/resources/opportunities beyond the classroom could be part of the solution?
- What other partners across campus can help you in this work?

*Discussion from these questions sparked* ideas for changes we implemented at Augsburg University. Small, isolated changes alone might not be enough. We encourage you, too, to take a multi-faceted approach.

## 5. *Acknowledgments*

This work is supported, in part, by grants from the National Science Foundation (Award Nos. 1625142, 1154096, and 1565060). Any opinions, findings, and conclusions or recommendations expressed in this material are those of the author(s) and do not necessarily reflect the views of the National Science Foundation.

## Glossary

| Terminology | Definition or meaning |
|---|---|
| **ABC** | An acronym used to describe the group of students who receive a letter grade of A, B, or C. This is typically a threshold for passing a course in order to take subsequent courses. |
| **BA** | Bachelor of Arts |
| **BIPOC** | Refers to any individual who identifies as Black, Indigenous, or is a person of color. Black and Indigenous individuals are centered in this description to underscore the generational effects of enslavement and genocide on these peoples. |
| **BS** | Bachelor of Science |
| **Corequisite** | A course taken at the same time as another often required course. The corequisite course (e.g., co-calculus) typically supports students regarding pre-requisite and foundational concepts. |
| **DEI** | An acronym that commonly stands for Diversity, Equity, and Inclusion. |
| **DFW** | An acronym used to describe the group of students who receive a letter grade of D, F, or withdraw from a course. This is typically a threshold for not passing a course and needing to retake the course. |
| **Formative assessment** | The goal of formative assessment is to monitor student learning to provide ongoing feedback that can be used by instructors to improve their teaching and by students to improve their learning. Formative assessments are generally *low-stakes* forms of assessment, which means that they contribute to a low (or no) percentage of the course grade. Examples of formative assessments include asking students to: draw a concept map in class to represent their understanding of a topic, submit one or two sentences identifying the main point of a lecture, or turn in a research proposal for early feedback. |
| **Gap-gazing** | Placing primary focus and resources on identifying what a student (or group of students) does not know or cannot do. Gap-gazing is often associated with an examination of outcome data disaggregated by student groups to identify "gaps." |
| **Hispanic-Serving institution (HSI)** | A Hispanic-Serving Institution (HSI) is defined as an institution of higher education that has an enrollment of undergraduate full-time equivalent students that is at least 25 percent Hispanic students. |
| **Historically Black Colleges and** | The Higher Education Act of 1965, as amended, defines an HBCU as: "…any historically black college or university that was established prior to 1964, whose principal mission was, and is, the education of black Americans, and that is accredited by a nationally recognized |
| **Universities (HBCU)** | accrediting agency or association determined by the Secretary [of Education] to be a reliable authority as to the quality of training offered or is, according to such an agency or association, making reasonable progress toward accreditation." |
| **Learning Assistant (LA)** | Learning assistants (Las) are generally described as paid undergraduate students who completed a particular course and performed well enough in the course to help facilitate discussions or activities with the guidance of a faculty member. The duties of a learning assistant can include grading, tutoring, assisting groups, etc. |
| **Minority-Serving Institution (MSI)** | Minority-Serving Institutions (MSIs) are institutions of higher education that serve minority populations, which includes Historically Black Colleges and Universities, Hispanic-Serving Institutions, Tribal Colleges and Universities, Alaska Native Serving Institutions, and Asian American and Pacific Islander Serving Institutions. |
| **Pell-grant** | Pell Grants are a form of federal financial aid which are usually awarded only to undergraduate students who display exceptional financial need. As such they are often used as indicators of family socio-economic status. |
| **Predominately White Institutions (PWI)** | Predominately White Institutions (PWIs) are not an official designation but are often used in literature and publications to indicate a higher educational institution that has an enrollment of undergraduate full-time equivalent students that is at least 50 white students. |
| **Role model** | A role model is a successful exemplar and can impact a person's motivation and achievement. |
| **STEM** | An acronym for Science, Technology, Engineering, and Mathematics. |
| **Summative Assessment** | The goal of summative assessment is to *evaluate student learning* at the end of an instructional unit by comparing it against some standard or benchmark. Summative assessments are often *high stakes*, which means that they individually contribute to a high percentage of the course grade. Examples of summative assessments include: a midterm exam, a final project, a paper, or a senior recital. |
| **Underrepresented and racially minoritized (URM)** | A term used to describe a group of individuals who are underrepresented in STEM and emphasizes the systemic nature in which they are minoritized through racial structures and discourses. |
| **Underrepresented students** | A term used to describe a group of individuals who are underrepresented in STEM, commonly used for racial groups (Black, Hispanic, Indigenous) that are underrepresented but without specially |
| | emphasizing race and may include women, veterans, and dis/abled individuals. |
| **Underserved students** | A term used to describe students who are underserved by the STEM discipline. This terminology places emphasis on the system and structures that are not serving particular students as opposed to common terms like “minority” which puts focus on the individual. |
| **White supremacy** | Refers to a political, economic, and cultural system in which white individuals overwhelmingly control power and material resources, white superiority and entitlement are widespread, and relations of white dominance and non-white subordination are reenacted daily. |